\documentclass[smallextended]{svjour3}

\usepackage{amsmath}
\usepackage[charter]{mathdesign}

\usepackage[scaled=.95,type1]{cabin} 

\DeclareTextFontCommand{\textbfit}{%
  \fontseries\bfdefault 
  \itshape
}

\usepackage[
scr=boondoxo,
]{mathalpha}
           
\DeclareFontFamily{U}{dutchcal}{\skewchar\font=45 }
\DeclareFontShape{U}{dutchcal}{m}{n}{<-> s*[1.0] dutchcal-r}{}
\DeclareFontShape{U}{dutchcal}{b}{n}{<-> s*[1.0] dutchcal-b}{}
\DeclareMathAlphabet{\mathlcal}{U}{dutchcal}{m}{n}
\SetMathAlphabet{\mathlcal}{bold}{U}{dutchcal}{b}{n}

\newcommand{\kos}[1]{\mathscr{#1}}

\newcommand{\mon}[1]{\mathlcal{#1}}
\newcommand{\cat}[1]{\textit{#1}}

\newcommand{\obj}[1]{\operatorname{Obj}{\left(#1\right)}}
\newcommand{\spec}[1]{\textit{Spec}\left({#1}\right)}

\DeclareSymbolFont{cmocal}{OMS}{cmsy}{m}{n}
\SetSymbolFont{cmocal}{bold}{OMS}{cmsy}{b}{n}
\DeclareSymbolFontAlphabet{\mathocal}{cmocal}

\usepackage{xcolor}
\usepackage[shortlabels]{enumitem} 
\usepackage[parfill]{parskip}
\usepackage{shuffle}

\usepackage[all,cmtip,color,matrix,arrow, 2cell]{xy}
\UseAllTwocells

\makeatletter
\def\labelbox#1{%
  \hbox{%
    \setbox\z@=\hbox{$\m@th\labelstyle{\,#1\,}$}%
    \setbox\tw@=\hbox{$\m@th\labelstyle\,$}%
    \dimen@=\ht\z@ \advance\dimen@ by \wd\tw@ \ht\z@=\dimen@
    \dimen@=\dp\z@ \advance\dimen@ by \wd\tw@ \dp\z@=\dimen@
    \box\z@
  }%
}
\makeatother

\usepackage[colorlinks=true]{hyperref}

\usepackage{etoolbox}
\tracingpatches
\makeatletter
\patchcmd{\@citeo}{\hskip0.1em}{\kern0.1em}{}{}
\patchcmd{\@citex}{\hskip0.1em}{\kern0.1em}{}{}
\makeatother

\makeatletter
\newcommand*{\Xbar}[1]{}%
\DeclareRobustCommand*{\Xbar}[1]{%
  \mathpalette\@Xbar{#1}%
}
\newcommand*{\@Xbar}[2]{%
  \sbox0{$#1\mathrm{#2}\m@th$}%
  \sbox2{$#1#2\m@th$}%
  \rlap{%
    \hbox to\wd2{%
      \hfill
      $\overline{%
        \vrule width 0pt height\ht0 %
        \kern\wd0 %
      }$%
    }%
  }%
  \copy2 %
}
\newcommand*{\Frozenbar}[1]{}%
\DeclareRobustCommand*{\Frozenbar}[1]{%
  \mathpalette\@Frozenbar{#1}%
}
\newcommand*{\@Frozenbar}[2]{%
  \sbox0{$#1\mathrm{W}\m@th$}%
  \sbox2{$#1#2\m@th$}%
  \rlap{%
    \hbox to\wd2{%
      \hfill
      $\overline{%
        \vrule width 0pt height\ht0 %
        \kern\wd0 %
      }$%
    }%
  }%
  \copy2 %
}
\makeatother

\usepackage{tocbibind}

\DeclareMathAlphabet{\pzc}{OT1}{pzc}{m}{it}

\def\CC{{\mathocal C}}
\def\CD{{\mathocal D}}

\def\CF{{\mathocal F}}
\def\CG{{\mathocal G}}

\def\CK{{\mathocal K}}
\def\CL{{\mathocal L}}

\def\CR{{\mathocal R}}
\def\CS{{\mathocal S}}
\def\CT{{\mathocal T}}

\def\pr{\prime}
\def\ppr{{\prime\!\prime}}

\usepackage{graphicx}
\makeatletter
\DeclareRobustCommand*\uell{\mathpalette\@uell\relax}
\newcommand*\@uell[2]{
  \setbox0=\hbox{$#1\ell$}
  \setbox1=\hbox{\rotatebox{10}{$#1\ell$}}
  \dimen0=\wd0 \advance\dimen0 by -\wd1 \divide\dimen0 by 2
  \mathord{\lower 0.1ex \hbox{\kern\dimen0\unhbox1\kern\dimen0}}
}

\makeatletter
\def\widebreve{\mathpalette\wide@breve}
\def\wide@breve#1#2{\sbox\z@{$#1#2$}%
     \mathop{\vbox{\m@th\ialign{##\crcr
\kern0.06em\brevefill#1{0.6\wd\z@}\crcr\noalign{\nointerlineskip}%
                    $\hss#1#2\hss$\crcr}}}\limits}
\def\brevefill#1#2{$\m@th\sbox\tw@{$#1($}%
  \hss\resizebox{#2}{\wd\tw@}{\rotatebox[origin=c]{90}{\upshape(}}\hss$}
\makeatletter

\usepackage{tikz}
\usepackage{tikz-cd}
\usetikzlibrary{calc,decorations.pathmorphing,shapes}
\usetikzlibrary{arrows,arrows.meta} 
\tikzcdset{arrow style=tikz, diagrams={>=stealth}}
\tikzset{shorten <>/.style={shorten >=#1,shorten <=#1}}

\DeclareMathOperator{\Hom}{Hom}

\DeclareMathOperator{\End}{End}
\DeclareMathOperator{\Aut}{Aut}

\usepackage{mathtools}

\def\slot{{\color{blue}-}}

\makeatletter
\newcommand{\superimpose}[2]{%
  {\ooalign{$#1\@firstoftwo#2$\cr\hfil$#1\@secondoftwo#2$\hfil\cr}}}
\makeatother

\makeatletter
\DeclareRobustCommand{\cev}[1]{%
  \mathpalette\do@cev{#1}%
}
\newcommand{\do@cev}[2]{%
  \fix@cev{#1}{+}%
  \reflectbox{$\m@th#1\vec{\reflectbox{$\,\fix@cev{#1}{-}\m@th#1#2\fix@cev{#1}{+}$}}$}%
  \fix@cev{#1}{-}%
}
\newcommand{\fix@cev}[2]{%
  \ifx#1\displaystyle
    \mkern#23mu
  \else
    \ifx#1\textstyle
      \mkern#23mu
    \else
      \ifx#1\scriptstyle
        \mkern#22mu
      \else
        \mkern#22mu
      \fi
    \fi
  \fi
}

\makeatother

\usepackage{circledsteps}

\newcommand\Tsr[2][]{
\scalebox{.6}{
\ifmmode
\Circled[fill color=black,inner color=none,#1]{$#2$}
\else
\Circled[fill color=none,inner color=black,#1]{$#2$}
\fi
}
}

\newcommand\tsr[2][]{
\scalebox{.6}{
\ifmmode
\Circled[fill color=black,inner color=none,#1]{$#2$}
\else
\Circled[fill color=none,inner color=black,#1]{$#2$}
\fi
}
}

\def\id{\textit{id}}

\newlength{\LETTERheight}
\AtBeginDocument{\settoheight{\LETTERheight}{I}}

\def\unit{\textit{1} }
\def\op{{\textit{op}}}

\usepackage{scalerel}

\def\cp{\scaleobj{.92}{\vartriangle}}
\def\pc{\scaleobj{.92}{\mathrel\triangledown}}

\def\tensor{\scaleobj{.92}{{\,\circledast\,}}}
\def\ctimes{\scaleobj{1.2}{{\,\bullet\,}}}

\def\acts{{\,\diamond\,}}
\def\bcts{{\,\bar\diamond\,}}

\usepackage[allcommands,tikz]{overarrows}
\NewOverArrowCommand[tikz]{vecar}{
path={(0,0)--(0.9,0)},
arrows={->[scale=0.75]},
thinner,
}
\NewOverArrowCommand[tikz]{cevar}{
path={(0,0)--(0.9,0)},
arrows={<[scale=0.75]-},
thinner,
}
\NewOverArrowCommand[tikz]{hatar}{
path={(0,0)--(0.41,0.1)--(0.9,0)},
arrows={->[scale=0.75]},
thinner,
}
\NewOverArrowCommand[tikz]{tahar}{
path={(0,0)--(0.49,0.1)--(0.9,0)},
arrows={<[scale=0.75]-},
thinner,
}

\NewOverArrowCommand[tikz]{what}{
path={(0,0)--(0.5,0.1)--(1,0)},
arrows={-},
thinner,
}

\def\bcts{\bar{\acts}}

\makeatletter
\def\para{\@startsection{subsubsection}{3}{\z@}%
    {-13pt plus-8pt minus-4pt}{\z@}{\normalsize\itshape}}
\makeatletter

\DeclareMathOperator{\Isom}{Isom}

\DeclareMathOperator{\Gal}{Gal}
\DeclareMathOperator{\GL}{GL}

\begin{document}

\title{Prekosmic Grothendieck/Galois Categories}

\author{
Jaehyeok Lee
}
\institute{
Ph.D. thesis, Department of Mathematics, POSTECH, 2025
}

\maketitle

\abstract{
  We establish a generalized version of the duality 
  between groups and the categories of their representations on sets.
  Given an abstract symmetric monoidal category
  $\kos{K}$ called Galois prekosmos,
  we define pre-Galois objects in $\kos{K}$
  and study the categories of their representations internal to $\kos{K}$.
  The motivating example of $\kos{K}$
  is the cartesian monoidal category $\kos{S\!e\!t}$ of sets,
  and pre-Galois objects in $\kos{S\!e\!t}$ are groups.
  We present an axiomatic definition of pre-Galois $\kos{K}$-categories,
  which is a complete abstract characterization
  of the categories of representations of pre-Galois objects in $\kos{K}$.
  The category of covering spaces over
  a well-connected topological space
  is a prototype of a pre-Galois $\kos{S\!e\!t}$-category.
  We establish a perfect correspondence between
  pre-Galois objects in $\kos{K}$
  and 
  pre-Galois $\kos{K}$-categories
  pointed with pre-fiber functors.
  
  We also establish a generalized version of the duality 
  between flat affine group schemes
  and the categories of their linear representations.
  Given an abstract symmetric monoidal category
  $\kos{K}$ called Grothendieck prekosmos,
  we define what are pre-Grothendieck objects in $\kos{K}$
  and study the categories of their representations internal to $\kos{K}$.
  The motivating example of $\kos{K}$
  is the symmetric monoidal category $\kos{V\!\!e\!c}_k$
  of vector spaces over a field $k$,
  and pre-Grothendieck objects in $\kos{V\!\!e\!c}_k$
  are affine group $k$-schemes.
  We present an axiomatic definition of pre-Grothendieck $\kos{K}$-categories,
  which is a complete abstract characterization of the categories
  of representations of pre-Grothendieck objects in $\kos{K}$.
  The indization of a neutral Tannakian category over a field $k$
  is a prototype of a pre-Grothendieck $\kos{V\!\!e\!c}_k$-category.
  We establish a perfect correspondence
  between
  pre-Grothendieck objects in $\kos{K}$
  and
  pre-Grothendieck $\kos{K}$-categories
  pointed with pre-fiber functors.

  This is the author's Ph.D. thesis.
}

\newpage
\tableofcontents

\newpage
\section{Introduction}

Following the idea and vision of Grothendieck,
we study the relation, or duality, between
group-like objects and the categories of their internal representations,
defined in contexts of great generality,
by providing a complete abstract characterization of such categories.
We begin by introducing two dualities 
that were originally proposed by Grothendieck.

A Galois category is a category
which is equivalent to
the category $G\text{-}\cat{FSet}$ of finite left $G$-sets
for some profinite group $G$.
Grothendieck \cite{SGA1}
gave an axiomatic definition
and established the concept of Galois categories.
This was a consequence of his attempt to unify the Galois theory of covering spaces and that of field extensions.
In the axiomatic definition,
we require a Galois category to admit a fiber functor.
Each profinite group $G$
determines a Galois category $G\text{-}\cat{FSet}$
and the functor 
$G\text{-}\cat{FSet}\to \cat{FSet}$
of forgetting $G$-actions
is a fiber functor.
Conversely,
each Galois category together with a choice of a fiber functor
determines a profinite group,
which is the group of natural automorphisms 
of the chosen fiber functor.
These correspondences form a perfect duality
between profinite groups and the categories of their representations on finite sets.
Such duality guided Grothendieck towards his invention of \'{e}tale fundamental groups of schemes.

A Tannakian category over a field $k$
is a $k$-linear tensor category
which is equivalent to
the $k$-linear tensor category 
of representations of some affine gerbe
over the fpqc site
$\cat{Aff}_k$
of affine $k$-schemes.
Grothendieck suggested the concept of Tannakian categories,
to understand various cohomology groups of a scheme
as different realizations of its motive.
The theory of Tannakian categories was first developed by his student
Saavedra Rivano \cite{Saavedra1972}
and was later rectified by Deligne \cite{Deligne1990}.
By restricting our attention to a special kind of Tannakian categories
we obtain another duality,
between affine group schemes over $k$
and the categories of their linear representations on finite dimensional $k$-vector spaces.

A neutral Tannakian category over a field $k$
is a $k$-linear tensor category
which is equivalent to the $k$-linear tensor category
$\cat{FRep}_k(\CG)$
of finite dimensional linear representations
of some affine group scheme $\CG$ over $k$.
In the axiomatic definition,
we require
a neutral Tannakian category
to admit a neutral fiber functor.
Each affine group scheme $\CG$ over $k$ determines a neutral Tannakian category
$\cat{FRep}_k(\CG)$ over $k$
and the functor
$\cat{FRep}_k(\CG)\to \cat{FVec}_k$
of forgetting $\CG$-actions is a neutral fiber functor.
Conversely,
each neutral Tannakian category over $k$
together with a choice of a neutral fiber functor
determines an affine group scheme over $k$,
which represents the presheaf of groups on $\cat{Aff}_k$
of linear tensor natural automorphisms of the chosen neutral fiber functor.
These correspondences form a perfect duality
between affine group schemes over $k$
and neutral Tannakian categories over $k$
pointed with neutral fiber functors.

In this thesis we establish two kinds of dualities,
each of which is analogous to one of the dualities introduced above,
in a general situation
where groups and flat affine group schemes as well as their representations
are defined internal to an abstract symmetric monoidal category.

To explain in more detail, we divide the content of the thesis into two contexts.
We name one of them as \emph{Galois context}
and the other as \emph{Grothendieck context}.
The underlying mechanism that lies in each context
is Grothendieck's yoga of `three' operations.
Two contexts are completely parallel,
in the sense that they are categorically dual to each other.

We begin the story of Galois context of the thesis.

\begin{definition}
  A \emph{Galois prekosmos} $\kos{K}$
  is a symmetric monoidal category
  $\kos{K}=(\CK,\otimes,\kappa)$
  whose underlying category $\CK$ has reflexive coequalizers.
\end{definition}

The motivating example of a Galois prekosmos is 
the cartesian monoidal category
$\kos{S\!e\!t}=(\cat{Set},\times,\{*\})$
of sets.
Some examples of Galois prekosmoi
are the symmetric monoidal categories
\begin{equation}\label{eq Intro example of kosmoi}
  \kos{S\!e\!t}
  ,\text{ }
  \kos{s\!S\!e\!t}
  ,\text{ }
  \kos{M\!\!o\!\!d}\!_R
  ,\text{ }
  \kos{C\!h}\!_R
  ,\text{ }
  \kos{B\!a\!\!n}
  ,\text{ }
  \kos{C\!G\!W\!H}
  ,\text{ }
  \kos{S\!p}\!^{\Sigma}
\end{equation}
of
(small) sets,
simplicial sets,
modules over a commutative ring $R$,
(co)chain complexes over a commutative ring $R$,
Banach spaces with linear contractions,
compactly generated weakly Hausdorff spaces,
symmetric spectra.
Every elementary topos,
hence every Grothendieck topos, is also an example of Galois prekosmos.

In Galois context of the thesis,
we specify a Galois prekosmos
$\kos{K}=(\CK,\otimes,\kappa)$
and establish a generalized version of the duality between
groups and the categories of their representations on sets.
We denote 
$\kos{E\!n\!s}(\kos{K})=(\cat{Ens}(\kos{K}),\otimes,\kappa)$
as the cartesian monoidal category
of cocommutative comonoids in $\kos{K}$.
A \emph{pre-Galois object} in $\kos{K}$
is a group object $\pi$ in $\kos{E\!n\!s}(\kos{K})$
such that the functor
$\pi\otimes\slot:\CK\to\CK$
preserves reflexive coequalizers.
A \emph{representation} of $\pi$ in $\kos{K}$
is a pair of an object $x$ in $\CK$
and a left $\pi$-action morphism
$\pi\otimes x\to x$ in $\CK$.
When $\kos{K}=\kos{S\!e\!t}$,
a pre-Galois object $\pi$ in $\kos{S\!e\!t}$
is a group
and representations of $\pi$ in $\kos{S\!e\!t}$
are left $\pi$-sets.

Representations
of $\pi$ in $\kos{K}$
also form a Galois prekosmos $\kos{R\!e\!p}(\pi)$.
To collect the properties
of $\kos{R\!e\!p}(\pi)$,
we introduce the $2$-category
$\mathbb{GAL}^{\cat{pre}}$
of Galois prekosmoi.
Let $\kos{T}$ be another Galois prekosmos.
A $1$-morphism 
$\kos{f}:\kos{T}\to\kos{K}$
in $\mathbb{GAL}^{\cat{pre}}$
is called a \emph{Galois morphism},
which is an adjunction
$\kos{f}=(\kos{f}_{!}:\kos{T}\rightleftarrows \kos{K}:\kos{f}^*)$
whose right adjoint $\kos{f}^*$ is compatible with symmetric monoidal structures.
A $2$-morphism 
in $\mathbb{GAL}^{\cat{pre}}$
is called a \emph{Galois transformation},
which is a comonoidal natural transformation
between right adjoints.

We have a canonical Galois morphism
$\kos{t}_{\pi}:\kos{R\!e\!p}(\pi)\to \kos{K}$
whose right adjoint is the functor of constructing trivial representations.
Thus we introduce Galois prekosmoi relative to $\kos{K}$.
A \emph{Galois $\kos{K}$-prekosmos}
is a pair $\mathfrak{T}=(\kos{T},\kos{t})$
of a Galois prekosmos $\kos{T}$
and a Galois morphism
$\kos{t}:\kos{T}\to\kos{K}$.
We can see $\kos{K}$ itself as a Galois $\kos{K}$-prekosmos
which we denote as $\text{$\mathfrak{K}$}=(\kos{K},\id_{\kos{K}})$.

We define the $2$-category
$\mathbb{GAL}^{\cat{pre}}_{\kos{K}}$
of Galois $\kos{K}$-prekosmoi
as the slice $2$-category of
$\mathbb{GAL}^{\cat{pre}}$ over $\kos{K}$.
A $1$-morphism in $\mathbb{GAL}^{\cat{pre}}_{\kos{K}}$
is called a \emph{Galois $\kos{K}$-morphism}
and a $2$-morphism in $\mathbb{GAL}^{\cat{pre}}_{\kos{K}}$
is called a \emph{Galois $\kos{K}$-transformation}.

\begin{definition}
  Let $\kos{K}$ be a Galois prekosmos.
  A \emph{pre-Galois $\kos{K}$-category}
  is a Galois $\kos{K}$-prekosmos $\mathfrak{T}=(\kos{T},\kos{t})$
  which is equivalent to the
  Galois $\kos{K}$-prekosmos 
  $\mathfrak{Rep}(\pi)=(\kos{R\!e\!p}(\pi),\kos{t}_{\pi})$
  of representations of $\pi$ in $\kos{K}$
  for some pre-Galois object $\pi$ in $\kos{K}$.
\end{definition}

We have a Galois $\kos{K}$-morphism
$\varpi_{\pi}:\mathfrak{K}\to\mathfrak{Rep}(\pi)$
whose right adjoint is the functor of forgetting
actions of $\pi$.
This is a prototype of a \emph{surjective pre-fiber functor}.

Given a Galois $\kos{K}$-prekosmos $\mathfrak{T}=(\kos{T},\kos{t})$,
we define a \emph{pre-fiber functor} for $\mathfrak{T}$
as a Galois $\kos{K}$-morphism
$\varpi:\mathfrak{K}\to\mathfrak{T}$ 
which satisfies some (but not necessarily all) of the axiomatic properties of
$\varpi_{\pi}:\mathfrak{K}\to\mathfrak{Rep}(\pi)$.
We say a pre-fiber functor $\varpi$ for $\mathfrak{T}$
is \emph{surjective}
if in addition
the right adjoint is conservative and preserves reflexive coequalizers.

We present an axiomatic definition of pre-Galois $\kos{K}$-categories.

\begin{definition}[Axiomatic]
  \label{def Intro axiomatic preGalKCat}
  Let $\kos{K}$ be a Galois prekosmos.
  A \emph{pre-Galois $\kos{K}$-category}
  is a Galois $\kos{K}$-prekosmos
  $\mathfrak{T}=(\kos{T},\kos{t})$
  which admits a surjective pre-fiber functor
  $\varpi:\mathfrak{K}\to\mathfrak{T}$.
\end{definition}

We show the following,
which implies that the two definitions are equivalent.

\begin{theorem}\label{mainthm Intro GAL Justifying axiomaticdef}
  Let $\kos{K}$ be a Galois prekosmos
  and let $\mathfrak{T}=(\kos{T},\kos{t})$
  be a pre-Galois $\kos{K}$-category 
  defined as in Definition~\ref{def Intro axiomatic preGalKCat}.
  For each pre-fiber functor $\varpi:\mathfrak{K}\to\mathfrak{T}$,
  the following are true.
  \begin{enumerate}
    \item 
    There exists a pre-Galois object
    $\pi$ in $\kos{K}$
    which represents the presheaf of groups of
    invertible Galois $\kos{K}$-transformations
    from $\varpi$ to $\varpi$.
    \begin{equation*}
      \vcenter{\hbox{
        \xymatrix{
          \Hom_{\cat{Ens}(\kos{K})}(\slot,\pi)
          \ar@2{->}[r]^-{\cong}
          &\underline{\Aut}_{\mathbb{GAL}^{\cat{pre}}_{\kos{K}}}(\varpi)
          :\cat{Ens}(\kos{K})^{\op}\to\cat{Grp}
        }
      }}
    \end{equation*}

    \item
    The given pre-fiber functor
    $\varpi:\mathfrak{K}\to \mathfrak{T}$
    is surjective,
    and factors through as an equivalence of Galois $\kos{K}$-prekosmoi
    $\widebreve{\varpi}:\mathfrak{Rep}(\pi)\xrightarrow{\simeq} \mathfrak{T}$.
    \begin{equation*}
      \vcenter{\hbox{
        \xymatrix@C=40pt{
          \mathfrak{T}
          &\mathfrak{Rep}(\pi)
          \ar[l]_-{\widebreve{\varpi}}^-{\simeq}
          \\
          \text{ }
          &\mathfrak{K}
          \ar@/^1pc/[ul]^-{\varpi}
          \ar[u]_-{\varpi_{\pi}}
        }
      }}
    \end{equation*}
  \end{enumerate}
\end{theorem}

We establish a perfect correspondence
between pre-Galois objects $\pi$ in $\kos{K}$
and pre-Galois $\kos{K}$-categories $\mathfrak{T}$
pointed with pre-fiber functors $\varpi$.

\begin{theorem} \label{mainthm Intro GAL 2catduality}
  Let $\kos{K}$ be a Galois prekosmos.
  We have a biequivalence
  \begin{equation*}
    \vcenter{\hbox{
      \xymatrix@C=30pt{
        \mathbb{GALOBJ}^{\cat{pre}}(\kos{K})
        \ar@<0.5ex>[r]^-{\simeq}
        &\mathbb{GALCAT}^{\cat{pre}}_*(\kos{K})
        \ar@<0.5ex>[l]^-{\simeq}
      }
    }}
  \end{equation*}
  between the $(2,1)$-category of pre-Galois objects in $\kos{K}$
  and the $(2,1)$-category of pre-Galois $\kos{K}$-categories
  pointed with pre-fiber functors.
\end{theorem}

A \emph{right $\pi$-torsor over $\kappa$}
is an object $p$ in $\cat{Ens}(\kos{K})$
equipped with a right $\pi$-action morphism
$p\otimes \pi\to p$ in $\cat{Ens}(\kos{K})$
such that
the functor $p\otimes\slot:\CK\to \CK$
is conservative and preserves reflexive coequalizers, and
the induced morphism $p\otimes \pi\to p\otimes p$
in $\cat{Ens}(\kos{K})$ is an isomorphism.
When $\kos{K}=\kos{S\!e\!t}$,
$\pi$ is a group and right $\pi$-torsors over $\kappa=\{*\}$
are the usual right torsors of the group $\pi$.

We also establish the following correspondence
between torsors and pre-fiber functors.

\begin{theorem}\label{mainthm Intro GAL torsor and prefib}
  Let $\kos{K}$ be a Galois prekosmos.
  For each pre-Galois $\kos{K}$-category 
  $\mathfrak{T}$ pointed with a pre-fiber functor
  $\varpi$,
  denote $\pi$ as the pre-Galois object in $\kos{K}$
  described in Theorem~\ref{mainthm Intro GAL Justifying axiomaticdef}.
  We have an adjoint equivalence
  \begin{equation*}
    \xymatrix{
      \cat{Fib}(\mathfrak{T})^{\cat{pre}}
      \ar@<0.5ex>[r]^-{\simeq}
      &(\cat{Tors-}\pi)_{\kappa}
      \ar@<0.5ex>[l]^-{\simeq}
    }
  \end{equation*}
  between the groupoid of pre-fiber functors for $\mathfrak{T}$
  and the groupoid of right $\pi$-torsors over $\kappa$.
\end{theorem}

We move on to the Grothendieck context of the thesis,
whose story is completely parallel to that of Galois context.

\begin{definition}
  A \emph{Grothendieck prekosmos} $\kos{K}$
  is a symmetric monoidal category
  $\kos{K}=(\CK,\otimes,\kappa)$
  whose underlying category $\CK$ has coreflexive equalizers.
\end{definition}

The motivating example of a Grothendieck prekosmos is 
the symmetric monoidal category
$\kos{V\!\!e\!c}_k=(\cat{Vec}_k,\otimes_k,k)$
of vector spaces over a field $k$.
All of the symmetric monoidal categories 
mentioned in (\ref{eq Intro example of kosmoi})
are also examples of Grothendieck prekosmoi.

In Grothendieck context,
we specify a Grothendieck prekosmos $\kos{K}$
and establish a generalized version of the duality between
flat affine group schemes and the categories of their linear representations.
We denote
$\kos{C\!o\!\!m\!\!m}(\kos{K})=(\cat{Comm}(\kos{K}),\otimes,\kappa)$
as the cocartesian monoidal category
of commutative monoids in $\kos{K}$,
and
$\kos{A\!f\!f}(\kos{K})=(\cat{Aff}(\kos{K}),\times,\spec{\kappa})$
as the opposite cartesian monoidal category.
For each commutative monoid $b$ in $\kos{K}$,
we denote $\cat{Spec}(b)$
as the corresponding object in $\cat{Aff}(\kos{K})$.

A \emph{pre-Grothendieck object} in $\kos{K}$
is a group object $\spec{\pi}$ in $\kos{A\!f\!f}(\kos{K})$
such that the functor $\slot\otimes \pi:\CK\to\CK$
preserves coreflexive equalizers.
A \emph{representation} of
$\spec{\pi}$ in $\kos{K}$
is a pair of an object $x$ in $\CK$
and a right $\pi$-coaction morphism $x\to x\otimes \pi$ in $\CK$.
When $\kos{K}=\kos{V\!\!e\!c}_k$ for some field $k$,
a pre-Grothendieck object $\cat{Spec}(\pi)$
in $\kos{V\!\!e\!c}_k$
is an affine group $k$-scheme
and representations of $\cat{Spec}(\pi)$ in $\kos{V\!\!e\!c}_k$
are linear representations of $\cat{Spec}(\pi)$
on $k$-vector spaces.

Representations
of $\cat{Spec}(\pi)$ in $\kos{K}$
also form a Grothendieck prekosmos $\kos{R\!e\!p}(\pi)$.
We introduce the $2$-category
$\mathbb{GRO}^{\cat{pre}}$
of Grothendieck prekosmoi.
Let $\kos{T}$ be another Grothendieck prekosmos.
A $1$-morphism 
$\kos{f}:\kos{T}\to\kos{K}$
in $\mathbb{GRO}^{\cat{pre}}$
is called a \emph{Grothendieck morphism},
which is an adjunction
$\kos{f}=(\kos{f}^*:\kos{K}\rightleftarrows \kos{T}:\kos{f}_*)$
whose left adjoint $\kos{f}^*$
is compatible with symmetric monoidal structures.
A $2$-morphism 
in $\mathbb{GRO}^{\cat{pre}}$
is called a \emph{Grothendieck transformation},
which is a monoidal natural transformation
between left adjoints.

We have a canonical Grothendieck morphism
$\kos{t}_{\pi}:\kos{R\!e\!p}(\pi)\to \kos{K}$
whose left adjoint is the functor of constructing trivial representations.
We define the $2$-category
$\mathbb{GRO}^{\cat{pre}}_{\kos{K}}$
of \emph{Grothendieck $\kos{K}$-prekosmoi},
\emph{Grothendieck $\kos{K}$-morphisms}
and \emph{Grothendieck $\kos{K}$-transformations}
as the slice $2$-category of
$\mathbb{GRO}^{\cat{pre}}$ over $\kos{K}$.

\begin{definition}
  Let $\kos{K}$ be a Grothendieck prekosmos.
  A \emph{pre-Grothendieck $\kos{K}$-category}
  is defined as a Grothendieck $\kos{K}$-prekosmos $\mathfrak{T}=(\kos{T},\kos{t})$
  which is equivalent to the
  Grothendieck $\kos{K}$-prekosmos 
  $\mathfrak{Rep}(\pi)=(\kos{R\!e\!p}(\pi),\kos{t}_{\pi})$
  of representations of $\cat{Spec}(\pi)$ in $\kos{K}$
  for some pre-Grothendieck object $\cat{Spec}(\pi)$ in $\kos{K}$.
\end{definition}

We have a Grothendieck $\kos{K}$-morphism
$\varpi_{\pi}:\mathfrak{K}\to\mathfrak{Rep}(\pi)$
whose left adjoint is the functor of forgetting
actions of $\cat{Spec}(\pi)$.
This is a prototype of a \emph{surjective pre-fiber functor}.

Given a Grothendieck $\kos{K}$-prekosmos
$\mathfrak{T}=(\kos{T},\kos{t})$,
we define a \emph{pre-fiber functor} for $\mathfrak{T}$
as a Grothendieck $\kos{K}$-morphism
$\varpi:\mathfrak{K}\to\mathfrak{T}$ 
which satisfies some of the axiomatic properties of
$\varpi_{\pi}:\mathfrak{K}\to\mathfrak{Rep}(\pi)$.
We say a pre-fiber functor $\varpi$ for $\mathfrak{T}$
is \emph{surjective}
if in addition
the left adjoint is conservative and preserves coreflexive equalizers.

We present an axiomatic definition of pre-Grothendieck $\kos{K}$-categories,
and show that the two definitions are equivalent.

\begin{definition}[Axiomatic]
  \label{def Intro axiomatic preGroKCat}
  Let $\kos{K}$ be a Grothendieck prekosmos.
  We define a \emph{pre-Grothendieck $\kos{K}$-category}
  as a Grothendieck $\kos{K}$-prekosmos
  $\mathfrak{T}=(\kos{T},\kos{t})$
  which admits a surjective pre-fiber functor
  $\varpi:\mathfrak{K}\to\mathfrak{T}$.
\end{definition}

\begin{theorem}\label{mainthm Intro GRO Justifying axiomaticdef}
  Let $\kos{K}$ be a Grothendieck prekosmos
  and let $\mathfrak{T}=(\kos{T},\kos{t})$
  be a pre-Grothendieck $\kos{K}$-category 
  defined as in Definition~\ref{def Intro axiomatic preGroKCat}.
  For each pre-fiber functor $\varpi:\mathfrak{K}\to\mathfrak{T}$,
  the following are true.
  \begin{enumerate}
    \item 
    We have a pre-Grothendieck object
    $\cat{Spec}(\pi)$ in $\kos{K}$
    which represents the presheaf of groups of
    invertible Grothendieck $\kos{K}$-transformations
    from $\varpi$ to $\varpi$.
    \begin{equation*}
      \vcenter{\hbox{
        \xymatrix{
          \Hom_{\cat{Aff}(\kos{K})}(\slot,\cat{Spec}(\pi))
          \ar@2{->}[r]^-{\cong}
          &\underline{\Aut}_{\mathbb{GRO}^{\cat{pre}}_{\kos{K}}}(\varpi)
          :\cat{Aff}(\kos{K})^{\op}\to\cat{Grp}
        }
      }}
    \end{equation*}

    \item
    The given pre-fiber functor
    $\varpi:\mathfrak{K}\to \mathfrak{T}$
    is surjective,
    and factors through as an equivalence of Grothendieck $\kos{K}$-prekosmoi
    $\widebreve{\varpi}:\mathfrak{Rep}(\pi)\xrightarrow{\simeq} \mathfrak{T}$.
    \begin{equation*}
      \vcenter{\hbox{
        \xymatrix@C=40pt{
          \mathfrak{T}
          &\mathfrak{Rep}(\pi)
          \ar[l]_-{\widebreve{\varpi}}^-{\simeq}
          \\
          \text{ }
          &\mathfrak{K}
          \ar@/^1pc/[ul]^-{\varpi}
          \ar[u]_-{\varpi_{\pi}}
        }
      }}
    \end{equation*}
  \end{enumerate}
\end{theorem}

We establish a perfect correspondence
between pre-Grothendieck objects
in $\kos{K}$
and 
pre-Grothendieck $\kos{K}$-categories
pointed with pre-fiber functors.

\begin{theorem} \label{mainthm Intro GRO 2catduality}
  Let $\kos{K}$ be a Grothendieck prekosmos.
  We have a biequivalence
  \begin{equation*}
    \vcenter{\hbox{
      \xymatrix@C=30pt{
        \mathbb{GROOBJ}^{\cat{pre}}(\kos{K})
        \ar@<0.5ex>[r]^-{\simeq}
        &\mathbb{GROCAT}^{\cat{pre}}_*(\kos{K})
        \ar@<0.5ex>[l]^-{\simeq}
      }
    }}
  \end{equation*}
  between the $(2,1)$-category of pre-Grothendieck objects in $\kos{K}$
  and the $(2,1)$-category of pre-Grothendieck $\kos{K}$-categories
  pointed with pre-fiber functors.
\end{theorem}

A \emph{left $\cat{Spec}(\pi)$-torsor over $\cat{Spec}(\kappa)$}
is an object $\cat{Spec}(p)$ in $\cat{Aff}(\kos{K})$
equipped with a left $\cat{Spec}(\pi)$-action morphism
$\cat{Spec}(\pi)\times \cat{Spec}(p)\to \cat{Spec}(p)$ in $\cat{Aff}(\kos{K})$
such that the functor $\slot\otimes p:\CK\to \CK$
is conservative and preserves coreflexive equalizers, and
the induced morphism
$\cat{Spec}(\pi)\times \cat{Spec}(p)\to \cat{Spec}(p)\times \cat{Spec}(p)$
in $\cat{Aff}(\kos{K})$
is an isomorphism.

When $\kos{K}=\kos{V\!\!e\!c}_k$ for some field $k$,
$\cat{Spec}(\pi)$ is an affine group scheme over $k$
and a left $\cat{Spec}(\pi)$-torsor over $\cat{Spec}(\kappa)$
is the usual left torsor
over $\cat{Spec}(k)$
with respect to the fpqc topology.

We also establish the following correspondence
between torsors and pre-fiber functors.

\begin{theorem}\label{mainthm Intro GRO torsor and prefib}
  Let $\kos{K}$ be a Grothendieck prekosmos.
  For each pre-Grothendieck $\kos{K}$-category 
  $\mathfrak{T}$ pointed with a pre-fiber functor
  $\varpi$,
  denote $\cat{Spec}(\pi)$
  as the pre-Grothendieck object in $\kos{K}$
  described in Theorem~\ref{mainthm Intro GRO Justifying axiomaticdef}.
  We have an adjoint equivalence
  \begin{equation*}
    \xymatrix{
      \cat{Fib}(\mathfrak{T})^{\cat{pre}}
      \ar@<0.5ex>[r]^-{\simeq}
      &(\cat{Spec}(\pi)\cat{-Tors})_{\kappa}^{\op}
      \ar@<0.5ex>[l]^-{\simeq}
    }
  \end{equation*}
  between the groupoid of pre-fiber functors for $\mathfrak{T}$
  and the opposite groupoid of left $\cat{Spec}(\pi)$-torsors over $\cat{Spec}(\kappa)$.
\end{theorem}

\subsection{Motivations}
\label{subsec Motivations}

\subsubsection{Fundamental group, coverings and path torsors}
\label{subsubsec pi1,covering,torsor}
Let $X$ be a
topological space
and let $p$ be a point in $X$.
For each covering $E$ over $X$,
its fiber $E_p$ over $p$
is equipped with the monodromy action
of the fundamental group $\pi_1(X,p)$ of $X$ based at $p$.
Assume that $X$ is connected, locally path-connected and semi-locally simply connected.
Then the functor
$\cat{Fib}_p:E\mapsto E_p$
of taking fibers over $p$
is represented by a universal covering,
and the automorphism group
$\Aut(\cat{Fib}_p)$ of this functor
is isomorphic to $\pi_1(X,p)$.
The induced action of $\Aut(\cat{Fib}_p)$ on
the fiber $E_p$
recovers the monodromy action,
and the functor
$\cat{Fib}_p:E\mapsto E_p$
lifts to an equivalence
between the category of coverings over $X$
and the category of left $\pi_1(X,p)$-sets.
Let $q$ be another point in $X$.
The set
$\Isom(\cat{Fib}_p,\cat{Fib}_q)$
of isomorphisms of functors
is a torsor with respect to the right action of
$\Aut(\cat{Fib}_p)\cong \pi_1(X,p)$,
and is isomorphic as torsors
to the path torsor $\pi_1(X;p,q)$.

\subsubsection{Grothendieck's formulation of Galois theory}
\label{subsubsec Grothendieck'sGaloisThy}
Let $X$ and $p$ be as in \textsection~\ref{subsubsec pi1,covering,torsor}
and restrict our attention
to finite coverings $E^{\cat{fin}}$ over $X$.
The functor
$\cat{Fib}_p^{\cat{fin}}:E^{\cat{fin}}\mapsto E^{\cat{fin}}_p$
of taking fibers over $p$ is then pro-represented by
finite Galois coverings over $X$.
The automorphism group
$\Aut(\cat{Fib}_p^{\cat{fin}})$
is isomorphic to the profinite completion
$\widehat{\pi_1(X,p)}$
of the fundamental group,
and the functor
$\cat{Fib}_p^{\cat{fin}}:E^{\cat{fin}}\mapsto E^{\cat{fin}}_p$
lifts to an equivalence between
the category of finite coverings over $X$
and the category of finite sets equipped with continuous left actions of
$\widehat{\pi_1(X,p)}$.

Meanwhile, let $k$ be a field.
A finite \'{e}tale $k$-algebra
is a commutative $k$-algebra
which is isomorphic to a finite product
of finite separable extensions of $k$.
Choose an algebraic closure $\bar{k}$ of $k$
and denote $k^{\cat{sep}}$
as the separable closure of $k$ in $\bar{k}$.
For each finite \'{e}tale $k$-algebra $A$,
the set
$\Hom_k(A,k^{\cat{sep}})$
of morphisms of $k$-algebras
is finite and
has a canonical continuous left
action of the absolute Galois group
$\Gal(k^{\cat{sep}}/k)$.
The contravariant functor
$A\mapsto \Hom_k(A,k^{\cat{sep}})$
is pro-represented by finite Galois extensions over $k$,
and the automorphism group of this functor
is isomorphic to 
$\Gal(k^{\cat{sep}}/k)$.
The functor
$A\mapsto \Hom_k(A,k^{\cat{sep}})$
lifts to an equivalence between
the opposite category 
of the category of finite \'{e}tale $k$-algebras
and the category of
finite continuous left $\Gal(k^{\cat{sep}}/k)$-sets.
This is often referred to as 
\emph{Grothendieck's formulation of Galois theory of field extensions}.

A \emph{Galois category} is a category
$\CC$
that is equivalent to
the category $G\text{-}\cat{FSet}$
of finite continuous left $G$-sets
for some profinite group $G$.
In the axiomatic definition
given by Grothendieck,
we require a Galois category $\CC$
to admit a functor 
$\CF:\CC\to \cat{FSet}$
to the category of finite sets,
which satisfies axiomatic properties that
the functors
$E^{\cat{fin}}\mapsto E^{\cat{fin}}_p$
and
$A\mapsto \Hom_k(A,k^{\cat{sep}})$
have.
Such functor $\CF$ is called a 
\emph{fiber functor} for $\CC$.
Each profinite group $G$ determines the Galois category
$G\text{-}\cat{FSet}$ of finite left $G$-sets,
and the functor $\cat{forget}:G\text{-}\cat{FSet}\to \cat{FSet}$
of forgetting $G$-actions is a fiber functor.
Conversely, each Galois category $\CC$
together with a choice of a fiber functor $\CF$
determines the profinite group
$\Aut(\CF)$ of natural automorphisms of $\CF$.

There is a complete analogy, or a \emph{duality},
between profinite groups and the categories of their representations on finite sets.
Given a profinite group $G$,
we have an isomorphism
$G\xrightarrow{\cong}\Aut(\cat{forget})$ of profinite groups.
Given a Galois category $\CC$,
each fiber functor $\CF$ for $\CC$
lifts to an equivalence
of categories $\CC\xrightarrow{\simeq}\Aut(\CF)\text{-}\cat{FSet}$.
The correspondences
$G\mapsto (G\text{-}\cat{FSet},\cat{forget})$
and
$(\CC,\CF)\mapsto \Aut(\CF)$
define an equivalence
between
the $(2,1)$-category of profinite groups
and the $(2,1)$-category of Galois categories pointed with fiber functors.

\subsubsection{Neutral Tannakian categories}

Let $k$ be a field
 and denote $\cat{Aff}_k$ as the fpqc site of affine $k$-schemes.
A fibered category $\Pi$ over $\cat{Aff}_k$ is called a \emph{gerbe}
if it is a stack in groupoids 
which satisfies the following properties:
(1) it is locally nonempty,
and
(2) any two objects in the same fiber are locally isomorphic.
A gerbe $\Pi$ over $\cat{Aff}_k$
is called \emph{affine} if for each object $\pzc{x}$ in the fiber $\Pi_R$
over $\spec{R}$,
the presheaf of groups $\underline{\Aut}(\pzc{x})$
is represented by an affine group scheme over $R$.
Denote $\cat{Proj}$ as the stack over $\cat{Aff}_k$
whose fiber at $\spec{R}$ is the category $\cat{Proj}_R$
of finitely generated projective $R$-modules.
A \emph{representation} of an affine gerbe $\Pi$ over $\cat{Aff}_k$
is a morphism $\Pi\to \cat{Proj}$ of fibered categories over $\cat{Aff}_k$.
A \emph{Tannakian category} over $k$
is a $k$-linear tensor category
$\pzc{T}$ which is equivalent to
the $k$-linear tensor category 
of representations of some affine gerbe $\Pi$ over $\cat{Aff}_k$.

Let $\CG$ be an affine group scheme over $k$.
A \emph{linear representation} of $\CG$ over $k$
is a pair of a $k$-vector space $V$
and a morphism 
$\CG\Rightarrow \GL_V$ of presheaves of groups
on $\cat{Aff}_k$.
It is called \emph{finite dimensional} if $V$ is a finite dimensional $k$-vector space.
A \emph{neutral Tannakian category} over $k$
is a $k$-linear tensor category $\pzc{T}$
which is equivalent to the $k$-linear tensor category
$\cat{FRep}_k(\CG)$
of finite dimensional linear representations
of some affine group scheme $\CG$ over $k$.

We explain why neutral Tannakian categories are special kind of Tannakian categories.
Each affine group scheme $\CG$ over $k$
determines an affine gerbe $\cat{Tors}(\CG)$ over $\cat{Aff}_k$,
whose fiber $\cat{Tors}(\CG)_R$
over $\cat{Spec}(R)$
is the category of $\CG$-torsor over $\cat{Spec}(R)$
for the fpqc topology.
We can regard $\CG$ itself as a $\CG$-torsor over $\cat{Spec}(k)$
and thus the fiber $\cat{Tors}(\CG)_k$ over $\cat{Spec}(k)$ is nonempty.
Each representation 
$\Lambda:\cat{Tors}(\CG)\to\cat{Proj}$
of $\cat{Tors}(\CG)$
induces a finite dimensional linear representation
of $\CG$,
whose underlying $k$-vector space
is the image $\Lambda_k(\CG)$
of $\CG\in \obj{\cat{Tors}(\CG)_k}$ under $\Lambda$.
The correspondence 
$\Lambda\mapsto \Lambda_k(\CG)$
defines an equivalence
from
the $k$-linear tensor category
of representations of $\cat{Tors}(\CG)$
to
the $k$-linear tensor category
$\cat{FRep}_k(\CG)$
of finite dimensional linear representations of $\CG$.

In the axiomatic definition,
we require a neutral Tannakian category $\pzc{T}$
to be a rigid abelian category,
where the term \emph{rigid} means that 
every object $X$ in $\pzc{T}$ has a dual object $X^{\vee}$.
We also require $\pzc{T}$ to admit a \emph{neutral fiber functor},
which is an exact faithful $k$-linear tensor functor
$\pzc{w}:\pzc{T}\to\cat{FVec}_k$
to the category of finite dimensional $k$-vector spaces.
Each neutral fiber functor $\pzc{w}$ for $\pzc{T}$
defines a presheaf of groups
$\underline{\Aut}^{\otimes}(\pzc{w})$
on $\cat{Aff}_k$,
which sends an affine $k$-scheme $\cat{Spec}(R)$
to the group of $R$-linear tensor natural automorphisms
of the functor $\pzc{w}_R:\pzc{T}\xrightarrow{\pzc{w}}\cat{FVec}_k\xrightarrow{\slot\otimes_kR}\cat{Mod}_R$.

Each affine group scheme $\CG$ over $k$
determines the neutral Tannakian category
$\cat{FRep}_k(\CG)$ over $k$,
and the functor 
$\cat{forget}:\cat{FRep}_k(\CG)\to \cat{FVec}_k$
of forgetting  $\CG$-actions is a neutral fiber functor.
We can recover the given $\CG$
from the pair $(\cat{FRep}_k(\CG),\cat{forget})$
as we have a canonical isomorphism
$\CG\xrightarrow{\cong} \underline{\Aut}^{\otimes}(\cat{forget})$
of presheaves of groups on $\cat{Aff}_k$.
Conversely,
let $\pzc{T}$ be a neutral Tannakian category over $k$
and let
$\pzc{w}$
be a neutral fiber functor for $\pzc{T}$.
Saavedra Rivano \cite{Saavedra1972}
showed that there exists an affine group scheme 
$\CG(\pzc{T},\pzc{w})$ over $k$
which represents the presheaf of groups 
$\underline{\Aut}^{\otimes}(\pzc{w})$
on $\cat{Aff}_k$, and
the given neutral fiber functor $\pzc{w}:\pzc{T}\to \cat{FVec}_k$
lifts to an equivalence of 
$k$-linear tensor categories
$\pzc{T}\xrightarrow{\simeq}\cat{FRep}_k(\CG(\pzc{T},\pzc{w}))$.
The correspondences
$\CG\mapsto (\cat{FRep}_k(\CG),\cat{forget})$
and
$(\pzc{T},\pzc{w})\mapsto \CG(\pzc{T},\pzc{w})$
provide a duality between affine group schemes over $k$
and neutral Tannakian categories over $k$
pointed with neutral fiber functors.

Let $(\pzc{T},\pzc{w})$ be a neutral Tannakian category over $k$
pointed with a neutral fiber functor
and consider the affine group scheme $\CG(\pzc{T},\pzc{w})$
over $k$ which represents $\underline{\Aut}^{\otimes}(\pzc{w})$.
Different neutral fiber functors for $\pzc{T}$
correspond to fpqc torsors of $\CG(\pzc{T},\pzc{w})$
over $\cat{Spec}(k)$.
This is analogous to the situation in \textsection~\ref{subsubsec pi1,covering,torsor}.
Given another neutral fiber functor
$\pzc{w}^{\pr}$ for $\pzc{T}$,
the presheaf
$\underline{\Isom}^{\otimes}(\pzc{w}^{\pr},\pzc{w})$
on $\cat{Aff}_k$
of linear tensor natural isomorphisms from $\pzc{w}^{\pr}$ to $\pzc{w}$
is represented by a left $\CG$-torsor over $\cat{Spec}(k)$.
This defines an equivalence
from the groupoid $\cat{Fib}(\pzc{T})_k$
of neutral fiber functors for $\pzc{T}$
to the opposite groupoid $(\CG(\pzc{T},\pzc{w})\text{-}\cat{Tors})_k^{\cat{op}}$
of left $\CG(\pzc{T},\pzc{w})$-torsors over $\spec{k}$.

\subsection{Outline of the thesis}
\label{subsec Outline}

The Galois context of the thesis
is throughout
\textsection~\ref{sec ColaxKtensorCat}-\textsection~\ref{sec PreGalCat}.
In \textsection~\ref{sec ColaxKtensorCat},
we collect technical tools that we crucially use in Galois context.
We fix a symmetric monoidal category $\kos{K}$
and study in general about symmetric monoidal categories
$\kos{T}$ equipped with colax symmetric monoidal functors
$\mon{t}:\kos{K}\to \kos{T}$.
In \textsection~\ref{sec PreGalCat},
we present the main story of Galois context of the thesis.
We fix a Galois prekosmos $\kos{K}$
and study about pre-Galois objects in $\kos{K}$
and pre-Galois $\kos{K}$-categories.
Theorem~\ref{mainthm Intro GAL Justifying axiomaticdef},
Theorem~\ref{mainthm Intro GAL 2catduality},
Theorem~\ref{mainthm Intro GAL torsor and prefib}
are presented as
Theorem~\ref{thm preGalKCat mainThm1},
Theorem~\ref{thm GALOBJbiequivGALCAT},
Theorem~\ref{thm TorsFib adjequiv}.

The Grothendieck context of the thesis
is throughout
\textsection~\ref{sec LaxKtensorCat}-\textsection~\ref{sec PreGroCat}.
In \textsection~\ref{sec LaxKtensorCat},
we collect technical tools that we crucially use in Grothendieck context.
We fix a symmetric monoidal category $\kos{K}$
and study in general about symmetric monoidal categories
$\kos{T}$ equipped with lax symmetric monoidal functors
$\mon{t}:\kos{K}\to \kos{T}$.
In \textsection~\ref{sec PreGroCat},
we present the main story of Grothendieck context of the thesis.
We fix a Grothendieck prekosmos $\kos{K}$
and study about pre-Grothendieck objects in $\kos{K}$
and pre-Grothendieck $\kos{K}$-categories.
Theorem~\ref{mainthm Intro GRO Justifying axiomaticdef},
Theorem~\ref{mainthm Intro GRO 2catduality},
Theorem~\ref{mainthm Intro GRO torsor and prefib}
are presented as
Theorem~\ref{thm preGroKCat mainThm1},
Theorem~\ref{thm GROOBJbiequivGROCAT},
Theorem~\ref{thm Gro TorsFib adjequiv}.

\subsection{Notations and preliminaries}
\label{subsec Preliminaries}

Let
$\CL:\CC\rightleftarrows \CD:\CR$
be an adjoint pair of functors
between categories $\CC$, $\CD$.
We say the adjunction
$\CL\dashv\CR$ is \emph{reflective}
if the right adjoint $\CR$ is fully faithful.
In this case, 
the essential image of the right adjoint
$\CR$
is a reflective full subcategory of $\CC$.
Dually, we say $\CL\dashv\CR$
is \emph{coreflective}
if the left adjoint $\CL$ is fully faithful.
In this case, 
the essential image of the left adjoint
$\CL$
is a coreflective full subcategory of $\CD$.

We refer \cite{Johnson2021} for
basic $2$-categorical notions that we use.
This includes 
symmetric monoidal categories,
Mac Lane's coherence theorem,
reversed monoidal categories,
colax/lax symmetric monoidal functors,
comonoidal/monoidal natural transformations,
$2$-categories,
bijective correspondences of mates,
slice $2$-categories,
colax/lax coslice $2$-categories.

We denote a symmetric monoidal category $\kos{K}$
as a triple
$\kos{K}=(\CK,\otimes,\kappa)$.
Object $x$, $y$, $z$, $w$ in $\CK$ are written in small letters
and symmetric monoidal coherence isomorphisms are denoted as
$a_{x,y,z}:x\otimes (y\otimes z)\xrightarrow{\cong}(x\otimes y)\otimes z$,
$\imath_x:x\xrightarrow{\cong}\kappa\otimes x$,
$\jmath_x:x\xrightarrow{\cong}x\otimes \kappa$,
$s_{x,y}:x\otimes y\xrightarrow{\cong}y\otimes x$.
Given any symmetric monoidal category
$\kos{T}=(\CT,\tensor,\unit)$ other than $\kos{K}$,
we try to avoid using small letters when referring to an object in $\CT$.

Let $\kos{T}=(\CT,\tensor,\unit)$,
$\text{$\kos{S}$}=(\CS,\ctimes,\pzc{1})$
be symmetric monoidal categories.
We denote a colax symmetric monoidal functor as
$\phi:\kos{S}\to\kos{T}$
and use the same notation for its underlying functor
$\phi:\CS\to \CT$.
The colax symmetric monoidal coherence morphisms are 
denoted as
\begin{equation*}
  \phi_{\text{$\pzc{X}$},\text{$\pzc{Y}$}}
  :\phi(\text{$\pzc{X}$}\ctimes \text{$\pzc{Y}$})
  \to \phi(\text{$\pzc{X}$})\tensor \phi(\text{$\pzc{Y}$})
  ,
  \qquad
  \phi_{\text{$\pzc{1}$}}:\phi(\text{$\pzc{1}$})\to \unit
  ,
  \qquad
  \text{$\pzc{X}$}, \text{$\pzc{Y}$}\in \obj{\CS}
  .
\end{equation*}
We use these notations throughout Galois context of the thesis.
We denote
$\mathbb{SMC}_{\cat{colax}}$
as the $2$-category 
of symmetric monoidal categories,
colax symmetric monoidal functors
and comonoidal natural transformations.

Similarly, we denote a lax symmetric monoidal functor
as $\phi:\kos{S}\to \kos{T}$
whose underlying functor is also denoted as $\phi:\CS\to \CT$.
The lax symmetric monoidal coherence morphisms are denoted as
\begin{equation*}
  \phi_{\text{$\pzc{X}$},\text{$\pzc{Y}$}}
  :\phi(\text{$\pzc{X}$})\tensor \phi(\text{$\pzc{Y}$})
  \to \phi(\text{$\pzc{X}$}\ctimes \text{$\pzc{Y}$})
  ,
  \qquad
  \phi_{\text{$\pzc{1}$}}:
  \unit\to \phi(\text{$\pzc{1}$})
  ,
  \qquad
  \text{$\pzc{X}$}, \text{$\pzc{Y}$}\in \obj{\CS}
  .
\end{equation*}
These notations are used throughout Grothendieck context of the thesis.
We denote
$\mathbb{SMC}_{\cat{lax}}$
as the $2$-category 
of symmetric monoidal categories,
lax symmetric monoidal functors
and monoidal natural transformations.

\subsection{Acknowledgement}
I have been fascinated by the vision of my thesis advisor Jae-Suk Park that there should a single framework, in the name of kosmic Grothendieck-Galois theory, 
which not only vastly generalizes and unifies all the pre-existing dualities between "group-like objects" and the "abstract" categories of their internal representations 
in the broadest possible terms as the practice of the "Grothendieckian" yoga for the $\otimes$-categories and the fibre functors, but also serves as a bridge that connects different fields 
of mathematics.  This thesis contains  a small portion of the results obtained by an intensive collaborations with him for the last six years.  
We are preparing a series of papers "Kosmic Grothendieck-Galois theory I, II \& III" 
on the foundation of kosmic Grothendieck-Galois theory with some applications, and
this thesis is essentially a detailed exposition of the first four sections of a prequel "Pre-Kosmic  Grothendieck-Galois theory",
where both the relative and equivariant versions of pre-kosmic  Grothendieck/Galois categories as well
as the pre-kosmic descent theory along suitable changes of base pre-kosmoi will be presented.

\newpage
\section{Colax $\kos{K}$-tensor categories}
\label{sec ColaxKtensorCat}
Throughout \textsection~\ref{sec ColaxKtensorCat},
we fix a symmetric monoidal category
$\kos{K}=(\CK,\otimes,\kappa)$.
We study about
symmetric monoidal categories
$\kos{T}=(\CT,\tensor,\unit)$
equipped with colax symmetric monoidal functors
$\mon{t}:\kos{K}\to\kos{T}$,
and about the induced actions of $\kos{K}$ on $\kos{T}$.

\begin{definition} \label{def ColaxKtensorCat}
  We define the $2$-category
  \begin{equation*}
    \mathbb{SMC}_{\cat{colax}}^{\kos{K}/\!\!/}
  \end{equation*}
  of
  colax $\kos{K}$-tensor categories,
  colax $\kos{K}$-tensor functors
  and
  comonoidal $\kos{K}$-tensor natural transformations
  as the colax coslice $2$-category 
  of $\mathbb{SMC}_{\cat{colax}}$ under $\kos{K}$.
\end{definition}

Let us explain Definition~\ref{def ColaxKtensorCat}
in detail.
A \emph{colax $\kos{K}$-tensor category}
$(\kos{T},\mon{t})$ is a pair
of a symmetric monoidal category
$\kos{T}=(\CT,\tensor,\unit)$
and a colax symmetric monoidal functor
$\mon{t}:\kos{K}\to\kos{T}$.
We denote the comonoidal coherence morphisms of $\mon{t}$ as
\begin{equation}
  \label{eq ColaxKtensorCat coherence}
  \mon{t}_{x,y}:
  \mon{t}(x\otimes y)
  \to
  \mon{t}(x)\tensor \mon{t}(y)
  ,
  \quad
  x,y\in\obj{\CK}
  ,
  \qquad
  \mon{t}_{\kappa}:
  \mon{t}(\kappa)\to \unit
  . 
\end{equation}
Let $(\kos{S},\mon{s})$ be another colax $\kos{K}$-tensor category
where $\text{$\kos{S}$}=(\CS,\ctimes,\pzc{1})$ is the underlying symmetric monoidal category.
A \emph{colax $\kos{K}$-tensor functor}
$(\phi,\what{\phi}):(\kos{S},\mon{s})\to (\kos{T},\mon{t})$
is a pair of a colax symmetric monoidal functor
$\phi:\kos{S}\to\kos{T}$
and a comonoidal natural transformation
\begin{equation*}
  \vcenter{\hbox{
    \xymatrix@C=40pt{
      \text{ }
      &\kos{S}
      \ar[d]^-{\phi}
      \\
      \kos{K}
      \ar@/^0.7pc/[ur]^-{\mon{s}}
      \ar[r]_-{\mon{t}}
      \xtwocell[r]{}<>{<-2.5>{\text{ }\text{ }\what{\phi}}}
      &\kos{T}
    }
  }}
  \qquad\quad
  \what{\phi}:\phi\mon{s}\Rightarrow\mon{t}:\kos{K}\to \kos{S}
  .
\end{equation*}
Let $(\psi,\what{\psi}):(\kos{S},\mon{s})\to (\kos{T},\mon{t})$
be another colax $\kos{K}$-tensor functor.
A \emph{comonoidal $\kos{K}$-tensor natural transformation}
$\vartheta:(\phi,\what{\phi})\Rightarrow (\psi,\what{\psi})
:(\kos{S},\mon{s})\to (\kos{T},\mon{t})$
is a comonoidal natural transformation
$\vartheta:\phi\Rightarrow\psi:\kos{S}\to\kos{T}$
which satisfies the relation
\begin{equation}
  \label{eq ColaxKtensorCat KtensorNat}
  \vcenter{\hbox{
    \xymatrix@C=20pt{
      \phi\mon{s}
      \ar@2{->}[rr]^-{\vartheta\mon{s}}
      \ar@2{->}[dr]_-{\what{\phi}}
      &\text{ }
      &\psi\mon{s}
      \ar@2{->}[dl]^-{\what{\psi}}
      \\
      \text{ }
      &\mon{t}
      &\text{ }
    }
  }}
  \qquad\quad
  \what{\phi}=\what{\psi}\circ (\vartheta\mon{s}):\phi\mon{s}\Rightarrow\mon{t}
  .
\end{equation}
We may also describe the relation
(\ref{eq ColaxKtensorCat KtensorNat})
as the following diagram.
\begin{equation*}
  \vcenter{\hbox{
    \xymatrix@C=40pt{
      \text{ }
      &\kos{S}
      \ar[d]^-{\phi}
      \\
      \kos{K}
      \ar@/^0.7pc/[ur]^-{\mon{s}}
      \ar[r]_-{\mon{t}}
      \xtwocell[r]{}<>{<-2.5>{\text{ }\text{ }\what{\phi}}}
      &\kos{T}
    }
  }}
  \quad=
  \vcenter{\hbox{
    \xymatrix@C=15pt{
      \text{ }
      &\text{ }
      &\text{ }
      &\kos{S}
      \ar@/^1.1pc/[d]^-{\phi}
      \ar@/_1.1pc/[d]_-{\psi}
      \xtwocell[d]{}<>{<0>{\vartheta}}
      \\
      \kos{K}
      \ar@/^0.8pc/[urrr]^-{\mon{s}}
      \ar[rrr]_-{\mon{t}}
      \xtwocell[rr]{}<>{<-2>{\text{ }\text{ }\what{\psi}}}
      &\text{ }
      &\text{ }
      &\kos{T}
    }
  }}
\end{equation*}

We introduce a special kind of colax $\kos{K}$-tensor categories.

\begin{definition}
  A colax $\kos{K}$-tensor category
  $(\kos{T},\mon{t})$
  is called a \emph{strong $\kos{K}$-tensor category}
  if $\mon{t}:\kos{K}\to \kos{T}$
  is a strong $\kos{K}$-tensor functor,
  i.e., the comonoidal coherence morphisms
  in (\ref{eq ColaxKtensorCat coherence})
  are isomorphisms.
\end{definition}

For instance, the pair $(\kos{K},\id_{\kos{K}})$
of $\kos{K}$ and the identity functor $\id_{\kos{K}}:\kos{K}\to \kos{K}$
is a strong $\kos{K}$-tensor category.

We also introduce a special kind of colax $\kos{K}$-tensor functors.

\begin{definition}
  Let 
  $(\kos{T},\mon{t})$, $(\kos{S},\mon{s})$
  be colax $\kos{K}$-tensor categories.
  A colax $\kos{K}$-tensor functor
  $(\phi,\what{\phi}):(\kos{S},\mon{s})\to (\kos{T},\mon{t})$
  is called a \emph{strong $\kos{K}$-tensor functor}
  if
  $\phi:\kos{S}\to\kos{T}$ is a strong symmetric monoidal functor
  and 
  $\what{\phi}:\phi\mon{s}\Rightarrow\mon{t}$
  is a comonoidal natural isomorphism.
\end{definition}

Let $(\kos{T},\mon{t})$ be a colax $\kos{K}$-tensor category.
The pair of the colax symmetric monoidal functor
$\mon{t}:\kos{K}\to \kos{T}$
and the identity natural transformation
$I_{\mon{t}}$ of $\mon{t}$
is a colax $\kos{K}$-tensor functor
$(\text{$\mon{t}$},I_{\text{$\mon{t}$}}):(\text{$\kos{K}$},\id_{\kos{K}})\to (\kos{T},\mon{t})$.
It is a strong $\kos{K}$-tensor functor
if and only if $(\kos{T},\mon{t})$ is a strong $\kos{K}$-tensor category.

We introduce the action of $\kos{K}$
associated to each colax $\kos{K}$-tensor category.

\begin{definition}
  Let $(\kos{T},\mon{t})$
  be a colax $\kos{K}$-tensor category
  where $\kos{T}=(\CT,\tensor,\unit)$
  is the underlying symmetric monoidal category.
  We define the functor
  \begin{equation*}
    \bar{\acts}\!_{\mon{t}}:\CK\times \CT\to \CT,
    \qquad
    z\bar{\acts}\!_{\mon{t}} X
    :=
    X\tensor \mon{t}(z)
    ,
    \qquad
    z\in \obj{\CK}
    ,
    \qquad
    X\in \obj{\CT}
  \end{equation*}
  and call it as the \emph{associated (colax) $\kos{K}$-action on $\kos{T}$}.
\end{definition}

We introduce how a colax $\kos{K}$-tensor functor
between colax $\kos{K}$-tensor categories
interacts with the associated $\kos{K}$-actions.

\begin{definition}
  Let $(\kos{T},\mon{t})$, $(\kos{S},\mon{s})$
  be colax $\kos{K}$-tensor categories
  where $\kos{T}=(\CT,\tensor,\unit)$,
  $\text{$\kos{S}$}=(\CS,\ctimes,\pzc{1})$
  are underlying symmetric monoidal categories.
  For each colax $\kos{K}$-tensor functor
  $(\phi,\what{\phi}):(\kos{S},\mon{s})\to (\kos{T},\mon{t})$,
  we define the natural transformation
  \begin{equation*}
    \vecar{\phi}_{\slot,\slot}:
    \phi(\slot\bar{\acts}\!_{\mon{s}}\slot)
    \Rightarrow
    \slot\bar{\acts}\!_{\mon{t}}\phi(\slot)
    :\CK\times \CS\to \CT
  \end{equation*}
  which we call as the \emph{associated (colax) $\kos{K}$-equivariance of $(\phi,\what{\phi})$},
  whose component at $z\in \obj{\CK}$, $\pzc{X}\in \obj{\CS}$ is
  described below.
  \begin{equation*}
    \xymatrix@C=40pt{
      \phi(z\bcts\!_{\mon{s}} \pzc{X})
      \ar[rr]^-{\vecar{\phi}_{z,\pzc{X}}}
      \ar@{=}[d]
      &\text{ }
      &z\bcts\!_{\mon{t}} \phi(\pzc{X})
      \ar@{=}[d]
      \\
      \phi(\text{$\pzc{X}$}\ctimes \mon{s}(z))
      \ar[r]^-{\phi_{\pzc{X},\mon{s}(z)}}
      &\phi(\text{$\pzc{X}$})\tensor \phi\mon{s}(z)
      \ar[r]^-{I_{\phi(\pzc{X})}\tensor \what{\phi}_z}
      &\phi(\text{$\pzc{X}$})\tensor \mon{t}(z)
    }
  \end{equation*}
\end{definition}

In particular, the $\kos{K}$-equivariance
$\vecar{\phi}$
associated to a strong $\kos{K}$-tensor functor
$(\phi,\what{\phi})$
between colax $\kos{K}$-tensor categories
is a natural isomorphism.

\begin{lemma} \label{lem ColaxKtensorCat whatphi via vecarphi}
  Let 
  $(\phi,\what{\phi}):(\kos{S},\mon{s})\to (\kos{T},\mon{t})$
  be a colax $\kos{K}$-tensor functor
  between colax $\kos{K}$-tensor categories.
  The given comonoidal natural transformation
  $\what{\phi}:\phi\mon{s}\Rightarrow\mon{t}$
  is described in terms of the associated $\kos{K}$-equivariance
  $\vecar{\phi}$ as follows.
  Let $z\in\obj{\CK}$.
  \begin{equation}\label{eq ColaxKtensorCat whatphi via vecarphi}
    \what{\phi}_z:
    \xymatrix{
      \phi\mon{s}(z)
      \ar[r]^-{\phi(\imath_{\mon{s}(z)})}_-{\cong}
      &\phi(\text{$\pzc{1}$}\ctimes \mon{s}(z))
      \ar[r]^-{\vecar{\phi}_{z,\pzc{1}}}
      &\phi(\text{$\pzc{1}$})\tensor \mon{t}(z)
      \ar[r]^-{\phi_{\text{$\pzc{1}$}}\tensor I_{\mon{t}(z)}}
      &\unit\tensor \mon{t}(z)
      \ar[r]^-{\imath_{\mon{t}(z)}^{-1}}_-{\cong}
      &\mon{t}(z)
    }
  \end{equation}
\end{lemma}
\begin{proof}
  We can check the relation (\ref{eq ColaxKtensorCat whatphi via vecarphi}) as follows.
  \begin{equation*}
    \vcenter{\hbox{
      \xymatrix@C=40pt{
        \phi\mon{s}(z)
        \ar[d]_-{\phi(\imath_{\mon{s}(z)})}^-{\cong}
        \ar@{=}[rr]
        &\text{ }
        &\phi\mon{s}(z)
        \ar@{=}[dddd]
        \\
        \phi(\text{$\pzc{1}$}\ctimes \mon{s}(z))
        \ar[dd]_-{\vecar{\phi}_{z,\text{$\pzc{1}$}}}
        \ar@/^0.5pc/[dr]|-{\phi_{\text{$\pzc{1}$},\mon{s}(z)}}
        &\text{ }
        &\text{ }
        \\
        \text{ }
        &\phi(\text{$\pzc{1}$})\tensor \phi\mon{s}(z)
        \ar[d]^-{\phi_{\text{$\pzc{1}$}}\tensor I_{\phi\mon{s}(z)}}
        \ar@/^0.5pc/[dl]|-{I_{\phi(\text{$\pzc{1}$})}\tensor \what{\phi}_z}
        &\text{ }
        \\
        \phi(\text{$\pzc{1}$})\tensor \mon{t}(z)
        \ar[d]_-{\phi_{\text{$\pzc{1}$}}\tensor I_{\mon{t}(z)}}
        &\unit\tensor \phi\mon{s}(z)
        \ar@/^0.5pc/[dl]|-{I_{\text{$\pzc{1}$}}\tensor \what{\phi}_z}
        \ar@/_0.5pc/[dr]_-{\imath^{-1}_{\phi\mon{s}(z)}}^-{\cong}
        &\text{ }
        \\
        \unit\tensor \mon{t}(z)
        \ar[d]_-{\imath_{\mon{t}(z)}^{-1}}^-{\cong}
        &\text{ }
        &\phi\mon{s}(z)
        \ar[d]^-{\what{\phi}_z}
        \\
        \mon{t}(z)
        \ar@{=}[rr]
        &\text{ }
        &\mon{t}(z)
      }
    }}
  \end{equation*}
  This completes the proof of Lemma~\ref{lem ColaxKtensorCat whatphi via vecarphi}.
\qed\end{proof}

The following lemma describes the $\kos{K}$-equivariances
associated to compositions of colax $\kos{K}$-tensor functors.

\begin{lemma} \label{lem ColaxKtensorCat compositionKequiv}
  Let $(\kos{T},\mon{t})$, $(\kos{T}^{\pr},\mon{t}^{\pr})$, $(\kos{S},\mon{s})$
  be colax $\kos{K}$-tensor categories
  and let
  \begin{equation*}
    \vcenter{\hbox{
      \xymatrix@C=30pt{
        (\kos{S},\mon{s})
        \ar[r]^-{(\phi,\what{\phi})}
        &(\kos{T},\mon{t})
        \ar[r]^-{(\psi,\what{\psi})}
        &(\kos{T}^{\pr},\mon{t}^{\pr})
      }
    }}
  \end{equation*}
  be colax $\kos{K}$-tensor functors.
  The $\kos{K}$-equivariance
  associated to the composition colax $\kos{K}$-tensor functor
  $(\psi\phi,\what{\psi\phi}):(\kos{S},\mon{s})\to (\kos{T}^{\pr},\mon{t}^{\pr})$
  is given as follows.
  Let $z\in \obj{\CK}$ and $\pzc{X}\in \obj{\CS}$.
  \begin{equation}\label{eq ColaxKtensorCat compositionKequiv}
    \vecar{\psi\phi}_{z,\text{$\pzc{X}$}}:
    \xymatrix@C=30pt{
      \psi\phi(z\bar{\acts}\!_{\text{$\mon{s}$}}\text{$\pzc{X}$})
      \ar[r]^-{\psi(\vecar{\phi}_{z,\text{$\pzc{X}$}})}
      &\psi(z\bar{\acts}\!_{\text{$\mon{t}$}}\phi(\text{$\pzc{X}$}))
      \ar[r]^-{\vecar{\psi}_{z,\text{$\phi(\pzc{X})$}}}
      &z\bar{\acts}\!_{\text{$\mon{t}^{\pr}$}}\psi\phi(\text{$\pzc{X}$})
    }
  \end{equation}
\end{lemma}
\begin{proof}
  The comonoidal natural transformation
  $\what{\psi\phi}$
  of the composition colax $\kos{K}$-tensor functor is given by
  $\what{\psi\phi}:\!\!
  \xymatrix@C=17pt{
    \psi\phi\mon{s}
    \ar@2{->}[r]^-{\psi\what{\phi}}
    &\psi\mon{t}
    \ar@2{->}[r]^-{\what{\psi}}
    &\mon{t}^{\pr}
    .
  }$
  We can check the relation (\ref{eq ColaxKtensorCat compositionKequiv}) as follows.
  \begin{equation*}
    \vcenter{\hbox{
      \xymatrix@C=60pt{
        \psi\phi(\text{$\pzc{X}$}\ctimes \mon{s}(z))
        \ar[dd]_-{\psi(\vecar{\phi}_{z,\text{$\pzc{X}$}})}
        \ar@{=}[r]
        &\psi\phi(\text{$\pzc{X}$}\ctimes \mon{s}(z))
        \ar[d]^-{\psi(\phi_{\text{$\pzc{X}$},\mon{s}(z)})}
        \ar@{=}[r]
        &\psi\phi(\text{$\pzc{X}$}\ctimes \mon{s}(z))
        \ar[dddd]^-{\vecar{\psi\phi}_{z,\text{$\pzc{X}$}}}
        \ar@/^1.8pc/[ddl]|-{(\psi\phi)_{\text{$\pzc{X}$},\mon{s}(z)}}
        \\
        \text{ }
        &\psi(\phi(\text{$\pzc{X}$})\tensor \phi\mon{s}(z))
        \ar@/_0.5pc/[dl]|-{\psi(I_{\phi(\text{$\pzc{X}$})}\tensor \what{\phi}_z)}
        \ar[d]^-{\psi_{\phi(\text{$\pzc{X}$}),\phi\mon{s}(z)}}
        &\text{ }
        \\
        \psi(\phi(\text{$\pzc{X}$})\tensor \mon{t}(z))
        \ar[dd]_-{\vecar{\psi}_{z,\phi(\text{$\pzc{X}$})}}
        \ar@/_0.5pc/[dr]|-{\psi_{\phi(\text{$\pzc{X}$}),\mon{t}(z)}}
        &\psi\phi(\text{$\pzc{X}$})\tensor\!^{\pr} \psi\phi\mon{s}(z)
        \ar[d]^-{I_{\psi\phi(\text{$\pzc{X}$})}\tensor\!^{\pr}\psi(\what{\phi}_z)}
        \ar@/^1.8pc/[ddr]|-{I_{\psi\phi(\text{$\pzc{X}$})}\tensor\!^{\pr}\what{\psi\phi}_z}
        &\text{ }
        \\
        \text{ }
        &\psi\phi(\text{$\pzc{X}$})\tensor\!^{\pr}\psi\mon{t}(z)
        \ar[d]^-{I_{\psi\phi(\text{$\pzc{X}$})}\tensor\!^{\pr}\what{\psi}_z}
        &\text{ }
        \\
        \psi\phi(\text{$\pzc{X}$})\tensor\!^{\pr}\mon{t}^{\pr}(z)
        \ar@{=}[r]
        &\psi\phi(\text{$\pzc{X}$})\tensor\!^{\pr}\mon{t}^{\pr}(z)
        \ar@{=}[r]
        &\psi\phi(\text{$\pzc{X}$})\tensor\!^{\pr}\mon{t}^{\pr}(z)
      }
    }}
  \end{equation*}
  This completes the proof of Lemma~\ref{lem ColaxKtensorCat compositionKequiv}.
\qed\end{proof}

The following lemma states that
a comonoidal $\kos{K}$-tensor natural transformation
between colax $\kos{K}$-tensor functors
is precisely
a comonoidal natural transformation
that is compatible with associated $\kos{K}$-equivariances.

\begin{lemma}
  \label{lem ColaxKtensorCat KtensorNat preserves Kaction}
  Let $(\phi,\what{\phi})$,
  $(\psi,\what{\psi}):(\kos{S},\mon{s})\to (\kos{T},\mon{t})$
  be colax $\kos{K}$-tensor functors
  between colax $\kos{K}$-tensor categories
  and let
  $\vartheta:\phi\Rightarrow\psi:\kos{S}\to \kos{T}$
  be a comonoidal natural transformation.
  Then
  $\vartheta:(\phi,\what{\phi})\Rightarrow (\psi,\what{\psi})
  :(\kos{S},\mon{s})\to (\kos{T},\mon{t})$
  is a comonoidal $\kos{K}$-tensor natural transformation
  if and only if
  the following relation
  holds for all
  $z\in \obj{\CK}$, $\pzc{X}\in\obj{\CS}$.
  \begin{equation}\label{eq ColaxKtensorCat KtensorNat preserves Kaction}
    \vcenter{\hbox{
      \xymatrix@C=50pt{
        \phi(z\bar{\acts}\!_{\text{$\mon{s}$}}\text{$\pzc{X}$})
        \ar[d]_-{\vartheta_{z\bar{\acts}\!_{\text{$\mon{s}$}}\text{$\pzc{X}$}}}
        \ar[r]^-{\vecar{\phi}_{z,\pzc{X}}}
        &z\bar{\acts}\!_{\text{$\mon{t}$}}\phi(\text{$\pzc{X}$})
        \ar[d]^-{I_z\bar{\acts}\!_{\text{$\mon{t}$}}\vartheta_{\text{$\pzc{X}$}}}
        \\
        \psi(z\bar{\acts}\!_{\text{$\mon{s}$}}\text{$\pzc{X}$})
        \ar[r]^-{\vecar{\psi}_{z,\pzc{X}}}
        &z\bar{\acts}\!_{\text{$\mon{t}$}}\psi(\text{$\pzc{X}$})
      }
    }}
  \end{equation}
\end{lemma}
\begin{proof}
  We first prove the only if part.
  Suppose that $\vartheta:(\phi,\what{\phi})\Rightarrow (\psi,\what{\psi})$
  is a comonoidal $\kos{K}$-tensor natural transformation.
  We obtain the relation (\ref{eq ColaxKtensorCat KtensorNat preserves Kaction})
  as follows.
  \begin{equation*}
    \vcenter{\hbox{
      \xymatrix@C=60pt{
        \phi(\text{$\pzc{X}$}\ctimes \mon{s}(z))
        \ar[dd]_-{\vartheta_{\text{$\pzc{X}$}\ctimes \mon{s}(z)}}
        \ar@{=}[r]
        &\phi(\text{$\pzc{X}$}\ctimes \mon{s}(z))
        \ar[d]^-{\phi_{\text{$\pzc{X}$},\mon{s}(z)}}
        \ar@{=}[r]
        &\phi(\text{$\pzc{X}$}\ctimes \mon{s}(z))
        \ar[dd]^-{\vecar{\phi}_{z,\text{$\pzc{X}$}}}
        \\
        \text{ }
        &\phi(\text{$\pzc{X}$})\tensor \phi\mon{s}(z)
        \ar[d]^-{\vartheta_{\text{$\pzc{X}$}}\tensor I_{\phi\mon{s}(z)}}
        \ar@/^0.5pc/[dr]|-{I_{\phi(\text{$\pzc{X}$})}\tensor \what{\phi}_z}
        \ar@<-3.5ex>@/_2.5pc/[dd]|(0.3){\vartheta_{\text{$\pzc{X}$}}\tensor \vartheta_{\mon{s}(z)}}
        &\text{ }
        \\
        \psi(\text{$\pzc{X}$}\ctimes \mon{s}(z))
        \ar[dd]_-{\vecar{\psi}_{z,\text{$\pzc{X}$}}}
        \ar@/_0.5pc/[dr]|-{\psi_{\text{$\pzc{X}$},\mon{s}(z)}}
        &\psi(\text{$\pzc{X}$})\tensor \phi\mon{s}(z)
        \ar[d]^-{I_{\psi(\text{$\pzc{X}$})}\tensor \vartheta_{\mon{s}(z)}}
        \ar@/^1pc/[ddr]|-{I_{\psi(\text{$\pzc{X}$})}\tensor \what{\phi}_z}
        &\phi(\text{$\pzc{X}$})\tensor \mon{t}(z)
        \ar[dd]^-{\vartheta_{\text{$\pzc{X}$}}\tensor I_{\mon{t}(z)}}
        \\
        \text{ }
        &\psi(\text{$\pzc{X}$})\tensor \psi\mon{s}(z)
        \ar[d]^-{I_{\psi(\text{$\pzc{X}$})}\tensor \what{\psi}_z}
        &\text{ }
        \\
        \psi(\text{$\pzc{X}$})\tensor \mon{t}(z)
        \ar@{=}[r]
        &\psi(\text{$\pzc{X}$})\tensor \mon{t}(z)
        \ar@{=}[r]
        &\psi(\text{$\pzc{X}$})\tensor \mon{t}(z)
      }
    }}
  \end{equation*}
  Next we prove the if part.
  Assume that $\vartheta:\phi\Rightarrow\psi$
  satisfies the relation (\ref{eq ColaxKtensorCat KtensorNat preserves Kaction})
  for all $z$, $\pzc{X}$.
  Using the descriptions
  of $\what{\phi}$, $\what{\psi}$
  given in (\ref{eq ColaxKtensorCat whatphi via vecarphi}),
  we can check that
  $\what{\phi}=\what{\psi}\circ (\vartheta\mon{s})$
  as follows.
  \begin{equation*}
    \vcenter{\hbox{
      \xymatrix{
        \phi\mon{s}(z)
        \ar[d]_-{\vartheta_{\mon{s}(z)}}
        \ar@{=}[rrrr]
        \ar@/^0.5pc/[dr]^(0.6){\phi(\imath_{\mon{s}(z)})}_-{\cong}
        &\text{ }
        &\text{ }
        &\text{ }
        &\phi\mon{s}(z)
        \ar[ddddd]^-{\what{\phi}_z}
        \\
        \psi\mon{s}(z)
        \ar@/_0.5pc/[dr]_-{\psi(\imath_{\mon{s}(z)})}^-{\cong}
        \ar[dddd]_-{\what{\psi}_z}
        &\phi(\text{$\pzc{1}$}\ctimes \mon{s}(z))
        \ar[d]^-{\vartheta_{\text{$\pzc{1}$}\ctimes \mon{s}(z)}}
        \ar@/^0.5pc/[dr]|-{\vecar{\phi}_{z,\text{$\pzc{1}$}}}
        &\text{ }
        &\text{ }
        &\text{ }
        \\
        \text{ }
        &\psi(\text{$\pzc{1}$}\ctimes \mon{s}(z))
        \ar@/_0.5pc/[dr]|-{\vecar{\psi}_{z,\text{$\pzc{1}$}}}
        &\phi(\text{$\pzc{1}$})\tensor \mon{t}(z)
        \ar[d]^-{\vartheta_{\text{$\pzc{1}$}}\tensor I_{\mon{t}(z)}}
        \ar@/^1.5pc/[ddr]|-{\phi_{\text{$\pzc{1}$}}\tensor I_{\mon{t}(z)}}
        &\text{ }
        &\text{ }
        \\
        \text{ }
        &\text{ }
        &\psi(\text{$\pzc{1}$})\tensor \mon{t}(z)
        \ar@/_0.5pc/[dr]|-{\psi_{\text{$\pzc{1}$}}\tensor I_{\mon{t}(z)}}
        &\text{ }
        &\text{ }
        \\
        \text{ }
        &\text{ }
        &\text{ }
        &\unit\tensor \mon{t}(z)
        \ar@/_0.5pc/[dr]^-{\imath_{\mon{t}(z)}^{-1}}_-{\cong}
        &\text{ }
        \\
        \mon{t}(z)
        \ar@{=}[rrrr]
        &\text{ }
        &\text{ }
        &\text{ }
        &\mon{t}(z)
      }
    }}
  \end{equation*}
  Therefore
  $\vartheta:(\phi,\what{\phi})\Rightarrow (\psi,\what{\psi})$
  is a comonoidal $\kos{K}$-tensor natural transformation.
  This completes the proof of Lemma~\ref{lem ColaxKtensorCat KtensorNat preserves Kaction}.
\qed\end{proof}

\subsection{Cocommutative comonoids as colax $\kos{K}$-tensor functors}
\label{subsec Ens(T)}
Let $\kos{T}=(\CT,\tensor,\unit)$
be a symmetric monoidal category.
We denote
\begin{equation*}
  \kos{E\!n\!s}(\kos{T})
  =(\cat{Ens}(\kos{T}),\tensor,\unit)
\end{equation*}
as the cartesian monoidal category
of cocommutative comonoids in $\kos{T}$.
An object in $\cat{Ens}(\kos{T})$
is a triple
$(C,\cp_C,e_C)$ 
where $C$ is an object in $\CT$ and 
$\cp_C:C\to C\tensor C$,
$e_C:C\to \unit$
are coproduct, counit morphisms in $\CT$
satisfying the coassociativity, counital, cocommutativity relations.
We often omit coproduct, counit morphisms
and simply denote an object $(C,\cp_C,e_C)$ in $\cat{Ens}(\kos{T})$ as $C$.
A morphism $C^{\pr}\to C$ in $\cat{Ens}(\kos{T})$
is a morphism in $\CT$
which is compatible with coproducts, counits morphisms of $C^{\pr}$, $C$.
The category $\cat{Ens}(\kos{T})$ has finite products.
The object $\unit$ equipped with
$\cp_{\unit}=\imath_{\unit}=\jmath_{\unit}:\unit\xrightarrow{\cong}\unit\tensor \unit$
and
$e_{\unit}=I_{\unit}:\unit\xrightarrow{\cong}\unit$
is a terminal object in $\cat{Ens}(\kos{T})$.
The product of objects $C$, $C^{\pr}$ in $\cat{Ens}(\kos{T})$
is $C\tensor C^{\pr}$ whose coproduct, counit morphisms are described below.
\begin{equation} \label{eq Ens(T) CC' definition}
  \begin{aligned}
    \cp_{C\tensor C^{\pr}}
    &:
    \xymatrix@C=45pt{
      C\tensor C^{\pr}
      \ar[r]^-{\cp_C\tensor \cp_{C^{\pr}}}
      &C\tensor C\tensor C^{\pr}\tensor C^{\pr}
      \ar[r]^-{I_C\tensor s_{C,C^{\pr}}\tensor I_{C^{\pr}}}_-{\cong}
      &C\tensor C^{\pr}\tensor C\tensor C^{\pr}
    }
    \\
    e_{C\tensor C^{\pr}}
    &:
    \xymatrix@C=30pt{
      C\tensor C^{\pr}
      \ar[r]^-{e_C\tensor e_{C^{\pr}}}
      &\unit\tensor \unit
      \ar[r]^-{\cp_{\unit}^{-1}}_-{\cong}
      &\unit
    }
  \end{aligned}
\end{equation}
The symmetric monoidal coherence isomorphisms
$a$, $\imath$, $\jmath$, $s$
of
$\kos{E\!n\!s}(\kos{T})$
are given by those of $\kos{T}$.

For each object $C$ in $\cat{Ens}(\kos{T})$,
we have a colax symmetric monoidal comonad 
\begin{equation*}
  \langle C\tensor\rangle=(C\tensor,\delta^{C\tensor},\varepsilon^{C\tensor})
\end{equation*}
on $\kos{T}$.  
The underlying colax symmetric monoidal endofunctor
$C\tensor:\kos{T}\to \kos{T}$
is the endofunctor
$C\tensor=C\tensor\slot:\CT\to \CT$
equipped with the following comonoidal coherence morphisms
\begin{equation}\label{eq Ens(T) Ctensor definition}
  \begin{aligned}
    (C\tensor)_{X,Y}
    &:
    \xymatrix@C=40pt{
      C\tensor X\tensor Y
      \ar[r]^-{\cp_C\tensor I_{X\tensor Y}}
      &C\tensor C\tensor X\tensor Y
      \ar[r]^-{I_C\tensor s_{C,X}\tensor I_Y}_-{\cong}
      &C\tensor X\tensor C\tensor Y
    }
    \\
    (C\tensor)_{\unit}
    &:
    \xymatrix@C=30pt{
      C\tensor \unit
      \ar[r]^-{e_C\tensor I_{\unit}}
      &\unit\tensor \unit
      \ar[r]^-{\cp_{\unit}^{-1}}_-{\cong}
      &\unit
    }
    \qquad
    X,Y\in\obj{\CT}
  \end{aligned}
\end{equation}
and the comonoidal natural transformations
$\delta^{C\tensor}$, $\varepsilon^{C\tensor}$
are described below.
\begin{equation}\label{eq Ens(T) deltaCepC definition}
  \begin{aligned}
    \delta^{C\tensor}_X
    &:
    \xymatrix@C=30pt{
      C\tensor X
      \ar[r]^-{\cp_C\tensor I_X}
      &(C\tensor C)\tensor X
      \ar[r]^-{a^{-1}_{C,C,X}}_-{\cong}
      &C\tensor (C\tensor X)
    }
    \\
    \varepsilon^{C\tensor}_X
    &:
    \xymatrix@C=30pt{
      C\tensor X
      \ar[r]^-{e_C\tensor I_X}
      &\unit\tensor X
      \ar[r]^-{\imath_X^{-1}}_-{\cong}
      &X
    }
    \qquad\quad
    X\in\obj{\CT}
  \end{aligned}
\end{equation}

\begin{proposition} \label{prop Ens(T) colaxK(K,T)Ens(T) adjunction}
  Let $(\kos{T},\mon{t})$ be a strong $\kos{K}$-tensor category.
  We have a reflective adjunction
  \begin{equation*}
    \CL:
    \mathbb{SMC}_{\cat{colax}}^{\kos{K}/\!\!/}\bigl((\kos{K},\id_{\kos{K}}),(\kos{T},\mon{t})\bigr)
    \rightleftarrows
    \cat{Ens}(\kos{T})
    :\iota
  \end{equation*}
  between
  the category of cocommutative comonoids in $\kos{T}$
  and
  the category of colax $\kos{K}$-tensor functors
  $(\kos{K},\id_{\kos{K}})\to (\kos{T},\mon{t})$.
  \begin{itemize}
    \item 
    The right adjoint $\iota$
    sends each object
    $C$ in $\cat{Ens}(\kos{T})$
    to the colax $\kos{K}$-tensor functor
    \begin{equation*}
      \iota(C)
      :=
      (C\tensor\text{$\mon{t}$},\what{C\tensor\mon{t}}):
      (\kos{K},\id_{\kos{K}})
      \to
      (\kos{T},\mon{t})
      .
    \end{equation*}
    The underlying colax symmetric monoidal functor is
    $C\tensor \text{$\mon{t}$}:
    \kos{K}
    \xrightarrow[]{\text{$\mon{t}$}}
    \kos{T}
    \xrightarrow[]{C\tensor}
    \kos{T}$
    and the component of the comonoidal natural transformation
    $\what{C\tensor\mon{t}}$ 
    at each object $x$ in $\CK$ is
    \begin{equation*}
      \vcenter{\hbox{
        \xymatrix@C=40pt{
          \text{ }
          &\kos{K}
          \ar[d]^-{C\tensor \mon{t}}
          \\
          \kos{K}
          \ar@/^0.7pc/[ur]^-{\id_{\kos{K}}}
          \ar[r]_-{\mon{t}}
          \xtwocell[r]{}<>{<-2.5>{\quad\text{ }\text{ }\what{C\tensor \mon{t}}}}
          &\kos{T}
        }
      }}
      \qquad
      \what{C\tensor \mon{t}}_x
      =\varepsilon^{C\tensor}_{\mon{t}(x)}
      :\!\!
      \xymatrix@C=25pt{
        C\tensor \mon{t}(x)
        \ar[r]^-{e_C\tensor I_{\mon{t}(x)}}
        &\unit\tensor \mon{t}(x)
        \ar[r]^-{\imath^{-1}_{\mon{t}(x)}}_-{\cong}
        &\mon{t}(x)
        .
      }
    \end{equation*}

    \item 
    The left adjoint $\CL$ sends each colax $\kos{K}$-tensor functor
    $(\phi,\what{\phi}):(\kos{K},\id_{\kos{K}})\to (\kos{T},\mon{t})$
    to the object
    $\CL(\phi,\what{\phi}):=\phi(\kappa)$
    in $\cat{Ens}(\kos{T})$,
    where
    \begin{equation*}
      \cp_{\phi(\kappa)}:\!
      \xymatrix{
        \phi(\kappa)
        \ar[r]^-{\phi(\cp_{\kappa})}_-{\cong}
        &\phi(\kappa\otimes\kappa)
        \ar[r]^-{\phi_{\kappa,\kappa}}
        &\phi(\kappa)\tensor \phi(\kappa)
        ,
      }
      \quad
      e_{\phi(\kappa)}:\!
      \xymatrix{
        \phi(\kappa)
        \ar[r]^-{\phi_{\kappa}}
        &\unit
        .
      }
    \end{equation*}

    \item 
    The component of the adjunction counit at each object
    $C$ in $\cat{Ens}(\kos{T})$ is
    \begin{equation*}
      \xymatrix@C=30pt{
        C\tensor \mon{t}(\kappa)
        \ar[r]^-{I_C\tensor \mon{t}_{\kappa}}_-{\cong}
        &C\tensor \unit
        \ar[r]^-{\jmath_C^{-1}}_-{\cong}
        &C
        .
      }
    \end{equation*}

    \item
    The component of the adjunction unit at each
    colax $\kos{K}$-tensor functor
    $(\phi,\what{\phi}):(\kos{K},\id_{\kos{K}})\to (\kos{T},\mon{t})$
    is the comonoidal $\kos{K}$-tensor natural transformation
    \begin{equation*}
      \hatar{\phi}
      :
      (\phi,\what{\phi})
      \Rightarrow
      (\phi(\kappa)\tensor\text{$\mon{t}$},\what{\phi(\kappa)\tensor\mon{t}})
      :(\kos{K},\id_{\kos{K}})\to (\kos{T},\mon{t})
    \end{equation*}
    whose component at $x\in\obj{\CK}$ is
    \begin{equation}
      \label{eq Ens(T) colaxK(K,T)Ens(T) adjunction}
      \hatar{\phi}_x:
      \xymatrix@C=25pt{
        \phi(x)
        \ar[r]^-{\phi(\imath_x)}_-{\cong}
        &\phi(\kappa\otimes x)
        \ar[r]^-{\phi_{\kappa,x}}
        \ar@/_1pc/@<-1ex>[rr]|-{\vecar{\phi}_{x,\kappa}}
        &\phi(\kappa)\tensor \phi(x)
        \ar[r]^-{I_{\phi(\kappa)}\tensor \what{\phi}_x}
        &\phi(\kappa)\tensor \mon{t}(x)
        .
      }
    \end{equation}
    In particular, the component of $\hatar{\phi}$ at $\kappa$ is
    \begin{equation}
      \label{eq2 Ens(T) colaxK(K,T)Ens(T) adjunction}
      \hatar{\phi}_{\kappa}:
      \xymatrix@C=30pt{
        \phi(\kappa)
        \ar[r]^-{\jmath_{\phi(\kappa)}}_-{\cong}
        &\phi(\kappa)\tensor \unit
        \ar[r]^-{I_{\phi(\kappa)}\tensor \mon{t}_{\kappa}^{-1}}_-{\cong}
        &\phi(\kappa)\tensor \mon{t}(\kappa).          
      }
    \end{equation}
  \end{itemize}
\end{proposition}
\begin{proof}
  We first check that the right adjoint functor $\iota$ is well-defined.
  For each object $C$ in $\cat{Ens}(\kos{T})$,
  the colax $\kos{K}$-tensor functor
  $\iota(C)=(C\tensor\text{$\mon{t}$},\what{C\tensor \mon{t}})
  :(\kos{K},\id_{\kos{K}})\to (\kos{T},\mon{t})$
  is well-defined.
  For each morphism $f:C^{\pr}\to C$ in $\cat{Ens}(\kos{T})$,
  the comonoidal natural transformation
  $f\tensor \mon{t}:C^{\pr}\tensor\mon{t}\Rightarrow C\tensor \mon{t}$
  becomes a comonoidal $\kos{K}$-tensor natural transformation
  $\iota(f):\iota(C^{\pr})\Rightarrow\iota(C)$
  as we can see from the diagram below.
  Let $x$ be an object in $\CK$.
  \begin{equation*}
    \vcenter{\hbox{
      \xymatrix@C=20pt{
        C^{\pr}\tensor \mon{t}
        \ar@2{->}[rr]^-{f\tensor \mon{t}}
        \ar@2{->}[dr]_-{\what{C^{\pr}\tensor \mon{t}}}
        &\text{ }
        &C\tensor \mon{t}
        \ar@2{->}[dl]^-{\what{C\tensor \mon{t}}}
        \\
        \text{ }
        &\mon{t}
        &\text{ }
      }
    }}
    \qquad
    \vcenter{\hbox{
      \xymatrix@C=15pt{
        C^{\pr}\tensor \mon{t}(x)
        \ar[d]_-{f\tensor I_{\mon{t}(x)}}
        \ar@{=}[rr]
        &\text{ }
        &C^{\pr}\tensor \mon{t}(x)
        \ar@/_1pc/[ddl]|-{e_{C^{\pr}}\tensor I_{\mon{t}(x)}}
        \ar[ddd]^-{\what{C^{\pr}\tensor\mon{t}}_x}
        \\
        C\tensor \mon{t}(x)
        \ar[dd]_-{\what{C\tensor\mon{t}}_x}
        \ar@/^0.5pc/[dr]|-{e_C\tensor I_{\mon{t}(x)}}
        &\text{ }
        &\text{ }
        \\
        \text{ }
        &\unit\tensor \mon{t}(x)
        \ar@/^0.5pc/[dr]_(0.4){\imath^{-1}_{\mon{t}(x)}}^(0.4){\cong}
        &\text{ }
        \\
        \mon{t}(x)
        \ar@{=}[rr]
        &\text{ }
        &\mon{t}(x)
      }
    }}
  \end{equation*}
  This shows that the right adjoint functor $\iota$ is well-defined.

  One can check that 
  for each colax $\kos{K}$-tensor functor
  $(\phi,\what{\phi}):(\kos{K},\id_{\kos{K}})\to (\kos{T},\mon{t})$
  the object 
  $\CL(\phi,\what{\phi})=\phi(\kappa)$
  in $\cat{Ens}(\kos{T})$
  is well-defined,
  and for each comonoidal $\kos{K}$-tensor natural transformation
  $\vartheta:
  (\phi,\what{\phi})\Rightarrow(\psi,\what{\psi})
  :(\kos{K},\id_{\kos{K}})\to (\kos{T},\mon{t})$
  its component at $\kappa$ is a morphism
  $\CL(\vartheta)=\vartheta_{\kappa}:\phi(\kappa)\to \psi(\kappa)$
  in $\cat{Ens}(\kos{T})$.
  Thus the left adjoint $\CL$ is also well-defined.

  Let $C$ be an object in $\cat{Ens}(\kos{T})$.
  The isomorphism
  $C\tensor \mon{t}(\kappa)
  \xrightarrow[\cong]{I_C\tensor \mon{t}_{\kappa}}
  C\tensor \unit
  \xrightarrow[\cong]{\jmath_C^{-1}}
  C$
  in $\cat{Ens}(\kos{T})$ is well-defined
  and is natural in variable $C$.
  This shows that the adjunction counit is well-defined.

  We show that the adjunction unit is well-defined in two steps.
  First step is to show that 
  for each colax $\kos{K}$-tensor functor
  $(\phi,\what{\phi}):(\kos{K},\id_{\kos{K}})\to (\kos{T},\mon{t})$,
  we have the comonoidal $\kos{K}$-tensor natural transformation
  $\hatar{\phi}:(\phi,\what{\phi})\Rightarrow(\phi(\kappa)\tensor \mon{t},\what{\phi(\kappa)\tensor \mon{t}})$
  as we claimed in (\ref{eq Ens(T) colaxK(K,T)Ens(T) adjunction}).
  We leave for the readers to check that
  $\hatar{\phi}:\phi\Rightarrow\phi(\kappa)\tensor \mon{t}:\kos{K}\to \kos{T}$
  is a comonoidal natural transformation.
  We have the relation
  \begin{equation*}
    \vcenter{\hbox{
      \xymatrix{
        \phi
        \ar@2{->}[rr]^-{\hatar{\phi}}
        \ar@2{->}[dr]_-{\what{\phi}}
        &\text{ }
        &\phi(\kappa)\tensor \mon{t}
        \ar@2{->}[dl]^-{\what{\phi(\kappa)\tensor\mon{t}}}
        \\
        \text{ }
        &\mon{t}
        &\text{ }
      }
    }}
  \end{equation*}
  as we can see from the diagram below.
  Let $x$ be an object in $\CK$.
  \begin{equation*}
    \xymatrix@C=70pt{
      \phi(x)
      \ar[ddd]_-{\hatar{\phi}_x}
      \ar@{=}[r]
      &\phi(x)
      \ar[d]^-{\phi(\imath_x)}_-{\cong}
      \ar@{=}[r]
      &\phi(x)
      \ar@{=}[dddd]
      \\
      \text{ }
      &\phi(\kappa\otimes x)
      \ar[d]^-{\phi_{\kappa,x}}
      &\text{ }
      \\
      \text{ }
      &\phi(\kappa)\tensor \phi(x)
      \ar@/_0.5pc/[dl]|-{I_{\phi(\kappa)}\tensor \what{\phi}_x}
      \ar[d]^-{\phi_{\kappa}\tensor I_{\phi(x)}}
      &\text{ }
      \\
      \phi(\kappa)\tensor \mon{t}(x)
      \ar[dd]_-{\what{\phi(\kappa)\tensor\mon{t}}_x}
      \ar@/_0.5pc/[dr]|-{\phi_{\kappa}\tensor I_{\mon{t}(x)}}
      &\unit\tensor \phi(x)
      \ar[d]^-{I_{\unit}\tensor \what{\phi}_x}
      \ar@/^0.5pc/[dr]^-{\imath^{-1}_{\phi(x)}}_-{\cong}
      &\text{ }
      \\
      \text{ }
      &\unit\tensor \mon{t}(x)
      \ar[d]^-{\imath^{-1}_{\mon{t}(x)}}_-{\cong}
      &\phi(x)
      \ar[d]^-{\what{\phi}_x}
      \\
      \mon{t}(x)
      \ar@{=}[r]
      &\mon{t}(x)
      \ar@{=}[r]
      &\mon{t}(x)
    }
  \end{equation*}
  This shows that 
  the comonoidal $\kos{K}$-tensor natural transformation
  \begin{equation*}
    \hatar{\phi}:(\phi,\what{\phi})\Rightarrow
    (\phi(\kappa)\tensor \mon{t},\what{\phi(\kappa)\tensor \mon{t}})
  \end{equation*}
  is well-defined.
  Second step is to show that
  $\hatar{\phi}$
  is natural in variable $(\phi,\what{\phi})$.
  Let
  $(\psi,\what{\psi}):(\kos{K},\id_{\kos{K}})\to (\kos{T},\mon{t})$
  be another colax $\kos{K}$-tensor functor
  and let
  $\vartheta:(\phi,\what{\phi})\Rightarrow (\psi,\what{\psi})$
  be a comonoidal $\kos{K}$-tensor natural transformation.
  Recall that we have $\what{\phi}=\what{\psi}\circ \vartheta
  :\phi\Rightarrow\mon{t}$.
  We need to verify the relation
  \begin{equation*}
    \vcenter{\hbox{
      \xymatrix@R=30pt@C=40pt{
        \phi
        \ar@2{->}[d]_-{\vartheta}
        \ar@2{->}[r]^-{\hatar{\phi}}
        &\phi(\kappa)\tensor \mon{t}
        \ar@2{->}[d]^-{\vartheta_{\kappa}\tensor \mon{t}}
        \\
        \psi
        \ar@2{->}[r]^-{\hatar{\psi}}
        &\psi(\kappa)\tensor \mon{t}
      }
    }}
    \qquad\quad
    (\vartheta_{\kappa}\tensor \mon{t})\circ \hatar{\phi}
    =
    \hatar{\psi}\circ\vartheta:
    \phi\Rightarrow\psi(\kappa)\tensor \mon{t}.
  \end{equation*}
  We can check this as follows.
  Let $x$ be an object in $\CK$.
  \begin{equation*}
    \vcenter{\hbox{
      \xymatrix@C=60pt{
        \phi(x)
        \ar[d]^-{\vartheta_x}
        \ar@{=}[r]
        &\phi(x)
        \ar[d]^-{\phi(\imath_x)}_-{\cong}
        \ar@/^3pc/@<5ex>[ddd]^-{\hatar{\phi}_x}
        \\
        \psi(x)
        \ar@/_3pc/@<-2ex>[ddd]_-{\hatar{\psi}_x}
        \ar[d]^-{\psi(\imath_x)}_-{\cong}
        &\phi(\kappa\otimes x)
        \ar@/^0.5pc/[dl]|-{\vartheta_{\kappa\otimes x}}
        \ar[d]^-{\phi_{\kappa,x}}
        \\
        \psi(\kappa\otimes x)
        \ar[d]^-{\psi_{\kappa,x}}
        &\phi(\kappa)\tensor\phi(x)
        \ar@/^0.5pc/[dl]|-{\vartheta_{\kappa}\tensor\vartheta_x}
        \ar[d]^-{I_{\phi(\kappa)}\tensor \what{\phi}_x}
        \\
        \psi(\kappa)\tensor \psi(x)
        \ar[d]^-{I_{\psi(\kappa)}\tensor \what{\psi}_x}
        &\phi(\kappa)\tensor \mon{t}(x)
        \ar[d]^-{\vartheta_{\kappa}\tensor I_{\mon{t}(x)}}
        \\
        \psi(\kappa)\tensor \mon{t}(x)
        \ar@{=}[r]
        &\psi(\kappa)\tensor \mon{t}(x)
      }
    }}
  \end{equation*}
  This shows that the adjunction unit is well-defined.

  We are left to show that the adjunction unit, counit
  satisfy the triangle identities.
  For each colax $\kos{K}$-tensor functor
  $(\phi,\what{\phi}):(\kos{K},\id_{\kos{K}})\to (\kos{T},\mon{t})$
  we verify one of the triangle identities,
  which is the description (\ref{eq2 Ens(T) colaxK(K,T)Ens(T) adjunction})
  of $\hatar{\phi}_{\kappa}$ that we claimed.
  \begin{equation*}
    \vcenter{\hbox{
      \xymatrix@C=30pt{
        \phi(\kappa)
        \ar@/_1.5pc/@{=}[ddr]
        \ar[r]^-{\hatar{\phi}_{\kappa}}
        &\phi(\kappa)\tensor \mon{t}(\kappa)
        \ar[d]^-{I_{\phi(\kappa)}\tensor \mon{t}_{\kappa}}_-{\cong}
        \\
        \text{ }
        &\phi(\kappa)\tensor \unit
        \ar[d]^-{\jmath^{-1}_{\phi(\kappa)}}_-{\cong}
        \\
        \text{ }
        &\phi(\kappa)
      }
    }}
    \quad
    \vcenter{\hbox{
      \xymatrix@C=30pt{
        \phi(\kappa)
        \ar[d]^-{\phi(\cp_{\kappa})}_-{\cong}
        \ar@/_2pc/@<-2ex>[ddd]_-{\hatar{\phi}_{\kappa}}
        \ar@{=}[rr]
        &\text{ }
        &\phi(\kappa)
        \ar@/_1pc/[dddl]^-{\jmath_{\phi(\kappa)}}_-{\cong}
        \ar@{=}[ddd]
        \\
        \phi(\kappa\otimes \kappa)
        \ar[d]^-{\phi_{\kappa,\kappa}}
        &\text{ }
        &\text{ }
        \\
        \phi(\kappa)\tensor \phi(\kappa)
        \ar[d]^-{I_{\phi(\kappa)}\tensor \what{\phi}_{\kappa}}
        \ar@/^0.5pc/[dr]^-{I_{\phi(\kappa)}\tensor \phi_{\kappa}}
        &\text{ }
        &\text{ }
        \\
        \phi(\kappa)\tensor \mon{t}(\kappa)
        \ar[r]_-{I_{\phi(\kappa)}\tensor \mon{t}_{\kappa}}^-{\cong}
        &\phi(\kappa)\tensor \unit
        \ar[r]_-{\jmath_{\phi(\kappa)}^{-1}}^-{\cong}
        &\phi(\kappa)
      }
    }}
  \end{equation*}
  Let $C$ be an object in $\cat{Ens}(\kos{T})$.
  For each object $x$ in $\CK$,
  we can describe $\hatar{C\tensor\mon{t}}_x$ as
  \begin{equation*}
    \vcenter{\hbox{
      \xymatrix@C=40pt{
        C\tensor \mon{t}(x)
        \ar[dddd]_-{\hatar{C\tensor \mon{t}}_x}
        \ar@{=}[rr]
        &\text{ }
        &C\tensor \mon{t}(x)
        \ar[d]^-{I_C\tensor \mon{t}(\imath_x)}_-{\cong}
        \\
        \text{ }
        &\text{ }
        &C\tensor \mon{t}(\kappa\otimes x)
        \ar@/_1.5pc/[ddl]_-{(C\tensor\mon{t})_{\kappa,x}}
        \ar[d]^-{I_C\tensor \mon{t}_{\kappa,x}}_-{\cong}
        \\
        \text{ }
        &\text{ }
        &C\tensor (\mon{t}(\kappa)\tensor \mon{t}(x))
        \ar@/^0.5pc/[dl]|-{(C\tensor)_{\mon{t}(\kappa),\mon{t}(x)}}
        \ar[dd]^-{a_{C,\mon{t}(\kappa),\mon{t}(x)}}_-{\cong}
        \\
        \text{ }
        &(C\tensor \mon{t}(\kappa))\tensor (C\tensor \mon{t}(x))
        \ar@/_0.5pc/[dl]_-{I_{C\tensor \mon{t}(x)}\tensor \what{C\tensor\mon{t}}_x}
        \ar@/^0.5pc/[dr]|-{I_{C\tensor \mon{t}(x)}\tensor \varepsilon^{C\tensor}_{\mon{t}(x)}}
        &\text{ }
        \\
        (C\tensor \mon{t}(\kappa))\tensor \mon{t}(x)
        \ar@{=}[rr]
        &\text{ }
        &(C\tensor \mon{t}(\kappa))\tensor \mon{t}(x)
      }
    }}
  \end{equation*}
  and using the above relation,
  we verify the other triangle identity as follows.
  \begin{equation*}
    \vcenter{\hbox{
      \xymatrix@C=20pt{
        C\tensor \mon{t}
        \ar@2{->}[r]^-{\hatar{C\tensor \mon{t}}}
        \ar@/_1.5pc/@{=}[ddr]
        &(C\tensor \mon{t}(\kappa))\tensor \mon{t}
        \ar@2{->}[d]^-{(I_{C}\tensor \mon{t}_{\kappa})\tensor I_{\mon{t}}}_-{\cong}
        \\
        \text{ }
        &(C\tensor \unit)\tensor \mon{t}
        \ar@2{->}[d]^-{\jmath^{-1}_{C}\tensor I_{\mon{t}}}_-{\cong}
        \\
        \text{ }
        &C\tensor \mon{t}
      }
    }}
    \quad
    \vcenter{\hbox{
      \xymatrix@C=40pt{
        C\tensor \mon{t}(x)
        \ar@/_3pc/@<-2ex>[ddd]|-{\hatar{C\tensor\mon{t}}_x}
        \ar[d]^-{I_C\tensor \mon{t}(\imath_x)}_-{\cong}
        \ar@{=}[r]
        &C\tensor \mon{t}(x)
        \ar[ddd]^-{I_C\tensor \imath_{\mon{t}(x)}}_-{\cong}
        \ar@{=}@/^3pc/@<3ex>[ddddd]
        \\
        C\tensor \mon{t}(\kappa\otimes x)
        \ar[d]^-{I_C\tensor \mon{t}_{\kappa,x}}_-{\cong}
        &\text{ }
        \\
        C\tensor (\mon{t}(\kappa)\tensor \mon{t}(x))
        \ar[d]^-{a_{C,\mon{t}(\kappa),\mon{t}(x)}}_-{\cong}
        \ar@/^1pc/[dr]^(0.5){I_C\tensor (\mon{t}_{\kappa}\tensor I_{\mon{t}(x)})}_-{\cong}
        &\text{ }
        \\
        (C\tensor \mon{t}(\kappa))\tensor \mon{t}(x)
        \ar[d]_-{(I_C\tensor \mon{t}_{\kappa})\tensor I_{\mon{t}(x)}}^-{\cong}
        &C\tensor (\unit \tensor \mon{t}(x))
        \ar@/^1pc/[dl]^(0.5){a_{C,\unit,\mon{t}(x)}}_(0.5){\cong}
        \ar[dd]^-{I_C\tensor \imath_{\mon{t}(x)}^{-1}}_-{\cong}
        \\
        (C\tensor \unit)\tensor \mon{t}(x)
        \ar[d]_-{\jmath_C^{-1}\tensor I_{\mon{t}(x)}}^-{\cong}
        &\text{ }
        \\
        C\tensor \mon{t}(x)
        \ar@{=}[r]
        &C\tensor \mon{t}(x)
      }
    }}
  \end{equation*}
  This completes the proof of Proposition~\ref{prop Ens(T) colaxK(K,T)Ens(T) adjunction}.
\qed\end{proof}


\begin{definition}
  Let $(\kos{T},\mon{t})$ be a strong $\kos{K}$-tensor category
  and let 
  $(\phi,\what{\phi}):(\kos{K},\id_{\kos{K}})\to (\kos{T},\mon{t})$
  be a colax $\kos{K}$-tensor functor.
  We define the \emph{reflection} of $(\phi,\what{\phi})$
  as the comonoidal $\kos{K}$-tensor natural transformation
  \begin{equation*}
    \hatar{\phi}:
    (\phi,\what{\phi})
    \Rightarrow
    (\phi(\kappa)\tensor\text{$\mon{t}$},\what{\phi(\kappa)\tensor\mon{t}})
    :(\kos{K},\id_{\kos{K}})\to (\kos{T},\mon{t})
  \end{equation*}
  introduced in Proposition~\ref{prop Ens(T) colaxK(K,T)Ens(T) adjunction}.
  We say $(\phi,\what{\phi})$  is \emph{reflective}
  if its reflection
  $\hatar{\phi}$
  is a comonoidal $\kos{K}$-tensor natural isomorphism.
\end{definition}

\begin{corollary} \label{cor Ens(T) strongKcat reflective adjunction}
  Let $(\kos{T},\mon{t})$ be a strong $\kos{K}$-tensor category
  and recall the reflective adjunction
  $\CL\dashv \iota$ in Proposition~\ref{prop Ens(T) colaxK(K,T)Ens(T) adjunction}.
  The essential image of the right adjoint $\iota$
  is the reflective full subcategory
  \begin{equation*}
    \mathbb{SMC}_{\cat{colax}}^{\kos{K}/\!\!/}\bigl((\kos{K},\id_{\kos{K}}),(\kos{T},\mon{t})\bigr)_{\cat{rfl}}
  \end{equation*}
  of reflective colax $\kos{K}$-tensor functors
  from $(\kos{K},\id_{\kos{K}})$ to $(\kos{T},\mon{t})$,
  and the adjunction $\CL\dashv \iota$ restricts to an adjoint equivalence of categories
  \begin{equation*}
    \CL:
    \mathbb{SMC}_{\cat{colax}}^{\kos{K}/\!\!/}\bigl((\kos{K},\id_{\kos{K}}),(\kos{T},\mon{t})\bigr)_{\cat{rfl}}
    \simeq
    \cat{Ens}(\kos{T})
    :\iota
    .
  \end{equation*}
\end{corollary}

Let
$(\phi,\what{\phi}):(\kos{K},\id_{\kos{K}})\to (\kos{T},\mon{t})$
be a colax $\kos{K}$-tensor functor
to a strong $\kos{K}$-tensor category $(\kos{T},\mon{t})$.
From the definition
(\ref{eq Ens(T) colaxK(K,T)Ens(T) adjunction}) of the reflection $\hatar{\phi}$,
we see that $(\phi,\what{\phi})$ is reflective
if the associated $\kos{K}$-equivariance
$\vecar{\phi}$ is a natural isomorphism.
The following lemma states that the converse is also true.

\begin{lemma} \label{lem Ens(T) equivariance and reflection}
  Let $(\kos{T},\mon{t})$ be a strong $\kos{K}$-tensor category
  and let $(\phi,\what{\phi}):(\kos{K},\id_{\kos{K}})\to (\kos{T},\mon{t})$
  be a colax $\kos{K}$-tensor functor.
  Then
  $(\phi,\what{\phi})$ is reflective
  if and only if
  the associated $\kos{K}$-equivariance $\vecar{\phi}$
  is a natural isomorphism.
\end{lemma}
\begin{proof}
  We explained the if part.
  To verify the only if part,
  it suffices to show that the following diagram
  commutes for all $x$, $y\in\obj{\CK}$.
  \begin{equation*}
    \vcenter{\hbox{
      \xymatrix@C=30pt{
        \phi(x\otimes y)
        \ar[r]^-{\hatar{\phi}_{x\otimes y}}
        \ar[d]_-{\vecar{\phi}_{y,x}}
        &\phi(\kappa)\tensor \mon{t}(x\otimes y)
        \ar[r]^-{I_{\phi(\kappa)}\tensor \mon{t}_{x,y}}_-{\cong}
        &\phi(\kappa)\tensor (\mon{t}(x)\tensor \mon{t}(y))
        \ar[d]^-{a_{\phi(\kappa),\mon{t}(x),\mon{t}(y)}}_-{\cong}
        \\
        \phi(x)\tensor \mon{t}(y)
        \ar[rr]^-{\hatar{\phi}_x\tensor I_{\mon{t}(y)}}
        &\text{ }
        &(\phi(\kappa)\tensor \mon{t}(x))\tensor \mon{t}(y)
      }
    }}
  \end{equation*}
  We can check this as follows.
  \begin{equation*}
    \vcenter{\hbox{
      \xymatrix@C=3pt{
        \phi(x\otimes y)
        \ar[d]^-{\phi(\imath_{x\otimes y})}_-{\cong}
        \ar@{.>}@/_2pc/@<-2.5ex>[ddd]|-{\hatar{\phi}_{x\otimes y}}
        \ar@{=}[r]
        &\phi(x\otimes y)
        \ar[dd]^-{\phi(\imath_x\otimes I_y)}_-{\cong}
        \ar@{=}[r]
        &\phi(x\otimes y)
        \ar[d]^-{\phi_{x,y}}
        \ar@{=}[r]
        &\phi(x\otimes y)
        \ar[dd]^-{\vecar{\phi}_{y,x}}
        \\
        \phi(\kappa\otimes (x\otimes y))
        \ar[d]^-{\phi_{\kappa,x\otimes y}}
        \ar@/^0.5pc/[dr]^-{\phi(a_{\kappa,x,y})}_-{\cong}
        &\text{ }
        &\phi(x)\tensor \phi(y)
        \ar[dd]^-{\phi(\imath_x)\tensor I_{\phi(y)}}_-{\cong}
        \ar@/^0.5pc/[dr]|-{I_{\phi(x)}\tensor \what{\phi}_y}
        &\text{ }
        \\
        \phi(\kappa)\tensor \phi(x\otimes y)
        \ar[d]|-{I_{\phi(\kappa)}\tensor \what{\phi}_{x\otimes y}}
        \ar@/^0.5pc/[dr]|-{I_{\phi(\kappa)}\tensor \phi_{x,y}}
        &\phi((\kappa\otimes x)\otimes y)
        \ar@/^0.5pc/[dr]|-{\phi_{\kappa\otimes x,y}}
        &\text{ }
        &\phi(x)\tensor \mon{t}(y)
        \ar@/^2pc/[dddl]|(0.3){\hatar{\phi}_x\tensor I_{\mon{t}(y)}}
        \\
        \phi(\kappa)\tensor \mon{t}(x\otimes y)
        \ar[d]_-{I_{\phi(\kappa)}\tensor \mon{t}_{x,y}}^-{\cong}
        &\phi(\kappa)\!\tensor\! (\phi(x)\!\tensor\! \phi(y))
        \ar@/^0.5pc/[dl]|-{I_{\phi(\kappa)}\tensor (\what{\phi}_x\tensor \what{\phi}_y)}
        \ar@/_0.5pc/[dr]_-{a_{\phi(\kappa),\phi(x),\phi(y)}}^-{\cong}
        &\phi(\kappa\otimes x)\tensor \phi(y)
        \ar[d]^-{\phi_{\kappa,x}\otimes I_{\phi(y)}}
        &\text{ }
        \\
        \phi(\kappa)\tensor (\mon{t}(x)\tensor \mon{t}(y))
        \ar[d]_-{a_{\phi(\kappa),\mon{t}(x),\mon{t}(y)}}^-{\cong}
        &\text{ }
        &(\phi(\kappa)\!\tensor\! \phi(x))\!\tensor\! \phi(y)
        \ar[d]_-{(I_{\phi(\kappa)}\tensor \what{\phi}_x)\tensor \what{\phi}_y}
        &\text{ }
        \\
        (\phi(\kappa)\tensor \mon{t}(x))\tensor \mon{t}(y)
        \ar@{=}[rr]
        &\text{ }
        &(\phi(\kappa)\tensor \mon{t}(x))\tensor \mon{t}(y)
        &\text{ }
      }
    }}
  \end{equation*}
  This completes the proof of Lemma~\ref{lem Ens(T) equivariance and reflection}.
\qed\end{proof}

For the rest of this subsection,
we focus on the case
$(\kos{T},\mon{t})=(\kos{K},\id_{\kos{K}})$.
Denote
\begin{equation} \label{eq Ens(T) Ens(K)def}
  \kos{E\!n\!s}(\kos{K})
  =(\cat{Ens}(\kos{K}),\otimes,\kappa)
\end{equation}
as the cartesian monoidal category
of cocommutative comonoids in $\kos{K}$.
An object in $\cat{Ens}(\kos{K})$
is a triple
$(c,\cp_c,e_c)$ 
which we often abbreviate as $c$.
The object $\kappa$ equipped with
$\cp_{\kappa}=\imath_{\kappa}=\jmath_{\kappa}:\kappa\xrightarrow{\cong}\kappa\otimes \kappa$
and 
$e_{\kappa}=I_{\kappa}:\kappa\xrightarrow{\cong}\kappa$
is a terminal object in $\cat{Ens}(\kos{K})$.
The product of objects $c$, $c^{\pr}$ in $\cat{Ens}(\kos{K})$
is $c\otimes c^{\pr}$ whose coproduct, counit morphisms
$\cp_{c\otimes c^{\pr}}$, $e_{c\otimes c^{\pr}}$ are
defined as in (\ref{eq Ens(T) CC' definition}).

Let $c$ be an object in $\cat{Ens}(\kos{K})$.
We have a colax symmetric monoidal comonad
\begin{equation*}
  \langle c\otimes \rangle=(c\otimes,\delta^{c\otimes},\varepsilon^{c\otimes})
\end{equation*}
on $\kos{K}$.  
The underlying colax symmetric monoidal endofunctor
$c\otimes:\kos{K}\to \kos{K}$
is the endofunctor
$c\otimes= c\otimes \slot:\CK\to \CK$
equipped with comonoidal coherence morphisms
\begin{equation*}
  \xymatrix{
    c\otimes (x\otimes y)
    \ar[r]^-{(c\otimes)_{x,y}}
    &(c\otimes x)\otimes (c\otimes y)
  }
  \qquad
  x,y\in\obj{\CK}
  \qquad
  \xymatrix{
    c\otimes \kappa
    \ar[r]^-{(c\otimes)_{\kappa}}
    &\kappa
  }
\end{equation*}
defined as in (\ref{eq Ens(T) Ctensor definition}),
and the comonoidal natural transformations
\begin{equation*}
  \xymatrix{
    c\otimes x
    \ar[r]^-{\delta^{c\otimes}_x}
    &c\otimes (c\otimes x)
  }
  \qquad
  \xymatrix{
    c\otimes x
    \ar[r]^-{\varepsilon^{c\otimes}_x}
    &x
  }
  \qquad
  x\in\obj{\CK}
\end{equation*}
are defined as in (\ref{eq Ens(T) deltaCepC definition}).

Let $(\phi,\what{\phi})$,
$(\psi,\what{\psi}):(\kos{K},\id_{\kos{K}})\to (\kos{K},\id_{\kos{K}})$
be reflective colax $\kos{K}$-tensor endofunctors
on $(\kos{K},\id_{\kos{K}})$.
By Lemma~\ref{lem Ens(T) equivariance and reflection},
the associated $\kos{K}$-equivariances $\vecar{\phi}$, $\vecar{\psi}$
are natural isomorphisms.
By Lemma~\ref{lem ColaxKtensorCat compositionKequiv},
the $\kos{K}$-equivariance
$\vecar{\psi\phi}$
associated to the composition $(\psi\phi,\what{\psi\phi})$
is also a natural isomorphism.
In particular,
the composition $(\psi\phi,\what{\psi\phi})$ is also reflective.
We denote
\begin{equation}\label{eq Ens(T) colaxKTend(K,idK) def}
  \kos{E\!n\!\!d}_{\cat{colax}}^{\kos{K}/\!\!/}(\kos{K},\id_{\kos{K}})_{\cat{rfl}}
  =
  \bigl(
    \mathbb{SMC}_{\cat{colax}}^{\kos{K}/\!\!/}\bigl((\kos{K},\id_{\kos{K}}),(\kos{K},\id_{\kos{K}})\bigr)_{\cat{rfl}}
    ,
    \circ
    ,
    (\id_{\kos{K}},I)
  \bigr)
\end{equation}
as the monoidal category of
reflective colax $\kos{K}$-tensor endofunctors
on $(\kos{K},\id_{\kos{K}})$.
The monoidal product is given by composition $\circ$
and the unit object
$(\id_{\kos{K}},I)$
is the identity strong $\kos{K}$-tensor endofunctor
of $(\kos{K},\id_{\kos{K}})$.

\begin{proposition} \label{prop Ens(T) KtensorEnd(K)=Ens(K)}
  The adjoint equivalence of categories
  \begin{equation*}
    \CL:
    \mathbb{SMC}_{\cat{colax}}^{\kos{K}/\!\!/}\bigl((\kos{K},\id_{\kos{K}}),(\kos{K},\id_{\kos{K}})\bigr)_{\cat{rfl}}
    \simeq
    \cat{Ens}(\kos{K})
    :\iota
  \end{equation*}
  described in Corollary~\ref{cor Ens(T) strongKcat reflective adjunction}
  becomes an adjoint equivalence of monoidal categories
  \begin{equation*}
    \CL:
    \kos{E\!n\!\!d}_{\cat{colax}}^{\kos{K}/\!\!/}(\kos{K},\id_{\kos{K}})_{\cat{rfl}}
    \simeq
    \kos{E\!n\!s}(\kos{K})
    :\iota
    .
  \end{equation*}
  Let $c$, $c^{\pr}$ be objects in $\cat{Ens}(\kos{K})$
  and let
  $(\phi,\what{\phi})$,
  $(\psi,\what{\psi})$ 
  be reflective colax $\kos{K}$-tensor endofunctors on $(\kos{K},\id_{\kos{K}})$.
  \begin{itemize}
    \item 
    The monoidal coherence isomorphisms of the right adjoint $\iota$
    are the following comonoidal $\kos{K}$-tensor natural isomorphisms.
    Let $x\in\obj{\CK}$.
    \begin{equation*}
      \begin{aligned}
        a_{c,c^{\pr},\slot}
        &:
        \xymatrix@C=15pt{
          \iota(c)\circ\iota(c^{\pr})
          \ar@2{->}[r]^-{\cong}
          &\iota(c\otimes c^{\pr})
          ,
        }
        &\qquad
        \imath_{\slot}
        &:
        \xymatrix@C=15pt{
          (\id_{\kos{K}},I)
          \ar@2{->}[r]^-{\cong}
          &\iota(\kappa)
          ,
        }
        \\
        a_{c,c^{\pr},x}
        &:
        \xymatrix@C=15pt{
          c\otimes (c^{\pr}\otimes x)
          \ar[r]^-{\cong}
          &(c\otimes c^{\pr})\otimes x
          ,
        }
        &\quad
        \imath_x
        &:
        \xymatrix@C=15pt{
          x
          \ar[r]^-{\cong}
          &\kappa\otimes x
          .
        }
      \end{aligned}
    \end{equation*}

    \item 
    The monoidal coherence isomorphisms of the left adjoint $\CL$
    are the following isomorphisms in $\cat{Ens}(\kos{K})$.
    \begin{equation*}
      \vcenter{\hbox{
        \xymatrix{
          \CL(\psi\phi,\what{\psi\phi})
          \ar@{=}[d]
          \ar[r]^-{\cong}
          &\CL(\psi,\what{\psi})\otimes \CL(\phi,\what{\phi})
          \ar@{=}[d]
          \\
          \psi\phi(\kappa)
          \ar[r]^-{\hatar{\psi}_{\phi(\kappa)}}_-{\cong}
          &\psi(\kappa)\otimes \phi(\kappa)
          ,
        }
      }}
      \qquad\qquad
      \CL(\id_{\kos{K}},I)=\kappa
      .
    \end{equation*}
  \end{itemize}
  In particular, we have the following relation
  for each object $x$ in $\CK$.
  \begin{equation}\label{eq Ens(T) KtensorEnd(K)=Ens(K)}
    \vcenter{\hbox{
      \xymatrix@C=60pt{
        \psi\phi(x)
        \ar[r]^-{\hatar{\psi\phi}_x}_-{\cong}
        \ar[d]_-{(\hatar{\psi}\hatar{\phi})_x}^-{\cong}
        &\psi\phi(\kappa)\otimes x
        \ar[d]^-{\hatar{\psi}_{\phi(\kappa)}\otimes I_x}_-{\cong}
        \\
        \psi(\kappa)\otimes (\phi(\kappa)\otimes x)
        \ar[r]^-{a_{\psi(\kappa),\phi(\kappa),x}}_-{\cong}
        &(\psi(\kappa)\otimes\phi(\kappa))\otimes x
      }
    }}
  \end{equation}
\end{proposition}
\begin{proof}
  We leave for the readers to check that
  $a_{c,c^{\pr},\slot}$
  and
  $\imath_{\slot}$
  are comonoidal $\kos{K}$-tensor natural transformations,
  and that $a_{c,c^{\pr},\slot}$ is natural in variables $c$, $c^{\pr}$.
  To conclude that the right adjoint $\iota$
  is a strong monoidal functor
  with monoidal coherence isomorphisms
  $a_{c,c^{\pr},\slot}$ and $\imath_{\slot}$,
  we need to check the following monoidal coherence relations.
  Let $c^{\ppr}$ be another object in $\cat{Ens}(\kos{K})$.
  \begin{equation*}
    \vcenter{\hbox{
      \xymatrix{
        \iota(c)\iota(c^{\pr})\iota(c^{\ppr})
        \ar@2{->}[d]_-{a_{c,c^{\pr},\slot}\iota(c^{\ppr})}^-{\cong}
        \ar@{=}[r]
        &\iota(c)\iota(c^{\pr})\iota(c^{\ppr})
        \ar@2{->}[d]^-{\iota(c)a_{c^{\pr},c^{\ppr},\slot}}_-{\cong}
        \\
        \iota(c)\iota(c^{\pr}\otimes c^{\ppr})
        \ar@2{->}[d]_-{a_{c,c^{\pr}\otimes c^{\ppr},\slot}}^-{\cong}
        &\iota(c\otimes c^{\pr})\iota(c^{\ppr})
        \ar@2{->}[d]^-{a_{c\otimes c^{\pr},c^{\ppr},\slot}}_-{\cong}
        \\
        \iota(c\!\otimes\! (c^{\pr}\!\otimes\! c^{\ppr}))
        \ar@2{->}[r]^-{\iota(a_{c,c^{\pr},c^{\ppr}})}_-{\cong}
        &\iota((c\!\otimes\! c^{\pr})\!\otimes\! c^{\ppr})
      }
    }}
    \vcenter{\hbox{
      \xymatrix{
        \iota(c)
        \ar@2{->}[r]^-{\imath_{\slot}\iota(c)}_-{\cong}
        \ar@/_1pc/@2{->}[dr]_-{\iota(\imath_c)}^-{\cong}
        &\iota(\kappa)\iota(c)
        \ar@2{->}[d]^-{a_{\kappa,c,\slot}}_-{\cong}
        \\
        \text{ }
        &\iota(\kappa\!\otimes\! c)
      }
    }}
    \vcenter{\hbox{
      \xymatrix{
        \iota(c)
        \ar@2{->}[r]^-{\iota(c)\imath_{\slot}}_-{\cong}
        \ar@/_1pc/@2{->}[dr]_-{\iota(\jmath_c)}^-{\cong}
        &\iota(c)\iota(\kappa)
        \ar@2{->}[d]^-{a_{c,\kappa,\slot}}_-{\cong}
        \\
        \text{ }
        &\iota(c\!\otimes\! \kappa)
      }
    }}
  \end{equation*}
  These relations are obtained by the Mac Lane coherence theorem.
  This show that the right adjoint $\iota$
  is a strong monoidal functor as we claimed.
  The left adjoint $\CL$ has a unique strong monoidal functor structure
  such that the adjoint equivalence of categories $\CL\dashv \iota$
  becomes an equivalence of monoidal categories.
  The monoidal coherence isomorphisms of $\CL$ are given as follows.
  \begin{equation*}
    \vcenter{\hbox{
      \xymatrix@C=5pt{
        \psi\phi(\kappa)
        \ar[dd]_-{(\hatar{\psi}\hatar{\phi})_{\kappa}}^-{\cong}
        \ar@{=}[rr]
        &\text{ }
        &\psi\phi(\kappa)
        \ar@/^0.5pc/[dl]_-{\psi(\hatar{\phi}_{\kappa})}^-{\cong}
        \ar@{=}[dd]
        \\
        \text{ }
        &\psi(\phi(\kappa)\otimes \kappa)
        \ar@/^0.5pc/[dl]_-{\hatar{\psi}_{\phi(\kappa)\otimes \kappa}}^-{\cong}
        \ar@/^0.5pc/[dr]_-{\psi(\jmath_{\phi(\kappa)}^{-1})}^-{\cong}
        &\text{ }
        \\
        \psi(\kappa)\otimes (\phi(\kappa)\otimes \kappa)
        \ar[d]_-{a_{\psi(\kappa),\phi(\kappa),\kappa}}^-{\cong}
        \ar@/^1pc/[ddrr]_-{I_{\psi(\kappa)}\otimes \jmath_{\phi(\kappa)}^{-1}}^-{\cong}
        &\text{ }
        &\psi\phi(\kappa)
        \ar[dd]^-{\hatar{\psi}_{\phi(\kappa)}}_-{\cong}
        \\
        (\psi(\kappa)\otimes \phi(\kappa))\otimes \kappa
        \ar[d]_-{\jmath_{\psi(\kappa)\otimes \phi(\kappa)}^{-1}}^-{\cong}
        &\text{ }
        &\text{ }
        \\
        \psi(\kappa)\otimes \phi(\kappa)
        \ar@{=}[rr]
        &\text{ }
        &\psi(\kappa)\otimes \phi(\kappa)
      }
    }}
    \vcenter{\hbox{
      \xymatrix{
        \CL(\id_{\kos{K}},I)
        \ar[d]_-{\CL(\imath_{\slot})}^-{\cong}
        \ar@{=}[r]
        &\kappa
        \ar[d]^-{\imath_{\kappa}}_-{\cong}
        \ar@{=}@/^2pc/@<2ex>[dd]
        \\
        \CL\iota(\kappa)
        \ar[d]_-{\jmath_{\kappa}^{-1}}^-{\cong}
        \ar@{=}[r]
        &\kappa\otimes\kappa
        \ar[d]^-{\jmath_{\kappa}^{-1}}_-{\cong}
        \\
        \kappa
        \ar@{=}[r]
        &\kappa
      }
    }}
  \end{equation*}
  This shows that we have the adjoint equivalence $\CL\dashv\iota$
  of monoidal categories as we claimed.
  In particular, we obtain the relation (\ref{eq Ens(T) KtensorEnd(K)=Ens(K)})
  which is a consequence of the fact that 
  $\CL\dashv\iota$
  is an adjoint equivalence of monoidal categories.
  \begin{equation*}
    \vcenter{\hbox{
      \xymatrix@C=30pt{
        (\psi,\what{\psi})\circ (\phi,\what{\phi})
        \ar@2{->}[d]_-{\hatar{\psi}\hatar{\phi}}^-{\cong}
        \ar@{=}[r]
        &(\psi\phi,\what{\psi\phi})
        \ar@2{->}[r]^-{\hatar{\psi\phi}}_-{\cong}
        &\iota\CL(\psi\phi,\what{\psi\phi})
        \ar@2{->}[d]^-{\iota(\hatar{\psi}_{\phi(\kappa)})}_-{\cong}
        \\
        \iota\CL(\psi,\what{\psi})\circ \iota\CL(\phi,\what{\phi})
        \ar@2{->}[rr]^-{a_{\CL(\psi,\what{\psi}),\CL(\phi,\what{\phi}),\slot}}_-{\cong}
        &\text{ }
        &\iota(\CL(\psi,\what{\psi})\otimes \CL(\phi,\what{\phi}))
      }
    }}
  \end{equation*}
  This completes the proof of Proposition~\ref{prop Ens(T) KtensorEnd(K)=Ens(K)}.
\qed\end{proof}

Since $\kos{E\!n\!s}(\kos{K})$ is a cartesian monoidal category,
each object $c$ in $\cat{Ens}(\kos{K})$
has a unique structure of a comonoid in $\kos{E\!n\!s}(\kos{K})$
with coproduct $\cp_c:c\to c\otimes c$
and counit $e_c:c\to \kappa$.
Recall the adjoint equivalence $\CL\dashv \iota$ of monoidal categories
in Proposition~\ref{prop Ens(T) KtensorEnd(K)=Ens(K)}.
The right adjoint $\iota$
sends the comonoid $(c,\cp_c,e_c)$ in $\kos{E\!n\!s}(\kos{K})$
to the colax $\kos{K}$-tensor comonad
$\langle c\otimes, \what{c\otimes}\rangle
=
(c\otimes ,\what{c\otimes},\delta^{c\otimes},\varepsilon^{c\otimes})$
on $(\kos{K},\id_{\kos{K}})$
where
\begin{equation*}
  \vcenter{\hbox{
    \xymatrix@C=40pt{
      \text{ }
      &\kos{K}
      \ar[d]^-{c\otimes}
      \\
      \kos{K}
      \ar@/^0.7pc/[ur]^-{\id_{\kos{K}}}
      \ar[r]_-{\id_{\kos{K}}}
      \xtwocell[r]{}<>{<-2.5>{\quad\what{c\otimes}}}
      &\kos{K}
    }
  }}
  \qquad
  \what{c\otimes}_x=\varepsilon^{c\otimes}_x:
  \xymatrix@C=18pt{
    c\otimes x
    \ar[r]^-{e_c\otimes I_x}
    &\kappa\otimes x
    \ar[r]^-{\imath_x^{-1}}_-{\cong}
    &x
    ,
  }
  \quad
  x\in\obj{\CK}.
\end{equation*}

\begin{proposition} \label{prop Ens(T) KtensorEnd(K) is cartesian}
  The monoidal category
  $\kos{E\!n\!\!d}_{\cat{colax}}^{\kos{K}/\!\!/}(\kos{K},\id_{\kos{K}})_{\cat{rfl}}$
  is cartesian.
  \begin{enumerate}
    \item 
    Every reflective colax $\kos{K}$-tensor endofunctor
    $(\phi,\what{\phi})$ on $(\kos{K},\id_{\kos{K}})$
    has a unique colax $\kos{K}$-tensor comonad structure
    \begin{equation*}
      \langle\phi,\what{\phi}\rangle
      =
      (\phi,\what{\phi},\delta^{\phi},\varepsilon^{\phi})
      \qquad\quad
      \begin{aligned}
        \delta^{\phi}
        &:
        (\phi,\what{\phi})\Rightarrow(\phi\phi,\what{\phi\phi})
        \\
        \varepsilon^{\phi}
        &:
        (\phi,\what{\phi})\Rightarrow (\id_{\kos{K}},I)
      \end{aligned}
    \end{equation*}
    where the components of $\delta^{\phi}$, $\varepsilon^{\phi}$ at
    $x\in\obj{\CK}$ are given below.
    \begin{equation}\label{eq Ens(T) KtensorEnd(K) is cartesian}
      \begin{aligned}
       \delta^{\phi}_x
        &:\!
        \xymatrix@C=25pt{
          \phi(x)
          \ar[r]^-{\phi(\imath_x)}_-{\cong}
          &\phi(\kappa\otimes x)
          \ar[r]^-{\phi_{\kappa,x}}
          &\phi(\kappa)\otimes \phi(x)
          \ar[r]^-{(\hatar{\phi}_{\phi(x)})^{-1}}_-{\cong}
          &\phi\phi(x)
        }
        \\
        \varepsilon^{\phi}_x
        &:\!
        \xymatrix@C=20pt{
          \phi(x)
          \ar[r]^-{\what{\phi}_x}
          &x
        }
      \end{aligned}
    \end{equation}
    \item 
    Every comonoidal $\kos{K}$-tensor natural transformation
    $\vartheta:(\phi,\what{\phi})\Rightarrow (\psi,\what{\psi})$
    between reflective colax $\kos{K}$-tensor endofunctors
    $(\phi,\what{\phi})$, $(\psi,\what{\psi})$
    on $(\kos{K},\id_{\kos{K}})$
    becomes a morphism 
    $\vartheta:\langle\phi,\what{\phi}\rangle\Rightarrow \langle\psi,\what{\psi}\rangle$
    of colax $\kos{K}$-tensor comonads.
  \end{enumerate}
\end{proposition}
\begin{proof}
  As we have the adjoint equivalence $\CL\dashv\iota$
  of monoidal categories described in Proposition~\ref{prop Ens(T) KtensorEnd(K)=Ens(K)},
  we conclude that 
  $\kos{E\!n\!\!d}_{\cat{colax}}^{\kos{K}/\!\!/}(\kos{K},\id_{\kos{K}})_{\cat{rfl}}$
  is cartesian monoidal.
  Therefore statements 1 and 2 are true,
  but we still need to verify
  the explicit descriptions 
  (\ref{eq Ens(T) KtensorEnd(K) is cartesian})
  of $\delta^{\phi}$ and $\varepsilon^{\phi}$.
  For each reflective colax $\kos{K}$-tensor endofunctor $(\phi,\what{\phi})$
  on $(\kos{K},\id_{\kos{K}})$,
  the reflection $\hatar{\phi}$
  becomes an isomorphism
  \begin{equation*}
    \xymatrix@C=15pt{
      \hatar{\phi}:
      \langle \phi,\what{\phi}\rangle
      \ar@2{->}[r]^-{\cong}
      &\langle \phi(\kappa)\otimes,\what{\phi(\kappa)\otimes}\rangle
    }
  \end{equation*}
  of reflective colax $\kos{K}$-tensor comonads on $(\kos{K},\id_{\kos{K}})$.
  Thus we have the following relations
  for each object $x$ in $\CK$.
  \begin{equation*}
    \vcenter{\hbox{
      \xymatrix@C=20pt{
        \phi(x)
        \ar[d]_-{\delta^{\phi}_x}
        \ar[rr]^-{\hatar{\phi}_x}_-{\cong}
        &\text{ }
        &\phi(\kappa)\otimes x
        \ar[d]^-{\delta^{\phi(\kappa)\otimes}_x}
        &\phi(x)
        \ar[rr]^-{\hatar{\phi}_x}_-{\cong}
        \ar[dr]_-{\varepsilon^{\phi}_x}
        &\text{ }
        &\phi(\kappa)\otimes x
        \ar[dl]^-{\varepsilon^{\phi(\kappa)\otimes}_x}
        \\
        \phi\phi(x)
        \ar[rr]^-{(\hatar{\phi}\hatar{\phi})_x}_-{\cong}
        &\text{ }
        &\phi(\kappa)\otimes (\phi(\kappa)\otimes x)
        &\text{ }
        &x
        &\text{ }
      }
    }}
  \end{equation*}
  We can check that $\varepsilon^{\phi}_x=\what{\phi}_x$ as follows.
  \begin{equation*}
    \qquad
    \varepsilon^{\phi}_x=
    \vcenter{\hbox{
      \xymatrix@C=30pt{
        \phi(x)
        \ar[ddd]_-{\hatar{\phi}_x}^-{\cong}
        \ar@{=}[r]
        &\phi(x)
        \ar[d]^-{\phi(\imath_x)}_-{\cong}
        \ar@{=}[r]
        &\phi(x)
        \ar@{=}[ddd]
        \\
        \text{ }
        &\phi(\kappa\otimes x)
        \ar[d]^-{\phi_{\kappa,x}}
        &\text{ }
        \\
        \text{ }
        &\phi(\kappa)\otimes \phi(x)
        \ar@/^1pc/[dl]_-{I_{\phi(\kappa)}\otimes \what{\phi}_x}
        \ar@/_1pc/[dr]^-{\varepsilon^{\phi(\kappa)\otimes}_{\phi(x)}}
        &\text{ }
        \\
        \phi(\kappa)\otimes x
        \ar[d]_-{\varepsilon^{\phi(\kappa)\otimes}_x}
        &\text{ }
        &\phi(x)
        \ar[d]^-{\what{\phi}_x}
        \\
        x
        \ar@{=}[rr]
        &\text{ }
        &x
      }
    }}
  \end{equation*}
  We are left to show that $\delta^{\phi}_x$ is as claimed in (\ref{eq Ens(T) KtensorEnd(K) is cartesian}).
  We have the relation
  \begin{equation*}
    (\dagger):
    \vcenter{\hbox{
      \xymatrix@C=20pt{
        \phi(\kappa\otimes x)
        \ar[d]_-{\phi_{\kappa,x}}
        \ar@{=}[rr]
        &\text{ }
        &\phi(\kappa\otimes x)
        \ar[d]^-{\phi_{\kappa,x}}
        \\
        \phi(\kappa)\otimes \phi(x)
        \ar[d]_-{I_{\phi(\kappa)}\otimes \what{\phi}_x}
        \ar@/_1pc/[ddr]|-{\delta^{\phi(\kappa)\otimes}_{\phi(x)}}
        \ar@{}[rr]|(0.4){(\dagger1)}
        &\text{ }
        &\phi(\kappa)\otimes \phi(x)
        \ar@/_1pc/[dl]_-{I_{\phi(\kappa)}\otimes \phi(\imath_x)}^-{\cong}
        \ar[ddd]^-{I_{\phi(\kappa)}\otimes \hatar{\phi}_x}_-{\cong}
        \\
        \phi(\kappa)\otimes x
        \ar[dd]_-{\delta^{\phi(\kappa)\otimes}_x}
        &\phi(\kappa)\otimes \phi(\kappa\otimes x)
        \ar[d]^-{I_{\phi(\kappa)}\otimes \phi_{\kappa,x}}
        &\text{ }
        \\
        \text{ }
        &\phi(\kappa)\otimes (\phi(\kappa)\otimes \phi(x))
        \ar@/_0.5pc/[dr]|-{I_{\phi(\kappa)}\otimes (I_{\phi(\kappa)}\otimes \what{\phi}_x)}
        &\text{ }
        \\
        \phi(\kappa)\otimes (\phi(\kappa)\otimes x)
        \ar@{=}[rr]
        &\text{ }
        &\phi(\kappa)\otimes (\phi(\kappa)\otimes x)
      }
    }}
  \end{equation*}
  where the diagram $(\dagger1)$ is verified below.
  \begin{equation*}
    (\dagger1):
    \!\!\!\!\!\!\!\!
    \!\!\!\!\!\!\!\!
    \!\!\!\!\!\!\!\!
    \vcenter{\hbox{
      \xymatrix@C=-5pt{
        \phi(\kappa\otimes x)
        \ar@<-1ex>[d]_-{\phi_{\kappa,x}}
        \ar@{=}[r]
        &\phi(\kappa\otimes x)
        \ar[d]^-{\phi(\cp_{\kappa}\otimes I_x)}_-{\cong}
        \ar@{=}[rr]
        &\text{ }
        &\phi(\kappa\otimes x)
        \ar@/_1pc/[ddl]_-{\phi(I_\kappa\otimes \imath_x)}^-{\cong}
        \ar[d]^-{\phi_{\kappa,x}}
        \\
        \phi(\kappa)\otimes \phi(x)
        \ar@/_0.5pc/[dr]_-{\phi(\cp_{\kappa})\otimes I_{\phi(x)}}^-{\cong}
        \ar@{.>}@<-0ex>@/_3pc/[ddrr]|-{\cp_{\phi(\kappa)}\otimes I_{\phi(x)}}
        \ar@<-1ex>[ddd]_-{\delta^{\phi(\kappa)\otimes}_{\phi(x)}}
        &\phi((\kappa\otimes \kappa)\otimes x)
        \ar[d]^-{\phi_{\kappa\otimes\kappa,x}}
        \ar@/^0.5pc/[dr]^-{\phi(a_{\kappa,\kappa,x}^{-1})}_-{\cong}
        &\text{ }
        &\phi(\kappa)\otimes \phi(x)
        \ar[dd]^-{I_{\phi(\kappa)}\otimes \phi(\imath_x)}_-{\cong}
        \\
        \text{ }
        &\phi(\kappa\otimes \kappa)\!\otimes\! \phi(x)
        \ar@/_0.5pc/[dr]|-{\phi_{\kappa,\kappa}\otimes I_{\phi(x)}}
        &\phi(\kappa\otimes (\kappa\otimes x))
        \ar@/^0.5pc/[dr]^-{\phi_{\kappa,\kappa\otimes x}}
        &\text{ }
        \\
        \text{ }
        &\text{ }
        &(\phi(\kappa)\!\otimes\! \phi(\kappa))\!\otimes\! \phi(x)
        \ar@/_0.5pc/[dr]_(0.35){a^{-1}_{\phi(\kappa),\phi(\kappa),\phi(x)}}^-{\cong}
        &\phi(\kappa)\otimes \phi(\kappa\otimes x)
        \ar[d]^-{I_{\phi(\kappa)}\otimes \phi_{\kappa,x}}
        \\
        \phi(\kappa)\!\otimes\! (\phi(\kappa)\!\otimes\! \phi(x))
        \ar@{=}[rrr]
        &\text{ }
        &\text{ }
        &\phi(\kappa)\!\otimes\! (\phi(\kappa)\!\otimes\! \phi(x))
      }
    }}
  \end{equation*}
  Using the relation $(\dagger)$,
  we obtain the description of $\delta^{\phi}_x$
  as we claimed.
  \begin{equation*}
    \vcenter{\hbox{
      \xymatrix@C=20pt{
        \phi(x)
        \ar[ddddd]_-{\delta^{\phi}_x}
        \ar@{=}[r]
        &\phi(x)
        \ar[ddd]_-{\hatar{\phi}_x}^-{\cong}
        \ar@{=}[rr]
        &\text{ }
        &\phi(x)
        \ar[d]^-{\phi(\imath_x)}_-{\cong}
        \\
        \text{ }
        &\text{ }
        &\text{ }
        &\phi(\kappa\otimes x)
        \ar[d]^-{\phi_{\kappa,x}}
        \ar@/_0.5pc/[dl]_-{\phi_{\kappa,x}}
        \\
        \text{ }
        &\text{ }
        &\phi(\kappa)\otimes \phi(x)
        \ar@/_0.5pc/[dl]_-{I_{\phi(\kappa)}\otimes \what{\phi}_x}
        \ar@{}[d]|-{(\dagger)}
        &\phi(\kappa)\otimes \phi(x)
        \ar@/^1pc/[ddll]^(0.6){I_{\phi(\kappa)}\otimes \hatar{\phi}_x}_(0.6){\cong}
        \ar[ddd]^-{\hatar{\phi}_{\phi(x)}^{-1}}_-{\cong}
        \\
        \text{ }
        &\phi(\kappa)\otimes x
        \ar[d]_-{\delta^{\phi(\kappa)\otimes}_x}
        &\text{ }
        &\text{ }
        \\
        \text{ }
        &\phi(\kappa)\otimes (\phi(\kappa)\otimes x)
        \ar[d]_-{(\hatar{\phi}\hatar{\phi})_x^{-1}}^-{\cong}
        &\text{ }
        &\text{ }
        \\
        \phi\phi(x)
        \ar@{=}[r]
        &\phi\phi(x)
        \ar@{=}[rr]
        &\text{ }
        &\phi\phi(x)
      }
    }}
  \end{equation*}
  This completes the proof of Proposition~\ref{prop Ens(T) KtensorEnd(K) is cartesian}.
\qed\end{proof}

\subsection{Right-strong $\kos{K}$-tensor adjunctions}
\label{subsec RSKtensoradj}

In this subsection we study about adjunctions
between colax $\kos{K}$-tensor categories.
We begin by introducing adjunctions
between symmetric monoidal categories
in colax context.

\begin{definition} \label{def RSKtensoradj RSSMadj}
  A \emph{right-strong symmetric monoidal (RSSM) adjunction}
  is an adjunction
  internal to the $2$-category
  $\mathbb{SMC}_{\cat{colax}}$.
\end{definition}

Let us explain Definition~\ref{def RSKtensoradj RSSMadj}
in detail.
Let $\kos{T}=(\CT,\tensor,\unit)$,
$\text{$\kos{S}$}=(\CS,\ctimes,\pzc{1})$
be symmetric monoidal categories.
A RSSM adjunction
\begin{equation*}
  \kos{f}:\kos{S}\to \kos{T}
\end{equation*}
is a data of a pair of colax symmetric monoidal functors
$\kos{f}_!:\kos{S}\to \kos{T}$,
$\kos{f}^*:\kos{T}\to \kos{S}$
and an adjunction
$\kos{f}_!:\CS\rightleftarrows \CT:\kos{f}^*$
such that the adjunction unit
$\eta:\id_{\kos{S}}\Rightarrow \kos{f}^*\kos{f}_!$
and the adjunction counit
$\epsilon:\kos{f}_!\kos{f}^*\Rightarrow \id_{\kos{T}}$
are comonoidal natural transformations.
We often denote a RSSM adjunction
$\kos{f}:\kos{S}\to\kos{T}$
as follows.
\begin{equation*}
  \vcenter{\hbox{
    \xymatrix{
      \kos{T}
      \ar@/^1pc/[d]^-{\kos{f}^*}
      \\
      \kos{S}
      \ar@/^1pc/[u]^-{\kos{f}_!}
    }
  }}
\end{equation*}
In this case,
the right adjoint $\kos{f}^*$
is a strong symmetric monoidal functor.
\begin{equation}\label{eq RSKtensoradj f*coherence}
  \text{$\kos{f}$}^*_{X,Y}:
  \text{$\kos{f}$}^*(X\tensor Y)
  \xrightarrow{\cong}
  \text{$\kos{f}$}^*(X)\ctimes \text{$\kos{f}$}^*(Y)
  \qquad
  X,Y\in\obj{\CT}
  \qquad
  \text{$\kos{f}$}^*_{\unit}:
  \text{$\kos{f}$}^*(\unit)
  \xrightarrow{\cong}
  \pzc{1}
\end{equation}
The inverses
$(\kos{f}^*_{X,Y})^{-1}:
\text{$\kos{f}$}^*(X)\ctimes \text{$\kos{f}$}^*(Y)
\xrightarrow{\cong}
\text{$\kos{f}$}^*(X\tensor Y)$
and
$(\kos{f}^*_{\unit})^{-1}:
\text{$\pzc{1}$}\xrightarrow{\cong}\kos{f}^*(\unit)$
of the coherence morphisms
of $\kos{f}^*$ are morphisms in $\CS$
that correspond to the following morphisms
in $\CT$ under the adjunction $\kos{f}_!\dashv\kos{f}^*$.
\begin{equation*}
  \xymatrix@C=35pt{
    \kos{f}_!(\kos{f}^*(X)\ctimes \kos{f}^*(Y))
    \ar[r]^-{\kos{f}_{!\kos{f}^*(X),\kos{f}^*(Y)}}
    &\kos{f}_!\kos{f}^*(X)\tensor \kos{f}_!\kos{f}^*(Y)
    \ar[r]^-{\epsilon_X\tensor \epsilon_Y}
    &X\tensor Y
  }
  \qquad
  \xymatrix{
    \kos{f}_!(\pzc{1})
    \ar[r]^-{\kos{f}_{!\pzc{1}}}
    &\unit
  }
\end{equation*}

Conversely,
suppose we are given an adjunction
$\kos{f}_!:\CS\rightleftarrows \CT:\kos{f}^*$
between the underlying categories $\CT$, $\CS$
and the right adjoint
has a strong symmetric monoidal functor structure
$\kos{f}^*:\kos{T}\to \kos{S}$.
Then the left adjoint has a unique colax symmetric monoidal functor structure
$\kos{f}_!:\kos{S}\to \kos{T}$
such that the given adjunction
$\kos{f}_!\dashv\kos{f}^*$
becomes a RSSM adjunction
$\kos{f}:\kos{S}\to \kos{T}$.
The comonoidal coherence morphisms
\begin{equation}\label{eq RSKtensoradj f!coherence}
  \text{$\kos{f}$}_{!\pzc{X},\pzc{Y}}
  :
  \text{$\kos{f}$}_!(\pzc{X}\ctimes \pzc{Y})
  \to
  \text{$\kos{f}$}_!(\pzc{X})\tensor \text{$\kos{f}$}_!(\pzc{Y})
  \qquad
  \pzc{X},\pzc{Y}\in \obj{\CS}
  \qquad
  \text{$\kos{f}$}_{!\pzc{1}}:
  \text{$\kos{f}$}_!(\pzc{1})\to \unit
\end{equation}
of $\kos{f}_!$
are morphisms in $\CT$
that correspond to the following morphisms in $\CS$
under the adjunction $\kos{f}_!\dashv\kos{f}^*$.
\begin{equation*}
  \xymatrix@C=35pt{
    \text{$\pzc{X}$}\ctimes \text{$\pzc{Y}$}
    \ar[r]^-{\eta_{\text{$\pzc{X}$}}\ctimes \eta_{\text{$\pzc{Y}$}}}
    &\text{$\kos{f}$}^*\text{$\kos{f}$}_!(\text{$\pzc{X}$})\ctimes \text{$\kos{f}$}^*\text{$\kos{f}$}_!(\text{$\pzc{Y}$})
    \ar[r]^-{(\text{$\kos{f}$}^*_{\text{$\kos{f}$}_!(\text{$\pzc{X}$}),\text{$\kos{f}$}_!(\text{$\pzc{Y}$})})^{-1}}_-{\cong}
    &\text{$\kos{f}$}^*(\text{$\kos{f}$}_!(\text{$\pzc{X}$})\tensor \text{$\kos{f}$}_!(\text{$\pzc{Y}$}))
  }
  \qquad
  \xymatrix{
    \text{$\pzc{1}$}
    \ar[r]^-{(\kos{f}^*_{\unit})^{-1}}_-{\cong}
    &\kos{f}^*(\unit)
  }
\end{equation*}

\begin{definition}
  Let $\kos{T}$, $\kos{S}$ be symmetric monoidal categories.
  A \emph{morphism} 
  $\vartheta:\kos{f}\Rightarrow\kos{g}$
  of RSSM adjunctions
  $\kos{f}$, $\kos{g}:\kos{S}\to\kos{T}$
  is a comonoidal natural transformation
  $\vartheta:\kos{f}^*\Rightarrow \kos{g}^*:\kos{T}\to \kos{S}$
  between right adjoints.
\end{definition}

We denote the $2$-category of
symmetric monoidal categories,
RSSM adjunctions
and morphisms of RSSM adjunctions
as
\begin{equation} \label{eq RSKtensoradj RSSM 2-cat}
  \mathbb{ADJ}_{\cat{right}}(\mathbb{SMC}_{\cat{colax}}).
\end{equation}
We use the subscript $_{\cat{right}}$
to indicate that $2$-cells in
$\mathbb{ADJ}_{\cat{right}}(\mathbb{SMC}_{\cat{colax}})$
are defined as $2$-cells in
$\mathbb{SMC}_{\cat{colax}}$
between right adjoints.
The $2$-category structure of
$\mathbb{ADJ}_{\cat{right}}(\mathbb{SMC}_{\cat{colax}})$
is defined as one expects.
For instance,
the vertical composition of $2$-cells
$\kos{f}\Rightarrow\kos{g}\Rightarrow\kos{h}$
in
$\mathbb{ADJ}_{\cat{right}}(\mathbb{SMC}_{\cat{colax}})$
is given by
the vertical composition of $2$-cells
$\kos{f}^*\Rightarrow\kos{g}^*\Rightarrow\kos{h}^*$
between right adjoints.

We introduce an important property of RSSM adjunctions.

\begin{definition} \label{def RSKtensoradj projformula}
  Let $\kos{f}:\kos{S}\to \kos{T}$
  be a RSSM adjunction between 
  symmetric monoidal categories $\kos{T}$, $\kos{S}$.
  We define a canonical natural transformation
  \begin{equation*}
    \varphi_{\slot,\slot}
    :\kos{f}_!(\slot\ctimes \kos{f}^*(\slot))\Rightarrow \kos{f}_!(\slot)\tensor \slot
    :\CS\times \CT\to \CT
  \end{equation*}
  whose component at 
  $\text{$\pzc{Z}$}\in\obj{\CS}$,
  $X\in\obj{\CT}$
  is given by
  \begin{equation*}
    \varphi_{\text{$\pzc{Z}$},X}
    :
    \xymatrix@C=30pt{
      \kos{f}_!(\text{$\pzc{Z}$}\ctimes \kos{f}^*(X))
      \ar[r]^-{\kos{f}_{!\pzc{Z},\kos{f}^*(X)}}
      &\kos{f}_!(\pzc{Z})\tensor \kos{f}_!\kos{f}^*(X)
      \ar[r]^-{I_{\kos{f}_!(\pzc{Z})}\tensor \epsilon_X}
      &\kos{f}_!(\pzc{Z})\tensor X
      .
    }
  \end{equation*}
  We say a RSSM adjunction $\kos{f}:\kos{S}\to \kos{T}$
  \emph{satisfies the projection formula}
  if the natural transformation
  $\varphi_{\slot,\slot}$
  associated to $\kos{f}$ is a natural isomorphism.
\end{definition}

Next we study about adjunctions
between colax $\kos{K}$-tensor categories.

\begin{definition} \label{def RSKtensoradj RSKtensoradjunction}
  A \emph{right-strong $\kos{K}$-tensor adjunction}
  is an adjunction
  internal to the $2$-category
  $\mathbb{SMC}^{\kos{K}/\!\!/}_{\cat{colax}}$.
\end{definition}

Let us explain Definition~\ref{def RSKtensoradj RSKtensoradjunction} in detail.
Let $(\kos{T},\mon{t})$, $(\kos{S},\mon{s})$
be colax $\kos{K}$-tensor categories.
A right-strong $\kos{K}$-tensor adjunction
\begin{equation*}
  (\kos{f},\what{\kos{f}}):(\kos{S},\mon{s})\to (\kos{T},\mon{t})
\end{equation*}
is a data of a pair of colax $\kos{K}$-tensor functors
$(\kos{f}_!,\what{\kos{f}}_!):(\kos{S},\mon{s})\to (\kos{T},\mon{t})$,
$(\kos{f}^*,\what{\kos{f}}^*):(\kos{T},\mon{t})\to (\kos{S},\mon{s})$
and an adjunction
$\kos{f}_!:\CS\rightleftarrows \CT:\kos{f}^*$
such that the adjunction unit
$\eta:(\id_{\kos{S}},I)\Rightarrow (\kos{f}^*\kos{f}_!,\what{\kos{f}^*\kos{f}_!})$
and the adjunction counit 
$\epsilon:(\kos{f}_!\kos{f}^*,\what{\kos{f}_!\kos{f}^*})\Rightarrow (\id_{\kos{T}},I)$
are comonoidal $\kos{K}$-tensor natural transformations.
We often denote a right-strong $\kos{K}$-tensor adjunction
$(\kos{f},\what{\kos{f}}):(\kos{S},\mon{s})\to (\kos{T},\mon{t})$
as follows.
\begin{equation*}
  \vcenter{\hbox{
    \xymatrix{
      (\kos{T},\mon{t})
      \ar@/^1pc/[d]^-{(\kos{f}^*,\what{\kos{f}}^*)}
      \\
      (\kos{S},\mon{s})
      \ar@/^1pc/[u]^-{(\kos{f}_!,\what{\kos{f}}_!)}
    }
  }}
\end{equation*}
In this case, the right adjoint 
$(\kos{f}^*,\what{\kos{f}}^*):(\kos{T},\mon{t})\to (\kos{S},\mon{s})$
is a strong $\kos{K}$-tensor functor.
We explained the coherence isomorphisms
of the underlying strong symmetric monoidal functor
$\kos{f}^*:\kos{T}\to \kos{S}$
in (\ref{eq RSKtensoradj f*coherence}).
The inverse of the comonoidal natural isomorphism
$\xymatrix@C=15pt{\what{\kos{f}}^*:\kos{f}^*\mon{t}\ar@2{->}[r]^-{\cong}&\mon{s}}$
is explicitly given below.
\begin{equation*}
  \xymatrix@C=30pt{
    (\what{\kos{f}}^*)^{-1}:
    \mon{s}
    \ar@2{->}[r]^-{\eta\mon{s}}
    &\kos{f}^*\kos{f}_!\mon{s}
    \ar@2{->}[r]^-{\kos{f}^*\what{\kos{f}}_!}
    &\kos{f}^*\mon{t}
  }
\end{equation*}

Conversely,
suppose we are given an adjunction
$\kos{f}_!:\CS\rightleftarrows \CT:\kos{f}^*$
between the underlying categories $\CT$, $\CS$
and the right adjoint has a strong $\kos{K}$-tensor functor structure
$(\kos{f}^*,\what{\kos{f}}^*):(\kos{T},\mon{t})\to (\kos{S},\mon{s})$.
Then the left adjoint has a unique colax $\kos{K}$-tensor functor structure
$(\kos{f}_!,\what{\kos{f}}_!):
(\kos{S},\mon{s})\to (\kos{T},\mon{t})$
such that the given adjunction
$\kos{f}_!\dashv\kos{f}^*$
becomes a right-strong $\kos{K}$-tensor adjunction
$(\kos{f},\what{\kos{f}}):(\kos{S},\mon{s})\to (\kos{T},\mon{t})$.
We explained the coherence morphisms of $\kos{f}_!$
in (\ref{eq RSKtensoradj f!coherence}).
The comonoidal natural transformation
$\what{\kos{f}}_!:\kos{f}_!\mon{s}\Rightarrow\mon{t}$
is described below.
\begin{equation} \label{eq RSKtensoradj leftadj uniqueKTstr}
  \xymatrix@C=30pt{
    \what{\kos{f}}_!:
    \kos{f}_!\mon{s}
    \ar@2{->}[r]^-{\kos{f}_!(\what{\kos{f}}^*)^{-1}}_-{\cong}
    &\kos{f}_!\kos{f}^*\mon{t}
    \ar@2{->}[r]^-{\epsilon\mon{t}}
    &\mon{t}
  }
\end{equation}

\begin{definition}
  Let $(\kos{T},\mon{t})$,
  $(\kos{S},\mon{s})$
  be colax $\kos{K}$-tensor categories.
  A \emph{morphism} 
  $\vartheta:(\kos{f},\what{\kos{f}})\Rightarrow(\kos{g},\what{\kos{g}})$
  of right-strong $\kos{K}$-tensor adjunctions
  $(\kos{f},\what{\kos{f}})$,
  $(\kos{g},\what{\kos{g}}):(\kos{S},\mon{s})\to(\kos{T},\mon{t})$
  is a comonoidal $\kos{K}$-tensor natural transformation
  $\vartheta:(\kos{f}^*,\what{\kos{f}}^*)\Rightarrow (\kos{g}^*,\what{\kos{g}}^*)
  :(\kos{T},\mon{t})\to (\kos{S},\mon{s})$
  between right adjoints.
\end{definition}

We denote the $2$-category of
colax $\kos{K}$-tensor categories,
right-strong $\kos{K}$-tensor adjunctions
and morphisms of right-strong $\kos{K}$-tensor adjunctions
as
\begin{equation} \label{eq RSKtensoradj RSKTadj 2-cat}
  \mathbb{ADJ}_{\cat{right}}\bigl(\mathbb{SMC}^{\kos{K}/\!\!/}_{\cat{colax}}\bigr).
\end{equation}
Again, we use the subscript $_{\cat{right}}$
to indicate that $2$-cells in
$\mathbb{ADJ}_{\cat{right}}\bigl(\mathbb{SMC}^{\kos{K}/\!\!/}_{\cat{colax}}\bigr)$
are defined as $2$-cells in
$\mathbb{SMC}^{\kos{K}/\!\!/}_{\cat{colax}}$
between right adjoints.
The $2$-category structure of
$\mathbb{ADJ}_{\cat{right}}\bigl(\mathbb{SMC}^{\kos{K}/\!\!/}_{\cat{colax}}\bigr)$
is defined as one expects.
For instance,
the vertical composition of $2$-cells in
$\mathbb{ADJ}_{\cat{right}}\bigl(\mathbb{SMC}^{\kos{K}/\!\!/}_{\cat{colax}}\bigr)$
is given by
the vertical composition of $2$-cells
between right adjoints.

\begin{definition}
  Let $(\kos{T},\mon{t})$, $(\kos{S},\mon{s})$
  be colax $\kos{K}$-tensor categories.
  We say a right-strong $\kos{K}$-tensor adjunction
  $(\kos{f},\what{\kos{f}}):(\kos{S},\mon{s})\to (\kos{T},\mon{t})$
  \emph{satisfies the projection formula}
  if the underlying RSSM adjunction
  $\kos{f}:\kos{S}\to \kos{T}$
  satisfies the projection formula.
\end{definition}

\begin{remark} \label{rem RSKtensoradj projformula as Tequiv}
  Let $\kos{T}$, $\kos{S}$ be symmetric monoidal categories
  and let $\kos{f}:\kos{S}\to \kos{T}$ be a RSSM adjunction.
  We explain how we can enhance the given RSSM adjunction
  $\kos{f}$
  to a right-strong $\kos{T}$-tensor adjunction.
  \begin{equation*}
    \vcenter{\hbox{
      \xymatrix{
        \kos{T}
        \ar@/^1pc/[d]^-{\kos{f}^*}
        \\
        \kos{S}
        \ar@/^1pc/[u]^-{\kos{f}_!}
      }
    }}
  \end{equation*}
  We can see $\kos{T}$ as a strong $\kos{T}$-tensor category
  which we denote as $(\kos{T},\id_{\kos{T}})$.
  We have a strong $\kos{T}$-tensor category
  $(\kos{S},\kos{f}^*)$
  as the right adjoint $\kos{f}^*:\kos{T}\to \kos{S}$
  is a strong symmetric monoidal functor.
  Then we obtain a right-strong $\kos{T}$-tensor adjunction
  \begin{equation*}
    \vcenter{\hbox{
      \xymatrix{
        (\kos{T},\id_{\kos{T}})
        \ar@/^1pc/[d]^-{(\kos{f}^*,I)}
        \\
        (\kos{S},\kos{f}^*)
        \ar@/^1pc/[u]^-{(\kos{f}_!,\epsilon)}
      }
    }}
  \end{equation*}
  as the right adjoint
  has a strong $\kos{T}$-tensor functor structure
  $(\kos{f}^*,I):(\kos{T},\id_{\kos{T}})\to (\kos{S},\kos{f}^*)$.
  The left adjoint colax $\kos{T}$-tensor functor
  is $\kos{f}_!$ equipped with
  $\epsilon:\kos{f}_!\kos{f}^*\Rightarrow \id_{\kos{T}}$.
  \begin{equation*}
    \vcenter{\hbox{
      \xymatrix@R=30pt@C=40pt{
        \text{ }
        &\kos{S}
        \ar[d]^-{\kos{f}_!}
        \\
        \kos{T}
        \ar@/^0.7pc/[ur]^-{\kos{f}^*}
        \ar[r]_-{\id_{\kos{T}}}
        \xtwocell[r]{}<>{<-3>{\epsilon}}
        &\kos{T}
      }
    }}
  \end{equation*}
  Moreover, the associated $\kos{T}$-equivariance
  $\vecar{\kos{f}}_!^{\kos{T}}$
  of the left adjoint 
  is equal to the natural transformation
  $\varphi_{\slot,\slot}$
  defined in Definition~\ref{def RSKtensoradj projformula}.
  Let $X\in\obj{\CT}$, $\pzc{Z}\in\obj{\CS}$.
  \begin{equation*}
    \vcenter{\hbox{
      \xymatrix@C=50pt{
        \kos{f}_!(X\bcts_{\kos{f}^*}\pzc{Z})
        \ar@{=}[d]
        \ar[r]^-{(\vecar{\kos{f}}_!^{\kos{T}})_{X,\pzc{Z}}}
        &X\bcts_{\id_{\kos{T}}}\kos{f}_!(\pzc{Z})
        \ar@{=}[d]
        \\
        \kos{f}_!(\pzc{Z}\ctimes \kos{f}^*(X))
        \ar[r]^-{\varphi_{\pzc{Z},X}}
        &\kos{f}_!(\pzc{Z})\tensor X
      }
    }}
  \end{equation*}
\end{remark}

Recall the category
$\cat{Ens}(\kos{K})$
of cocommutative comonoids in $\kos{K}$
introduced in (\ref{eq Ens(T) Ens(K)def}).
For each object $c$ in $\cat{Ens}(\kos{K})$,
we have a symmetric monoidal category
\begin{equation*}
  \kos{K}_c=(\CK_c,\otimes_c,\kos{c}^*(\kappa))
\end{equation*}
and a RSSM adjunction
$\kos{c}:\kos{K}_c\to \kos{K}$
satisfying the projection formula.
\begin{equation*}
  \vcenter{\hbox{
    \xymatrix{
      \kos{K}
      \ar@/^1pc/[d]^-{\kos{c}^*}
      \\
      \kos{K}_c
      \ar@/^1pc/[u]^-{\kos{c}_!}
    }
  }}
\end{equation*}
\begin{itemize}
  \item 
  The underlying category $\CK_c$ of $\kos{K}_c$
  is the coKleisli category associated to the comonad
  $\langle c\otimes\rangle=(c\otimes,\delta^{c\otimes},\varepsilon^{c\otimes})$
  on $\CK$.
  An object in $\CK_c$ is denoted as $\kos{c}^*(x)$ where $x$ is an object in $\CK$.
  A morphism $\kos{c}^*(x)\xrightarrow{l} \kos{c}^*(y)$ in $\CK_c$
  is a morphism $c\otimes x\xrightarrow{l} y$ in $\CK$.
  The composition of morphisms
  $\kos{c}^*(x)\xrightarrow{l}\kos{c}^*(y)\xrightarrow{l^{\pr}}\kos{c}^*(z)$
  in $\CK_c$ is given by
  $c\otimes x
  \xrightarrow{\delta^{c\otimes}_x}
  c\otimes (c\otimes x)
  \xrightarrow{I_c\otimes l}
  c\otimes y
  \xrightarrow{l^{\pr}}
  z$
  and the identity morphism of $\kos{c}^*(x)$ in $\CK_c$ is
  $I_{\kos{c}^*(x)}=\varepsilon^{c\otimes}_x:c\otimes x\to x$.

  \item
  We have an adjunction
  $\kos{c}_!:\CK_c\rightleftarrows \CK:\kos{c}^*$.
  The right adjoint $\kos{c}^*$ sends each object $x$ in $\CK$ to $\kos{c}^*(x)\in\obj{\CK_c}$.
  The left adjoint $\kos{c}_!$ sends each object $\kos{c}^*(z)$ in $\CK_c$
  to $\kos{c}_!\kos{c}^*(z)=c\otimes z\in\obj{\CK}$.
  The component of the adjunction counit at $x\in\obj{\CK}$
  is $\varepsilon^{c\otimes}_x:\kos{c}_!\kos{c}^*(x)=c\otimes x\to x$.
  The component $\kos{c}^*(z)\to \kos{c}^*\kos{c}_!\kos{c}^*(z)$
  of the adjunction unit at $\kos{c}^*(z)\in\obj{\CK_c}$
  is $I_{c\otimes z}:c\otimes z\xrightarrow{\cong}c\otimes z$.
  
  \item
  The symmetric monoidal category structure of
  $\kos{K}_c$
  is induced from that of $\kos{K}$
  via the right adjoint $\kos{c}^*$.
  For instance, the monoidal product in $\kos{K}_c$ is given by
  $\kos{c}^*(x)\otimes_c\kos{c}^*(y)=\kos{c}^*(x\otimes y)$
  where $x$, $y\in\obj{\CK}$.
  The right adjoint is then a strong symmetric monoidal functor
  $\kos{c}^*:\kos{K}\to \kos{K}_c$
  whose coherence isomorphisms are identity morphisms.
  The adjunction $\kos{c}_!\dashv \kos{c}^*$
  becomes a RSSM adjunction $\kos{c}:\kos{K}_c\to \kos{K}$,
  and the colax symmetric monoidal coherences of
  $\kos{c}_!:\kos{K}_c\to \kos{K}$ are described below.
  \begin{equation*}
    \vcenter{\hbox{
      \xymatrix@C=30pt{
        \kos{c}_!(\kos{c}^*(x)\otimes_c\kos{c}^*(y))
        \ar@{=}[d]
        \ar[r]^-{(\kos{c}_!)_{\kos{c}^*(x),\kos{c}^*(y)}}
        &\kos{c}_!\kos{c}^*(x)\otimes \kos{c}_!\kos{c}^*(y)
        \ar@{=}[d]
        &\kos{c}_!\kos{c}^*(\kappa)
        \ar@{=}[d]
        \ar[r]^-{(\kos{c}_!)_{\kos{c}^*(\kappa)}}
        &\kappa
        \ar@{=}[d]
        \\
        c\otimes (x\otimes y)
        \ar[r]^-{(c\otimes)_{x,y}}
        &(c\otimes x)\otimes (c\otimes y)
        &c\otimes \kappa
        \ar[r]^-{(c\otimes)_{\kappa}}
        &\kappa
      }
    }}
  \end{equation*}

  \item
  The RSSM adjunction $\kos{c}:\kos{K}_c\to \kos{K}$
  satisfies the projection formula,
  as we can see from the diagram below.
  Let $x\in\obj{\CK}$ and $\kos{c}^*(z)\in \obj{\CK_c}$.
  \begin{equation*}
    \vcenter{\hbox{
      \xymatrix@C=40pt{
        \kos{c}_!(\kos{c}^*(z)\otimes_c\kos{c}^*(x))
        \ar[d]^-{\kos{c}_{!\kos{c}^*(z),\kos{c}^*(x)}}
        \ar@/_2pc/@<-5ex>[dd]_-{\varphi_{\kos{c}^*(z),x}}
        \ar@{=}[r]
        &c\otimes (z\otimes x)
        \ar[d]^-{(c\otimes)_{z,x}}
        \ar@<6ex>@/^2pc/[dd]^-{a_{c,z,x}}_-{\cong}
        \\
        \kos{c}_!\kos{c}^*(z)\otimes \kos{c}_!\kos{c}^*(x)
        \ar[d]^-{I_{\kos{c}_!\kos{c}^*(z)}\otimes \varepsilon^{c\otimes}_x}
        \ar@{=}[r]
        &(c\otimes z)\otimes (c\otimes x)
        \ar[d]^-{I_{c\otimes z}\otimes \varepsilon^{c\otimes}_x}
        \\
        \kos{c}_!\kos{c}^*(z)\otimes x
        \ar@{=}[r]
        &(c\otimes z)\otimes x
      }
    }}
  \end{equation*}
\end{itemize}
Thus we obtain a right-strong $\kos{K}$-tensor adjunction
\begin{equation*}
  \vcenter{\hbox{
    \xymatrix{
      (\kos{K},\id_{\kos{K}})
      \ar@/^1pc/[d]^-{(\kos{c}^*,\what{\kos{c}}^*)}
      \\
      (\kos{K}_c,\kos{c}^*)
      \ar@/^1pc/[u]^-{(\kos{c}_!,\what{\kos{c}}_!)}
    }
  }}
\end{equation*}
as explained in Remark~\ref{rem RSKtensoradj projformula as Tequiv}.
\begin{itemize}
  \item 
  The comonoidal natural isomorphism
  $\what{\kos{c}}^*:\kos{c}^*=\kos{c}^*$
  is the identity natural transformation of $\kos{c}^*$,
  and the comonoidal natural transformation
  $\what{\kos{c}}_!:\kos{c}_!\kos{c}^*\Rightarrow \id_{\kos{K}}$
  is given by $\varepsilon^{c\otimes}:c\otimes \Rightarrow \id_{\kos{K}}$.
  One can check that
  \begin{equation*}
    (\kos{c}_!\kos{c}^*,\what{\kos{c}_!\kos{c}^*})
    =(c\otimes,\what{c\otimes})
    :(\kos{K},\id_{\kos{K}})\to (\kos{K},\id_{\kos{K}})
    .
  \end{equation*}

  \item
  The associated $\kos{K}$-equivariance $\vecar{\kos{c}}_!$
  of the left adjoint
  $(\kos{c}_!,\what{\kos{c}}_!)$
  is a natural isomorphism:
  \begin{equation*}
    \vecar{\kos{c}}_{!x,\kos{c}^*(z)}
    =
    \varphi_{\kos{c}^*(z),x}
    =
    a_{c,z,x}:
    \xymatrix{
      c\otimes (z\otimes x)
      \ar[r]^-{\cong}
      &(c\otimes z)\otimes x
      .
    }
  \end{equation*}
\end{itemize}
Let $f:c^{\pr}\to c$ be a morphism in $\cat{Ens}(\kos{K})$.
We have a strong $\kos{K}$-tensor functor
\begin{equation} \label{eq RSKtensoradj functorbetweencoKelisliCat}
  (\kos{f}^*,\what{\kos{f}}^*):
  (\kos{K}_c,\kos{c}^*)\to (\kos{K}_{c^{\pr}},\kos{c}^{\pr*})
\end{equation}
which satisfies the relation
\begin{equation*}
  \vcenter{\hbox{
    \xymatrix@C=0pt{
      \text{ }
      &(\kos{K},\id_{\kos{K}})
      \ar[dl]_-{(\kos{c}^*,\what{\kos{c}}^*)}
      \ar[dr]^-{(\kos{c}^{\pr*},\what{\kos{c}}^{\pr*})}
      \\
      (\kos{K}_c,\kos{c}^*)
      \ar[rr]^-{(\kos{f}^*,\what{\kos{f}}^*)}
      &\text{ }
      &(\kos{K}_{c^{\pr}},\kos{c}^{\pr*})
    }
  }}
  \qquad\quad
  (\kos{f}^*,\what{\kos{f}}^*)(\kos{c}^*,\what{\kos{c}}^*)
  =
  (\kos{c}^{\pr*},\what{\kos{c}}^{\pr*}).
\end{equation*}
\begin{itemize}
  \item 
  The functor $\kos{f}^*$ sends each object $\kos{c}^*(x)$ in $\CK_c$
  to the object
  $\kos{f}^*\kos{c}^*(x)=\kos{c}^{\pr*}(x)$ in $\CK_{c^{\pr}}$,
  and each morphism $\kos{c}^*(x)\xrightarrow{l} \kos{c}^*(y)$ in $\CK_c$
  to the morphism $\kos{c}^{\pr*}(x)\xrightarrow{\kos{f}^*(l)} \kos{c}^{\pr*}(y)$ in $\CK_{c^{\pr}}$
  where $\kos{f}^*(l):c^{\pr}\otimes x\xrightarrow{f\otimes I_x}c\otimes x\xrightarrow{l} y$.

  \item
  The symmetric monoidal coherence isomorphisms of $\kos{f}^*:\kos{K}_c\to \kos{K}_{c^{\pr}}$
  are identity morphisms,
  and the comonoidal natural isomorphism
  $\what{\kos{f}}^*:\kos{f}^*\kos{c}^*=\kos{c}^{\pr*}$
  is the identity natural transformation.
\end{itemize}

\begin{definition} \label{def RSKtensoradj Hom,Isom}
  Let $(\kos{T},\mon{t})$ be a colax $\kos{K}$-tensor category
  and let $(\phi,\what{\phi})$, $(\psi,\what{\psi})
  :(\kos{T},\mon{t})\to (\kos{K},\id_{\kos{K}})$
  be colax $\kos{K}$-tensor functors.
  \begin{enumerate}
    \item 
    We define the \emph{presheaf of comonoidal $\kos{K}$-tensor natural transformations
    from $(\phi,\what{\phi})$ to $(\psi,\what{\psi})$}
    as the presheaf on $\cat{Ens}(\kos{K})$
    \begin{equation*}
      \underline{\Hom}_{\mathbb{SMC}_{\cat{colax}}^{\kos{K}/\!\!/}}
      \bigl(
        (\phi,\what{\phi})
        ,
        (\psi,\what{\psi})
      \bigr)
      :
      \cat{Ens}(\kos{K})^{\op}
      \to
      \cat{Set}
    \end{equation*}
    which sends each object $c$ in $\cat{Ens}(\kos{K})$
    to the set 
    of comonoidal $\kos{K}$-tensor natural transformations
    $(\kos{c}^*\phi,\what{\kos{c}^*\phi})
    \Rightarrow
    (\kos{c}^*\psi,\what{\kos{c}^*\psi})
    :(\kos{T},\mon{t})\to (\kos{K}_c,\kos{c}^*)$.
    
    \item
    We define \emph{the presheaf of comonoidal $\kos{K}$-tensor natural isomorphisms
    from $(\phi,\what{\phi})$ to $(\psi,\what{\psi})$}
    as the presheaf on $\cat{Ens}(\kos{K})$
    \begin{equation*}
      \underline{\Isom}_{\mathbb{SMC}_{\cat{colax}}^{\kos{K}/\!\!/}}
      \bigl(
        (\phi,\what{\phi})
        ,
        (\psi,\what{\psi})
      \bigr)
      :
      \cat{Ens}(\kos{K})^{\op}
      \to
      \cat{Set}
    \end{equation*}
    which sends each object $c$ in $\cat{Ens}(\kos{K})$
    to the set 
    of comonoidal $\kos{K}$-tensor natural isomorphisms
    $\xymatrix@C=12pt{(\kos{c}^*\phi,\what{\kos{c}^*\phi})
    \ar@2{->}[r]^-{\cong}
    &(\kos{c}^*\psi,\what{\kos{c}^*\psi})
    :(\kos{T},\mon{t})\to (\kos{K}_c,\kos{c}^*).}$
  \end{enumerate}
\end{definition}

For the rest of this subsection,
we study about right-strong $\kos{K}$-tensor adjunctions
from $(\kos{K},\id_{\kos{K}})$
to a strong $\kos{K}$-tensor category $(\kos{T},\mon{t})$.

\begin{definition}\label{def RSKtensoradj reflective}
  Let $(\kos{T},\mon{t})$ be a strong $\kos{K}$-tensor category
  and let $(\omega,\what{\omega}):(\kos{K},\id_{\kos{K}})\to (\kos{T},\mon{t})$
  be a right-strong $\kos{K}$-tensor adjunction.
  \begin{equation*}
    \vcenter{\hbox{
      \xymatrix{
        (\kos{T},\mon{t})
        \ar@/^1pc/[d]^-{(\omega^*,\what{\omega}^*)}
        \\
        (\kos{K},\id_{\kos{K}})
        \ar@/^1pc/[u]^-{(\omega_!,\what{\omega}_!)}
      }
    }}
  \end{equation*}
  \begin{enumerate}
    \item
    We say $(\omega,\what{\omega})$
    is \emph{reflective} if the left adjoint
    $(\omega_!,\what{\omega}_!):(\kos{K},\id_{\kos{K}})\to (\kos{T},\mon{t})$
    is a reflective colax $\kos{K}$-tensor functor.
    This amounts to saying that
    the associated $\kos{K}$-equivariance
    $\vecar{\omega}_!$
    of the left adjoint
    is a natural isomorphism:
    see Lemma~\ref{lem Ens(T) equivariance and reflection}.

    \item 
    We define the natural transformation
    \begin{equation*}
      \varphi_{\slot}:
      \omega_!\omega^*\Rightarrow \omega_!(\kappa)\tensor\slot:\CT\to \CT
    \end{equation*}
    whose component at $X\in\obj{\CT}$ is
    \begin{equation} \label{eq RSKtensoradj modifiedvarphi def}
      \varphi_X:
      \xymatrix@C=30pt{
        \omega_!\omega^*(X)
        \ar[r]^-{\omega_!(\imath_{\omega^*(X)})}_-{\cong}
        &\omega_!(\kappa\otimes \omega^*(X))
        \ar[r]^-{\varphi_{\kappa,X}}
        &\omega_!(\kappa)\tensor X
        .
      }
    \end{equation}
  \end{enumerate}
\end{definition}

\begin{lemma} \label{lem RSKtensoradj phi-modifiedphi}
  Let $(\kos{T},\mon{t})$
  be a strong $\kos{K}$-tensor category
  and let
  $(\omega,\what{\omega}):(\kos{K},\id_{\kos{K}})\to (\kos{T},\mon{t})$
  be a reflective right-strong $\kos{K}$-tensor adjunction.
  \begin{equation*}
    \vcenter{\hbox{
      \xymatrix{
        (\kos{T},\mon{t})
        \ar@/^1pc/[d]^-{(\omega^*,\what{\omega}^*)}
        \\
        (\kos{K},\id_{\kos{K}})
        \ar@/^1pc/[u]^-{\text{ reflective }(\omega_!,\what{\omega}_!)}
      }
    }}
  \end{equation*}
  For each object $X$ in $\CT$, the following are equivalent:
  \begin{enumerate}
    \item 
    $\varphi_{\slot,X}
    :\omega_!(\slot\otimes \omega^*(X))
    \Rightarrow\omega_!(\slot)\tensor X:\CK\to \CT$
    is a natural isomorphism;

    \item
    $\varphi_X:\omega_!\omega^*(X)\to \omega_!(\kappa)\tensor X$
    is an isomorphism in $\CT$.
  \end{enumerate}
  In particular,
  $(\omega,\what{\omega})$ satisfies the projection formula
  if and only if
  $\varphi_{\slot}:\omega_!\omega^*\Rightarrow \omega_!(\kappa)\tensor\slot:\CT\to \CT$
  is a natural isomorphism.
\end{lemma}
\begin{proof}
  Statement 1 implies statemet 2
  as we can see from the definition (\ref{eq RSKtensoradj modifiedvarphi def})
  of $\varphi_X:\omega_!\omega^*(X)\to \omega_!(\kappa)\tensor X$.
  To prove the converse,
  it suffices to show that
  for each object $z$ in $\CK$,
  the following diagram commutes.
  \begin{equation} \label{eq RSKtensoradj phi}
    \vcenter{\hbox{
      \xymatrix@C=40pt{
        \omega_!(z\otimes \omega^*(X))
        \ar[d]_-{\varphi_{z,X}}
        \ar[r]^-{\omega_!(s_{z,\omega^*(X)})}_-{\cong}
        &\omega_!(\omega^*(X)\otimes z)
        \ar[r]^-{(\vecar{\omega}_!)_{z,\omega^*(X)}}_-{\cong}
        &\omega_!\omega^*(X)\tensor \mon{t}(z)
        \ar[d]^-{\varphi_X\tensor I_{\mon{t}(z)}}
        \\
        \omega_!(z)\tensor X
        \ar[r]^-{(\hatar{\omega}_!)_z\tensor I_X}_-{\cong}
        &\omega_!(\kappa)\tensor \mon{t}(z)\tensor X
        \ar[r]^-{I_{\omega_!(\kappa)}\tensor s_{\mon{t}(z),X}}_-{\cong}
        &\omega_!(\kappa)\tensor X\tensor \mon{t}(z)
      }
    }}
  \end{equation}
  We can check this as follows.
  First, we have the following relation
  \begin{equation*}
    (\dagger):
    \!\!\!
    \!\!\!
    \!\!\!
    \!\!\!
    \vcenter{\hbox{
      \xymatrix@C=15pt{
        \omega_!(\omega^*(X)\otimes z)
        \ar[dd]_-{(\vecar{\omega}_!)_{z,\omega^*(X)}}^-{\cong}
        \ar@{=}[r]
        &\omega_!(\omega^*(X)\otimes z)
        \ar[d]^-{(\omega_!)_{\omega^*(X),z}}
        \ar@{=}[r]
        &\omega_!(\omega^*(X)\otimes z)
        \ar[d]^-{\omega_!(\imath_{\omega^*(X)}\otimes I_z)}_-{\cong}
        \\
        \text{ }
        &\omega_!\omega^*(X)\tensor \omega_!(z)
        \ar@/^0.5pc/[dl]|-{I_{\omega_!\omega^*(X)\tensor \what{\omega}_{!z}}}
        \ar[ddd]_-{\varphi_X\tensor I_{\omega_!(z)}}
        \ar@/_0.5pc/[dr]_-{\omega_!(\imath_{\omega^*(X)})\tensor I_{\omega_!(z)}}^-{\cong}
        &\omega_!(\kappa\otimes \omega^*(X)\otimes z)
        \ar[d]^-{\omega_{!\kappa\otimes \omega^*(X),z}}
        \ar@/^2.5pc/@<9ex>[dd]|-{\omega_{!\kappa,\omega^*(X),z}}
        \\
        \omega_!\omega^*(X)\tensor \mon{t}(z)
        \ar[ddd]_-{\varphi_X\tensor I_{\mon{t}(z)}}
        &\text{ }
        &\omega_!(\kappa\otimes \omega^*(X))\tensor \omega_!(z)
        \ar[d]^-{\omega_{!\kappa,\omega^*(X)}\tensor I_{\omega_!(z)}}
        \\
        \text{ }
        &\text{ }
        &\omega_!(\kappa)\tensor \omega_!\omega^*(X)\tensor \omega_!(z)
        \ar@/_0.5pc/[dl]|-{I_{\omega_!(\kappa)}\tensor \epsilon_X\tensor \omega_!(z)}
        \ar[dd]^-{I_{\omega_!(\kappa)}\tensor \epsilon_X\tensor \what{\omega}_{!z}}
        \\
        \text{ }
        &\omega_!(\kappa)\tensor X\tensor \omega_!(z)
        \ar@/_0.5pc/[dr]|-{I_{\omega_!(\kappa)\tensor X}\tensor \what{\omega}_{!z}}
        &\text{ }
        \\
        \omega_!(\kappa)\tensor X\tensor \mon{t}(z)
        \ar@{=}[rr]
        &\text{ }
        &\omega_!(\kappa)\tensor X\tensor \mon{t}(z)
      }
    }}
  \end{equation*}
  where the notation $\omega_{!\kappa,\omega^*(X),z}$
  means the colax symmetric monoidal coherence of $\omega_!:\kos{K}\to\kos{T}$
  with respect to three input objects
  $\kappa$, $\omega^*(X)$, $z$ in $\CK$.
  We also have the following relation
  \begin{equation*}
    (\dagger\!\dagger):
    \!\!\!
    \!\!\!
    \!\!\!
    \!\!\!
    \vcenter{\hbox{
      \xymatrix@C=15pt{
        \omega_!(z\otimes \omega^*(X))
        \ar[dd]_-{\varphi_{z,X}}
        \ar@{=}[rr]
        &\text{ }
        &\omega_!(z\otimes \omega^*(X))
        \ar[d]^-{\omega_!(\imath_z\otimes I_{\omega^*(X)})}_-{\cong}
        \\
        \text{ }
        &\text{ }
        &\omega_!(\kappa\otimes z\otimes \omega^*(X))
        \ar[d]^-{\omega_{!\kappa\otimes z,\omega^*(X)}}
        \ar@/_2.5pc/[ddl]|-{\varphi_{\kappa\otimes z,X}}
        \ar@/^2.5pc/@<9ex>[dd]|-{\omega_{!\kappa,z,\omega^*(X)}}
        \\
        \omega_!(z)\tensor X
        \ar[ddd]_-{\hatar{\omega}_{!z}\tensor I_X}^-{\cong}
        \ar@/_0.5pc/[dr]_-{\omega_!(\imath_z)\tensor I_X}^-{\cong}
        &\text{ }
        &\omega_!(\kappa\otimes z)\tensor \omega_!\omega^*(X)
        \ar[d]^-{\omega_{!\kappa,z}\tensor I_{\omega_!\omega^*(X)}}
        \ar@/_0.5pc/[dl]|-{I_{\omega_!(\kappa\otimes z)}\tensor \epsilon_X}
        \\
        \text{ }
        &\omega_!(\kappa\otimes z)\tensor X
        \ar[d]^-{\omega_{!\kappa,z}\tensor I_X}
        &\omega_!(\kappa)\tensor \omega_!(z)\tensor \omega_!\omega^*(X)
        \ar[dd]^-{I_{\omega_!(\kappa)}\tensor \what{\omega}_{!z}\tensor \epsilon_X}
        \ar@/^0.5pc/[dl]|-{I_{\omega_!(\kappa)\tensor \omega_!(z)}\tensor \epsilon_X}
        \\
        \text{ }
        &\omega_!(\kappa)\tensor \omega_!(z)\tensor X
        \ar@/_0.5pc/[dr]|-{I_{\omega_!(\kappa)}\tensor \what{\omega}_{!z}\tensor I_X}
        &\text{ }
        \\
        \omega_!(\kappa)\tensor \mon{t}(z)\tensor X
        \ar@{=}[rr]
        &\text{ }
        &\omega_!(\kappa)\tensor \mon{t}(z)\tensor X
      }
    }}
  \end{equation*}
  where the notation $\omega_{!\kappa,z,\omega^*(X)}$
  means the colax symmetric monoidal coherence of $\omega_!:\kos{K}\to\kos{T}$
  with respect to three input objects
  $\kappa$, $z$, $\omega^*(X)$ in $\CK$.
  Using these relations,
  we verify the diagram (\ref{eq RSKtensoradj phi}) as follows.
  \begin{equation*}
    \vcenter{\hbox{
      \xymatrix@C=0pt{
        \omega_!(z\otimes \omega^*(X))
        \ar[d]^-{\omega_!(s_{z,\omega^*(X)})}_-{\cong}
        \ar@{=}[rr]
        &\text{ }
        &\omega_!(z\otimes \omega^*(X))
        \ar@<-1ex>[d]_-{\omega_!(\imath_z\otimes I_{\omega^*(X)})}^-{\cong}
        \ar@{=}[r]
        &\omega_!(z\otimes \omega^*(X))
        \ar[d]^-{\varphi_{z,X}}
        \\
        \omega_!(\omega^*(X)\otimes z)
        \ar[dd]^(0.66){(\vecar{\omega}_!)_{z,\omega^*(X)}}_(0.66){\cong}
        \ar@/_0.5pc/[dr]^-{\omega_!(\imath_{\omega^*(X)}\otimes I_z)}_-{\cong}
        &\text{ }
        &\omega_!(\kappa\otimes z\otimes \omega^*(X))
        \ar@/_0.5pc/[dl]_-{\omega_!(I_{\kappa}\otimes s_{z,\omega^*(X)})}^-{\cong}
        \ar@<-1ex>[d]^-{\omega_{!\kappa,z,\omega^*(X)}}
        \ar@{}[r]|-{(\dagger\!\dagger)}
        &\omega_!(z)\tensor X
        \ar[dd]^-{\hatar{\omega}_{!z}\tensor I_X}_-{\cong}
        \\
        \text{ }
        &\omega_!(\kappa\otimes \omega^*(X)\otimes z)
        \ar@/_1.5pc/[dr]|-{\omega_{!\kappa,\omega^*(X),z}}
        &\omega_!(\kappa)\!\tensor\! \omega_!(z)\!\tensor\! \omega_!\omega^*(X)
        \ar@/^0.7pc/[dr]|-{I_{\omega_!(\kappa)}\tensor \what{\omega}_{!z}\tensor \epsilon_X}
        \ar@<-1ex>[d]|-{I_{\omega_!(\kappa)}\tensor s_{\omega_!(z),\omega_!\omega^*(X)}}
        &\text{ }
        \\
        \omega_!\omega^*(X)\tensor \mon{t}(z)
        \ar[d]^-{\varphi_X\tensor I_{\mon{t}(z)}}
        \ar@{}[rr]|-{(\dagger)}
        &\text{ }
        &\omega_!(\kappa)\!\tensor\! \omega_!\omega^*(X)\!\tensor\! \omega_!(z)
        \ar@<-1ex>[d]|-{I_{\omega_!(\kappa)}\tensor \epsilon_X\tensor \what{\omega}_{!z}}
        &\omega_!(\kappa)\tensor \mon{t}(z)\tensor X
        \ar[d]^-{I_{\omega_!(\kappa)}\tensor s_{\mon{t}(z),X}}_-{\cong}
        \\
        \omega_!(\kappa)\tensor X\tensor \mon{t}(z)
        \ar@{=}[rr]
        &\text{ }
        &\omega_!(\kappa)\tensor X\tensor \mon{t}(z)
        \ar@{=}[r]
        &\omega_!(\kappa)\tensor X\tensor \mon{t}(z)
      }
    }}
  \end{equation*}
  This completes the proof of Lemma~\ref{lem RSKtensoradj phi-modifiedphi}.
\qed\end{proof}

The following proposition is the main result of this subsection.

\begin{proposition} \label{prop RSKtensoradj HomRep,Hom=Isom}
  Let $(\kos{T},\mon{t})$
  be a strong $\kos{K}$-tensor category.
  Suppose we are given a reflective right-strong $\kos{K}$-tensor adjunction
  $(\omega,\what{\omega})
  :(\kos{K},\id_{\kos{K}})\to (\kos{T},\mon{t})$
  and a strong $\kos{K}$-tensor functor
  $(\omega^{\pr*},\what{\omega}^{\pr*})
  :(\kos{T},\mon{t})\to (\kos{K},\id_{\kos{K}})$.
  \begin{equation*}
    \vcenter{\hbox{
      \xymatrix{
        (\kos{T},\mon{t})
        \ar@/^1pc/[d]^-{(\omega^*,\what{\omega}^*)}
        \\
        (\kos{K},\id_{\kos{K}})
        \ar@/^1pc/[u]^-{\text{ reflective }(\omega_!,\what{\omega}_!)}
      }
    }}
    \qquad\qquad
    \vcenter{\hbox{
      \xymatrix{
        (\kos{T},\mon{t})
        \ar[d]^-{(\omega^{\pr*},\what{\omega}^{\pr*})}
        \\
        (\kos{K},\id_{\kos{K}})
      }
    }}
  \end{equation*}
  \begin{enumerate}
    \item 
    The presheaf of comonoidal $\kos{K}$-tensor natural transformations
    from
    $(\omega^*,\what{\omega}^*)$
    to
    $(\omega^{\pr*},\what{\omega}^{\pr*})$
    is represented by the object
    $\omega^{\pr*}\omega_!(\kappa)$ in $\cat{Ens}(\kos{K})$.
    \begin{equation*}
      \hspace{-0.5cm}
      \xymatrix@C=12pt{
        \Hom_{\cat{Ens}(\kos{K})}\bigl(\slot, \omega^{\pr*}\omega_!(\kappa)\bigr)
        \ar@2{->}[r]^-{\cong}
        &\underline{\Hom}_{\mathbb{SMC}_{\cat{colax}}^{\kos{K}/\!\!/}}
        \!\bigl(
          (\omega^*,\what{\omega}^*)
          ,
          (\omega^{\pr*},\what{\omega}^{\pr*})
        \bigr)
        \!:\!
        \cat{Ens}(\kos{K})^{\op}
        \to
        \cat{Set}
      }
    \end{equation*}
    The universal element is the comonoidal $\kos{K}$-tensor natural transformation
    \begin{equation*}
      \begin{aligned}
        \xi
        &:
        \xymatrix{
          (\kos{p}^*\omega^*,\what{\kos{p}^*\omega^*})
          \Rightarrow
          (\kos{p}^*\omega^{\pr*},\what{\kos{p}^*\omega^{\pr*}})
          ,
        }
        \qquad
        p:=\omega^{\pr*}\omega_!(\kappa)
      \end{aligned}
    \end{equation*}
    whose component at each object $X$ in $\CT$ is
    \begin{equation*}
      \begin{aligned}
        \xi_X
        &:
        \kos{p}^*\omega^*(X)\to \kos{p}^*\omega^{\pr*}(X)
        ,
        \\
        \xi_X
        &:
        \xymatrix@C=40pt{
          \omega^{\pr*}\omega_!(\kappa)\otimes \omega^*(X)
          \ar[r]^-{(\hatar{\omega^{\pr*}\omega_!}_{\omega^*(X)})^{-1}}_-{\cong}
          &\omega^{\pr*}\omega_!\omega^*(X)
          \ar[r]^-{\omega^{\pr*}(\epsilon_X)}
          &\omega^{\pr*}(X)
          .
        }
      \end{aligned}
    \end{equation*}

    \item
    Let $X$ be an object in $\CT$.
    The following are equivalent:
    \begin{enumerate}
      \item 
      $\omega^{\pr*}(\varphi_X):
      \omega^{\pr*}\omega_!\omega^*(X)
      \to
      \omega^{\pr*}(\omega_!(\kappa)\tensor X)$
      is an isomorphism in $\CT$.
      
      \item 
      The component
      $\xi_X:\kos{p}^*\omega^*(X)\to \kos{p}^*\omega^{\pr*}(X)$
      of the universal element $\xi$ at $X$
      is an isomorphism in the coKleisli category $\CK_p$.
      
      \item
      For each object $c$ in $\cat{Ens}(\kos{K})$
      and for each comonoidal $\kos{K}$-tensor natural transformation
      $\vartheta:
      (\kos{c}^*\omega^*,\what{\kos{c}^*\omega^*})
      \Rightarrow
      (\kos{c}^*\omega^{\pr*},\what{\kos{c}^*\omega^{\pr*}})
      :(\kos{T},\mon{t})\to (\kos{K}_c,\kos{c}^*)$,
      the component
      $\vartheta_X:\kos{c}^*\omega^*(X)\to \kos{c}^*\omega^{\pr*}(X)$
      of $\vartheta$ at $X$ is an isomorphism
      in the coKleisli category $\CK_c$.
    \end{enumerate}
    In this case, if we denote
    \begin{equation*}
      \check{\xi}_X:
      \xymatrix@C=30pt{
        \omega^{\pr*}\omega_!\omega^{\pr*}(X)
        \ar[r]^-{\hatar{\omega^{\pr*}\omega_!}_{\omega^{\pr*}(X)}}_-{\cong}
        &\omega^{\pr*}\omega_!(\kappa)\otimes \omega^{\pr*}(X)
        \ar[r]^-{\xi_X^{-1}}
        &\omega^*(X)
      }
    \end{equation*}
    then the inverse $\vartheta_X^{-1}$ is explicitly described as
    \begin{equation*}
      \vartheta_X^{-1}:\!\!
      \xymatrix@C=28pt{
        c\otimes \omega^{\pr*}(X)
        \ar[r]^-{I_c\otimes \eta_{\omega^{\pr*}(X)}}
        &c\otimes \omega^*\omega_!\omega^{\pr*}(X)
        \ar[r]^-{\vartheta_{\omega_!\omega^{\pr*}(X)}}
        &\omega^{\pr*}\omega_!\omega^{\pr*}(X)
        \ar[r]^-{\check{\xi}_X}
        &\omega^*(X)
        .
      }
    \end{equation*}

    \item
    We have
    \begin{equation*}
      \underline{\Hom}_{\mathbb{SMC}_{\cat{colax}}^{\kos{K}/\!\!/}}
      \bigl(
        (\omega^*,\what{\omega}^*)
        ,
        (\omega^{\pr*},\what{\omega}^{\pr*})
      \bigr)
      =
      \underline{\Isom}_{\mathbb{SMC}_{\cat{colax}}^{\kos{K}/\!\!/}}
      \bigl(
        (\omega^*,\what{\omega}^*)
        ,
        (\omega^{\pr*},\what{\omega}^{\pr*})
      \bigr)
    \end{equation*}
    if and only if
    $\omega^{\pr*}(\varphi_{\slot}):
    \omega^{\pr*}\omega_!\omega^*
    \Rightarrow
    \omega^{\pr*}(\omega_!(\kappa)\tensor \slot)
    :\CT\to \CK$
    is a natural isomorphism.
  \end{enumerate}
\end{proposition}
\begin{proof}
  We first prove statement 1.
  For each object $c$ in $\cat{Ens}(\kos{K})$,
  we have bijections
  \begin{equation} \label{eq1 RSKtensoradj HomRep,Hom=Isom}
    \begin{aligned}
      \Hom_{\cat{Ens}(\kos{K})}\bigl(c,\text{ }\omega^{\pr*}\omega_!(\kappa)\bigr)
      &\cong
      \Hom_{\cat{Ens}(\kos{K})}\bigl(\kos{c}_!\kos{c}^*(\kappa),\text{ }\omega^{\pr*}\omega_!(\kappa)\bigr)
      \\
      &\cong
      \mathbb{SMC}_{\cat{colax}}^{\kos{K}/\!\!/}\bigl((\kos{c}_!\kos{c}^*,\what{\kos{c}_!\kos{c}^*}),\text{ }(\omega^{\pr*}\omega_!,\what{\omega^{\pr*}\omega_!})\bigr)
      \\
      &\cong
      \mathbb{SMC}_{\cat{colax}}^{\kos{K}/\!\!/}\bigl((\kos{c}^*\omega^*,\what{\kos{c}^*\omega^*}),\text{ }(\kos{c}^*\omega^{\pr*},\what{\kos{c}^*\omega^{\pr*}})\bigr)
      .
    \end{aligned}
  \end{equation}
  The first bijection in (\ref{eq1 RSKtensoradj HomRep,Hom=Isom})
  is obtained by precomposing the inverse of the isomorphism
  $\jmath_c:c\xrightarrow{\cong}c\otimes \kappa=\kos{c}_!\kos{c}^*(\kappa)$
  in $\cat{Ens}(\kos{K})$.
  Recall the adjoint equivalence of categories
  \begin{equation*}
    \CL:
    \mathbb{SMC}_{\cat{colax}}^{\kos{K}/\!\!/}\bigl((\kos{K},\id_{\kos{K}}),(\kos{K},\id_{\kos{K}})\bigr)_{\cat{rfl}}
    \simeq
    \cat{Ens}(\kos{K})
    :\iota
  \end{equation*}
  described in Corollary~\ref{cor Ens(T) strongKcat reflective adjunction}.
  As the colax $\kos{K}$-tensor functors
  \begin{equation*}
    (\kos{c}_!\kos{c}^*,\what{\kos{c}_!\kos{c}^*}),
    (\omega^{\pr*}\omega_!,\what{\omega^{\pr*}\omega_!})
    :(\kos{K},\id_{\kos{K}})\to (\kos{K},\id_{\kos{K}})
  \end{equation*}
  are both reflective,
  we obtain the second bijection in (\ref{eq1 RSKtensoradj HomRep,Hom=Isom}).
  The last bijection in (\ref{eq1 RSKtensoradj HomRep,Hom=Isom})
  is the bijective correspondence of mates
  in the $2$-category $\mathbb{SMC}_{\cat{colax}}^{\kos{K}/\!\!/}$.
  \begin{equation*}
    \vcenter{\hbox{
      \xymatrix@R=30pt@C=40pt{
        (\kos{K},\id_{\kos{K}})
        \ar[r]^-{(\omega_!,\what{\omega}_!)}
        \ar[d]_-{(\kos{c}^*,\what{\kos{c}}^*)}
        &(\kos{T},\mon{t})
        \ar[d]^-{(\omega^{\pr*},\what{\omega}^{\pr*})}
        \\
        (\kos{K}_c,\kos{c}^*)
        \ar[r]_-{(\kos{c}_!,\what{\kos{c}}_!)}
        &(\kos{K},\id_{\kos{K}})
        \xtwocell[ul]{}<>{<0>{\vartheta_!\text{ }\text{ }}}
      }
    }}
    \quad\leftrightsquigarrow\quad
    \vcenter{\hbox{
      \xymatrix@R=30pt@C=40pt{
        (\kos{K},\id_{\kos{K}})
        \ar[d]_-{(\kos{c}^*,\what{\kos{c}}^*)}
        &(\kos{T},\kos{t})
        \ar[l]_-{(\omega^*,\what{\omega}^*)}
        \ar[d]^-{(\omega^{\pr*},\what{\omega}^{\pr*})}
        \\
        (\kos{K}_c,\kos{c}^*)
        \xtwocell[ur]{}<>{<0>{\text{ }\vartheta}}
        &(\kos{K},\id_{\kos{K}})
        \ar[l]^-{(\kos{c}^*,\what{\kos{c}}^*)}
      }
    }}
  \end{equation*}
  The mate pair $\vartheta_!$ and $\vartheta$
  are related as follows.
  \begin{equation*}
    \begin{aligned}
      \vartheta_X
      &:
      \xymatrix@C=30pt{
        c\otimes \omega^*(X)
        \ar[r]^-{(\vartheta_!)_{\omega^*(X)}}
        &\omega^{\pr*}\omega_!\omega^*(X)
        \ar[r]^-{\omega^{\pr*}(\epsilon_X)}
        &\omega^{\pr*}(X)
      }
      ,
      \qquad
      X\in\obj{\CT}
      ,
      \\
      (\vartheta_!)_x
      &:
      \xymatrix@C=30pt{
        c\otimes x
        \ar[r]^-{I_c\otimes \eta_x}
        &c\otimes \omega^*\omega_!(x)
        \ar[r]^-{\vartheta_{\omega_!(x)}}
        &\omega^{\pr*}\omega_!(x)
      }   
      ,
      \qquad
      x\in\obj{\CK}
      .
    \end{aligned}
  \end{equation*}
  To complete the proof of statement 1,
  it suffices to show that the bijection (\ref{eq1 RSKtensoradj HomRep,Hom=Isom})
  sends each morphism $g:c\to \omega^{\pr*}\omega_!(\kappa)$ in $\cat{Ens}(\kos{K})$
  to the comonoidal $\kos{K}$-tensor natural transformation
  $\vartheta:(\kos{c}^*\omega^*,\what{\kos{c}^*\omega^*})\Rightarrow (\kos{c}^*\omega^{\pr*},\what{\kos{c}^*\omega^{\pr*}})$
  whose component at each object $X$ in $\CT$ is
  \begin{equation*}
    \vartheta_X:
    \xymatrix@C=30pt{
      c\otimes \omega^*(X)
      \ar[r]^-{g\otimes I_{\omega^*(X)}}
      &\omega^{\pr*}\omega_!(\kappa)\otimes \omega^*(X)
      \ar[r]^-{\xi_X}
      &\omega^{\pr*}(X)
      .
    }
  \end{equation*}
  We can check this as follows.
  The component of the corresponding colax $\kos{K}$-tensor natural transformation
  $\vartheta_!:(\kos{c}_!\kos{c}^*,\what{\kos{c}_!\kos{c}^*})\Rightarrow (\omega^{\pr*}\omega_!,\what{\omega^{\pr*}\omega_!})$
  at each object $x$ in $\CK$ is
  \begin{equation*}
    (\vartheta_!)_x:
      \xymatrix@C=30pt{ 
      c\otimes x
      \ar[r]^-{g\otimes I_x}
      &\omega^{\pr*}\omega_!(\kappa)\otimes x
      \ar[r]^-{(\hatar{\omega^{\pr*}\omega_!}_x)^{-1}}_-{\cong}
      &\omega^{\pr*}\omega_!(x)
      .
    }
  \end{equation*}
  Therefore the component of $\vartheta$ at $X$ is given as follows.
  \begin{equation*}
    \vcenter{\hbox{
      \xymatrix@C=35pt{
        c\otimes \omega^*(X)
        \ar[dd]^-{(\vartheta_!)_{\omega^*(X)}}
        \ar@/_2pc/@<-3ex>[ddd]_-{\vartheta_X}
        \ar@{=}[r]
        &c\otimes \omega^*(X)
        \ar[d]^-{g\otimes I_{\omega^*(X)}}
        \\
        \text{ }
        &\omega^{\pr*}\omega_!(\kappa)\otimes \omega^*(X)
        \ar@/^0.5pc/[dl]^(0.45){(\hatar{\omega^{\pr*}\omega_!}_{\omega^*(X)})^{-1}}_-{\cong}
        \ar[dd]^-{\xi_X}
        \\
        \omega^{\pr*}\omega_!\omega^*(X)
        \ar[d]^-{\omega^{\pr*}(\epsilon_X)}
        &\text{ }
        \\
        \omega^{\pr*}(X)
        \ar@{=}[r]
        &\omega^{\pr*}(X)
      }
    }}
  \end{equation*}
  This completes the proof of statement 1.
  We mention here that under the bijection (\ref{eq1 RSKtensoradj HomRep,Hom=Isom}),
  each comonoidal $\kos{K}$-tensor natural transformation
  $\vartheta:(\kos{c}^*\omega^*,\what{\kos{c}^*\omega^*})\Rightarrow (\kos{c}^*\omega^{\pr*},\what{\kos{c}^*\omega^{\pr*}})$
  corresponds to the following morphism in $\cat{Ens}(\kos{K})$.
  \begin{equation*}
    g:
    \xymatrix{
      c
      \ar[r]^-{\jmath_c}_-{\cong}
      &c\otimes \kappa
      \ar[r]^-{I_c\otimes \eta_{\kappa}}
      \ar@/_1pc/@<-1ex>[rr]|-{(\vartheta_!)_{\kappa}}
      &c\otimes \omega^*\omega_!(\kappa)
      \ar[r]^-{\vartheta_{\omega_!(\kappa)}}
      &\omega^{\pr*}\omega_!(\kappa)
    }
  \end{equation*}

  We are left to prove statement 2,
  as it implies statment 3.
  To do so, we first show that the morphism $\kos{p}_!(\xi_X)$ in $\CK$
  is explicitly given as follows.
  \begin{equation}\label{eq2 RSKtensoradj HomRep,Hom=Isom}
    \vcenter{\hbox{
      \xymatrix@C=40pt{
        \omega^{\pr*}\omega_!(\kappa)\otimes \omega^*(X)
        \ar[r]^-{(\hatar{\omega^{\pr*}\omega_!}_{\omega^*(X)})^{-1}}_-{\cong}
        \ar@/_1.2pc/[drr]_-{\kos{p}_!(\xi_X)}
        &\omega^{\pr*}\omega_!\omega^*(X)
        \ar[r]^-{\omega^{\pr*}(\varphi_X)}
        &\omega^{\pr*}(\omega_!(\kappa)\tensor X)
        \ar[d]^-{\omega^{\pr*}_{\omega_!(\kappa),X}}_-{\cong}
        \\
        \text{ }
        &\text{ }
        &\omega^{\pr*}\omega_!(\kappa)\otimes \omega^{\pr*}(X)
      }
    }}
  \end{equation}
  Let us denote $(\phi,\what{\phi})=(\omega^{\pr*}\omega_!,\what{\omega^{\pr*}\omega_!})
  :(\kos{K},\id_{\kos{K}})\to (\kos{K},\id_{\kos{K}})$
  and
  $p=\omega^{\pr*}\omega_!(\kappa)$.
  Then we can check the relation (\ref{eq2 RSKtensoradj HomRep,Hom=Isom})
  as
  \begin{equation*}
    \vcenter{\hbox{
      \xymatrix@C=30pt{
        p\otimes\omega^*(X)
        \ar@/_3pc/@<-3ex>[dddd]_-{\kos{p}_!(\xi_X)}
        \ar[dd]_-{\delta^{p\otimes}_{\omega^{\pr*}(X)}}
        \ar@{=}[rr]
        &\text{ }
        &p\otimes\omega^*(X)
        \ar[d]^-{(\hatar{\phi}_{\omega^*(X)})^{-1}}_-{\cong}
        \\
        \text{ }
        &\text{ }
        &\omega^{\pr*}\omega_!\omega^*(X)
        \ar@/_1pc/[dl]_-{\delta^{\phi}_{\omega^*(X)}}
        \ar[dd]^-{\omega^{\pr*}(\varphi_X)}
        \\
        p\otimes (p\otimes \omega^*(X))
        \ar[dd]_-{I_p\otimes \xi_X}
        \ar@/_0.5pc/[dr]_-{I_p\otimes (\hatar{\phi}_{\omega^*(X)})^{-1}}^-{\cong}
        &\omega^{\pr*}\omega_!\omega^{\pr*}\omega_!\omega^*(X)
        \ar[d]^-{\hatar{\phi}_{\omega^{\pr*}\omega_!\omega^*(X)}}_-{\cong}
        \ar@{}[dr]|-{(\dagger)}
        &\text{ }
        \\
        \text{ }
        &p\otimes \omega^{\pr*}\omega_!\omega^*(X)
        \ar[d]^-{I_p\otimes \omega^{\pr*}(\epsilon_X)}
        &\omega^{\pr*}(\omega_!(\kappa)\tensor X)
        \ar[d]^-{\omega^{\pr*}_{\omega_!(\kappa),X}}_-{\cong}
        \\
        p\otimes \omega^{\pr*}(X)
        \ar@{=}[r]
        &p\otimes \omega^{\pr*}(X)
        \ar@{=}[r]
        &p\otimes \omega^{\pr*}(X)
      }
    }}
  \end{equation*}
  where the diagram $(\dagger)$ is verified below.
  \begin{equation*}
    (\dagger):
    \!\!\!\!
    \!\!\!\!
    \vcenter{\hbox{
      \xymatrix@C=15pt{
        \omega^{\pr*}\omega_!\omega^*(X)
        \ar[ddd]_-{\omega^{\pr*}(\varphi_X)}
        \ar@{=}[r]
        &\omega^{\pr*}\omega_!\omega^*(X)
        \ar[d]^-{\omega^{\pr*}\omega_!(\imath_{\omega^*(X)})}_-{\cong}
        \ar@{=}[r]
        &\omega^{\pr*}\omega_!\omega^*(X)
        \ar[d]^-{\delta^{\omega^{\pr*}\omega_!}_{\omega^*(X)}}
        \\
        \text{ }
        &\omega^{\pr*}\omega_!(\kappa\otimes \omega^*(X))
        \ar[d]^-{\omega^{\pr*}(\omega_{!\kappa,\omega^*(X)})}
        \ar@<0.3ex>@/^1.5pc/[ddr]|-{(\omega^{\pr*}\omega_!)_{\kappa,\omega^*(X)}}
        &\omega^{\pr*}\omega_!\omega^{\pr*}\omega_!\omega^*(X)
        \ar[dd]^-{\hatar{\omega^{\pr*}\omega_!}_{\omega^{\pr*}\omega_!\omega^*(X)}}_-{\cong}
        \\
        \text{ }
        &\omega^{\pr*}(\omega_!(\kappa)\tensor \omega_!\omega^*(X))
        \ar@/^1pc/[dl]^(0.5){\omega^{\pr*}(I_{\omega_!(\kappa)}\tensor \epsilon_X)}
        \ar@/_1pc/[dr]_-{\omega^{\pr*}_{\omega_!(\kappa),\omega_!\omega^*(X)}}^-{\cong}
        &\text{ }
        \\
        \omega^{\pr*}(\omega_!(\kappa)\tensor X)
        \ar[d]_{\omega^{\pr*}_{\omega_!(\kappa),X}}^-{\cong}
        &\text{ }
        &\omega^{\pr*}\omega_!(\kappa)\otimes \omega^{\pr*}\omega_!\omega^*(X)
        \ar[d]^-{I_{\omega^{\pr*}\omega_!(\kappa)}\otimes \omega^{\pr*}(\epsilon_X)}
        \\
        \omega^{\pr*}\omega_!(\kappa)\otimes \omega^{\pr*}(X)
        \ar@{=}[rr]
        &\text{ }
        &\omega^{\pr*}\omega_!(\kappa)\otimes \omega^{\pr*}(X)
      }
    }}
  \end{equation*}
  From the explicit description
  (\ref{eq2 RSKtensoradj HomRep,Hom=Isom}) of
  $\kos{p}_!(\xi_X)$ 
  and the fact that the functor $\kos{p}_!$ is conservative,
  we deduce that
  $\xi_X$ is an isomorphism in $\CK_p$
  if and only if
  $\kos{p}_!(\xi_X)$ is an isomorphism in $\CK$
  if and only if
  $\omega^{\pr*}(\varphi_X)$ is an isomorphism in $\CK$.
  This shows that statements \emph{(a)} and \emph{(b)} are equivalent.

  As statement \emph{(c)} impiles statement \emph{(b)},
  we are left to show that \emph{(a)}$\Leftrightarrow$\emph{(b)} implies \emph{(c)}.
  We begin by verifying that the morphism $\check{\xi}_X:\omega^{\pr*}\omega_!\omega^{\pr*}(X)\to \omega^*(X)$
  in $\CK$
  satisfies the following relations.
  \begin{equation} \label{eq3 RSKtensoradj HomRep,Hom=Isom}
    \vcenter{\hbox{
      \xymatrix@C=20pt{
        \omega^{\pr*}\omega_!\omega^*(X)
        \ar[r]^-{\delta^{\omega^{\pr*}\omega_!}_{\omega^*(X)}}
        \ar@/_1.5pc/[ddr]_-{\varepsilon^{\omega^{\pr*}\omega_!}_{\omega^*(X)}}
        &\omega^{\pr*}\omega_!\omega^{\pr*}\omega_!\omega^*(X)
        \ar[d]^-{\omega^{\pr*}\omega_!\omega^{\pr*}(\epsilon_X)}
        \\
        \text{ }        
        &\omega^{\pr*}\omega_!\omega^{\pr*}(X)
        \ar[d]^-{\check{\xi}_X}
        \\
        \text{ }
        &\omega^*(X)
      }
    }}
    \quad
    \vcenter{\hbox{
      \xymatrix@C=20pt{
        \omega^{\pr*}\omega_!\omega^{\pr*}(X)
        \ar[r]^-{\delta^{\omega^{\pr*}\omega_!}_{\omega^{\pr*}(X)}}
        \ar@/_1.5pc/[ddr]_-{\varepsilon^{\omega^{\pr*}\omega_!}_{\omega^{\pr*}(X)}}
        &\omega^{\pr*}\omega_!\omega^{\pr*}\omega_!\omega^{\pr*}(X)
        \ar[d]^-{\omega^{\pr*}\omega_!(\check{\xi}_X)}
        \\
        \text{ }        
        &\omega^{\pr*}\omega_!\omega^*(X)
        \ar[d]^-{\omega^{\pr*}(\epsilon_X)}
        \\
        \text{ }
        &\omega^{\pr*}(X)
      }
    }}
  \end{equation}
  Again, let us denote $(\phi,\what{\phi})=(\omega^{\pr*}\omega_!,\what{\omega^{\pr*}\omega_!})
  :(\kos{K},\id_{\kos{K}})\to (\kos{K},\id_{\kos{K}})$
  and
  $p=\omega^{\pr*}\omega_!(\kappa)$.
  We verify the left relation in (\ref{eq3 RSKtensoradj HomRep,Hom=Isom}) as follows.
  \begin{equation*}
    \vcenter{\hbox{
      \xymatrix@C=40pt{
        \omega^{\pr*}\omega_!\omega^*(X)
        \ar[ddd]_-{\delta^{\phi}_{\omega^*(X)}}
        \ar@{=}[rr]
        &\text{ }
        &\omega^{\pr*}\omega_!\omega^*(X)
        \ar@/_0.5pc/[dl]^-{\hatar{\phi}_{\omega^*(X)}}_-{\cong}
        \ar[ddddd]^-{\varepsilon^{\phi}_{\omega^*(X)}}
        \\
        \text{ }
        &p\otimes \omega^*(X)
        \ar[d]^-{\delta^{p\otimes}_{\omega^*(X)}}
        \ar@/^2pc/[ddddr]^(0.35){\varepsilon^{p\otimes}_{\omega^*(X)}}
        &\text{ }
        \\
        \text{ }
        &p\otimes (p\otimes \omega^*(X))
        \ar@/^1pc/[dl]_-{(\text{ }\!\hatar{\phi}\hatar{\phi}\text{ }\!)_{\omega^*(X)}^{-1}}^-{\cong}
        \ar[d]^-{I_p\otimes \xi_X}
        &\text{ }
        \\
        \omega^{\pr*}\omega_!\omega^{\pr*}\omega_!\omega^*(X)
        \ar[d]_-{\phi\omega^{\pr*}(\epsilon_X)}
        &p\otimes \omega^{\pr*}(X)
        \ar@/_1pc/[ddr]^-{\xi_X^{-1}}
        \ar@/^1pc/[dl]_-{(\hatar{\phi}_{\omega^{\pr*}(X)})^{-1}}^-{\cong}
        &\text{ }
        \\
        \omega^{\pr*}\omega_!\omega^{\pr*}(X)
        \ar[d]_-{\check{\xi}_X}
        &\text{ }
        &\text{ }
        \\
        \omega^*(X)
        \ar@{=}[rr]
        &\text{ }
        &\omega^*(X)
      }
    }}
  \end{equation*}
  We verify the right relation in (\ref{eq3 RSKtensoradj HomRep,Hom=Isom}) as follows.
  \begin{equation*}
    \vcenter{\hbox{
      \xymatrix@C=40pt{
        \omega^{\pr*}\omega_!\omega^{\pr*}(X)
        \ar[ddd]_-{\delta^{\phi}_{\omega^{\pr*}(X)}}
        \ar@{=}[rr]
        &\text{ }
        &\omega^{\pr*}\omega_!\omega^{\pr*}(X)
        \ar@/_0.5pc/[dl]^-{\hatar{\phi}_{\omega^{\pr*}(X)}}_-{\cong}
        \ar[ddddd]^-{\varepsilon^{\phi}_{\omega^{\pr*}(X)}}
        \\
        \text{ }
        &p\otimes \omega^{\pr*}(X)
        \ar[d]^-{\delta^{p\otimes}_{\omega^{\pr*}(X)}}
        \ar@/^2pc/[ddddr]^(0.35){\varepsilon^{p\otimes}_{\omega^{\pr*}(X)}}
        &\text{ }
        \\
        \text{ }
        &p\otimes (p\otimes \omega^{\pr*}(X))
        \ar@/^1pc/[dl]_-{(\text{ }\!\hatar{\phi}\hatar{\phi}\text{ }\!)_{\omega^{\pr*}(X)}^{-1}}^-{\cong}
        \ar[d]^-{I_p\otimes \xi^{-1}_X}
        &\text{ }
        \\
        \omega^{\pr*}\omega_!\omega^{\pr*}\omega_!\omega^{\pr*}(X)
        \ar[d]_-{\phi(\check{\xi}_X)}
        &p\otimes \omega^*(X)
        \ar@/_1pc/[ddr]^-{\xi_X}
        \ar@/^1pc/[dl]_-{(\hatar{\phi}_{\omega^*(X)})^{-1}}^-{\cong}
        &\text{ }
        \\
        \omega^{\pr*}\omega_!\omega^*(X)
        \ar[d]_-{\omega^{\pr*}(\epsilon_X)}
        &\text{ }
        &\text{ }
        \\
        \omega^{\pr*}(X)
        \ar@{=}[rr]
        &\text{ }
        &\omega^{\pr*}(X)
      }
    }}
  \end{equation*}
  We are ready to show that $\vartheta_X$ and $\vartheta_X^{-1}$ are inverse to each other.
  We can describe both $\vartheta_X$ and $\vartheta_X^{-1}$
  using $\vartheta_!$ as follows.
  \begin{equation*}
    \begin{aligned}
      \vartheta_X
      &:
      \xymatrix@C=30pt{
        c\otimes \omega^*(X)
        \ar[r]^-{(\vartheta_!)_{\omega^*(X)}}
        &\omega^{\pr*}\omega_!\omega^*(X)
        \ar[r]^-{\omega^{\pr*}(\epsilon_X)}
        &\omega^{\pr*}(X)
      }
      \\
      \vartheta^{-1}_X
      &:
      \xymatrix@C=30pt{
        c\otimes \omega^{\pr*}(X)
        \ar[r]^-{(\vartheta_!)_{\omega^{\pr*}(X)}}
        &\omega^{\pr*}\omega_!\omega^{\pr*}(X)
        \ar[r]^-{\check{\xi}_X}
        &\omega^*(X)
      }   
    \end{aligned}
  \end{equation*}
  We have
  \begin{equation*}
    \vcenter{\hbox{
      \xymatrix@C=18pt{
        c\otimes \omega^*(X)
        \ar[d]_-{\delta^{c\otimes}_{\omega^*(X)}}
        \ar@{=}[rr]
        &\text{ }
        &c\otimes \omega^*(X)
        \ar[dd]^-{(\vartheta_!)_{\omega^*(X)}}
        \ar@/^3pc/@<4ex>[ddddd]|-{\varepsilon^{c\otimes}_{\omega^*(X)}}
        \\
        c\otimes (c\otimes \omega^*(X))
        \ar[dd]_-{I_c\otimes \vartheta_X}
        \ar@/^1pc/[dr]^-{I_c\otimes (\vartheta_!)_{\omega^*(X)}}
        &\text{ }
        &\text{ }
        \\
        \text{ }
        &c\otimes \omega^{\pr*}\omega_!\omega^*(X)
        \ar@/_0.5pc/[dl]_-{I_c\otimes \omega^{\pr*}(\epsilon_X)}
        \ar[d]^-{(\vartheta_!)_{\omega^{\pr*}\omega_!\omega^*(X)}}
        &\omega^{\pr*}\omega_!\omega^*(X)
        \ar@/^1pc/[dl]^-{\delta^{\omega^{\pr*}\omega_!}_{\omega^*(X)}}
        \ar[ddd]^-{\varepsilon^{\omega^{\pr*}\omega_!}_{\omega^*(X)}}
        \\
        c\otimes \omega^{\pr*}(X)
        \ar[dd]_-{\vartheta_X^{-1}}
        \ar@/_0.5pc/[dr]_-{(\vartheta_!)_{\omega^{\pr*}(X)}}
        &\omega^{\pr*}\omega_!\omega^{\pr*}\omega_!\omega^*(X)
        \ar[d]^-{\omega^{\pr*}\omega_!\omega^{\pr*}(\epsilon_X)}
        &\text{ }
        \\
        \text{ }
        &\omega^{\pr*}\omega_!\omega^{\pr*}(X)
        \ar[d]^-{\check{\xi}_X}
        &\text{ }
        \\
        \omega^*(X)
        \ar@{=}[r]
        &\omega^*(X)
        \ar@{=}[r]
        &\omega^*(X)
      }
    }}
  \end{equation*}
  as well as
  \begin{equation*}
    \vcenter{\hbox{
      \xymatrix@C=18pt{
        c\otimes \omega^{\pr*}(X)
        \ar[d]_-{\delta^{c\otimes}_{\omega^{\pr*}(X)}}
        \ar@{=}[rr]
        &\text{ }
        &c\otimes \omega^{\pr*}(X)
        \ar[dd]^-{(\vartheta_!)_{\omega^{\pr*}(X)}}
        \ar@/^3pc/@<4ex>[ddddd]|-{\varepsilon^{c\otimes}_{\omega^{\pr*}(X)}}
        \\
        c\otimes (c\otimes \omega^{\pr*}(X))
        \ar[dd]_-{I_c\otimes \vartheta_X^{-1}}
        \ar@/^1pc/[dr]^-{I_c\otimes (\vartheta_!)_{\omega^{\pr*}(X)}}
        &\text{ }
        &\text{ }
        \\
        \text{ }
        &c\otimes \omega^{\pr*}\omega_!\omega^{\pr*}(X)
        \ar@/_0.5pc/[dl]_-{I_c\otimes \check{\xi}_X}
        \ar[d]^-{(\vartheta_!)_{\omega^{\pr*}\omega_!\omega^{\pr*}(X)}}
        &\omega^{\pr*}\omega_!\omega^{\pr*}(X)
        \ar@/^1pc/[dl]^-{\delta^{\omega^{\pr*}\omega_!}_{\omega^{\pr*}(X)}}
        \ar[ddd]^-{\varepsilon^{\omega^{\pr*}\omega_!}_{\omega^{\pr*}(X)}}
        \\
        c\otimes \omega^*(X)
        \ar[dd]_-{\vartheta_X}
        \ar@/_0.5pc/[dr]_-{(\vartheta_!)_{\omega^*(X)}}
        &\omega^{\pr*}\omega_!\omega^{\pr*}\omega_!\omega^{\pr*}(X)
        \ar[d]^-{\omega^{\pr*}\omega_!(\check{\xi}_X)}
        \\
        \text{ }
        &\omega^{\pr*}\omega_!\omega^*(X)
        \ar[d]^-{\omega^{\pr*}(\epsilon_X)}
        \\
        \omega^{\pr*}(X)
        \ar@{=}[r]
        &\omega^{\pr*}(X)
        \ar@{=}[r]
        &\omega^{\pr*}(X)
      }
    }}
  \end{equation*}
  where we used both relations
  in (\ref{eq3 RSKtensoradj HomRep,Hom=Isom}).
  This shows that $\vartheta_X$ and $\vartheta_X^{-1}$
  are inverse to each other as morphisms in the coKleisli category
  $\CK_c$.
  This completes the proof of Proposition~\ref{prop RSKtensoradj HomRep,Hom=Isom}.
\qed\end{proof}

\begin{lemma}
  \label{lem RSKtensoradj antipode}
  Let $(\kos{T},\mon{t})$
  be a strong $\kos{K}$-tensor category
  and let 
  $(\omega,\what{\omega})$,
  $(\omega^{\pr},\what{\omega}^{\pr})
  :(\kos{K},\id_{\kos{K}})\to (\kos{T},\mon{t})$
  be reflective right-strong $\kos{K}$-tensor adjunctions.
  \begin{equation*}
    \vcenter{\hbox{
      \xymatrix{
        (\kos{T},\mon{t})
        \ar@/^1pc/[d]^-{(\omega^*,\what{\omega}^*)}
        \\
        (\kos{K},\id_{\kos{K}})
        \ar@/^1pc/[u]^-{\text{ reflective }(\omega_!,\what{\omega}_!)}
      }
    }}
    \qquad
    \vcenter{\hbox{
      \xymatrix{
        (\kos{T},\mon{t})
        \ar@/^1pc/[d]^-{(\omega^{\pr*},\what{\omega}^{\pr*})}
        \\
        (\kos{K},\id_{\kos{K}})
        \ar@/^1pc/[u]^-{\text{ reflective }(\omega^{\pr}_!,\what{\omega}^{\pr}_!)}
      }
    }}
  \end{equation*}
  Assume that the natural transformations
  \begin{equation*}
    \omega^{\pr*}(\varphi_{\slot}):
    \omega^{\pr*}\omega_!\omega^*
    \Rightarrow
    \omega^{\pr*}(\omega_!(\kappa)\tensor \slot)
    ,
    \qquad
    \omega^*(\varphi^{\pr}_{\slot}):
    \omega^*\omega^{\pr}_!\omega^{\pr*}
    \Rightarrow
    \omega^*(\omega^{\pr}_!(\kappa)\tensor \slot)
  \end{equation*}
  are natural isomorphisms.
  \begin{enumerate}
    \item 
    We have a comonoidal $\kos{K}$-tensor natural isomorphism
    \begin{equation*}
      \varsigma^{\omega,\omega^{\pr}}:
      \xymatrix@C=18pt{
        (\omega^{\pr*}\omega_!,\what{\omega^{\pr*}\omega_!})
        \ar@2{->}[r]^-{\cong}
        &(\omega^*\omega^{\pr}_!,\what{\omega^*\omega^{\pr}_!})
        :(\kos{K},\id_{\kos{K}})\to (\kos{K},\id_{\kos{K}})
      }
    \end{equation*}
    whose component at each object $x$ in $\CK$ is
    \begin{equation*}
      \varsigma^{\omega,\omega^{\pr}}_x:
      \xymatrix@C=30pt{
        \omega^{\pr*}\omega_!(x)
        \ar[r]^-{\omega^{\pr*}\omega_!(\eta^{\pr}_x)}
        &\omega^{\pr*}\omega_!\omega^{\pr*}\omega^{\pr}_!(x)
        \ar[r]^-{\check{\xi}_{\omega^{\pr}_!(x)}}
        &\omega^*\omega^{\pr}_!(x)
        .
      }
    \end{equation*}
    
    \item
    The following diagram of presheaves on $\cat{Ens}(\kos{K})$ strictly commutes.
    \begin{equation*}
      \xymatrix@C=20pt{
        \Hom_{\cat{Ens}(\kos{K})}\bigl(\slot, \omega^{\pr*}\omega_!(\kappa)\bigr)
        \ar@2{->}[r]^-{\cong}
        \ar@2{->}[d]_-{(\varsigma^{\omega,\omega^{\pr}}_{\kappa})\circ \text{ }\slot}^-{\cong}
        &\underline{\Isom}_{\mathbb{SMC}_{\cat{colax}}^{\kos{K}/\!\!/}}
        \bigl(
          (\omega^*,\what{\omega}^*)
          ,
          (\omega^{\pr*},\what{\omega}^{\pr*})
        \bigr)
        \ar@2{->}[d]^-{\text{taking inverse}}_-{\cong}
        &\vartheta
        \ar@{|->}[d]^-{\cong}
        \\
        \Hom_{\cat{Ens}(\kos{K})}\bigl(\slot, \omega^*\omega^{\pr}_!(\kappa)\bigr)
        \ar@2{->}[r]^-{\cong}
        &\underline{\Isom}_{\mathbb{SMC}_{\cat{colax}}^{\kos{K}/\!\!/}}
        \bigl(
          (\omega^{\pr*},\what{\omega}^{\pr*})
          ,
          (\omega^*,\what{\omega}^*)
        \bigr)
        &\vartheta^{-1}
      }
    \end{equation*}
  \end{enumerate}
\end{lemma}
\begin{proof}
  Let us denote 
  $(\phi,\what{\phi})=(\omega^{\pr*}\omega_!,\what{\omega^{\pr*}\omega_!})
  :(\kos{K},\id_{\kos{K}})\to (\kos{K},\id_{\kos{K}})$
  and $p=\omega^{\pr*}\omega_!(\kappa)$.
  By Proposition~\ref{prop RSKtensoradj HomRep,Hom=Isom},
  the universal element
  \begin{equation*}
    \xymatrix{
      \xi:
      (\kos{p}^*\omega^*,\what{\kos{p}^*\omega^*})
      \ar@2{->}[r]^-{\cong}
      &(\kos{p}^*\omega^{\pr*},\what{\kos{p}^*\omega^{\pr*}})
      :(\kos{T},\mon{t})\to (\kos{K}_p,\kos{p}^*)
    }
  \end{equation*}
  is a comonoidal $\kos{K}$-tensor natural isomorphism.
  Therefore
  the following are also comonoidal $\kos{K}$-tensor natural transformations.
  \begin{equation*}
    \begin{aligned}
      \xi
      &:\!
      \xymatrix@C=20pt{
        (p\otimes \omega^*,\what{p\otimes \omega^*})
        \ar@2{->}[r]^-{\kos{p}_!\xi}_-{\cong}
        &(p\otimes \omega^{\pr*},\what{p\otimes \omega^{\pr*}})
        \ar@2{->}[r]^-{\varepsilon^{p\otimes}\omega^*}
        &(\omega^{\pr*},\what{\omega}^{\pr*})
        :(\kos{T},\mon{t})\to (\kos{K},\id_{\kos{K}})
      }
      \\
      \xi^{-1}
      &:\!
      \xymatrix@C=20pt{
        (p\otimes \omega^{\pr*},\what{p\otimes \omega^{\pr*}})
        \ar@2{->}[r]^-{(\kos{p}_!\xi)^{-1}}_-{\cong}
        &(p\otimes \omega^*,\what{p\otimes \omega^*})
        \ar@2{->}[r]^-{\varepsilon^{p\otimes}\omega^*}
        &(\omega^*,\what{\omega}^*)
        :(\kos{T},\mon{t})\to (\kos{K},\id_{\kos{K}})
      }
      \\
      \check{\xi}
      &:\!
      \xymatrix@C=14pt{
        (\omega^{\pr*}\omega_!\omega^{\pr*},\what{\omega^{\pr*}\omega_!\omega^{\pr*}})
        \ar@2{->}[r]^-{\hatar{\phi}\omega^{\pr*}}_-{\cong}
        &(p\otimes \omega^{\pr*},\what{p\otimes \omega^{\pr*}})
        \ar@2{->}[r]^-{\xi^{-1}}
        &(\omega^*,\what{\omega}^*)
        :(\kos{T},\mon{t})\to (\kos{K},\id_{\kos{K}})
      }
    \end{aligned}
  \end{equation*}
  We conclude that
  \begin{equation*}
    \xymatrix@C=50pt{
      \varsigma^{\omega,\omega^{\pr}}:
      (\omega^{\pr*}\omega_!,\what{\omega^{\pr*}\omega_!})
      \ar@2{->}[r]^-{(\check{\xi}\omega^{\pr}_!)\circ (\omega^{\pr*}\omega_!\eta^{\pr})}
      &(\omega^*\omega^{\pr}_!,\what{\omega^*\omega^{\pr}_!})
      :(\kos{K},\id_{\kos{K}})\to (\kos{K},\id_{\kos{K}})
    }
  \end{equation*}
  is a comonoidal $\kos{K}$-tensor natural transformation.
  Next, we verify that the diagram of presheaves commutes.
  It suffices to show that
  $\varsigma^{\omega,\omega^{\pr}}_{\kappa}
  :\omega^{\pr*}\omega_!(\kappa)\to \omega^*\omega^{\pr}_!(\kappa)$
  is the unique morphism 
  in $\cat{Ens}(\kos{K})$ corresponding to 
  the inverse $\xi^{-1}$ of the universal element.
  We can check this as follows.
  \begin{equation*}
    \vcenter{\hbox{
      \xymatrix@C=40pt{
        \omega^{\pr*}\omega_!(\kappa)
        \ar[d]^-{\jmath_{\omega^{\pr*}\omega_!(\kappa)}}_-{\cong}
        \ar@{=}[r]
        &\omega^{\pr*}\omega_!(\kappa)
        \ar[ddd]^-{\omega^{\pr*}\omega_!(\eta^{\pr}_{\kappa})}
        \ar@/^2pc/@<5ex>[dddd]^-{\varsigma^{\omega,\omega^{\pr}}_{\kappa}}
        \\
        \omega^{\pr*}\omega_!(\kappa)\otimes \kappa
        \ar[d]^-{I_{\omega^{\pr*}\omega_!(\kappa)}\otimes \eta^{\pr}_{\kappa}}
        \\
        \omega^{\pr*}\omega_!(\kappa)\otimes \omega^{\pr*}\omega^{\pr}_!(\kappa)
        \ar@/_1pc/[dr]^(0.5){(\hatar{\omega^{\pr*}\omega_!}_{\omega^{\pr*}\omega^{\pr}_!(\kappa)})^{-1}}_-{\cong}
        \ar[dd]^-{\xi^{-1}_{\omega^{\pr}_!(\kappa)}}
        \\
        \text{ }
        &\omega^{\pr*}\omega_!\omega^{\pr*}\omega^{\pr}_!(\kappa)
        \ar[d]^-{\check{\xi}_{\omega^{\pr}_!(\kappa)}}
        \\
        \omega^*\omega^{\pr}_!(\kappa)
        \ar@{=}[r]
        &\omega^*\omega^{\pr}_!(\kappa)
      }
    }}
  \end{equation*}
  This shows that the diagram of presheaves commutes.
  As a consequence,
  we obtain that
  $\varsigma^{\omega,\omega^{\pr}}:
  (\omega^{\pr*}\omega_!,\what{\omega^{\pr*}\omega_!})
  \Rightarrow
  (\omega^*\omega^{\pr}_!,\what{\omega^*\omega^{\pr}_!})$
  is a comonoidal $\kos{K}$-tensor natural isomorphism,
  whose inverse is
  $\varsigma^{\omega^{\pr},\omega}$.
  This completes the proof of Lemma~\ref{lem RSKtensoradj antipode}.
\qed\end{proof}

\begin{lemma}
  \label{lem RSKtensoradj composition}
  Let $(\kos{T},\mon{t})$
  be a strong $\kos{K}$-tensor category.
  Suppose we are given reflective right-strong $\kos{K}$-tensor adjunctions
  $(\omega,\what{\omega})$,
  $(\omega^{\pr},\what{\omega}^{\pr})
  :(\kos{K},\id_{\kos{K}})\to (\kos{T},\mon{t})$
  and a strong $\kos{K}$-tensor functor
  $(\omega^{\ppr*},\what{\omega}^{\ppr*})
  :(\kos{T},\mon{t})\to (\kos{K},\id_{\kos{K}})$.
  \begin{equation*}
    \vcenter{\hbox{
      \xymatrix{
        (\kos{T},\mon{t})
        \ar@/^1pc/[d]^-{(\omega^*,\what{\omega}^*)}
        \\
        (\kos{K},\id_{\kos{K}})
        \ar@/^1pc/[u]^-{\text{ reflective }(\omega_!,\what{\omega}_!)}
      }
    }}
    \qquad
    \vcenter{\hbox{
      \xymatrix{
        (\kos{T},\mon{t})
        \ar@/^1pc/[d]^-{(\omega^{\pr*},\what{\omega}^{\pr*})}
        \\
        (\kos{K},\id_{\kos{K}})
        \ar@/^1pc/[u]^-{\text{ reflective }(\omega^{\pr}_!,\what{\omega}^{\pr}_!)}
      }
    }}
    \qquad
    \vcenter{\hbox{
      \xymatrix{
        (\kos{T},\mon{t})
        \ar[d]^-{(\omega^{\ppr*},\what{\omega}^{\ppr*})}
        \\
        (\kos{K},\id_{\kos{K}})
      }
    }}
  \end{equation*}
  Then the morphism of presheaves
  \begin{equation*}
    \begin{aligned}
      \underline{\Hom}_{\mathbb{SMC}_{\cat{colax}}^{\kos{K}/\!\!/}}
      \bigl(
        (\omega^{\pr*},\what{\omega}^{\pr*})
        ,&
        (\omega^{\ppr*},\what{\omega}^{\ppr*})
      \bigr)
      \text{ }\!\times\text{ }
      \underline{\Hom}_{\mathbb{SMC}_{\cat{colax}}^{\kos{K}/\!\!/}}
      \bigl(
        (\omega^*,\what{\omega}^*)
        ,
        (\omega^{\pr*},\what{\omega}^{\pr*})
      \bigr)
      \\
      &\Rightarrow\text{ }\!
      \underline{\Hom}_{\mathbb{SMC}_{\cat{colax}}^{\kos{K}/\!\!/}}
      \bigl(
        (\omega^*,\what{\omega}^*)
        ,
        (\omega^{\ppr*},\what{\omega}^{\ppr*})
      \bigr)
      :\cat{Ens}(\kos{K})^{\op}\to\cat{Set}
    \end{aligned}
  \end{equation*}
  obtained by composing comonoidal $\kos{K}$-tensor natural transformations
  corresponds to the following morphism 
  between the representing objects in $\cat{Ens}(\kos{K})$.
  \begin{equation*}
    \xymatrix@C=50pt{
      \omega^{\ppr*}\omega^{\pr}_!(\kappa)
      \otimes
      \omega^{\pr*}\omega_!(\kappa)
      \ar[r]^-{(\hatar{\omega^{\ppr*}\omega^{\pr}_!}_{\omega^{\pr*}\omega_!(\kappa)})^{-1}}_-{\cong}
      &\omega^{\ppr*}\omega^{\pr}_!\omega^{\pr*}\omega_!(\kappa)
      \ar[r]^-{\omega^{\ppr*}(\epsilon^{\pr}_{\omega_!(\kappa)})}
      &\omega^{\ppr*}\omega_!(\kappa)
    }
  \end{equation*}
\end{lemma}
\begin{proof}
  Let $c$ be an object in $\cat{Ens}(\kos{K})$.
  Suppose we are given comonoidal $\kos{K}$-tensor natural transformations
  \begin{equation*}
    \vcenter{\hbox{
      \xymatrix{
        (\kos{c}^*\omega^*,\what{\kos{c}^*\omega^*})
        \ar@2{->}[r]^-{\vartheta}
        &(\kos{c}^*\omega^{\pr*},\what{\kos{c}^*\omega^{\pr*}})
        \ar@2{->}[r]^-{\vartheta^{\pr}}
        &(\kos{c}^*\omega^{\ppr*},\what{\kos{c}^*\omega^{\ppr*}})
      }
    }}
  \end{equation*}
  and consider their corresponding mates
  $\vartheta_!:c\otimes\Rightarrow \omega^{\pr*}\omega_!$
  and
  $\vartheta^{\pr}_!:c\otimes\Rightarrow \omega^{\ppr*}\omega^{\pr}_!$.
  Then the corresponding mate of the composition
  $\tilde{\vartheta}=\vartheta^{\pr}\circ\vartheta:
  (\kos{c}^*\omega^*,\what{\kos{c}^*\omega^*})
  \Rightarrow
  (\kos{c}^*\omega^{\ppr*},\what{\kos{c}^*\omega^{\ppr*}})$
  is given by
  \begin{equation*}
    \vcenter{\hbox{
      \xymatrix{
        \tilde{\vartheta}_!:
        c\otimes
        \ar@2{->}[r]^-{\delta^{c\otimes}}
        &c\otimes (c\otimes)
        \ar@2{->}[r]^-{\vartheta^{\pr}_!\vartheta_!}
        &\omega^{\ppr*}\omega^{\pr}_!\omega^{\pr*}\omega_!
        \ar@2{->}[r]^-{\omega^{\ppr*}\epsilon^{\pr}\omega_!}
        &\omega^{\ppr*}\omega_!
      }
    }}
  \end{equation*}
  as we can see from the diagram below.
  Let $x$ be an object in $\CK$.
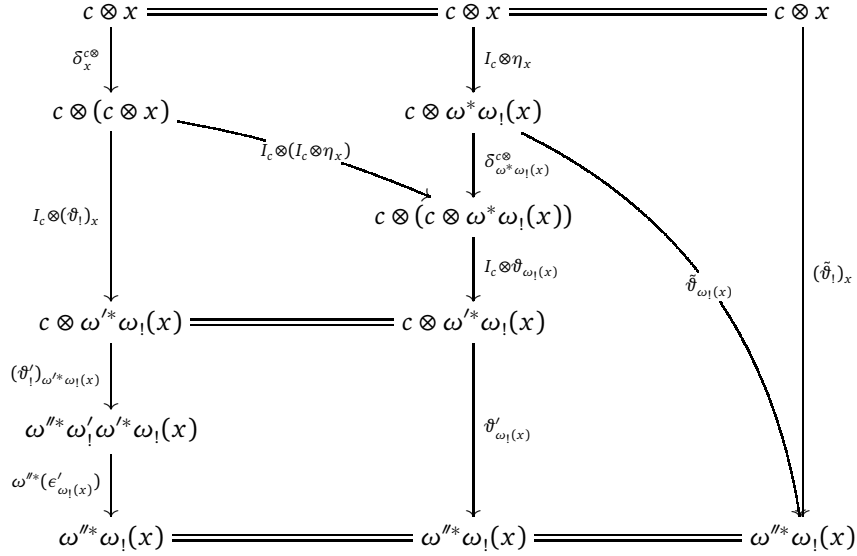
  \begin{equation*}
    \vcenter{\hbox{
      \xymatrix@C=60pt{
        c\otimes x
        \ar[d]_-{\delta^{c\otimes}_x}
        \ar@{=}[r]
        &c\otimes x
        \ar[d]^-{I_c\otimes \eta_x}
        \ar@{=}[r]
        &c\otimes x
        \ar[ddddd]^-{(\tilde{\vartheta}_!)_x}
        \\
        c\otimes (c\otimes x)
        \ar[dd]_-{I_c\otimes (\vartheta_!)_x}
        \ar@/^0.5pc/[dr]|-{I_c\otimes (I_c\otimes \eta_x)}
        &c\otimes \omega^*\omega_!(x)
        \ar[d]^-{\delta^{c\otimes}_{\omega^*\omega_!(x)}}
        \ar@/^2.5pc/[ddddr]|-{\tilde{\vartheta}_{\omega_!(x)}}
        &\text{ }
        \\
        \text{ }
        &c\otimes (c\otimes \omega^*\omega_!(x))
        \ar[d]^-{I_c\otimes \vartheta_{\omega_!(x)}}
        &\text{ }
        \\
        c\otimes \omega^{\pr*}\omega_!(x)
        \ar[d]_-{(\vartheta^{\pr}_!)_{\omega^{\pr*}\omega_!(x)}}
        \ar@{=}[r]
        &c\otimes \omega^{\pr*}\omega_!(x)
        \ar[dd]^-{\vartheta^{\pr}_{\omega_!(x)}}
        &\text{ }
        \\
        \omega^{\ppr*}\omega^{\pr}_!\omega^{\pr*}\omega_!(x)
        \ar[d]_-{\omega^{\ppr*}(\epsilon^{\pr}_{\omega_!(x)})}
        &\text{ }
        &\text{ }
        \\
        \omega^{\ppr*}\omega_!(x)
        \ar@{=}[r]
        &\omega^{\ppr*}\omega_!(x)
        \ar@{=}[r]
        &\omega^{\ppr*}\omega_!(x)
      }
    }}
  \end{equation*}
  Let us denote 
  $c\xrightarrow{g} \omega^{\pr*}\omega_!(\kappa)$,
  $c\xrightarrow{g^{\pr}} \omega^{\ppr*}\omega^{\pr}_!(\kappa)$,
  $c\xrightarrow{\tilde{g}} \omega^{\ppr*}\omega_!(\kappa)$
  as the morphisms in $\cat{Ens}(\kos{K})$
  corresponding to
  $\vartheta$, $\vartheta^{\pr}$, $\tilde{\vartheta}$.
  Then we have the following relation.
  \begin{equation*}
    \vcenter
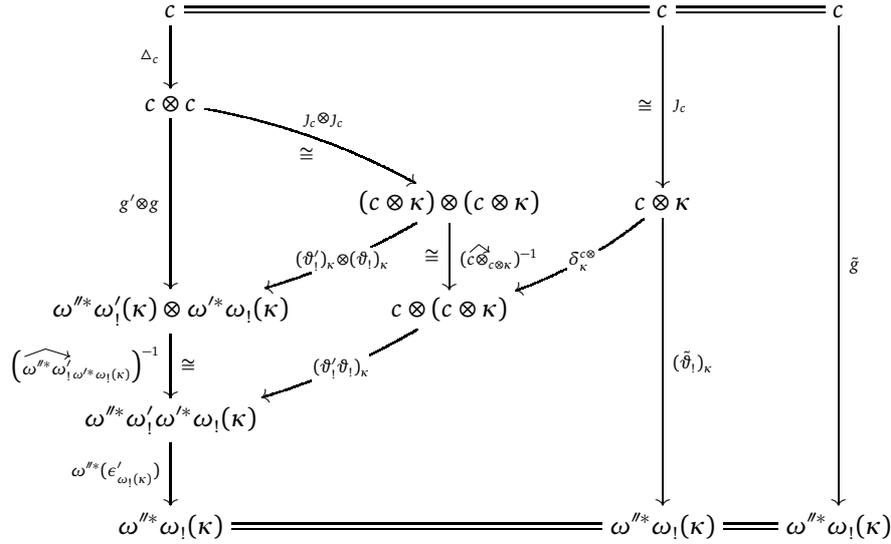
{\hbox{
      \xymatrix@C=20pt{
        c
        \ar[d]_-{\cp_c}
        \ar@{=}[rr]
        &\text{ }
        &c
        \ar[dd]^-{\jmath_c}_-{\cong}
        \ar@{=}[r]
        &c
        \ar[ddddd]^-{\tilde{g}}
        \\
        c\otimes c
        \ar[dd]_-{g^{\pr}\otimes g}
        \ar@/^0.5pc/[dr]^-{\jmath_c\otimes \jmath_c}_-{\cong}
        &\text{ }
        &\text{ }
        &\text{ }
        \\
        \text{ }
        &(c\otimes \kappa)\otimes (c\otimes \kappa)
        \ar@/^0.5pc/[dl]|-{(\vartheta^{\pr}_!)_{\kappa}\otimes (\vartheta_!)_{\kappa}}
        \ar[d]^-{(\hatar{c\otimes}_{c\otimes \kappa})^{-1}}_-{\cong}
        &c\otimes \kappa
        \ar[ddd]^-{(\tilde{\vartheta}_!)_{\kappa}}
        \ar@/^0.5pc/[dl]|-{\delta^{c\otimes}_{\kappa}}
        &\text{ }
        \\
        \omega^{\ppr*}\omega^{\pr}_!(\kappa)\otimes \omega^{\pr*}\omega_!(\kappa)
        \ar[d]_-{\bigl(\hatar{\omega^{\ppr*}\omega^{\pr}_!}_{\omega^{\pr*}\omega_!(\kappa)}\bigr)^{-1}}^-{\cong}
        &c\otimes (c\otimes \kappa)
        \ar@/^0.5pc/[dl]|-{(\vartheta^{\pr}_!\vartheta_!)_{\kappa}}
        &\text{ }
        &\text{ }
        \\
        \omega^{\ppr*}\omega^{\pr}_!\omega^{\pr*}\omega_!(\kappa)
        \ar[d]_-{\omega^{\ppr*}(\epsilon^{\pr}_{\omega_!(\kappa)})}
        &\text{ }
        &\text{ }
        &\text{ }
        \\
        \omega^{\ppr*}\omega_!(\kappa)
        \ar@{=}[rr]
        &\text{ }
        &\omega^{\ppr*}\omega_!(\kappa)
        \ar@{=}[r]
        &\omega^{\ppr*}\omega_!(\kappa)
      }
    }}
  \end{equation*}
  This completes the proof of Lemma~\ref{lem RSKtensoradj composition}.
\qed\end{proof}

\subsection{Colax $\kos{K}$-tensor monad}


\begin{definition}\label{def ColaxKTmonad}
  A \emph{colax $\kos{K}$-tensor monad}
  is a monad internal to the $2$-category $\mathbb{SMC}_{\cat{colax}}^{\kos{K}/\!\!/}$.
\end{definition}

Let us explain Definition~\ref{def ColaxKTmonad} in detail.
Let $(\kos{S},\mon{s})$ be a colax $\kos{K}$-tensor category
where
$\text{$\kos{S}$}=(\CS,\ctimes,\pzc{1})$
is the underlying symmetric monoidal category.
A colax $\kos{K}$-tensor monad on $(\kos{S},\mon{s})$
is a tuple
\begin{equation*}
  \langle \phi,\what{\phi}\rangle
  =(\phi,\what{\phi},\mu,\eta)  
\end{equation*}
of a colax $\kos{K}$-tensor endofunctor
$(\phi,\what{\phi}):(\kos{S},\mon{s})\to (\kos{S},\mon{s})$
and comonoidal $\kos{K}$-tensor natural transformations
$\mu:(\phi\phi,\what{\phi\phi})
\Rightarrow (\phi,\what{\phi})$,
$\eta:(\id_{\kos{S}},I)\Rightarrow (\phi,\what{\phi})$
which satisfy the associativity, unital relations.
Equivalently,
a colax $\kos{K}$-tensor monad on $(\kos{S},\mon{s})$
is a monoid object in the monoidal category
\begin{equation*}
  \kos{E\!n\!\!d}_{\cat{colax}}^{\kos{K}/\!\!/}(\kos{S},\mon{s})
  =
  \Bigl(
    \mathbb{SMC}_{\cat{colax}}^{\kos{K}/\!\!/}\bigl((\kos{S},\mon{s}),(\kos{S},\mon{s})\bigr)
    ,
    \circ
    ,
    (\id_{\kos{S}},I)
  \Bigr)
\end{equation*}
of colax $\kos{K}$-tensor endofunctors on $(\kos{S},\mon{s})$.

\begin{definition}
  \label{def ColaxKTmonad HopfKTmonad}
  Let $(\kos{S},\mon{s})$ be a colax $\kos{K}$-tensor category.
  \begin{enumerate}
    \item 
    The \emph{fusion operator} associated to
    a colax $\kos{K}$-tensor monad
    $\langle\phi,\what{\phi}\rangle$ on $(\kos{S},\mon{s})$
    is the natural transformation
    $\chi_{\slot,\slot}:\phi(\slot\ctimes \phi(\slot))\Rightarrow \phi(\slot)\ctimes\phi(\slot):\CS\times\CS\to \CS$
    whose component at $\pzc{X}$, $\pzc{Y}\in\obj{\CS}$ is
    \begin{equation*}
      \chi_{\pzc{X},\pzc{Y}}:
      \xymatrix@C=35pt{
        \phi(\pzc{X}\ctimes \phi(\pzc{Y}))
        \ar[r]^-{\phi_{\pzc{X}, \phi(\pzc{Y})}}
        &\phi(\pzc{X})\ctimes \phi\phi(\pzc{Y})
        \ar[r]^-{I_{\phi(\pzc{X})}\ctimes \text{ }\mu_{\pzc{Y}}}
        &\phi(\pzc{X})\ctimes \phi(\pzc{Y})
        .
      }
    \end{equation*}

    \item 
    A \emph{Hopf $\kos{K}$-tensor monad} on $(\kos{S},\mon{s})$
    is a colax $\kos{K}$-tensor monad $\langle\phi,\what{\phi}\rangle$ on $(\kos{S},\mon{s})$
    whose associated fusion operator $\chi_{\slot,\slot}$ is a natural isomorphism.
  \end{enumerate}
\end{definition}


Let $\langle \phi,\what{\phi}\rangle=(\phi,\what{\phi},\mu,\eta)$
be a colax $\kos{K}$-tensor monad on
a colax $\kos{K}$-tensor category
$(\kos{S},\mon{s})$.
We have the \emph{Eilenberg-Moore colax $\kos{K}$-tensor category}
\begin{equation}\label{eq ColaxKTmonad defE-Mcat}
  (\text{$\kos{S}$},\mon{s})^{\langle\phi,\what{\phi}\rangle}
  =(\text{$\kos{S}$}^{\langle\phi,\what{\phi}\rangle},\mon{s}^{\langle\phi,\what{\phi}\rangle})
\end{equation}
where we denote
$\text{$\kos{S}$}^{\langle\phi,\what{\phi}\rangle}
=(\CS^{\langle\phi,\what{\phi}\rangle},\ctimes\!^{\langle\phi,\what{\phi}\rangle},\pzc{1}^{\langle\phi,\what{\phi}\rangle})$
as the underlying symmetric monoidal category.
\begin{itemize}
  \item 
  The category $\CS^{\langle\phi,\what{\phi}\rangle}$
  is the Eilenberg-Moore category
  associated to the monad $\langle\phi\rangle=(\phi,\mu,\eta)$ on $\CS$.
  An object in
  $\CS^{\langle\phi,\what{\phi}\rangle}$
  is a pair $X=(\pzc{X},\gamma_{\pzc{X}})$
  of an object $\pzc{X}$ in $\CS$
  and a morphism $\gamma_{\pzc{X}}:\phi(\pzc{X})\to \pzc{X}$ in $\CS$
  which satisfies the following relations.
  \begin{equation*}
    \vcenter{\hbox{
      \xymatrix@C=30pt{
        \phi\phi(\pzc{X})
        \ar[r]^-{\mu_{\pzc{X}}}
        \ar[d]_-{\phi(\gamma_{\pzc{X}})}
        &\phi(\pzc{X})
        \ar[d]^-{\gamma_{\pzc{X}}}
        &\text{ }
        &\pzc{X}
        \ar[r]^-{\eta_{\pzc{X}}}
        \ar@/_1pc/@{=}[dr]
        &\phi(\pzc{X})
        \ar[d]^-{\gamma_{\pzc{X}}}
        \\
        \phi(\pzc{X})
        \ar[r]^-{\gamma_{\pzc{X}}}
        &\pzc{X}
        &\text{ }
        &\text{ }
        &\pzc{X}
      }
    }}
  \end{equation*} 
  Let $Y=(\pzc{Y},\gamma_{\pzc{Y}})$
  be another object in $\CS^{\langle\phi,\what{\phi}\rangle}$.
  A morphism $X\to Y$ in $\CS^{\langle\phi,\what{\phi}\rangle}$
  is a morphism $\pzc{X}\to \pzc{Y}$ in $\CS$ which is compatible with
  $\gamma_{\pzc{X}}$, $\gamma_{\pzc{Y}}$.

  \item
  The monoidal product 
  of objects
  $X$, $Y$
  in $\CS^{\langle\phi,\what{\phi}\rangle}$
  is the pair
  $X\ctimes\!^{\langle\phi,\what{\phi}\rangle} Y
  =(\pzc{X}\ctimes \pzc{Y},\gamma_{\pzc{X}\ctimes \pzc{Y}})$
  where
  $\gamma_{\pzc{X}\ctimes \pzc{Y}}:
  \phi(\pzc{X}\ctimes \pzc{Y})
  \xrightarrow{\phi_{\pzc{X},\pzc{Y}}}
  \phi(\pzc{X})\ctimes \phi(\pzc{Y})
  \xrightarrow{\gamma_{\pzc{X}}\ctimes \gamma_{\pzc{Y}}}
  \pzc{X}\ctimes \pzc{Y}$
  and the unit object 
  in $\CS^{\langle\phi,\what{\phi}\rangle}$
  is the pair
  $\pzc{1}^{\langle\phi,\what{\phi}\rangle}
  =(\pzc{1},\gamma_{\pzc{1}})$ where
  $\gamma_{\pzc{1}}:
  \phi(\pzc{1})\xrightarrow{\phi_{\pzc{1}}}\pzc{1}$.
  The symmetric monoidal coherence morphisms
  $a$, $\imath$, $\jmath$, $s$ 
  of $\kos{S}^{\langle\phi,\what{\phi}\rangle}$
  are given by those of $\kos{S}$.

  \item
  The colax symmetric monoidal functor
  $\text{$\mon{s}$}^{\langle\phi,\what{\phi}\rangle}:
  \kos{K}\to \kos{S}^{\langle\phi,\what{\phi}\rangle}$
  sends each object $x$ in $\CK$ to the pair
  $\mon{s}^{\langle\phi,\what{\phi}\rangle}(x)
  =(\mon{s}(x),\what{\phi}_x)$
  where
  $\what{\phi}_x:\phi\mon{s}(x)\to \mon{s}(x)$.
  The colax symmetric monoidal coherence morphisms
  of $\mon{s}^{\langle\phi,\what{\phi}\rangle}$ are given by those of $\mon{s}$.
  In particular,
  $(\text{$\kos{S}$},\mon{s})^{\langle\phi,\what{\phi}\rangle}
  =(\text{$\kos{S}$}^{\langle\phi,\what{\phi}\rangle},\mon{s}^{\langle\phi,\what{\phi}\rangle})$
  is a strong $\kos{K}$-tensor category
  if and only if
  $(\kos{S},\mon{s})$
  is a strong $\kos{K}$-tensor category.
\end{itemize}
We also have a right-strong $\kos{K}$-tensor adjunction
$(\omega,\what{\omega}):(\kos{S},\mon{s})\to 
(\kos{S},\mon{s})^{\langle\phi,\what{\phi}\rangle}$.
\begin{equation}\label{eq ColaxKTmonad E-McatRSKTadj}
  \vcenter{\hbox{
    \xymatrix{
      (\kos{S},\mon{s})^{\langle\phi,\what{\phi}\rangle}
      \ar@/^1pc/[d]^-{(\omega^*,\what{\omega}^*)}
      \\
      (\kos{S},\mon{s})
      \ar@/^1pc/[u]^-{(\omega_!,\what{\omega}_!)}
    }
  }}
\end{equation}
\begin{itemize}
  \item 
  The right adjoint $\omega^*$ is
  the forgetful functor
  which sends each object
  $X=(\pzc{X},\gamma_{\pzc{X}})$ in $\CS^{\langle\phi,\what{\phi}\rangle}$
  to the underlying object $\omega^*(X)=\pzc{X}$ in $\CS$.
  The symmetric monoidal coherence isomorphisms of
  $\omega^*:\text{$\kos{S}$}^{\langle\phi,\what{\phi}\rangle}\to \text{$\kos{S}$}$
  are identity morphisms,
  and the comonoidal natural isomorphism
  $\what{\omega}^*:\text{$\omega$}^*\mon{s}^{\langle\phi,\what{\phi}\rangle}=\mon{s}$
  is the identity natural transformation.
  
  \item
  The left adjoint $\omega_!$ sends each object
  $\pzc{X}$ in $\CS$
  to the object
  $\omega_!(\pzc{X})=(\phi(\pzc{X}),\mu_{\pzc{X}})$
  in $\CS^{\langle\phi,\what{\phi}\rangle}$
  where $\mu_{\pzc{X}}:\phi\phi(\pzc{X})\to\phi(\pzc{X})$.
  The colax symmetric monoidal coherence morphisms of
  $\omega_!:\kos{S}\to \kos{S}^{\langle\phi,\what{\phi}\rangle}$ are
  given by those of $\phi:\kos{S}\to\kos{S}$,
  and the comonoidal natural transformation
  $\what{\omega}_!:\omega_!\mon{s}\Rightarrow\mon{s}^{\langle\phi,\what{\phi}\rangle}$
  is given by
  $\what{\phi}:\phi\mon{s}\Rightarrow\mon{s}$.

  \item
  The component of the adjunction counit at each object
  $X=(\pzc{X},\gamma_{\pzc{X}})$ in $\CS^{\langle\phi,\what{\phi}\rangle}$
  is $\gamma_{\pzc{X}}:\omega_!\omega^*(X)\to X$.
  The component of the adjunction unit at each object
  $\pzc{X}$ in $\CS$ is
  $\eta_{\pzc{X}}:\pzc{X}\to \phi(\pzc{X})=\omega^*\omega_!(\pzc{X})$.
  The adjunction colax $\kos{K}$-tensor monad on $(\kos{S},\mon{s})$
  induced from the right-strong $\kos{K}$-tensor adjunction
  $(\omega,\what{\omega}):(\kos{S},\mon{s})
  \to (\kos{S},\mon{s})^{\langle\phi,\what{\phi}\rangle}$
  is
  $\langle\omega^*\omega_!,\what{\omega^*\omega_!}\rangle
  =\langle\phi,\what{\phi}\rangle$.
\end{itemize}

Moreover, one can check that the following are true.
\begin{enumerate}
  \item 
  The $\kos{K}$-equivariances
  $\vecar{\phi}$, $\vecar{\omega}_!$
  associated to
  $(\phi,\what{\phi}):(\kos{S},\mon{s})\to (\kos{S},\mon{s})$
  and
  $(\omega_!,\what{\omega}_!):
  (\kos{S},\mon{s})\to
  (\kos{S},\mon{s})^{\langle\phi,\what{\phi}\rangle}$
  satisfy the relation
  $\vecar{\phi}=\omega^*\vecar{\omega}_!$.
  In particular,
  $\vecar{\phi}$ is a natural isomorphism
  if and only if $\vecar{\omega}_!$ is a natural isomorphism.
  
  \item
  The right-strong $\kos{K}$-tensor adjunction
  $(\omega,\what{\omega}):
  (\kos{S},\mon{s})\to
  (\kos{S},\mon{s})^{\langle\phi,\what{\phi}\rangle}$
  satisfies the projection formula
  if and only if
  for every object
  $X=(\pzc{X},\gamma_{\pzc{X}})$ in $\CS^{\langle\phi,\what{\phi}\rangle}$
  and every object $\pzc{Z}$ in $\CS$,
  \begin{equation*}
    \varphi_{\pzc{Z},X}:
    \xymatrix{
      \phi(\pzc{Z}\ctimes \pzc{X})
      \ar[r]^-{\phi_{\pzc{Z},\pzc{X}}}
      &\phi(\pzc{Z})\ctimes \phi(\pzc{X})
      \ar[r]^-{I_{\phi(\pzc{Z})}\ctimes \gamma_{\pzc{X}}}
      &\phi(\pzc{Z})\ctimes \pzc{X}
    }
  \end{equation*}
  is an isomorphism in $\CS$.

  \item
  The given colax $\kos{K}$-tensor monad
  $\langle\phi,\what{\phi}\rangle$ on $(\kos{S},\mon{s})$
  is a Hopf $\kos{K}$-tensor monad
  if and only if 
  for every pair of objects $\pzc{X}$, $\pzc{Z}$ in $\CS$,
  \begin{equation*}
    \chi_{\pzc{Z},\pzc{X}}=
    \varphi_{\pzc{Z},\omega_!(\pzc{X})}:
    \xymatrix@C=30pt{
      \phi(\pzc{Z}\ctimes \phi(\pzc{X}))
      \ar[r]^-{\phi_{\pzc{Z},\phi(\pzc{X})}}
      &\phi(\pzc{Z})\ctimes \phi\phi(\pzc{X})
      \ar[r]^-{I_{\phi(\pzc{Z})}\ctimes \mu_{\pzc{X}}}
      &\phi(\pzc{Z})\ctimes \phi(\pzc{X})
    }
  \end{equation*}
  is an isomorphism in $\CS$.
\end{enumerate}

\subsection{Representations of a group object in $\kos{E\!n\!s}(\kos{K})$}
\label{subsec RepGrpEns(K)}

Recall the cartesian monoidal category
$\kos{E\!n\!s}(\kos{K})=(\cat{Ens}(\kos{K}),\otimes,\kappa)$
of cocommutative comonoids in $\kos{K}$
introduced in (\ref{eq Ens(T) Ens(K)def}).
We denote
\begin{equation*}
  \cat{Mon}(\kos{E\!n\!s}(\kos{K}))
\end{equation*}
as the category of monoid objects in $\kos{E\!n\!s}(\kos{K})$.
A monoid object in $\kos{E\!n\!s}(\kos{K})$
is a triple $(\pi,\pc_{\pi},u_{\pi})$
where $\pi$ is an object in $\cat{Ens}(\kos{K})$
and 
$\pi\otimes \pi\xrightarrow{\pc_{\pi}} \pi$,
$\kappa\xrightarrow{u_{\pi}} \pi$
are product, unit morphisms
in $\cat{Ens}(\kos{K})$
satisfying the associativity, unital relations.
Equivalently,
a monoid object $(\pi,\pc_{\pi},u_{\pi})$
in $\kos{E\!n\!s}(\kos{K})$ is
a cocommutative bimonoid $(\pi,\cp_{\pi},e_{\pi},\pc_{\pi},u_{\pi})$ in $\kos{K}$.
We often denote a monoid object $(\pi,\pc_{\pi},u_{\pi})$
in $\kos{E\!n\!s}(\kos{K})$ as $\pi$.
A morphism $\pi^{\pr}\to \pi$ of monoid objects in $\kos{E\!n\!s}(\kos{K})$
is a morphism in $\cat{Ens}(\kos{K})$
which is compatible with product, unit morphisms of $\pi^{\pr}$, $\pi$.

Let $\pi$ be a monoid object in $\kos{E\!n\!s}(\kos{K})$.
We have the following representable presheaf
on $\cat{Ens}(\kos{K})$
which factors through the category $\cat{Mon}$ of monoids.
\begin{equation*}
  \Hom_{\cat{Ens}(\kos{K})}(\slot,\pi):
  \cat{Ens}(\kos{K})^{\op}\to\cat{Mon}
\end{equation*}
For each object $c$ in $\cat{Ens}(\kos{K})$,
the set
$\Hom_{\cat{Ens}(\kos{K})}(c,\pi)$
is a monoid equipped with convolution product
$(c\xrightarrow{f}\pi,c\xrightarrow{g}\pi)
\mapsto
f\star g:c\xrightarrow{\cp_c}c\otimes c\xrightarrow{f\otimes g}\pi\otimes \pi\xrightarrow{\pc_{\pi}}\pi$
and identity object
$c\xrightarrow{e_c}\kappa\xrightarrow{u_{\pi}}\pi$.
Using the Yoneda lemma,
one can check that the correspondence
$\pi\mapsto \Hom_{\cat{Ens}(\kos{K})}(\slot,\pi)$
defines an equivalence of categories
from
$\cat{Mon}(\kos{E\!n\!s}(\kos{K}))$
to the category of
representable presheaves
on $\cat{Ens}(\kos{K})$
which factors through $\cat{Mon}$.

\begin{remark} \label{rem RepGrpEns(K) pitensor conservative}
  Let $\pi$ be a monoid object in $\kos{E\!n\!s}(\kos{K})$.
  Then the functor $\pi\otimes\slot:\CK\to \CK$ is conservative.
  Indeed, if $x\xrightarrow{f} y$ is a morphism in $\CK$
  such that $\pi\otimes x\xrightarrow{I_{\pi}\otimes f} \pi\otimes y$ is an isomorphism,
  then $f$ has an inverse
  $f^{-1}:
  y
  \xrightarrow[\cong]{\imath_y}
  \kappa\otimes y
  \xrightarrow{u_{\pi}\otimes I_y}
  \pi\otimes y
  \xrightarrow[\cong]{(I_{\pi}\otimes f)^{-1}}
  \pi\otimes x
  \xrightarrow{e_{\pi}\otimes I_x}
  \kappa\otimes x
  \xrightarrow[\cong]{\imath_x^{-1}}
  x$.
\end{remark}

In (\ref{eq Ens(T) colaxKTend(K,idK) def})
we introduced the monoidal category
$\kos{E\!n\!\!d}_{\cat{colax}}^{\kos{K}/\!\!/}(\kos{K},\id_{\kos{K}})_{\cat{rfl}}$
of reflective colax $\kos{K}$-tensor endofunctors on $(\kos{K},\id_{\kos{K}})$.
Recall the definition of colax $\kos{K}$-tensor monads
in Definition~\ref{def ColaxKTmonad}.

\begin{definition}
  A colax $\kos{K}$-tensor monad
  $\langle\phi,\what{\phi}\rangle=(\phi,\what{\phi},\mu,\eta)$ on
  $(\kos{K},\id_{\kos{K}})$
  is called \emph{reflective}
  if the underlying colax $\kos{K}$-tensor endofunctor
  $(\phi,\what{\phi})$ on $(\kos{K},\id_{\kos{K}})$ is reflective.
  We denote
  \begin{equation*}
    \cat{Monad}_{\cat{colax}}^{\kos{K}/\!\!/}(\kos{K},\id_{\kos{K}})_{\cat{rfl}}
  \end{equation*}
  as the category of reflective colax $\kos{K}$-tensor monads on $(\kos{K},\id_{\kos{K}})$,
  which is the category of monoid objects in
  $\kos{E\!n\!\!d}_{\cat{colax}}^{\kos{K}/\!\!/}(\kos{K},\id_{\kos{K}})_{\cat{rfl}}$.
\end{definition}

Recall the adjoint equivalence of monoidal categories
\begin{equation*}
  \CL:
  \kos{E\!n\!\!d}_{\cat{colax}}^{\kos{K}/\!\!/}(\kos{K},\id_{\kos{K}})_{\cat{rfl}}
  \simeq
  \kos{E\!n\!s}(\kos{K})
  :\iota
\end{equation*}
in Proposition~\ref{prop Ens(T) KtensorEnd(K)=Ens(K)}.
By considering the category of monoid objects on both sides,
we obtain the following corollary.

\begin{corollary}
  \label{cor RepGrpEns(K) rflCKTmonad=Mon(Ens(K))}
  We have an adjoint equivalence of categories
  \begin{equation*}
    \CL:
    \cat{Monad}_{\cat{colax}}^{\kos{K}/\!\!/}(\kos{K},\id_{\kos{K}})_{\cat{rfl}}
    \simeq
    \cat{Mon}(\kos{E\!n\!s}(\kos{K}))
    :\iota
    .
  \end{equation*}
  \begin{itemize}
    \item 
    The right adjoint
    $\iota$
    sends each monoid object $\pi$ in $\kos{E\!n\!s}(\kos{K})$
    to the reflective colax $\kos{K}$-tensor monad
    $\iota(\pi)
    =\langle\pi\otimes,\what{\pi\otimes}\rangle
    =(\pi\otimes,\what{\pi\otimes},\mu^{\pi\otimes},\eta^{\pi\otimes})$
    on $(\kos{K},\id_{\kos{K}})$.
    The components of $\mu^{\pi\otimes}$, $\eta^{\pi\otimes}$
    at each object $x$ in $\CK$ are described below.
    \begin{equation*}
      \hspace*{-0.5cm}
      \mu^{\pi\otimes}_x
      :\!\!
      \xymatrix@C=18pt{
        \pi\otimes (\pi\otimes x)
        \ar[r]^-{a_{\pi,\pi,x}}_-{\cong}
        &(\pi\otimes \pi)\otimes x
        \ar[r]^-{\pc_{\pi}\otimes I_x}
        &\pi\otimes x
      }
      \quad
      \eta^{\pi\otimes}_x
      :\!\!
      \xymatrix@C=18pt{
        x
        \ar[r]^-{\imath_x}_-{\cong}
        &\kappa\otimes x
        \ar[r]^-{u_{\pi}\otimes I_x}
        &\pi\otimes x
      }
    \end{equation*}
   
    \item 
    The left adjoint
    sends each reflective colax $\kos{K}$-tensor monad
    $\!\langle\phi,\what{\phi}\rangle\!=\!(\phi,\what{\phi},\mu,\eta)$
    on $(\kos{K},\id_{\kos{K}})$
    to the monoid object
    $\CL\langle \phi,\what{\phi}\rangle
    =\phi(\kappa)$ in $\kos{E\!n\!s}(\kos{K})$
    whose product, unit morphisms are
    \begin{equation}\label{eq RepGrpEns(K) rflCKTmonad=Mon(Ens(K))}
      \pc_{\phi(\kappa)}
      :
      \xymatrix{
        \phi(\kappa)\otimes \phi(\kappa)
        \ar[r]^-{(\hatar{\phi}_{\phi(\kappa)})^{-1}}_-{\cong}
        &\phi\phi(\kappa)
        \ar[r]^-{\mu_{\kappa}}
        &\phi(\kappa)
        ,
      }
      \qquad
      u_{\phi(\kappa)}
      :
      \xymatrix{
        \kappa
        \ar[r]^-{\eta_{\kappa}}
        &\phi(\kappa)
        .
      }
    \end{equation}
  \end{itemize}
\end{corollary}


A group object in $\kos{E\!n\!s}(\kos{K})$
is a tuple $(\pi,\pc_{\pi},u_{\pi},\varsigma_{\pi})$
where $(\pi,\pc_{\pi},u_{\pi})$ is a monoid object in $\kos{E\!n\!s}(\kos{K})$
and
$\varsigma_{\pi}:\pi\to \pi$
is the antipode morphism,
which is a morphism in $\cat{Ens}(\kos{K})$
that satisfies the relation
$\varsigma_{\pi}\star I_{\pi}=I_{\pi}\star \varsigma_{\pi}=u_{\pi}\circ e_{\pi}:\pi\to \pi$.
Equivalently,
a group object $(\pi,\pc_{\pi},u_{\pi},\varsigma_{\pi})$
in $\kos{E\!n\!s}(\kos{K})$ is
a cocommutative Hopf monoid $(\pi,\cp_{\pi},e_{\pi},\pc_{\pi},u_{\pi},\varsigma_{\pi})$ in $\kos{K}$.
One can check that the antipode morphism
$\varsigma_{\pi}$
is involutive (hence is an isomorphism),
and is an anti-morphism of monoid objects in $\kos{E\!n\!s}(\kos{K})$.
\begin{equation*}
  \vcenter{\hbox{
    \xymatrix@C=30pt{
      \pi
      \ar[r]^-{\varsigma_{\pi}}_-{\cong}
      \ar@/_1pc/@{=}[dr]
      &\pi
      \ar[d]^-{\varsigma_{\pi}}_-{\cong}
      &\pi\otimes \pi
      \ar[d]_-{\pc_{\pi}}
      \ar[r]^-{s_{\pi,\pi}}_-{\cong}
      &\pi\otimes \pi
      \ar[r]^-{\varsigma_{\pi}\otimes \varsigma_{\pi}}_-{\cong}
      &\pi\otimes \pi
      \ar[d]^-{\pc_{\pi}}
      &\kappa
      \ar[d]_-{u_{\pi}}
      \ar@/^1pc/[dr]^-{u_{\pi}}
      &\text{ }
      \\
      \text{ }
      &\pi
      &\pi
      \ar[rr]^-{\varsigma_{\pi}}_-{\cong}
      &\text{ }
      &\pi
      &\pi
      \ar[r]^-{\varsigma_{\pi}}_-{\cong}
      &\pi
    }
  }}
\end{equation*}
We also denote a group object $(\pi,\pc_{\pi},u_{\pi},\varsigma_{\pi})$
in $\kos{E\!n\!s}(\kos{K})$ as $\pi$.
A morphism $\pi^{\pr}\to \pi$
of group objects in $\kos{E\!n\!s}(\kos{K})$
is a morphism as monoid objects in $\kos{E\!n\!s}(\kos{K})$,
i.e., it is a morphism $\pi^{\pr}\to \pi$ in $\cat{Ens}(\kos{K})$
which is compatible with product, unit morphisms of $\pi^{\pr}$, $\pi$.

Let $\pi$ be a group object in $\kos{E\!n\!s}(\kos{K})$.
We have the following representable presheaf
on $\cat{Ens}(\kos{K})$
which factors through the category $\cat{Grp}$ of groups.
\begin{equation} \label{eq RepGrpEns(K) Hom(-,pi)}
  \Hom_{\cat{Ens}(\kos{K})}(\slot,\pi):
  \cat{Ens}(\kos{K})^{\op}\to\cat{Grp}
\end{equation}
For each object $c$ in $\cat{Ens}(\kos{K})$,
the set
$\Hom_{\cat{Ens}(\kos{K})}(c,\pi)$
is a group
where the inverse of 
$c\xrightarrow{f}\pi$
with respect to the convolution product $\star$
is $c\xrightarrow{f}\pi\xrightarrow[\cong]{\varsigma_{\pi}}\pi$.
Using the Yoneda lemma,
one can check that the correspondence
$\pi\mapsto \Hom_{\cat{Ens}(\kos{K})}(\slot,\pi)$
defines an equivalence of categories
from
the category of group objects in $\kos{E\!n\!s}(\kos{K})$
to the category of
representable presheaves
on $\cat{Ens}(\kos{K})$
which factors through $\cat{Grp}$.

Recall the definition of Hopf $\kos{K}$-tensor monads
in Definition~\ref{def ColaxKTmonad HopfKTmonad}.

\begin{definition}
  A Hopf $\kos{K}$-tensor monad
  $\langle\phi,\what{\phi}\rangle=(\phi,\what{\phi},\mu,\eta)$ on
  $(\kos{K},\id_{\kos{K}})$
  is called \emph{reflective}
  if the underlying colax $\kos{K}$-tensor endofunctor
  $(\phi,\what{\phi})$ on $(\kos{K},\id_{\kos{K}})$ is reflective.
\end{definition}

Each group object $\pi$ in $\kos{E\!n\!s}(\kos{K})$
determines a reflective Hopf $\kos{K}$-tensor monad
$\iota(\pi)=\langle\pi\otimes,\what{\pi\otimes}\rangle$
on $(\kos{K},\id_{\kos{K}})$.
The component of the associated fusion operator
at $x$, $y\in\obj{\CK}$ is
\begin{equation*}
  \chi_{x,y}:
  \xymatrix@C=30pt{
    \pi\otimes x\otimes \pi\otimes y
    \ar[r]^-{(\pi\otimes)_{x,\pi\otimes y}}
    &\pi\otimes x\otimes \pi\otimes \pi\otimes y
    \ar[r]^-{I_{\pi\otimes x}\otimes \pc_{\pi}\otimes I_y}
    &\pi\otimes x\otimes \pi\otimes y
  }
\end{equation*}
which is an isomorphism in $\CK$
whose inverse is described below.
\begin{equation*}
  \xymatrix@C=40pt{
    \pi\otimes x\otimes \pi\otimes y
    \ar[d]_-{\chi_{x,y}^{-1}}
    \ar[r]^-{(\pi\otimes)_{x,\pi\otimes y}}
    &\pi\otimes x\otimes \pi\otimes \pi\otimes y
    \ar[d]^-{I_{\pi\otimes x}\otimes \varsigma_{\pi}\otimes I_{\pi\otimes y}}_-{\cong}
    \\
    \pi\otimes x\otimes \pi\otimes y
    &\pi\otimes x\otimes \pi\otimes \pi\otimes y
    \ar[l]_-{I_{\pi\otimes x}\otimes \pc_{\pi}\otimes I_y}
  }
\end{equation*}

Let $\pi$ be a group object in $\kos{E\!n\!s}(\kos{K})$.
We define the strong $\kos{K}$-tensor category
of representations of $\pi$ in $\kos{K}$ as follows.
Recall the Eilenberg-Moore category 
associated to a colax $\kos{K}$-tensor monad
introduced in (\ref{eq ColaxKTmonad defE-Mcat}).
  
\begin{definition}\label{def RepGrpEns(K) Rep(pi)Ktensorcat}
  Let $\pi$ be a group object in $\kos{E\!n\!s}(\kos{K})$.
  We define the strong $\kos{K}$-tensor category
  \begin{equation*}
    (\kos{R\!e\!p}(\pi),\kos{t}^*_{\pi})    
  \end{equation*}
  of representations of $\pi$ in $\kos{K}$
  as the Eilenberg-Moore strong $\kos{K}$-tensor category
  associated to the reflective Hopf $\kos{K}$-tensor monad
  $\langle \pi\otimes,\what{\pi\otimes}\rangle$ on $(\kos{K},\id_{\kos{K}})$.
  The underlying symmetric monoidal category is denoted as
  \begin{equation*}
    \kos{R\!e\!p}(\pi)=(\cat{Rep}(\pi),\tensor\!_{\pi},\unit\!_{\pi})
  \end{equation*}
  and the strong symmetric monoidal functor
  $\kos{t}_{\pi}^*:\kos{K}\to \kos{R\!e\!p}(\pi)$
  is the functor of constructing trivial representations.
\end{definition}

Let us explain Definition~\ref{def RepGrpEns(K) Rep(pi)Ktensorcat} in detail.
Let $\pi$ be a group object in $\kos{E\!n\!s}(\kos{K})$.
An object in
$\cat{Rep}(\pi)$
is a pair $X=(x,\gamma_x)$
of an object $x$ in $\CK$
and a morphism $\gamma_x:\pi\otimes x\to x$ in $\CK$
satisfying the left $\pi$-action relations.
Let $Y=(y,\gamma_y)$ be another object in $\cat{Rep}(\pi)$.
A morphism $X\to Y$ in $\cat{Rep}(\pi)$
is a morphism $x\to y$ in $\CK$ which is compatible with
$\gamma_x$, $\gamma_y$.
The monoidal product 
of objects $X$, $Y$ in $\cat{Rep}(\pi)$
is the pair
$X\tensor\!_{\pi}Y=(x\otimes y,\gamma_{x\otimes y})$
where 
$\gamma_{x\otimes y}:
\pi\otimes (x\otimes y)
\xrightarrow{(\pi\otimes)_{x,y}}
(\pi\otimes x)\otimes (\pi\otimes y)
\xrightarrow{\gamma_x\otimes \gamma_y}
x\otimes y$
and the unit object in $\cat{Rep}(\pi)$
is the pair
$\unit\!_{\pi}=(\kappa,\gamma_{\kappa})$ where
$\gamma_{\kappa}=\varepsilon^{\pi\otimes}_{\kappa}:\pi\otimes \kappa\xrightarrow{e_{\pi}\otimes I_{\kappa}}\kappa\otimes\kappa\xrightarrow[\cong]{\cp_{\kappa}^{-1}}\kappa$.
The symmetric monoidal coherence morphisms
$a$, $\imath$, $\jmath$, $s$ 
of $\kos{R\!e\!p}(\pi)$ are given by those of $\kos{K}$.
The functor
$\kos{t}_{\pi}^*$
sends each object $z$ in $\CK$ to the pair
\begin{equation*}
  \kos{t}_{\pi}^*(z)
  =
  (z,\what{\pi\otimes}_z:\pi\otimes z\to z),
  \qquad
  \what{\pi\otimes}_z
  =
  \varepsilon^{\pi\otimes}_z:
  \xymatrix{
    \pi\otimes z
    \ar[r]^-{e_{\pi}\otimes I_z}
    &\kappa\otimes z
    \ar[r]^-{\imath_z^{-1}}_-{\cong}
    &z
    .
  }
\end{equation*}
The symmetric monoidal coherence isomorphisms of
$\kos{t}_{\pi}^*:\kos{K}\to\kos{R\!e\!p}(\pi)$
are given by identity morphisms.

\begin{remark}
  Let $\pi$ be a group object in $\kos{E\!n\!s}(\kos{K})$.
  We explain why we define representations of $\pi$ in $\kos{K}$
  as in Definition~\ref{def RepGrpEns(K) Rep(pi)Ktensorcat}.
  For each object $x$ in $\CK$,
  we denote 
  \begin{equation*}
    \underline{\Aut}(x):\cat{Ens}(\kos{K})^{\op}\to\cat{Grp}
  \end{equation*}
  as the presheaf of groups of automorphisms of $x$.
  It sends each object $c$ in $\cat{Ens}(\kos{K})$
  to the group $\Aut_{\CK_c}(\kos{c}^*(x))$
  of automorphisms of $\kos{c}^*(x)$ in the coKleisli category $\CK_c$,
  and 
  each morphism $c^{\pr}\xrightarrow{f} c$ in $\cat{Ens}(\kos{K})$
  to the morphism of groups
  $\Aut_{\CK_c}(\kos{c}^*(x))
  \to
  \Aut_{\CK_{c^{\pr}}}(\kos{c}^{\pr*}(x))$
  induced by the functor
  $\kos{f}^*:\kos{K}_c\to\kos{K}_{c^{\pr}}$
  introduced in (\ref{eq RSKtensoradj functorbetweencoKelisliCat}).
  \begin{itemize}
    \item 
    Let $X=(x,\gamma_x:\pi\otimes x\to x)$
    be a representation of $\pi$ in $\kos{K}$.
    Then we have a morphism of presheaves of groups
    \begin{equation}\label{eq RepGrpEns(K) mor of presheaves}
      r^x:
      \Hom_{\cat{Ens}(\kos{K})}(\slot,\pi)
      \Rightarrow
      \underline{\Aut}(x):
      \cat{Ens}(\kos{K})^{\op}\to\cat{Grp}
    \end{equation}
    which sends each morphism $c\xrightarrow{g}\pi$ in $\cat{Ens}(\kos{K})$
    to the automorphism $r^x_c(g):\kos{c}^*(x)\xrightarrow{\cong}\kos{c}^*(x)$
    in $\CK_c$ where
    $r^x_c(g):c\otimes x\xrightarrow{g\otimes I_x}\pi\otimes x\xrightarrow{\gamma_x}x$.

    \item
    Conversely, suppose we are given a pair
    $(x,r^x)$ of an object $x$ in $\CK$
    and a morphism of presheaves $r^x$ as in (\ref{eq RepGrpEns(K) mor of presheaves}).
    Then we have a representation $X=(x,\gamma_x)$ of $\pi$ in $\kos{K}$
    where $\gamma_x=r^x_{\pi}(I_{\pi}):\pi\otimes x\to x$.
  \end{itemize}
  Using the Yoneda Lemma,
  one can show that the correspondences
  $X=(x,\gamma_x)\mapsto (x,r^x)$
  and 
  $(x,r^x)\mapsto X=(x,\gamma_x)$
  described above are inverse to each other.
\end{remark}

\begin{lemma} \label{lem RepGrpEns(K) Rep(pi)}
  Let $\pi$ be a group object in $\kos{E\!n\!s}(\kos{K})$
  and recall the strong $\kos{K}$-tensor category
  $(\kos{R\!e\!p}(\pi),\kos{t}^*_{\pi})$
  of representations of $\pi$ in $\kos{K}$
  defined in Definition~\ref{def RepGrpEns(K) Rep(pi)Ktensorcat}.
  \begin{enumerate}
    \item 
    The functor $\kos{t}^*_{\pi}$ of constructing trivial representations
    is fully faithful.
    
    \item 
    We have a right-strong $\kos{K}$-tensor adjunction
    \begin{equation*}
      \vcenter{\hbox{
        \xymatrix{
          (\kos{R\!e\!p}(\pi),\kos{t}_{\pi}^*)
          \ar@/^1pc/[d]^-{(\omega_{\pi}^*,\what{\omega}_{\pi}^*)}
          \\
          (\kos{K},\id_{\kos{K}})
          \ar@/^1pc/[u]^-{(\omega_{\pi!},\what{\omega}_{\pi!})}
        }
      }}
      \qquad
      (\omega_{\pi},\what{\omega}_{\pi})
      :(\kos{K},\id_{\kos{K}})\to 
      (\kos{R\!e\!p}(\pi),\kos{t}_{\pi}^*)
    \end{equation*}
    satisfying the projection formula,
    such that the left adjoint
    $(\omega_{\pi!},\what{\omega}_{\pi!})$
    is a reflective colax $\kos{K}$-tensor functor
    and the underlying functor $\omega_{\pi!}$ is conservative.
  \end{enumerate}
\end{lemma}
\begin{proof}
  One can easily check that the functor
  $\kos{t}^*_{\pi}:\CK\to \cat{Rep}(\pi)$
  of constructing trivial representations
  is fully faithful.
  The right-strong $\kos{K}$-tensor adjunction
  $(\omega_{\pi},\what{\omega}_{\pi})$
  is defined as in (\ref{eq ColaxKTmonad E-McatRSKTadj}).
  \begin{itemize}
    \item 
    The right adjoint $\omega^*_{\pi}$ is
    the forgetful functor
    which sends each object
    $X=(x,\gamma_x)$ in $\cat{Rep}(\pi)$
    to the underlying object $\omega^*_{\pi}(X)=x$ in $\CK$.
    The symmetric monoidal coherence isomorphisms of
    $\omega^*_{\pi}:\kos{R\!e\!p}(\pi)\to \kos{K}$
    are identity morphisms,
    and the comonoidal natural isomorphism
    $\what{\omega}_{\pi}^*:\omega_{\pi}^*\kos{t}_{\pi}^*=\id_{\kos{K}}$
    is the identity natural transformation.
    
    \item
    The left adjoint $\omega_{\pi!}$ sends each object
    $x$ in $\CK$
    to the object
    $\omega_{\pi!}(x)=(\pi\otimes x,\mu^{\pi\otimes}_x)$
    in $\cat{Rep}(\pi)$
    where
    $\mu^{\pi\otimes}_x:\pi\otimes (\pi\otimes x)\xrightarrow[\cong]{a_{\pi,\pi,x}}(\pi\otimes \pi)\otimes x\xrightarrow{\pc_{\pi}\otimes I_x}\pi\otimes x$.
    The colax symmetric monoidal coherence morphisms of
    $\omega_{\pi!}:\kos{K}\to \kos{R\!e\!p}(\pi)$ are
    given by those of $\pi\otimes:\kos{K}\to\kos{K}$,
    and the comonoidal natural transformation
    $\what{\omega}_{\pi!}:\omega_{\pi!}\Rightarrow\kos{t}_{\pi}^*$
    is given by
    $\what{\pi\otimes}:\pi\otimes\Rightarrow\id_{\kos{K}}$.
  
    \item
    The component of the adjunction counit at each object
    $X=(x,\gamma_x)$ in $\cat{Rep}(\pi)$
    is $\gamma_x:\omega_{\pi!}\omega_{\pi}^*(X)\to X$.
    The component of the adjunction unit at each object
    $x$ in $\CK$ is
    $\eta^{\pi\otimes}_x:x\xrightarrow[\cong]{\imath_x}\kappa\otimes x\xrightarrow{u_{\pi}\otimes I_x}\pi\otimes x=\omega_{\pi}^*\omega_{\pi!}(x)$.
  \end{itemize}
  The right-strong $\kos{K}$-tensor adjunction
  $(\omega_{\pi},\what{\omega}_{\pi})$
  satisfies the projection formula.
  For each object
  $X=(x,\gamma_x)$ in $\cat{Rep}(\pi)$
  and each object $z$ in $\CK$,
  \begin{equation*}
    \varphi_{z,X}:
    \xymatrix@C=30pt{
      \pi\otimes z\otimes x
      \ar[r]^-{(\pi\otimes)_{z,x}}
      &\pi\otimes z\otimes \pi\otimes x
      \ar[r]^-{I_{\pi\otimes z}\otimes \gamma_x}
      &\pi\otimes z\otimes x
    }
  \end{equation*}
  is an isomorphism in $\CK$,
  whose inverse is
  \begin{equation*}
    \xymatrix@C=35pt{
      \varphi_{z,X}^{-1}:
      \pi\otimes z\otimes x
      \ar[r]^-{(\pi\otimes)_{z,x}}
      &\pi\otimes z\otimes \pi\otimes x
      \ar[r]^-{I_{\pi\otimes z}\otimes \varsigma_{\pi}\otimes I_x}_-{\cong}
      &\pi\otimes z\otimes \pi\otimes x
      \ar[r]^-{I_{\pi\otimes z}\otimes \gamma_x}
      &\pi\otimes z\otimes x
      .
    }
  \end{equation*}
  The colax $\kos{K}$-tensor endofunctor
  $(\pi\otimes,\what{\pi\otimes})$ on $(\kos{K},\id_{\kos{K}})$ is reflective
  and the underlying functor
  $\pi\otimes:\CK\to \CK$ 
  is conservative as explained in Remark~\ref{rem RepGrpEns(K) pitensor conservative}.
  Since we have
  $(\omega^*_{\pi}\omega_{\pi!},\what{\omega^*_{\pi}\omega_{\pi!}})
  =(\pi\otimes,\what{\pi\otimes})$
  as colax $\kos{K}$-tensor endofunctor on $(\kos{K},\id_{\kos{K}})$,
  we obtain that the left adjoint
  $(\omega_{\pi!},\what{\omega}_{\pi!})
  :(\kos{K},\id_{\kos{K}})\to 
  (\kos{R\!e\!p}(\pi),\kos{t}_{\pi}^*)$
  is also reflective,
  and the underlying functor
  $\omega_{\pi!}:\CK\to \cat{Rep}(\pi)$
  is also conservative.
  This completes the proof of Lemma~\ref{lem RepGrpEns(K) Rep(pi)}.
\qed\end{proof}

We define as follows,
which is analogous to Definition~\ref{def RSKtensoradj Hom,Isom}.

\begin{definition} \label{def RepGrpEns(K) End,Aut}
  Let $(\kos{T},\mon{t})$ be a colax $\kos{K}$-tensor category
  and let $(\phi,\what{\phi}):(\kos{T},\mon{t})\to (\kos{K},\id_{\kos{K}})$
  be a colax $\kos{K}$-tensor functor.
  \begin{enumerate}
    \item 
    We define the \emph{presheaf of monoids of comonoidal $\kos{K}$-tensor natural endomorphisms
    of $(\phi,\what{\phi})$}
    as the presheaf of monoids on $\cat{Ens}(\kos{K})$
    \begin{equation*}
      \underline{\End}_{\mathbb{SMC}_{\cat{colax}}^{\kos{K}/\!\!/}}
      (\phi,\what{\phi})
      :
      \cat{Ens}(\kos{K})^{\op}
      \to
      \cat{Mon}
    \end{equation*}
    which sends each object $c$ in $\cat{Ens}(\kos{K})$
    to the monoid of comonoidal $\kos{K}$-tensor natural endomorphisms
    $(\kos{c}^*\phi,\what{\kos{c}^*\phi})
    \Rightarrow
    (\kos{c}^*\phi,\what{\kos{c}^*\phi})
    :(\kos{T},\mon{t})\to (\kos{K}_c,\kos{c}^*)$.
    
    \item 
    We define the \emph{presheaf of groups of comonoidal $\kos{K}$-tensor natural automorphisms
    of $(\phi,\what{\phi})$}
    as the presheaf of groups on $\cat{Ens}(\kos{K})$
    \begin{equation*}
      \underline{\Aut}_{\mathbb{SMC}_{\cat{colax}}^{\kos{K}/\!\!/}}
      (\phi,\what{\phi})
      :
      \cat{Ens}(\kos{K})^{\op}
      \to
      \cat{Grp}
    \end{equation*}
    which sends each object $c$ in $\cat{Ens}(\kos{K})$
    to the group of comonoidal $\kos{K}$-tensor natural automorphisms
    $\xymatrix@C=12pt{(\kos{c}^*\phi,\what{\kos{c}^*\phi})
    \ar@2{->}[r]^-{\cong}
    &(\kos{c}^*\phi,\what{\kos{c}^*\phi})
    :(\kos{T},\mon{t})\to (\kos{K}_c,\kos{c}^*).}$
  \end{enumerate}
\end{definition}

The following is the main result of this subsection,
which we obtain by applying the results in \textsection~\ref{subsec RSKtensoradj}.

\begin{proposition} \label{prop RepGrpEns(K) comparison functor}
  Let $(\kos{T},\mon{t})$
  be a strong $\kos{K}$-tensor category
  and let 
  $(\omega,\what{\omega})
  :(\kos{K},\id_{\kos{K}})\to (\kos{T},\mon{t})$
  be a reflective right-strong $\kos{K}$-tensor adjunction.
  \begin{equation*}
    \vcenter{\hbox{
      \xymatrix{
        (\kos{T},\mon{t})
        \ar@/^1pc/[d]^-{(\omega^*,\what{\omega}^*)}
        \\
        (\kos{K},\id_{\kos{K}})
        \ar@/^1pc/[u]^-{\text{ reflective }(\omega_!,\what{\omega}_!)}
      }
    }}
  \end{equation*}
  Assume further that
  $\omega^*(\varphi_{\slot}):
  \omega^*\omega_!\omega^*(\slot)
  \Rightarrow
  \omega^*(\omega_!(\kappa)\tensor \slot)
  :\CT\to \CK$
  is a natural isomorphism.
  \begin{enumerate}
    \item 
    We have a group object $\pi:=\omega^*\omega_!(\kappa)$ in $\kos{E\!n\!s}(\kos{K})$
    which represents the presheaf of groups
    of comonoidal $\kos{K}$-tensor natural
    automorphisms of $(\omega^*,\what{\omega}^*)$.
    \begin{equation*}
      \xymatrix@C=18pt{
        \Hom_{\cat{Ens}(\kos{K})}\bigl(\slot, \omega^*\omega_!(\kappa)\bigr)
        \ar@2{->}[r]^-{\cong}
        &\underline{\Aut}_{\mathbb{SMC}_{\cat{colax}}^{\kos{K}/\!\!/}}
        (\omega^*,\what{\omega}^*)
        :
        \cat{Ens}(\kos{K})^{\op}
        \to
        \cat{Grp}
      }
    \end{equation*}
    The product, unit, antipode morphisms of $\omega^*\omega_!(\kappa)$
    are explicitly given below.
    Recall the comonoidal $\kos{K}$-tensor natural automorphism
    $\varsigma^{\omega,\omega}$
    of $(\omega^*\omega_!,\what{\omega^*\omega_!})$
    defined in Lemma~\ref{lem RSKtensoradj antipode}.
    \begin{equation}\label{eq RepGrpEns(K) comparison functor}
      \begin{aligned}
        \pc_{\omega^*\omega_!(\kappa)}
        &:
        \xymatrix@C=40pt{
          \omega^*\omega_!(\kappa)
          \otimes
          \omega^*\omega_!(\kappa)
          \ar[r]^-{(\hatar{\omega^*\omega_!}_{\omega^*\omega_!(\kappa)})^{-1}}_-{\cong}
          &\omega^*\omega_!\omega^*\omega_!(\kappa)
          \ar[r]^-{\omega^*(\epsilon_{\omega_!(\kappa)})}
          &\omega^*\omega_!(\kappa)
        }
        \\
        u_{\omega^*\omega_!(\kappa)}
        &:
        \xymatrix{
          \kappa
          \ar[r]^-{\eta_{\kappa}}
          &\omega^*\omega_!(\kappa)
        }
        \\
        \varsigma_{\omega^*\omega_!(\kappa)}
        &:
        \xymatrix{
          \omega^*\omega_!(\kappa)
          \ar[r]^-{\varsigma^{\omega,\omega}_{\kappa}}_-{\cong}
          &\omega^*\omega_!(\kappa)
        }
      \end{aligned}
    \end{equation}
    
    \item 
    The strong $\kos{K}$-tensor functor
    $(\omega^*,\what{\omega}^*):
    (\kos{T},\mon{t})\to (\kos{K},\id_{\kos{K}})$
    factors through as a strong $\kos{K}$-tensor functor
    $(\widebreve{\omega}\!^*,\what{\widebreve{\omega}}^*)
    :(\kos{T},\mon{t})\to (\kos{R\!e\!p}(\pi),\kos{t}^*_{\pi})$.
    \begin{equation*}
      \vcenter{\hbox{
        \xymatrix@C=40pt{
          (\kos{T},\mon{t})
          \ar[r]^-{(\widebreve{\omega}\!^*,\what{\widebreve{\omega}}^*)}
          \ar@/_1pc/[dr]_-{(\omega^*,\what{\omega}^*)}
          &(\kos{R\!e\!p}(\pi),\kos{t}^*_{\pi})
          \ar[d]^-{(\omega_{\pi}^*,\what{\omega}_{\pi}^*)}
          \\
          \text{ }
          &(\kos{K},\id_{\kos{K}})
        }
      }}
      \qquad\qquad
      \pi=\omega^*\omega_!(\kappa)
    \end{equation*}
  \end{enumerate}
\end{proposition}
\begin{proof}
  Consider the case
  $(\omega^{\pr*},\what{\omega}^{\pr*})=(\omega^*,\what{\omega}^*)$
  in Proposition~\ref{prop RSKtensoradj HomRep,Hom=Isom}.
  Statements 1,3 of Proposition~\ref{prop RSKtensoradj HomRep,Hom=Isom}
  imply that we have an equality
  of presheaves on $\cat{Ens}(\kos{K})$
  \begin{equation*}
    \underline{\End}_{\mathbb{SMC}_{\cat{colax}}^{\kos{K}/\!\!/}}
    (\omega^*,\what{\omega}^*)
    =
    \underline{\Aut}_{\mathbb{SMC}_{\cat{colax}}^{\kos{K}/\!\!/}}
    (\omega^*,\what{\omega}^*)
    :\cat{Ens}(\kos{K})^{\op}\to\cat{Set}
  \end{equation*}
  which is represented by the object
  $\pi:=\omega^*\omega_!(\kappa)$ in $\cat{Ens}(\kos{K})$.
  The universal element is
  the comonoidal $\kos{K}$-tensor natural automorphism
  \begin{equation*}
    \xi:
    \xymatrix@C=15pt{
      (\pi^*\omega^*,\what{\pi^*\omega^*})
      \ar@2{->}[r]^-{\cong}
      &(\pi^*\omega^*,\what{\pi^*\omega^*})
      :(\kos{T},\mon{t})\to (\kos{K}_{\pi},\pi^*)
    }
  \end{equation*}
  whose component
  $\xi_X:\pi^*\omega^*(X)\to \pi^*\omega^*(X)$
  at each object $X$ in $\CT$ is
  \begin{equation*}
    \xi_X:
    \xymatrix@C=40pt{
      \omega^*\omega_!(\kappa)\otimes \omega^*(X)
      \ar[r]^-{(\hatar{\omega^*\omega_!}_{\omega^*(X)})^{-1}}_-{\cong}
      &\omega^*\omega_!\omega^*(X)
      \ar[r]^-{\omega^*(\epsilon_X)}
      &\omega^*(X)
      .
    }
  \end{equation*}
  Thus the representing object
  $\pi=\omega^*\omega_!(\kappa)$
  has a unique structure 
  of a group object in $\kos{E\!n\!s}(\kos{K})$
  so that the representation
  \begin{equation*}
    \xymatrix@C=15pt{
      \Hom_{\cat{Ens}(\kos{K})}\bigl(\slot,\omega^*\omega_!(\kappa)\bigr)
      \ar@2{->}[r]^-{\cong}
      &\underline{\Aut}_{\mathbb{SMC}_{\cat{colax}}^{\kos{K}/\!\!/}}(\omega^*,\what{\omega}^*)
      :\cat{Ens}(\kos{K})^{\op}\to \cat{Grp}
    }
  \end{equation*}
  becomes an isomorphism between presheaves of groups
  on $\cat{Ens}(\kos{K})$.
  One can check that the identity natural transformations of 
  $(\omega^*,\what{\omega^*})$
  corresponds to the unit morphism
  $u_{\omega^*\omega_!(\kappa)}:\kappa\xrightarrow{\eta_{\kappa}}\omega^*\omega_!(\kappa)$
  as described in (\ref{eq RepGrpEns(K) comparison functor}).
  By considering the case
  $(\omega^{\pr},\what{\omega}^{\pr})=(\omega,\what{\omega})$
  in Lemma~\ref{lem RSKtensoradj antipode}
  and the case
  $(\omega^{\pr},\what{\omega}^{\pr})=(\omega,\what{\omega})$,
  $(\omega^{\ppr*},\what{\omega}^{\ppr*})=(\omega^*,\what{\omega}^*)$
  in Lemma~\ref{lem RSKtensoradj composition},
  we obtain that the product, antipode morphisms of $\pi=\omega^*\omega_!(\kappa)$
  are also given as described in (\ref{eq RepGrpEns(K) comparison functor}).
  
  The given right-strong $\kos{K}$-tensor adjunction
  $(\omega,\what{\omega}):(\kos{K},\id_{\kos{K}})\to(\kos{T},\mon{t})$
  induces a colax $\kos{K}$-tensor monad
  $\langle\omega^*\omega_!,\what{\omega^*\omega_!}\rangle
  =(\omega^*\omega_!,\what{\omega^*\omega_!},\mu,\eta)$
  on $(\kos{K},\id_{\kos{K}})$
  where 
  \begin{equation*}
    \mu=\omega^*\epsilon\omega_!:(\omega^*\omega_!\omega^*\omega_!,\what{\omega^*\omega_!\omega^*\omega_!})\Rightarrow(\omega^*\omega_!,\what{\omega^*\omega_!}).
  \end{equation*}
  Recall the adjoint equivalence of categories
  $\CL\dashv\iota$
  in Corollary~\ref{cor RepGrpEns(K) rflCKTmonad=Mon(Ens(K))}.
  We have an equality
  $\CL\langle\omega^*\omega_!,\what{\omega^*\omega_!}\rangle
  =\omega^*\omega_!(\kappa)$
  of monoid objects in $\kos{E\!n\!s}(\kos{K})$
  as the product, unit morphisms of $\omega^*\omega_!(\kappa)$
  in (\ref{eq RepGrpEns(K) comparison functor})
  is the same as those described in (\ref{eq RepGrpEns(K) rflCKTmonad=Mon(Ens(K))}).
  The component of the adjunction unit
  of $\CL\dashv \iota$ in Corollary~\ref{cor RepGrpEns(K) rflCKTmonad=Mon(Ens(K))}
  at $\langle\omega^*\omega_!,\what{\omega^*\omega_!}\rangle$
  is the following isomorphism of colax $\kos{K}$-tensor monads on $(\kos{K},\id_{\kos{K}})$.
  \begin{equation}\label{eq2 RepGrpEns(K) comparison functor}
    \hatar{\omega^*\omega_!}:
    \xymatrix{
      \langle\omega^*\omega_!,\what{\omega^*\omega_!}\rangle
      \ar@2{->}[r]^-{\cong}
      &\langle\omega^*\omega_!(\kappa)\otimes,\what{\omega^*\omega_!(\kappa)\otimes}\rangle
    }
  \end{equation}
  Using the above isomorphism of colax $\kos{K}$-tensor monads,
  we can show that the strong $\kos{K}$-tensor functor
  $(\omega^*,\what{\omega}^*)$
  factors through as a strong $\kos{K}$-tensor functor
  $(\widebreve{\omega}\!^*,\what{\widebreve{\omega}}^*)
  :(\kos{T},\mon{t})\to (\kos{R\!e\!p}(\pi),\kos{t}^*_{\pi})$.
  The underlying functor $\widebreve{\omega}\!^*$ sends each object $X$ in $\CT$
  to the object
  $\widebreve{\omega}\!^*(X)=(\omega^*(X),\gamma_{\omega^*(X)})$
  in $\cat{Rep}(\pi)$
  where
  \begin{equation*}
    \gamma_{\omega^*(X)}=\xi_X:
    \!\!
    \xymatrix@C=40pt{
      \omega^*\omega_!(\kappa)\otimes \omega^*(X)
      \ar[r]^-{(\hatar{\omega^*\omega_!}_{\omega^*(X)})^{-1}}_-{\cong}
      &\omega^*\omega_!\omega^*(X)
      \ar[r]^-{\omega^*(\epsilon_X)}
      &\omega^*(X)
      .
    }
  \end{equation*}
  We can check that $\widebreve{\omega}\!^*(X)$ is a well-defined object in $\cat{Rep}(\pi)$
  as follows.
  Let us denote $\pi=\omega^*\omega_!(\kappa)$
  and $(\phi,\what{\phi})=(\omega^*\omega_!,\what{\omega^*\omega_!})$
  as well as $\mu=\omega^*\epsilon\omega_!:\phi\phi\Rightarrow\phi$.
  \begin{equation*}
    \vcenter{\hbox{
      \xymatrix@C=40pt{
        \pi\otimes (\pi\otimes \omega^*(X))
        \ar[dd]_-{\mu^{\pi\otimes}_{\omega^*(X)}}
        \ar@{=}[r]
        &\pi\otimes (\pi\otimes \omega^*(X))
        \ar[d]^-{I_{\pi}\otimes (\hatar{\phi}_{\omega^*(X)})^{-1}}_-{\cong}
        \ar@{=}[r]
        &\pi\otimes (\pi\otimes \omega^*(X))
        \ar[dd]^-{I_{\pi}\otimes \gamma_{\omega^*(X)}}
        \\
        \text{ }
        &\pi\otimes \phi\omega^*(X)
        \ar@/^0.5pc/[dr]^-{I_{\pi}\otimes \omega^*(\epsilon_X)}
        \ar[d]^-{(\hatar{\phi}_{\phi\omega^*(X)})^{-1}}_-{\cong}
        &\text{ }
        \\
        \pi\otimes \omega^*(X)
        \ar[dd]_-{\gamma_{\omega^*(X)}}
        \ar@/_0.5pc/[dr]_-{(\hatar{\phi}_{\omega^*(X)})^{-1}}^-{\cong}
        &\phi\phi\omega^*(X)
        \ar[d]^-{\mu_{\omega^*(X)}}
        &\pi\otimes \omega^*(X)
        \ar@/^0.5pc/[dl]^-{(\hatar{\phi}_{\omega^*(X)})^{-1}}_-{\cong}
        \ar[dd]^-{\gamma_{\omega^*(X)}}
        \\
        \text{ }
        &\phi\omega^*(X)
        \ar[d]^-{\omega^*(\epsilon_X)}
        &\text{ }
        \\
        \omega^*(X)
        \ar@{=}[r]
        &\omega^*(X)
        \ar@{=}[r]
        &\omega^*(X)
      }
    }}
  \end{equation*}
  \begin{equation*}
    \vcenter{\hbox{
      \xymatrix{
        \omega^*(X)
        \ar[d]_-{\eta^{\pi\otimes}_{\omega^*(X)}}
        \ar@{=}[r]
        &\omega^*(X)
        \ar[dd]^-{\eta_{\omega^*(X)}}
        \ar@{=}[r]
        &\omega^*(X)
        \ar@{=}[ddd]
        \\
        \pi\otimes\omega^*(X)
        \ar[dd]_-{\gamma_{\omega^*(X)}}
        \ar@/_0.5pc/[dr]^-{(\hatar{\phi}_{\omega^*(X)})^{-1}}_-{\cong}
        &\text{ }
        &\text{ }
        \\
        \text{ }
        &\omega^*\omega_!\omega^*(X)
        \ar@/_0.5pc/[dr]^-{\omega^*(\epsilon_X)}
        &\text{ }
        \\
        \omega^*(X)
        \ar@{=}[rr]
        &\text{ }
        &\omega^*(X)
      }
    }}
  \end{equation*}
  Thus the underlying functor $\widebreve{\omega}\!^*$ is well-defined.
  Let $Y$ be another object in $\CT$.
  The symmetric monoidal coherernce isomorphism
  $\omega^*_{X,Y}:\omega^*(X\tensor Y)\xrightarrow{\cong}\omega^*(X)\otimes \omega^*(Y)$
  of $\omega^*:\kos{T}\to\kos{K}$
  becomes an isomorphism
  $\widebreve{\omega}\!^*_{X,Y}:
  \widebreve{\omega}\!^*(X\tensor Y)\xrightarrow{\cong}\widebreve{\omega}\!^*(X)\tensor\!_{\pi} \widebreve{\omega}^*(Y)$
  in $\cat{Rep}(\pi)$
  as we can see from the diagrams below.
  We begin our calculation as
  \begin{equation*}
    \vcenter{\hbox{
      \xymatrix@C=30pt{
        \pi\otimes \omega^*(X\tensor Y)
        \ar[ddd]_-{\gamma_{\omega^*(X\tensor Y)}}
        \ar@{=}[r]
        &\pi\otimes \omega^*(X\tensor Y)
        \ar[d]^-{(\hatar{\phi}_{\omega^*(X\tensor Y)})^{-1}}_-{\cong}
        \ar@{=}[r]
        &\pi\otimes \omega^*(X\tensor Y)
        \ar[d]^-{I_{\pi}\otimes \omega^*_{X,Y}}_-{\cong}
        \\
        \text{ }
        &\phi\omega^*(X\tensor Y)
        \ar[d]|-{\omega^*((\omega_!\omega^*)_{X,Y})}
        \ar@/_2pc/[ddl]|-{\omega^*(\epsilon_{X\tensor Y})}
        \ar@/^0.5pc/[dr]^(0.45){\phi(\omega^*_{X,Y})}_(0.45){\cong}
        &\pi\otimes (\omega^*(X)\otimes \omega^*(Y))
        \ar[d]^-{(\hatar{\phi}_{\omega^*(X)\otimes \omega^*(Y)})^{-1}}_-{\cong}
        \\
        \text{ }
        &\omega^*(\omega_!\omega^*(X)\tensor \omega_!\omega^*(Y))
        \ar@/^0.5pc/[dl]^-{\omega^*(\epsilon_X\tensor \epsilon_Y)}
        \ar@/_0.5pc/[dr]_-{\omega^*_{\omega_!\omega^*(X),\omega_!\omega^*(Y)}}^-{\cong}
        &\phi(\omega^*(X)\otimes \omega^*(Y))
        \ar[d]^-{\phi_{\omega^*(X),\omega^*(Y)}}
        \\
        \omega^*(X\tensor Y)
        \ar[d]_-{\omega^*_{X,Y}}^-{\cong}
        &\text{ }
        &\phi\omega^*(X)\otimes \phi\omega^*(Y)
        \ar[d]^-{\omega^*(\epsilon_X)\otimes \omega^*(\epsilon_Y)}
        \\
        \omega^*(X)\otimes \omega^*(Y)
        \ar@{=}[rr]
        &\text{ }
        &\omega^*(X)\otimes \omega^*(Y)
      }
    }}
  \end{equation*}
  and we finish as follows.
  \begin{equation*}
    \vcenter{\hbox{
      \xymatrix@C=10pt{
        \pi\otimes \omega^*(X\tensor Y)
        \ar[d]_-{I_{\pi}\otimes \omega^*_{X,Y}}^-{\cong}
        \ar@{=}[rr]
        &\text{ }
        &\pi\otimes \omega^*(X\tensor Y)
        \ar[d]^-{I_{\pi}\otimes \omega^*_{X,Y}}_-{\cong}
        \\
        \pi\otimes (\omega^*(X)\otimes \omega^*(Y))
        \ar[d]_-{(\hatar{\phi}_{\omega^*(X)\otimes \omega^*(Y)})^{-1}}^-{\cong}
        \ar@{=}[r]
        &\pi\otimes (\omega^*(X)\otimes \omega^*(Y))
        \ar[d]^-{(\pi\otimes)_{\omega^*(X),\omega^*(Y)}}
        \ar@{=}[r]
        &\pi\otimes (\omega^*(X)\otimes \omega^*(Y))
        \ar[ddd]^-{\gamma_{\omega^*(X)\otimes\omega^*(Y)}}
        \\
        \phi(\omega^*(X)\otimes \omega^*(Y))
        \ar[d]_-{\phi_{\omega^*(X),\omega^*(Y)}}
        &(\pi\otimes \omega^*(X))\otimes (\pi\otimes \omega^*(Y))
        \ar@/^1pc/[dl]^-{\quad(\hatar{\phi}_{\omega^*(X)})^{-1}\otimes (\hatar{\phi}_{\omega^*(Y)})^{-1}}_-{\cong}
        \ar@/^1pc/[ddr]|-{\gamma_{\omega^*(X)}\otimes \gamma_{\omega^*(Y)}}
        &\text{ }
        \\
        \phi\omega^*(X)\otimes \phi\omega^*(Y)
        \ar[d]_-{\omega^*(\epsilon_X)\otimes \omega^*(\epsilon_Y)}
        &\text{ }
        &\text{ }
        \\
        \omega^*(X)\otimes \omega^*(Y)
        \ar@{=}[rr]
        &\text{ }
        &\omega^*(X)\otimes \omega^*(Y)
      }
    }}
  \end{equation*}
  The symmetric monoidal coherence isomorphism
  $\omega^*_{\unit}:\omega^*(\unit)\xrightarrow{\cong}\kappa$
  of $\omega^*:\kos{T}\to\kos{K}$
  also becomes an isomorphism
  $\widebreve{\omega}\!^*_{\unit}:\widebreve{\omega}\!^*(\unit)\xrightarrow{\cong}\unit\!_{\pi}$
  in $\cat{Rep}(\pi)$
  as we can see from the diagram below.
  \begin{equation*}
    \vcenter{\hbox{
      \xymatrix@C=40pt{
        \pi\otimes \omega^*(\unit)
        \ar[dd]_-{\gamma_{\omega^*(\unit)}}
        \ar@{=}[r]
        &\pi\otimes \omega^*(\unit)
        \ar[d]^-{(\hatar{\phi}_{\omega^*(\unit)})^{-1}}_-{\cong}
        \ar@{=}[rr]
        &\text{ }
        &\pi\otimes \omega^*(\unit)
        \ar[d]^-{I_{\pi}\otimes \omega^*_{\unit}}_-{\cong}
        \\
        \text{ }
        &\phi\omega^*(\unit)
        \ar@/^0.5pc/[dl]|-{\omega^*(\epsilon_{\unit})=\omega^*((\omega_!\omega^*)_{\unit})}
        \ar[dd]^-{(\phi\omega^*)_{\unit}}
        \ar@/^0.5pc/[dr]^-{\phi(\omega^*_{\unit})}_-{\cong}
        &\text{ }
        &\pi\otimes \kappa
        \ar@/_0.5pc/[dl]_-{(\hatar{\phi}_{\kappa})^{-1}}^-{\cong}
        \ar[dd]^-{(\pi\otimes)_{\kappa}=\gamma_{\kappa}}
        \\
        \omega^*(\unit)
        \ar[d]_-{\omega^*_{\unit}}^-{\cong}
        &\text{ }
        &\phi(\kappa)
        \ar[d]^-{\phi_{\kappa}}
        &\text{ }
        \\
        \kappa
        \ar@{=}[r]
        &\kappa
        \ar@{=}[r]
        &\kappa
        \ar@{=}[r]
        &\kappa
      }
    }}
  \end{equation*}
  Finally, the comonoidal natural isomorphism
  $\what{\omega}^*:\omega^*\text{$\mon{t}$}\cong\id_{\kos{K}}:\kos{K}\to\kos{K}$
  becomes a comonoidal natural isomorphism
  $\what{\widebreve{\omega}}^*:\widebreve{\omega}\!^*\text{$\mon{t}$}\cong \kos{t}^*_{\pi}:\kos{K}\to\kos{R\!e\!p}(\pi)$
  as we can see from the following diagram.
  Let $z$ be an object in $\CK$.
  \begin{equation*}
    \vcenter{\hbox{
      \xymatrix@C=40pt{
        \pi\otimes \omega^*\mon{t}(z)
        \ar[dd]_-{\gamma_{\omega^*\mon{t}(z)}}
        \ar@{=}[r]
        &\pi\otimes \omega^*\mon{t}(z)
        \ar[d]^-{(\hatar{\phi}_{\omega^*\mon{t}(z)})^{-1}}_-{\cong}
        \ar@{=}[rr]
        &\text{ }
        &\pi\otimes \omega^*\mon{t}(z)
        \ar[d]^-{I_{\pi}\otimes \what{\omega}^*_z}_-{\cong}
        \\
        \text{ }
        &\phi\omega^*\mon{t}(z)
        \ar@/^0.5pc/[dl]|-{\omega^*(\epsilon_{\mon{t}(z)})}
        \ar[dd]^-{\what{\phi\omega^*}_z}
        \ar@/^0.5pc/[dr]^-{\phi(\what{\omega}^*_z)}_-{\cong}
        &\text{ }
        &\pi\otimes z
        \ar@/_0.5pc/[dl]_-{(\hatar{\phi}_z)^{-1}}^-{\cong}
        \ar[dd]|-{\what{\pi\otimes}_z=\varepsilon^{\pi\otimes}_z=\gamma_z}
        \\
        \omega^*\mon{t}(z)
        \ar[d]_-{\what{\omega}^*_z}^-{\cong}
        &\text{ }
        &\phi(z)
        \ar[d]^-{\what{\phi}_z}
        &\text{ }
        \\
        z
        \ar@{=}[r]
        &z
        \ar@{=}[r]
        &z
        \ar@{=}[r]
        &z
      }
    }}
  \end{equation*}
  This completes the proof of Proposition~\ref{prop RepGrpEns(K) comparison functor}.
\qed\end{proof}

\newpage

\newpage
\section{Prekosmic Galois categories}
\label{sec PreGalCat}

Recall the definition of
right-strong symmetric monoidal (RSSM) adjunctions
in Definition~\ref{def RSKtensoradj RSSMadj}
as well as the $2$-category
$\mathbb{ADJ}_{\cat{right}}(\mathbb{SMC}_{\cat{colax}})$
of 
symmetric monoidal categories,
RSSM adjunctions
and morphisms of RSSM adjunctions
introduced in (\ref{eq RSKtensoradj RSSM 2-cat}).

\subsection{Galois prekosmoi}
\label{subsec Galprekosmoi}

\begin{definition} \label{def Galprekosmoi GALpre}
  A \emph{Galois prekosmos}
  is a symmetric monoidal category
  \begin{equation*}
    \kos{K}=(\CK,\otimes,\kappa)
  \end{equation*}
  whose underlying category $\CK$
  has reflexive coequalizers.
  Let $\kos{T}=(\CT,\tensor,\unit)$
  be another Galois prekosmos.
  A \emph{Galois morphism}
  from $\kos{T}$ to $\kos{K}$
  is a right-strong symmetric monoidal (RSSM) adjunction
  $\kos{f}:\kos{T}\to \kos{K}$.
  \begin{equation*}
    \vcenter{\hbox{
      \xymatrix{
        \kos{K}
        \ar@/^1pc/[d]^-{\kos{f}^*}
        \\
        \kos{T}
        \ar@/^1pc/[u]^-{\kos{f}_!}
      }
    }}
  \end{equation*}
  Let
  $\kos{g}:\kos{T}\to \kos{K}$
  be another Galois morphism.
  A \emph{Galois transformation}
  $\vartheta:\kos{f}\Rightarrow\kos{g}:\kos{T}\to\kos{K}$
  is a comonoidal natural transformation
  $\vartheta:\kos{f}^*\Rightarrow\kos{g}^*:\kos{K}\to\kos{T}$
  between right adjoints.
  We define the $2$-category
  $\mathbb{GAL}^{\cat{pre}}$
  of Galois prekosmoi,
  Galois morphisms and Galois transformations
  as the full sub-$2$-category of
  $\mathbb{ADJ}_{\cat{right}}(\mathbb{SMC}_{\cat{colax}})$.
  \begin{equation*}
    \mathbb{GAL}^{\cat{pre}}
    \hookrightarrow
    \mathbb{ADJ}_{\cat{right}}(\mathbb{SMC}_{\cat{colax}})
  \end{equation*}
\end{definition}

We introduced the $2$-category structure of
$\mathbb{GAL}^{\cat{pre}}$
in Definition~\ref{def Galprekosmoi GALpre}.
Let $\kos{T}$, $\kos{K}$ be Galois prekosmoi
and let
$\kos{f}$, $\kos{g}$, $\kos{h}:\kos{T}\to\kos{K}$
be Galois morphisms.
The vertical composition of Galois transformations
$\kos{f}\Rightarrow\kos{g}\Rightarrow\kos{h}$
is defined as the vertical composition
$\kos{f}^*\Rightarrow\kos{g}^*\Rightarrow\kos{h}^*$
between right adjoints.

\begin{definition}\label{def Galprekosmoi Galmor property}
  We say a Galois morphism
  $\kos{f}$
  between Galois prekosmoi is
  \begin{itemize}
    \item 
    \emph{connected}
    if the right adjoint $\kos{f}^*$
    is fully faithful.
  
    \item 
    \emph{locally connected}
    if the RSSM adjunction $\kos{f}$
    satisfies the projection formula:
    see Definition~\ref{def RSKtensoradj projformula}.
  
    \item
    \emph{\'{e}tale}
    if the RSSM adjunction $\kos{f}$ satisfies the projection formula
    and the left adjoint $\kos{f}_!$ is conservative;
  
    \item 
    \emph{surjective}
    if the right adjoint
    $\kos{f}^*$
    is conservative and preserves reflexive coequalizers.
  \end{itemize}  
\end{definition}

\begin{remark}
  The terminologies
  are borrowed from geometric morphisms of topoi.
  Let $\kos{f}$ be a Galois morphism whose domain and codomain are
  topoi, considered as cartesian Galois prekosmoi.
  Then $\kos{f}^*$ admits a right adjoint $\kos{f}_*$
  if and only if $\kos{f}$ is an essential geometric morphism.
  In this case,
  $\kos{f}$ is
  connected|locally-connected|\'{e}tale|surjective
  as a Galois morphism
  if and only if it is
  connected|locally-connected|\'{e}tale|surjective
  as a geometric morphism.
\end{remark}

\begin{lemma}
  Let $\kos{K}=(\CK,\otimes,\kappa)$ be a Galois prekosmos.
  Then the cartesian monoidal category
  $\kos{E\!n\!s}(\kos{K})=(\cat{Ens}(\kos{K}),\otimes,\kappa)$
  of cocommutative comonoids in $\kos{K}$
  is also a Galois prekosmos.
\end{lemma}


For the rest of this section \textsection~\ref{sec PreGalCat},
we fix a Galois prekosmos
$\kos{K}=(\CK,\otimes,\kappa)$.

\begin{definition} \label{def Galprekosmoi GalpreK}
  We define the $2$-category
  \begin{equation*}
    \mathbb{GAL}^{\cat{pre}}_{\kos{K}}
  \end{equation*}
  of Galois $\kos{K}$-prekosmoi,
  Galois $\kos{K}$-morphisms
  and Galois $\kos{K}$-transformations
  as the slice $2$-category of $\mathbb{GAL}^{\cat{pre}}$ over $\kos{K}$.
\end{definition}

Let us explain Definition~\ref{def Galprekosmoi GalpreK}
in detail.
A \emph{Galois $\kos{K}$-prekosmos}
$\mathfrak{T}=(\kos{T},\kos{t})$ is a pair
of a Galois prekosmos
$\kos{T}=(\CT,\tensor,\unit)$
and a Galois morphism 
\begin{equation*}
  \vcenter{\hbox{
    \xymatrix{
      \kos{K}
      \ar@/^1pc/[d]^-{\kos{t}^*}
      \\
      \kos{T}
      \ar@/^1pc/[u]^-{\kos{t}_!}
    }
  }}
  \qquad\quad
  \kos{t}:\kos{T}\to\kos{K}.
\end{equation*}
We can see $\kos{K}$ itself as a Galois $\kos{K}$-prekosmos
which we denote as
$\text{$\mathfrak{K}$}=(\kos{K},\id_{\kos{K}})$.
Let $\mathfrak{S}=(\kos{S},\kos{s})$
be another Galois $\kos{K}$-prekosmos
and denote
$\text{$\kos{S}$}=(\CS,\ctimes,\pzc{1})$
as the underlying Galois prekosmos.
A \emph{Galois $\kos{K}$-morphism}
$\mathfrak{f}=(\kos{f},\what{\kos{f}}^*):
\mathfrak{S}\to \mathfrak{T}$
is a pair of a Galois morphism
$\kos{f}:\kos{S}\to\kos{T}$
and an invertible Galois transformation
\begin{equation*}
  \xymatrix@C=15pt{
    \what{\kos{f}}^*:\kos{t}\kos{f}\ar@2{->}[r]^-{\cong}& \kos{s}:\kos{S}\to \kos{K}
  }  
\end{equation*}
which is a comonoidal natural isomorphism between right adjoints
\begin{equation*}
  \vcenter{\hbox{
    \xymatrix@R=30pt@C=40pt{
      \text{ }
      &\kos{T}
      \ar[d]^-{\kos{f}^*}
      \\
      \kos{K}
      \ar@/^0.7pc/[ur]^-{\kos{t}^*}
      \ar[r]_-{\kos{s}^*}
      \xtwocell[r]{}<>{<-3>{\text{ }\text{ }\what{\kos{f}}^*}}
      &\kos{S}
    }
  }}
  \qquad\quad
  \xymatrix@C=15pt{
    \what{\kos{f}}^*:
    \kos{f}^*\kos{t}^*
    \ar@2{->}[r]^-{\cong}
    &\kos{s}^*
    :\kos{K}\to \kos{S}
    .
  }
\end{equation*}
Let $\mathfrak{g}=(\kos{g},\what{\kos{g}}^*):
\mathfrak{S}\to\mathfrak{T}$
be another Galois $\kos{K}$-morphism.
A \emph{Galois $\kos{K}$-transformation}
$\vartheta:\mathfrak{f}\Rightarrow\mathfrak{g}:\mathfrak{S}\to\mathfrak{T}$
is a Galois transformation
$\vartheta:\kos{f}\Rightarrow\kos{g}:\kos{S}\to\kos{T}$
which satisfies the relation
\begin{equation} \label{eq Galprekosmoi GalKtransformation}
  \vcenter{\hbox{
    \xymatrix@C=15pt{
      \kos{f}^*\kos{t}^*
      \ar@2{->}[rr]^-{\vartheta\kos{t}^*}
      \ar@2{->}[dr]_-{\what{\kos{f}}^*}^-{\cong}
      &\text{ }
      &\kos{g}^*\kos{t}^*
      \ar@2{->}[dl]^-{\what{\kos{g}}^*}_-{\cong}
      \\
      \text{ }
      &\kos{s}^*
    }
  }}
  \qquad\quad
  \what{\kos{f}}^*
  =
  \what{\kos{g}}^*\circ (\vartheta^*\kos{t}^*)
  :\kos{K}\to \kos{S}
  .
\end{equation}
We can also describe relation (\ref{eq Galprekosmoi GalKtransformation})
as follows.
\begin{equation*}
  \vcenter{\hbox{
    \xymatrix@R=30pt@C=40pt{
      \text{ }
      &\kos{T}
      \ar[d]^-{\kos{f}^*}
      \\
      \kos{K}
      \ar@/^0.7pc/[ur]^-{\kos{t}^*}
      \ar[r]_-{\kos{s}^*}
      \xtwocell[r]{}<>{<-3>{\text{ }\text{ }\what{\kos{f}}^*}}
      &\kos{S}
    }
  }}
  \quad=
  \vcenter{\hbox{
    \xymatrix@R=30pt@C=25pt{
      \text{ }
      &\text{ }
      &\kos{T}
      \ar@/^1pc/[d]^-{\kos{f}^*}
      \ar@/_1pc/[d]_-{\kos{g}^*}
      \xtwocell[d]{}<>{<0>{\vartheta}}
      \\
      \kos{K}
      \ar@/^1pc/[urr]^-{\kos{t}^*}
      \ar[rr]_-{\kos{s}^*}
      \xtwocell[rr]{}<>{<-2>{\text{ }\text{ }\what{\kos{g}}^*}}
      &\text{ }
      &\kos{S}
    }
  }}
\end{equation*}

Recall the $2$-category
$\mathbb{ADJ}_{\cat{right}}\bigl(\mathbb{SMC}_{\cat{colax}}^{\kos{K}/\!\!/}\bigr)$
of colax $\kos{K}$-tensor categories,
right-strong $\kos{K}$-tensor adjunctions,
and morphisms of right-strong $\kos{K}$-tensor adjunctions
introduced in (\ref{eq RSKtensoradj RSKTadj 2-cat}).
The $2$-category
$\mathbb{GAL}^{\cat{pre}}_{\kos{K}}$
of Galois $\kos{K}$-prekosmoi
is a full sub-$2$-category of
$\mathbb{ADJ}_{\cat{right}}(\mathbb{SMC}_{\cat{colax}}^{\kos{K}/\!\!/})$.
\begin{equation} \label{eq Galprekosmoi GalpreK RSKTadj}
  \mathbb{GAL}^{\cat{pre}}_{\kos{K}}
  \hookrightarrow
  \mathbb{ADJ}_{\cat{right}}\bigl(\mathbb{SMC}_{\cat{colax}}^{\kos{K}/\!\!/}\bigr)
\end{equation}
\begin{itemize}
  \item 
  Let $\mathfrak{T}=(\kos{T},\kos{t})$
  be a Galois $\kos{K}$-prekosmos.
  Then the pair $(\kos{T},\kos{t}^*)$ is a
  strong $\kos{K}$-tensor category.

  \item 
  Let 
  $\mathfrak{S}=(\kos{S},\kos{s})$
  be another Galois $\kos{K}$-prekosmos
  and let
  $\mathfrak{f}=(\kos{f},\what{\kos{f}}^*):
  \mathfrak{S}\to \mathfrak{T}$
  be a Galois $\kos{K}$-morphism.
  Then the pair
  $(\kos{f}^*,\what{\kos{f}}^*):
  (\kos{T},\kos{t}^*)\to (\kos{S},\kos{s}^*)$
  is a strong $\kos{K}$-tensor functor.
  Thus the left adjoint has a unique colax $\kos{K}$-tensor functor structure
  $(\kos{f}_!,\what{\kos{f}}_!):
  (\kos{S},\kos{s}^*)\to (\kos{T},\kos{t}^*)$
  such that the given Galois morphism
  $\kos{f}:\kos{S}\to \kos{T}$
  becomes a right-strong $\kos{K}$-tensor adjunction
  \begin{equation} \label{eq Galprekosmoi GalKmor associated RSKTadj}
    \vcenter{\hbox{
      \xymatrix{
        (\kos{T},\kos{t}^*)
        \ar@/^1pc/[d]^-{(\kos{f}^*,\what{\kos{f}}^*)}
        \\
        (\kos{S},\kos{s}^*)
        \ar@/^1pc/[u]^-{(\kos{f}_!,\what{\kos{f}}_!)}
      }
    }}
    \qquad\qquad
    (\kos{f},\what{\kos{f}}):
    (\kos{S},\kos{s}^*)\to (\kos{T},\kos{t}^*)
    .
  \end{equation}
  The comonoidal natural transformation
  $\what{\kos{f}}_!$ is 
  given as in (\ref{eq RSKtensoradj leftadj uniqueKTstr}).
  \begin{equation*}
    \vcenter{\hbox{
      \xymatrix@R=30pt@C=40pt{
        \text{ }
        &\kos{S}
        \ar[d]^-{\kos{f}_!}
        \\
        \kos{K}
        \ar@/^0.7pc/[ur]^-{\kos{s}^*}
        \ar[r]_-{\kos{t}^*}
        \xtwocell[r]{}<>{<-3>{\text{ }\text{ }\what{\kos{f}}_!}}
        &\kos{T}
      }
    }}
    \qquad\quad
    \xymatrix@C=30pt{
      \what{\kos{f}}_!:
      \kos{f}_!\kos{s}^*
      \ar@2{->}[r]^-{\kos{f}_!\what{\kos{f}}^{*-1}}_-{\cong}
      &\kos{f}_!\kos{f}^*\kos{t}^*
      \ar@2{->}[r]^-{\epsilon^{\kos{f}}\kos{t}^*}
      &\kos{t}^*
    }
  \end{equation*}

  \item 
  Let $\mathfrak{g}=(\kos{g},\what{\kos{g}}^*):\mathfrak{S}\to\mathfrak{T}$
  be another Galois $\kos{K}$-morphism.
  Then a Galois $\kos{K}$-transformation
  $\vartheta:\mathfrak{f}\Rightarrow\mathfrak{g}:\mathfrak{S}\to\mathfrak{T}$
  is precisely 
  a comonoidal $\kos{K}$-tensor natural transformation
  $\vartheta:(\kos{f}^*,\what{\kos{f}}^*)\Rightarrow (\kos{g}^*,\what{\kos{g}}^*)
  :(\kos{T},\kos{t}^*)\to (\kos{S},\kos{s}^*)$.
\end{itemize}
From now on, we denote a Galois $\kos{K}$-morphism
$\mathfrak{f}:\mathfrak{S}\to \mathfrak{T}$
between Galois $\kos{K}$-prekosmoi
$\mathfrak{T}=(\kos{T},\kos{t})$,
$\mathfrak{S}=(\kos{S},\kos{s})$
as
\begin{equation*}
  \mathfrak{f}=(\kos{f},\what{\kos{f}}):\mathfrak{S}\to\mathfrak{T}
\end{equation*}
where
$(\kos{f},\what{\kos{f}}):
(\kos{S},\kos{s}^*)\to (\kos{T},\kos{t}^*)$
is the associated right-strong $\kos{K}$-tensor adjunction
described in (\ref{eq Galprekosmoi GalKmor associated RSKTadj}).

\begin{definition}
  A Galois $\kos{K}$-prekosmoi 
  $\mathfrak{T}=(\kos{T},\kos{t})$
  is called
  \begin{equation*}
    \text{ 
      \emph{connected|locally connected|\'{e}tale|surjective}
    }
  \end{equation*}
  if the Galois morphism
  $\kos{t}:\kos{T}\to \kos{K}$ has such property.
\end{definition}

\begin{definition}
  Let
  $\mathfrak{T}=(\kos{T},\kos{t})$,
  $\mathfrak{S}=(\kos{S},\kos{s})$
  be Galois $\kos{K}$-prekosmoi.
  A Galois $\kos{K}$-morphism
  $\mathfrak{f}=(\kos{f},\what{\kos{f}}):\mathfrak{S}\to\mathfrak{T}$
  is called \emph{connected|locally connected|\'{e}tale|surjective}
  if the underlying Galois morphism $\kos{f}:\kos{S}\to\kos{T}$ has such property.
\end{definition}

\begin{definition}
  Let $\mathfrak{T}=(\kos{T},\kos{t})$
  be a Galois $\kos{K}$-prekosmos
  and recall that we denote
  $\text{$\mathfrak{K}$}=(\kos{K},\id_{\kos{K}})$.
  We say a Galois $\kos{K}$-morphism
  \begin{equation*}
    \varpi=(\omega,\what{\omega})
    :\mathfrak{K}\to\mathfrak{T}
  \end{equation*}
  is \emph{reflective}
  if the associated right-strong $\kos{K}$-tensor adjunction
  $(\omega,\what{\omega}):(\kos{K},\id_{\kos{K}})\to (\kos{T},\kos{t}^*)$
  is reflective:
  see Definition~\ref{def RSKtensoradj reflective}.
  This amounts to saying that
  the left adjoint
  $(\omega_!,\what{\omega}_!):
  (\kos{K},\id_{\kos{K}})\to (\kos{T},\kos{t}^*)$
  is a reflective colax $\kos{K}$-tensor functor,
  or equivalently,
  the associated $\kos{K}$-equivariance
  $\vecar{\omega}_!$
  is a natural isomorphism:
  see Lemma~\ref{lem Ens(T) equivariance and reflection}.
\end{definition}


\begin{definition}
  Let $\mathfrak{T}=(\kos{T},\kos{t})$
  be a Galois $\kos{K}$-prekosmos.
  \begin{enumerate}
    \item 
    A \emph{pre-fiber functor} for $\mathfrak{T}$
    is a reflective \'{e}tale Galois $\kos{K}$-morphism
    \begin{equation*}
      \varpi=(\omega,\what{\omega}):\mathfrak{K}\to \mathfrak{T}.
    \end{equation*}
    Equivalently, it is a Galois $\kos{K}$-morphism
    $\varpi=(\omega,\what{\omega}):\mathfrak{K}\to\mathfrak{T}$
    which has the following properties:
    \begin{itemize}
      \item 
      the Galois morphism $\omega:\kos{K}\to \kos{T}$
      satisfies the projection formula;
      
      \item 
      the left adjoint $(\omega_!,\what{\omega}_!):(\kos{K},\id_{\kos{K}})\to (\kos{T},\kos{t}^*)$
      is a reflective colax $\kos{K}$-tensor functor,
      and the underlying functor
      $\omega_!$ is conservative.
    \end{itemize}
    
    \item
    A pre-fiber functor $\varpi=(\omega,\what{\omega}):\mathfrak{K}\to\mathfrak{T}$
    is called \emph{surjective}
    if it is surjective as a Galois $\kos{K}$-morphism.
    This amounts to saying that
    the right adjoint $\omega^*$ is conservative and preserves reflexive coequalizers.
  \end{enumerate}
  A \emph{morphism} of pre-fiber functors $\varpi$, $\varpi^{\pr}$ for $\mathfrak{T}$
  is a Galois $\kos{K}$-transformation 
  $\vartheta:\varpi\Rightarrow\varpi^{\pr}:\mathfrak{K}\to\mathfrak{T}$.
  We denote $\cat{Fib}(\mathfrak{T})^{\cat{pre}}$
  as the category of pre-fiber functors for $\mathfrak{T}$.
\end{definition}


\begin{definition}
  Let $\mathfrak{T}=(\kos{T},\kos{t})$ be a Galois $\kos{K}$-prekosmos
  and let $\varpi=(\omega,\what{\omega})$,
  $\varpi^{\pr}=(\omega^{\pr},\what{\omega}^{\pr}):\mathfrak{K}\to\mathfrak{T}$
  be pre-fiber functors.
  Consider the strong $\kos{K}$-tensor functors
  $(\omega^*,\what{\omega}^*)$,
  $(\omega^{\pr*},\what{\omega}^{\pr*})
  :(\kos{T},\kos{t}^*)\to (\kos{K},\id_{\kos{K}})$
  and recall Definition~\ref{def RSKtensoradj Hom,Isom}
  as well as Definition~\ref{def RepGrpEns(K) End,Aut}.
  \begin{enumerate}
    \item 
    We define the \emph{presheaf of Galois $\kos{K}$-transformations
    from $\varpi$ to $\varpi^{\pr}$}
    as
    \begin{equation*}
      \underline{\Hom}_{\mathbb{GAL}^{\cat{pre}}_{\kos{K}}}(\varpi,\varpi^{\pr})
      :=
      \underline{\Hom}_{\mathbb{SMC}^{\kos{K}/\!\!/}_{\cat{colax}}}
      \bigl((\omega^*,\what{\omega}^*),(\omega^{\pr*},\what{\omega}^{\pr*})\bigr)
      :\cat{Ens}(\kos{K})^{\op}\to\cat{Set}
    \end{equation*}
    and the \emph{presheaf of monoids of Galois $\kos{K}$-transformations
    from $\varpi$ to $\varpi$} as
    \begin{equation*}
      \underline{\End}_{\mathbb{GAL}^{\cat{pre}}_{\kos{K}}}(\varpi)
      :=
      \underline{\End}_{\mathbb{SMC}^{\kos{K}/\!\!/}_{\cat{colax}}}
      (\omega^*,\what{\omega}^*)
      :\cat{Ens}(\kos{K})^{\op}\to\cat{Mon}
      . 
    \end{equation*}

    \item 
    We define the \emph{presheaf of invertible Galois $\kos{K}$-transformations
    from $\varpi$ to $\varpi^{\pr}$}
    as
    \begin{equation*}
      \underline{\Isom}_{\mathbb{GAL}^{\cat{pre}}_{\kos{K}}}(\varpi,\varpi^{\pr})
      :=
      \underline{\Isom}_{\mathbb{SMC}^{\kos{K}/\!\!/}_{\cat{colax}}}
      \bigl((\omega^*,\what{\omega}^*),(\omega^{\pr*},\what{\omega}^{\pr*})\bigr)
      :\cat{Ens}(\kos{K})^{\op}\to\cat{Set}
    \end{equation*}
    and the \emph{presheaf of groups of invertible Galois $\kos{K}$-transformations
    from $\varpi$ to $\varpi$} as
    \begin{equation*}
      \underline{\Aut}_{\mathbb{GAL}^{\cat{pre}}_{\kos{K}}}(\varpi)
      :=
      \underline{\Aut}_{\mathbb{SMC}^{\kos{K}/\!\!/}_{\cat{colax}}}
      (\omega^*,\what{\omega}^*)
      :\cat{Ens}(\kos{K})^{\op}\to\cat{Grp}
      . 
    \end{equation*}
  \end{enumerate}
\end{definition}

\begin{remark}
  Let $\mathfrak{T}=(\kos{T},\kos{t})$ be a Galois $\kos{K}$-prekosmos
  and let $\varpi=(\omega,\what{\omega})$,
  $\varpi^{\pr}=(\omega^{\pr},\what{\omega}^{\pr}):\mathfrak{K}\to\mathfrak{T}$
  be pre-fiber functors.
  For each object $c$ in $\cat{Ens}(\kos{K})$,
  an element $\vartheta$ in 
  $\underline{\Hom}_{\mathbb{GAL}^{\cat{pre}}_{\kos{K}}}(\varpi,\varpi^{\pr})(c)$
  is a comonoidal $\kos{K}$-tensor natural transformation
  \begin{equation*}
    \vartheta:(\kos{c}^*\omega^*,\what{\kos{c}^*\omega^*})\Rightarrow
    (\kos{c}^*\omega^{\pr*},\what{\kos{c}^*\omega^{\pr*}})
    :(\kos{T},\kos{t}^*)\to (\kos{K}_c,\kos{c}^*).
  \end{equation*}
  We cannot say that $\vartheta$ is a Galois $\kos{K}$-transformation,
  since the symmetric monoidal category $\kos{K}_c$ is in general not 
  a Galois prekosmos.
  Nontheless, 
  we are going to call 
  $\underline{\Hom}_{\mathbb{GAL}^{\cat{pre}}_{\kos{K}}}(\varpi,\varpi^{\pr})$
  as the presheaf of Galois $\kos{K}$-transformations
  from $\varpi$ to $\varpi^{\pr}$.
\end{remark}

The following is a consequence of results in \textsection~\ref{subsec RSKtensoradj}.

\begin{corollary} \label{cor Galprekosmoi summary}
  Let $\mathfrak{T}=(\kos{T},\kos{t})$ be a Galois $\kos{K}$-prekosmos
  and let $\varpi=(\omega,\what{\omega})$,
  $\varpi^{\pr}=(\omega^{\pr},\what{\omega}^{\pr}):\mathfrak{K}\to\mathfrak{T}$
  be pre-fiber functors.
  \begin{enumerate}
    \item 
    The presheaf of Galois $\kos{K}$-transformations from $\varpi$ to $\varpi^{\pr}$
    is represented by the object $\omega^{\pr*}\omega_!(\kappa)$ in $\cat{Ens}(\kos{K})$.
    \begin{equation*}
      \vcenter{\hbox{
        \xymatrix{
          \Hom_{\cat{Ens}(\kos{K})}\bigl(\slot,\omega^{\pr*}\omega_!(\kappa)\bigr)
          \ar@2{->}[r]^-{\cong}
          &\underline{\Hom}_{\mathbb{GAL}^{\cat{pre}}_{\kos{K}}}(\varpi,\varpi^{\pr})
          :\cat{Ens}(\kos{K})^{\op}\to\cat{Set}
        }
      }}
    \end{equation*}

    \item
    We have
    $\underline{\Hom}_{\mathbb{GAL}^{\cat{pre}}_{\kos{K}}}(\varpi,\varpi^{\pr})
    =
    \underline{\Isom}_{\mathbb{GAL}^{\cat{pre}}_{\kos{K}}}(\varpi,\varpi^{\pr})
    :\cat{Ens}(\kos{K})^{\op}\to\cat{Set}$.

    \item
    We have a comonoidal $\kos{K}$-tensor natural isomorphism
    \begin{equation*}
      \varsigma^{\omega,\omega^{\pr}}:
      \xymatrix@C=18pt{
        (\omega^{\pr*}\omega_!,\what{\omega^{\pr*}\omega_!})
        \ar@2{->}[r]^-{\cong}
        &(\omega^*\omega^{\pr}_!,\what{\omega^*\omega^{\pr}_!})
        :(\kos{K},\id_{\kos{K}})\to (\kos{K},\id_{\kos{K}})
      }
    \end{equation*}
    whose component at $\kappa$ is an isomorphism
    $\varsigma^{\omega,\omega^{\pr}}_{\kappa}:
    \omega^{\pr*}\omega_!(\kappa)
    \xrightarrow[]{\cong}
    \omega^*\omega^{\pr}_!(\kappa)$
    in $\cat{Ens}(\kos{K})$,
    and the following diagram of presheaves strictly commutes.
    \begin{equation*}
      \vcenter{\hbox{
        \xymatrix{
          \Hom_{\cat{Ens}(\kos{K})}\bigl(\slot,\omega^{\pr*}\omega_!(\kappa)\bigr)
          \ar@2{->}[d]_-{\varsigma^{\omega,\omega^{\pr}}_{\kappa}\circ(\slot)}^-{\cong}
          \ar@2{->}[r]^-{\cong}
          &\underline{\Isom}_{\mathbb{GAL}^{\cat{pre}}_{\kos{K}}}(\varpi,\varpi^{\pr})
          \ar@2{->}[d]^-{\cat{inverse}}_-{\cong}
          \\
          \Hom_{\cat{Ens}(\kos{K})}\bigl(\slot,\omega^*\omega^{\pr}_!(\kappa)\bigr)
          \ar@2{->}[r]^-{\cong}
          &\underline{\Isom}_{\mathbb{GAL}^{\cat{pre}}_{\kos{K}}}(\varpi^{\pr},\varpi)
        }
      }}
    \end{equation*}
  \end{enumerate}
  In particular, the category
  $\cat{Fib}(\mathfrak{T})^{\cat{pre}}$
  of pre-fiber functors for $\mathfrak{T}$ is a groupoid.
\end{corollary}
\begin{proof}
  The right-strong $\kos{K}$-tensor adjunctions
  $(\omega,\what{\omega})$,
  $(\omega^{\pr},\what{\omega}^{\pr})$
  associated to $\varpi$, $\varpi^{\pr}$
  are reflective and satisfy the projection formula.
  In particular, the associated natural transformations
  $\varphi_{\slot}$,
  $\varphi^{\pr}_{\slot}$
  defined as in Definition~\ref{def RSKtensoradj reflective}
  are natural isomorphisms.
  By 
  Proposition~\ref{prop RSKtensoradj HomRep,Hom=Isom}
  and 
  Lemma~\ref{lem RSKtensoradj antipode},
  we immediately obtain statements 1-3.
  Moreover, statement 2 implies that
  $\cat{Fib}(\mathfrak{T})_{\kappa}$
  is a groupoid.
\qed\end{proof}

\subsection{Pre-Galois objects in a Galois prekosmos $\kos{K}$}
\label{subsec preGalobj}

Let $\kos{K}=(\CK,\otimes,\kappa)$ be a Galois prekosmos.
Recall the definition of a group object $\pi$ in $\kos{E\!n\!s}(\kos{K})$
introduced in \textsection\!~\ref{subsec RepGrpEns(K)}

\begin{definition}
  A \emph{pre-Galois object} in $\kos{K}$
  is a group object $\pi$ in $\kos{E\!n\!s}(\kos{K})$
  such that the functor $\pi\otimes\slot:\CK\to\CK$
  preserves reflexive coequalizers.
\end{definition}

Let $\pi$ be a pre-Galois object in $\kos{K}$.
We denote the product, unit, antipode morphisms of $\pi$ as
$\pi\otimes \pi\xrightarrow{\pc_{\pi}} \pi$,
$\kappa\xrightarrow{u_{\pi}} \pi$,
$\pi\xrightarrow[\cong]{\varsigma_{\pi}}\pi$.
The functor $\pi\otimes\slot:\CK\to \CK$
is conservative
as explained in Remark~\ref{rem RepGrpEns(K) pitensor conservative}.

A \emph{morphism} $\pi^{\pr}\to\pi$
of pre-Galois objects in $\kos{K}$
is a morphism $\pi^{\pr}\to\pi$ of group objects in $\kos{E\!n\!s}(\kos{K})$.

A \emph{representation} of a pre-Galois object $\pi$ in $\kos{K}$
is a representation of $\pi$ as a group object in $\kos{E\!n\!s}(\kos{K})$:
it is a pair $(x,\gamma_x)$
of an object $x$ in $\CK$
and a left $\pi$-action morphism $\gamma_x:\pi\otimes x\to x$
in $\CK$.
Recall the strong $\kos{K}$-tensor category
$(\kos{R\!e\!p}(\kos{K}),\kos{t}_{\pi}^*)$
of representations of $\pi$ in $\kos{K}$
defined in Definition~\ref{def RepGrpEns(K) Rep(pi)Ktensorcat}.

\begin{lemma} \label{lem preGalobj Rep(pi)}
  Let $\pi$ be a pre-Galois object in $\kos{K}$
  and recall Lemma~\ref{lem RepGrpEns(K) Rep(pi)}.
  \begin{enumerate}
    \item 
    We have a connected Galois $\kos{K}$-prekosmos
    $\mathfrak{Rep}(\pi)=(\kos{R\!e\!p}(\pi),\kos{t}_{\pi})$
    of representations of $\pi$.

    \item 
    We have a surjective pre-fiber functor
    $\varpi_{\pi}=(\omega_{\pi},\what{\omega}_{\pi}):\mathfrak{K}\to \mathfrak{Rep}(\pi)$
    whose right adjoint
    $\omega_{\pi}^*$ is the functor of forgetting the actions of $\pi$.
  \end{enumerate}
\end{lemma}
\begin{proof}
  Recall the symmetric monoidal category
  $\kos{R\!e\!p}(\pi)=(\cat{Rep}(\pi),\tensor\!_{\pi},\unit\!_{\pi})$
  of representations of $\pi$.
  The category
  $\cat{Rep}(\pi)$ has reflexive coequalizers
  since the functor $\pi\otimes\slot:\CK\to \CK$
  preserves reflexive coequalizers.
  Thus $\kos{R\!e\!p}(\pi)$ is a Galois prekosmos.
  The fully faithful strong symmetric monoidal functor
  $\kos{t}_{\pi}^*:\kos{K}\to\kos{R\!e\!p}(\pi)$
  of constructing trivial representations
  admits a left adjoint
  $\kos{t}_{\pi!}$
  which sends each representation $X=(x,\gamma_x)$ of $\pi$
  to the following reflexive coequalizer in $\CK$.
  \begin{equation*}
    \vcenter{\hbox{
      \xymatrix@C=30pt{
        \pi\otimes x
        \ar@<0.5ex>[r]^-{\gamma_x}
        \ar@<-0.5ex>[r]_-{e_{\pi}\otimes I_x}
        &x
        \ar@/_2pc/@<-1ex>[l]|-{u_{\pi}\otimes I_x}
        \ar@{->>}[r]
        &\kos{t}_{\pi!}(X)
      }
    }}
  \end{equation*}
  Thus we obtain a connected Galois morphism
  \begin{equation*}
    \vcenter{\hbox{
      \xymatrix{
        \kos{K}
        \ar@/^1pc/[d]^{\kos{t}^*_{\pi}}
        \\
        \kos{R\!e\!p}(\pi)
        \ar@/^1pc/[u]^{\kos{t}_{\pi!}}
      }
    }}
    \qquad\quad
    \kos{t}_{\pi}:\kos{R\!e\!p}(\pi)\to \kos{K}
  \end{equation*}
  and the pair
  $\mathfrak{Rep}(\pi):=(\kos{R\!e\!p}(\pi),\kos{t}_{\pi})$
  is a connected Galois $\kos{K}$-prekosmos.
  Recall that in Lemma~\ref{lem RepGrpEns(K) Rep(pi)},
  we introduced the reflective right-strong $\kos{K}$-tensor adjunction
  \begin{equation*}
    \vcenter{\hbox{
      \xymatrix{
        (\kos{R\!e\!p}(\pi),\kos{t}_{\pi}^*)
        \ar@/^1pc/[d]^-{(\omega_{\pi}^*,\what{\omega}_{\pi}^*\text{ }\!)}
        \\
        (\kos{K},\id_{\kos{K}})
        \ar@/^1pc/[u]^-{(\omega_{\pi!},\what{\omega}_{\pi!})}
      }
    }}
    \qquad\quad
    (\omega_{\pi},\what{\omega}_{\pi}):
    (\kos{K},\id_{\kos{K}})\to (\kos{R\!e\!p}(\pi),\kos{t}^*_{\pi})
  \end{equation*}
  satisfying the projection formula,
  such that the left adjoint
  $(\omega_{\pi!},\what{\omega}_{\pi!})$
  is a reflective colax $\kos{K}$-tensor functor
  and the underlying functor $\omega_{\pi!}$ is conservative.
  Thus 
  \begin{equation*}
    \varpi_{\pi}=(\omega_{\pi},\what{\omega}_{\pi})
    :\mathfrak{K}\to\mathfrak{Rep}(\pi)
  \end{equation*}
  is a pre-fiber functor for $\mathfrak{Rep}(\pi)$.
  The forgetful functor
  $\omega^*$ is conservative,
  and it preserves reflexive coequalizers
  since the functor 
  $\pi\otimes\slot:\CK\to \CK$
  preserves reflexive coequalizers.
  We conclude that
  $\varpi_{\pi}:\mathfrak{K}\to\mathfrak{Rep}(\pi)$
  is a surjective pre-fiber functor.
  This completes the proof of Lemma~\ref{lem preGalobj Rep(pi)}.
\qed\end{proof}

Let $\pi$ be a pre-Galois object in $\kos{K}$.
In (\ref{eq RepGrpEns(K) Hom(-,pi)})
we explained the presheaf of groups
$\Hom_{\cat{Ens}(\kos{K})}(\slot,\pi)$
on $\cat{Ens}(\kos{K})$,
where the group operation $\star$
is the convolution product.
The group $\Hom_{\cat{Ens}(\kos{K})}(\kappa,\pi)$
acts on the presheaf of groups
$\Hom_{\cat{Ens}(\kos{K})}(\slot,\pi)$
by conjugation.
The conjugation action of each element
$\kappa\xrightarrow{\theta}\pi$
corresponds to the automorphism 
$\sigma_{\theta}:\pi\xrightarrow{\cong}\pi$
where
\begin{equation*}
  \sigma_{\theta}:
  \xymatrix@C=50pt{
    \pi
    \ar[r]^-{(\imath_p\otimes I_{\kappa})\circ \jmath_p}_-{\cong}
    &\kappa\otimes \pi\otimes \kappa
    \ar[r]^-{\theta\otimes I_{\pi}\otimes \theta}
    &\pi\otimes \pi\otimes \pi
    \ar[r]^-{\pc_{\pi}\circ (\pc_{\pi}\otimes \varsigma_{\pi})}
    &\pi
    .
  }
\end{equation*}
We call $\sigma_{\theta}$ as the \emph{inner automorphism} of $\pi$ associated to $\theta$.
We have
the relation
$\sigma_{\theta\star\tilde{\theta}}
:\pi\xrightarrow[\cong]{\sigma_{\tilde{\theta}}}\pi\xrightarrow[\cong]{\sigma_{\theta}}\pi$
for every pair
$\kappa\xrightarrow{\theta}\pi$,
$\kappa\xrightarrow{\tilde{\theta}}\pi$
and we have 
$\sigma_{u_{\pi}}=I_{\pi}$
where $\kappa\xrightarrow{u_{\pi}} \pi$.

Let $\pi^{\pr}\xrightarrow{f} \pi$ be a morphism of pre-Galois objects in $\kos{K}$.
For each $\kappa\xrightarrow{\theta^{\pr}}\pi^{\pr}$,
let us denote $f\theta^{\pr}:\kappa\xrightarrow{\theta^{\pr}} \pi^{\pr}\xrightarrow{f} \pi$.
Then the associated inner automorphisms
$\sigma^{\pr}_{\theta^{\pr}}:\pi^{\pr}\xrightarrow{\cong}\pi^{\pr}$,
$\sigma_{f\theta^{\pr}}:\pi\xrightarrow{\cong}\pi$
satisfy the relation
\begin{equation*}
  \vcenter{\hbox{
    \xymatrix@C=30pt{
      \pi^{\pr}
      \ar[d]_-{\sigma^{\pr}_{\theta^{\pr}}}^-{\cong}
      \ar[r]^-{f}
      &\pi
      \ar[d]^-{\sigma_{f\theta^{\pr}}}_-{\cong}
      \\
      \pi^{\pr}
      \ar[r]^-{f}
      &\pi
    }
  }}
  \qquad\quad
  \sigma_{f\theta^{\pr}}\circ f
  =
  f\circ \sigma^{\pr}_{\theta^{\pr}}
  .
\end{equation*}

\begin{definition}\label{def preGalobj GALOBJpre2cat}
  We define the $(2,1)$-category
  \begin{equation*}
    \mathbb{GALOBJ}^{\cat{pre}}(\kos{K})
  \end{equation*}
  of pre-Galois objects in $\kos{K}$ as follows.
  \begin{itemize}
    \item 
    A $0$-cell is a pre-Galois object $\pi$ in $\kos{K}$.

    \item 
    A $1$-cell $\pi^{\pr}\xrightarrow{f} \pi$
    is a morphism of pre-Galois objects in $\kos{K}$,
    i.e., a morphism of group objects in $\kos{E\!n\!s}(\kos{K})$.

    \item 
    A $2$-cell $\theta:f_1\Rightarrow f_2$
    between $1$-cells $f_1$, $f_2:\pi^{\pr}\to \pi$
    is a morphism $\kappa\xrightarrow{\theta} \pi$
    in $\cat{Ens}(\kos{K})$
    whose associated inner automorphism
    $\sigma_{\theta}:\pi\xrightarrow{\cong}\pi$
    satisfies the relation
    \begin{equation*}
      \vcenter{\hbox{
        \xymatrix@C=40pt{
          \pi^{\pr}
          \ar@/^1pc/[r]^-{f_1}
          \ar@/_1pc/[r]_-{f_2}
          \xtwocell[r]{}<>{<0>{\text{ }\text{ }\theta}}
          &\pi
        }
      }}
      \qquad\quad
      \vcenter{\hbox{
        \xymatrix@R=10pt@C=30pt{
          \text{ }
          &\pi
          \ar[dd]^-{\sigma_{\theta}}_-{\cong}
          \\
          \pi^{\pr}
          \ar[ur]^-{f_1}
          \ar[dr]_-{f_2}
          &\text{ }
          \\
          \text{ }
          &\pi
        }
      }}
      \qquad\quad
      \sigma_{\theta}\circ f_1=f_2.
    \end{equation*}
    This is equivalent to saying that
    the relation
    $\theta\star f_1=f_2\star \theta$ holds.
    \begin{equation*}
      \vcenter{\hbox{
        \xymatrix@C=30pt{
          \kappa\otimes \pi^{\pr}
          \ar[d]_-{\theta\otimes f_1}
          &\pi^{\pr}
          \ar[l]_-{\imath_{\pi^{\pr}}}^-{\cong}
          \ar[r]^-{\jmath_{\pi^{\pr}}}_-{\cong}
          &\pi^{\pr}\otimes \kappa
          \ar[d]^-{f_2\otimes \theta}
          \\
          \pi\otimes \pi
          \ar[r]^-{\pc_{\pi}}
          &\pi
          &\pi\otimes \pi
          \ar[l]_-{\pc_{\pi}}
        }
      }}
    \end{equation*}

    \item 
    Identity $1$-cell of the $0$-cell $\pi$ is the identity morphism $I_{\pi}:\pi\to \pi$.
    Composition of $1$-cells is the usual composition of morphisms of pre-Galois objects in $\kos{K}$.
    Identity $2$-cell of the $1$-cell $\pi^{\pr}\xrightarrow{f} \pi$
    is $\kappa\xrightarrow{u_{\pi}} \pi$.

    \item 
    Vertical composition of $2$-cells
    $\xymatrix@C=18pt{f_1\ar@2{->}[r]^-{\theta_1}&f_2\ar@2{->}[r]^-{\theta_2}&f_3}$
    between $1$-cells $f_1$, $f_2$, $f_3:\pi^{\pr}\to \pi$
    is the convolution product $\theta_2\star\theta_1:f_1\Rightarrow f_3$.
    \begin{equation*}
      \vcenter{\hbox{
        \xymatrix@C=50pt{
          \pi^{\pr}
          \ar@/^2.5pc/[r]^-{f_1}
          \ar[r]|-{f_2}
          \ar@/_2.5pc/[r]_-{f_3}
          \xtwocell[r]{}<>{<-3>{\text{ }\text{ }\theta_1}}
          \xtwocell[r]{}<>{<3>{\text{ }\text{ }\theta_2}}
          &\pi
        }
      }}
      \qquad\qquad
      \vcenter{\hbox{
        \xymatrix@R=20pt@C=50pt{
          \text{ }
          &\pi
          \ar[d]^-{\sigma_{\theta_1}}_-{\cong}
          \ar@/^2pc/@<3ex>[dd]^-{\sigma_{\theta_2\star \theta_1}}_-{\cong}
          \\
          \pi^{\pr}
          \ar[ur]^-{f_1}
          \ar[r]|-{f_2}
          \ar[dr]_-{f_3}
          &\pi
          \ar[d]^-{\sigma_{\theta_2}}_-{\cong}
          \\
          \text{ }
          &\pi
        }
      }}
    \end{equation*}

    \item 
    Every $2$-cell $\theta:f_1\Rightarrow f_2:\pi^{\pr}\to \pi$
    is invertible, whose inverse is
    $\kappa\xrightarrow{\theta}\pi\xrightarrow[\cong]{\varsigma_{\pi}}\pi$.

    \item 
    Horizontal composition of $2$-cells
    $\theta:f_1\Rightarrow f_2:\pi^{\pr}\to \pi$
    and
    $\theta^{\pr}:f^{\pr}_1\Rightarrow f^{\pr}_2:\pi^{\ppr}\to \pi^{\pr}$
    is defined as
    $\theta\centerdot\theta^{\pr}
    :=\theta\star f_1\theta^{\pr}
    =f_2\theta^{\pr}\star \theta$.
    \begin{equation*}
      \vcenter{\hbox{
        \xymatrix@C=40pt{
          \pi^{\ppr}
          \ar@/^1pc/[r]^-{f_1^{\pr}}
          \ar@/_1pc/[r]_-{f_2^{\pr}}
          \xtwocell[r]{}<>{<0>{\text{ }\text{ }\theta^{\pr}}}
          &\pi^{\pr}
          \ar@/^1pc/[r]^-{f_1}
          \ar@/_1pc/[r]_-{f_2}
          \xtwocell[r]{}<>{<0>{\text{ }\text{ }\theta}}
          &\pi
        }
      }}
      \qquad\quad
      \vcenter{\hbox{
        \xymatrix@C=50pt{
          \pi^{\ppr}
          \ar@/^1.2pc/[r]^-{f_1\circ f_1^{\pr}}
          \ar@/_1.2pc/[r]_-{f_2\circ f_2^{\pr}}
          \xtwocell[r]{}<>{<0>{\quad\theta\centerdot \theta^{\pr}}}
          &\pi
        }
      }}
    \end{equation*}
    \begin{equation*}
      \hspace{-0.5cm}
      \vcenter{\hbox{
        \xymatrix@R=12pt@C=18pt{
          \text{ }
          &\text{ }
          &\pi
          \ar[dddd]^-{\sigma_{\theta\centerdot\theta^{\pr}}}_-{\cong}
          \\
          \text{ }
          &\pi^{\pr}
          \ar[ur]^-{f_1}
          &\text{ }
          \\
          \pi^{\ppr}
          \ar[ur]^-{f^{\pr}_1}
          \ar[dr]_-{f^{\pr}_2}
          &\text{ }
          &\text{ }
          \\
          \text{ }
          &\pi^{\pr}
          \ar[dr]_-{f_2}
          &\text{ }
          \\
          \text{ }
          &\text{ }
          &\pi
        }
      }}
      =\text{ }\text{ }
      \vcenter{\hbox{
        \xymatrix@R=12pt@C=18pt{
          \text{ }
          &\text{ }
          &\pi
          \ar[dd]^-{\sigma_{f_1\theta^{\pr}}}_-{\cong}
          \\
          \text{ }
          &\pi^{\pr}
          \ar[ur]^-{f_1}
          \ar[dd]^-{\sigma^{\pr}_{\theta^{\pr}}}_-{\cong}
          &\text{ }
          \\
          \pi^{\ppr}
          \ar[ur]^-{f^{\pr}_1}
          \ar[dr]_-{f^{\pr}_2}
          &\text{ }
          &\pi
          \ar[dd]^-{\sigma_{\theta}}_-{\cong}
          \\
          \text{ }
          &\pi^{\pr}
          \ar[ur]|-{f_1}
          \ar[dr]_-{f_2}
          &\text{ }
          \\
          \text{ }
          &\text{ }
          &\pi
        }
      }}
      =\text{ }\text{ }
      \vcenter{\hbox{
        \xymatrix@R=12pt@C=18pt{
          \text{ }
          &\text{ }
          &\pi
          \ar[dd]^-{\sigma_{\theta}}_-{\cong}
          \\
          \text{ }
          &\pi^{\pr}
          \ar[ur]^-{f_1}
          \ar[dr]|-{f_2}
          \ar[dd]^-{\sigma^{\pr}_{\theta^{\pr}}}_-{\cong}
          &\text{ }
          \\
          \pi^{\ppr}
          \ar[ur]^-{f^{\pr}_1}
          \ar[dr]_-{f^{\pr}_2}
          &\text{ }
          &\pi
          \ar[dd]^-{\sigma_{f_2\theta^{\pr}}}_-{\cong}
          \\
          \text{ }
          &\pi^{\pr}
          \ar[dr]_-{f_2}
          &\text{ }
          \\
          \text{ }
          &\text{ }
          &\pi
        }
      }}
    \end{equation*}
    
    \item 
    Associators and left, right unitors are identities.
  \end{itemize}
\end{definition}

\subsection{Pre-Galois $\kos{K}$-category}
\label{subsec preGalKCat}

Let $\kos{K}=(\CK,\otimes,\kappa)$ be a Galois prekosmos.

\begin{definition}
  A \emph{pre-Galois $\kos{K}$-category}
  is a Galois $\kos{K}$-prekosmos $\mathfrak{T}=(\kos{T},\kos{t})$
  which admits a surjective pre-fiber functor
  $\varpi=(\omega,\what{\omega}):\mathfrak{K}\to \mathfrak{T}$.
\end{definition}

Let $\pi$ be a pre-Galois object in $\kos{K}$
and recall Lemma~\ref{lem preGalobj Rep(pi)}.
The connected Galois $\kos{K}$-prekosmos
$\mathfrak{Rep}(\pi)=(\kos{R\!e\!p}(\pi),\kos{t}_{\pi})$
of representations of $\pi$
is a pre-Galois $\kos{K}$-category,
as it admits a surjective pre-fiber functor
$\varpi_{\pi}:\mathfrak{K}\to\mathfrak{Rep}(\pi)$.

\begin{lemma} \label{lem preGalKCat prefib is surjective}
  Let $\mathfrak{T}=(\kos{T},\kos{t})$ be a pre-Galois $\kos{K}$-category.
  Then every pre-fiber functor $\mathfrak{K}\to\mathfrak{T}$ is surjective.
\end{lemma}
\begin{proof}
  From the definition of a pre-Galois $\kos{K}$-category $\mathfrak{T}$,
  there exists some
  surjective pre-fiber functor
  $\varpi=(\omega,\what{\omega}):\mathfrak{K}\to \mathfrak{T}$.
  Let $\varpi^{\pr}=(\omega^{\pr},\what{\omega}^{\pr}):\mathfrak{K}\to\mathfrak{T}$
  be an arbitrary pre-fiber functor.
  We claim that $\varpi$ and $\varpi^{\pr}$ are locally isomorphic
  and conclude that $\varpi^{\pr}$ is also surjective.
  Let us explain in detail.
  By Corollary~\ref{cor Galprekosmoi summary}
  the object $p:=\omega^{\pr*}\omega_!(\kappa)$ in $\cat{Ens}(\kos{K})$
  represents the presheaf of invertible Galois $\kos{K}$-transformations
  from $\varpi$ to $\varpi^{\pr}$,
  and the universal element is a comonoidal $\kos{K}$-tensor natural isomorphism
  \begin{equation} \label{eq preGalKCat prefib is surjective}
    \xi:
    \xymatrix{
      (\kos{p}^*\omega^*,\what{\kos{p}^*\omega^*})
      \ar@2{->}[r]^-{\cong}
      &(\kos{p}^*\omega^{\pr*},\what{\kos{p}^*\omega^{\pr*}})
      :(\kos{T},\kos{t}^*)\to (\kos{K}_p,\kos{p}^*).
    }
  \end{equation}
  Note that the functor
  $p\otimes:\CK\to \CK$
  is conservative and preserves reflexive coequalizers.
  This is because
  we have natural isomorphisms
  \begin{equation*}
    \xymatrix@C=40pt{
      p\otimes
      \ar@2{->}[r]^-{\bigl(\hatar{\omega^{\pr*}\omega_!}\bigr)^{-1}}_-{\cong}
      &\omega^{\pr*}\omega_!
      \ar@2{->}[r]^-{\varsigma^{\omega,\omega^{\pr}}}_-{\cong}
      &\omega^*\omega^{\pr}_!
      :\CK\to\CK
    }
  \end{equation*}
  where $\varsigma^{\omega,\omega^{\pr}}$
  is mentioned in Corollary~\ref{cor Galprekosmoi summary},
  and the functor
  $\omega^*\omega^{\pr}_!$
  is conservative and preserves reflexive coequalizers.
  Thus the functor
  $\kos{p}^*:\kos{K}\to\kos{K}_p$
  is conservative and preserves reflexive coequalizers.
  As we have the natural isomorphism
  (\ref{eq preGalKCat prefib is surjective})
  the composition
  $\kos{p}^*\omega^{\pr*}$
  is conservative and preserves reflexive coequalizers.
  We conclude that
  $\omega^{\pr*}$
  is conservative and preserves reflexive coequalizers.
  This shows that the pre-fiber functor
  $\varpi^{\pr}:\mathfrak{K}\to\mathfrak{T}$
  is surjective.
  This completes the proof of Lemma~\ref{lem preGalKCat prefib is surjective}.
\qed\end{proof}

\begin{theorem} \label{thm preGalKCat mainThm1}
  Let $\mathfrak{T}=(\kos{T},\kos{t})$ be a pre-Galois $\kos{K}$-category
  and let $\varpi:\mathfrak{K}\to\mathfrak{T}$ be a pre-fiber functor for $\mathfrak{T}$.
  \begin{enumerate}
    \item 
    We have a pre-Galois object
    $\omega^*\omega_!(\kappa)$ in $\kos{K}$
    which represents the presheaf of groups of
    invertible Galois $\kos{K}$-transformations
    from $\varpi$ to $\varpi$.
    \begin{equation*}
      \vcenter{\hbox{
        \xymatrix{
          \Hom_{\cat{Ens}(\kos{K})}\bigl(\slot,\omega^*\omega_!(\kappa)\bigr)
          \ar@2{->}[r]^-{\cong}
          &\underline{\Aut}_{\mathbb{GAL}^{\cat{pre}}_{\kos{K}}}(\varpi)
          :\cat{Ens}(\kos{K})^{\op}\to\cat{Grp}
        }
      }}
    \end{equation*}

    \item
    The given pre-fiber functor
    $\varpi:\mathfrak{K}\to \mathfrak{T}$
    is surjective,
    and factors through as an equivalence of Galois $\kos{K}$-prekosmoi
    $\widebreve{\varpi}:\mathfrak{Rep}(\omega^*\omega_!(\kappa))\xrightarrow{\simeq} \mathfrak{T}$.
    \begin{equation*}
      \vcenter{\hbox{
        \xymatrix@C=40pt{
          \mathfrak{T}
          &\mathfrak{Rep}(\pi)
          \ar[l]_-{\widebreve{\varpi}}^-{\simeq}
          \\
          \text{ }
          &\mathfrak{K}
          \ar@/^1pc/[ul]^-{\varpi}
          \ar[u]_-{\varpi_{\pi}}
        }
      }}
      \qquad\quad
      \pi=\omega^*\omega_!(\kappa)
    \end{equation*}
  \end{enumerate}
\end{theorem}
\begin{proof}
  In Proposition~\ref{prop RepGrpEns(K) comparison functor}
  we showed that we have a group object $\omega^*\omega_!(\kappa)$ in $\kos{E\!n\!s}(\kos{K})$
  which represents the presheaf of groups 
  $\underline{\Aut}_{\mathbb{GAL}^{\cat{pre}}_{\kos{K}}}(\varpi)$.
  The product, unit and antipode morphisms are explicitly given as follows.
  \begin{equation*}
    \begin{aligned}
      \pc_{\omega^*\omega_!(\kappa)}
      &:
      \xymatrix@C=40pt{
        \omega^*\omega_!(\kappa)
        \otimes
        \omega^*\omega_!(\kappa)
        \ar[r]^-{(\hatar{\omega^*\omega_!}_{\omega^*\omega_!(\kappa)})^{-1}}_-{\cong}
        &\omega^*\omega_!\omega^*\omega_!(\kappa)
        \ar[r]^-{\omega^*(\epsilon_{\omega_!(\kappa)})}
        &\omega^*\omega_!(\kappa)
      }
      \\
      u_{\omega^*\omega_!(\kappa)}
      &:
      \xymatrix{
        \kappa
        \ar[r]^-{\eta_{\kappa}}
        &\omega^*\omega_!(\kappa)
      }
      \\
      \varsigma_{\omega^*\omega_!(\kappa)}
      &:
      \xymatrix{
        \omega^*\omega_!(\kappa)
        \ar[r]^-{\varsigma^{\omega,\omega}_{\kappa}}_-{\cong}
        &\omega^*\omega_!(\kappa)
      }
    \end{aligned}
  \end{equation*}
  By Lemma~\ref{lem preGalKCat prefib is surjective},
  the given pre-fiber functor $\varpi:\mathfrak{K}\to \mathfrak{T}$
  is surjective.
  Therefore
  $\omega^*\omega_!:\kos{K}\to \kos{K}$
  is conservative and preserves reflexive coequalizers.
  As we have the natural isomorphism
  $\xymatrix@C=15pt{
    \hatar{\omega^*\omega_!}:
    \omega^*\omega_!
    \ar@2{->}[r]^-{\cong}
    &\omega^*\omega_!(\kappa)\otimes \slot
  }$
  the functor
  $\omega^*\omega_!(\kappa)\otimes\slot:\CK\to \CK$
  also preserves reflexive coequalizers.
  Thus $\omega^*\omega_!(\kappa)$ is a pre-Galois object in $\kos{K}$.
  This proves statement 1.

  Next we prove statement 2.
  Let us denote $\pi:=\omega^*\omega_!(\kappa)$.
  By Proposition~\ref{prop RepGrpEns(K) comparison functor},
  the strong $\kos{K}$-tensor functor
  $(\omega^*,\what{\omega}^*)
  :(\kos{T},\kos{t}^*)\to (\kos{K},\id_{\kos{K}})$
  factors through as a strong $\kos{K}$-tensor functor
  $(\widebreve{\omega}\!^*,\what{\widebreve{\omega}}^*)
  :(\kos{T},\kos{t}^*)\to (\kos{R\!e\!p}(\pi),\kos{t}_{\pi}^*)$.
  \begin{equation*}
    \vcenter{\hbox{
      \xymatrix@C=40pt{
        (\kos{T},\kos{t}^*)
        \ar[r]^-{(\widebreve{\omega}\!^*,\what{\widebreve{\omega}}^*)}
        \ar@/_1pc/[dr]_-{(\omega^*,\what{\omega}^*)}
        &(\kos{R\!e\!p}(\pi),\kos{t}_{\pi}^*)
        \ar[d]^-{(\omega_{\pi}^*,\what{\omega}_{\pi}^*)}
        \\
        \text{ }
        &(\kos{K},\id_{\kos{K}})
      }
    }}
  \end{equation*}  
  To conclude statement 2,
  it suffices to show that
  $\widebreve{\omega}\!^*$
  is an equivalence of categories.
  This follows from the crude monadicity theorem,
  which we explain in detail.
  We first construct the left adjoint
  $\widebreve{\omega}\!_!$
  of $\widebreve{\omega}\!^*$.
  The functor $\widebreve{\omega}\!_!$
  sends each object $(x,\gamma_x)$ in $\cat{Rep}(\pi)$
  to the following reflexive coequalizer in $\CT$.
  \begin{equation}\label{eq defbreveomega!}
    \vcenter{\hbox{
      \xymatrix@C=35pt{
        \omega_!\omega^*\omega_!(x)
        \ar@<0.5ex>[rr]^-{\epsilon_{\omega_!(x)}}
        \ar@<-0.5ex>[rr]_-{\omega_!\bigl(\gamma_x\text{ }\!\circ\text{ }\! \hatar{\omega^*\omega_!}_x\bigr)}
        &\text{ }
        &\omega_!(x)
        \ar@{->>}[r]^-{\widebreve{\textit{cq}}_{(x,\gamma_x)}}
        \ar@<-1ex>@/_2pc/[ll]_-{\omega_!(\eta_x)}
        &\widebreve{\omega}\!_!(x,\gamma_x)
      }
    }}
  \end{equation}
  The component 
  $\widebreve{\eta}\!_{(x,\gamma_x)}
  :(x,\gamma_x)\to \widebreve{\omega}\!^*\!\widebreve{\omega}\!_!(x,\gamma_x)$
  of the adjunction unit $\widebreve{\eta}$
  at each object $(x,\gamma_x)$ in $\cat{Rep}(\pi)$
  is given by
  \begin{equation*}
    \xymatrix@C=40pt{
      \widebreve{\eta}\!_{(x,\gamma_x)}:
      x
      \ar[r]^-{\eta_x}
      &\omega^*\omega_!(x)
      \ar@{->>}[r]^-{\omega^*(\widebreve{\textit{cq}}_{(x,\gamma_x)})}
      &\omega^*\!\widebreve{\omega}\!_!(x,\gamma_x)
      .
    }
  \end{equation*}
  We can check that
  $\widebreve{\eta}\!_{(x,\gamma_x)}$
  is indeed a morphism in $\cat{Rep}(\pi)$ as follows.
  Let us denote
  $\tilde{\gamma}_x=\gamma_x\circ \hatar{\omega^*\omega_!}_x:\omega^*\omega_!(x)\to x$
  and $X=\widebreve{\omega}\!_!(x,\gamma_x)$.
  Note that we have
  $\omega^*(\epsilon_X)
  =\gamma_X
  \circ
  \hatar{\omega^*\omega_!}_{\omega^*(X)}
  :\omega^*\omega_!\omega^*(X)
  \to \omega^*(X)$.
  \begin{equation*}
    \vcenter{\hbox{
      \xymatrix@C=15pt{
        \omega^*\omega_!(x)
        \ar[d]_-{\tilde{\gamma}_x}
        \ar@{=}[r]
        &\omega^*\omega_!(x)
        \ar[d]_-{\eta_{\omega^*\omega_!(x)}}
        \ar@{=}[r]
        &\omega^*\omega_!(x)
        \ar@{=}[dd]
        \ar@{=}[r]
        &\omega^*\omega_!(x)
        \ar[d]^-{\omega^*\omega_!(\eta_x)}
        \ar@{=}[r]
        &\omega^*\omega_!(x)
        \ar[dd]|-{\omega^*\omega_!(\widebreve{\eta}\!_{(x,\gamma_x)})}
        \\
        x
        \ar[dd]_-{\widebreve{\eta}\!_{(x,\gamma_x)}}
        \ar@/_1pc/[dr]_-{\eta_x}
        &\omega^*\omega_!\omega^*\omega_!(x)
        \ar[d]_-{\omega^*\omega_!(\tilde{\gamma}_x)}
        \ar@/^0.5pc/[dr]|-{\omega^*(\epsilon_{\omega_!(x)})}
        &\text{ }
        &\omega^*\omega_!\omega^*\omega_!(x)
        \ar@/_0.5pc/[dl]|-{\omega^*(\epsilon_{\omega_!(x)})}
        \ar@{->>}@/_0.5pc/[dr]|-{\omega^*\omega_!\omega^*(\widebreve{\textit{cq}}\!_{(x,\gamma_x)})}
        &\text{ }
        \\
        \text{ }
        &\omega^*\omega_!(x)
        \ar@/_0.5pc/@{->>}[dr]|-{\omega^*(\widebreve{\textit{cq}}\!_{(x,\gamma_x)})}
        &\omega^*\omega_!(x)
        \ar@{->>}[d]^-{\omega^*(\widebreve{\textit{cq}}\!_{(x,\gamma_x)})}
        &\text{ }
        &\omega^*\omega_!\omega^*(X)
        \ar[d]^-{\omega^*(\epsilon_X)}
        \\
        \omega^*(X)
        \ar@{=}[rr]
        &\text{ }
        &\omega^*(X)
        \ar@{=}[rr]
        &\text{ }
        &\omega^*(X)
      }
    }}
  \end{equation*}
  The component 
  $\widebreve{\epsilon}\!_Y:
  \widebreve{\omega}\!_!\!\widebreve{\omega}\!^*(Y)\to Y$
  of the adjunction counit $\widebreve{\epsilon}$
  at each object $Y$ in $\CT$ is the
  unique morphism in $\CT$ which satisfies the relation below.
  Note that we have
  $\omega^*(\epsilon_Y)=\gamma_{\omega^*(Y)}\circ \hatar{\omega^*\omega_!}_{\omega^*(Y)}:\omega^*\omega_!\omega^*(Y)\to \omega^*(Y)$.
  \begin{equation}\label{eq defbreveepsilon}
    \vcenter{\hbox{
      \xymatrix@C=40pt{
        \omega_!\omega^*\omega_!\omega^*(Y)
        \ar@<0.5ex>[r]^-{\epsilon_{\omega_!\omega^*(Y)}}
        \ar@<-0.5ex>[r]_-{\omega_!\omega^*(\epsilon_Y)}
        &\omega_!\omega^*(Y)
        \ar@{->>}[r]^-{\widebreve{\textit{cq}}_{\widebreve{\omega}\!^*(Y)}}
        \ar@/_1pc/[dr]_-{\epsilon_Y}
        &\widebreve{\omega}\!_!\!\widebreve{\omega}\!^*(Y)
        \ar@{.>}[d]^-{\widebreve{\epsilon}\!_Y}_-{\exists!}
        \\
        \text{ }
        &\text{ }
        &Y
      }
    }}
  \end{equation}
  We need to show that $\widebreve{\eta}$, $\widebreve{\epsilon}$
  satisfy the triangle identities.
  Let $(x,\gamma_x)$ be an object in $\cat{Rep}(\pi)$.
  We obtain one of the triangle identities
  by right-cancelling the epimorphism 
  $\widebreve{\textit{cq}}\!_{(x,\gamma_x)}$ 
  in the diagram below.
  \begin{equation*}
    \vcenter{\hbox{
      \xymatrix{
        \omega_!(x)
        \ar@{->>}[d]_-{\widebreve{\textit{cq}}\!_{(x,\gamma_x)}}
        \ar@{=}[r]
        &\omega_!(x)
        \ar[dd]_-{\omega_!(\widebreve{\eta}\!_{(x,\gamma_x)})}
        \ar@{=}[rr]
        &\text{ }
        &\omega_!(x)
        \ar@/_0.5pc/[dl]|-{\omega_!(\eta_x)}
        \ar@{=}[dd]
        \\
        \widebreve{\omega}\!_!(x,\gamma_x)
        \ar[dd]_-{\widebreve{\omega}_!(\widebreve{\eta}\!_{(x,\gamma_x)})}
        &\text{ }
        &\omega_!\omega^*\omega_!(x)
        \ar@{->>}@/_0.5pc/[dl]|-{\omega_!\omega^*(\widebreve{\textit{cq}}_{(x,\gamma_x)})}
        \ar@/_0.5pc/[dr]|-{\epsilon_{\omega_!(x)}}
        &\text{ }
        \\
        \text{ }
        &\omega_!\omega^*\!\widebreve{\omega}\!_!(x,\gamma_x)
        \ar@{->>}@/_0.5pc/[dl]|-{\widebreve{\textit{cq}}\!_{\widebreve{\omega}\!^*\!\widebreve{\omega}\!_!(x,\gamma_x)}}
        \ar[dd]^-{\epsilon_{\widebreve{\omega}\!_!(x,\gamma_x)}}
        &\text{ }
        &\omega_!(x)
        \ar@{->>}[dd]^-{\widebreve{\textit{cq}}_{(x,\gamma_x)}}
        \\
        \widebreve{\omega}\!_!\!\widebreve{\omega}\!^*\!\widebreve{\omega}\!_!(x,\gamma_x)
        \ar[d]_-{\widebreve{\epsilon}_{\widebreve{\omega}\!_!(x,\gamma_x)}}
        &\text{ }
        &\text{ }
        &\text{ }
        \\
        \widebreve{\omega}\!_!(x,\gamma_x)
        \ar@{=}[r]
        &\widebreve{\omega}\!_!(x,\gamma_x)
        \ar@{=}[rr]
        &\text{ }
        &\widebreve{\omega}\!_!(x,\gamma_x)
      }
    }}
  \end{equation*}  
  Let $Y$ be an object in $\CT$.
  We obtain the other triangle identity as follows.
  \begin{equation*}
    \vcenter{\hbox{
      \xymatrix{
        \omega^*(Y)
        \ar[dd]_-{\widebreve{\eta}\!_{\widebreve{\omega}\!^*(Y)}}
        \ar@{=}[rr]
        &\text{ }
        &\omega^*(Y)
        \ar@/_0.5pc/[dl]|-{\eta_{\omega^*(Y)}}
        \ar@{=}[ddd]
        \\
        \text{ }
        &\omega^*\omega_!\omega^*(Y)
        \ar@{->>}@/_0.5pc/[dl]|-{\omega^*(\widebreve{\textit{cq}}\!_{\widebreve{\omega}\!^*(Y)})}
        \ar@/_1pc/[ddr]|-{\omega^*(\epsilon_Y)}
        &\text{ }
        \\
        \omega^*\!\widebreve{\omega}\!_!\!\widebreve{\omega}\!^*(Y)
        \ar[d]_-{\omega^*(\widebreve{\epsilon}\!_Y)}
        &\text{ }
        &\text{ }
        \\
        \omega^*(Y)
        \ar@{=}[rr]
        &\text{ }
        &\omega^*(Y)
      }
    }}
  \end{equation*}
  This shows that $\widebreve{\omega}\!^*$ admits a left adjoint $\widebreve{\omega}\!_!$.
  Next we show that the adjunction unit, counit are natural isomorphisms.
  For each object $(x,\gamma_x)$ in $\cat{Rep}(\pi)$
  let us denote
  $\tilde{\gamma}_x=\gamma_x\circ \hatar{\omega^*\omega_!}_x:\omega^*\omega_!(x)\to x$
  and $\phi=\omega^*\omega_!:\kos{K}\to \kos{K}$.
  Consider the following diagram.
  \begin{equation*}
    \vcenter{\hbox{
      \xymatrix@R=35pt@C=60pt{
        \pi\otimes \pi\otimes x
        \ar@<0.5ex>[r]^-{\pc_{\pi}\otimes I_x}
        \ar@<-0.5ex>[r]_-{I_{\pi}\otimes \gamma_x}
        \ar[d]_-{(\hatar{\phi}\hatar{\phi})_x^{-1}}^-{\cong}
        &\pi\otimes x
        \ar[r]^-{\gamma_x}
        \ar[d]^-{(\hatar{\phi})^{-1}_x}_-{\cong}
        \ar@/_1.5pc/@<-1ex>[l]_-{u_{\pi}\otimes I_{\pi\otimes x}}
        &x
        \ar@{.>}[d]^-{\exists!\text{ }\widebreve{\eta}\!_{(x,\gamma_x)}}_-{\cong}
        \ar@/_1.5pc/@<-1ex>[l]_-{u_{\pi}\otimes I_x}
        \\
        \omega^*\omega_!\omega^*\omega_!(x)
        \ar@<0.5ex>[r]^-{\omega^*(\epsilon_{\omega_!(x)})}
        \ar@<-0.5ex>[r]_-{\omega^*\omega_!(\tilde{\gamma}_x)}
        &\omega^*\omega_!(x)
        \ar@{->>}[r]^-{\omega^*(\widebreve{\textit{cq}}\!_{(x,\gamma_x)})}
        &\omega^*\!\widebreve{\omega}\!_!(x,\gamma_x)
      }
    }}
  \end{equation*}
  The top horizontal morphisms form a split coequalizer diagram.
  The bottom horizontal morphisms also form a coequalizer diagram,
  as the functor $\omega^*$ preserves the reflexive coequalizer diagram
  (\ref{eq defbreveomega!}).
  One can check that $\widebreve{\eta}\!_{(x,\gamma_x)}$
  is the unique morphism which makes the above diagram commutative.
  Thus $\widebreve{\eta}\!_{(x,\gamma_x)}$ is an isomorphism.
  This shows that $\widebreve{\eta}$ is a natural isomorphism.
  Let $Y$ be an object in $\CT$.
  We claim that $\widebreve{\epsilon}_Y$ is an isomorphism.
  If we apply the functor $\omega^*$ to the diagram (\ref{eq defbreveepsilon})
  we obtain the following.
  \begin{equation*}
    \vcenter{\hbox{
      \xymatrix@C=60pt{
        \omega^*\omega_!\omega^*\omega_!\omega^*(Y)
        \ar@<0.5ex>[r]^-{\omega^*(\epsilon_{\omega_!\omega^*(Y)})}
        \ar@<-0.5ex>[r]_-{\omega^*\omega_!\omega^*(\epsilon_Y)}
        \ar@{=}[d]
        &\omega^*\omega_!\omega^*(Y)
        \ar@{->>}[r]^-{\omega^*(\widebreve{\textit{cq}}_{\widebreve{\omega}\!^*(Y)})}
        \ar@{=}[d]
        &\omega^*\!\widebreve{\omega}\!_!\!\widebreve{\omega}\!^*(Y)
        \ar[d]^-{\omega^*(\widebreve{\epsilon}\!_Y)}_-{\cong}
        \\
        \omega^*\omega_!\omega^*\omega_!\omega^*(Y)
        \ar@<0.5ex>[r]^-{\omega^*(\epsilon_{\omega_!\omega^*(Y)})}
        \ar@<-0.5ex>[r]_-{\omega^*\omega_!\omega^*(\epsilon_Y)}
        &\omega^*\omega_!\omega^*(Y)
        \ar[r]^-{\omega^*(\epsilon_Y)}
        \ar@/^2pc/@<1ex>[l]|-{\eta_{\omega^*\omega_!\omega^*(Y)}}
        &\omega^*(Y)
        \ar@/^2pc/@<1ex>[l]|-{\eta_{\omega^*(Y)}}
      }
    }}
  \end{equation*}
  The top horizontal morphisms form a coequalizer diagram
  as the functor $\omega^*$ preserves reflexive coequalizers.
  The bottom horizontal morphisms form a split coequalizer diagram.
  Thus we obtain that
  $\omega^*(\widebreve{\epsilon}\!_Y)$
  is an isomoprhism.
  As the functor $\omega^*$
  is conservative,
  we conclude that $\widebreve{\epsilon}\!_Y$ is an isomorphism
  for every object $Y$ in $\CT$.
  This shows that $\widebreve{\epsilon}$ is also a natural isomorphism,
  and that $\widebreve{\omega}\!_!\dashv \widebreve{\omega}\!^*$
  is an adjoint equivalence of categories.
  This completes the proof of Theorem~\ref{thm preGalKCat mainThm1}.
\qed\end{proof}

\begin{corollary}
  Every pre-Galois $\kos{K}$-category
  $\mathfrak{T}=(\kos{T},\kos{t})$
  is connected as a Galois $\kos{K}$-prekosmos.
\end{corollary}
\begin{proof}
  This is an immediate consequence of 
  Lemma~\ref{lem preGalobj Rep(pi)}
  and 
  Theorem~\ref{thm preGalKCat mainThm1}.
\qed\end{proof}

\subsection{Pre-fundamental functors}
\label{subsec Galprefundamentalfunct}

Let $\kos{K}=(\CK,\otimes,\kappa)$ be a Galois prekosmos.

In this subsection, we introduce morphisms between pre-Galois $\kos{K}$-categories
and define 
the $(2,1)$-category $\mathbb{GALCAT}^{\cat{pre}}_*(\kos{K})$
of pointed pre-Galois $\kos{K}$-categories.
The main goal is to establish an equivalence
of $2$-categories
between 
$\mathbb{GALOBJ}^{\cat{pre}}(\kos{K})$
and 
$\mathbb{GALCAT}^{\cat{pre}}_*(\kos{K})$.

\begin{definition}
  Let $\mathfrak{T}=(\kos{T},\kos{t})$, $\mathfrak{T}^{\pr}=(\kos{T}^{\pr},\kos{t}^{\pr})$
  be pre-Galois $\kos{K}$-categories.
  A \emph{pre-fundamental functor} from $\mathfrak{T}^{\pr}$ to $\mathfrak{T}$
  is a Galois $\kos{K}$-morphism
  $\mathfrak{f}:\mathfrak{T}^{\pr}\to\mathfrak{T}$
  which has the following property:
  there exists some pre-fiber functor
  $\varpi^{\pr}:\mathfrak{K}\to\mathfrak{T}^{\pr}$
  such that the composition 
  $\mathfrak{f}\varpi^{\pr}:\mathfrak{K}\xrightarrow{\varpi^{\pr}}\mathfrak{T}^{\pr}\xrightarrow{\mathfrak{f}}\mathfrak{T}$
  is also a pre-fiber functor.
\end{definition}

\begin{lemma} \label{lem Galprefundamentalfunct Hom=Isom}
  Let $\mathfrak{T}=(\kos{T},\kos{t})$, $\mathfrak{T}^{\pr}=(\kos{T}^{\pr},\kos{t}^{\pr})$
  be pre-Galois $\kos{K}$-categories
  and let
  $\mathfrak{f}_1$, $\mathfrak{f}_2:\mathfrak{T}^{\pr}\to\mathfrak{T}$
  be pre-fundamental functors.
  Then every Galois $\kos{K}$-transformation
  $\vartheta:\mathfrak{f}_1\Rightarrow\mathfrak{f}_2$
  is invertible.
\end{lemma}
\begin{proof}
  Choose a pre-fiber functor $\varpi^{\pr}:\mathfrak{K}\to\mathfrak{T}^{\pr}$
  such that the composition
  $\mathfrak{f}_1\varpi^{\pr}:\mathfrak{K}\to\mathfrak{T}$
  is also a pre-fiber functor.
  Consider an arbitrary Galois $\kos{K}$-transformation
  $\vartheta:\mathfrak{f}_1\Rightarrow\mathfrak{f}_2:\mathfrak{T}^{\pr}\to\mathfrak{T}$
  which is a comonoidal $\kos{K}$-tensor natural transformation
  \begin{equation*}
    \vartheta:(\kos{f}_1^*,\what{\kos{f}_1^*})\Rightarrow (\kos{f}_2^*,\what{\kos{f}_2^*})
    :(\kos{T},\kos{t}^*)\to (\kos{T}^{\pr},\kos{t}^{\pr*})
    .
  \end{equation*}
  By Lemma~\ref{lem preGalKCat prefib is surjective},
  the pre-fiber functor $\varpi^{\pr}$ is surjective.
  In particular, the right adjoint $\omega^{\pr*}$ is conservative.
  Therefore to conclude that $\vartheta$ is invertible,
  it suffices to show that
  \begin{equation*}
    \omega^{\pr*}\vartheta:
    (\omega^{\pr*}\kos{f}_1^*,\what{\omega^{\pr*}\kos{f}_1^*})\Rightarrow (\omega^{\pr*}\kos{f}_2^*,\what{\omega^{\pr*}\kos{f}_2^*})
    :(\kos{T},\kos{t}^*)\to (\kos{K},\id_{\kos{K}})
  \end{equation*}
  is invertible.
  This follows from Proposition~\ref{prop RSKtensoradj HomRep,Hom=Isom},
  since 
  $(\kos{f}_1\omega^{\pr},\what{\kos{f}_1\omega^{\pr}}):(\kos{K},\id_{\kos{K}})\to (\kos{T},\kos{t}^*)$
  is a reflective right-strong $\kos{K}$-tensor adjunction
  satisfying the projection formula,
  and
  $(\omega^{\pr*}\kos{f}_2^*,\what{\omega^{\pr*}\kos{f}_2^*}):(\kos{T},\kos{t}^*)\to (\kos{K},\id_{\kos{K}})$
  is a strong $\kos{K}$-tensor functor.
  \begin{equation*}
    \vcenter{\hbox{
      \xymatrix{
        (\kos{T},\kos{t}^*)
        \ar@/^1pc/[d]^-{(\omega^{\pr*}\kos{f}_1^*\text{ }\!,\text{ }\!\what{\omega^{\pr*}\kos{f}_1^*})}
        \\
        (\kos{K},\id_{\kos{K}})
        \ar@/^1pc/[u]^-{(\kos{f}_{1!}\omega_!^{\pr}\text{ }\!,\text{ }\!\what{\kos{f}_{1!}\omega_!^{\pr}})}
      }
    }}
    \qquad\qquad\quad
    \vcenter{\hbox{
      \xymatrix{
        (\kos{T},\kos{t}^*)
        \ar[d]^-{(\omega^{\pr*}\kos{f}_2^*\text{ }\!,\text{ }\!\what{\omega^{\pr*}\kos{f}_2^*})}
        \\
        (\kos{K},\id_{\kos{K}})
      }
    }}
  \end{equation*}
  This completes the proof of Lemma~\ref{lem Galprefundamentalfunct Hom=Isom}.
\qed\end{proof}

\begin{definition}
  We define the $(2,1)$-category 
  \begin{equation*}
    \mathbb{GALCAT}^{\cat{pre}}_*(\kos{K})
  \end{equation*}
  of pre-Galois $\kos{K}$-categories pointed with pre-fiber functors
  as follows.
  \begin{itemize}
    \item 
    A $0$-cell is a pair $(\mathfrak{T},\varpi)$
    of a pre-Galois $\kos{K}$-category $\mathfrak{T}$
    and a pre-fiber functor $\varpi:\mathfrak{K}\to\mathfrak{T}$.

    \item 
    A $1$-cell $(\mathfrak{f},\alpha):(\mathfrak{T}^{\pr},\varpi^{\pr})\to (\mathfrak{T},\varpi)$
    is a pair of a pre-fundamental functor $\mathfrak{f}:\mathfrak{T}^{\pr}\to\mathfrak{T}$
    and an invertible Galois $\kos{K}$-transformation
    $\xymatrix@C=15pt{\alpha:\varpi\ar@2{->}[r]^-{\cong}&\mathfrak{f}\varpi^{\pr}.}$
    \begin{equation*}
      \vcenter{\hbox{
        \xymatrix@C=15pt{
          \mathfrak{T}
          &\text{ }
          &\mathfrak{T}^{\pr}
          \ar[ll]_-{\mathfrak{f}}
          \\
          \text{ }
          &\mathfrak{K}
          \ar[ur]_-{\varpi^{\pr}}
          \ar[ul]^-{\varpi}
          \xtwocell[u]{}<>{<0>{\text{ }\text{ }\alpha}}
        }
      }}
      \qquad\quad
      \alpha:
      \xymatrix@C=18pt{
        \omega^*
        \ar@2{->}[r]^-{\cong}
        &\omega^{\pr*}\kos{f}^*
        :\kos{T}\to\kos{K}
      }
      \qquad
      \vcenter{\hbox{
        \xymatrix{
          \kos{K}
          \ar@{=}[d]
          &\kos{T}
          \ar[l]_-{\omega^*}
          \ar[d]^-{\kos{f}^*}
          \\
          \kos{K}
          \xtwocell[ur]{}<>{<0>{\text{ }\text{ }\alpha}}
          &\kos{T}^{\pr}
          \ar[l]^-{\omega^{\pr*}}
        }
      }}
    \end{equation*}
    We denote the mate associated to $\alpha$ as
    \begin{equation*}
      \vcenter{\hbox{
        \xymatrix{
          \kos{K}
          \ar@{=}[d]
          \ar[r]^-{\omega_!}
          &\kos{T}
          \ar[d]^-{\kos{f}^*}
          \\
          \kos{K}
          \ar[r]_-{\omega^{\pr}_!}
          &\kos{T}^{\pr}
          \xtwocell[ul]{}<>{<0>{\alpha_!\text{ }\text{ }}}
        }
      }}
      \qquad
      \alpha_!:
      \xymatrix{
        \omega^{\pr}_!
        \ar@2{->}[r]^-{\omega^{\pr}_!\eta}
        &\omega^{\pr}_!\omega^*\omega_!
        \ar@2{->}[r]^-{\omega^{\pr}_!\alpha\omega_!}_-{\cong}
        &\omega^{\pr}_!\omega^{\pr*}\kos{f}^*\omega_!
        \ar@2{->}[r]^-{\epsilon^{\pr}\kos{f}^*\omega_!}
        &\kos{f}^*\omega_!
        .
      }
    \end{equation*}

    \item
    A $2$-cell $\vartheta:(\mathfrak{f}_1,\alpha_1)\Rightarrow (\mathfrak{f}_2,\alpha_2)$
    between $1$-cells
    $(\mathfrak{f}_1,\alpha_1)$, $(\mathfrak{f}_2,\alpha_2):(\mathfrak{T}^{\pr},\varpi^{\pr})\to (\mathfrak{T},\varpi)$
    is a Galois $\kos{K}$-transformation
    $\vartheta:\mathfrak{f}_1\Rightarrow\mathfrak{f}_2$
    between pre-fundamental functors.
    It is automatically an invertible Galois $\kos{K}$-transformation
    by Lemma~\ref{lem Galprefundamentalfunct Hom=Isom}.
    \begin{equation*}
      \vcenter{\hbox{
        \xymatrix@C=35pt{
          (\mathfrak{T},\varpi)
          \xtwocell[r]{}<>{<0>{\text{ }\text{ }\vartheta}}
          &(\mathfrak{T}^{\pr},\varpi^{\pr})
          \ar@/_1pc/@<-0.5ex>[l]_-{(\mathfrak{f}_1,\alpha_1)}
          \ar@/^1pc/@<0.5ex>[l]^-{(\mathfrak{f}_2,\alpha_2)}
        }
      }}
      \qquad\quad
      \vartheta:
      \xymatrix@C=15pt{
        \mathfrak{f}_1
        \ar@2{->}[r]^-{\cong}
        &\mathfrak{f}_2
      }
    \end{equation*}

    \item
    Identity $1$-cell of a $0$-cell $(\mathfrak{T},\varpi)$
    is a pair of an identity fundamental functor $I_{\mathfrak{T}}:\mathfrak{T}\to\mathfrak{T}$
    and the identity Galois $\kos{K}$-transformation
    $\varpi=I_{\mathfrak{T}}\varpi$.
    Identity $2$-cell of a $1$-cell
    $(\mathfrak{f},\alpha):(\mathfrak{T}^{\pr},\varpi^{\pr})\to (\mathfrak{T},\varpi)$
    is the identity Galois $\kos{K}$-transformation
    $I_{\mathfrak{f}}:\mathfrak{f}=\mathfrak{f}$.

    \item 
    Composition of $1$-cells is given as follows.
    \begin{equation*}
      \vcenter{\hbox{
        \xymatrix@C=20pt{
          (\mathfrak{T},\varpi)
          &(\mathfrak{T}^{\pr},\varpi^{\pr})
          \ar[l]_-{(\mathfrak{f},\alpha)}
          &(\mathfrak{T}^{\ppr},\varpi^{\ppr})
          \ar[l]_-{(\mathfrak{f}^{\pr},\alpha^{\pr})}
          \ar@/^1pc/@<0.5ex>[ll]^-{(\mathfrak{f},\alpha)\circ(\mathfrak{f}^{\pr},\alpha^{\pr})=\bigl(\mathfrak{f}\mathfrak{f}^{\pr},(\mathfrak{f}\alpha^{\pr})\circ\alpha\bigr)}
        }
      }}
      \quad
      \vcenter{\hbox{
        \xymatrix@C=15pt{
          (\mathfrak{f}\alpha^{\pr})\circ\alpha:
          \varpi
          \ar@2{->}[r]^-{\alpha}_-{\cong}
          &\mathfrak{f}\varpi^{\pr}
          \ar@2{->}[r]^-{\mathfrak{f}\alpha^{\pr}}_-{\cong}
          &\mathfrak{f}\mathfrak{f}^{\pr}\varpi^{\ppr}
        }
      }}
    \end{equation*}

    \item
    Vertical composition of $2$-cells
    is given by vertical composition as Galois $\kos{K}$-transformations.
    \begin{equation*}
      \vcenter{\hbox{
        \xymatrix@C=50pt{
          (\mathfrak{T},\varpi)
          \xtwocell[r]{}<>{<-3>{\text{ }\text{ }\vartheta_1}}
          \xtwocell[r]{}<>{<3>{\text{ }\text{ }\vartheta_2}}
          &(\mathfrak{T}^{\pr},\varpi^{\pr})
          \ar@/_2pc/@<-0.5ex>[l]_-{(\mathfrak{f}_1,\alpha_1)}
          \ar[l]|-{(\mathfrak{f}_2,\alpha_2)}
          \ar@/^2pc/@<0.5ex>[l]^-{(\mathfrak{f}_3,\alpha_3)}
        }
      }}
      \qquad\quad
      \begin{aligned}
        \vartheta_2\circ\vartheta_1
        &:
        \xymatrix@C=18pt{
          \mathfrak{f}_1
          \ar@2{->}[r]^-{\vartheta_1}_-{\cong}
          &\mathfrak{f}_2
          \ar@2{->}[r]^-{\vartheta_2}_-{\cong}
          &\mathfrak{f}_3
        }
      \end{aligned}
    \end{equation*}

    \item 
    Horizontal composition of $2$-cells 
    is given by 
    horizontal composition as Galois $\kos{K}$-transformations.
    \begin{equation*}
      \vcenter{\hbox{
        \xymatrix@C=35pt{
          (\mathfrak{T},\varpi)
          \xtwocell[r]{}<>{<0>{\text{ }\text{ }\vartheta}}
          &(\mathfrak{T}^{\pr},\varpi^{\pr})
          \ar@/_1pc/@<-0.5ex>[l]_-{(\mathfrak{f}_1,\alpha_1)}
          \ar@/^1pc/@<0.5ex>[l]^-{(\mathfrak{f}_2,\alpha_2)}
          \xtwocell[r]{}<>{<0>{\text{ }\text{ }\text{ }\vartheta^{\pr}}}
          &(\mathfrak{T}^{\ppr},\varpi^{\ppr})
          \ar@/_1pc/@<-0.5ex>[l]_-{(\mathfrak{f}^{\pr}_1,\alpha^{\pr}_1)}
          \ar@/^1pc/@<0.5ex>[l]^-{(\mathfrak{f}^{\pr}_2,\alpha^{\pr}_2)}
        }
      }}
    \end{equation*}
    \begin{equation*}
      \vcenter{\hbox{
        \xymatrix@C=50pt{
          (\mathfrak{T},\varpi)
          \xtwocell[r]{}<>{<0>{\quad\vartheta\vartheta^{\pr}}}
          &(\mathfrak{T}^{\ppr},\varpi^{\ppr})
          \ar@/_1pc/@<-0.5ex>[l]_-{(\mathfrak{f}_1,\alpha_1)\circ(\mathfrak{f}^{\pr}_1,\alpha^{\pr}_1)}
          \ar@/^1pc/@<0.5ex>[l]^-{(\mathfrak{f}_2,\alpha_2)\circ(\mathfrak{f}^{\pr}_2,\alpha^{\pr}_2)}
        }
      }}
      \qquad\quad
      \begin{aligned}
        \vartheta\vartheta^{\pr}
        &:
        \mathfrak{f}_1\mathfrak{f}_1^{\pr}
        \Rightarrow
        \mathfrak{f}_2\mathfrak{f}_2^{\pr}
        :\mathfrak{T}^{\ppr}\to\mathfrak{T}
        \\
        \vartheta^{\pr}\vartheta
        &:
        \kos{f}_1^{\pr*}\kos{f}_1^*
        \Rightarrow
        \kos{f}_2^{\pr*}\kos{f}_2^*
        :\kos{T}\to\kos{T}^{\ppr}
      \end{aligned}
    \end{equation*}
  \end{itemize}
\end{definition}

\begin{proposition}\label{prop inducedfundamentalfunctor}
  Let $\kos{K}$ be a Galois prekosmos $\kos{K}$.
  We have a $2$-functor
  \begin{equation*}
    \vcenter{\hbox{
      \xymatrix@C=30pt{
        \mathbb{GALOBJ}^{\cat{pre}}(\kos{K})
        \ar[r]
        &\mathbb{GALCAT}^{\cat{pre}}_*(\kos{K}).
      }
    }}
  \end{equation*}
\end{proposition}
\begin{proof}
  First, we describe how the $2$-functor
  $\mathbb{GALOBJ}^{\cat{pre}}(\kos{K})
  \to \mathbb{GALCAT}^{\cat{pre}}_*(\kos{K})$
  sends $0$-cells to $0$-cells.
  For each pre-Galois object $\pi$ in $\kos{K}$,
  we have the pointed pre-Galois $\kos{K}$-category
  $(\mathfrak{Rep}(\pi),\varpi_{\pi})$
  which we described in Lemma~\ref{lem preGalobj Rep(pi)}.
  
  Next, we describe how the $2$-functor
  $\mathbb{GALOBJ}^{\cat{pre}}(\kos{K})
  \to \mathbb{GALCAT}^{\cat{pre}}_*(\kos{K})$
  sends
  $1$-cells to $1$-cells.
  Let $\pi^{\pr}\xrightarrow{f}\pi$
  be a morphism of pre-Galois objects in $\kos{K}$.
  We claim that there exists a Galois $\kos{K}$-morphism
  $\mathfrak{f}=(\kos{f},\what{\kos{f}}):
  \mathfrak{Rep}(\pi^{\pr})\to \mathfrak{Rep}(\pi)$
  and an invertible Galois $\kos{K}$-transformation
  $\xymatrix@C=15pt{
    I:
    \varpi_{\pi}
    \ar@2{->}[r]^-{\cong}
    &\mathfrak{f}\varpi_{\pi^{\pr}}
  }$.
  \begin{equation*}
    \vcenter{\hbox{
      \xymatrix@C=15pt{
        \mathfrak{Rep}(\pi)
        &\text{ }
        &\mathfrak{Rep}(\pi^{\pr})
        \ar[ll]_-{\mathfrak{f}}
        \\
        \text{ }
        &\mathfrak{K}
        \ar[ul]^-{\varpi_{\pi}}
        \ar[ur]_-{\varpi_{\pi^{\pr}}}
        \xtwocell[u]{}<>{<0>{I}}
      }
    }}
  \end{equation*}
  Let us explain the underlying Galois morphism
  $\kos{f}:\kos{R\!e\!p}(\pi^{\pr})\to \kos{R\!e\!p}(\pi)$.
  The right adjoint $\kos{f}^*:\kos{R\!e\!p}(\pi)\to \kos{R\!e\!p}(\pi^{\pr})$
  sends each object $X=(x,\gamma_x)$ in $\cat{Rep}(\pi)$ to
  \begin{equation*}
    \kos{f}^*(X)=
    (x,\xymatrix@C=15pt{
      \pi^{\pr}\otimes x
      \ar[r]^-{f\otimes I_x}
      &\pi\otimes x
      \ar[r]^-{\gamma_x}
      &x
    }\!).
  \end{equation*}
  The comonoidal coherence isomorphisms of $\kos{f}^*$ are given by identity morphisms.
  We have
  $\kos{f}^*_{\unit_{\pi}}:\kos{f}^*(\unit\!_{\pi})=\unit\!_{\pi^{\pr}}$
  and for another representation $Y=(y,\gamma_y)$ of $\pi$
  we have
  $\kos{f}^*_{X,Y}:
  \kos{f}^*(X\tensor\!_{\pi}Y)
  =\kos{f}^*(X)\tensor\!_{\pi^{\pr}} \kos{f}^*(Y)$.
  The left adjoint $\kos{f}_!$ sends
  each object $Z=(z,\gamma^{\pr}_z)$ in $\cat{Rep}(\pi^{\pr})$ to
  $\kos{f}_!(Z)=(\pi\otimes_{\pi^{\pr}}z,\gamma_{\pi\otimes_{\pi^{\pr}}z})$
  whose underlying object in $\CK$ is the reflexive coequalizer in $\CK$
  \begin{equation} \label{eq inducedfundamentalfunctor f!def}
    \vcenter{\hbox{
      \xymatrix@C=30pt{
        \pi\otimes \pi^{\pr}\otimes z
        \ar@<0.5ex>[rr]^-{I_{\pi}\otimes \gamma^{\pr}_z}
        \ar@<-0.5ex>[rr]_-{(\pc_{\pi}\otimes I_z)\circ (I_{\pi}\otimes f\otimes I_z)}
        &\text{ }
        &\pi\otimes z
        \ar@/_2pc/@<-1ex>[ll]|-{I_{\pi}\otimes u_{\pi^{\pr}}\otimes I_z}
        \ar@{->>}[r]^-{\textit{cq}^f_Z}
        &\pi\otimes_{\pi^{\pr}}z
      }
    }}
  \end{equation}
  and the left $\pi$-action
  $\gamma_{\pi\otimes_{\pi^{\pr}}z}$
  is the unique morphism in $\CK$ satisfying the relation below.
  \begin{equation*}
    \vcenter{\hbox{
      \xymatrix@C=50pt{
        \pi\otimes \pi\otimes z
        \ar@{->>}[r]^-{I_{\pi}\otimes \textit{cq}^f_Z}
        \ar[d]_-{\pc_{\pi}\otimes I_z}
        &\pi\otimes (\pi\otimes_{\pi^{\pr}}z)
        \ar@{.>}[d]^-{\gamma_{\pi\otimes_{\pi^{\pr}}z}}_-{\exists!}
        \\
        \pi\otimes z
        \ar@{->>}[r]^-{\textit{cq}^f_Z}
        &\pi\otimes_{\pi^{\pr}}z
      }
    }}
  \end{equation*}
  We can check that the morphism
  $\gamma_{\pi\otimes_{\pi^{\pr}}z}$
  is well-defined as follows.
  Note that the functor $\pi\otimes\slot:\CK\to \CK$
  preserves the defining reflexive coequalizer diagram
  (\ref{eq inducedfundamentalfunctor f!def}) of $\pi\otimes_{\pi^{\pr}}z$.
  \begin{equation*}
    \vcenter{\hbox{
      \xymatrix@C=60pt{
        \pi\pi\pi^{\pr}z
        \ar[d]_-{I_{\pi\pi}fI_z}
        \ar@{=}[r]
        &\pi\pi\pi^{\pr}z
        \ar[d]^-{\pc_{\pi}I_{\pi^{\pr}z}}
        \ar@{=}[r]
        &\pi\pi\pi^{\pr}z
        \ar[dd]^-{I_{\pi\pi}\gamma^{\pr}_z}
        \\
        \pi\pi\pi z
        \ar[d]_-{I_{\pi}\pc_{\pi}I_z}
        \ar@/^0.5pc/[dr]|-{\pc_{\pi}I_{\pi z}}
        &\pi\pi^{\pr}z
        \ar[d]^-{I_{\pi}fI_z}
        \ar@/^1pc/[ddr]|-{I_{\pi}\gamma^{\pr}_z}
        &\text{ }
        \\
        \pi\pi z
        \ar[d]_-{\pc_{\pi}I_z}
        &\pi\pi z
        \ar@/^0.5pc/[dl]|-{\pc_{\pi}I_z}
        &\pi\pi z
        \ar[d]^-{\pc_{\pi}I_z}
        \\
        \pi z
        \ar@{->>}[d]_-{\textit{cq}^f_Z}
        &\text{ }
        &\pi z
        \ar@{->>}[d]^-{\textit{cq}^f_Z}
        \\
        \pi\otimes_{\pi^{\pr}}z
        \ar@{=}[rr]
        &\text{ }
        &\pi\otimes_{\pi^{\pr}}z
      }
    }}
  \end{equation*}
  We obtain that
  $\gamma_{\pi\otimes_{\pi^{\pr}}z}$
  satisfies the left $\pi$-action relations,
  by right-cancelling the epimorphism
  $I_{\pi\otimes \pi}\otimes \textit{cq}^f_Z$
  in the following diagram
  \begin{equation*}
    \vcenter{\hbox{
      \xymatrix@C=10pt{
        \pi\otimes \pi\otimes \pi\otimes z
        \ar@{->>}[d]_-{I_{\pi\otimes \pi}\otimes \textit{cq}^f_Z}
        \ar@{=}[rr]
        &\text{ }
        &\pi\otimes \pi\otimes \pi\otimes z
        \ar@/_0.5pc/[dl]|-{I_{\pi}\otimes \pc_{\pi}\otimes I_x}
        \ar[d]^-{\pc_{\pi}\otimes I_{\pi\otimes z}}
        \ar@{=}[r]
        &\pi\otimes \pi\otimes \pi\otimes z
        \ar@{->>}[d]^-{I_{\pi\otimes \pi}\otimes \textit{cq}^f_Z}
        \\
        \pi\otimes \pi\otimes (\pi\otimes_{\pi^{\pr}}z)
        \ar[d]_-{I_{\pi}\otimes \gamma_{\pi\otimes_{\pi^{\pr}}z}}
        &\pi\otimes \pi\otimes z
        \ar@{->>}@/_0.5pc/[dl]|-{I_{\pi}\otimes \textit{cq}^f_Z}
        \ar[d]|-{\pc_{\pi}\otimes I_z}
        &\pi\otimes \pi\otimes z
        \ar@/^0.5pc/[dl]|-{\pc_{\pi}\otimes I_z}
        \ar@{->>}@/_0.5pc/[dr]|-{I_{\pi}\otimes \textit{cq}^f_Z}
        &\pi\otimes \pi\otimes (\pi\otimes_{\pi^{\pr}}z)
        \ar[d]^-{\pc_{\pi}\otimes I_{\pi\otimes_{\pi^{\pr}}z}}
        \\
        \pi\otimes (\pi\otimes_{\pi^{\pr}}z)
        \ar[d]_-{\gamma_{\pi\otimes_{\pi^{\pr}}z}}
        &\pi\otimes z
        \ar@/^0.5pc/@{->>}[dl]|-{\textit{cq}^f_Z}
        &\text{ }
        &\pi\otimes (\pi\otimes_{\pi^{\pr}}z)
        \ar[d]^-{\gamma_{\pi\otimes_{\pi^{\pr}}z}}
        \\
        \pi\otimes_{\pi^{\pr}}z
        \ar@{=}[rrr]
        &\text{ }
        &\text{ }
        &\pi\otimes_{\pi^{\pr}}z
      }
    }}
  \end{equation*}
  and right-cancelling the epimorphism
  $\textit{cq}^f_Z$
  in the diagram below.
  \begin{equation*}
    \vcenter{\hbox{
      \xymatrix{
        \pi\otimes z
        \ar@{->>}[d]_-{\textit{cq}^f_Z}
        \ar@{=}[rr]
        &\text{ }
        &\pi\otimes z
        \ar@/_0.5pc/[dl]|-{u_{\pi}\otimes I_{\pi\otimes z}}
        \ar@{=}[dd]
        \\
        \pi\otimes_{\pi^{\pr}}z
        \ar[d]_-{u_{\pi}\otimes I_{\pi\otimes_{\pi^{\pr}}z}}
        &\pi\otimes \pi\otimes z
        \ar@{->>}@/_0.5pc/[dl]|-{I_{\pi}\otimes \textit{cq}^f_Z}
        \ar@/^0.5pc/[dr]|-{\pc_{\pi}\otimes I_z}
        &\text{ }
        \\
        \pi\otimes (\pi\otimes_{\pi^{\pr}}z)
        \ar[d]_-{\gamma_{\pi\otimes_{\pi^{\pr}}z}}
        &\text{ }
        &\pi\otimes z
        \ar@{->>}[d]^-{\textit{cq}^f_Z}
        \\
        \pi\otimes_{\pi^{\pr}}z
        \ar@{=}[rr]
        &\text{ }
        &\pi\otimes_{\pi^{\pr}}z
      }
    }}
  \end{equation*}
  This shows that the left adjoint functor
  $\kos{f}_!:\cat{Rep}(\pi^{\pr})\to \cat{Rep}(\pi)$
  is well-defined.
  We explain the adjunction unit $\eta^f$ and counit $\epsilon^f$
  associated to $\kos{f}_!\dashv\kos{f}^*$.
  The component of the adjunction unit $\eta^f$ at
  each object $Z=(z,\gamma^{\pr}_z)$ in $\cat{Rep}(\pi^{\pr})$ is
  the morphism $\eta^f_Z:Z\to \kos{f}^*\kos{f}_!(Z)$ in $\cat{Rep}(\pi^{\pr})$
  whose underlying morphism in $\CK$ is
  \begin{equation*}
    \xymatrix{
      \eta^f_Z:
      z
      \ar[r]^-{\imath_z}_-{\cong}
      &\kappa\otimes z
      \ar[r]^-{u_{\pi}\otimes I_z}
      &\pi\otimes z
      \ar@{->>}[r]^-{\textit{cq}^f_Z}
      &\pi\otimes_{\pi^{\pr}}z
      .
    }
  \end{equation*}
  We can check that $\eta^f_Z$ is indeed a morphism in $\cat{Rep}(\pi^{\pr})$
  as follows.
  \begin{equation*}
    \vcenter
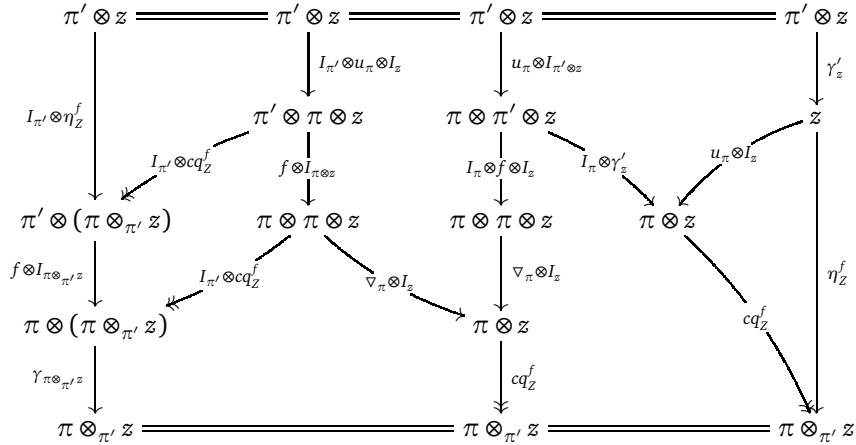
{\hbox{
      \xymatrix{
        \pi^{\pr}\otimes z
        \ar[dd]_-{I_{\pi^{\pr}}\otimes \eta^f_Z}
        \ar@{=}[r]
        &\pi^{\pr}\otimes z
        \ar[d]^-{I_{\pi^{\pr}}\otimes u_{\pi}\otimes I_z}
        \ar@{=}[r]
        &\pi^{\pr}\otimes z
        \ar[d]^-{u_{\pi}\otimes I_{\pi^{\pr}\otimes z}}
        \ar@{=}[rr]
        &\text{ }
        &\pi^{\pr}\otimes z
        \ar[d]^-{\gamma^{\pr}_z}
        \\
        \text{ }
        &\pi^{\pr}\otimes \pi\otimes z
        \ar@{->>}@/_0.5pc/[dl]|-{I_{\pi^{\pr}}\otimes \textit{cq}^f_Z}
        \ar[d]|-{f\otimes I_{\pi\otimes z}}
        &\pi\otimes \pi^{\pr}\otimes z
        \ar[d]|-{I_{\pi}\otimes f\otimes I_z}
        \ar@/^0.5pc/[dr]|-{I_{\pi}\otimes \gamma^{\pr}_z}
        &\text{ }
        &z
        \ar@/_0.5pc/[dl]|-{u_{\pi}\otimes I_z}
        \ar[ddd]^-{\eta^f_Z}
        \\
        \pi^{\pr}\otimes (\pi\otimes_{\pi^{\pr}}z)
        \ar[d]_-{f\otimes I_{\pi\otimes_{\pi^{\pr}}z}}
        &\pi\otimes \pi\otimes z
        \ar@{->>}@/^0.5pc/[dl]|-{I_{\pi^{\pr}}\otimes \textit{cq}^f_Z}
        \ar@/_0.5pc/[dr]|-{\pc_{\pi}\otimes I_z}
        &\pi\otimes \pi\otimes z
        \ar[d]^-{\pc_{\pi}\otimes I_z}
        &\pi\otimes z
        \ar@{->>}@/^0.5pc/[ddr]|-{\textit{cq}^f_Z}
        &\text{ }
        \\
        \pi\otimes (\pi\otimes_{\pi^{\pr}}z)
        \ar[d]_-{\gamma_{\pi\otimes_{\pi^{\pr}}z}}
        &\text{ }
        &\pi\otimes z
        \ar@{->>}[d]^-{\textit{cq}^f_Z}
        &\text{ }
        &\text{ }
        \\
        \pi\otimes_{\pi^{\pr}}z
        \ar@{=}[rr]
        &\text{ }
        &\pi\otimes_{\pi^{\pr}}z
        \ar@{=}[rr]
        &\text{ }
        &\pi\otimes_{\pi^{\pr}}z
      }
    }}
  \end{equation*}
  We leave for the readers to check that
  $\eta^f_Z$ is natural in variable $Z$.
  The component of the adjunction counit $\epsilon^f$
  at each object $X=(x,\gamma_x)$ in $\cat{Rep}(\pi)$
  is the unique morphism
  $\epsilon^f_X:\kos{f}_!\kos{f}^*(X)\to X$
  in $\cat{Rep}(\pi)$
  whose underlying morphism in $\CK$ satisfies the following relation.
  \begin{equation*}
    \vcenter{\hbox{
      \xymatrix@C=40pt{
        \pi\otimes x
        \ar@{->>}[r]^-{\textit{cq}^f_{\kos{f}^*(X)}}
        \ar@/_1pc/[dr]_-{\gamma_x}
        &\pi\otimes_{\pi^{\pr}}x
        \ar@{.>}[d]^-{\epsilon^f_X}_-{\exists!}
        \\
        \text{ }
        &x
      }
    }}
  \end{equation*}
  One can check that the morphism $\epsilon^f_X$ in $\CK$ is well-defined.
  The morphism $\epsilon^f_X$ is indeed a morphism in $\cat{Rep}(\pi)$
  which can be checked by right-cancelling the epimorphism
  $I_{\pi}\otimes \textit{cq}^f_{\kos{f}^*(X)}$
  in the diagram below.
  \begin{equation*}
    \vcenter
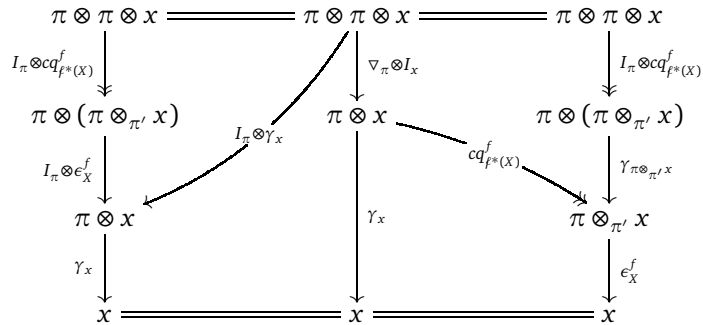
{\hbox{
      \xymatrix@C=40pt{
        \pi\otimes \pi\otimes x
        \ar@{->>}[d]_-{I_{\pi}\otimes \textit{cq}^f_{\kos{f}^*(X)}}
        \ar@{=}[r]
        &\pi\otimes \pi\otimes x
        \ar@/^1pc/[ddl]|-{I_{\pi}\otimes \gamma_x}
        \ar[d]^-{\pc_{\pi}\otimes I_x}
        \ar@{=}[r]
        &\pi\otimes \pi\otimes x
        \ar@{->>}[d]^-{I_{\pi}\otimes \textit{cq}^f_{\kos{f}^*(X)}}
        \\
        \pi\otimes (\pi\otimes_{\pi^{\pr}}x)
        \ar[d]_-{I_{\pi}\otimes \epsilon^f_X}
        &\pi\otimes x
        \ar@{->>}@/^0.5pc/[dr]|-{\textit{cq}^f_{\kos{f}^*(X)}}
        \ar[dd]^-{\gamma_x}
        &\pi\otimes (\pi\otimes_{\pi^{\pr}}x)
        \ar[d]^-{\gamma_{\pi\otimes_{\pi^{\pr}}x}}
        \\
        \pi\otimes x
        \ar[d]_-{\gamma_x}
        &\text{ }
        &\pi\otimes_{\pi^{\pr}}x
        \ar[d]^-{\epsilon^f_X}
        \\
        x
        \ar@{=}[r]
        &x
        \ar@{=}[r]
        &x
      }
    }}
  \end{equation*}
  We leave for the readers to check that
  $\epsilon^f_X$ is natural in variable $X$.
  We need to verify that the adjunction unit, counit satisfy the triangle identities.
  For each object $X=(x,\gamma_x)$ in $\cat{Rep}(\pi)$,
  we verify one of the triangle identities as follows.
  \begin{equation*}
    \vcenter{\hbox{
      \xymatrix{
        \kos{f}^*(X)
        \ar[r]^-{\eta^f_{\kos{f}^*(X)}}
        \ar@{=}@/_1pc/[dr]
        &\kos{f}^*\kos{f}_!\kos{f}^*(X)
        \ar[d]^-{\kos{f}^*(\epsilon^f_X)}
        \\
        \text{ }
        &\kos{f}^*(X)
      }
    }}
    \qquad\quad
    \vcenter{\hbox{
      \xymatrix@C=30pt{
        x
        \ar[dd]_-{\eta^f_{\kos{f}^*(X)}}
        \ar@{=}[r]
        &x
        \ar[d]_-{u_{\pi}\otimes I_x}
        \ar@{=}[r]
        &x
        \ar@{=}[ddd]
        \\
        \text{ }
        &\pi\otimes x
        \ar@{->>}[d]_-{\textit{cq}^f_{\kos{f}^*(X)}}
        \ar@/^1pc/[ddr]^(0.35){\gamma_x}
        &\text{ }
        \\
        \pi\otimes_{\pi^{\pr}}x
        \ar[d]_-{\kos{f}^*(\epsilon^f_X)}
        \ar@{=}[r]
        &\pi\otimes_{\pi^{\pr}}x
        \ar[d]_-{\epsilon^f_X}
        &\text{ }
        \\
        x
        \ar@{=}[r]
        &x
        \ar@{=}[r]
        &x
      }
    }}
  \end{equation*}
  For each object $Z=(z,\gamma^{\pr}_z)$ in $\cat{Rep}(\pi^{\pr})$,
  we obtain the other triangle identity
  by right-cancelling the epimorphism 
  $\textit{cq}^f_Z$ in the diagram below.
  \begin{equation*}
    \vcenter{\hbox{
      \xymatrix{
        \kos{f}_!(Z)
        \ar[r]^-{\kos{f}_!(\eta^f_Z)}
        \ar@{=}@/_1pc/[dr]
        &\kos{f}_!\kos{f}^*\kos{f}_!(Z)
        \ar[d]^-{\epsilon^f_{\kos{f}_!(Z)}}
        \\
        \text{ }
        &\kos{f}_!(Z)
      }
    }}
    \vcenter{\hbox{
      \xymatrix@C=10pt{
        \pi\otimes z
        \ar@{->>}[d]_-{\textit{cq}^f_Z}
        \ar@{=}[r]
        &\pi\otimes z
        \ar[dd]_-{I_\pi\otimes \eta^f_Z}
        \ar@{=}[rr]
        &\text{ }
        &\pi\otimes z
        \ar@/_0.5pc/[dl]|-{I_{\pi}\otimes u_{\pi}\otimes I_z}
        \ar@{=}[dd]
        \\
        \pi\otimes_{\pi^{\pr}}z
        \ar[dd]_-{\kos{f}_!(\eta^f_Z)}
        &\text{ }
        &\pi\otimes \pi\otimes z
        \ar@{->>}@/_0.5pc/[dl]|-{I_{\pi}\otimes \textit{cq}^f_Z}
        \ar@/_0.5pc/[dr]|-{\pc_{\pi}\otimes I_z}
        &\text{ }
        \\
        \text{ }
        &\pi\otimes (\pi\otimes_{\pi^{\pr}}z)
        \ar@/_0.5pc/@{->>}[dl]|-{\textit{cq}^f_{\kos{f}^*\kos{f}_!(Z)}}
        \ar[dd]^-{\gamma_{\pi\otimes_{\pi^{\pr}}z}}
        &\text{ }
        &\pi\otimes z
        \ar@{->>}[dd]^-{\textit{cq}^f_Z}
        \\
        \pi\otimes_{\pi^{\pr}}(\pi\otimes_{\pi^{\pr}}z)
        \ar[d]_-{\epsilon^f_{\kos{f}_!(Z)}}
        &\text{ }
        &\text{ }
        &\text{ }
        \\
        \pi\otimes_{\pi^{\pr}}z
        \ar@{=}[r]
        &\pi\otimes_{\pi^{\pr}}z
        \ar@{=}[rr]
        &\text{ }
        &\pi\otimes_{\pi^{\pr}}z
      }
    }}
  \end{equation*}
  This shows that we have the Galois morphism
  $\kos{f}:\kos{R\!e\!p}(\pi^{\pr})\to \kos{R\!e\!p}(\pi)$
  as we claimed.
  Moreover, we have the invertible Galois morphism
  $\what{\kos{f}}^*:\kos{t}_{\pi}\kos{f}\Rightarrow\kos{t}_{\pi^{\pr}}$
  which is the identity natural transformation
  \begin{equation*}
    \vcenter{\hbox{
      \xymatrix@R=30pt@C=40pt{
        \text{ }
        &\kos{R\!e\!p}(\pi)
        \ar[d]^-{\kos{f}^*}
        \\
        \kos{K}
        \ar@/^0.7pc/[ur]^-{\kos{t}^*_{\pi}}
        \ar[r]_-{\kos{t}^*_{\pi^{\pr}}}
        \xtwocell[r]{}<>{<-3>{\text{ }\text{ }\what{\kos{f}}^*}}
        &\kos{R\!e\!p}(\pi^{\pr})
      }
    }}
    \qquad\qquad
    \what{\kos{f}}^*:
    \kos{f}^*\kos{t}^*_{\pi}=\kos{t}^*_{\pi^{\pr}}.
  \end{equation*}
  Thus we have a well-defined Galois $\kos{K}$-morphism
  $\mathfrak{f}=(\kos{f},\what{\kos{f}}):\mathfrak{Rep}(\pi^{\pr})\to \mathfrak{Rep}(\pi)$.
  Finally, the invertible Galois $\kos{K}$-transformation
  $\xymatrix@C=15pt{
    I:
    \varpi_{\pi}
    \ar@2{->}[r]^-{\cong}
    &\mathfrak{f}\varpi_{\pi^{\pr}}
  }$
  is also given by identity natural transformation.
  \begin{equation*}
    \vcenter{\hbox{
      \xymatrix@C=15pt{
        \mathfrak{Rep}(\pi)
        &\text{ }
        &\mathfrak{Rep}(\pi^{\pr})
        \ar[ll]_-{\mathfrak{f}}
        \\
        \text{ }
        &\mathfrak{K}
        \ar[ul]^-{\varpi_{\pi}}
        \ar[ur]_-{\varpi_{\pi^{\pr}}}
        \xtwocell[u]{}<>{<0>{I}}
      }
    }}
    \qquad\quad
    I:
    (\omega_{\pi}^*,\what{\omega}_{\pi}^*)
    =
    (\omega_{\pi^{\pr}}^*\kos{f}^*,\what{\omega_{\pi^{\pr}}^*\kos{f}^*})
    .
  \end{equation*}

  Now We describe how the $2$-functor
  $\mathbb{GALOBJ}^{\cat{pre}}(\kos{K})
  \to \mathbb{GALCAT}^{\cat{pre}}_*(\kos{K})$
  sends $2$-cells to $2$-cells.
  Let $\theta:f_1\Rightarrow f_2$ be a $2$-cell between
  morphisms $f_1$, $f_2:\pi^{\pr}\to\pi$
  of pre-Galois objects $\pi$, $\pi^{\pr}$ in $\kos{K}$.
  We claim that there is an invertible Galois $\kos{K}$-transformation
  \begin{equation*}
    \vcenter{\hbox{
      \xymatrix@C=40pt{
        (\mathfrak{Rep}(\pi),\varpi_{\pi})
        \xtwocell[r]{}<>{<0>{\text{ }\text{ }\vartheta}}
        &(\mathfrak{Rep}(\pi^{\pr}),\varpi_{\pi^{\pr}})
        \ar@/_1pc/[l]_-{(\mathfrak{f}_1,I)}
        \ar@/^1pc/[l]^-{(\mathfrak{f}_2,I)}
      }
    }}
    \qquad\quad
    \vartheta:
    \xymatrix@C=18pt{
      \mathfrak{f}_1
      \ar@2{->}[r]^-{\cong}
      &\mathfrak{f}_2
    }
  \end{equation*}
  whose component $\vartheta_X:\kos{f}_1^*(X)\xrightarrow[]{\cong} \kos{f}_2^*(X)$
  at each object $X=(x,\gamma_x)$ in $\cat{Rep}(\pi)$ is
  \begin{equation*}
    \vartheta_X:
    \xymatrix{
      x
      \ar[r]^-{\imath_x}_-{\cong}
      &\kappa\otimes x
      \ar[r]^-{\theta\otimes I_x}
      &\pi\otimes x
      \ar[r]^-{\gamma_x}
      &x
      .
    }
  \end{equation*}
  We can check that $\vartheta_X$ is a morphism in $\cat{Rep}(\pi^{\pr})$
  as follows.
  \begin{equation*}
    \vcenter
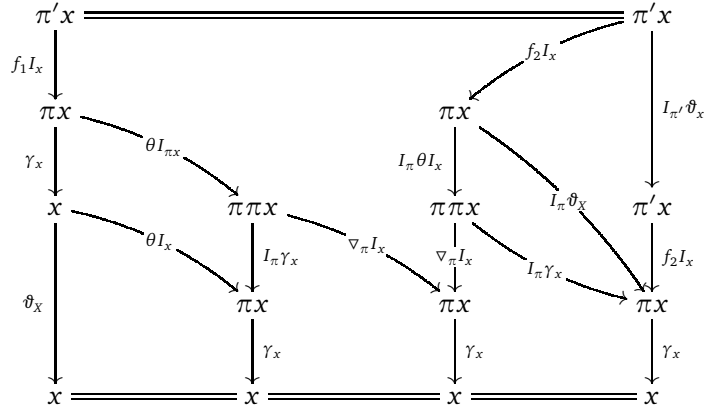
{\hbox{
      \xymatrix@C=50pt{
        \pi^{\pr}x
        \ar[d]_-{f_1I_x}
        \ar@{=}[rrr]
        &\text{ }
        &\text{ }
        &\pi^{\pr}x
        \ar@/_0.5pc/[dl]|-{f_2I_x}
        \ar[dd]^-{I_{\pi^{\pr}}\vartheta_x}
        \\
        \pi x
        \ar[d]_-{\gamma_x}
        \ar@/^0.5pc/[dr]|-{\theta I_{\pi x}}
        &\text{ }
        &\pi x
        \ar[d]_-{I_{\pi}\theta I_x}
        \ar@/^0.5pc/[ddr]|-{I_{\pi}\vartheta_X}
        &\text{ }
        \\
        x
        \ar[dd]_-{\vartheta_X}
        \ar@/^0.5pc/[dr]|-{\theta I_x}
        &\pi \pi x
        \ar[d]^-{I_{\pi}\gamma_x}
        \ar@/^0.5pc/[dr]|-{\pc_{\pi}I_x}
        &\pi \pi x
        \ar[d]|-{\pc_{\pi}I_x}
        \ar@/_0.5pc/[dr]|-{I_{\pi}\gamma_x}
        &\pi^{\pr}x
        \ar[d]^-{f_2I_x}
        \\
        \text{ }
        &\pi x
        \ar[d]^-{\gamma_x}
        &\pi x
        \ar[d]^-{\gamma_x}
        &\pi x
        \ar[d]^-{\gamma_x}
        \\
        x
        \ar@{=}[r]
        &x
        \ar@{=}[r]
        &x
        \ar@{=}[r]
        &x
      }
    }}
  \end{equation*}
  Let $Y=(y,\gamma_y)$ be another object in $\cat{Rep}(\pi)$
  and let $w$ be an object in $\CK$.
  The following diagrams show that $\vartheta$ is a comonoidal $\kos{K}$-tensor natural transformation.
  \begin{equation*}
    \vcenter
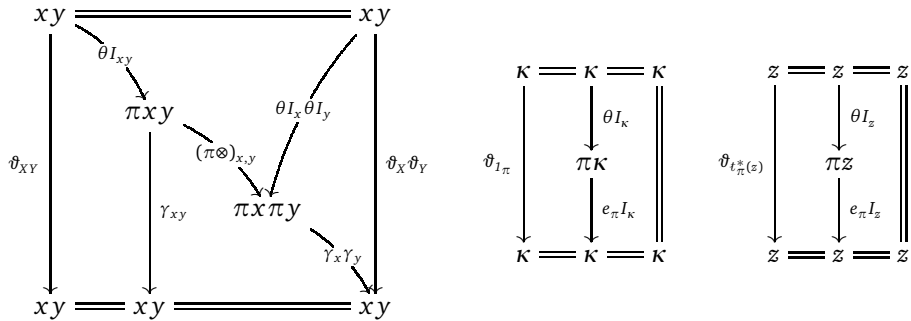
{\hbox{
      \xymatrix@C=15pt{
        xy
        \ar[ddd]_-{\vartheta_{XY}}
        \ar@/^0.5pc/[dr]|-{\theta I_{xy}}
        \ar@{=}[rrr]
        &\text{ }
        &\text{ }
        &xy
        \ar@/_0.5pc/[ddl]|-{\theta I_x \theta I_y}
        \ar[ddd]^-{\vartheta_X\vartheta_Y}
        \\
        \text{ }
        &\pi xy
        \ar[dd]^-{\gamma_{xy}}
        \ar@/^0.5pc/[dr]|-{(\pi\otimes)_{x,y}}
        &\text{ }
        &\text{ }
        \\
        \text{ }
        &\text{ }
        &\pi x\pi y
        \ar@/^0.5pc/[dr]|-{\gamma_x\gamma_y}
        &\text{ }
        \\
        xy
        \ar@{=}[r]
        &xy
        \ar@{=}[rr]
        &\text{ }
        &xy
      }
    }}
    \quad
    \vcenter{\hbox{
      \xymatrix@C=10pt{
        \kappa
        \ar[dd]_-{\vartheta_{\unit_{\pi}}}
        \ar@{=}[r]
        &\kappa
        \ar[d]^-{\theta I_{\kappa}}
        \ar@{=}[r]
        &\kappa
        \ar@{=}[dd]
        \\
        \text{ }
        &\pi\kappa
        \ar[d]^-{e_{\pi}I_{\kappa}}
        &\text{ }
        \\
        \kappa
        \ar@{=}[r]
        &\kappa
        \ar@{=}[r]
        &\kappa
      }
    }}
    \quad
    \vcenter{\hbox{
      \xymatrix@C=10pt{
        z
        \ar[dd]_-{\vartheta_{\kos{t}^*_{\pi}(z)}}
        \ar@{=}[r]
        &z
        \ar[d]^-{\theta I_z}
        \ar@{=}[r]
        &z
        \ar@{=}[dd]
        \\
        \text{ }
        &\pi z
        \ar[d]^-{e_{\pi}I_z}
        &\text{ }
        \\
        z
        \ar@{=}[r]
        &z
        \ar@{=}[r]
        &z
      }
    }}
  \end{equation*}
  One can also explicitly check that
  the inverse of $\vartheta_X$ is given by
  \begin{equation*}
    \vartheta_X^{-1}:
    \xymatrix{
      x
      \ar[r]^-{\imath_x}_-{\cong}
      &\kappa\otimes x
      \ar[r]^-{\theta\otimes I_x}
      &\pi\otimes x
      \ar[r]^-{\varsigma_{\pi}\otimes I_x}_-{\cong}
      &\pi\otimes x
      \ar[r]^-{\gamma_x}
      &x
      .
    }
  \end{equation*}
  This shows that the invertible Galois $\kos{K}$-transformation
  $\vartheta$ is well-defined.

  It is immediate that the correspondence
  $\mathbb{GALOBJ}^{\cat{pre}}(\kos{K})
  \to \mathbb{GALCAT}^{\cat{pre}}_*(\kos{K})$
  preserves identity $1$-cells, identity $2$-cells
  and horizontal composition of $1$-cells.
  We can check that the correspondence
  $\mathbb{GALOBJ}^{\cat{pre}}(\kos{K})
  \to \mathbb{GALCAT}^{\cat{pre}}_*(\kos{K})$
  preserves vertical composition of $2$-cells as follows.
  Let $f_1$, $f_2$, $f_3:\pi^{\pr}\to \pi$
  be morphisms of pre-Galois objects in $\kos{K}$.
  let $\theta_1:f_1\Rightarrow f_2$, $\theta_2:f_2\Rightarrow f_3$
  be $2$-cells.
  If we denote $\vartheta_1$, $\vartheta_2$, $\vartheta$
  as the Galois $\kos{K}$-transformations
  corresponding to $\theta_1$, $\theta_2$, $\theta_2\star\theta_1$,
  then we have the following.
  \begin{equation*}
    \vcenter
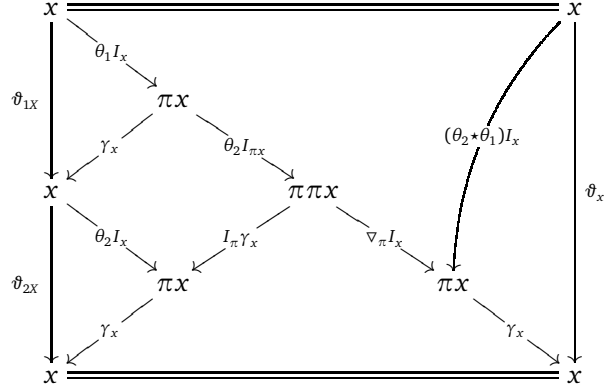
{\hbox{
      \xymatrix@C=30pt{
        x
        \ar[dd]_-{\vartheta_{1X}}
        \ar@{=}[rrrr]
        \ar[dr]|-{\theta_1 I_x}
        &\text{ }
        &\text{ }
        &\text{ }
        &x
        \ar@/_1pc/[dddl]|-{(\theta_2\star\theta_1)I_x}
        \ar[dddd]^-{\vartheta_x}
        \\
        \text{ }
        &\pi x
        \ar[dl]|-{\gamma_x}
        \ar[dr]|-{\theta_2 I_{\pi x}}
        &\text{ }
        &\text{ }
        &\text{ }
        \\
        x
        \ar[dd]_-{\vartheta_{2X}}
        \ar[dr]|-{\theta_2 I_x}
        &\text{ }
        &\pi \pi x
        \ar[dl]|-{I_{\pi}\gamma_x}
        \ar[dr]|-{\pc_{\pi}I_x}
        &\text{ }
        &\text{ }
        \\
        \text{ }
        &\pi x
        \ar[dl]|-{\gamma_x}
        &\text{ }
        &\pi x
        \ar[dr]|-{\gamma_x}
        &\text{ }
        \\
        x
        \ar@{=}[rrrr]
        &\text{ }
        &\text{ }
        &\text{ }
        &x
      }
    }}
  \end{equation*}
  This shows that the correspondence
  $\mathbb{GALOBJ}^{\cat{pre}}(\kos{K})
  \to \mathbb{GALCAT}^{\cat{pre}}_*(\kos{K})$
  preserves vertical composition of $2$-cells as follows.
  We can also check that the correspondence
  preserves horizontal composition of $2$-cells as follows.
  Consider the following situation
  \begin{equation*}
    \vcenter
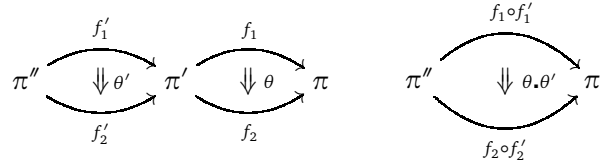
{\hbox{
      \xymatrix@C=40pt{
        \pi^{\ppr}
        \ar@/^1pc/[r]^-{f_1^{\pr}}
        \ar@/_1pc/[r]_-{f_2^{\pr}}
        \xtwocell[r]{}<>{<0>{\text{ }\text{ }\theta^{\pr}}}
        &\pi^{\pr}
        \ar@/^1pc/[r]^-{f_1}
        \ar@/_1pc/[r]_-{f_2}
        \xtwocell[r]{}<>{<0>{\text{ }\text{ }\theta}}
        &\pi
      }
    }}
    \qquad
    \vcenter{\hbox{
      \xymatrix@C=50pt{
        \pi^{\ppr}
        \ar@/^1.5pc/[r]^-{f_1\circ f_1^{\pr}}
        \ar@/_1.5pc/[r]_-{f_2\circ f_2^{\pr}}
        \xtwocell[r]{}<>{<0>{\quad\text{ }\text{ }\theta\centerdot \theta^{\pr}}}
        &\pi
      }
    }}
  \end{equation*}
  \begin{equation*}
    \theta\centerdot\theta^{\pr}
    =
    \theta\star f_1\theta^{\pr}
    =
    f_2\theta^{\pr}\star \theta:f_1\circ f^{\pr}_1\Rightarrow f_2\circ f^{\pr}_2.
  \end{equation*}
  and denote $\vartheta$, $\vartheta^{\pr}$, $\tilde{\vartheta}$
  as the Galois $\kos{K}$-transformations
  corresponding to $\theta$, $\theta^{\pr}$, $\theta\centerdot\theta^{\pr}$.
  We denote the horizontal composition of $\theta$ and $\theta^{\pr}$ as
  $\vartheta\vartheta^{\pr}:\mathfrak{f}_1\mathfrak{f}_1^{\pr}\Rightarrow \mathfrak{f}_2\mathfrak{f}_2^{\pr}$.
  We have the following relation for each object 
  $X=(x,\gamma_x)$ in $\cat{Rep}(\pi)$.
  \begin{equation*}
    \vcenter{\hbox{
      \xymatrix@C=40pt{
        x
        \ar[ddddd]_-{(\vartheta\vartheta^{\pr})_X}
        \ar@{=}[r]
        &x
        \ar[ddd]_-{\vartheta^{\pr}_{\kos{f}_1^*(X)}}
        \ar[dr]|-{\theta^{\pr} I_x}
        \ar@{=}[rrrr]
        &\text{ }
        &\text{ }
        &\text{ }
        &x
        \ar@/_2pc/[ddddl]|-{(\theta\star f_1\theta^{\pr})I_x=(\theta\centerdot \theta^{\pr})I_x}
        \ar[ddddd]^-{\tilde{\vartheta}_X}
        \\
        \text{ }
        &\text{ }
        &\pi^{\pr} x
        \ar[d]^-{f_1 I_x}
        &\text{ }
        &\text{ }
        &\text{ }
        \\
        \text{ }
        &\text{ }
        &\pi x
        \ar[dl]|-{\gamma_x}
        \ar[dr]|-{\theta I_{\pi x}}
        &\text{ }
        &\text{ }
        &\text{ }
        \\
        \text{ }
        &x
        \ar[dd]_-{\kos{f}^{\pr*}_2(\vartheta_X)}
        \ar[dr]|-{\theta I_x}
        &\text{ }
        &\pi \pi x
        \ar[dl]|-{I_{\pi}\gamma_x}
        \ar[dr]|-{\pc_{\pi}I_x}
        &\text{ }
        &\text{ }
        \\
        \text{ }
        &\text{ }
        &\pi x
        \ar[dl]|-{\gamma_x}
        &\text{ }
        &\pi x
        \ar[dr]|-{\gamma_x}
        &\text{ }
        \\
        x
        \ar@{=}[r]
        &x
        \ar@{=}[rrrr]
        &\text{ }
        &\text{ }
        &\text{ }
        &x
      }
    }}
  \end{equation*}
  We can also obtain the same relation as follows.
  \begin{equation*}
    \vcenter{\hbox{
      \xymatrix@C=40pt{
        x
        \ar[ddddd]_-{(\vartheta\vartheta^{\pr})_X}
        \ar@{=}[r]
        &x
        \ar[dd]_-{\kos{f}_1^*(\theta_X)}
        \ar[dr]|-{\theta I_x}
        \ar@{=}[rrrr]
        &\text{ }
        &\text{ }
        &\text{ }
        &x
        \ar@/_2pc/[ddddl]|-{(f_2\theta^{\pr}\star \theta)I_x=(\theta\centerdot\theta^{\pr})I_x}
        \ar[ddddd]^-{\tilde{\vartheta}_X}
        \\
        \text{ }
        &\text{ }
        &\pi x
        \ar[dl]|-{\gamma_x}
        \ar[dr]|-{\theta^{\pr}I_{\pi x}}
        &\text{ }
        &\text{ }
        &\text{ }
        \\
        \text{ }
        &x
        \ar[ddd]_-{\theta^{\pr}_{\kos{f}_2^*(X)}}
        \ar[dr]|-{\theta^{\pr} I_x}
        &\text{ }
        &\pi^{\pr}\pi x
        \ar[dl]|-{I_{\pi^{\pr}}\gamma_x}
        \ar[d]^-{f_2I_{\pi x}}
        &\text{ }
        &\text{ }
        \\
        \text{ }
        &\text{ }
        &\pi^{\pr} x
        \ar[d]_-{f_2 I_x}
        &\pi \pi x
        \ar[dl]|-{I_{\pi}\gamma_x}
        \ar[dr]|-{\pc_{\pi}I_x}
        &\text{ }
        &\text{ }
        \\
        \text{ }
        &\text{ }
        &\pi x
        \ar[dl]|-{\gamma_x}
        &\text{ }
        &\pi x
        \ar[dr]|-{\gamma_x}
        &\text{ }
        \\
        x
        \ar@{=}[r]
        &x
        \ar@{=}[rrrr]
        &\text{ }
        &\text{ }
        &\text{ }
        &x
      }
    }}
  \end{equation*}
  We conclude that
  we have a well-defined $2$-functor
  $\mathbb{GALOBJ}^{\cat{pre}}(\kos{K})
  \to \mathbb{GALCAT}^{\cat{pre}}_*(\kos{K})$.
  This completes the proof of Proposition~\ref{prop inducedfundamentalfunctor}.
\qed\end{proof}

\begin{proposition} \label{prop CATtoOBJ}
  Let $\kos{K}$ be a Galois prekosmos $\kos{K}$.
  We have a $2$-functor
  \begin{equation*}
    \vcenter{\hbox{
      \xymatrix@C=30pt{
        \mathbb{GALCAT}^{\cat{pre}}_*(\kos{K})
        \ar[r]
        &\mathbb{GALOBJ}^{\cat{pre}}(\kos{K}).
      }
    }}
  \end{equation*}
\end{proposition}
\begin{proof}
  First, we describe how the $2$-functor
  $\mathbb{GALCAT}^{\cat{pre}}_*(\kos{K})
  \to \mathbb{GALOBJ}^{\cat{pre}}(\kos{K})$
  sends $0$-cells to $0$-cells.
  For each pointed pre-Galois $\kos{K}$-category
  $(\mathfrak{T},\varpi)$
  we have a pre-Galois object $\omega^*\omega_!(\kappa)$ in $\kos{K}$
  as explained in Theorem~\ref{thm preGalKCat mainThm1}.
  
  Next, we describe how the $2$-functor
  $\mathbb{GALCAT}^{\cat{pre}}_*(\kos{K})
  \to \mathbb{GALOBJ}^{\cat{pre}}(\kos{K})$
  sends $1$-cells to $1$-cells.
  Let 
  $(\mathfrak{f},\alpha):(\mathfrak{T}^{\pr},\varpi^{\pr})\to (\mathfrak{T},\varpi)$
  be a $1$-cell between pointed pre-Galois $\kos{K}$-categories.
  \begin{equation*}
    \vcenter{\hbox{
      \xymatrix{
        \mathfrak{T}
        &\text{ }
        &\mathfrak{T}^{\pr}
        \ar[ll]_-{\mathfrak{f}}
        \\
        \text{ }
        &\mathfrak{K}
        \ar[ur]_-{\varpi^{\pr}}
        \ar[ul]^-{\varpi}
        \xtwocell[u]{}<>{<0>{\text{ }\text{ }\alpha}}
      }
    }}
    \qquad\quad
    \alpha:
    \xymatrix@C=18pt{
      \omega^*
      \ar@2{->}[r]^-{\cong}
      &\omega^{\pr*}\kos{f}^*
      :\kos{T}\to\kos{K}
    }
    \qquad
    \vcenter{\hbox{
      \xymatrix{
        \kos{K}
        \ar@{=}[d]
        &\kos{T}
        \ar[l]_-{\omega^*}
        \ar[d]^-{\kos{f}^*}
        \\
        \kos{K}
        \xtwocell[ur]{}<>{<0>{\text{ }\text{ }\alpha}}
        &\kos{T}^{\pr}
        \ar[l]^-{\omega^{\pr*}}
      }
    }}
  \end{equation*}
  We denote the mate associated to $\alpha$ as
  \begin{equation*}
    \vcenter{\hbox{
      \xymatrix{
        \kos{K}
        \ar@{=}[d]
        \ar[r]^-{\omega_!}
        &\kos{T}
        \ar[d]^-{\kos{f}^*}
        \\
        \kos{K}
        \ar[r]_-{\omega^{\pr}_!}
        &\kos{T}^{\pr}
        \xtwocell[ul]{}<>{<0>{\alpha_!\text{ }\text{ }}}
      }
    }}
    \qquad\quad
    \alpha_!:
    \xymatrix{
      \omega^{\pr}_!
      \ar@2{->}[r]^-{\omega^{\pr}_!\eta}
      &\omega^{\pr}_!\omega^*\omega_!
      \ar@2{->}[r]^-{\omega^{\pr}_!\alpha\omega_!}_-{\cong}
      &\omega^{\pr}_!\omega^{\pr*}\kos{f}^*\omega_!
      \ar@2{->}[r]^-{\epsilon^{\pr}\kos{f}^*\omega_!}
      &\kos{f}^*\omega_!
    }
  \end{equation*}
  which satisfies the following relations.
  \begin{equation}\label{eq alpha! relation}
    \vcenter{\hbox{
      \xymatrix@C=30pt{
        \omega^{\pr}_!\omega^*
        \ar@2{->}[d]_-{\omega^{\pr}_!\alpha}^-{\cong}
        \ar@2{->}[r]^-{\alpha_!\omega^*}
        &\kos{f}^*\omega_!\omega^*
        \ar@2{->}[d]^-{\kos{f}^*\epsilon}
        \\
        \omega^{\pr}_!\omega^{\pr*}\kos{f}^*
        \ar@2{->}[r]_-{\epsilon^{\pr}\kos{f}^*}
        &\kos{f}^*
      }
    }}
    \qquad\quad
    \vcenter{\hbox{
      \xymatrix@C=30pt{
        I_{\kos{K}}
        \ar@2{->}[d]_-{\eta^{\pr}}
        \ar@2{->}[r]^-{\eta}
        &\omega^*\omega_!
        \ar@2{->}[d]^-{\alpha\omega_!}_-{\cong}
        \\
        \omega^{\pr*}\omega^{\pr}_!
        \ar@2{->}[r]_-{\omega^{\pr*}\alpha_!}
        &\omega^{\pr*}\kos{f}^*\omega_!
      }
    }}
  \end{equation}
  If we denote
  \begin{equation*}
    \check{\alpha}:
    \xymatrix@C=30pt{
      \omega^{\pr*}\omega^{\pr}_!
      \ar@2{->}[r]^-{\omega^{\pr*}\alpha_!}
      &\omega^{\pr*}\kos{f}^*\omega_!
      \ar@2{->}[r]^-{\alpha^{-1}\omega_!}_-{\cong}
      &\omega^*\omega_!
    }
  \end{equation*}
  then we claim that
  $\check{\alpha}_{\kappa}:
  \omega^{\pr*}\omega^{\pr}_!(\kappa)
  \to \omega^*\omega_!(\kappa)$
  is a morphism of pre-Galois objects in $\kos{K}$.
  It suffices to show that
  $\check{\alpha}:\omega^{\pr*}\omega^{\pr}_!\Rightarrow \omega^*\omega_!$
  is a morphism of monads,
  as we have the adjoint equivalence 
  $\CL\dashv\iota$
  in Corollary~\ref{cor RepGrpEns(K) rflCKTmonad=Mon(Ens(K))}.
  We can check this as follows.
  We will use $\to$ instead of $\Rightarrow$.
  We have
  \begin{equation*}
    \vcenter{\hbox{
      \xymatrix{
        \omega^{\pr*}\omega^{\pr}_!\omega^{\pr*}\omega^{\pr}_!
        \ar[dddd]_-{\check{\alpha}\check{\alpha}}
        \ar@{=}[rrr]
        \ar@/^0.5pc/[dr]|-{\omega^{\pr*}\omega^{\pr}_!\omega^{\pr*}\alpha_!}
        &\text{ }
        &\text{ }
        &\omega^{\pr*}\omega^{\pr}_!\omega^{\pr*}\omega^{\pr}_!
        \ar[dd]^-{\omega^{\pr*}\epsilon^{\pr}\omega^{\pr}_!}
        \\
        \text{ }
        &\omega^{\pr*}\omega^{\pr}_!\omega^{\pr*}\kos{f}^*\omega_!
        \ar[d]_-{\omega^{\pr*}\omega^{\pr}_!\alpha^{-1}\omega_!}^-{\cong}
        \ar@/^2pc/[dddr]|-{\omega^{\pr*}\epsilon^{\pr}\kos{f}^*\omega_!}
        &\text{ }
        &\text{ }
        \\
        \text{ }
        &\omega^{\pr*}\omega^{\pr}_!\omega^*\omega_!
        \ar[d]_-{\omega^{\pr*}\alpha_!\omega^*\omega_!}
        &\text{ }
        &\omega^{\pr*}\omega^{\pr}_!
        \ar[ddd]^-{\check{\alpha}}
        \ar@/_1pc/[ddl]|-{\omega^{\pr*}\alpha_!}
        \\
        \text{ }
        &\omega^{\pr*}\kos{f}^*\omega_!\omega^*\omega_!
        \ar@/^0.5pc/[dl]^-{\alpha^{-1}\omega_!\omega^*\omega_!}_-{\cong}
        \ar@/_0.5pc/[dr]|-{\omega^{\pr*}\kos{f}^*\epsilon\omega_!}
        \ar@{}[r]|-{(\dagger)}
        &\text{ }
        &\text{ }
        \\
        \omega^*\omega_!\omega^*\omega_!
        \ar[d]_-{\omega^*\epsilon\omega_!}
        &\text{ }
        &\omega^{\pr*}\kos{f}^*\omega_!
        \ar@/_0.5pc/[dr]_-{\alpha^{-1}\omega_!}^-{\cong}
        &\text{ }
        \\
        \omega^*\omega_!
        \ar@{=}[rrr]
        &\text{ }
        &\text{ }
        &\omega^*\omega_!
      }
    }}
  \end{equation*}
  where the diagram $(\dagger)$
  is the left relation of (\ref{eq alpha! relation}).
  We also have
  \begin{equation*}
    \vcenter{\hbox{
      \xymatrix@C=40pt{
        I_{\kos{K}}
        \ar[d]_-{\eta^{\pr}}
        \ar@{=}[rr]
        &\text{ }
        \ar@{}[dd]|-{(\ddagger)}
        &I_{\kos{K}}
        \ar@2{->}[ddd]^-{\eta}
        \\
        \omega^{\pr*}\omega^{\pr}_!
        \ar[dd]_-{\check{\alpha}}
        \ar[dr]|-{\omega^{\pr*}\alpha_!}
        &\text{ }
        &\text{ }
        \\
        \text{ }
        &\omega^{\pr*}\kos{f}^*\omega_!
        \ar[dr]^-{\alpha^{-1}\omega_!}_-{\cong}
        &\text{ }
        \\
        \omega^*\omega_!
        \ar@{=}[rr]
        &\text{ }
        &\omega^*\omega_!
      }
    }}
  \end{equation*}
  where the diagram $(\ddagger)$
  is the right relation of (\ref{eq alpha! relation}).
  This shows that 
  $\check{\alpha}_{\kappa}:\omega^{\pr*}\omega^{\pr}_!(\kappa)\to \omega^*\omega_!(\kappa)$
  is a morphism of pre-Galois objects in $\kos{K}$.

  We also describe how the $2$-functor
  $\mathbb{GALCAT}^{\cat{pre}}_*(\kos{K})
  \to \mathbb{GALOBJ}^{\cat{pre}}(\kos{K})$
  sends $2$-cells to $2$-cells.
  Let 
  $\vartheta:(\mathfrak{f}_1,\alpha_1)\Rightarrow (\mathfrak{f}_2,\alpha_2)$
  be a $2$-cell between $1$-cells
  $(\mathfrak{f}_1,\alpha_1)$, $(\mathfrak{f}_2,\alpha_2):(\mathfrak{T}^{\pr},\varpi^{\pr})\to (\mathfrak{T},\varpi)$.
  If we denote
  \begin{equation*}
    \check{\vartheta}:
    \xymatrix@C=30pt{
      I_{\kos{K}}
      \ar@2{->}[r]^-{\eta}
      &\omega^*\omega_!
      \ar@2{->}[r]^-{\alpha_1\omega_!}_-{\cong}
      &\omega^{\pr*}\kos{f}_1^*\omega_!
      \ar@2{->}[r]^-{\omega^{\pr*}\vartheta\omega_!}_-{\cong}
      &\omega^{\pr*}\kos{f}_2^*\omega_!
      \ar@2{->}[r]^-{\alpha_2^{-1}\omega_!}_-{\cong}
      &\omega^*\omega_!
    }
  \end{equation*}
  then $\check{\vartheta}$ satisfies the following relation.
  \begin{equation} \label{eq checkvartheta relation}
    \vcenter{\hbox{
      \xymatrix@C=30pt{
        \omega^*
        \ar@2{-}[d]_-{\alpha_1}^-{\cong}
        \ar@2{->}[rr]^-{\check{\vartheta}\omega^*}
        &\text{ }
        &\omega^*\omega_!\omega^*
        \ar@2{->}[d]^-{\omega^*\epsilon}
        \\
        \omega^{\pr*}\kos{f}^*_1
        \ar@2{->}[r]^-{\omega^{\pr*}\vartheta}_-{\cong}
        &\omega^{\pr*}\kos{f}^*_2
        \ar@2{->}[r]^-{\alpha_2^{-1}}_-{\cong}
        &\omega^*
      }
    }}
  \end{equation}
  We claim that
  $\check{\vartheta}_{\kappa}:\kappa\to \omega^*\omega_!(\kappa)$
  is a $2$-cell
  $\check{\alpha}_{1\kappa}\Rightarrow \check{\alpha}_{2\kappa}
  :\omega^{\pr*}\omega^{\pr}_!(\kappa)\to \omega^*\omega_!(\kappa)$
  between morphisms of pre-Galois objects in $\kos{K}$.
  It suffices to show the relation
  \begin{equation*}
    \vcenter{\hbox{
      \xymatrix{
        \omega^{\pr*}\omega^{\pr}_!
        \ar@2{->}[d]_-{\check{\vartheta}\check{\alpha}_1}
        \ar@2{->}[r]^-{\check{\alpha}_2\check{\vartheta}}
        &\omega^*\omega_!\omega^*\omega_!
        \ar@2{->}[d]^-{\omega^*\epsilon\omega_!}
        \\
        \omega^*\omega_!\omega^*\omega_!
        \ar@2{->}[r]_-{\omega^*\epsilon\omega_!}
        &\omega^*\omega_!
      }
    }}
  \end{equation*}
  which we verify as follows.
  \begin{equation*}
    \vcenter{\hbox{
      \xymatrix@C=15pt{
        \omega^{\pr*}\omega^{\pr}_!
        \ar[dd]_-{\check{\alpha}_1}
        \ar@{=}[r]
        &\omega^{\pr*}\omega^{\pr}_!
        \ar[ddd]_-{\omega^{\pr*}\alpha_{1!}}
        \ar@{=}[r]
        &\omega^{\pr*}\omega^{\pr}_!
        \ar[d]^-{\omega^{\pr*}\omega^{\pr}_!\eta}
        \ar@{=}[rr]
        &\text{ }
        &\omega^{\pr*}\omega^{\pr}_!
        \ar[d]^-{\omega^{\pr*}\omega^{\pr}_!\check{\vartheta}}
        \\
        \text{ }
        &\text{ }
        &\omega^{\pr*}\omega^{\pr}_!\omega^*\omega_!
        \ar[d]^-{\omega^{\pr*}\omega_!\alpha_1\omega_!}_-{\cong}
        &\text{ }
        &\omega^{\pr*}\omega^{\pr}_!\omega^*\omega_!
        \ar@/^1pc/[ddll]_(0.4){\omega^{\pr*}\omega^{\pr}_!\alpha_2\omega_!}^-{\cong}
        \ar@{=}[d]
        \\
        \omega^*\omega_!
        \ar[dd]_-{\check{\vartheta}\omega^*\omega_!}
        \ar@/_0.5pc/[dr]_-{\alpha_1\omega_!}^-{\cong}
        &\text{ }
        &\omega^{\pr*}\omega^{\pr}_!\omega^{\pr*}\kos{f}_1^*\omega_!
        \ar@/_0.5pc/[dl]|-{\omega^{\pr*}\epsilon^{\pr}\kos{f}_1^*\omega_!}
        \ar[d]^-{\omega^{\pr*}\omega^{\pr}_!\omega^{\pr*}\vartheta\omega_!}_-{\cong}
        &\text{ }
        &\omega^{\pr*}\omega^{\pr}_!\omega^*\omega_!
        \ar@/^0.5pc/[dl]|-{\omega^{\pr*}\alpha_{2!}\omega^*\omega_!}
        \ar[dd]^-{\check{\alpha}_2\omega^*\omega_!}
        \\
        \text{ }
        &\omega^{\pr*}\kos{f}^*_1\omega_!
        \ar@/_0.5pc/[dr]_-{\omega^{\pr*}\vartheta\omega_!}^-{\cong}
        &\omega^{\pr*}\omega^{\pr}_!\omega^{\pr*}\kos{f}_2^*\omega_!
        \ar[d]^-{\omega^{\pr*}\epsilon^{\pr}\kos{f}^*_2\omega_!}
        &\omega^{\pr*}\kos{f}_2^*\omega_!\omega^*\omega_!
        \ar@/^1pc/[dl]^-{\omega^{\pr*}\kos{f}_2^*\epsilon\omega_!}
        \ar@/_1pc/[dr]^-{\alpha_2^{-1}\omega_!\omega^*\omega_!}_-{\cong}
        &\text{ }
        \\
        \omega^*\omega_!\omega^*\omega_!
        \ar[d]_-{\omega^*\epsilon\omega_!}
        &\text{ }
        &\omega^{\pr*}\kos{f}^*_2\omega_!
        \ar[d]^-{\alpha_2^{-1}\omega_!}_-{\cong}
        &\text{ }
        &\omega^*\omega_!\omega^*\omega_!
        \ar[d]^-{\omega^*\epsilon\omega_!}
        \\
        \omega^*\omega_!
        \ar@{=}[rr]
        &\text{ }
        &\omega^*\omega_!
        \ar@{=}[rr]
        &\text{ }
        &\omega^*\omega_!
      }
    }}
  \end{equation*}
  This shows that
  $\check{\vartheta}_{\kappa}:\kappa\to \omega^*\omega_!(\kappa)$
  is a $2$-cell
  $\check{\alpha}_{1\kappa}\Rightarrow \check{\alpha}_{2\kappa}
  :\omega^{\pr*}\omega^{\pr}_!(\kappa)\to \omega^*\omega_!(\kappa)$.

  We leave for the readers to check that
  the correspondence
  $\mathbb{GALCAT}^{\cat{pre}}_*(\kos{K})
  \to \mathbb{GALOBJ}^{\cat{pre}}(\kos{K})$
  preserves identity $1$-cells and identity $2$-cells.
  We can check that
  $\mathbb{GALCAT}^{\cat{pre}}_*(\kos{K})
  \to \mathbb{GALOBJ}^{\cat{pre}}(\kos{K})$
  preserves vertical composition of $2$-cells as follows.
  \begin{equation*}
    \vcenter{\hbox{
      \xymatrix@C=50pt{
        (\mathfrak{T},\varpi)
        \xtwocell[r]{}<>{<-3>{\text{ }\text{ }\vartheta_1}}
        \xtwocell[r]{}<>{<3>{\text{ }\text{ }\vartheta_2}}
        &(\mathfrak{T}^{\pr},\varpi^{\pr})
        \ar@/_2pc/@<-0.5ex>[l]_-{(\mathfrak{f}_1,\alpha_1)}
        \ar[l]|-{(\mathfrak{f}_2,\alpha_2)}
        \ar@/^2pc/@<0.5ex>[l]^-{(\mathfrak{f}_3,\alpha_3)}
      }
    }}
    \qquad\quad
    \begin{aligned}
      \vartheta:=\vartheta_2\circ\vartheta_1
      &:
      \xymatrix@C=18pt{
        (\mathfrak{f}_1,\alpha_1)
        \ar@2{->}[r]
        &(\mathfrak{f}_3,\alpha_3)
        ,
      }
      \\
      \vartheta:=\vartheta_2\circ\vartheta_1
      &:
      \xymatrix@C=18pt{
        \mathfrak{f}_1
        \ar@2{->}[r]^-{\vartheta_1}_-{\cong}
        &\mathfrak{f}_2
        \ar@2{->}[r]^-{\vartheta_2}_-{\cong}
        &\mathfrak{f}_3
        .
      }
    \end{aligned}
  \end{equation*}
  We need to check the relation
  $\xymatrix@C=10pt{
    \check{\vartheta}:
    I_{\kos{K}}
    \ar@2{->}[r]^-{\check{\vartheta}_2\check{\vartheta}_1}
    &\omega^*\omega_!\omega^*\omega_!
    \ar@2{->}[r]^-{\omega^*\epsilon\omega_!}
    &\omega^*\omega_!
  }$
  which we verify as follows.
  \begin{equation*}
    \vcenter
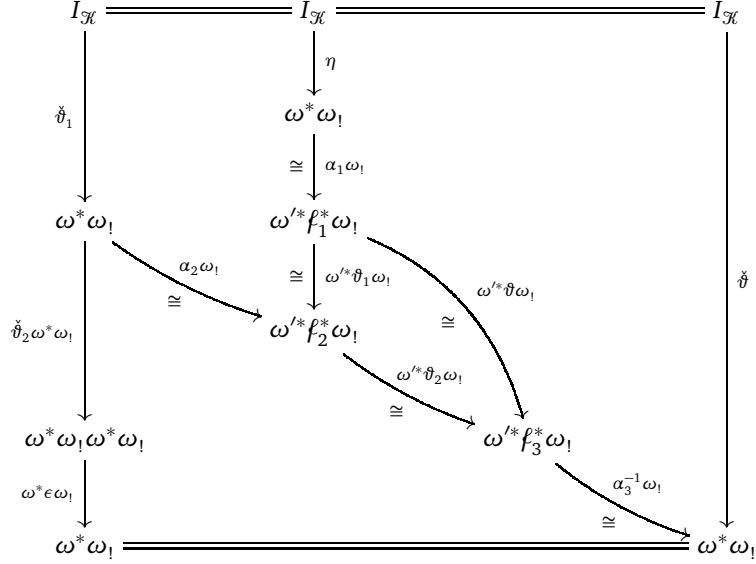
{\hbox{
      \xymatrix@C=40pt{
        I_{\kos{K}}
        \ar[dd]_-{\check{\vartheta}_1}
        \ar@{=}[r]
        &I_{\kos{K}}
        \ar[d]^-{\eta}
        \ar@{=}[rr]
        &\text{ }
        &I_{\kos{K}}
        \ar[ddddd]^-{\check{\vartheta}}
        \\
        \text{ }
        &\omega^*\omega_!
        \ar[d]^-{\alpha_1\omega_!}_-{\cong}
        &\text{ }
        &\text{ }
        \\
        \omega^*\omega_!
        \ar[dd]_-{\check{\vartheta}_2\omega^*\omega_!}
        \ar@/_0.5pc/[dr]^-{\alpha_2\omega_!}_-{\cong}
        &\omega^{\pr*}\kos{f}^*_1\omega_!
        \ar[d]^-{\omega^{\pr*}\vartheta_1\omega_!}_-{\cong}
        \ar@/^1.5pc/[ddr]^-{\omega^{\pr*}\vartheta\omega_!}_-{\cong}
        &\text{ }
        &\text{ }
        \\
        \text{ }
        &\omega^{\pr*}\kos{f}^*_2\omega_!
        \ar@/_0.5pc/[dr]^-{\omega^{\pr*}\vartheta_2\omega_!}_-{\cong}
        &\text{ }
        &\text{ }
        \\
        \omega^*\omega_!\omega^*\omega_!
        \ar[d]_-{\omega^*\epsilon\omega_!}
        &\text{ }
        &\omega^{\pr*}\kos{f}^*_3\omega_!
        \ar@/_0.5pc/[dr]^-{\alpha_3^{-1}\omega_!}_-{\cong}
        &\text{ }
        \\
        \omega^*\omega_!
        \ar@{=}[rrr]
        &\text{ }
        &\text{ }
        &\omega^*\omega_!
      }
    }}
  \end{equation*}
  We also check that the correspondence
  $\mathbb{GALCAT}^{\cat{pre}}_*(\kos{K})
  \to \mathbb{GALOBJ}^{\cat{pre}}(\kos{K})$
  preserves composition of $1$-cells as follows.

  \begin{equation*}
    \vcenter{\hbox{
      \xymatrix@C=15pt{
        (\mathfrak{T},\varpi)
        &(\mathfrak{T}^{\pr},\varpi^{\pr})
        \ar[l]_-{(\mathfrak{f},\alpha)}
        &(\mathfrak{T}^{\ppr},\varpi^{\ppr})
        \ar[l]_-{(\mathfrak{f}^{\pr},\alpha^{\pr})}
        \ar@/^1pc/@<0.5ex>[ll]^-{(\mathfrak{f},\alpha)\circ(\mathfrak{f}^{\pr},\alpha^{\pr})=\bigl(\mathfrak{f}\mathfrak{f}^{\pr},(\mathfrak{f}\alpha^{\pr})\circ\alpha\bigr)}
      }
    }}
    \begin{aligned}
      \beta:=
      (\mathfrak{f}\alpha^{\pr})\circ\alpha
      &:\!\!
      \vcenter{\hbox{
        \xymatrix@C=15pt{
          \varpi
          \ar@2{->}[r]^-{\alpha}_-{\cong}
          &\mathfrak{f}\varpi^{\pr}
          \ar@2{->}[r]^-{\mathfrak{f}\alpha^{\pr}}_-{\cong}
          &\mathfrak{f}\mathfrak{f}^{\pr}\varpi^{\ppr}
        }
      }}
      \\
      \beta:=
      (\alpha^{\pr}\kos{f}^*)\circ\alpha
      &:\!\!
      \vcenter{\hbox{
        \xymatrix@C=15pt{
          \omega^*
          \ar@2{->}[r]^-{\alpha}_-{\cong}
          &\omega^{\pr*}\kos{f}^*
          \ar@2{->}[r]^-{\alpha^{\pr}\kos{f}^*}_-{\cong}
          &\omega^{\ppr*}\kos{f}^{\pr*}\kos{f}^*
        }
      }}
    \end{aligned}
  \end{equation*}
  Then the associated mate of $\beta$ is given by
  $\xymatrix@C=15pt{
    \beta_!:
    \omega^{\ppr}_!
    \ar@2{->}[r]^-{\alpha^{\pr}_!}
    &\kos{f}^{\pr*}\omega^{\pr}_!
    \ar@2{->}[r]^-{\kos{f}^{\pr*}\alpha_!}
    &\kos{f}^{\pr*}\kos{f}^*\omega_!
  }$
  which we can see from the diagram below.
  \begin{equation*}
    \vcenter
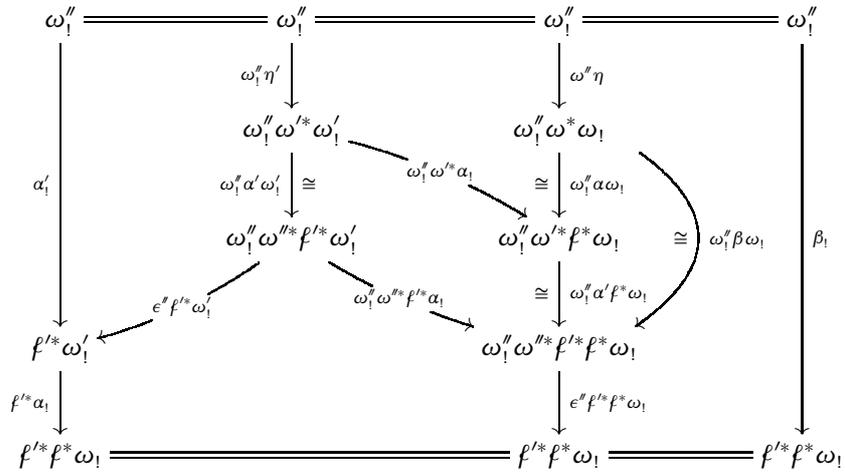
{\hbox{
      \xymatrix@C=40pt{
        \omega^{\ppr}_!
        \ar@{=}[r]
        \ar[ddd]_-{\alpha^{\pr}_!}
        &\omega^{\ppr}_!
        \ar[d]_-{\omega^{\ppr}_!\eta^{\pr}}
        \ar@{=}[r]
        &\omega^{\ppr}_!
        \ar[d]^-{\omega^{\ppr}\eta}
        \ar@{=}[r]
        &\omega^{\ppr}_!
        \ar[dddd]^-{\beta_!}
        \\
        \text{ }
        &\omega^{\ppr}_!\omega^{\pr*}\omega^{\pr}_!
        \ar[d]_-{\omega^{\ppr}_!\alpha^{\pr}\omega^{\pr}_!}^-{\cong}
        \ar@/^0.5pc/[dr]|-{\omega^{\ppr}_!\omega^{\pr*}\alpha_!}
        &\omega^{\ppr}_!\omega^*\omega_!
        \ar[d]^-{\omega^{\ppr}_!\alpha\omega_!}_-{\cong}
        \ar@<3.5ex>@/^3pc/[dd]^-{\omega^{\ppr}_!\beta\omega_!}_-{\cong}
        &\text{ }
        \\
        \text{ }
        &\omega^{\ppr}_!\omega^{\ppr*}\kos{f}^{\pr*}\omega^{\pr}_!
        \ar@/^0.5pc/[dl]|-{\epsilon^{\ppr}\kos{f}^{\pr*}\omega^{\pr}_!}
        \ar@/_0.5pc/[dr]|-{\omega^{\ppr}_!\omega^{\ppr*}\kos{f}^{\pr*}\alpha_!}
        &\omega^{\ppr}_!\omega^{\pr*}\kos{f}^*\omega_!
        \ar[d]^-{\omega^{\ppr}_!\alpha^{\pr}\kos{f}^*\omega_!}_-{\cong}
        &\text{ }
        \\
        \kos{f}^{\pr*}\omega^{\pr}_!
        \ar[d]_-{\kos{f}^{\pr*}\alpha_!}
        &\text{ }
        &\omega^{\ppr}_!\omega^{\ppr*}\kos{f}^{\pr*}\kos{f}^*\omega_!
        \ar[d]^-{\epsilon^{\ppr}\kos{f}^{\pr*}\kos{f}^*\omega_!}
        &\text{ }
        \\
        \kos{f}^{\pr*}\kos{f}^*\omega_!
        \ar@{=}[rr]
        &\text{ }
        &\kos{f}^{\pr*}\kos{f}^*\omega_!
        \ar@{=}[r]
        &\kos{f}^{\pr*}\kos{f}^*\omega_!
      }
    }}
  \end{equation*}
  Using the above description of $\beta_!$,
  we obtain the following relation.
  \begin{equation*}
    \vcenter
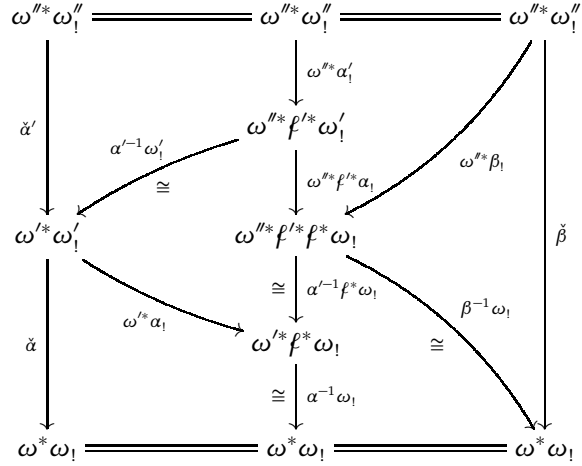
{\hbox{
      \xymatrix@C=50pt{
        \omega^{\ppr*}\omega^{\ppr}_!
        \ar[dd]_-{\check{\alpha}^{\pr}}
        \ar@{=}[r]
        &\omega^{\ppr*}\omega^{\ppr}_!
        \ar[d]^-{\omega^{\ppr*}\alpha^{\pr}_!}
        \ar@{=}[r]
        &\omega^{\ppr*}\omega^{\ppr}_!
        \ar[dddd]^-{\check{\beta}}
        \ar@/^1pc/[ddl]^-{\omega^{\ppr*}\beta_!}
        \\
        \text{ }
        &\omega^{\ppr*}\kos{f}^{\pr*}\omega^{\pr}_!
        \ar@/_0.5pc/[dl]_-{\alpha^{\pr-1}\omega^{\pr}_!}^-{\cong}
        \ar[d]^-{\omega^{\ppr*}\kos{f}^{\pr*}\alpha_!}
        &\text{ }
        \\
        \omega^{\pr*}\omega^{\pr}_!
        \ar[dd]_-{\check{\alpha}}
        \ar@/_0.5pc/[dr]_-{\omega^{\pr*}\alpha_!}
        &\omega^{\ppr*}\kos{f}^{\pr*}\kos{f}^*\omega_!
        \ar[d]^-{\alpha^{\pr-1}\kos{f}^*\omega_!}_-{\cong}
        \ar@/^1pc/[ddr]^-{\beta^{-1}\omega_!}_-{\cong}
        &\text{ }
        \\
        \text{ }
        &\omega^{\pr*}\kos{f}^*\omega_!
        \ar[d]^-{\alpha^{-1}\omega_!}_-{\cong}
        &\text{ }
        \\
        \omega^*\omega_!
        \ar@{=}[r]
        &\omega^*\omega_!
        \ar@{=}[r]
        &\omega^*\omega_!
      }
    }}
  \end{equation*}
  This shows that the correspondence
  $\mathbb{GALCAT}^{\cat{pre}}_*(\kos{K})
  \to \mathbb{GALOBJ}^{\cat{pre}}(\kos{K})$
  preserves composition of $1$-cells.

  Finally, we show that the correspondence
  $\mathbb{GALCAT}^{\cat{pre}}_*(\kos{K})
  \to \mathbb{GALOBJ}^{\cat{pre}}(\kos{K})$
  preserves horizontal composition of $2$-cells as follows.
  \begin{equation*}
    \vcenter
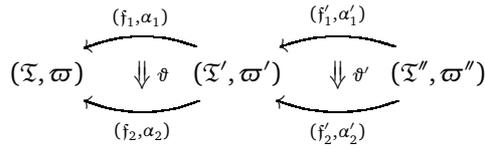
{\hbox{
      \xymatrix@C=35pt{
        (\mathfrak{T},\varpi)
        \xtwocell[r]{}<>{<0>{\text{ }\text{ }\vartheta}}
        &(\mathfrak{T}^{\pr},\varpi^{\pr})
        \ar@/_1pc/@<-0.5ex>[l]_-{(\mathfrak{f}_1,\alpha_1)}
        \ar@/^1pc/@<0.5ex>[l]^-{(\mathfrak{f}_2,\alpha_2)}
        \xtwocell[r]{}<>{<0>{\text{ }\text{ }\text{ }\vartheta^{\pr}}}
        &(\mathfrak{T}^{\ppr},\varpi^{\ppr})
        \ar@/_1pc/@<-0.5ex>[l]_-{(\mathfrak{f}^{\pr}_1,\alpha^{\pr}_1)}
        \ar@/^1pc/@<0.5ex>[l]^-{(\mathfrak{f}^{\pr}_2,\alpha^{\pr}_2)}
      }
    }}
  \end{equation*}
  \begin{equation*}
    \vartheta:\kos{f}_1^*\Rightarrow \kos{f}_2^*
    \quad
    \vartheta^{\pr}:\kos{f}^{\pr*}_1\Rightarrow \kos{f}^{\pr*}_2
    \quad
    \tilde{\vartheta}:=\vartheta^{\pr}\vartheta:\kos{f}^{\pr*}_1\kos{f}^*_1\Rightarrow \kos{f}^{\pr*}_2\kos{f}^*_2
  \end{equation*}
  We need to verify the relation
  $\check{\tilde{\vartheta}}=\check{\vartheta}\centerdot \check{\vartheta}^{\pr}$
  where
  \begin{equation*}
    \vcenter{\hbox{
      \xymatrix@C=30pt{
        \check{\vartheta}\centerdot \check{\vartheta}^{\pr}:
        I_{\kos{K}}
        \ar[r]^-{\check{\vartheta}\check{\vartheta}^{\pr}}
        &\omega^*\omega_!\omega^{\pr*}\omega^{\pr}_!
        \ar[r]^-{\omega^*\omega_!\check{\alpha}_1}
        &\omega^*\omega_!\omega^*\omega_!
        \ar[r]^-{\omega^{\pr}\epsilon\omega_!}
        &\omega^*\omega_!
        .
      }
    }}
  \end{equation*}
  We can check this as follows.
  We begin
  \begin{equation*}
    \vcenter{\hbox{
      \xymatrix@C=40pt{
        I_{\kos{K}}
        \ar[ddddd]_-{\check{\vartheta}\centerdot \check{\vartheta}^{\pr}}
        \ar@{=}[r]
        &I_{\kos{K}}
        \ar[d]_-{\check{\vartheta}^{\pr}}
        \ar@{=}[r]
        &I_{\kos{K}}
        \ar[d]^-{\check{\vartheta}^{\pr}}
        \\
        \text{ }
        &\omega^{\pr*}\omega^{\pr}_!
        \ar[d]_-{\check{\alpha}_1}
        \ar@{=}[r]
        &\omega^{\pr*}\omega^{\pr}_!
        \ar[dd]^-{\omega^{\pr*}\alpha_{1!}}
        \\
        \text{ }
        &\omega^*\omega_!
        \ar[dd]_-{\check{\vartheta}\omega^*\omega_!}
        \ar@/^0.5pc/[dr]^-{\alpha_1\omega_!}_-{\cong}
        &\text{ }
        \\
        \text{ }
        &\text{ }
        &\omega^{\pr*}\kos{f}^*_1\omega_!
        \ar[d]^-{\omega^{\pr*}\vartheta\omega_!}_-{\cong}
        \\
        \text{ }
        &\omega^*\omega_!\omega^*\omega_!
        \ar[d]_-{\omega^{\pr}\epsilon\omega_!}
        &\omega^{\pr*}\kos{f}^*_2\omega_!
        \ar[d]^-{\alpha_2^{-1}\omega_!}_-{\cong}
        \\
        \omega^*\omega_!
        \ar@{=}[r]
        &\omega^*\omega_!
        \ar@{=}[r]
        &\omega^*\omega_!
      }
    }}
  \end{equation*}
  and we finish.
  \begin{equation*}
    \vcenter{\hbox{
      \xymatrix@C=25pt{
        I_{\kos{K}}
        \ar[dddd]_-{\check{\vartheta}^{\pr}}
        \ar@{=}[rr]
        &\text{ }
        &I_{\kos{K}}
        \ar[d]^-{\eta^{\pr}}
        \ar@{=}[r]
        &I_{\kos{K}}
        \ar[d]^-{\eta}
        \ar@{=}[r]
        &I_{\kos{K}}
        \ar[ddddddd]^-{\check{\tilde{\theta}}}
        \\
        \text{ }
        &\text{ }
        &\omega^{\pr*}\omega^{\pr}_!
        \ar[d]^-{\alpha^{\pr}_1\omega^{\pr}_!}_-{\cong}
        \ar@/^0.5pc/[dr]|-{\omega^{\pr*}\alpha_{1!}}
        &\omega\omega_!
        \ar[d]^-{\alpha_1\omega_!}_-{\cong}
        &\text{ }
        \\
        \text{ }
        &\text{ }
        &\omega^{\ppr*}\kos{f}^{\pr*}_1\omega^{\pr}_!
        \ar@/_0.5pc/[dl]_-{\omega^{\ppr*}\vartheta^{\pr}\omega^{\pr}_!}^-{\cong}
        \ar[d]^-{\omega^{\ppr*}\kos{f}^{\pr*}_1\alpha_{1!}}
        &\omega^{\pr*}\kos{f}^*_1\omega_!
        \ar[d]^-{\alpha^{\pr}_1\kos{f}^*_1\omega_1}_-{\cong}
        &\text{ }
        \\
        \text{ }
        &\omega^{\ppr*}\kos{f}^{\pr*}_2\omega^{\pr}_!
        \ar@/_0.5pc/[dl]_-{\alpha^{\pr-1}_2\omega^{\pr}_!}^-{\cong}
        \ar[d]^-{\omega^{\ppr*}\kos{f}^{\pr*}_2\alpha_{1!}}
        &\omega^{\ppr*}\kos{f}^{\pr*}_1\kos{f}^*_1\omega_!
        \ar@/^0.5pc/[dl]^-{\omega^{\ppr*}\vartheta^{\pr}\kos{f}^*_1\omega_!}_-{\cong}
        \ar@/^1pc/@<1ex>[dd]^-{\omega^{\ppr*}\tilde{\vartheta}\omega_!}_-{\cong}
        \ar@{=}[r]
        &\omega^{\ppr*}\kos{f}^{\pr*}_1\kos{f}^*_1\omega_!
        &\text{ }
        \\
        \omega^{\pr*}\omega^{\pr}_!
        \ar[d]_-{\omega^{\pr*}\alpha_{1!}}
        &\omega^{\ppr*}\kos{f}^{\pr*}_2\kos{f}^*_1\omega_!
        \ar@/^0.5pc/[dl]_-{\alpha^{\pr-1}_2\kos{f}^*_1\omega_!}^-{\cong}
        \ar[d]^-{\omega^{\ppr*}\kos{f}^{\pr*}_2\vartheta\omega_!}_-{\cong}
        &\text{ }
        &\text{ }
        &\text{ }
        \\
        \omega^{\pr*}\kos{f}^*_1\omega_!
        \ar[d]_-{\omega^{\pr*}\vartheta\omega_!}^-{\cong}
        &\omega^{\ppr*}\kos{f}^{\pr*}_2\kos{f}^*_2\omega_!
        \ar[d]^-{\alpha^{\pr-1}_2\kos{f}^*_2\omega_!}_-{\cong}
        \ar@{=}[r]
        &\omega^{\ppr*}\kos{f}^{\pr*}_2\kos{f}^*_2\omega_!
        &\text{ }
        &\text{ }
        \\
        \omega^{\pr*}\kos{f}^*_2\omega_!
        \ar[d]_-{\alpha_2^{-1}\omega_!}^-{\cong}
        \ar@{=}[r]
        &\omega^{\pr*}\kos{f}^*_2\omega_!
        \ar[d]^-{\alpha_2^{-1}\omega_!}_-{\cong}
        &\text{ }
        &\text{ }
        &\text{ }
        \\
        \omega^*\omega_!
        \ar@{=}[r]
        &\omega^*\omega_!
        \ar@{=}[rrr]
        &\text{ }
        &\text{ }
        &\omega^*\omega_!
      }
    }}
  \end{equation*}
  This completes the proof of Proposition~\ref{prop CATtoOBJ}.
\qed\end{proof}

\begin{theorem} \label{thm GALOBJbiequivGALCAT}
  Let $\kos{K}$ be a Galois prekosmos $\kos{K}$.
  We have a biequivalence of $(2,1)$-categories
  \begin{equation*}
    \vcenter{\hbox{
      \xymatrix@C=30pt{
        \mathbb{GALOBJ}^{\cat{pre}}(\kos{K})
        \ar@<0.5ex>[r]^-{\simeq}
        &\mathbb{GALCAT}^{\cat{pre}}_*(\kos{K}).
        \ar@<0.5ex>[l]^-{\simeq}
      }
    }}
  \end{equation*}
\end{theorem}
\begin{proof}
  We are going to show that the $2$-functors described in 
  Proposition~\ref{prop inducedfundamentalfunctor}
  and
  Proposition~\ref{prop CATtoOBJ}
  are quasi-inverse to each other.
  The composition $2$-functor
  \begin{equation}\label{eq OBJ composition2funct}
    \mathbb{GALOBJ}^{\cat{pre}}(\kos{K})
    \to 
    \mathbb{GALCAT}^{\cat{pre}}_*(\kos{K})
    \to
    \mathbb{GALOBJ}^{\cat{pre}}(\kos{K})
  \end{equation}
  sends each $0$-cell $\pi$ to $\pi\otimes \kappa$,
  each $1$-cell $\pi^{\pr}\xrightarrow{f} \pi$ to 
  $\pi^{\pr}\otimes\kappa\xrightarrow{f\otimes I_{\kappa}} \pi\otimes \kappa$,
  and each $2$-cell $\theta:f_1\Rightarrow f_2:\pi^{\pr}\to \pi$
  to the $2$-cell 
  $\jmath_{\pi}\circ \theta:
  \xymatrix@C=15pt{
    \kappa
    \ar[r]^-{\theta}
    &\pi
    \ar[r]^-{\jmath_{\pi}}_-{\cong}
    &\pi\otimes \kappa
  }$
  between $1$-cells
  $f_1\otimes I_{\kappa}$,
  $f_2\otimes I_{\kappa}:\pi^{\pr}\otimes \kappa\to \pi\otimes \kappa$
  We claim that we have a $2$-natural isomorphism
  from the identity $2$-functor of 
  $\mathbb{NGALOBJ}(\kos{K})$
  to the composition $2$-functor (\ref{eq OBJ composition2funct}).
  The component of the $2$-natural isomorphism 
  at each $0$-cell $\pi$
  is the isomorphism $\jmath_{\pi}:\pi\cong \pi\otimes \kappa$.
  We need to show that the following diagram of functors strictly commute.
  Let $\pi^{\pr}$ be another $0$-cell.
  \begin{equation*}
    \vcenter{\hbox{
      \xymatrix@C=40pt{
        \mathbb{GALOBJ}^{\cat{pre}}(\kos{K})(\pi^{\pr},\pi)
        \ar[d]_-{(\slot)\otimes\kappa}
        \ar@{=}[r]
        &\mathbb{GALOBJ}^{\cat{pre}}(\kos{K})(\pi^{\pr},\pi)
        \ar[d]^-{\jmath_{\pi}\circ (-)}
        \\
        \mathbb{GALOBJ}^{\cat{pre}}(\kos{K})(\pi^{\pr}\otimes \kappa,\pi\otimes \kappa)
        \ar[r]^-{(-)\circ \jmath_{\pi^{\pr}}}
        &\mathbb{GALOBJ}^{\cat{pre}}(\kos{K})(\pi^{\pr},\pi\otimes \kappa)
      }
    }}
  \end{equation*}
  The two functors send
  each $1$-cell $f:\pi^{\pr}\to \pi$ to the same $1$-cell
  $\jmath_{\pi}\circ f=(f\otimes I_{\kappa})\circ \jmath_{\pi^{\pr}}:
  \pi^{\pr}\to \pi\otimes \kappa$
  and each $2$-cell $\sigma:f_1\Rightarrow f_2:\pi^{\pr}\to \pi$
  to the same $2$-cell
  \begin{equation*}
    \begin{aligned}
      \vcenter{\hbox{
        \xymatrix@C=20pt{
          \pi^{\pr}
          \ar[r]^-{\jmath_{\pi^{\pr}}}
          &\pi^{\pr}\otimes \kappa
          \ar@/^1pc/[rr]^-{f_1\otimes I_{\kappa}}
          \ar@/_1pc/[rr]_-{f_2\otimes I_{\kappa}}
          \xtwocell[rr]{}<>{<0>{\quad\jmath_{\pi}\circ \theta}}
          &\text{ }
          &\pi\otimes \kappa
        }
      }}
      \quad&=\quad
      \vcenter{\hbox{
        \xymatrix@C=30pt{
          \pi^{\pr}
          \ar@/^1pc/[r]^-{f_1}
          \ar@/_1pc/[r]_-{f_2}
          \xtwocell[r]{}<>{<0>{\text{ }\text{ }\theta}}
          &\pi
          \ar[r]^-{\jmath_{\pi}}
          &\pi\otimes \kappa
        }
      }}
      \\
      \bigl((f_2\otimes I_{\kappa})\circ\jmath_{\pi^{\pr}}\circ u_{\pi^{\pr}}\bigr)\star (\jmath_{\pi}\circ \theta)
      \quad&=\quad
      \jmath_{\pi}\circ \theta
    \end{aligned}
  \end{equation*}
  which we can see from the diagram below.
  \begin{equation*}
    \vcenter{\hbox{
      \xymatrix@C=40pt{
        \kappa
        \ar[ddd]|-{\bigl((f_2\otimes I_{\kappa})\circ\jmath_{\pi^{\pr}}\circ u_{\pi^{\pr}}\bigr)\star (\jmath_{\pi}\circ \theta)}
        \ar@{=}[r]
        &\kappa
        \ar[ddd]^-{(\jmath_{\pi}\circ u_{\pi})\star (\jmath_{\pi}\circ \theta)}
        \ar@{=}[r]
        &\kappa
        \ar[d]^-{u_{\pi}\otimes \theta}
        \ar@{=}[r]
        &\kappa
        \ar[dd]^-{\theta}
        \\
        \text{ }
        &\text{ }
        &\pi\otimes \pi
        \ar[d]^-{\jmath_{\pi}\otimes \jmath_{\pi}}_-{\cong}
        \ar@/^0.5pc/[dr]|-{\pc_{\pi}}
        &\text{ }
        \\
        \text{ }
        &\text{ }
        &(\pi\otimes \kappa)\otimes (\pi\otimes \kappa)
        \ar[d]^-{\pc_{\pi\otimes \kappa}}
        &\pi
        \ar[d]^-{\jmath_{\pi}}_-{\cong}
        \\
        \pi\otimes \kappa
        \ar@{=}[r]
        &\pi\otimes \kappa
        \ar@{=}[r]
        &\pi\otimes \kappa
        \ar@{=}[r]
        &\pi\otimes \kappa
      }
    }}
  \end{equation*}
  Thus we obtain a $2$-natural isomorphism
  from the identity $2$-functor of 
  $\mathbb{GALOBJ}^{\cat{pre}}(\kos{K})$
  to the composition $2$-functor (\ref{eq OBJ composition2funct}).

  Consider the composition $2$-functor
  \begin{equation}\label{eq CAT composition2funct}
    \mathbb{GALCAT}^{\cat{pre}}_*(\kos{K})\to
    \mathbb{GALOBJ}^{\cat{pre}}(\kos{K})\to 
    \mathbb{GALCAT}^{\cat{pre}}_*(\kos{K})
    .
  \end{equation}
  \begin{itemize}
    \item 
    It sends each $0$-cell $(\mathfrak{T},\varpi)$
    to the $0$-cell
    \begin{equation*}
      \bigl(\mathfrak{Rep}(\omega^*\omega_!(\kappa)),\varpi_{\omega^*\omega_!(\kappa)}\bigr)
      .
    \end{equation*}

    \item
    It sends 
    each $1$-cell
    $(\mathfrak{f},\alpha):(\mathfrak{T}^{\pr},\varpi^{\pr})\to (\mathfrak{T},\varpi)$
    to the $1$-cell
    \begin{equation*}
      (\mathfrak{f}(\check{\alpha}_{\kappa}),I):
      \bigl(\mathfrak{Rep}(\omega^{\pr*}\omega^{\pr}_!(\kappa)),\varpi_{\omega^{\pr*}\omega^{\pr}_!(\kappa)}\bigr)
      \to \bigl(\mathfrak{Rep}(\omega^*\omega_!(\kappa)),\varpi_{\omega^*\omega_!(\kappa)}\bigr)
    \end{equation*}
    where
    $\mathfrak{f}(\check{\alpha}_{\kappa}):
    \mathfrak{Rep}(\omega^{\pr*}\omega^{\pr}_!(\kappa))
    \to \mathfrak{Rep}(\omega^*\omega_!(\kappa))$
    is the fundamental functor
    obtained from the $1$-cell
    $\check{\alpha}_{\kappa}:\omega^{\pr*}\omega^{\pr}_!(\kappa)\to \omega^*\omega_!(\kappa)$.

    \item
    It sends each $2$-cell
    $\vartheta:(\mathfrak{f}_1,\alpha_1)\Rightarrow (\mathfrak{f}_2,\alpha_2)
    :(\mathfrak{T}^{\pr},\varpi^{\pr})\to (\mathfrak{T},\varpi)$
    to the $2$-cell
    \begin{equation*}
      \vartheta(\check{\vartheta}_{\kappa}):
      \bigl(\mathfrak{f}(\check{\alpha}_{1\kappa}),I\bigr)
      \Rightarrow
      \bigl(\mathfrak{f}(\check{\alpha}_{2\kappa}),I\bigr)
    \end{equation*}
    where 
    $\vartheta(\check{\vartheta}_{\kappa})$
    is the $2$-cell
    obtained from 
    $\check{\vartheta}_{\kappa}:
    \check{\alpha}_{1\kappa}
    \Rightarrow
    \check{\alpha}_{2\kappa}
    :\omega^{\pr*}\omega^{\pr}_!(\kappa)\to \omega^*\omega_!(\kappa)$.
  \end{itemize}
  We have the following
  weak $2$-natural transformation
  from the composition $2$-functor
  (\ref{eq CAT composition2funct})
  to the identity $2$-functor of
  $\mathbb{GALCAT}^{\cat{pre}}_*(\kos{K})$.
  By the term `weak $2$-natural transformation'
  we mean that we associate 
  to each $1$-cell
  $(\mathfrak{f},\alpha)$
  a constraint $2$-cell $\widebreve{\alpha}$
  which we explain below.
  \begin{itemize}
    \item 
    The component of the weak $2$-natural transformation
    at each $0$-cell $(\mathfrak{T},\varpi)$
    is the $1$-cell
    \begin{equation} \label{eq component of invertiblestrongtransformation}
      \vcenter{\hbox{
        \xymatrix@C=30pt{
          \mathfrak{T}
          &\mathfrak{Rep}(\omega^*\omega_!(\kappa))
          \ar[l]_-{\widebreve{\varpi}}^-{\simeq}
          \\
          \text{ }
          &\mathfrak{K}
          \ar@/^1pc/[ul]^-{\varpi}
          \ar[u]_-{\varpi_{\omega^*\omega_!(\kappa)}}
          \xtwocell[u]{}<>{<-5>{I}}
        }
      }}
      \quad
      \xymatrix@C=15pt{
        (\widebreve{\varpi},I):
        \bigl(\mathfrak{Rep}(\omega^*\omega_!(\kappa)),\varpi_{\omega^*\omega_!(\kappa)}\bigr)
        \ar[r]^-{\simeq}
        &(\mathfrak{T},\varpi)
      }
    \end{equation}
    where 
    $\widebreve{\varpi}:\mathfrak{Rep}(\omega^*\omega_!(\kappa))\xrightarrow{\simeq}\mathfrak{T}$
    is the equivalence of Galois $\kos{K}$-prekosmoi
    described in Theorem~\ref{thm preGalKCat mainThm1},
    and the invertible Galois $\kos{K}$-transformation
    $I:\varpi\Rightarrow\widebreve{\varpi}\varpi_{\omega^*\omega_!(\kappa)}$
    is the identity comonoidal $\kos{K}$-tensor natural isomorphism
    \begin{equation*}
      I:(\omega^*,\what{\omega^*})=(\omega^*_{\omega^*\omega_!(\kappa)}\widebreve{\omega}\!^*,\what{\omega^*_{\omega^*\omega_!(\kappa)}\widebreve{\omega}\!^*}).
    \end{equation*}
    
    \item
    Let 
    $(\mathfrak{f},\alpha):(\mathfrak{T}^{\pr},\varpi^{\pr})\to (\mathfrak{T},\varpi)$
    be a $1$-cell.
    Then
    the invertible Galois $\kos{K}$-transformation
    $\alpha:\widebreve{\varpi}\Rightarrow\mathfrak{f}^*\widebreve{\varpi}\!^{\pr}$
    becomes an invertible $2$-cell
    \begin{equation*}
      \vcenter{\hbox{
        \xymatrix@C=40pt{
          (\mathfrak{T}^{\pr},\varpi^{\pr})
          \ar[d]_-{(\mathfrak{f},\alpha)}
          &\bigl(\mathfrak{Rep}(\omega^{\pr*}\omega^{\pr}_!(\kappa)),\varpi_{\omega^{\pr*}\omega^{\pr}_!(\kappa)}\bigr)
          \ar[l]_-{(\widebreve{\varpi}^{\pr},I)}^-{\simeq}
          \ar[d]^-{(\mathfrak{f}(\check{\alpha}_{\kappa}),I)}
          \xtwocell[dl]{}<>{<0>{\widebreve{\alpha}\text{ }}}
          \\
          (\mathfrak{T},\varpi)
          &\bigl(\mathfrak{Rep}(\omega^*\omega_!(\kappa)),\varpi_{\omega^*\omega_!(\kappa)}\bigr)
          \ar[l]^-{(\widebreve{\varpi},I)}_-{\simeq}
        }
      }}
      \qquad
      \begin{aligned}
        \widebreve{\alpha}
        &:\!\!
        \xymatrix@C=15pt{
          \widebreve{\varpi}\mathfrak{f}(\check{\alpha}_{\kappa})
          \ar@2{->}[r]^-{\cong}
          &\mathfrak{f}\widebreve{\varpi}\!^{\pr}
        }
        \\
        \widebreve{\alpha}
        &:\!\!
        \xymatrix@C=15pt{
          \kos{f}(\check{\alpha}_{\kappa})^*\widebreve{\omega}\!^*
          \ar@2{->}[r]^-{\cong}
          &\widebreve{\omega}\!^{\pr*}\kos{f}^*
        }
      \end{aligned}
    \end{equation*}
    which means that for each object $X$ in $\CT$,
    the isomorphism
    $\alpha_X:\omega^*(X)\xrightarrow{\cong} \omega^{\pr*}\kos{f}^*(X)$
    is compatible left $\omega^{\pr*}\omega^{\pr}_!(\kappa)$-actions.
    We can check this as follows.
    \begin{equation*}
      \hspace*{-1cm}
      \vcenter{\hbox{
        \xymatrix@C=25pt{
          \omega^{\pr*}\omega^{\pr}_!(\kappa)\otimes \omega^*(X)
          \ar[d]_-{I_{\omega^{\pr*}\omega^{\pr}_!(\kappa)}\otimes \alpha_X}^-{\cong}
          \ar@{=}[r]
          &\omega^{\pr*}\omega^{\pr}_!(\kappa)\otimes \omega^*(X)
          \ar[d]^-{\bigl(\hatar{\omega^{\pr*}\omega^{\pr}_!}_{\omega^*(X)}\bigr)^{-1}}_-{\cong}
          \ar@{=}[r]
          &\omega^{\pr*}\omega^{\pr}_!(\kappa)\otimes \omega^*(X)
          \ar[d]^-{\check{\alpha}_{\kappa}\otimes I_{\omega^*(X)}}
          \\
          \omega^{\pr*}\omega^{\pr}_!(\kappa)\otimes \omega^{\pr*}\kos{f}^*(X)
          \ar[dd]_-{\bigl(\hatar{\omega^{\pr*}\omega^{\pr}_!}_{\omega^{\pr*}\kos{f}^*(X)}\bigr)^{-1}}^-{\cong}
          &\omega^{\pr*}\omega^{\pr}_!\omega^*(X)
          \ar@/^0.5pc/[ddl]_-{\omega^{\pr*}\omega^{\pr}_!(\alpha_X)}^-{\cong}
          \ar[dd]|-{\omega^{\pr*}((\alpha_!)_{!\omega^*(X)})}
          \ar@/^0.5pc/[dr]|-{\check{\alpha}_{\omega^*(X)}}
          &\omega^*\omega_!(\kappa)\otimes \omega^*(X)
          \ar[d]^-{\bigl(\hatar{\omega^*\omega_!}_{\omega^*(X)}\bigr)^{-1}}_-{\cong}
          \\
          \text{ }
          &\text{ }
          &\omega^*\omega_!\omega^*(X)
          \ar@/^0.5pc/[dl]^-{\alpha_{\omega_!\omega^*(X)}}_-{\cong}
          \ar[d]^-{\omega^*(\epsilon_X)}
          \\
          \omega^{\pr*}\omega^{\pr}_!\omega^{\pr*}\kos{f}^*(X)
          \ar[d]_-{\omega^{\pr*}(\epsilon^{\pr}_{\kos{f}^*(X)})}
          &\omega^{\pr*}\kos{f}^*\omega_!\omega^*(X)
          \ar[d]^-{\omega^{\pr*}\kos{f}^*(\epsilon_X)}
          &\omega^*(X)
          \ar[d]^-{\alpha_X}_-{\cong}
          \\
          \omega^{\pr*}\kos{f}^*(X)
          \ar@{=}[r]
          &\omega^{\pr*}\kos{f}^*(X)
          \ar@{=}[r]
          &\omega^{\pr*}\kos{f}^*(X)
        }
      }}
    \end{equation*}
  \end{itemize}
  To each $1$-cell $(\mathfrak{f},\alpha)$,
  we associate the constraint $2$-cell $\widebreve{\alpha}$ 
  which satisfies the following relations.
  \begin{itemize}
    \item 
    For each $2$-cell
    $\vartheta:(\mathfrak{f}_1,\alpha_1)\Rightarrow (\mathfrak{f}_2,\alpha_2)$,
    the following diagram commutes.
    \begin{equation*}
      \vcenter{\hbox{
        \xymatrix@C=40pt{
          \kos{f}(\check{\alpha}_{1\kappa})^*\widebreve{\omega}\!^*
          \ar@2{->}[d]_-{\vartheta(\check{\vartheta}_{\kappa})\widebreve{\omega}\!^*}^-{\cong}
          \ar@2{->}[r]^-{\widebreve{\alpha}_1}_-{\cong}
          &\widebreve{\omega}\!^{\pr*}\kos{f}^*_1
          \ar@2{->}[d]^-{\widebreve{\omega}\!^{\pr*}\vartheta}
          \\
          \kos{f}(\check{\alpha}_{2\kappa})^*\widebreve{\omega}\!^*
          \ar@2{->}[r]^-{\widebreve{\alpha}_2}_-{\cong}
          &\widebreve{\omega}\!^{\pr*}\kos{f}^*_2
        }
      }}
    \end{equation*}
    We can check this as follows.
    Let $X$ be an object in $\CT$.
    \begin{equation*}
      \vcenter{\hbox{
        \xymatrix{
          \omega^*(X)
          \ar[ddd]_-{\vartheta(\check{\vartheta}_{\kappa})_{\omega^*(X)}}
          \ar@{=}[r]
          &\omega^*(X)
          \ar[d]_-{\check{\vartheta}_{\kappa}\otimes I_{\omega^*(X)}}
          \ar@{=}[r]
          &\omega^*(X)
          \ar@/^1pc/[ddl]|-{\check{\vartheta}_{\omega^*(X)}}
          \ar[ddd]^-{\alpha_{1\omega^*(X)}}_-{\cong}
          \\
          \text{ }
          &\omega^*\omega_!(\kappa)\otimes \omega^*(X)
          \ar[d]_-{\bigl(\hatar{\omega^*\omega_!}_{\omega^*(X)}\bigr)}^-{\cong}
          &\text{ }
          \\
          \text{ }
          &\omega^*\omega_!\omega^*(X)
          \ar@/^0.5pc/[dl]^-{\omega^*(\epsilon_X)}
          \ar@{}[dd]|-{(\dagger)}
          &\text{ }
          \\
          \omega^*(X)
          \ar[d]_-{\alpha_{2\omega^*(X)}}^-{\cong}
          &\text{ }
          &\omega^{\pr*}\kos{f}^*_1(X)
          \ar[d]^-{\omega^{\pr*}\vartheta_X}_-{\cong}
          \\
          \omega^{\pr*}\kos{f}^*_2(X)
          \ar@{=}[rr]
          &\text{ }
          &\omega^{\pr*}\kos{f}^*_2(X)
        }
      }}
    \end{equation*}
    We used the relation
    (\ref{eq checkvartheta relation})
    in the diagram $(\dagger)$.
    
    \item
    The constraint 
    associated to the identity $1$-cell
    is the identity $2$-cell.
    \begin{equation*}
      \vcenter{\hbox{
        \xymatrix@C=40pt{
          (\mathfrak{T},\varpi)
          \ar[d]_-{(\id,I)}
          &\bigl(\mathfrak{Rep}(\omega^{\pr*}\omega^{\pr}_!(\kappa)),\varpi_{\omega^{\pr*}\omega^{\pr}_!(\kappa)}\bigr)
          \ar[l]_-{(\widebreve{\varpi},I)}^-{\simeq}
          \ar[d]^-{(\id,I)}
          \xtwocell[dl]{}<>{<0>{I\text{ }}}
          \\
          (\mathfrak{T},\varpi)
          &\bigl(\mathfrak{Rep}(\omega^*\omega_!(\kappa)),\varpi_{\omega^*\omega_!(\kappa)}\bigr)
          \ar[l]^-{(\widebreve{\varpi},I)}_-{\simeq}
        }
      }}
    \end{equation*}

    \item
    Whenever we have a composible pair of $1$-cells
    \begin{equation*}
      \vcenter{\hbox{
        \xymatrix@C=15pt{
          (\mathfrak{T},\varpi)
          &(\mathfrak{T}^{\pr},\varpi^{\pr})
          \ar[l]_-{(\mathfrak{f},\alpha)}
          &(\mathfrak{T}^{\ppr},\varpi^{\ppr})
          \ar[l]_-{(\mathfrak{f}^{\pr},\alpha^{\pr})}
          \ar@/^1pc/@<0.5ex>[ll]^-{(\mathfrak{f},\alpha)\circ(\mathfrak{f}^{\pr},\alpha^{\pr})=\bigl(\mathfrak{f}\mathfrak{f}^{\pr},(\mathfrak{f}\alpha^{\pr})\circ\alpha\bigr)}
        }
      }}
    \end{equation*}
    \begin{equation*}
      \begin{aligned}
        \beta:=
        (\mathfrak{f}\alpha^{\pr})\circ\alpha
        &:\!\!
        \vcenter{\hbox{
          \xymatrix@C=15pt{
            \varpi
            \ar@2{->}[r]^-{\alpha}_-{\cong}
            &\mathfrak{f}\varpi^{\pr}
            \ar@2{->}[r]^-{\mathfrak{f}\alpha^{\pr}}_-{\cong}
            &\mathfrak{f}\mathfrak{f}^{\pr}\varpi^{\ppr}
          }
        }}
        \\
        \beta:=
        (\alpha^{\pr}\kos{f}^*)\circ\alpha
        &:\!\!
        \vcenter{\hbox{
          \xymatrix@C=15pt{
            \omega^*
            \ar@2{->}[r]^-{\alpha}_-{\cong}
            &\omega^{\pr*}\kos{f}^*
            \ar@2{->}[r]^-{\alpha^{\pr}\kos{f}^*}_-{\cong}
            &\omega^{\ppr*}\kos{f}^{\pr*}\kos{f}^*
          }
        }}
      \end{aligned}
    \end{equation*}
    we have
    \begin{equation*}
      \vcenter{\hbox{
        \xymatrix@C=25pt{
          (\mathfrak{T}^{\ppr},\varpi^{\ppr})
          \ar[d]_-{(\mathfrak{f}^{\pr},\alpha^{\pr})}
          &\bigl(\mathfrak{Rep}(\omega^{\ppr*}\omega^{\ppr}_!(\kappa)),\varpi^{\ppr}\bigr)
          \ar[l]_-{(\widebreve{\varpi}\!^{\ppr},I)}^-{\simeq}
          \ar[d]^-{(\mathfrak{f}(\check{\alpha}^{\pr}_{\kappa}),I)}
          \xtwocell[dl]{}<>{<1>{\widebreve{\alpha}\!^{\pr}\text{ }\text{ }\text{ }}}
          \\
          (\mathfrak{T}^{\pr},\varpi^{\pr})
          \ar[d]_-{(\mathfrak{f},\alpha)}
          &\bigl(\mathfrak{Rep}(\omega^{\pr*}\omega^{\pr}_!(\kappa)),\varpi^{\pr}\bigr)
          \ar[l]_-{(\widebreve{\varpi}\!^{\pr},I)}^-{\simeq}
          \ar[d]^-{(\mathfrak{f}(\check{\alpha}_{\kappa}),I)}
          \xtwocell[dl]{}<>{<0>{\widebreve{\alpha}\text{ }}}
          \\
          (\mathfrak{T},\varpi)
          &\bigl(\mathfrak{Rep}(\omega^*\omega_!(\kappa)),\varpi\bigr)
          \ar[l]^-{(\widebreve{\varpi},I)}_-{\simeq}
        }
      }}
      \qquad\qquad
    \end{equation*}
    \begin{equation*}
      \qquad\qquad=
      \vcenter{\hbox{
        \xymatrix@C=25pt{
          (\mathfrak{T}^{\ppr},\varpi^{\ppr})
          \ar[d]_-{(\mathfrak{f}\mathfrak{f}^{\pr},\beta)}
          &\bigl(\mathfrak{Rep}(\omega^{\ppr*}\omega^{\ppr}_!(\kappa)),\varpi^{\ppr}\bigr)
          \ar[l]_-{(\widebreve{\varpi}\!^{\ppr},I)}^-{\simeq}
          \ar[d]^-{(\mathfrak{f}(\check{\beta}_{\kappa}),I)}
          \xtwocell[dl]{}<>{<0>{\widebreve{\beta}\text{ }}}
          \\
          (\mathfrak{T},\varpi)
          &\bigl(\mathfrak{Rep}(\omega^*\omega_!(\kappa)),\varpi\bigr)
          \ar[l]^-{(\widebreve{\varpi},I)}_-{\simeq}
        }
      }}
    \end{equation*}
    \begin{equation*}
      \widebreve{\beta}
      :\!\!
      \xymatrix@C=40pt{
        \kos{f}(\check{\alpha}^{\pr}_{\kappa})^*
        \kos{f}(\check{\alpha}_{\kappa})^*
        \widebreve{\omega}\!^*
        \ar@2{->}[r]^{\kos{f}(\check{\alpha}^{\pr}_{\kappa})^*\widebreve{\alpha}}_-{\cong}
        &\kos{f}(\check{\alpha}^{\pr}_{\kappa})^*
        \widebreve{\omega}\!^{\pr*}\kos{f}^*
        \ar@2{->}[r]^-{\widebreve{\alpha}\!^{\pr}\kos{f}^*}_-{\cong}
        &\widebreve{\omega}^{\ppr*}\kos{f}^{\pr*}\kos{f}^*
        .
      }
    \end{equation*}
    This simply follows from the definition
    $\beta:\!\!
    \xymatrix@C=15pt{
      \omega^*
      \ar@2{->}[r]^-{\alpha}_-{\cong}
      &\omega^{\pr*}\kos{f}^*
      \ar@2{->}[r]^-{\alpha^{\pr}\kos{f}^*}_-{\cong}
      &\omega^{\ppr*}\kos{f}^{\pr*}\kos{f}^*.
    }$
  \end{itemize}
  Moreover, this weak $2$-transformation
  from the composition $2$-functor
  (\ref{eq CAT composition2funct})
  to the identity $2$-functor on
  $\mathbb{GALCAT}^{\cat{pre}}_*(\kos{K})$
  is also invertible in weak sense.
  This means that every component $1$-cell has a quasi-inverse
  $1$-cell.
  The $1$-cell
  (\ref{eq component of invertiblestrongtransformation})
  has the following quasi-inverse $1$-cell
  \begin{equation*}
    \vcenter{\hbox{
      \xymatrix@C=25pt{
        \mathfrak{T}
        \ar[r]^-{\widetilde{\varpi}}_-{\simeq}
        &\mathfrak{Rep}(\omega^*\omega_!(\kappa))
        \xtwocell[d]{}<>{<5>{\widebreve{\eta}\varpi}}
        \\
        \text{ }
        &\mathfrak{K}
        \ar@/^1pc/[ul]^-{\varpi}
        \ar[u]_-{\varpi_{\omega^*\omega_!(\kappa)}}
      }
    }}
    \xymatrix@C=15pt{
      (\widetilde{\varpi},\widebreve{\eta}\varpi_{\omega^*\omega_!(\kappa)}):
      (\mathfrak{T},\varpi)
      \ar[r]^-{\simeq}
      &\bigl(\mathfrak{Rep}(\omega^*\omega_!(\kappa)),\varpi_{\omega^*\omega_!(\kappa)}\bigr)
    }
  \end{equation*}
  where $\widetilde{\varpi}:\mathfrak{T}\xrightarrow{\simeq}\mathfrak{Rep}(\omega^*\omega_!(\kappa))$
  is the opposite Galois $\kos{K}$-morphism
  of $\widebreve{\varpi}$
  defined as
  \begin{equation*}
    \vcenter{\hbox{
      \xymatrix{
        (\kos{R\!e\!p}(\omega^*\omega_!(\kappa)),\kos{t}^*_{\omega^*\omega_!(\kappa)})
        \ar@/^1pc/[d]^-{(\widetilde{\omega}^*,\what{\widetilde{\omega}}^*)=(\widebreve{\omega}\!_!,\what{\widebreve{\omega}}_!)}
        \\
        (\kos{T},\kos{t}^*)
        \ar@/^1pc/[u]^-{(\widebreve{\omega}\!^*,\what{\widebreve{\omega}}^*)=(\widetilde{\omega}_!,\what{\widetilde{\omega}}_!)}
      }
    }}
    \qquad
    \begin{aligned}
      \widetilde{\eta}
      &:\!\!
      \xymatrix@C=15pt{
        I_{\kos{T}}
        \ar@2{->}[r]^-{\widebreve{\epsilon}^{-1}}_-{\cong}
        &\widebreve{\omega}\!_!\widebreve{\omega}\!^*
        =\widetilde{\omega}^*\widetilde{\omega}_!
      }
      \\
      \widetilde{\epsilon}
      &:\!\!
      \xymatrix@C=15pt{
        \widetilde{\omega}_!\widetilde{\omega}^*
        =\widebreve{\omega}\!^*\widebreve{\omega}\!_!
        \ar@2{->}[r]^-{\widebreve{\eta}^{-1}}_-{\cong}
        &I_{\kos{R\!e\!p}(\omega^*\omega_!(\kappa))}
      }
    \end{aligned}
  \end{equation*}
  and the invertible Galois $\kos{K}$-transformation
  $\widebreve{\eta}\varpi_{\omega^*\omega_!(\kappa)}
  :\varpi_{\omega^*\omega_!(\kappa)}
  \Rightarrow\widetilde{\varpi}\varpi$
  is the comonoidal $\kos{K}$-tensor natural isomorphism
  \begin{equation*}
    \omega_{\omega^*\omega_!(\kappa)}^*\widebreve{\eta}
    :
    \xymatrix{
      \omega_{\omega^*\omega_!(\kappa)}^*
      \ar@2{->}[r]^-{\cong}
      &\omega_{\omega^*\omega_!(\kappa)}^*
      \widebreve{\omega}\!^*\widebreve{\omega}\!_!
      =
      \omega^*\widetilde{\omega}^*
      .
    }
  \end{equation*}
  We conclude that the $2$-functors defined in 
  Proposition~\ref{prop inducedfundamentalfunctor}
  and
  Proposition~\ref{prop CATtoOBJ}
  are quasi-inverse to each other.
  Thus we established the biequivalence between
  $\mathbb{GALOBJ}^{\cat{pre}}(\kos{K})$
  and
  $\mathbb{GALCAT}^{\cat{pre}}_*(\kos{K})$
  as we claimed.
  This completes the proof of Theorem~\ref{thm GALOBJbiequivGALCAT}.
\qed\end{proof}

\subsection{Torsors of a pre-Galois object}
\label{subsec Torsors GAL}
Let $\kos{K}=(\CK,\otimes,\kappa)$ be a Galois prekosmos.

\begin{definition}
  Let $\pi$ be a pre-Galois object in $\kos{K}$.
  A \emph{right $\pi$-torsor over $\kappa$}
  is a pair $(p,\lambda_p)$
  of an object $p$ in $\cat{Ens}(\kos{K})$
  and morphism $p\otimes \pi\xrightarrow{\lambda_p} p$ in $\cat{Ens}(\kos{K})$
  such that
  \begin{itemize}
    \item 
    $p\otimes\slot:\CK\to \CK$ is conservative and preserves reflexive coequalizers;
  
    \item
    $p\otimes \pi\xrightarrow{\lambda_p} p$ satisfies the right $\pi$-action relations;
  
    \item
    the composition
    $\tau_p:p\otimes \pi\xrightarrow[]{\cp_p\otimes I_{\pi}} p\otimes p\otimes \pi\xrightarrow[]{I_p\otimes \lambda_p} p\otimes p$
    is an isomorphism.
  \end{itemize}
  We define the \emph{division morphism} of $(p,\lambda_p)$ as follows.
  \begin{equation*}
    \mon{d}_p:
    \xymatrix{
      p\otimes p
      \ar[r]^-{\tau_p^{-1}}_-{\cong}
      &p\otimes \pi
      \ar[r]^-{e_p\otimes I_{\pi}}
      &\kappa\otimes \pi
      \ar[r]^-{\imath_{\pi}^{-1}}_-{\cong}
      &\pi
    }
  \end{equation*}
\end{definition}

We often omit the right $\pi$-action $\lambda_p$
and simply denote a right $\pi$-torsor $(p,\lambda_p)$ over $\kappa$
as $p$.
A \emph{morphism}
$p\to p^{\pr}$ of right $\pi$-torsors over $\kappa$
is a morphism in $\cat{Ens}(\kos{K})$
which is compatible with right $\pi$-actions
$\lambda_p$, $\lambda_{p^{\pr}}$.
We denote
$(\cat{Tors-}\pi)_{\kappa}$
as the category of right $\pi$-torsors over $\kappa$.

\begin{remark}
  Let $\pi$ be a pre-Galois object in $\kos{K}$.
  The pair $(\pi,\pc_{\pi})$ is a right $\pi$-torsor over $\kappa$.
  The composition 
  \begin{equation*}
    \varphi_{\pi}:=\tau_{\pi}:
    \pi\otimes \pi
    \xrightarrow{\cp_{\pi}\otimes I_{\pi}}
    \pi\otimes \pi\otimes \pi
    \xrightarrow{I_{\pi}\otimes \pc_{\pi}}
    \pi\otimes \pi
  \end{equation*}
  is an isomorphism
  whose inverse is
  $\tau_{\pi}^{-1}:
  \pi\otimes \pi
  \xrightarrow{\cp_{\pi}\otimes I_{\pi}}
  \pi\otimes \pi\otimes \pi
  \xrightarrow[\cong]{I_{\pi}\otimes\varsigma_{\pi}\otimes I_{\pi}}
  \pi\otimes \pi\otimes \pi
  \xrightarrow{I_{\pi}\otimes \pc_{\pi}}
  \pi\otimes \pi$
  The division morphism is given by
  $\mon{d}_{\pi}:
  \pi\otimes \pi
  \xrightarrow[\cong]{\varsigma_{\pi}\otimes I_{\pi}}
  \pi\otimes \pi
  \xrightarrow{\pc_{\pi}}
  \pi$.
\end{remark}

\begin{lemma} \label{lem Torsors tau relations}
  Let $\pi$ be a pre-Galois object in $\kos{K}$
  and let $p$ be a right $\pi$-torsor over $\kappa$.
  The isomorphism
  $\tau_p:p\otimes \pi\xrightarrow{\cong}p\otimes p$
  satisfies the following relations.
  \begin{equation*}
    \vcenter{\hbox{
      \xymatrix@C=30pt{
        p\otimes \pi
        \ar[r]^-{\cp_p\otimes I_{\pi}}
        \ar[d]_-{\tau_p}^-{\cong}
        &p\otimes p\otimes \pi
        \ar[d]^-{I_p\otimes \tau_p}_-{\cong}
        &p\otimes \pi
        \ar[d]_-{\tau_p}^-{\cong}
        \ar@/^1.2pc/[drr]^-{\lambda_p}
        &\text{ }
        &\text{ }
        \\
        p\otimes p
        \ar[r]^-{\cp_p\otimes I_p}
        &p\otimes p\otimes p
        &p\otimes p
        \ar[r]^-{e_p\otimes I_p}
        &\kappa\otimes p
        \ar[r]^-{\imath_p^{-1}}_-{\cong}
        &p
        \\
        p\otimes \pi\otimes \pi
        \ar[r]^-{I_p\otimes \pc_{\pi}}
        \ar[d]_-{\tau_p\otimes I_{\pi}}^-{\cong}
        &p\otimes \pi
        \ar[d]^-{\tau_p}_-{\cong}
        &p
        \ar[r]^-{\jmath_p}_-{\cong}
        \ar@/_1pc/[drr]_-{\cp_p}
        &p\otimes \kappa
        \ar[r]^-{I_p\otimes u_{\pi}}
        &p\otimes \pi
        \ar[d]^-{\tau_p}_-{\cong}
        \\
        p\otimes p\otimes \pi
        \ar[r]^-{I_p\otimes \lambda_p}
        &p\otimes p
        &\text{ }
        &\text{ }
        &p\otimes p
      }
    }}
  \end{equation*}
\end{lemma}
\begin{proof}
  We omit the symbol $\otimes$ and
  the coherence isomorphism $a$.
  \begin{equation*}
    \vcenter{\hbox{
      \xymatrix@C=35pt{
        p \pi
        \ar[dd]_-{\tau_p}^-{\cong}
        \ar@{=}[rr]
        &\text{ }
        &p \pi
        \ar[d]^-{\cp_p I_{\pi}}
        \ar@/_0.5pc/[dl]|-{\cp_p I_{\pi}}
        &p \pi
        \ar[dd]_-{\tau_p}^-{\cong}
        \ar@{=}[rr]
        &\text{ }
        &p \pi
        \ar[d]^-{\lambda_p}
        \ar@/_0.5pc/[dl]|-{\cp_p I_{\pi}}
        \ar@/^1pc/[ddl]^(0.6){\imath_p I_{\pi}}_-{\cong}
        \\
        \text{ }
        &p p \pi
        \ar@/_0.5pc/[dl]|-{I_p \lambda_p}
        \ar[d]|-{\cp_p I_{p \pi}}
        &p p \pi
        \ar@/^0.5pc/[dl]|-{I_p \cp_p I_{\pi}}
        \ar[dd]^-{I_p \tau_p}_-{\cong}
        &\text{ }
        &p p \pi
        \ar@/_0.5pc/[dl]|-{I_p \lambda_p}
        \ar[d]|-{e_p I_{p \pi}}
        &p
        \ar[dd]^(0.6){\imath_p}_(0.6){\cong}
        \\
        p p
        \ar[d]_-{\cp_p I_p}
        &p p p \pi
        \ar@/^0.5pc/[dl]|-{I_{p p} \lambda_p}
        &\text{ }
        &p p
        \ar[d]_-{e_p I_p}
        &\kappa p \pi
        \ar@/^0.5pc/[dl]|-{I_{\kappa} \lambda_p}
        &\text{ }
        \\
        p p p
        \ar@{=}[rr]
        &\text{ }
        &p p p
        &\kappa p
        \ar@{=}[rr]
        &\text{ }
        &\kappa p
      }
    }}
  \end{equation*}
  
  \begin{equation*}
    \vcenter{\hbox{
      \xymatrix@C=35pt{
        p \pi \pi
        \ar[dd]_-{\tau_p I_{\pi}}^-{\cong}
        \ar@{=}[rr]
        &\text{ }
        &p \pi \pi
        \ar@/_0.5pc/[dl]|-{\cp_p I_{\pi \pi}}
        \ar[d]^-{I_p \pc_{\pi}}
        &p \kappa
        \ar[d]_-{I_p u_{\pi}}
        \ar@{=}[rr]
        \ar@/^0.5pc/[dr]|-{\cp_p I_{\kappa}}
        &\text{ }
        &p \kappa
        \ar[dd]^(0.4){\jmath_p^{-1}}_(0.4){\cong}
        \\
        \text{ }
        &p p \pi \pi
        \ar@/_0.5pc/[dl]|-{I_p \lambda_p I_{\pi}}
        \ar[d]|-{I_{p p} \pc_{\pi}}
        &p \pi
        \ar@/^0.5pc/[dl]|-{\cp_p I_{\pi}}
        \ar[dd]^-{\tau_p}_-{\cong}
        &p \pi
        \ar[dd]_-{\tau_p}^-{\cong}
        \ar@/_0.5pc/[dr]|-{\cp_p I_{\pi}}
        &p p \kappa
        \ar[d]|-{I_{p p} u_{\pi}}
        \ar@<-0.3ex>@/^1pc/[ddr]^(0.4){I_p \jmath_p^{-1}}_-{\cong}
        &\text{ }
        \\
        p p \pi
        \ar[d]_-{I_p \lambda_p}
        &p p \pi
        \ar@/^0.5pc/[dl]|-{I_p \lambda_p}
        &\text{ }
        &\text{ }
        &p p \pi
        \ar@/_0.5pc/[dr]|-{I_p \lambda_p}
        &p
        \ar[d]^-{\cp_p}
        \\
        p p
        \ar@{=}[rr]
        &\text{ }
        &p p
        &p p
        \ar@{=}[rr]
        &\text{ }
        &p p
      }
    }}
  \end{equation*}
  This completes the proof of Lemma~\ref{lem Torsors tau relations}.
\qed\end{proof}

\begin{lemma} \label{lem Torsors division relations}
  Let $\pi$ be a pre-Galois object in $\kos{K}$
  and let $p$
  be a right $\pi$-torsor over $\kappa$.
  The division morphism
  $\mon{d}_p:p\otimes p\to \pi$
  satisfies the following relations.
  \begin{equation*}
    \vcenter{\hbox{
      \xymatrix@C=30pt{
        p\otimes p\otimes \pi
        \ar[d]_-{I_p\otimes \lambda_p}
        \ar[r]^-{\mon{d}_p\otimes I_{\pi}}
        &\pi\otimes \pi
        \ar[d]^-{\pc_{\pi}}
        &p
        \ar[r]^-{e_p}
        \ar[d]_-{\cp_p}
        &\kappa
        \ar[d]^-{u_{\pi}}
        &p\otimes p
        \ar[d]_-{\cp_p\otimes I_p}
        \ar@/^1.2pc/[dr]^-{(\tau_p)^{-1}}_-{\cong}
        &\text{ }
        \\
        p\otimes p
        \ar[r]^-{\mon{d}_p}
        &\pi
        &p\otimes p
        \ar[r]^-{\mon{d}_p}
        &\pi
        &p\otimes p\otimes p
        \ar[r]^-{I_p\otimes \mon{d}_p}
        &p\otimes \pi
      }
    }}
  \end{equation*}
\end{lemma}
\begin{proof}
  We omit the symbol $\otimes$
  and the coherence isomorphisms $a$, $\imath$, $\jmath$.
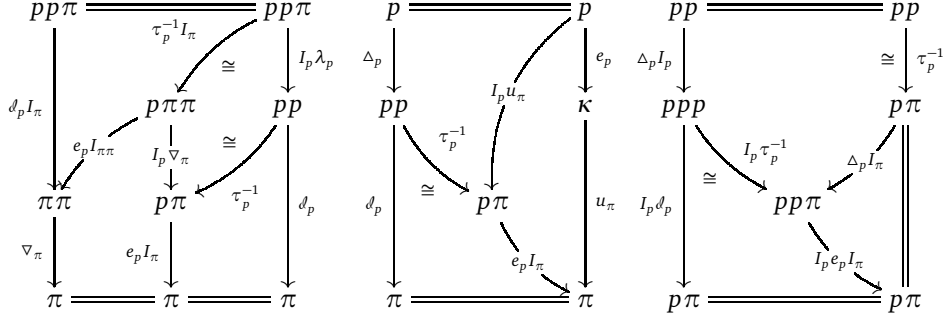
  \begin{equation*}
    \vcenter{\hbox{
      \xymatrix@C=19pt{
        pp\pi
        \ar@{=}[rr]
        \ar[dd]_-{\mon{d}_pI_{\pi}}
        &\text{ }
        &pp\pi
        \ar@/_0.5pc/[dl]_(0.6){\tau_p^{-1}I_{\pi}}^-{\cong}
        \ar[d]^-{I_p\lambda_p}
        &p
        \ar[d]_-{\cp_p}
        \ar@{=}[rr]
        &\text{ }
        &p
        \ar@/_1pc/[ddl]|-{I_pu_{\pi}}
        \ar[d]^-{e_p}
        &pp
        \ar[d]_-{\cp_pI_p}
        \ar@{=}[rr]
        &\text{ }
        &pp
        \ar[d]^-{\tau_p^{-1}}_-{\cong}
        \\
        \text{ }
        &p\pi\pi
        \ar[d]|-{I_p\pc_{\pi}}
        \ar@/_0.5pc/[dl]|-{e_pI_{\pi\pi}}
        &pp
        \ar@/^0.5pc/[dl]^(0.6){\tau_p^{-1}}_(0.45){\cong}
        \ar[dd]^-{\mon{d}_p}
        &pp
        \ar[dd]_-{\mon{d}_p}
        \ar@/_0.5pc/[dr]^(0.45){\tau_p^{-1}}_(0.6){\cong}
        &\text{ }
        &\kappa
        \ar[dd]^-{u_{\pi}}
        &ppp
        \ar@/_0.5pc/[dr]^(0.6){I_p\tau_p^{-1}}_-{\cong}
        \ar[dd]_-{I_p\mon{d}_p}
        &\text{ }
        &p\pi
        \ar@/^0.5pc/[dl]|-{\cp_pI_{\pi}}
        \ar@{=}[dd]
        \\
        \pi\pi
        \ar[d]_-{\pc_{\pi}}
        &p\pi
        \ar[d]_-{e_pI_{\pi}}
        &\text{ }
        &\text{ }
        &p\pi
        \ar@/_0.5pc/[dr]|-{e_pI_{\pi}}
        &\text{ }
        &\text{ }
        &pp\pi
        \ar@/_0.5pc/[dr]|-{I_pe_pI_{\pi}}
        &\text{ }
        \\
        \pi
        \ar@{=}[r]
        &\pi
        \ar@{=}[r]
        &\pi
        &\pi
        \ar@{=}[rr]
        &\text{ }
        &\pi
        &p\pi
        \ar@{=}[rr]
        &\text{ }
        &p\pi
      }
    }}
  \end{equation*}
  This completes the proof of Lemma~\ref{lem Torsors division relations}.
\qed\end{proof}

\begin{lemma} \label{lem Torsors Hom=Isom}
  Let $\pi$ be a pre-Galois object in $\kos{K}$.
  Every morphism of right $\pi$-torsors over $\kappa$
  is an isomorphism, i.e.,
  the category
  $(\cat{Tors-}\pi)_{\kappa}$
  of right $\pi$-torsors over $\kappa$
  is a groupoid.
\end{lemma}
\begin{proof}
  Let $f:p\to p^{\pr}$ be a morphism of right $\pi$-torsors over $\kappa$.
  Then we have
  \begin{equation*}
    \vcenter
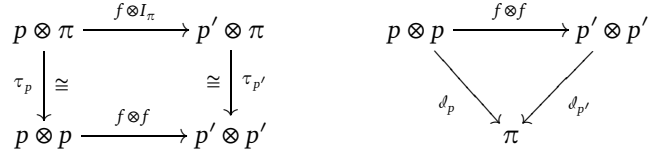
{\hbox{
      \xymatrix@C=15pt{
        p\otimes \pi
        \ar[rr]^-{f\otimes I_{\pi}}
        \ar[d]_-{\tau_p}^-{\cong}
        &\text{ }
        &p^{\pr}\otimes \pi
        \ar[d]^-{\tau_{p^{\pr}}}_-{\cong}
        &\text{ }
        &p\otimes p
        \ar[dr]_-{\mon{d}_p}
        \ar[rr]^-{f\otimes f}
        &\text{ }
        &p^{\pr}\otimes p^{\pr}
        \ar[dl]^-{\mon{d}_{p^{\pr}}}
        \\
        p\otimes p
        \ar[rr]^-{f\otimes f}
        &\text{ }
        &p^{\pr}\otimes p^{\pr}
        &\text{ }
        &\text{ }
        &\pi
        &\text{ }
      }
    }}
  \end{equation*}
  as we can see from the diagrams below.
  We omit the symbol $\otimes$ and coherence isomorphisms $a$, $\imath$, $\jmath$.
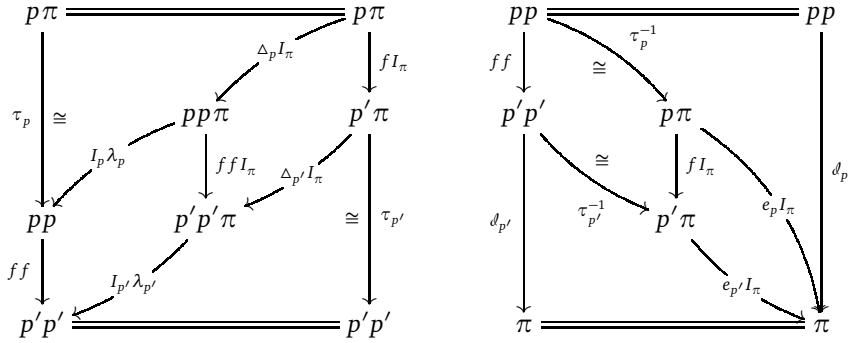
  \begin{equation*}
    \vcenter{\hbox{
      \xymatrix@C=35pt{
        p\pi
        \ar[dd]_-{\tau_p}^-{\cong}
        \ar@{=}[rr]
        &\text{ }
        &p\pi
        \ar@/_0.5pc/[dl]|-{\cp_pI_{\pi}}
        \ar[d]^-{fI_{\pi}}
        &pp
        \ar[d]_-{ff}
        \ar@/^0.5pc/[dr]^(0.6){\tau_p^{-1}}_-{\cong}
        \ar@{=}[rr]
        &\text{ }
        &pp
        \ar[ddd]^-{\mon{d}_p}
        \\
        \text{ }
        &pp\pi
        \ar@/_0.5pc/[dl]|-{I_p\lambda_p}
        \ar[d]^-{ffI_{\pi}}
        &p^{\pr}\pi
        \ar@/^0.5pc/[dl]|-{\cp_{p^{\pr}}I_{\pi}}
        \ar[dd]^-{\tau_{p^{\pr}}}_-{\cong}
        &p^{\pr}p^{\pr}
        \ar[dd]_-{\mon{d}_{p^{\pr}}}
        \ar@/_0.5pc/[dr]_(0.6){\tau_{p^{\pr}}^{-1}}^-{\cong}
        &p\pi
        \ar[d]^-{fI_{\pi}}
        \ar@/^1pc/[ddr]|-{e_pI_{\pi}}
        &\text{ }
        \\
        pp
        \ar[d]_-{ff}
        &p^{\pr}p^{\pr}\pi
        \ar@/^0.5pc/[dl]|-{I_{p^{\pr}}\lambda_{p^{\pr}}}
        &\text{ }
        &\text{ }
        &p^{\pr}\pi
        \ar@/_0.5pc/[dr]|-{e_{p^{\pr}}I_{\pi}}
        &\text{ }
        \\
        p^{\pr}p^{\pr}
        \ar@{=}[rr]
        &\text{ }
        &p^{\pr}p^{\pr}
        &\pi
        \ar@{=}[rr]
        &\text{ }
        &\pi
      }
    }}
  \end{equation*}
  To conclude that $f:p\to p^{\pr}$ is an isomorphism,
  it suffies to show that
  $I_p\otimes f:p\otimes p\to p\otimes p^{\pr}$
  is an isomorphism.
  This is because the functor $p\otimes\slot:\CK\to \CK$ is conservative.
  We claim that the inverse of $I_p\otimes f$ is given by
  \begin{equation*}
    (I_p\otimes f)^{-1}:
    \xymatrix@C=30pt{
      p\otimes p^{\pr}
      \ar[r]^-{\cp_p\otimes I_{p^{\pr}}}
      &p\otimes p\otimes p^{\pr}
      \ar[r]^-{I_p\otimes f\otimes I_{p^{\pr}}}
      &p\otimes p^{\pr}\otimes p^{\pr}
      \ar[r]^-{I_p\otimes \mon{d}_{p^{\pr}}}
      &p\otimes \pi
      \ar[r]^-{\tau_p}_-{\cong}
      &p\otimes p
      .
    }
  \end{equation*}
  We can check this as follows.
  \begin{equation*}
    \vcenter{\hbox{
      \xymatrix@C=30pt{
        pp
        \ar@{=}[rrrr]
        \ar[d]_-{I_pf}
        &\text{ }
        &\text{ }
        &\text{ }
        &pp
        \ar@/_0.5pc/[dll]|-{\cp_pI_p}
        \ar@{=}[ddddd]
        \ar@/_0.5pc/[ddddl]_(0.45){\tau_p^{-1}}^(0.45){\cong}
        \\
        pp^{\pr}
        \ar[dddd]_-{(I_pf)^{-1}}
        \ar@/_0.5pc/[dr]|-{\cp_pI_{p^{\pr}}}
        &\text{ }
        &ppp
        \ar@/_0.5pc/[dl]_-{I_{pp}f}
        \ar[dd]_-{I_pff}
        \ar@/^0.5pc/[dddr]|-{I_p\mon{d}_p}
        &\text{ }
        &\text{ }
        \\
        \text{ }
        &ppp^{\pr}
        \ar@/_0.5pc/[dr]_-{I_pfI_{p^{\pr}}}
        &\text{ }
        &\text{ }
        &\text{ }
        \\
        \text{ }
        &\text{ }
        &pp^{\pr}p^{\pr}
        \ar@/_0.5pc/[dr]_-{I_p\mon{d}_{p^{\pr}}}
        &\text{ }
        &\text{ }
        \\
        \text{ }
        &\text{ }
        &\text{ }
        &p\pi
        \ar@/_0.5pc/[dr]_-{\tau_p}^-{\cong}
        &\text{ }
        \\
        pp
        \ar@{=}[rrrr]
        &\text{ }
        &\text{ }
        &\text{ }
        &pp
      }
    }}
  \end{equation*}
  \begin{equation*}
    \vcenter{\hbox{
      \xymatrix@C=50pt{
        pp^{\pr}
        \ar@/^0.5pc/[dr]|-{\cp_pI_{p^{\pr}}}
        \ar[ddddd]_-{(I_pf)^{-1}}
        \ar@{=}[rrrr]
        &\text{ }
        &\text{ }
        &\text{ }
        &pp^{\pr}
        \ar@{=}[dddddd]
        \ar@/_0.5pc/[dl]|-{\cp_pI_{p^{\pr}}}
        \\
        \text{ }
        &ppp^{\pr}
        \ar[d]_-{I_pfI_{p^{\pr}}}
        \ar@/^0.5pc/[dr]|-{\cp_pI_{pp^{\pr}}}
        &\text{ }
        &ppp^{\pr}
        \ar@/_0.5pc/[dl]|-{I_p\cp_pI_{p^{\pr}}}
        \ar[dd]^-{I_pfI_{p^{\pr}}}
        \ar@/^1.5pc/[dddddr]|(0.35){I_pe_pI_{p^{\pr}}}
        &\text{ }
        \\
        \text{ }
        &pp^{\pr}p^{\pr}
        \ar[d]_-{I_p\mon{d}_{p^{\pr}}}
        \ar@/^0.5pc/[dr]|-{\cp_pI_{p^{\pr}p^{\pr}}}
        &pppp^{\pr}
        \ar[d]^-{I_{pp}fI_{p^{\pr}}}
        &\text{ }
        &\text{ }
        \\
        \text{ }
        &p\pi
        \ar@/_1pc/[ddl]_-{\tau_p}^-{\cong}
        \ar[d]_-{\cp_pI_{\pi}}
        &ppp^{\pr}p^{\pr}
        \ar[d]^-{I_pfI_{p^{\pr}p^{\pr}}}
        \ar@/^0.5pc/[dl]|-{I_{pp}\mon{d}_{p^{\pr}}}
        &pp^{\pr}p^{\pr}
        \ar@/^0.5pc/[dl]|-{I_p\cp_{p^{\pr}}I_{p^{\pr}}}
        \ar[dd]^-{I_p\tau_{p^{\pr}}^{-1}}_-{\cong}
        \ar@/^0.5pc/[dddr]|-{I_pe_{p^{\pr}}I_{p^{\pr}}}
        &\text{ }
        \\
        \text{ }
        &pp\pi
        \ar@/^0.5pc/[dl]^-{I_p\lambda_p}
        \ar@/_0.5pc/[dr]_-{I_pfI_{\pi}}
        &pp^{\pr}p^{\pr}p^{\pr}
        \ar[d]^-{I_{pp^{\pr}}\mon{d}_{p^{\pr}}}
        &\text{ }
        &\text{ }
        \\
        pp
        \ar[d]_-{I_pf}
        &\text{ }
        &pp^{\pr}\pi
        \ar@{=}[r]
        &pp^{\pr}\pi
        \ar@/_0.5pc/[dr]|-{I_p\lambda_{p^{\pr}}}
        &\text{ }
        \\
        pp^{\pr}
        \ar@{=}[rrrr]
        &\text{ }
        &\text{ }
        &\text{ }
        &pp^{\pr}
      }
    }}
  \end{equation*}
  This completes the proof of Lemma~\ref{lem Torsors Hom=Isom}.
\qed\end{proof}

Let $\pi$ be a pre-Galois object in $\kos{K}$.
A \emph{left $\pi$-torsor over $\kappa$}
is a pair $(p,\gamma_p)$
of an object $p$ in $\cat{Ens}(\kos{K})$
and a morphism $\gamma_p:\pi\otimes p\to p$ in $\cat{Ens}(\kos{K})$
such that
\begin{itemize}
  \item 
  $p\otimes\slot:\CK\to \CK$ is conservative and preserves reflexive coequalizers;

  \item 
  $\gamma_p:\pi\otimes p\to p$ satisfies the left $\pi$-action relations;

  \item 
  the composition
  $\pi\otimes p\xrightarrow[]{I_{\pi}\otimes \cp_p} \pi\otimes p\otimes p\xrightarrow[]{\gamma_p\otimes I_p} p\otimes p$
  is an isomorphism.
\end{itemize}
We often omit $\gamma_p$
and simply denote a left $\pi$-torsor $(p,\gamma_p)$ over $\kappa$
as $p$.

\begin{remark}
  Let $\pi$ be a pre-Galois object in $\kos{K}$.
  Suppose $p$ is a right $\pi$-torsor over $\kappa$ with right $\pi$-action $\lambda_p:p\otimes \pi\to p$.
  Then $p$ equipped with a left $\pi$-action
  \begin{equation*}
    \xymatrix{
      \gamma_p:
      \pi\otimes p
      \ar[r]^-{\varsigma_{\pi}\otimes I_p}_-{\cong}
      &\pi\otimes p
      \ar[r]^-{s_{\pi,p}}_-{\cong}
      &p\otimes \pi
      \ar[r]^-{\lambda_p}
      &p
    }
  \end{equation*}
  becomes a left $\pi$-torsor over $\kappa$.
  Conversely, if $p$ is a left $\pi$-torsor over $\kappa$
  with left $\pi$-action $\gamma_p:\pi\otimes p\to p$,
  then $p$ equipped with
  \begin{equation*}
    \lambda_p:
    \xymatrix{
      p\otimes \pi
      \ar[r]^-{I_p\otimes \varsigma_{\pi}}_-{\cong}
      &p\otimes \pi
      \ar[r]^-{s_{p,\pi}}_-{\cong}
      &\pi\otimes p
      \ar[r]^-{\gamma_p}
      &p
    }
  \end{equation*}
  becomes a right $\pi$-torsor over $\kappa$.
\end{remark}

Let $\pi^{\pr}$ be another pre-Galois object in $\kos{K}$.
A \emph{$(\pi^{\pr},\pi)$-bitorsor over $\kappa$}
is a triple $(p,\gamma^{\pr}_p,\lambda_p)$
where
$(p,\gamma^{\pr}_p)$
is a left $\pi^{\pr}$-torsor over $\kappa$,
$(p,\lambda_p)$
is a right $\pi$-torsor over $\kappa$
and 
the left $\pi^{\pr}$-action
$\gamma^{\pr}_p$
is compatible with
the right $\pi$-action
$\lambda_p$.
\begin{equation} \label{eq Torsors bitorsor actions}
  \vcenter{\hbox{
    \xymatrix@C=35pt{
      \pi^{\pr}\otimes (p\otimes \pi)
      \ar[d]_-{I_{\pi^{\pr}}\otimes \lambda_p}
      \ar[r]^-{a_{\pi^{\pr},p,\pi}}_-{\cong}
      &(\pi^{\pr}\otimes p)\otimes \pi
      \ar[r]^-{\gamma^{\pr}_p\otimes I_{\pi}}
      &p\otimes \pi
      \ar[d]^-{\lambda_p}
      \\
      \pi^{\pr}\otimes p
      \ar[rr]^-{\gamma^{\pr}_p}
      &\text{ }
      &p
    }
  }}
\end{equation}
We often omit
$\gamma^{\pr}_p$
as well as
$\lambda_p$
and simply denote a $(\pi^{\pr},\pi)$-torsor
$(p,\gamma^{\pr}_p,\lambda_p)$ over $\kappa$
as $p$.
The
\emph{opposite $(\pi,\pi^{\pr})$-bitorsor
of $p$ over $\kappa$}
is a triple 
$\mathring{p}=(p,\gamma_p,\lambda^{\pr}_p)$
whose left $\pi$-action $\gamma_p$
and right $\pi^{\pr}$-action $\lambda^{\pr}_p$ are
\begin{equation*}
  \begin{aligned}
    \gamma_p
    &:
    \xymatrix{
      \pi\otimes p
      \ar[r]^-{\varsigma_{\pi}\otimes I_p}_-{\cong}
      &\pi\otimes p
      \ar[r]^-{s_{\pi,p}}_-{\cong}
      &p\otimes \pi
      \ar[r]^-{\lambda_p}
      &p
      ,
    }
    \\
    \lambda^{\pr}_p
    &:
    \xymatrix{
      p\otimes \pi^{\pr}
      \ar[r]^-{I_p\otimes \varsigma_{\pi^{\pr}}}_-{\cong}
      &p\otimes \pi^{\pr}
      \ar[r]^-{s_{p,\pi^{\pr}}}_-{\cong}
      &\pi^{\pr}\otimes p
      \ar[r]^-{\gamma^{\pr}_p}
      &p
      .
    }
  \end{aligned}
\end{equation*}

\begin{lemma} \label{lem Torsors pip definition}
  Let $\pi$ be a pre-Galois object in $\kos{K}$
  and let $(p,\lambda_p)$ be a right $\pi$-torsor over $\kappa$.
  We define an object $\pi^p$ in $\cat{Ens}(\kos{K})$
  as the following reflexive coequalizer.
  \begin{equation} \label{eq Torsors pip definition}
    \vcenter{\hbox{
      \xymatrix@C=50pt{
        p\otimes p\otimes p
        \ar@<0.5ex>[rr]^-{I_{p\otimes p}\otimes e_p}
        \ar@<-0.5ex>[rr]_-{(\lambda_p\otimes I_p)\circ (I_p\otimes \mon{d}_p\otimes I_p)\circ (I_{p\otimes p}\otimes \cp_p)}
        &\text{ }
        &p\otimes p
        \ar@{->>}[r]^-{\textit{cq}_p}
        \ar@<-1ex>@/_2pc/[ll]|-{I_p\otimes \cp_p}
        &\pi^p
      }
    }}
  \end{equation}
  \begin{enumerate}
    \item 
    We have a unique morphism
    $\gamma^p_p:\pi^p\otimes p\to p$
    in $\cat{Ens}(\kos{K})$
    which satisfies the following relation.
    \begin{equation*}
      \vcenter{\hbox{
        \xymatrix@C=40pt{
          p\otimes p\otimes p
          \ar[d]_-{I_p\otimes \mon{d}_p}
          \ar@{->>}[r]^-{\textit{cq}_p\otimes I_p}
          &\pi^p\otimes p
          \ar@{.>}[d]^-{\gamma^p_p}_-{\exists!}
          \\
          p\otimes \pi
          \ar[r]^-{\lambda_p}
          &p
        }
      }}
    \end{equation*}
    Moreover,
    $\gamma^p_p$
    is compatible with
    $\lambda_p$
    as described in (\ref{eq Torsors bitorsor actions}).

    \item 
    We have an isomorphism
    $\Phi:p\otimes p \xrightarrow{\cong}p\otimes \pi^p$
    in $\cat{Ens}(\kos{K})$
    defined by
    \begin{equation*}
      \Phi:
      \xymatrix@C=30pt{
        p\otimes p
        \ar[r]^-{\cp_p\otimes I_p}
        &p\otimes p\otimes p
        \ar[r]^-{I_p\otimes s_{p,p}}_-{\cong}
        &p\otimes p\otimes p
        \ar@{->>}[r]^-{I_p\otimes \textit{cq}_p}
        &p\otimes \pi^p
      }
    \end{equation*}
    whose inverse is 
    \begin{equation*}
      \Phi^{-1}:
      \xymatrix@C=30pt{
        p\otimes \pi^p
        \ar[r]^-{\cp_p\otimes I_{\pi^p}}
        &p\otimes p\otimes \pi^p
        \ar[r]^-{I_p\otimes s_{p,\pi^p}}_-{\cong}
        &p\otimes \pi^p\otimes p
        \ar[r]^-{I_p\otimes \gamma^p_p}
        &p\otimes p
        .
      }
    \end{equation*}
    Moreover, the isomorphism
    $\Phi$
    satisfies the following relation.
    \begin{equation} \label{eq2 Torsors pip definition}
      \vcenter{\hbox{
        \xymatrix@C=30pt{
          p\otimes p
          \ar[d]_-{\Phi}^-{\cong}
          \ar[r]^-{\cp_{p\otimes p}}
          &p\otimes p\otimes p\otimes p
          \ar[r]^-{I_{p\otimes p}\otimes s_{p,p}}_-{\cong}
          &p\otimes p\otimes p\otimes p
          \ar@{->>}[r]^-{I_{p\otimes p}\otimes \textit{cq}_p}
          &p\otimes p\otimes \pi^p
          \ar[d]^-{\Phi\otimes I_{\pi^p}}_-{\cong}
          \\
          p\otimes \pi^p
          \ar[rrr]^-{I_p\otimes \cp_{\pi^p}}
          &\text{ }
          &\text{ }
          &p\otimes \pi^p\otimes \pi^p
        }
      }}
    \end{equation}
  \end{enumerate}
\end{lemma}
\begin{proof}
  We omit the symbol $\otimes$ and the
  coherence isomorphisms $a$, $\imath$, $\jmath$
  when presenting the details.
  First we show that the object $\pi^p$ in $\cat{Ens}(\kos{K})$ is well-defined.
  We need to show that 
  $I_p\otimes \cp_p:p\otimes p\to p\otimes p\otimes p$
  is a common section of the parallel pair of morphisms in (\ref{eq Torsors pip definition}).
  The morphism $I_p\otimes \cp_p$ is obviously a section of $I_{p\otimes p}\otimes e_p:p\otimes p\otimes p\to p\otimes p$.
  We show that $I_p\otimes \cp_p$ is also a section of the other morphism as follows.
  \begin{equation*}
    \vcenter{\hbox{
      \xymatrix@C=50pt{
        pp
        \ar[d]_-{I_p\cp_p}
        \ar@{=}[rr]
        &\text{ }
        &pp
        \ar@{=}[dd]
        \ar@/_0.5pc/[dl]|-{I_p\cp_p}
        \\
        ppp
        \ar[d]_-{I_{pp}\cp_p}
        &ppp
        \ar@/_0.5pc/[dl]|-{I_p\cp_pI_p}
        \ar[d]^{I_pe_pI_p}
        &\text{ }
        \\
        pppp
        \ar[d]_-{I_p\mon{d}_pI_p}
        &pp
        \ar@/^0.5pc/[dl]|-{I_pu_{\pi}I_p}
        \ar@{=}[r]
        &pp
        \ar@{=}[dd]
        \\
        p\pi p
        \ar[d]_-{\lambda_pI_p}
        &\text{ }
        &\text{ }
        \\
        pp
        \ar@{=}[rr]
        &\text{ }
        &pp
      }
    }}
  \end{equation*}
  This shows that the object $\pi^p$ in $\cat{Ens}(\kos{K})$ is well-defined.
  Next we show that the morphism $\gamma^p_p:\pi^p\otimes p\to p$ is well-defined.
  As the functor $\slot\otimes p:\CK\to \CK$ preserves 
  the reflexive coequalizer (\ref{eq Torsors pip definition}),
  we need to check that the morphism
  $\lambda_p\circ (I_p\otimes \mon{d}_p):p\otimes p\otimes p\to p$
  coequalizes the parallel pair of morphisms
  $I_{p\otimes p}\otimes e_p\otimes I_p$,
  $(\lambda_p\otimes I_{p\otimes p})\circ (I_p\otimes \mon{d}_p\otimes I_{p\otimes p})\circ (I_{p\otimes p}\otimes \cp_p\otimes I_p)
  :p\otimes p\otimes p\otimes p\to p\otimes p\otimes p$.
  We can check this as follows.
  \begin{equation*}
    \vcenter{\hbox{
      \xymatrix@C=70pt{
        pppp
        \ar[d]_-{I_{pp}\cp_pI_p}
        \ar@{=}[rr]
        &\text{ }
        &pppp
        \ar@/_1pc/[ddl]_-{I_{pp}\tau_p^{-1}}^-{\cong}
        \ar[ddd]^-{I_{pp}e_pI_p}
        \\
        ppppp
        \ar[d]_-{I_p\mon{d}_pI_{pp}}
        \ar@/^0.5pc/[dr]|-{I_{ppp}\mon{d}_p}
        &\text{ }
        &\text{ }
        \\
        p\pi pp
        \ar[d]_-{\lambda_pI_{pp}}
        \ar@/_0.5pc/[dr]|-{I_{p\pi}\mon{d}_p}
        &ppp\pi
        \ar@/^0.5pc/[dr]|-{I_{pp}\lambda_p}
        \ar[d]^-{I_p\mon{d}_pI_{\pi}}
        &\text{ }
        \\
        ppp
        \ar[d]_-{I_p\mon{d}_p}
        &p\pi\pi
        \ar@/^0.5pc/[dl]|-{\lambda_pI_{\pi}}
        \ar@/_0.5pc/[dr]|-{I_p\pc_{\pi}}
        &ppp
        \ar[d]^-{I_p\mon{d}_p}
        \\
        p\pi
        \ar[d]_-{\lambda_p}
        &\text{ }
        &p\pi
        \ar[d]^-{\lambda_p}
        \\
        p
        \ar@{=}[rr]
        &\text{ }
        &p
      }
    }}
  \end{equation*}
  This shows that the morphism $\gamma^p_p:\pi^p\otimes p\to p$ is well-defined.
  We also obtain that $\gamma^p_p$ and $\lambda_p$ are compatible,
  by right-cancelling the epimorphism
  $\textit{cq}_p\otimes I_{p\otimes \pi}$
  in the diagram below.
  \begin{equation*}
    \vcenter{\hbox{
      \xymatrix@C=30pt{
        ppp\pi
        \ar@{->>}[d]_-{\textit{cq}_pI_{p\pi}}
        \ar@{=}[rr]
        &\text{ }
        \ar@{=}[r]
        &ppp\pi
        \ar@/_0.5pc/[dl]|-{I_p\mon{d}_pI_{\pi}}
        \ar[d]^-{I_{pp}\lambda_p}
        \ar@{=}[r]
        &ppp\pi
        \ar@{->>}[d]^-{\textit{cq}_pI_{p\pi}}
        \\
        \pi^p p\pi
        \ar[d]_-{\gamma^p_pI_{\pi}}
        &p\pi\pi
        \ar@/_0.5pc/[dl]|-{\lambda_pI_{\pi}}
        \ar[d]^-{I_p\pc_{\pi}}
        &ppp
        \ar@/^0.5pc/[dl]|-{I_p\mon{d}_p}
        \ar@{->>}@/_0.5pc/[dr]|-{\textit{cq}_pI_p}
        &\pi^p p\pi
        \ar[d]^-{I_{\pi^p}\lambda_p}
        \\
        p\pi
        \ar[d]_-{\lambda_p}
        &p\pi
        \ar@/^0.5pc/[dl]|-{\lambda_p}
        &\text{ }
        &\pi^p p
        \ar[d]^-{\gamma^p_pI_p}
        \\
        p
        \ar@{=}[rrr]
        &\text{ }
        &\text{ }
        &p
      }
    }}
  \end{equation*}
  Next we show that
  $\Phi:p\otimes p\to p\otimes \pi^p$
  and 
  $\Phi^{-1}:p\otimes \pi^p\to p\otimes p$
  are inverse to each other.
  We verify the relation
  $\Phi^{-1}\circ \Phi=I_{p\otimes p}$ as follows.
  \begin{equation*}
    \vcenter
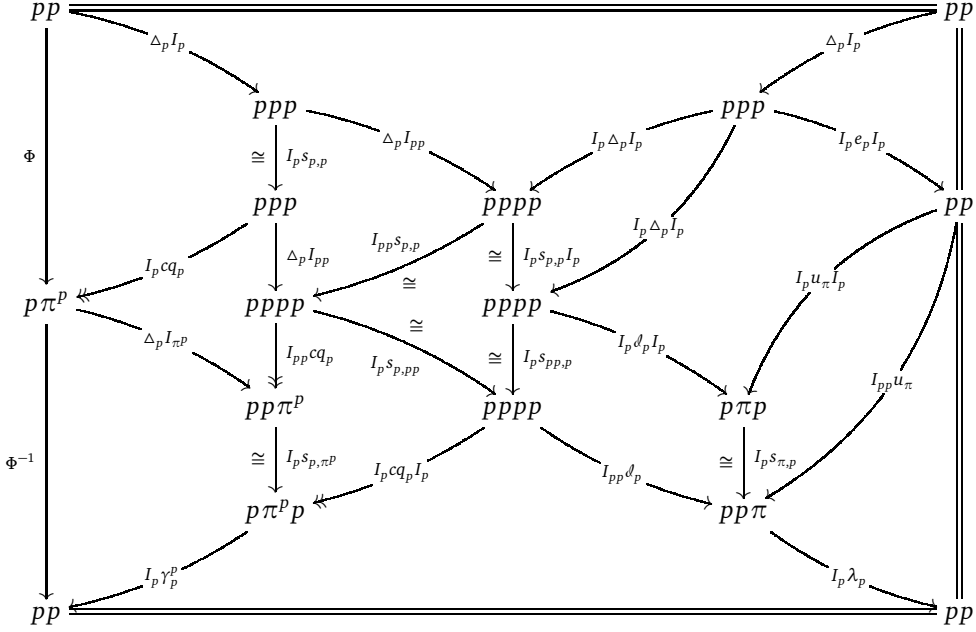
{\hbox{
      \xymatrix@C=60pt{
        pp
        \ar[ddd]_-{\Phi}
        \ar@/^0.5pc/[dr]|-{\cp_pI_p}
        \ar@{=}[rrrr]
        &\text{ }
        &\text{ }
        &\text{ }
        &pp
        \ar@/_0.5pc/[dl]|-{\cp_pI_p}
        \ar@{=}[dd]
        \\
        \text{ }
        &ppp
        \ar[d]^-{I_ps_{p,p}}_-{\cong}
        \ar@/^0.5pc/[dr]|-{\cp_pI_{pp}}
        &\text{ }
        &ppp
        \ar@/_0.5pc/[dl]|-{I_p\cp_pI_p}
        \ar@/^1.2pc/[ddl]|-{I_p\cp_pI_p}
        \ar@/^0.5pc/[dr]|-{I_pe_pI_p}
        &\text{ }
        \\
        \text{ }
        &ppp
        \ar@{->>}@/^0.5pc/[dl]|-{I_p\textit{cq}_p}
        \ar[d]^-{\cp_pI_{pp}}
        &pppp
        \ar@/^0.5pc/[dl]_(0.45){I_{pp}s_{p,p}}^-{\cong}
        \ar[d]^-{I_ps_{p,p}I_p}_-{\cong}
        &\text{ }
        &pp
        \ar@/_1.2pc/[ddl]|-{I_pu_{\pi}I_p}
        \ar@/^1.5pc/[dddl]|-{I_{pp}u_{\pi}}
        \ar@{=}[dddd]
        \\
        p\pi^p
        \ar[ddd]_-{\Phi^{-1}}
        \ar@/^0.5pc/[dr]|-{\cp_pI_{\pi^p}}
        &pppp
        \ar@{->>}[d]^-{I_{pp}\textit{cq}_p}
        \ar@/^0.5pc/[dr]_(0.55){I_ps_{p,pp}}^-{\cong}
        &pppp
        \ar[d]^-{I_ps_{pp,p}}_-{\cong}
        \ar@/^0.5pc/[dr]|-{I_p\mon{d}_pI_p}
        &\text{ }
        &\text{ }
        \\
        \text{ }
        &pp\pi^p
        \ar[d]^-{I_ps_{p,\pi^p}}_-{\cong}
        &pppp
        \ar@{->>}@/^0.5pc/[dl]|-{I_p\textit{cq}_pI_p}
        \ar@/_0.5pc/[dr]|-{I_{pp}\mon{d}_p}
        &p\pi p
        \ar[d]^{I_ps_{\pi,p}}_-{\cong}
        &\text{ }
        \\
        \text{ }
        &p\pi^p p
        \ar@/^0.5pc/[dl]|-{I_p\gamma^p_p}
        &\text{ }
        &pp\pi
        \ar@/_0.5pc/[dr]|-{I_p\lambda_p}
        &\text{ }
        \\
        pp
        \ar@{=}[rrrr]
        &\text{ }
        &\text{ }
        &\text{ }
        &pp
      }
    }}
  \end{equation*}
  We obtain the other relation
  $\Phi\circ \Phi^{-1}=I_{p\otimes \pi^p}$
  by right-cancelling the epimorphism
  $I_p\otimes \textit{cq}_p$
  in the following diagram,
  where we postpone the proof of the diagram $(\dagger)$.
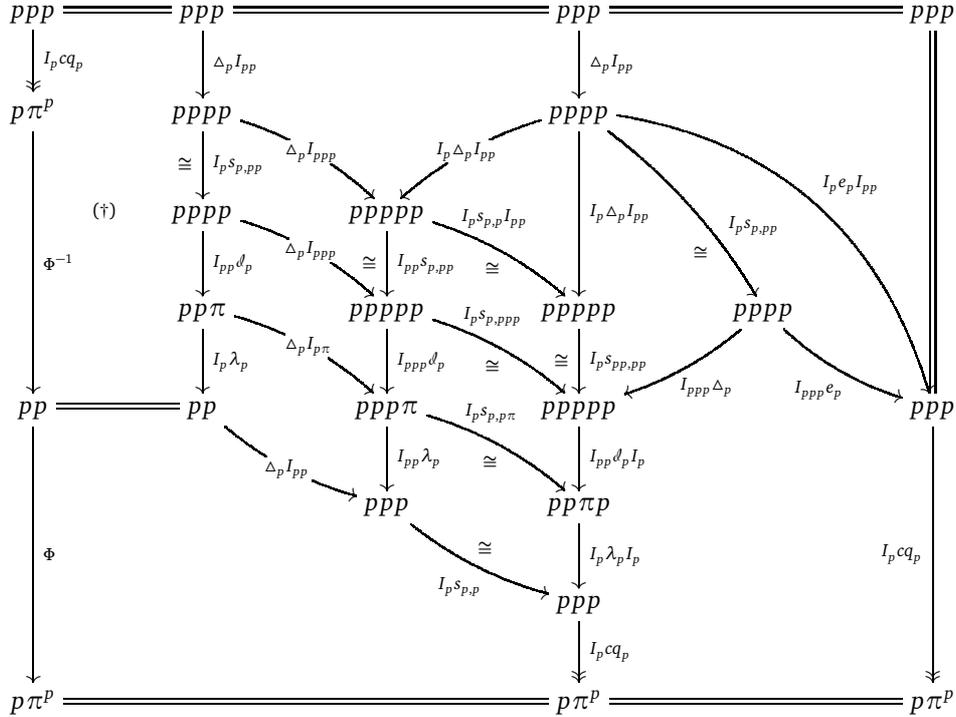
  \begin{equation*}
    \vcenter{\hbox{
      \xymatrix@C=37.5pt{
        ppp
        \ar@{->>}[d]^-{I_p\textit{cq}_p}
        \ar@{=}[r]
        &ppp
        \ar[d]^-{\cp_pI_{pp}}
        \ar@{=}[rr]
        &\text{ }
        &ppp
        \ar[d]^-{\cp_pI_{pp}}
        \ar@{=}[rr]
        &\text{ }
        &ppp
        \ar@{=}[dddd]
        \\
        p\pi^p
        \ar[ddd]^-{\Phi^{-1}}
        &pppp
        \ar[d]^-{I_ps_{p,pp}}_-{\cong}
        \ar@/^0.5pc/[dr]|-{\cp_pI_{ppp}}
        &\text{ }
        &pppp
        \ar@/_0.5pc/[dl]|-{I_p\cp_pI_{pp}}
        \ar[dd]^-{I_p\cp_pI_{pp}}
        \ar@/^0.5pc/[ddr]^(0.7){I_ps_{p,pp}}_(0.7){\cong}
        \ar@/^2.5pc/[dddrr]^-{I_pe_pI_{pp}}
        &\text{ }
        &\text{ }
        \\
        \text{ }
        \ar@{}[r]|-{(\dagger)}
        &pppp
        \ar[d]^-{I_{pp}\mon{d}_p}
        \ar@/^0.5pc/[dr]|-{\cp_pI_{ppp}}
        &ppppp
        \ar[d]^-{I_{pp}s_{p,pp}}_-{\cong}
        \ar@/^0.5pc/[dr]^(0.42){I_ps_{p,p}I_{pp}}_-{\cong}
        &\text{ }
        &\text{ }
        &\text{ }
        \\
        \text{ }
        &pp\pi
        \ar[d]^-{I_p\lambda_p}
        \ar@/^0.5pc/[dr]|-{\cp_pI_{p\pi}}
        &ppppp
        \ar[d]^-{I_{ppp}\mon{d}_p}
        \ar@/^0.5pc/[dr]^(0.42){I_ps_{p,ppp}}_-{\cong}
        &ppppp
        \ar[d]^-{I_ps_{pp,pp}}_-{\cong}
        &pppp
        \ar@/^0.5pc/[dl]^-{I_{ppp}\cp_p}
        \ar@/_0.5pc/[dr]_-{I_{ppp}e_p}
        &\text{ }
        \\
        pp
        \ar[ddd]^-{\Phi}
        \ar@{=}[r]
        &pp
        \ar@/_0.5pc/[dr]|-{\cp_pI_{pp}}
        &ppp\pi
        \ar[d]^-{I_{pp}\lambda_p}
        \ar@/^0.5pc/[dr]^(0.42){I_ps_{p,p\pi}}_-{\cong}
        &ppppp
        \ar[d]^-{I_{pp}\mon{d}_pI_p}
        &\text{ }
        &ppp
        \ar@{->>}[ddd]_-{I_p\textit{cq}_p}
        \\
        \text{ }
        &\text{ }
        &ppp
        \ar@/_0.5pc/[dr]_-{I_ps_{p,p}}^-{\cong}
        &pp\pi p
        \ar[d]^-{I_p\lambda_pI_p}
        \\
        \text{ }
        &\text{ }
        &\text{ }
        &ppp
        \ar@{->>}[d]^-{I_p\textit{cq}_p}
        &\text{ }
        &\text{ }
        \\
        p\pi^p
        \ar@{=}[rrr]
        &\text{ }
        &\text{ }
        &p\pi^p
        \ar@{=}[rr]
        &\text{ }
        &p\pi^p
      }
    }}
  \end{equation*}
  We verify the diagram $(\dagger)$ as follows.
  \begin{equation*}
    (\dagger):
    \vcenter
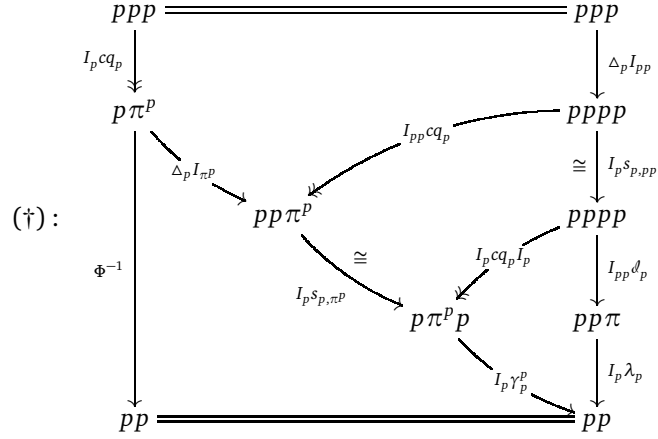
{\hbox{
      \xymatrix@C=30pt{
        ppp
        \ar@{->>}[d]_-{I_p\textit{cq}_p}
        \ar@{=}[rrr]
        &\text{ }
        &\text{ }
        &ppp
        \ar[d]^-{\cp_pI_{pp}}
        \\
        p\pi^p
        \ar[ddd]_-{\Phi^{-1}}
        \ar@/_0.5pc/[dr]|-{\cp_pI_{\pi^p}}
        &\text{ }
        &\text{ }
        &pppp
        \ar[d]^-{I_ps_{p,pp}}_-{\cong}
        \ar@{->>}@/_1pc/[dll]|-{I_{pp}\textit{cq}_p}
        \\
        \text{ }
        &pp\pi^p
        \ar@/_0.5pc/[dr]_-{I_ps_{p,\pi^p}}^-{\cong}
        &\text{ }
        &pppp
        \ar[d]^-{I_{pp}\mon{d}_p}
        \ar@{->>}@/_0.5pc/[dl]|-{I_p\textit{cq}_pI_p}
        \\
        \text{ }
        &\text{ }
        &p\pi^p p
        \ar@/_0.5pc/[dr]|-{I_p\gamma^p_p}
        &pp\pi
        \ar[d]^-{I_p\lambda_p}
        \\
        pp
        \ar@{=}[rrr]
        &\text{ }
        &\text{ }
        &pp
      }
    }}
  \end{equation*}
  This shows that $\Phi$ and $\Phi^{-1}$ are inverse to each other.
  Finally, we verify the relation (\ref{eq2 Torsors pip definition}) as follows.
  \begin{equation*}
    \vcenter
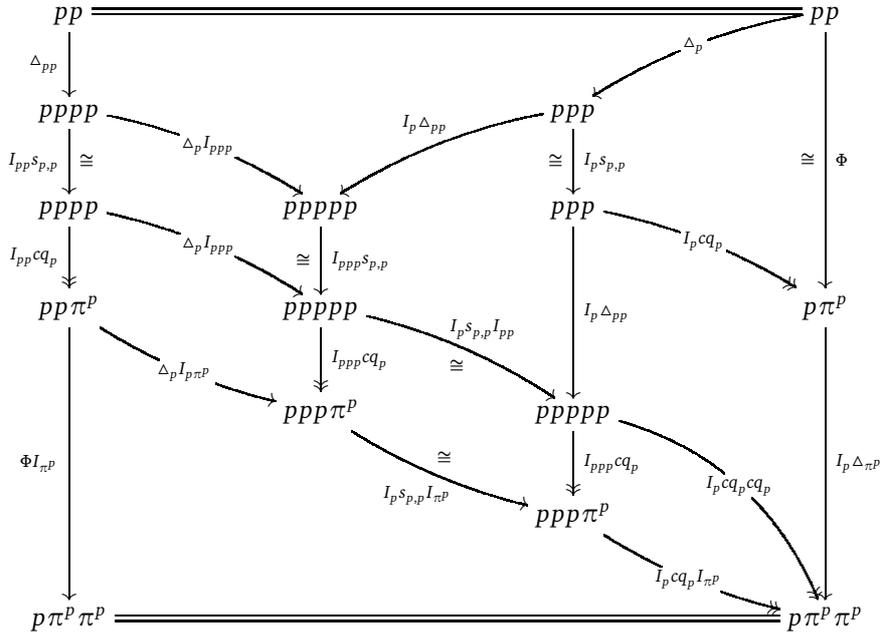
{\hbox{
      \xymatrix@C=60pt{
        pp
        \ar[d]_-{\cp_{pp}}
        \ar@{=}[rrr]
        &\text{ }
        &\text{ }
        &pp
        \ar@/_0.5pc/[dl]|-{\cp_p}
        \ar[ddd]^-{\Phi}_-{\cong}
        \\
        pppp
        \ar[d]_-{I_{pp}s_{p,p}}^-{\cong}
        \ar@/^0.5pc/[dr]|-{\cp_pI_{ppp}}
        &\text{ }
        &ppp
        \ar[d]^-{I_ps_{p,p}}_-{\cong}
        \ar@/_0.5pc/[dl]_-{I_p\cp_{pp}}
        &\text{ }
        \\
        pppp
        \ar@{->>}[d]_-{I_{pp}\textit{cq}_p}
        \ar@/^0.5pc/[dr]|-{\cp_pI_{ppp}}
        &ppppp
        \ar[d]^-{I_{ppp}s_{p,p}}_-{\cong}
        &ppp
        \ar[dd]^-{I_p\cp_{pp}}
        \ar@{->>}@/^0.5pc/[dr]|-{I_p\textit{cq}_p}
        &\text{ }
        \\
        pp\pi^p
        \ar[ddd]_-{\Phi I_{\pi^p}}
        \ar@/_0.5pc/[dr]|-{\cp_pI_{p\pi^p}}
        &ppppp
        \ar@{->>}[d]^-{I_{ppp}\textit{cq}_p}
        \ar@/^0.5pc/[dr]^-{I_ps_{p,p}I_{pp}}_-{\cong}
        &\text{ }
        &p\pi^p
        \ar[ddd]^-{I_p\cp_{\pi^p}}
        \\
        \text{ }
        &ppp\pi^p
        \ar@/_0.5pc/[dr]_-{I_ps_{p,p}I_{\pi^p}}^-{\cong}
        &ppppp
        \ar@{->>}[d]^-{I_{ppp}\textit{cq}_p}
        \ar@{->>}@/^1.5pc/[ddr]|-{I_p\textit{cq}_p\textit{cq}_p}
        &\text{ }
        \\
        \text{ }
        &\text{ }
        &ppp\pi^p
        \ar@/_0.5pc/@{->>}[dr]|-{I_p\textit{cq}_pI_{\pi^p}}
        &\text{ }
        \\
        p\pi^p\pi^p
        \ar@{=}[rrr]
        &\text{ }
        &\text{ }
        &p\pi^p\pi^p
      }
    }}
  \end{equation*}
  This completes the proof of Lemma~\ref{lem Torsors pip definition}.
\qed\end{proof}

\begin{proposition} \label{prop Torsors righttors become bitors}
  Let $\pi$ be a pre-Galois object in $\kos{K}$
  and let $(p,\lambda_p)$ be a right $\pi$-torsor over $\kappa$.
  Recall the object $\pi^p$ in $\cat{Ens}(\kos{K})$
  and the morphism $\gamma^p_p:\pi^p\otimes p\to p$
  in $\cat{Ens}(\kos{K})$
  introduced in Lemma~\ref{lem Torsors pip definition}.
  \begin{enumerate}
    \item 
    We have a pre-Galois object $\pi^p$ in $\kos{K}$,
    whose product, unit morphisms
    $\pc_{\pi^p}$, $u_{\pi^p}$
    are the unique morphisms satisfying the following relations.
    \begin{equation*}
      \vcenter{\hbox{
        \xymatrix@C=40pt{
          \pi^p\otimes p\otimes p
          \ar[d]_-{\gamma^p_p\otimes I_p}
          \ar@{->>}[r]^-{I_{\pi^p}\otimes \textit{cq}_p}
          &\pi^p\otimes \pi^p
          \ar@{.>}[d]^-{\pc_{\pi^p}}_-{\exists!}
          &p
          \ar[r]^-{e_p}
          \ar[d]_-{\cp_p}
          &\kappa
          \ar@{.>}[d]^-{u_{\pi^p}}_-{\exists!}
          \\
          p\otimes p
          \ar@{->>}[r]^-{\textit{cq}_p}
          &\pi^p
          &p\otimes p
          \ar@{->>}[r]^-{\textit{cq}_p}
          &\pi^p
        }
      }}
    \end{equation*}

    \item
    The given right $\pi$-torsor
    $(p,\lambda_p)$ over $\kappa$
    becomes a $(\pi^p,\pi)$-bitorsor
    $(p,\gamma^p_p,\lambda_p)$ over $\kappa$.
  \end{enumerate}
\end{proposition}
\begin{proof}
  We claim that $\pi^p$ is a pre-Galois object in $\kos{K}$
  as described in statement 1.
  We have a natural isomorphism
  $\Phi\otimes\slot:
  \xymatrix@C=15pt{
    p\otimes p\otimes \slot
    \ar@2{->}[r]^-{\cong}
    &p\otimes \pi^p\otimes \slot
    :\CK\to \CK
  }$
  where $\Phi:p\otimes p\xrightarrow{\cong}p\otimes \pi^p$
  is the isomorphism described in Lemma~\ref{lem Torsors pip definition}.
  As the functor $p\otimes \slot:\CK\to \CK$
  is conservative and preserves reflexive coequalizers,
  we obtain that the functor
  $\pi^p\otimes \slot:\CK\to \CK$
  is conservative and preserves reflexive coequalizers.

  We then show that the morphism
  $\pc_{\pi^p}$ is well-defined.
  As the functor $\pi^p\otimes\slot:\CK\to \CK$
  preserves the reflexive coequalizer (\ref{eq2 Torsors pip definition}),
  we need to show that the composition
  $\textit{cq}_p\circ (\gamma^p_p\otimes I_p):\pi^p\otimes p\otimes p\to \pi^p$
  coequalizes the parallel pair of morphisms
  $I_{\pi^p\otimes p\otimes p}\otimes e_p$,
  $(I_{\pi^p}\otimes \lambda_p\otimes I_p)\circ (I_{\pi^p\otimes p}\otimes \mon{d}_p\otimes I_p)\circ (I_{\pi^p\otimes p\otimes p}\otimes \cp_p)
  :\pi^p\otimes p\otimes p\otimes p\to \pi^p\otimes p\otimes p$.
  We can check this as follows.
  \begin{equation*}
    \vcenter{\hbox{
      \xymatrix@C=40pt{
        \pi^p ppp
        \ar[d]_-{I_{\pi^p pp}\cp_p}
        \ar@/^0.5pc/[dr]|-{\gamma^p_pI_{pp}}
        \ar@{=}[rr]
        &\text{ }
        &\pi^p ppp
        \ar[d]^-{I_{\pi^p pp}e_p}
        \\
        \pi^p pppp
        \ar[d]_-{I_{\pi^pp}\mon{d}_p I_p}
        \ar@/^0.5pc/[dr]|-{\gamma^p_pI_{ppp}}
        &ppp
        \ar[d]^-{I_{pp}\cp_p}
        \ar@/^1pc/[ddr]|-{I_{pp}e_p}
        &\pi^p pp
        \ar[dd]^-{\gamma^p_pI_p}
        \\
        \pi^p p\pi p
        \ar[d]_-{I_{\pi^p}\lambda_pI_p}
        \ar@/^0.5pc/[dr]|-{\gamma^p_pI_{\pi p}}
        &pppp
        \ar[d]^-{I_p\mon{d}_pI_p}
        &\text{ }
        \\
        \pi^p pp
        \ar[d]_-{\gamma^p_pI_p}
        &p\pi p
        \ar@/^0.5pc/[dl]|-{I_p\lambda_p}
        &pp
        \ar@{->>}[dd]^-{\textit{cq}_p}
        \\
        pp
        \ar@{->>}[d]_-{\textit{cq}_p}
        &\text{ }
        &\text{ }
        \\
        \pi^p
        \ar@{=}[rr]
        &\text{ }
        &\pi^p
      }
    }}
  \end{equation*}
  This shows that the morphism $\pc_{\pi^p}$ is well-defined.
  To conclude that $u_{\pi^p}$ is well-defined,
  it suffices to check that the morphism
  $\textit{cq}_p\circ \cp_p:p\to \pi^p$
  coequalizes the parallel pair of morphisms
  $I_p\otimes e_p$, $e_p\otimes I_p:p\otimes p\to p$.
  This is because we have the reflexive coequalizer diagram
  \begin{equation*}
    \vcenter{\hbox{
      \xymatrix@C=40pt{
        p\otimes p
        \ar@<0.5ex>[r]^-{I_p\otimes e_p}
        \ar@<-0.5ex>[r]_-{e_p\otimes I_p}
        &p
        \ar[r]^-{e_p}
        \ar@/_1pc/@<-2ex>[l]|-{\cp_p}
        &\kappa
        .
      }
    }}
  \end{equation*}
  We have
  \begin{equation*}
    \vcenter{\hbox{
      \xymatrix@C=40pt{
        pp
        \ar[ddd]_-{I_pe_p}
        \ar@{=}[r]
        &pp
        \ar[d]^-{\cp_pI_p}
        \ar@{=}[rr]
        &\text{ }
        &pp
        \ar@/^0.5pc/[dl]|-{I_p\cp_p}
        \ar[dd]^-{e_pI_p}
        \\
        \text{ }
        &ppp
        \ar@/_0.5pc/[dddl]|-{I_{pp}e_p}
        \ar[d]^-{I_{pp}\cp_p}
        &ppp
        \ar@/^0.5pc/[dl]|-{\cp_pI_{pp}}
        \ar[dd]^-{\tau_p^{-1}I_p}_-{\cong}
        \ar@/^0.5pc/[dddr]|-{e_pI_{pp}}
        &\text{ }
        \\
        \text{ }
        &pppp
        \ar@/_0.5pc/[dr]|-{I_p\mon{d}_pI_p}
        &\text{ }
        &p
        \ar[dd]^-{\cp_p}
        \\
        p
        \ar[d]_-{\cp_p}
        &\text{ }
        &p\pi p
        \ar@/_0.5pc/[dr]|-{\lambda_pI_p}
        &\text{ }
        \\
        pp
        \ar@{->>}[d]_-{\textit{cq}_p}
        &\text{ }
        &\text{ }
        &pp
        \ar@{->>}[d]^-{\textit{cq}_p}
        \\
        \pi^p
        \ar@{=}[rrr]
        &\text{ }
        &\text{ }
        &\pi^p
      }
    }}
  \end{equation*}
  which shows that $u_{\pi^p}$ is well-defined.
  We obtain the associativity relation of $\pc_{\pi^p}$
  by right-cancelling the epimorphisms
  $I_{\pi^p}\otimes \textit{cq}_p\otimes I_{p\otimes p}$
  and
  $I_{\pi^p\otimes \pi^p}\otimes \textit{cq}_p$
  in the diagram below.
  \begin{equation*}
    \vcenter{\hbox{
      \xymatrix@C=30pt{
        \pi^p pppp
        \ar@{->>}[d]_-{I_{\pi^p}\textit{cq}_pI_{pp}}
        \ar@{=}[rr]
        &\text{ }
        &\pi^p pppp
        \ar[d]^-{\gamma^p_pI_{ppp}}
        \ar@{=}[rr]
        \ar@/^0.5pc/[dr]|-{I_{\pi^pp}\mon{d}_pI_p}
        &\text{ }
        &\pi^p pppp
        \ar@{->>}[d]^-{I_{\pi^p}\textit{cq}_pI_{pp}}
        \\
        \pi^p \pi^p pp
        \ar@{->>}[d]_-{I_{\pi^p\pi^p}\textit{cq}_p}
        \ar@/^0.5pc/[dr]|-{\pc_{\pi^p}I_{pp}}
        &\text{ }
        &pppp
        \ar@{->>}@/_0.5pc/[dl]|-{\textit{cq}_pI_{pp}}
        \ar[d]^-{I_p\mon{d}_pI_p}
        &\pi^pp\pi p
        \ar[d]^-{I_{\pi^p}\lambda_pI_p}
        \ar@/^0.5pc/[dl]|-{\gamma^p_pI_{\pi p}}
        &\pi^p \pi^p pp
        \ar@{->>}[d]^-{I_{\pi^p\pi^p}\textit{cq}_p}
        \ar@/^0.5pc/[dl]|-{I_{\pi^p}\gamma^p_pI_p}
        \\
        \pi^p \pi^p \pi^p
        \ar[d]_-{\pc_{\pi^p}I_{\pi^p}}
        &\pi^p pp 
        \ar@{->>}@/^0.5pc/[dl]|-{I_{\pi^p}\textit{cq}_p}
        \ar@/_0.5pc/[dr]|-{\gamma^p_pI_p}
        &p\pi p
        \ar[d]^-{\lambda_pI_p}
        &\pi^p pp
        \ar@/^0.5pc/[dl]|-{\gamma^p_pI_p}
        \ar@{->>}@/_0.5pc/[dr]|-{I_{\pi^p}\textit{cq}_p}
        &\pi^p \pi^p \pi^p
        \ar[d]^-{I_{\pi^p}\pc_{\pi^p}}
        \\
        \pi^p \pi^p
        \ar[d]_-{\pc_{\pi^p}}
        &\text{ }
        &pp
        \ar@{->>}[d]^-{\textit{cq}_p}
        &\text{ }
        &\pi^p \pi^p
        \ar[d]^-{\pc_{\pi^p}}
        \\
        \pi^p
        \ar@{=}[rr]
        &\text{ }
        &\pi^p
        \ar@{=}[rr]
        &\text{ }
        &\pi^p
      }
    }}
  \end{equation*}
  In the above diagram, we used the fact that $\gamma^p_p$ and $\lambda_p$ are compatible
  which we showed in Lemma~\ref{lem Torsors pip definition}.
  We obtain the left unital relation of $\pc_{\pi^p}$, $u_{\pi^p}$
  by right-cancelling the epimorphisms
  $e_p\otimes I_{p\otimes p}$
  and
  $\textit{cq}_p$
  in the diagram below.
  \begin{equation*}
    \vcenter{\hbox{
      \xymatrix@C=40pt{
        ppp
        \ar[d]_-{e_pI_{pp}}
        \ar@{=}[r]
        &ppp
        \ar[d]^-{\cp_pI_{pp}}
        \ar@{=}[rr]
        &\text{ }
        &ppp
        \ar@/_1pc/[ddl]_-{\tau_p^{-1}I_p}^-{\cong}
        \ar[ddd]^-{e_pI_{pp}}
        \\
        pp
        \ar@{->>}[d]_-{\textit{cq}_p}
        \ar@/^0.5pc/[dr]|-{u_{\pi^p}I_{pp}}
        &pppp
        \ar@{->>}[d]^-{\textit{cq}_pI_{pp}}
        \ar@/^0.5pc/[dr]|-{I_p\mon{d}_pI_p}
        &\text{ }
        &\text{ }
        \\
        \pi^p
        \ar[d]_-{u_{\pi^p}I_{\pi^p}}
        &\pi^p pp
        \ar@/^0.5pc/@{->>}[dl]|-{I_{\pi^p}\textit{cq}_p}
        \ar@/_0.5pc/[dr]_-{\gamma^p_pI_p}
        &p\pi p
        \ar[d]^-{\lambda_pI_p}
        &\text{ }
        \\
        \pi^p\pi^p
        \ar[d]_-{\pc_{\pi^p}}
        &\text{ }
        &pp
        \ar@{=}[r]
        &pp
        \ar@{->>}[d]^-{\textit{cq}_p}
        \\
        \pi^p
        \ar@{=}[rrr]
        &\text{ }
        &\text{ }
        &\pi^p
      }
    }}
  \end{equation*}
  We obtain the right unital relation of $\pc_{\pi^p}$, $u_{\pi^p}$
  by right-cancelling the epimorphisms
  $I_{p\otimes p}\otimes e_p$ and $\textit{cq}_p$ 
  in the diagram below.
  \begin{equation*}
    \vcenter{\hbox{
      \xymatrix@C=40pt{
        ppp
        \ar[d]_-{I_{pp}e_p}
        \ar@{=}[rr]
        &\text{ }
        &ppp
        \ar[d]^-{I_{pp}\cp_p}
        \ar@{->>}@/_0.5pc/[dl]|-{\textit{cq}_pI_p}
        \ar@{=}[r]
        &ppp
        \ar[ddd]^-{I_{pp}e_p}
        \\
        pp
        \ar@{->>}[d]_-{\textit{cq}_p}
        &\pi^p p
        \ar@/_0.5pc/[dl]|-{I_{\pi^p}e_p}
        \ar[d]^-{I_{\pi^p}\cp_p}
        &pppp
        \ar[d]^-{I_p\mon{d}_pI_p}
        \ar@{->>}@/^0.5pc/[dl]|-{\textit{cq}_pI_{pp}}
        &\text{ }
        \\
        \pi^p
        \ar[d]_-{I_{\pi^p}u_{\pi^p}}
        &\pi^p pp
        \ar@/^0.5pc/@{->>}[dl]|-{I_{\pi^p}\textit{cq}_p}
        \ar@/_0.5pc/[dr]_-{\gamma^p_pI_p}
        &p \pi p
        \ar[d]^-{\lambda_pI_p}
        &\text{ }
        \\
        \pi^p\pi^p
        \ar[d]_-{\pc_{\pi^p}}
        &\text{ }
        &pp
        \ar@{->>}[d]^-{\textit{cq}_p}
        &pp
        \ar@{->>}[d]^-{\textit{cq}_p}
        \\
        \pi^p
        \ar@{=}[rr]
        &\text{ }
        &\pi^p
        \ar@{=}[r]
        &\pi^p
      }
    }}
  \end{equation*}
  To conclude that $\pi^p$ is a pre-Galois object in $\kos{K}$,
  the only thing we are left to show is that the composition
  $\varphi_{\pi^p}:
  \pi^p\otimes \pi^p
  \xrightarrow{\cp_{\pi^p}\otimes I_{\pi^p}}
  \pi^p\otimes \pi^p\otimes \pi^p
  \xrightarrow{I_{\pi^p}\otimes \pc_{\pi^p}}
  \pi^p\otimes \pi^p$
  is an isomorphism.
  As the functor $p\otimes\slot:\CK\to \CK$ is conservative,
  it suffices to show that
  \begin{equation*}
    I_p\otimes \varphi_{\pi^p}:
    \xymatrix@C=35pt{
      p\otimes \pi^p\otimes \pi^p
      \ar[r]^-{I_p\otimes \cp_{\pi^p}\otimes I_{\pi^p}}
      &p\otimes \pi^p\otimes \pi^p\otimes \pi^p
      \ar[r]^-{I_{p\otimes \pi^p}\otimes \pc_{\pi^p}}
      &p\otimes \pi^p\otimes \pi^p
    }
  \end{equation*}
  is an isomorphism.
  This is also equivalent to saying that the composition
  \begin{equation*}
    \Psi:
    \xymatrix@C=35pt{
      pp\pi^p
      \ar[r]^-{\cp_{pp}I_{\pi^p}}
      &pppp\pi^p
      \ar[r]^-{I_{pp}s_{p,p}I_{\pi^p}}_-{\cong}
      &pppp\pi^p
      \ar@{->>}[r]^-{I_{pp}\textit{cq}_pI_{\pi^p}}
      &pp\pi^p\pi^p
      \ar[r]^-{I_{pp}\pc_{\pi^p}}
      &pp\pi^p
    }
  \end{equation*}
  is an isomorphism,
  as we can see from the diagram below.
  \begin{equation*}
    \vcenter{\hbox{
      \xymatrix@C=60pt{
        pp\pi^p
        \ar[d]^-{\cp_{pp}I_{\pi^p}}
        \ar@/_3pc/@<-3ex>[dddd]_-{\Psi}
        \ar[rr]^-{\Phi I_{\pi^p}}_-{\cong}
        &\text{ }
        &p\pi^p\pi^p
        \ar[dddd]^-{I_p\varphi_{\pi^p}}
        \ar@/_1pc/[dddl]|-{I_p\cp_{\pi^p}I_{\pi^p}}
        \\
        pppp\pi^p
        \ar[d]^-{I_{pp}s_{p,p}I_{\pi^p}}_-{\cong}
        &\text{ }
        &\text{ }
        \\
        pppp\pi^p
        \ar@{->>}[d]^-{I_{pp}\textit{cq}_pI_{\pi^p}}
        &\text{ }
        &\text{ }
        \\
        pp\pi^p\pi^p
        \ar[d]^-{I_{pp}\pc_{\pi^p}}
        \ar[r]^-{\Phi I_{\pi^p\pi^p}}_-{\cong}
        &p\pi^p\pi^p\pi^p
        \ar@/^0.5pc/[dr]|-{I_{p\pi^p}\pc_{\pi^p}}
        &\text{ }
        \\
        pp\pi^p
        \ar[rr]^-{\Phi I_{\pi^p}}_-{\cong}
        &\text{ }
        &p\pi^p\pi^p
      }
    }}
  \end{equation*}
  In the above diagram, we used the relation (\ref{eq2 Torsors pip definition}).
  We claim that $\Psi$ is an isomorphism whose inverse is
  \begin{equation*}
    \Psi^{-1}:
    \xymatrix@C=35pt{
      pp\pi^p
      \ar[r]^-{\cp_{pp}I_{\pi^p}}
      &pppp\pi^p
      \ar[r]^-{I_{pp}\textit{cq}_pI_{\pi^p}}
      &pp\pi^p\pi^p
      \ar[r]^-{I_{pp}\pc_{\pi^p}}
      &pp\pi^p
      .
    }
  \end{equation*}
  If we denote
  \begin{equation*}
    \begin{aligned}
      \Upsilon
      &:
      \xymatrix@C=30pt{
        pp\pi^p
        \ar[r]^-{s_{p,p}I_{\pi^p}}_-{\cong}
        &pp\pi^p
        \ar@{->>}[r]^-{\textit{cq}_pI_{\pi^p}}
        &\pi^p\pi^p
        \ar[r]^-{\pc_{\pi^p}}
        &\pi^p
      }
      \\
      \widetilde{\Upsilon}
      &:
      \xymatrix@C=30pt{
        pp\pi^p
        \ar@{->>}[r]^-{\textit{cq}_pI_{\pi^p}}
        &\pi^p\pi^p
        \ar[r]^-{\pc_{\pi^p}}
        &\pi^p
      }
    \end{aligned}
  \end{equation*}
  then $\Psi$, $\Psi^{-1}$ are inverse to each other
  if and only if $\Upsilon$, $\widetilde{\Upsilon}$
  satisfies the relations below.
  \begin{equation} \label{eq Torsors righttors become bitors}
    \vcenter{\hbox{
      \xymatrix{
        pp\pi^p
        \ar[r]^-{\cp_{pp}I_{\pi^p}}
        \ar@/_1pc/[drr]_-{e_{pp}I_{\pi^p}}
        &pppp\pi^p
        \ar[r]^-{I_{pp}\Upsilon}
        &pp\pi^p
        \ar[d]^-{\widetilde{\Upsilon}}
        \\
        \text{ }
        &\text{ }
        &\pi^p
      }
    }}
    \qquad
    \vcenter{\hbox{
      \xymatrix{
        pp\pi^p
        \ar[r]^-{\cp_{pp}I_{\pi^p}}
        \ar@/_1pc/[drr]_-{e_{pp}I_{\pi^p}}
        &pppp\pi^p
        \ar[r]^-{I_{pp}\widetilde{\Upsilon}}
        &pp\pi^p
        \ar[d]^-{\Upsilon}
        \\
        \text{ }
        &\text{ }
        &\pi^p
      }
    }}
  \end{equation}
  When verifying the relations (\ref{eq Torsors righttors become bitors}),
  we are going to use the following relation.
  \begin{equation*}
    (\dagger):
    \vcenter{\hbox{
      \xymatrix@C=40pt{
        ppp
        \ar[d]_-{I_p\cp_pI_p}
        \ar@{=}[rrr]
        &\text{ }
        &\text{ }
        &ppp
        \ar[d]^-{I_pe_pI_p}
        \\
        pppp
        \ar@{->>}[dd]_-{\textit{cq}_p\textit{cq}_p}
        \ar@{=}[r]
        &pppp
        \ar@{->>}[d]^-{\textit{cq}_pI_{pp}}
        \ar@/^0.5pc/[dr]|-{I_p\mon{d}_pI_p}
        &\text{ }
        &pp
        \ar@/_0.5pc/[dl]|-{I_pu_{\pi}I_p}
        \ar@{=}[dd]
        \\
        \text{ }
        &\pi^p pp
        \ar@{->>}@/^0.5pc/[dl]|-{I_{\pi^p}\textit{cq}_p}
        \ar@/_0.5pc/[dr]|-{\gamma^p_pI_p}
        &p\pi p
        \ar[d]^-{\lambda_pI_p}
        &\text{ }
        \\
        \pi^p\pi^p
        \ar[d]_-{\pc_{\pi^p}}
        &\text{ }
        &pp
        \ar@{=}[r]
        &pp
        \ar@{->>}[d]^-{\textit{cq}_p}
        \\
        \pi^p
        \ar@{=}[rrr]
        &\text{ }
        &\text{ }
        &\pi^p
      }
    }}
  \end{equation*}
  We show one of the relations in (\ref{eq Torsors righttors become bitors})
  as follows.
  \begin{equation*}
    \vcenter{\hbox{
      \xymatrix@C=50pt{
        pp\pi^p
        \ar[dd]_-{\cp_{pp}I_{\pi^p}}
        \ar@{=}[r]
        &pp\pi^p
        \ar[d]^-{\cp_pI_{p\pi^p}}
        \ar@{=}[rr]
        &\text{ }
        &pp\pi^p
        \ar[dd]^-{I_pe_pI_{\pi^p}}
        \\
        \text{ }
        &ppp\pi^p
        \ar[d]^-{I_ps_{p,p}I_{\pi^p}}_-{\cong}
        \ar@/^1.5pc/[ddr]^-{I_{pp}e_pI_{\pi^p}}
        &\text{ }
        &\text{ }
        \\
        pppp\pi^p
        \ar[ddd]_-{I_{pp}\Upsilon}
        \ar@/^0.5pc/[dr]^(0.45){I_{pp}s_{p,p}I_{\pi^p}}_-{\cong}
        &ppp\pi^p
        \ar[d]^-{I_p\cp_p I_{p\pi^p}}
        \ar@/^0.5pc/[dr]|-{I_pe_pI_{p\pi^p}}
        &\text{ }
        &p\pi^p
        \ar@/_0.5pc/[dl]|-{\cp_pI_{\pi^p}}
        \ar[dd]^-{e_pI_{\pi^p}}
        \\
        \text{ }
        &pppp\pi^p
        \ar@{->>}[d]^-{I_{pp}\textit{cq}_pI_{\pi^p}}
        &pp\pi^p
        \ar@{->>}[ddd]^-{\textit{cq}_p I_{\pi^p}}
        &\text{ }
        \\
        \text{ }
        &pp\pi^p\pi^p
        \ar@/^0.5pc/[dl]|-{I_{pp}\pc_{\pi^p}}
        \ar@{->>}[d]^-{\textit{cq}_pI_{\pi^p\pi^p}}
        \ar@{}[r]|-{(\dagger)}
        &\text{ }
        &\pi^p
        \ar@/^1pc/[ddl]|-{u_{\pi^p}I_{\pi^p}}
        \ar@{=}[ddd]
        \\
        pp\pi^p
        \ar[dd]_-{\widetilde{\Upsilon}}
        \ar@{->>}@/_0.5pc/[dr]|-{\textit{cq}_pI_{\pi^p}}
        &\pi^p\pi^p\pi^p
        \ar[d]^-{I_{\pi^p}\pc_{\pi^p}}
        \ar@/^0.5pc/[dr]|-{\pc_{\pi^p}I_{\pi^p}}
        &\text{ }
        &\text{ }
        \\
        \text{ }
        &\pi^p\pi^p
        \ar@/_0.5pc/[dr]|-{\pc_{\pi^p}}
        &\pi^p\pi^p
        \ar[d]^-{\pc_{\pi^p}}
        &\text{ }
        \\
        \pi^p
        \ar@{=}[rr]
        &\text{ }
        &\pi^p
        \ar@{=}[r]
        &\pi^p
      }
    }}
  \end{equation*}
  We show the other relation in (\ref{eq Torsors righttors become bitors})
  as follows.
  \begin{equation*}
    \vcenter{\hbox{
      \xymatrix@C=20pt{
        pp\pi^p
        \ar[dd]_-{\cp_{pp}I_{\pi^p}}
        \ar@{=}[r]
        &pp\pi^p
        \ar[d]^-{I_p\cp_pI_{\pi^p}}
        \ar@{=}[rrrr]
        &\text{ }
        &\text{ }
        &\text{ }
        &pp\pi^p
        \ar[dd]^-{e_pI_{p\pi^p}}
        \\
        \text{ }
        &ppp\pi^p
        \ar[d]^-{s_{p,p}I_{p\pi^p}}_-{\cong}
        \ar@/^2pc/[dddrrr]^-{e_pI_{pp\pi^p}}
        &\text{ }
        &\text{ }
        &\text{ }
        &\text{ }
        \\
        pppp\pi^p
        \ar[dd]_-{I_{pp}\widetilde{\Upsilon}}
        \ar@/^0.5pc/[dr]_-{s_{p,p}I_{pp\pi^p}}^(0.4){\cong}
        &ppp\pi^p
        \ar[d]^-{I_p\cp_p I_{p\pi^p}}
        \ar@/^1pc/[ddrrr]|-{I_pe_pI_{p\pi^p}}
        &\text{ }
        &\text{ }
        &\text{ }
        &p\pi^p
        \ar@/_1pc/[ddl]|-{\cp_pI_{\pi^p}}
        \ar[ddd]^-{e_pI_{\pi^p}}
        \\
        \text{ }
        &pppp\pi^p
        \ar[dd]_-{I_{pp}\widetilde{\Upsilon}}
        \ar@{->>}@/^0.5pc/[dr]|-{I_{pp}\textit{cq}_pI_{\pi^p}}
        &\text{ }
        &\text{ }
        &\text{ }
        &\text{ }
        \\
        pp\pi^p
        \ar[ddd]_-{\Upsilon}
        \ar@/_0.5pc/[dr]_-{s_{p,p}I_{\pi^p}}^-{\cong}
        &\text{ }
        &pp\pi^p\pi^p
        \ar@/^0.5pc/[dl]|-{I_{pp}\pc_{\pi^p}}
        \ar@{->>}@/^0.5pc/[dr]|-{\textit{cq}_pI_{\pi^p\pi^p}}
        \ar@{}[rr]|-{(\dagger)}
        &\text{ }
        &pp\pi^p
        \ar@{->>}[dd]^-{\textit{cq}_pI_{\pi^p}}
        &\text{ }
        \\
        \text{ }
        &pp\pi^p
        \ar@{->>}@/_0.5pc/[dr]|-{\textit{cq}_pI_{\pi^p}}
        &\text{ }
        &\pi^p\pi^p\pi^p
        \ar@/^0.5pc/[dl]|-{I_{\pi^p}\pc_{\pi^p}}
        \ar@/^0.5pc/[dr]|-{\pc_{\pi^p}I_{\pi^p}}
        &\text{ }
        &\pi^p
        \ar@/^0.5pc/[dl]|-{u_{\pi^p}I_{\pi^p}}
        \ar@{=}[dd]
        \\
        \text{ }
        &\text{ }
        &\pi^p\pi^p
        \ar@/_0.5pc/[dr]|-{\pc_{\pi^p}}
        &\text{ }
        &\pi^p\pi^p
        \ar@/^0.5pc/[dl]|-{\pc_{\pi^p}}
        &\text{ }
        \\
        \pi^p
        \ar@{=}[rrr]
        &\text{ }
        &\text{ }
        &\pi^p
        \ar@{=}[rr]
        &\text{ }
        &\pi^p
      }
    }}
  \end{equation*}
  We conclude that $\pi^p$ is a pre-Galois object in $\kos{K}$
  as we claimed.
  We are left to show that $(p,\gamma^p_p,\lambda_p)$
  is a $(\pi^p,\pi)$-bitorsor over $\kappa$.
  We already showed in Lemma~\ref{lem Torsors pip definition}
  that $\gamma^p_p$
  and $\lambda_p$ are compatible.
  Thus we only need to show that $(p,\gamma^p_p)$
  is a left $\pi^p$-torsor over $\kappa$.
  We can check that
  $\gamma^p_p$ satisfies the left $\pi^p$-action relations
  by right-cancelling the epimorphism
  $I_{\pi^p}\otimes \textit{cq}_p\otimes I_p$
  in the following diagram
  \begin{equation*}
    \vcenter
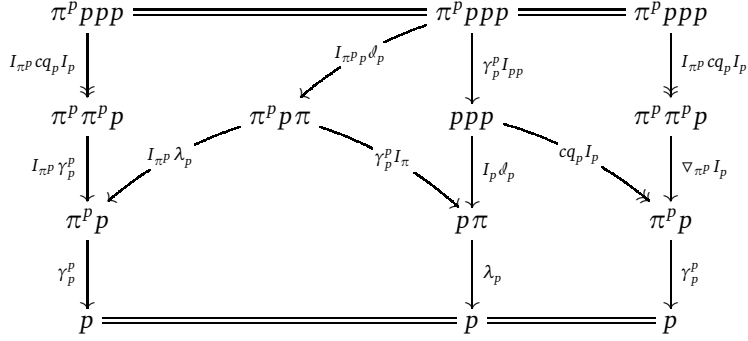
{\hbox{
      \xymatrix@C=40pt{
        \pi^p ppp
        \ar@{->>}[d]_-{I_{\pi^p}\textit{cq}_pI_p}
        \ar@{=}[rr]
        &\text{ }
        &\pi^p ppp
        \ar@/_0.5pc/[dl]|-{I_{\pi^p p}\mon{d}_p}
        \ar[d]^-{\gamma^p_p I_{pp}}
        \ar@{=}[r]
        &\pi^p ppp
        \ar@{->>}[d]^-{I_{\pi^p}\textit{cq}_pI_p}
        \\
        \pi^p\pi^p p
        \ar[d]_-{I_{\pi^p}\gamma^p_p}
        &\pi^p p\pi
        \ar@/_0.5pc/[dl]|-{I_{\pi^p}\lambda_p}
        \ar@/^0.5pc/[dr]|-{\gamma^p_p I_{\pi}}
        &ppp
        \ar@{->>}@/^0.5pc/[dr]|-{\textit{cq}_pI_p}
        \ar[d]^-{I_p\mon{d}_p}
        &\pi^p\pi^p p
        \ar[d]^-{\pc_{\pi^p}I_p}
        \\
        \pi^p p
        \ar[d]_-{\gamma^p_p}
        &\text{ }
        &p \pi
        \ar[d]^-{\lambda_p}
        &\pi^p p
        \ar[d]^-{\gamma^p_p}
        \\
        p
        \ar@{=}[rr]
        &\text{ }
        &p
        \ar@{=}[r]
        &p
      }
    }}
  \end{equation*}
  and by right-cancelling the epimorphism
  $e_p\otimes I_p$
  in the diagram below.
  \begin{equation*}
    \vcenter
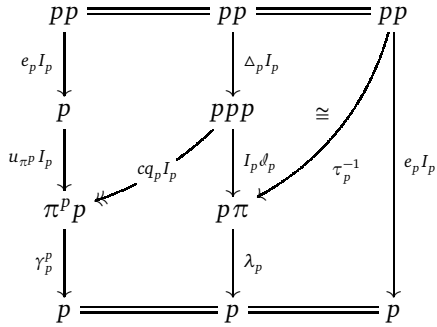
{\hbox{
      \xymatrix@C=40pt{
        pp
        \ar[d]_-{e_pI_p}
        \ar@{=}[r]
        &pp
        \ar[d]^-{\cp_pI_p}
        \ar@{=}[r]
        &pp
        \ar@/^1pc/[ddl]^(0.6){\tau_p^{-1}}_-{\cong}
        \ar[ddd]^-{e_pI_p}
        \\
        p
        \ar[d]_-{u_{\pi^p}I_p}
        &ppp
        \ar@{->>}@/^0.5pc/[dl]|-{\textit{cq}_pI_p}
        \ar[d]^-{I_p\mon{d}_p}
        &\text{ }
        \\
        \pi^p p
        \ar[d]_-{\gamma^p_p}
        &p\pi
        \ar[d]^-{\lambda_p}
        &\text{ }
        \\
        p
        \ar@{=}[r]
        &p
        \ar@{=}[r]
        &p
      }
    }}
  \end{equation*}
  Furthermore, we have
  \begin{equation*}
    \vcenter
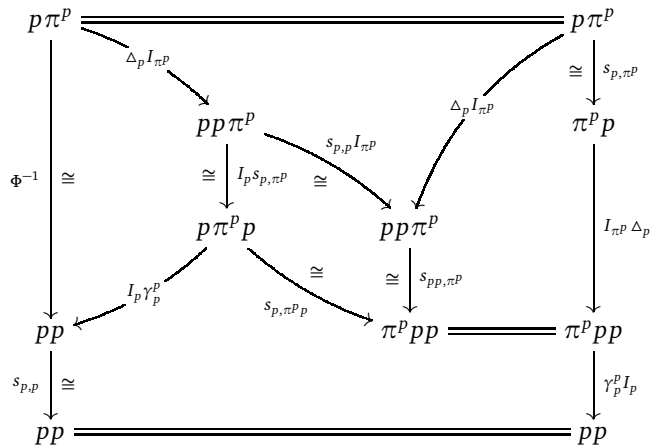
{\hbox{
      \xymatrix@C=40pt{
        p\pi^p
        \ar[ddd]_-{\Phi^{-1}}^-{\cong}
        \ar@/^0.5pc/[dr]|-{\cp_pI_{\pi^p}}
        \ar@{=}[rrr]
        &\text{ }
        &\text{ }
        &p\pi^p
        \ar[d]^-{s_{p,\pi^p}}_-{\cong}
        \ar@/_1pc/[ddl]|-{\cp_pI_{\pi^p}}
        \\
        \text{ }
        &pp\pi^p 
        \ar[d]^-{I_ps_{p,\pi^p}}_-{\cong}
        \ar@/^0.5pc/[dr]^-{s_{p,p}I_{\pi^p}}_-{\cong}
        &\text{ }
        &\pi^p p
        \ar[dd]^-{I_{\pi^p}\cp_p}
        \\
        \text{ }
        &p\pi^p p
        \ar@/^0.5pc/[dl]|-{I_p\gamma^p_p}
        \ar@/_0.5pc/[dr]_-{s_{p,\pi^pp}}^-{\cong}
        &pp\pi^p
        \ar[d]^-{s_{pp,\pi^p}}_-{\cong}
        &\text{ }
        \\
        pp
        \ar[d]_-{s_{p,p}}^-{\cong}
        &\text{ }
        &\pi^ppp
        \ar@{=}[r]
        &\pi^ppp
        \ar[d]^-{\gamma^p_pI_p}
        \\
        pp
        \ar@{=}[rrr]
        &\text{ }
        &\text{ }
        &pp
      }
    }}
  \end{equation*}
  which in particular implies that the composition
  $\pi^p\otimes p
  \xrightarrow{I_{\pi^p}\cp_p}
  \pi^p\otimes p\otimes p
  \xrightarrow{\gamma_p^p\otimes I_p}
  p\otimes p$
  is an isomorphism.
  This shows that $(p,\gamma^p_p)$ is a left $\pi^p$-torsor over $\kappa$.
  We conclude that $(p,\gamma^p_p,\lambda_p)$
  is a $(\pi^p,\pi)$-bitorsor over $\kappa$.
  This completes the proof of Proposition~\ref{prop Torsors righttors become bitors}.
\qed\end{proof}

We introduce a lemma which we will use in the next subsection.

\begin{lemma} \label{lem Torsors antipode formula}
  Let $\pi$ be a pre-Galois object in $\kos{K}$
  and let $p$
  be a right $\pi$-torsor over $\kappa$.
  \begin{enumerate}
    \item 
    The isomorphism $\tau_p:p\otimes \pi\xrightarrow{\cong}p\otimes p$
    satisfies the following relation.
    \begin{equation} \label{eq Torsors antipode formula}
      \vcenter{\hbox{
        \xymatrix@C=35pt{
          p\otimes \pi\otimes \pi
          \ar[dd]_-{\tau_p\otimes I_{\pi}}^-{\cong}
          \ar[r]^-{I_p\otimes \varphi_{\pi}}_-{\cong}
          &p\otimes \pi\otimes \pi
          \ar[r]^-{\tau_p\otimes I_{\pi}}_-{\cong}
          &p\otimes p\otimes \pi
          \ar@/_1pc/[dr]_-{(p\otimes)_{p,\pi}}
          \ar[r]^-{\cp_p\otimes I_{p\otimes \pi}}
          &p\otimes p\otimes p\otimes \pi
          \ar[d]^-{I_p\otimes s_{p,p}\otimes I_{\pi}}_-{\cong}
          \\
          \text{ }
          &\text{ }
          &\text{ }
          &p\otimes p\otimes p\otimes \pi
          \ar[d]^-{I_{p\otimes p}\otimes \lambda_p}
          \\
          p\otimes p\otimes \pi
          \ar[rrr]^-{I_p\otimes \tau_p}_-{\cong}
          &\text{ }
          &\text{ }
          &p\otimes p\otimes p
        }
      }}
    \end{equation}
  
    \item 
    The division morphism $\mon{d}_p:p\otimes p\to \pi$
    satisfies the following relation.
    \begin{equation} \label{eq2 Torsors antipode formula}
      \vcenter{\hbox{
        \xymatrix@C=40pt{
          p\otimes p
          \ar[d]_-{s_{p,p}}^-{\cong}
          \ar[r]^-{\mon{d}_p}
          &\pi
          \ar[d]^-{\varsigma_{\pi}}_-{\cong}
          \\
          p\otimes p
          \ar[r]^-{\mon{d}_p}
          &\pi
        }
      }}
    \end{equation}

    \item 
    The induced left $\pi$-action
    $\xymatrix{
      \gamma_p:\pi\otimes p
      \ar[r]^-{\varsigma_{\pi}\otimes I_p}_-{\cong}
      &\pi\otimes p
      \ar[r]^-{s_{\pi,p}}_-{\cong}
      &p\otimes \pi
      \ar[r]^-{\lambda_p}
      &p
    }$
    satisfies the following relation.
    \begin{equation}\label{eq3 Torsors antipode formula}
      \vcenter{\hbox{
        \xymatrix@C=30pt{
          p\otimes p
          \ar[d]_-{I_p\otimes \cp_p}
          \ar[r]^-{I_p\otimes e_p}
          &p\otimes \kappa
          \ar[r]^-{\jmath_p^{-1}}_-{\cong}
          &p
          \ar[d]^-{\cp_p}
          \\
          p\otimes p\otimes p
          \ar[r]^-{\tau_p^{-1}\otimes I_p}_-{\cong}
          &p\otimes \pi\otimes p
          \ar[r]^-{I_p\otimes \gamma_p}
          &p\otimes p
        }
      }}
    \end{equation}
  \end{enumerate}
\end{lemma}
\begin{proof}
  We omit the symbol $\otimes$ and the coherence isomorphisms $a$, $\imath$, $\jmath$.
  We can check the relation (\ref{eq Torsors antipode formula})
  as in the following diagram
  \begin{equation*}
    \vcenter{\hbox{
      \xymatrix@C=65pt{
        p\pi\pi
        \ar[dd]_-{I_p\varphi_{\pi}}^-{\cong}
        \ar@{=}[rr]
        &\text{ }
        &p\pi\pi
        \ar@/_0.5pc/[dl]|-{I_p\cp_{\pi}I_{\pi}}
        \ar[d]^-{\cp_pI_{\pi\pi}}
        \ar@{=}[r]
        &p\pi\pi
        \ar[dddd]^-{\tau_pI_{\pi}}_-{\cong}
        \\
        \text{ }
        &p\pi\pi\pi
        \ar@/_0.5pc/[dl]|-{I_{p\pi}\pc_{\pi}}
        \ar[d]^-{\cp_pI_{\pi\pi\pi}}
        &p p \pi\pi
        \ar@/^0.5pc/[dl]|-{I_{pp}\cp_{\pi}I_{\pi}}
        \ar[dd]^-{I_p\cp_{p\pi}I_{\pi}}
        \ar@/^1pc/[dddr]|-{I_p\lambda_pI_{\pi}}
        &\text{ }
        \\
        p\pi\pi
        \ar[dd]_-{\tau_p I_{\pi}}^-{\cong}
        \ar@/_0.5pc/[dr]|-{\cp_pI_{\pi\pi}}
        &p p \pi\pi\pi
        \ar[d]^-{I_{pp\pi}\pc_{\pi}}
        \ar@/^0.5pc/[dr]|-{I_p(p\otimes)_{\pi,\pi\pi}}
        &\text{ }
        &\text{ }
        \\
        \text{ }
        \ar@{}[dr]|-{(\dagger)}
        &p p \pi \pi
        \ar[d]^-{I_p(p\otimes)_{\pi,\pi}}
        &p p \pi p \pi\pi
        \ar@/^0.5pc/[dl]|-{I_{pp\pi p}\pc_{\pi}}
        \ar[dd]^-{I_p\lambda_p\lambda_pI_{\pi}}
        &\text{ }
        \\
        p p \pi
        \ar[d]_-{(p\otimes)_{p,\pi}}
        &p p \pi p \pi
        \ar@/^0.5pc/[dl]|-{I_p\lambda_pI_{p\pi}}
        \ar[dd]^-{I_p\lambda_p\lambda_p}
        &\text{ }
        &p p \pi
        \ar@/^0.5pc/[dl]|-{I_P\cp_pI_{\pi}}
        \ar[dd]^-{I_p\tau_p}_-{\cong}
        \\
        p p p \pi
        \ar[d]_-{I_{pp}\lambda_p}
        &\text{ }
        &p p p \pi
        \ar[d]_-{I_{pp}\lambda_p}
        &\text{ }
        \\
        p p p
        \ar@{=}[r]
        &p p p
        \ar@{=}[r]
        &p p p
        \ar@{=}[r]
        &p p p
      }
    }}
  \end{equation*}
  where the diagram $(\dagger)$ is verified below.
  \begin{equation*}
    (\dagger):
    \vcenter{\hbox{
      \xymatrix@C=50pt{
        p\pi\pi
        \ar[dd]_-{\tau_p I_{\pi}}^-{\cong}
        \ar@{=}[r]
        &p\pi\pi
        \ar[d]^-{\cp_pI_{\pi\pi}}
        \ar@{=}[rr]
        &\text{ }
        &p\pi\pi
        \ar[d]^-{\cp_pI_{\pi\pi}}
        \\
        \text{ }
        &pp\pi\pi
        \ar@/^0.5pc/[dl]|-{I_p\lambda_pI_{\pi}}
        \ar@/^0.5pc/[dr]|-{\cp_pI_{p\pi\pi}}
        \ar[ddd]_-{(p\otimes)_{p\pi,\pi}}
        &\text{ }
        &pp\pi\pi
        \ar@/_0.5pc/[dl]|-{I_p\cp_pI_{\pi\pi}}
        \ar[ddd]^-{I_p(p\otimes)_{\pi,\pi}}
        \ar@/^1.5pc/[ddl]|-{I_p\cp_pI_{\pi\pi}}
        \\
        pp\pi
        \ar[ddd]_-{(p\otimes)_{p,\pi}}
        &\text{ }
        &ppp\pi\pi
        \ar[d]^-{I_ps_{p,p}I_{\pi\pi}}_-{\cong}
        \ar@/_2pc/@<-2ex>[dd]_-{I_ps_{p,p\pi}I_{\pi}}^-{\cong}
        &\text{ }
        \\
        \text{ }
        &\text{ }
        &ppp\pi\pi
        \ar[d]^-{I_{pp}s_{p,\pi}I_{\pi}}_-{\cong}
        &\text{ }
        \\
        \text{ }
        &pp\pi p\pi
        \ar@{=}[r]
        &pp\pi p\pi
        \ar@{=}[r]
        &pp\pi p\pi
        \ar[d]^-{I_p\lambda_pI_{\pi}}
        \\
        ppp\pi
        \ar@{=}[rrr]
        &\text{ }
        &\text{ }
        &ppp\pi
      }
    }}
  \end{equation*}
  Next we show the relation (\ref{eq2 Torsors antipode formula}).
  \begin{equation*}
    \vcenter{\hbox{
      \xymatrix@C=40pt{
        pp
        \ar[ddd]_-{\mon{d}_p}
        \ar@{=}[rr]
        &\text{ }
        &pp
        \ar[d]^-{I_{pp}u_{\pi}}
        \ar@{=}[rrr]
        &\text{ }
        &\text{ }
        &pp
        \ar[dd]^-{s_{p,p}}_-{\cong}
        \\
        \text{ }
        &\text{ }
        &pp\pi
        \ar@/_1pc/[dddl]|-{\mon{d}_pI_{\pi}}
        \ar[dd]^-{\tau_p^{-1}I_{\pi}}_-{\cong}
        \ar@{=}[r]
        &pp\pi
        \ar[d]^-{(p\otimes)_{p,\pi}}
        \ar@/^1.8pc/[ddr]^(0.45){s_{p,p}I_{\pi}}_(0.45){\cong}
        &\text{ }
        &\text{ }
        \\
        \text{ }
        &\text{ }
        &\text{ }
        &ppp\pi
        \ar[d]^-{I_{pp}\lambda_p}
        \ar@/^0.5pc/[dr]|-{e_pI_{pp\pi}}
        &\text{ }
        &pp
        \ar@/^0.5pc/[dl]|-{I_{pp}u_{\pi}}
        \ar@{=}[dd]
        \\
        \pi
        \ar[dddd]_-{\varsigma_{\pi}}^-{\cong}
        \ar@/^0.5pc/[dr]|-{I_{\pi}u_{\pi}}
        &\text{ }
        &p\pi\pi
        \ar@/_0.5pc/[dl]|-{e_pI_{\pi\pi}}
        \ar[dd]^-{I_p\varphi_{\pi}^{-1}}_-{\cong}
        \ar@{}[r]|-{(\dagger\!\dagger)}
        &ppp
        \ar[d]^-{I_p\tau_p^{-1}}_-{\cong}
        \ar@/^0.5pc/[dr]|-{e_pI_{pp}}
        &pp\pi
        \ar[d]^-{I_p\lambda_p}
        &\text{ }
        \\
        \text{ }
        &\pi\pi
        \ar[dd]^-{\varphi_{\pi}^{-1}}_-{\cong}
        &\text{ }
        &pp\pi
        \ar@/^0.5pc/[dl]^-{\tau_p^{-1}I_{\pi}}_-{\cong}
        \ar[ddd]^-{e_{pp}I_{\pi}}
        &pp
        \ar@{=}[r]
        &pp
        \ar[ddd]^-{\mon{d}_p}
        \\
        \text{ }
        &\text{ }
        &p\pi\pi
        \ar@/^0.5pc/[dl]|-{e_pI_{\pi\pi}}
        \ar[dd]^-{e_{p\pi}I_{\pi}}
        &\text{ }
        &\text{ }
        &\text{ }
        \\
        \text{ }
        &\pi\pi
        \ar[d]^-{e_{\pi}I_{\pi}}
        &\text{ }
        &\text{ }
        &\text{ }
        &\text{ }
        \\
        \pi
        \ar@{=}[r]
        &\pi
        \ar@{=}[r]
        &\pi
        \ar@{=}[r]
        &\pi
        \ar@{=}[rr]
        &\text{ }
        &\pi
      }
    }}
  \end{equation*}
  We used the relation (\ref{eq Torsors antipode formula})
  in the diagram $(\dagger\!\dagger)$.
  Finally we obtain the relation (\ref{eq3 Torsors antipode formula})
  as follows
  \begin{equation*}
    \vcenter{\hbox{
      \xymatrix@C=60pt{
        pp
        \ar[d]_-{I_p\cp_p}
        \ar@{=}[r]
        &pp
        \ar[d]^-{\cp_pI_p}
        \ar@{=}[r]
        &pp
        \ar[d]^-{I_pe_p}
        \\
        ppp
        \ar[dd]_-{\tau_p^{-1}I_p}^-{\cong}
        \ar@/^0.5pc/[dr]|-{\cp_pI_{pp}}
        &ppp
        \ar[d]^-{I_{pp}\cp_p}
        \ar@/^1.5pc/[dddr]|-{I_{pp}e_p}
        &p
        \ar[ddd]^-{\cp_p}
        \\
        \text{ }
        &pppp
        \ar@/^0.5pc/[dl]|-{I_p\mon{d}_pI_p}
        \ar@{}[dd]|-{(\dagger\!\dagger\!\dagger)}
        &\text{ }
        \\
        p\pi p
        \ar[d]_-{I_p\gamma_p}
        &\text{ }
        &\text{ }
        \\
        pp
        \ar@{=}[rr]
        &\text{ }
        &pp
      }
    }}
  \end{equation*}
  where the diagram $(\dagger\!\dagger\!\dagger)$
  is verified below.
  \begin{equation*}
    (\dagger\!\dagger\!\dagger):
    \quad
    \vcenter{\hbox{
      \xymatrix@C=60pt{
        pp
        \ar[d]_-{I_p\cp_p}
        \ar@{=}[rr]
        &\text{ }
        &pp
        \ar[d]^-{I_p\cp_p}
        \ar@{=}[r]
        &pp
        \ar[dd]^-{s_{p,p}}_-{\cong}
        \ar@/^2pc/@<3ex>[ddddd]^-{I_pe_p}
        \\
        ppp
        \ar[d]_-{\mon{d}_pI_p}
        \ar[rr]^-{I_ps_{p,p}}_-{\cong}
        \ar@/_0.5pc/[dr]_-{s_{p,p}I_p}^-{\cong}
        &\text{ }
        &ppp
        \ar[dd]^-{s_{p,pp}}_-{\cong}
        &\text{ }
        \\
        \pi p
        \ar[ddd]_-{\gamma_p}
        \ar@/_0.5pc/[dr]_-{\varsigma_{\pi}I_p}^-{\cong}
        &ppp
        \ar[d]^-{\mon{d}_pI_p}
        \ar@/_0.5pc/[dr]_-{s_{pp,p}}^-{\cong}
        &\text{ }
        &pp
        \ar@/_0.5pc/[dl]|-{\cp_pI_p}
        \ar@/^1pc/[ddl]^-{\tau_p^{-1}}_-{\cong}
        \ar[ddd]^-{e_pI_p}
        \\
        \text{ }
        &\pi p
        \ar@/_0.5pc/[dr]_-{s_{\pi,p}}^-{\cong}
        &ppp
        \ar[d]^-{I_p\mon{d}_p}
        &\text{ }
        \\
        \text{ }
        &\text{ }
        &p\pi
        \ar@/_0.5pc/[dr]|-{\lambda_p}
        &\text{ }
        \\
        p
        \ar@{=}[rrr]
        &\text{ }
        &\text{ }
        &p
      }
    }}
  \end{equation*}
  This completes the proof of Lemma~\ref{lem Torsors antipode formula}.
\qed\end{proof}


\subsubsection{Pre-fiber functor twisted by a torsor}
\begin{lemma}\label{lem twistbyatorsor}
  Let $\pi$ be a pre-Galois object in $\kos{K}$
  and let $p$ be a right $\pi$-torsor over $\kappa$.
  We can twist the forgetful
  strong $\kos{K}$-tensor functor
  $(\omega_{\pi}^*,\what{\omega}_{\pi}^*):
  (\kos{R\!e\!p}(\pi),\kos{t}^*_{\pi})
  \to
  (\kos{K},\id_{\kos{K}})$
  by $p$ and obtain the strong $\kos{K}$-tensor functor
  \begin{equation*}
    \vcenter{\hbox{
      \xymatrix@R=30pt@C=40pt{
        \text{ }
        &\kos{R\!e\!p}(\pi)
        \ar[d]^-{\omega^{p*}}
        \\
        \kos{K}
        \ar@/^1pc/[ur]^-{\kos{t}_{\pi}^*}
        \ar[r]_-{\id_{\kos{K}}}
        \xtwocell[r]{}<>{<-3>{\quad\what{\omega}^{p*}}}
        &\kos{K}
      }
    }}
    \qquad\quad
    (\omega^{p*},\what{\omega}^{p*}):
    (\kos{R\!e\!p}(\pi),\kos{t}^*_{\pi})
    \to
    (\kos{K},\id_{\kos{K}})
    .
  \end{equation*}
  \begin{itemize}
    \item 
    The underlying functor $\omega^{p*}$ sends 
    each object $X=(x,\gamma_x)$ in $\cat{Rep}(\pi)$
    to the following reflexive coequalizer in $\CK$.
    \begin{equation}\label{eq twistbyatorsor}
      \vcenter{\hbox{
        \xymatrix@C=40pt{
          p\otimes \pi\otimes x
          \ar@<0.5ex>[r]^-{I_p\otimes \gamma_x}
          \ar@<-0.5ex>[r]_-{\lambda_p\otimes I_x}
          &p\otimes x
          \ar@{->>}[r]^-{\xi^p_X}
          \ar@<-1ex>@/_2pc/[l]|-{I_p\otimes u_{\pi}\otimes I_x}
          &\omega^{p*}(X)
        }
      }}
    \end{equation}
    For every object $z$ in $\CK$,
    the functor $\slot\otimes z:\CK\to \CK$
    preserves the reflexive coequalizer (\ref{eq twistbyatorsor}).
    
    \item 
    The strong symmetric monoidal coherence isomorphisms of
    $\omega^{p*}$ are unique morphisms satisfying the relations below.
    Let $Y=(y,\gamma_y)$
    be another object in $\cat{Rep}(\pi)$.
    \begin{equation*}
      \vcenter{\hbox{
        \xymatrix@C=40pt{
          p\otimes (x\otimes y)
          \ar[r]^-{(p\otimes)_{x,y}}
          \ar@{->>}[d]_-{\xi^p_{X\tensor\!_{\pi}Y}}
          &(p\otimes x)\otimes (p\otimes y)
          \ar@{->>}[d]^-{\xi^p_X\otimes \xi^p_Y}
          \\
          \omega^{p*}(X\tensor\!_{\pi}Y)
          \ar@{.>}[r]^-{\exists!\text{ }\omega^{p*}_{X,Y}}_-{\cong}
          &\omega^{p*}(X)\otimes \omega^{p*}(Y)
        }
      }}
      \qquad\qquad
      \vcenter{\hbox{
        \xymatrix@C=40pt{
          p\otimes \kappa
          \ar@{->>}[d]_-{\xi^p_{\unit\!_{\pi}}}
          \ar[r]^-{\jmath_p^{-1}}_-{\cong}
          &p
          \ar[d]^-{e_p}
          \\
          \omega^{p*}(\unit\!_{\pi})
          \ar@{.>}[r]^-{\exists!\text{ }\omega^{p*}_{\unit\!_{\pi}}}_-{\cong}
          &\kappa
        }
      }}
    \end{equation*}
  
    \item
    The comonoidal natural isomorphism
    $\what{\omega}^{p*}:\omega^{p*}\kos{t}_{\pi}^*\cong \id_{\kos{K}}$
    is described below.
    Let $z$ be an object in $\CK$.
    \begin{equation*}
      \vcenter{\hbox{
        \xymatrix@C=40pt{
          p\otimes z
          \ar@{->>}[r]^-{\xi^p_{\kos{t}_{\pi}^*(z)}}
          \ar@/_1pc/[dr]_-{e_p\otimes I_z}
          &\omega^{p*}\kos{t}_{\pi}^*(z)
          \ar@{.>}[d]^-{\exists!\text{ }\what{\omega}^{p*}_z}_-{\cong}
          \\
          \text{ }
          &z
        }
      }}
    \end{equation*}
    If we denote
    $X\otimes z
    =z\bar{\acts}X
    =X\tensor\!_{\pi}\kos{t}_{\pi}^*(z)$
    as the action of $z$ on $X$,
    then the associated $\kos{K}$-equivariance 
    $\vecar{\omega}^{p*}$ of
    $(\omega^{p*},\what{\omega}^{p*})$
    is described below.
    \begin{equation*}
      \vcenter{\hbox{
        \xymatrix@C=40pt{
          p\otimes (x\otimes z)
          \ar@{->>}[d]_-{\xi^p_{X\otimes z}}
          \ar[r]^-{a_{p,x,z}}_-{\cong}
          &(p\otimes x)\otimes z
          \ar@{->>}[d]^-{\xi^p_X\otimes I_z}
          \\
          \omega^{p*}(X\otimes z)
          \ar@{.>}[r]^-{\exists!\text{ }\vecar{\omega}^{p*}_{z,X}}_-{\cong}
          &\omega^{p*}(X)\otimes z
        }
      }}
    \end{equation*}
  \end{itemize}
\end{lemma}
\begin{proof}
  Let $X=(x,\gamma_x)$ be a representation of $\pi$.
  The reflexive coequalizer
  $\xi^p_X$
  in (\ref{eq twistbyatorsor}) exists in $\CK$.
  This is because $I_p\otimes u_{\pi}\otimes I_x:p\otimes x\to p\otimes \pi\otimes x$
  is a common section of
  $I_p\otimes \gamma_x$,
  $\lambda_p\otimes I_x: p\otimes \pi\otimes x\to p\otimes x$.
  This shows that the underlying functor $\omega^{p*}$ is well-defined.
  If we treat $\xi^p_X$
  as a morphism $\kos{p}^*(x)\to \kos{p}^*(\omega^{p*}(X))$
  in the coKleisli category $\kos{K}_p$,
  then we obtain a morphism
  \begin{equation*}
    \kos{p}_!(\xi^p_X):
    \xymatrix@C=30pt{
      p\otimes x
      \ar[r]^-{\cp_p\otimes I_x}
      &p\otimes p\otimes x
      \ar[r]^-{I_p\otimes \xi^p_X}
      &p\otimes \omega^{*p}(X)
    }
  \end{equation*}  
  in $\CK$.
  We claim that $\kos{p}_!(\xi^p_X)$ is an isomorphism in $\CK$.
  Consider the following diagram.
  \begin{equation}\label{eq p!xi isomoprhism}
    \vcenter{\hbox{
      \xymatrix@C=60pt{
        p\pi \pi x
        \ar[d]_-{\tau_pI_{\pi x}}^-{\cong}
        \ar@<0.5ex>[r]^-{I_{p\pi}\gamma_x}
        \ar@<-0.5ex>[r]_-{I_p\pc_{\pi}I_x}
        &p\pi x
        \ar[r]^-{I_p\gamma_x}
        \ar@/_1.5pc/@<-1ex>[l]_-{I_pu_{\pi}I_{\pi x}}
        \ar[d]^-{\tau_pI_x}_-{\cong}
        &px
        \ar@{.>}[d]^-{\kos{p}_!(\xi^p_X)}_-{\exists!}
        \ar@/_1.5pc/@<-1ex>[l]_-{I_pu_{\pi}I_x}
        \\
        pp\pi x
        \ar@<0.5ex>[r]^-{I_{pp}\gamma_x}
        \ar@<-0.5ex>[r]_-{I_p\lambda_pI_x}
        &ppx
        \ar@{->>}[r]^-{I_p\xi^p_X}
        &p\omega^{*p}(X)
      }
    }}
  \end{equation}
  The top horizontal morphisms in (\ref{eq p!xi isomoprhism})
  form a split coequalizer diagram.
  The bottom horizontal morphisms in (\ref{eq p!xi isomoprhism})
  form a reflexive coequalizer diagram,
  as the functor $p\otimes\slot:\CK\to\CK$
  preserves the defining reflexive coequalizer
  (\ref{eq twistbyatorsor})
  of $\omega^{p*}(X)$.
  The left and middle vertical isomorphisms
  are compatible with horizontal parallel pairs of morphisms.
  Thus we obtain a unique vertical isomorphism
  \begin{equation*}
    \xymatrix@C=30pt{
      p\otimes x
      \ar[r]^-{I_p\otimes u_{\pi}\otimes I_x}
      \ar@/_1pc/@<-1ex>[rr]|-{\cp_p\otimes I_x}
      &p\otimes \pi\otimes x
      \ar[r]^-{\tau_p\otimes I_x}_-{\cong}
      &p\otimes p\otimes x
      \ar@{->>}[r]^-{I_p\otimes \xi^p_X}
      &p\otimes \omega^{p*}(X)
    }
  \end{equation*}
  which is equal to $\kos{p}_!(\xi^p_X)$,
  and makes the diagram (\ref{eq p!xi isomoprhism}) commutative.
  This shows that $\kos{p}_!(\xi^p_X)$ is an isomorphism in $\CK$.

  Let $z$ be an object in $\CK$.
  We claim that the functor $\slot\otimes z:\CK\to \CK$
  preserves the defining reflexive coequalizer diagram
  (\ref{eq twistbyatorsor}) of $\omega^{p*}(X)$.
  As the functor $p\otimes\slot:\CK\to \CK$
  is conservative and preserves reflexive coequalizers,
  it suffices to show that the functor
  $p\otimes (\slot\otimes z):\CK\to \CK$
  preserves the reflexive coequalizer diagram
  (\ref{eq twistbyatorsor}).
  As we have
  $a_{p,\slot,z}:
  p\otimes (\slot\otimes z)
  \cong
  (p\otimes \slot)\otimes z
  :\CK\to \CK$
  it suffices to show that the functor
  $p\otimes (\slot\otimes z):\CK\to \CK$
  preserves the reflexive coequalizer diagram
  (\ref{eq twistbyatorsor}).
  The diagram (\ref{eq p!xi isomoprhism}) in particular tells us that
  we obtain a split coequalizer diagram after applying the functor
  $p\otimes\slot:\CK\to\CK$ to the diagram (\ref{eq twistbyatorsor}).
  Hence the coequalizer diagram in (\ref{eq p!xi isomoprhism})
  is preserved under the functor $\slot\otimes z:\CK\to \CK$.
  This implies that the functor
  $(p\otimes \slot)\otimes z:\CK\to \CK$
  preserves the coequalizer diagram (\ref{eq twistbyatorsor}).
  This shows that the functor
  $\slot\otimes z:\CK\to \CK$
  preserves the reflexive coequalizer
  (\ref{eq twistbyatorsor}).
  
  We show that $\omega^{p*}$ is a strong symmetric monoidal functor
  as we claimed.
  Let $X=(x,\gamma_x)$, $Y=(y,\gamma_y)$ be representations of $\pi$.
  The following diagrams imply that
  the morphisms $\omega^{p*}_{X,Y}$, $\omega^{p*}_{\unit}$ are well-defined.
  \begin{equation*}
    \vcenter
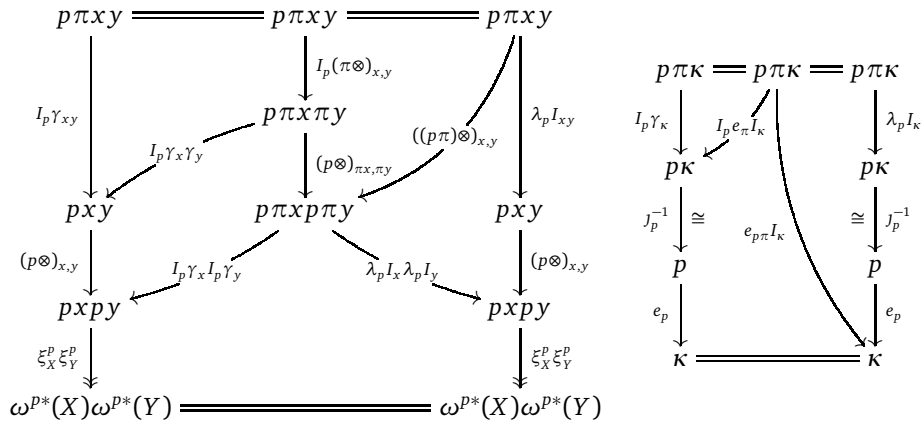
{\hbox{
      \xymatrix@C=25pt{
        p\pi xy
        \ar@{=}[r]
        \ar[dd]_-{I_p\gamma_{xy}}
        &p\pi xy
        \ar[d]^-{I_p(\pi\otimes)_{x,y}}
        \ar@{=}[r]
        &p\pi xy
        \ar@/^1.5pc/[ddl]|-{((p\pi)\otimes)_{x,y}}
        \ar[dd]^-{\lambda_pI_{xy}}
        \\
        \text{ }
        &p\pi x\pi y
        \ar@/_0.5pc/[dl]|-{I_p\gamma_x\gamma_y}
        \ar[d]^-{(p\otimes)_{\pi x,\pi y}}
        &\text{ }
        \\
        pxy
        \ar[d]_-{(p\otimes)_{x,y}}
        &p\pi xp\pi y
        \ar@/^0.5pc/[dl]|-{I_p\gamma_x I_p\gamma_y}
        \ar@/_0.5pc/[dr]|-{\lambda_pI_x\lambda_pI_y}
        &pxy
        \ar[d]^-{(p\otimes)_{x,y}}
        \\
        pxpy
        \ar@{->>}[d]_-{\xi^p_X\xi^p_Y}
        &\text{ }
        &pxpy
        \ar@{->>}[d]^-{\xi^p_X\xi^p_Y}
        \\
        \omega^{p*}(X)\omega^{p*}(Y)
        \ar@{=}[rr]
        &\text{ }
        &\omega^{p*}(X)\omega^{p*}(Y)
      }
    }}
    \text{ }
    \vcenter{\hbox{
      \xymatrix@C=12pt{
        p\pi \kappa
        \ar[d]_-{I_p\gamma_{\kappa}}
        \ar@{=}[r]
        &p\pi \kappa
        \ar@/^0.5pc/[dl]|-{I_pe_{\pi}I_{\kappa}}
        \ar@/_1pc/[dddr]_-{e_{p\pi}I_{\kappa}}
        \ar@{=}[r]
        &p\pi \kappa
        \ar[d]^-{\lambda_pI_{\kappa}}
        \\
        p\kappa
        \ar[d]_-{\jmath_p^{-1}}^-{\cong}
        &\text{ }
        &p\kappa
        \ar[d]^-{\jmath_p^{-1}}_-{\cong}
        \\
        p
        \ar[d]_-{e_p}
        &\text{ }
        &p
        \ar[d]^-{e_p}
        \\
        \kappa
        \ar@{=}[rr]
        &\text{ }
        &\kappa
      }
    }}
  \end{equation*}
  We need to show that the morphisms
  $\omega^{p*}_{X,Y}$, $\omega^{p*}_{\unit}$
  are isomorphisms.
  As the functor $p\otimes\slot:\CK\to \CK$
  is conservative,
  it suffices to show that 
  $I_p\otimes \omega^{p*}_{X,Y}$,
  $I_p\otimes \omega^{p*}_{\unit}$
  are isomorphisms.
  We obtain that $I_p\otimes \omega^{p*}_{X,Y}$
  is an isomorphism from the following relation
  \begin{equation*}
    \vcenter{\hbox{
      \xymatrix@C=45pt{
        p\otimes x\otimes y
        \ar[d]_-{I_p\otimes s_{x,y}}^-{\cong}
        \ar[rr]^-{\kos{p}_!(\xi^p_{X\tensor Y})}_-{\cong}
        &\text{ }
        &p\otimes \omega^{p*}(X\tensor Y)
        \ar[dd]^-{I_p\otimes \omega^{p*}_{X,Y}}
        \\
        p\otimes y\otimes x
        \ar[d]_-{\kos{p}_!(\xi^p_Y)\otimes I_x}^-{\cong}
        &\text{ }
        &\text{ }
        \\
        p\otimes \omega^{p*}(Y)\otimes x
        \ar[r]^-{I_p\otimes s_{\omega^{p*}(Y),x}}_-{\cong}
        &p\otimes x\otimes \omega^{p*}(Y)
        \ar[r]^-{\kos{p}_!(\xi^p_X)\otimes I_{\omega^{p*}(Y)}}_-{\cong}
        &p\otimes \omega^{p*}(X)\otimes \omega^{p*}(Y)
      }
    }}
  \end{equation*}
  which we verify below.
  \begin{equation*}
    \vcenter{\hbox{
      \xymatrix@C=25pt{
        pxy
        \ar[d]_-{I_ps_{x,y}}^-{\cong}
        \ar@{=}[rr]
        &\text{ }
        &pxy
        \ar@/^1.7pc/[dddl]|(0.5){(p\otimes)_{x,y}}
        \ar[dd]^-{\cp_pI_{xy}}
        \ar@{=}[r]
        &pxy
        \ar[ddd]^-{\kos{p}_!(\xi^p_{XY})}_-{\cong}
        \\
        pyx
        \ar[dd]_-{\kos{p}_!(\xi^p_Y) I_x}^-{\cong}
        \ar@/^0.5pc/[dr]|-{\cp_pI_{yx}}
        &\text{ }
        &\text{ }
        &\text{ }
        \\
        \text{ }
        &ppyx
        \ar@{->>}@/^0.5pc/[dl]|-{I_p\xi^p_YI_x}
        \ar[d]^-{I_ps_{py,x}}_-{\cong}
        &ppxy
        \ar[dd]^-{I_p(p\otimes)_{x,y}}
        \ar@{->>}@/^0.5pc/[dr]|-{I_p\xi^p_{XY}}
        &\text{ }
        \\
        p\omega^{p*}(Y) x
        \ar[d]_-{I_p s_{\omega^{p*}(Y),x}}^-{\cong}
        &pxpy
        \ar@{->>}@/^0.5pc/[dl]|-{I_{px}\xi^p_Y}
        \ar[dd]^-{\kos{p}_!(\xi^p_X) I_{py}}_-{\cong}
        \ar@/^0.5pc/[dr]|-{\cp_pI_{py}}
        &\text{ }
        &p\omega^{p*}(XY)
        \ar[ddd]^-{I_p\omega^{p*}_{X,Y}}
        \\
        p x \omega^{p*}(Y)
        \ar[dd]_-{\kos{p}_!(\xi^p_X) I_{\omega^{p*}(Y)}}^-{\cong}
        &\text{ }
        &ppxpy
        \ar@{->>}@/^0.5pc/[dl]|-{I_p\xi^p_XI_{py}}
        \ar@{->>}@/^1pc/[ddr]|-{I_p\xi^p_X\xi^p_Y}
        &\text{ }
        \\
        \text{ }
        &p\omega^{p*}(X)py
        \ar@{->>}[d]^-{I_{p\omega^{p*}(X)}\xi^p_Y}
        &\text{ }
        &\text{ }
        \\
        p \omega^{p*}(X)\omega^{p*}(Y)
        \ar@{=}[r]
        &p \omega^{p*}(X)\omega^{p*}(Y)
        \ar@{=}[rr]
        &\text{ }
        &p \omega^{p*}(X)\omega^{p*}(Y)
      }
    }}
  \end{equation*}
  We also obtain that $I_p\otimes \omega^{p*}_{\unit}$
  is an isomorphism
  from the diagrams below.
  \begin{equation*}
    \vcenter{\hbox{
      \xymatrix@C=40pt{
        p\otimes \kappa
        \ar[r]^-{\kos{p}_!(\xi^p_{\unit})}_-{\cong}
        \ar@/_1pc/@{=}[dr]
        &p\otimes \omega^{p*}(\unit)
        \ar[d]^-{I_p\otimes \omega^{p*}_{\unit}}
        \\
        \text{ }
        &p\otimes \kappa
      }
    }}
    \qquad\quad
    \vcenter{\hbox{
      \xymatrix{
        p\kappa
        \ar[dd]_-{\kos{p}_!(\xi^p_{\unit})}^-{\cong}
        \ar@{=}[rr]
        &\text{ }
        &p\kappa
        \ar@/_0.5pc/[dl]|-{\cp_pI_{\kappa}}
        \ar@{=}[ddd]
        \\
        \text{ }
        &pp\kappa
        \ar@{->>}@/_0.5pc/[dl]|-{I_p\xi^p_{\unit}}
        \ar@/_1pc/[ddr]|-{I_pe_pI_{\kappa}}
        &\text{ }
        \\
        p\omega^{p*}(\unit)
        \ar[d]_-{I_p\omega^{p*}_{\unit}}
        &\text{ }
        &\text{ }
        \\
        p\kappa
        \ar@{=}[rr]
        &\text{ }
        &p\kappa
      }
    }}
  \end{equation*}
  This shows that
  $\omega^{p*}_{X,Y}$, $\omega^{p*}_{\unit}$
  are isomorphisms.
  Next we verify the comonoidal coherence relations
  of $\omega^{p*}$.
  Let $Z=(z,\gamma_z)$
  be another representation of $\pi$.
  We have
  \begin{equation*}
    \vcenter{\hbox{
      \xymatrix@C=60pt{
        \omega^{p*}(X\tensor (Y\tensor Z))
        \ar[d]_-{\omega^{p*}_{X,Y\tensor Z}}^-{\cong}
        \ar[r]^-{\omega^{p*}(a_{X,Y,Z})}_-{\cong}
        &\omega^{p*}((X\tensor Y)\tensor Z)
        \ar[d]^-{\omega^{p*}_{X\tensor Y,Z}}_-{\cong}
        \\
        \omega^{p*}(X)\otimes \omega^{p*}(Y\tensor Z)
        \ar[d]_-{I_{\omega^{p*}(X)}\otimes \omega^{p*}_{Y,Z}}^-{\cong}
        &\omega^{p*}(X\tensor Y)\otimes \omega^{p*}(Z)
        \ar[d]^-{\omega^{p*}_{X,Y}\otimes I_{\omega^{p*}(Z)}}_-{\cong}
        \\
        \omega^{p*}(X)\otimes (\omega^{p*}(Y)\otimes \omega^{p*}(Z))
        \ar[r]^-{a_{\omega^{p*}(X),\omega^{p*}(Y),\omega^{p*}(Z)}}_-{\cong}
        &(\omega^{p*}(X)\otimes \omega^{p*}(Y))\otimes \omega^{p*}(Z)
      }
    }}
  \end{equation*}
  which we obtain by right-cancelling the epimorphism
  $\xi^p_{X(YZ)}$ in the diagram below.
  \begin{equation*}
    \vcenter{\hbox{
      \xymatrix@C=-5pt{
        p(x(yz))
        \ar@{->>}[d]_-{\xi^p_{X(YZ)}}
        \ar@{=}[r]
        &p(x(yz))
        \ar[d]^-{(p\otimes)_{x,yz}}
        \ar@{=}[r]
        &p(x(yz))
        \ar[d]^-{I_pa_{x,y,z}}_-{\cong}
        \ar@{=}[r]
        &p(x(yz))
        \ar@{->>}[d]^-{\xi^p_{X(YZ)}}
        \\
        \omega^{p*}(X(YZ))
        \ar[d]_-{\omega^{p*}_{X,YZ}}^-{\cong}
        &(pz)(p(yz))
        \ar@{->>}@/^0.5pc/[dl]|-{\xi^p_X\xi^p_{YZ}}
        \ar[d]^-{I_{px}(p\otimes)_{y,z}}
        \ar@{}[r]|-{(\dagger)}
        &p((xy)z)
        \ar[d]^-{(p\otimes)_{xy,z}}
        \ar@{->>}@/^0.5pc/[dr]|-{\xi^p_{(XY)Z}}
        &\omega^{p*}(X(YZ))
        \ar[d]^-{\omega^{p*}(a_{X,Y,Z})}_-{\cong}
        \\
        \omega^{p*}(X)\omega^{p*}(YZ)
        \ar[d]_-{I_{\omega^{p*}(X)}\omega^{p*}_{Y,Z}}^-{\cong}
        &(px)((py)(pz))
        \ar@{->>}@/^0.5pc/[dl]|-{\xi^p_X(\xi^p_Y\xi^p_Z)}
        \ar@/_0.5pc/[dr]_-{a_{px,py,pz}}^-{\cong}
        &(p(xy))(pz)
        \ar[d]^-{(p\otimes)_{x,y}I_{pz}}
        \ar@{->>}@/^0.5pc/[dr]|-{\xi^p_{XY}\xi^p_Z}
        &\omega^{p*}((XY)Z)
        \ar[d]^-{\omega^{p*}_{XY,Z}}_-{\cong}
        \\
        \omega^{p*}(X)(\omega^{p*}(Y)\omega^{p*}(Z))
        \ar[d]^-{a_{\omega^{p*}(X),\omega^{p*}(Y),\omega^{p*}(Z)}}_-{\cong}
        &\text{ }
        &((px)(py))(pz)
        \ar@{->>}@/_0.5pc/[dr]|-{(\xi^p_X\xi^p_Y)\xi^p_Z}
        &\omega^{p*}(XY)\omega^{p*}(Z)
        \ar[d]^-{\omega^{p*}_{X,Y}I_{\omega^{p*}(Z)}}_-{\cong}
        \\
        (\omega^{p*}(X)\omega^{p*}(Y))\omega^{p*}(Z)
        \ar@{=}[rrr]
        &\text{ }
        &\text{ }
        &(\omega^{p*}(X)\omega^{p*}(Y))\omega^{p*}(Z)
      }
    }}
  \end{equation*}
  In the diagram $(\dagger)$,
  we used the associativity of $\cp_p:p\to p\otimes p$.
  We also have
  \begin{equation*}
    \vcenter{\hbox{
      \xymatrix@C=40pt{
        \omega^{p*}(X)
        \ar[r]^-{\omega^{p*}(\imath_X)}_-{\cong}
        \ar[d]_-{\imath_{\omega^{p*}(X)}}^-{\cong}
        &\omega^{p*}(\unit\tensor X)
        \ar[d]^-{\omega^{p*}_{\unit,X}}_-{\cong}
        \\
        \kappa\otimes \omega^{p*}(X)
        &\omega^{p*}(\unit)\otimes \omega^{p*}(X)
        \ar[l]_-{\omega^{p*}_{\unit}\otimes I_{\omega^{p*}(X)}}^-{\cong}
      }
    }}
  \end{equation*}
  which we obtain by right-cancelling the epimorphism
  $\xi^p_X$
  in the diagram below.
  \begin{equation*}
    \vcenter{\hbox{
      \xymatrix@C=60pt{
        px
        \ar@{->>}[d]_-{\xi^p_X}
        \ar@{=}[r]
        &px
        \ar[d]^-{I_p\imath_x}_-{\cong}
        \ar@{=}[r]
        &px
        \ar@/^1.8pc/[dddl]_-{\imath_{px}}^-{\cong}
        \ar@{->>}[dd]^-{\xi^p_X}
        \\
        \omega^{p*}(X)
        \ar[d]_-{\omega^{p*}(\imath_X)}^-{\cong}
        &p\kappa x
        \ar@{->>}@/^0.5pc/[dl]|-{\xi^p_{\unit\tensor X}}
        \ar[d]^-{(p\otimes)_{\kappa,x}}
        &\text{ }
        \\
        \omega^{p*}(\unit X)
        \ar[d]_-{\omega^{p*}_{\unit,X}}^-{\cong}
        &p\kappa px
        \ar@{->>}@/^0.5pc/[dl]|-{\xi^p_{\unit}\xi^p_X}
        \ar[d]^-{e_pI_{\kappa px}}
        &\omega^{p*}(X)
        \ar[dd]^-{\imath_{\omega^{p*}(X)}}_-{\cong}
        \\
        \omega^{p*}(\unit)\omega^{p*}(X)
        \ar[d]_-{\omega^{p*}_{\unit} I_{\omega^{p*}(X)}}^-{\cong}
        &\kappa px
        \ar@{->>}[d]^-{I_{\kappa}\xi^p_X}
        &\text{ }
        \\
        \kappa\omega^{p*}(X)
        \ar@{=}[r]
        &\kappa\omega^{p*}(X)
        \ar@{=}[r]
        &\kappa\omega^{p*}(X)
      }
    }}
  \end{equation*}
  Similarly, we have
  \begin{equation*}
    \vcenter{\hbox{
      \xymatrix@C=40pt{
        \omega^{p*}(X)
        \ar[r]^-{\omega^{p*}(\jmath_X)}_-{\cong}
        \ar[d]_-{\jmath_{\omega^{p*}(X)}}^-{\cong}
        &\omega^{p*}(X\tensor \unit)
        \ar[d]^-{\omega^{p*}_{X,\unit}}_-{\cong}
        \\
        \omega^{p*}(X)\otimes \kappa
        &\omega^{p*}(X)\otimes \omega^{p*}(\unit)
        \ar[l]_-{I_{\omega^{p*}(X)}\otimes \omega^{p*}_{\unit}}^-{\cong}
      }
    }}
  \end{equation*}
  which we obtain by right-cancelling the epimorphism
  $\xi^p_X$
  in the diagram below.
  \begin{equation*}
    \vcenter{\hbox{
      \xymatrix@C=60pt{
        px
        \ar@{->>}[d]_-{\xi^p_X}
        \ar@{=}[r]
        &px
        \ar[d]^-{I_p\jmath_x}_-{\cong}
        \ar@{=}[r]
        &px
        \ar@/^1.8pc/[dddl]_-{\jmath_{px}}^-{\cong}
        \ar@{->>}[dd]^-{\xi^p_X}
        \\
        \omega^{p*}(X)
        \ar[d]_-{\omega^{p*}(\jmath_X)}^-{\cong}
        &px\kappa
        \ar@{->>}@/^0.5pc/[dl]|-{\xi^p_{X\unit}}
        \ar[d]^-{(p\otimes)_{x,\kappa}}
        &\text{ }
        \\
        \omega^{p*}(X\unit)
        \ar[d]_-{\omega^{p*}_{X,\unit}}^-{\cong}
        &pxp\kappa
        \ar@{->>}@/^0.5pc/[dl]|-{\xi^p_X\xi^p_{\unit}}
        \ar[d]^-{I_{px}e_pI_{\kappa}}
        &\omega^{p*}(X)
        \ar[dd]^-{\jmath_{\omega^{p*}(X)}}_-{\cong}
        \\
        \omega^{p*}(X)\omega^{p*}(\unit)
        \ar[d]_-{I_{\omega^{p*}(X)}\omega^{p*}_{\unit}}^-{\cong}
        &px\kappa
        \ar@{->>}[d]^-{\xi^p_XI_{\kappa}}
        &\text{ }
        \\
        \omega^{p*}(X)\kappa
        \ar@{=}[r]
        &\omega^{p*}(X)\kappa
        \ar@{=}[r]
        &\omega^{p*}(X)\kappa
      }
    }}
  \end{equation*}
  Finally we have
  \begin{equation*}
    \vcenter{\hbox{
      \xymatrix@C=50pt{
        \omega^{p*}(X\tensor Y)
        \ar[d]_-{\omega^{p*}_{X,Y}}^-{\cong}
        \ar[r]^-{\omega^{p*}(s_{X,Y})}_-{\cong}
        &\omega^{p*}(Y\tensor X)
        \ar[d]^-{\omega^{p*}_{Y,X}}_-{\cong}
        \\
        \omega^{p*}(X)\otimes \omega^{p*}(Y)
        \ar[r]^-{s_{\omega^{p*}(X),\omega^{p*}(Y)}}_-{\cong}
        &\omega^{p*}(Y)\otimes \omega^{p*}(X)
      }
    }}
  \end{equation*}
  which we obtain by right-cancelling the epimorphism
  $\xi^p_{X\tensor Y}$
  in the diagram below.
  \begin{equation*}
    \vcenter{\hbox{
      \xymatrix@C=35pt{
        pxy
        \ar@{->>}[d]_-{\xi^p_{XY}}
        \ar@{=}[r]
        &pxy
        \ar[d]^-{(p\otimes)_{x,y}}
        \ar@{=}[r]
        &pxy
        \ar[d]^-{I_ps_{x,y}}_-{\cong}
        \ar@{=}[r]
        &pxy
        \ar@{->>}[d]^-{\xi^p_{XY}}
        \\
        \omega^{p*}(XY)
        \ar[d]_-{\omega^{p*}_{X,Y}}^-{\cong}
        &pxpy
        \ar@{->>}@/^0.5pc/[dl]|-{\xi^p_X\xi^p_Y}
        \ar[d]^-{s_{px,py}}_-{\cong}
        \ar@{}[r]|-{(\dagger\!\dagger)}
        &pyx
        \ar@/^0.5pc/[dl]|-{(p\otimes)_{y,x}}
        \ar@{->>}@/_0.5pc/[dr]|-{\xi^p_{YX}}
        &\omega^{p*}(XY)
        \ar[d]^-{\omega^{p*}(s_{X,Y})}_-{\cong}
        \\
        \omega^{p*}(X)\omega^{p*}(Y)
        \ar[d]_-{s_{\omega^{p*}(X),\omega^{p*}(Y)}}^-{\cong}
        &pypx
        \ar@{->>}[d]^-{\xi^p_Y\xi^p_X}
        &\text{ }
        &\omega^{p*}(YX)
        \ar[d]^-{\omega^{p*}_{Y,X}}_-{\cong}
        \\
        \omega^{p*}(Y)\omega^{p*}(X)
        \ar@{=}[r]
        &\omega^{p*}(Y)\omega^{p*}(X)
        \ar@{=}[rr]
        &\text{ }
        &\omega^{p*}(Y)\omega^{p*}(X)
      }
    }}
  \end{equation*}
  In the diagram $(\dagger\!\dagger)$,
  we used the cocommutativity of $\cp_p:p\to p\otimes p$.
  We conclude that $\omega^{p*}:\kos{R\!e\!p}(\pi)\to \kos{K}$
  is a strong symmetric monoidal functor as we claimed.
  Now we show that
  $(\omega^{p*},\what{\omega}^{p*})$
  is a strong $\kos{K}$-tensor functor.
  Let $z$ be an object in $\CK$.
  The following diagram shows that 
  the natural isomorphism
  $\what{\omega}^{p*}$ is well-defined.
  \begin{equation*}
    \vcenter{\hbox{
      \xymatrix@C=50pt{
        p\otimes \pi\otimes z
        \ar[d]_-{\tau_p\otimes I_z}^-{\cong}
        \ar@<0.5ex>[r]^-{\lambda_p\otimes I_z}
        \ar@<-0.5ex>[r]_-{I_p\otimes e_{\pi}\otimes I_z}
        &p\otimes z
        \ar@{=}[d]
        \ar@{->>}[r]^-{\xi^p_{\kos{t}_{\pi}^*(z)}}
        &\omega^{p*}\kos{t}_{\pi}^*(z)
        \ar@{.>}[d]^-{\exists!\text{ }\what{\omega}^{p*}_z}_-{\cong}
        \\
        p\otimes p\otimes z
        \ar@<0.5ex>[r]^-{e_p\otimes I_{p\otimes z}}
        \ar@<-0.5ex>[r]_-{I_p\otimes e_p\otimes I_z}
        &p\otimes z
        \ar[r]^-{e_p\otimes I_z}
        &z
      }
    }}
  \end{equation*}
  We need to show that
  $\what{\omega}^{p*}$ is a comonoidal natural isomorphism.
  Let $w$ be another object in $\CK$.
  We have
  \begin{equation*}
    \vcenter{\hbox{
      \xymatrix@C=30pt{
        \omega^{p*}\kos{t}_{\pi}^*(z\otimes w)
        \ar[d]_-{\what{\omega}^{p*}_{z\otimes w}}^-{\cong}
        \ar@{=}[r]
        &\omega^{p*}(\kos{t}_{\pi}^*(z)\tensor \kos{t}_{\pi}^*(w))
        \ar[r]^-{\omega^{p*}_{\kos{t}^*_{\pi}(z),\kos{t}^*_{\pi}(w)}}_-{\cong}
        &\omega^{p*}\kos{t}_{\pi}^*(z)\otimes \omega^{p*}\kos{t}_{\pi}^*(w)
        \ar[d]^-{\what{\omega}^{p*}_z\otimes \what{\omega}^{p*}_w}_-{\cong}
        \\
        z\otimes w
        \ar@{=}[rr]
        &\text{}
        &z\otimes w
      }
    }}
  \end{equation*}
  which we obtain after right-cancelling the epimorphism
  $\xi^p_{\kos{t}^*_{\pi}(z)\tensor \kos{t}^*_{\pi}(w)}
  =\xi^p_{\kos{t}_{\pi}^*(z\otimes w)}$
  in the diagram below.
  \begin{equation*}
    \vcenter{\hbox{
      \xymatrix@C=50pt{
        p\otimes z\otimes w
        \ar@{->>}[d]_-{\xi^p_{\kos{t}_{\pi}^*(z)\tensor \kos{t}_{\pi}^*(w)}}
        \ar@{=}[r]
        &p\otimes z\otimes w
        \ar[d]^-{(p\otimes)_{z,w}}
        \ar@{=}[r]
        &p\otimes z\otimes w
        \ar@/^0.5pc/[dddl]|-{e_p\otimes I_{z\otimes w}}
        \ar@{->>}[d]^-{\xi^p_{\kos{t}_{\pi}^*(z\otimes w)}}
        \\
        \omega^{p*}(\kos{t}^*_{\pi}(z)\tensor \kos{t}^*_{\pi}(w))
        \ar[d]_-{\omega^{p*}_{\kos{t}^*_{\pi}(z),\kos{t}^*_{\pi}(w)}}^-{\cong}
        &p\otimes z\otimes p\otimes w
        \ar@{->>}@/^0.5pc/[dl]|-{\xi^p_{\kos{t}^*_{\pi}(z)}\otimes \xi^p_{\kos{t}^*_{\pi}(w)}}
        \ar[dd]|-{e_p\otimes I_z\otimes e_p\otimes I_w}
        &\omega^{p*}\kos{t}^*_{\pi}(z\otimes w)
        \ar[dd]^-{\what{\omega}^{p*}_{z\otimes w}}_-{\cong}
        \\
        \omega^{p*}\kos{t}^*_{\pi}(z)\otimes \omega^{p*}\kos{t}^*_{\pi}(w)
        \ar[d]_-{\what{\omega}^{p*}_z\otimes \what{\omega}^{p*}_w}^-{\cong}
        &\text{ }
        \\
        z\otimes w
        \ar@{=}[r]
        &z\otimes w
        \ar@{=}[r]
        &z\otimes w
      }
    }}
  \end{equation*}
  We conclude that
  $(\omega^{p*},\what{\omega}^{p*})
  :(\kos{R\!e\!p}(\pi),\kos{t}_{\pi}^*)
  \to (\kos{K},\id_{\kos{K}})$
  is a strong $\kos{K}$-tensor functor.
  We are left to verify the characterization
  of the $\kos{K}$-equivariance
  $\vecar{\omega}^{p*}$.
  Let $X$ be a representation of $\pi$ and let $z$ be an object in $\CK$.
  Let us denote $X\otimes z=X\tensor\!_{\pi}\kos{t}^*_{\pi}(z)$.
  The $\kos{K}$-equivariance $\vecar{\omega}^{p*}$ is given by
  \begin{equation*}
    \vecar{\omega}^{p*}_{z,X}:
    \xymatrix@C=40pt{
      \omega^{p*}(X\otimes z)
      \ar[r]^-{\omega^{p*}_{X,\kos{t}_{\pi}^*(w)}}_-{\cong}
      &\omega^{p*}(X)\otimes \omega^{p*}\kos{t}_{\pi}^*(w)
      \ar[r]^-{I_{\omega^{p*}(X)}\otimes \what{\omega}^{p*}_w}_-{\cong}
      &\omega^{p*}(X)\otimes z
      .
    }
  \end{equation*}
  We obtain the description of $\vecar{\omega}^{p*}$
  after right-cancelling the epimorphism
  $\xi^p_{X\otimes z}$ in the diagram below.
  \begin{equation*}
    \vcenter{\hbox{
      \xymatrix@C=60pt{
        p\otimes (x\otimes z)
        \ar@{->>}[d]_-{\xi^p_{X\otimes z}}
        \ar@{=}[r]
        &p\otimes (x\otimes z)
        \ar[d]^-{(p\otimes)_{x,z}}
        \ar@{=}[r]
        &p\otimes (x\otimes z)
        \ar[dd]^-{a_{p,x,z}}_-{\cong}
        \\
        \omega^{p*}(X\otimes z)
        \ar[dd]_-{\vecar{\omega}^{p*}_{z,X}}^-{\cong}
        \ar@/_0.5pc/[dr]^-{\omega^{p*}_{X,\kos{t}^*_{\pi}(z)}}_-{\cong}
        &(p\otimes x)\otimes (p\otimes z)
        \ar@{->>}[d]^-{\xi^p_X\otimes \xi^p_{\kos{t}^*_{\pi}(z)}}
        \ar@/^0.5pc/[dr]|-{I_{p\otimes x}\otimes (e_p\otimes I_z)}
        &\text{ }
        \\
        \text{ }
        &\omega^{p*}(X)\otimes \omega^{p*}\kos{t}_{\pi}^*(z)
        \ar@/_0.5pc/[dr]^-{I_{\omega^{p*}(X)}\otimes \what{\omega}^{p*}_z}_-{\cong}
        &(p\otimes x)\otimes z
        \ar@{->>}[d]^-{\xi^p_X\otimes I_z}
        \\
        \omega^{p*}(X)\otimes z
        \ar@{=}[rr]
        &\text{ }
        &\omega^{p*}(X)\otimes z
      }
    }}
  \end{equation*}
  This completes the proof of Lemma~\ref{lem twistbyatorsor}.
\qed\end{proof}

\begin{proposition} \label{prop twistbyatorsor}
  Let $\pi$, $\pi^{\pr}$ be pre-Galois objects in $\kos{K}$
  and $p$
  be a $(\pi^{\pr},\pi)$-bitorsor over $\kappa$.
  \begin{enumerate}
    \item 
    We have a pre-fiber functor
    \begin{equation*}
      \varpi^p=(\omega^p,\what{\omega}^p)
      :\mathfrak{K}\to \mathfrak{Rep}(\pi)      
    \end{equation*}
    which is obtained by twisting the pre-fiber functor
    $\varpi_{\pi}=(\omega_{\pi},\what{\omega}_{\pi})
    :\mathfrak{K}\to\mathfrak{Rep}(\pi)$.
    The right adjoint strong $\kos{K}$-tensor functor
    $(\omega^{p*},\what{\omega}^{p*})$
    of $\varpi^p$
    is described in Lemma~\ref{lem twistbyatorsor}.

    \item
    The twisted fiber functor
    $\varpi^p:\mathfrak{K}\to\mathfrak{Rep}(\pi)$
    factors through as an equivalence of Galois $\kos{K}$-prekosmoi
    $\widebreve{\varpi}\!^p:
    \mathfrak{Rep}(\pi^{\pr})
    \xrightarrow{\simeq}
    \mathfrak{Rep}(\pi)$.
    \begin{equation*}
      \vcenter{\hbox{
        \xymatrix@C=30pt{
          \mathfrak{Rep}(\pi)
          &\mathfrak{Rep}(\pi^{\pr})
          \ar[l]_-{\widebreve{\varpi}\!^p}^-{\simeq}
          \\
          \text{ }
          &\mathfrak{K}
          \ar[u]_-{\varpi_{\pi^{\pr}}}
          \ar@/^1pc/[ul]^-{\varpi^p}
        }
      }}
    \end{equation*}
  \end{enumerate}
\end{proposition}
\begin{proof}
  It suffices to show that the strong $\kos{K}$-tensor functor
  $(\omega^{p*},\what{\omega}^{p*})$
  described in Lemma~\ref{lem twistbyatorsor}
  factors through as an equivalence of strong $\kos{K}$-tensor categories
  $(\widebreve{\omega}\!^{p*},\what{\widebreve{\omega}}^{p*}):
  (\kos{R\!e\!p}(\pi),\kos{t}_{\pi}^*)
  \xrightarrow{\simeq}
  (\kos{R\!e\!p}(\pi^{\pr}),\kos{t}_{\pi^{\pr}}^*)$.
  \begin{equation*}
    \vcenter{\hbox{
      \xymatrix@C=40pt{
        (\kos{R\!e\!p}(\pi),\kos{t}_{\pi}^*)
        \ar[r]^-{(\widebreve{\omega}\!^{p*},\what{\widebreve{\omega}}^{p*})}_-{\simeq}
        \ar@/_1pc/[dr]_-{(\omega^{p*},\what{\omega}^{p*})}
        &(\kos{R\!e\!p}(\pi^{\pr}),\kos{t}_{\pi^{\pr}}^*)
        \ar[d]^-{(\omega_{\pi^{\pr}},\what{\omega}_{\pi^{\pr}})}
        \\
        \text{ }
        &(\kos{K},\id_{\kos{K}})
      }
    }}
  \end{equation*}
  The underlying functor $\widebreve{\omega}\!^{p*}$
  sends each object $X=(x,\gamma_x)$ in $\cat{Rep}(\pi)$ to
  \begin{equation*}
    \widebreve{\omega}\!^{p*}(X)=
    \bigl(
      \omega^{p*}(X)
      ,\text{ }
      \gamma^{\pr}_{\omega^{p*}(X)}
      :\pi^{\pr}\otimes \omega^{p*}(X)\to \omega^{p*}(X)
    \bigr)
  \end{equation*}
  whose left $\pi^{\pr}$-action $\gamma^{\pr}_{\omega^{p*}(X)}$ is
  the unique morphism satisfying the relation below.
  \begin{equation*}
    \vcenter{\hbox{
      \xymatrix@C=40pt{
        \pi^{\pr}\otimes p\otimes x
        \ar[d]_-{\gamma^{\pr}_p\otimes I_x}
        \ar@{->>}[r]^-{I_{\pi^{\pr}}\otimes \xi^p_X}
        &\pi^{\pr}\otimes \omega^{p*}(X)
        \ar@{.>}[d]^-{\gamma^{\pr}_{\omega^{p*}(X)}}_-{\exists!}
        \\
        p\otimes x
        \ar@{->>}[r]^-{\xi^p_X}
        &\omega^{p*}(X)
      }
    }}
  \end{equation*}
  The following diagram show that the morphism $\gamma^{\pr}_{\omega^{p*}(X)}$ is well-defined.
  \begin{equation*}
    \vcenter{\hbox{
      \xymatrix@C=40pt{
        \pi^{\pr}p\pi x
        \ar[d]_-{I_{\pi^{\pr}p}\gamma_x}
        \ar@{=}[r]
        &\pi^{\pr}p\pi x
        \ar[d]^-{\gamma^{\pr}_pI_{\pi x}}
        \ar@{=}[r]
        &\pi^{\pr}p\pi x
        \ar[d]^-{I_{\pi^{\pr}}\lambda_pI_x}
        \\
        \pi^{\pr}px
        \ar[d]_-{\gamma^{\pr}_pI_x}
        &p\pi x
        \ar@/^0.5pc/[dl]|-{I_p\gamma_x}
        \ar@/_0.5pc/[dr]|-{\lambda_pI_x}
        &\pi^{\pr}px
        \ar[d]^-{\gamma^{\pr}_pI_x}
        \\
        px
        \ar@{->>}[d]_-{\xi^p_X}
        &\text{ }
        &px
        \ar@{->>}[d]_-{\xi^p_X}
        \\
        \omega^{p*}(X)
        \ar@{=}[rr]
        &\text{ }
        &\omega^{p*}(X)
      }
    }}
  \end{equation*}
  To conclude that the underlying functor
  $\widebreve{\omega}\!^{p*}$
  is well-defined,
  we need to check that $\gamma^{\pr}_{\omega^{p*}(X)}$
  satisfies the left $\pi^{\pr}$-action relations.
  We obtain this by right-cancelling the epimorphisms
  $I_{\pi^{\pr}\otimes \pi^{\pr}}\otimes \xi^p_{\omega^{p*}(X)}$
  and
  $\xi^p_X$
  in each of the following diagrams.
  \begin{equation*}
    \vcenter{\hbox{
      \xymatrix{
        \pi^{\pr}\pi^{\pr}px
        \ar@{->>}[d]_-{I_{\pi^{\pr}\pi^{\pr}}\xi^p_{\omega^{p*}(X)}}
        \ar@{=}[r]
        &\pi^{\pr}\pi^{\pr}px
        \ar[d]^-{I_{\pi^{\pr}}\gamma^{\pr}_pI_x}
        \ar@{=}[r]
        &\pi^{\pr}\pi^{\pr}px
        \ar[d]^-{\pc_{\pi^{\pr}}I_x}
        \ar@{=}[r]
        &\pi^{\pr}\pi^{\pr}px
        \ar@{->>}[d]^-{I_{\pi^{\pr}\pi^{\pr}}\xi^p_{\omega^{p*}(X)}}
        \\
        \pi^{\pr}\pi^{\pr}\omega^{p*}(X)
        \ar[d]_-{I_{\pi^{\pr}}\gamma^{\pr}_{\omega^{p*}(X)}}
        &\pi^{\pr}px
        \ar@{->>}@/^0.5pc/[dl]|-{I_{\pi^{\pr}}\xi^p_X}
        \ar[d]^-{\gamma^{\pr}_pI_x}
        &\pi^{\pr}px
        \ar@/^0.5pc/[dl]|-{\gamma^{\pr}_pI_x}
        \ar@/_0.5pc/[dr]|-{I_{\pi^{\pr}}\xi^p_X}
        &\pi^{\pr}\pi^{\pr}\omega^{p*}(X)
        \ar[d]^-{\pc_{\pi^{\pr}}I_{\omega^{p*}(X)}}
        \\
        \pi^{\pr}\omega^{p*}(X)
        \ar[d]_-{\gamma^{\pr}_{\omega^{p*}(X)}}
        &px
        \ar@{->>}[d]^-{\xi^p_X}
        &\text{ }
        &\pi^{\pr}\omega^{p*}(X)
        \ar[d]^-{\gamma^{\pr}_{\omega^{p*}(X)}}
        \\
        \pi^{\pr}\omega^{p*}(X)
        \ar@{=}[r]
        &\pi^{\pr}\omega^{p*}(X)
        \ar@{=}[rr]
        &\text{ }
        &\pi^{\pr}\omega^{p*}(X)
      }
    }}
  \end{equation*}
  \begin{equation*}
    \vcenter{\hbox{
      \xymatrix{
        px
        \ar@{->>}[d]_-{\xi^p_X}
        \ar@{=}[r]
        &px
        \ar[d]^-{u_{\pi^{\pr}}I_{px}}
        \ar@{=}[r]
        &px
        \ar@{=}[dd]
        \\
        \omega^{p*}(X)
        \ar[d]_-{u_{\pi^{\pr}}I_{\omega^{p*}(X)}}
        &\pi^{\pr}px
        \ar@{->>}@/^0.5pc/[dl]|-{I_{\pi^{\pr}}\xi^p_X}
        \ar@/_0.5pc/[dr]|-{\gamma^{\pr}_pI_x}
        &\text{ }
        \\
        \pi^{\pr}\omega^{p*}(X)
        \ar[d]_-{\gamma^{\pr}_{\omega^{p*}(X)}}
        &\text{ }
        &px
        \ar@{->>}[d]^-{\xi^p_X}
        \\
        \omega^{p*}(X)
        \ar@{=}[rr]
        &\text{ }
        &\omega^{p*}(X)
      }
    }}
  \end{equation*}
  This shows that the underlying functor
  $\widebreve{\omega}\!^{p*}$ is well-defined.
  We also obtain that the functor
  $\omega^{p*}$ factors through the functor
  $\widebreve{\omega}\!^{p*}$.

  Next we show that
  $(\widebreve{\omega}\!^{p*},\what{\widebreve{\omega}}^{p*})$
  is a strong $\kos{K}$-tensor functor,
  and the strong $\kos{K}$-tensor functor
  $(\omega^{p*},\what{\omega}^{p*})$ factors through
  $(\widebreve{\omega}\!^{p*},\what{\widebreve{\omega}}^{p*})$.
  Let $X=(x,\gamma_x)$ and $Y=(y,\gamma_y)$
  be representations of $\pi$
  and let $z$ be an object in $\CK$.
  It suffices to show that 
  the isomorphisms
  \begin{equation*}
    \begin{aligned}
      \omega^{p*}_{X,Y}
      &:
      \omega^{p*}(X\tensor\!_{\pi} Y)
      \xrightarrow{\cong}
      \omega^{p*}(X)\otimes \omega^{p*}(Y)
      \\
      \omega^{p*}_{\unit\!_{\pi}}
      &:
      \omega^{p*}(\unit\!_{\pi})
      \xrightarrow{\cong}
      \kappa
      \\
      \what{\omega}^{p*}_z
      &:
      \omega^{p*}\kos{t}_{\pi}^*(z)
      \xrightarrow{\cong}
      z
    \end{aligned}
  \end{equation*}
  in $\CK$ become isomorphisms 
  \begin{equation*}
    \begin{aligned}
      \widebreve{\omega}\!^{p*}_{X,Y}
      &:
      \widebreve{\omega}\!^{p*}(X\tensor\!_{\pi} Y)
      \xrightarrow{\cong}
      \widebreve{\omega}\!^{p*}(X)\tensor\!_{\pi^{\pr}}
      \widebreve{\omega}\!^{p*}(Y)
      \\
      \widebreve{\omega}\!^{p*}_{\unit}
      &:
      \widebreve{\omega}\!^{p*}(\unit)
      \xrightarrow{\cong}
      \unit_{\pi^{\pr}}
      \\
      \what{\widebreve{\omega}}^{p*}_z
      &:
      \widebreve{\omega}\!^{p*}\kos{t}_{\pi}^*(z)
      \xrightarrow{\cong}
      \kos{t}_{\pi^{\pr}}^*(z)
    \end{aligned}
  \end{equation*}
  in $\kos{R\!e\!p}(\pi^{\pr})$.
  We obtain this by right-cancelling
  the epimorphisms
  $I_{\pi^{\pr}}\otimes \xi^p_{X\tensor\!_{\pi}Y}$,
  $I_{\pi^{\pr}}\otimes \xi^p_{\unit\!_{\pi}}$
  and
  $I_{\pi^{\pr}}\otimes \xi^p_{\kos{t}^*_{\pi}(z)}$
  in each of the following diagrams.
  \begin{equation*}
    \vcenter{\hbox{
      \xymatrix@C=15pt{
        \pi^{\pr}pxy
        \ar@{->>}[d]_-{I_{\pi^{\pr}}\xi^p_{XY}}
        \ar@{=}[r]
        &\pi^{\pr}pxy
        \ar[d]^-{I_{\pi^{\pr}}(p\otimes)_{x,y}}
        \ar@{=}[r]
        &\pi^{\pr}pxy
        \ar@/^1.5pc/[ddl]|-{((\pi^{\pr}p)\otimes)_{x,y}}
        \ar[dd]^-{\gamma^{\pr}_pI_{xy}}
        \ar@{=}[r]
        &\pi^{\pr}pxy
        \ar@{->>}[d]^-{I_{\pi^{\pr}}\xi^p_{XY}}
        \\
        \pi^{\pr}\omega^{p*}(XY)
        \ar[d]_-{I_{\pi^{\pr}}\omega^{p*}_{X,Y}}^-{\cong}
        &\pi^{\pr}pxpy
        \ar@{->>}@/_0.5pc/[dl]|-{I_{\pi^{\pr}}\xi^p_X\xi^p_Y}
        \ar[d]^-{(\pi^{\pr}\otimes)_{px,py}}
        &\text{ }
        &\pi^{\pr}\omega^{p*}(XY)
        \ar[dd]^-{\gamma^{\pr}_{\omega^{p*}(XY)}}
        \\
        \pi^{\pr}\omega^{p*}(X)\omega^{p*}(Y)
        \ar[dd]_-{\gamma^{\pr}_{\omega^{p*}(X)\omega^{p*}(Y)}}
        \ar@/_0.5pc/[dr]|-{(\pi^{\pr}\otimes)_{\omega^{p*}(X),\omega^{p*}(Y)}}
        &\pi^{\pr}px\pi^{\pr}py
        \ar@{->>}[d]^-{I_{\pi^{\pr}}\xi^p_X I_{\pi^{\pr}}\xi^p_Y}
        \ar@/^1pc/[dr]|-{\gamma^{\pr}_pI_x\gamma^{\pr}_pI_y}
        &pxy
        \ar[d]^-{(p\otimes)_{x,y}}
        \ar@{->>}@/^0.5pc/[dr]|-{\xi^p_{XY}}
        &\text{ }
        \\
        \text{ }
        &\pi^{\pr}\omega^{p*}(X)\pi^{\pr}\omega^{p*}(X)
        \ar[d]^-{\gamma^{\pr}_{\omega^{p*}(X)}\gamma^{\pr}_{\omega^{p*}(Y)}}
        &pxpy
        \ar@{->>}[d]^-{\xi^p_X\xi^p_Y}
        &\omega^{p*}(XY)
        \ar[d]^-{\omega^{p*}_{X,Y}}_-{\cong}
        \\
        \omega^{p*}(X)\omega^{p*}(Y)
        \ar@{=}[r]
        &\omega^{p*}(X)\omega^{p*}(Y)
        \ar@{=}[r]
        &\omega^{p*}(X)\omega^{p*}(Y)
        \ar@{=}[r]
        &\omega^{p*}(X)\omega^{p*}(Y)
      }
    }}
  \end{equation*}
  \begin{equation*}
    \vcenter{\hbox{
      \xymatrix@C=50pt{
        \pi^{\pr}p \kappa
        \ar@{->>}[d]_-{I_{\pi^{\pr}}\xi^p_{\unit}}
        \ar@{=}[r]
        &\pi^{\pr}p \kappa
        \ar@/^1pc/[ddl]|-{I_{\pi^{\pr}}e_pI_{\kappa}}
        \ar[ddd]^-{e_{\pi^{\pr}p}I_{\kappa}}
        \ar@{=}[r]
        &\pi^{\pr}p \kappa
        \ar[d]^-{\gamma^{\pr}_pI_{\kappa}}
        \ar@{=}[r]
        &\pi^{\pr}p \kappa
        \ar@{->>}[d]^-{I_{\pi^{\pr}}\xi^p_{\unit}}
        \\
        \pi^{\pr} \omega^{p*}(\unit)
        \ar[d]_-{I_{\pi^{\pr}}\widebreve{\omega}^{p*}_{\unit}}^-{\cong}
        &\text{ }
        &p\kappa
        \ar[dd]^-{e_pI_{\kappa}}
        \ar@{->>}@/^0.5pc/[dr]|-{\xi^p_{\unit}}
        &\pi^{\pr} \omega^{p*}(\unit)
        \ar[d]^-{\gamma^{\pr}_{\omega^{p*}(\unit)}}
        \\
        \pi^{\pr}\kappa
        \ar[d]_-{\gamma^{\pr}_{\kappa}}
        \ar@/^0.5pc/[dr]|-{e_{\pi^{\pr}}I_{\kappa}}
        &\text{ }
        &\text{ }
        &\omega^{p*}(\unit)
        \ar[d]^-{\omega^{p*}_{\unit}}_-{\cong}
        \\
        \kappa
        \ar@{=}[r]
        &\kappa
        \ar@{=}[r]
        &\kappa
        \ar@{=}[r]
        &\kappa
      }
    }}
  \end{equation*}
  \begin{equation*}
    \vcenter{\hbox{
      \xymatrix@C=50pt{
        \pi^{\pr}pz
        \ar@{->>}[d]_-{I_{\pi^{\pr}}\xi^p_{\kos{t}_{\pi}^*(z)}}
        \ar@{=}[r]
        &\pi^{\pr}pz
        \ar@/^1pc/[ddl]|-{I_{\pi^{\pr}}e_pI_z}
        \ar[ddd]^-{e_{\pi^{\pr}p}I_z}
        \ar@{=}[r]
        &\pi^{\pr}pz
        \ar[d]^-{\gamma^{\pr}_pI_z}
        \ar@{=}[r]
        &\pi^{\pr}pz
        \ar@{->>}[d]^-{I_{\pi^{\pr}}\xi^p_{\kos{t}_{\pi}^*(z)}}
        \\
        \pi^{\pr}\omega^{p*}\kos{t}_{\pi}^*(z)
        \ar[d]_-{I_{\pi^{\pr}}\widehat{\omega}^{p*}_z}^-{\cong}
        &\text{ }
        &pz
        \ar[dd]^-{e_pI_z}
        \ar@{->>}@/^0.5pc/[dr]|-{\xi^p_{\kos{t}^*_{\pi}(z)}}
        &\pi^{\pr}\omega^{p*}\kos{t}_{\pi}^*(z)
        \ar[d]^-{\gamma^{\pr}_{\omega^{p*}\kos{t}_{\pi}^*(z)}}
        \\
        \pi^{\pr}z
        \ar[d]_-{\gamma^{\pr}_z}
        \ar@/^0.5pc/[dr]|-{e_{\pi^{\pr}}I_z}
        &\text{ }
        &\text{ }
        &\omega^{p*}\kos{t}^*_{\pi}(z)
        \ar[d]^-{\what{\omega}^{p*}_z}_-{\cong}
        \\
        z
        \ar@{=}[r]
        &z
        \ar@{=}[r]
        &z
        \ar@{=}[r]
        &z
      }
    }}
  \end{equation*}
  This shows that 
  $(\omega^{p*},\what{\omega}^{p*})$
  factors through
  $(\widebreve{\omega}\!^{p*}, \what{\widebreve{\omega}}^{p*})$
  as strong $\kos{K}$-tensor functors.
  We are left to show that 
  $(\widebreve{\omega}\!^{p*}, \what{\widebreve{\omega}}^{p*})$
  is an equivalence of strong $\kos{K}$-tensor categories.
  Consider the opposite $(\pi,\pi^{\pr})$-bitorsor
  $\mathring{p}$ over $\kappa$.
  \begin{equation*}
    \mathring{p}=(p,\gamma_p,\lambda^{\pr}_p)
    \qquad\qquad
    \begin{aligned}
      \gamma_p
      &:\!\!
      \xymatrix@C=20pt{
        \pi\otimes p
        \ar[r]^-{\varsigma_{\pi}\otimes I_p}_-{\cong}
        &\pi\otimes p
        \ar[r]^-{s_{\pi,p}}_-{\cong}
        &p\otimes \pi
        \ar[r]^-{\lambda_p}
        &p
      }
      \\
      \lambda^{\pr}_p
      &:\!\!
      \xymatrix@C=20pt{
        p\otimes \pi^{\pr}
        \ar[r]^-{I_p\otimes \varsigma_{\pi^{\pr}}}_-{\cong}
        &p\otimes \pi^{\pr}
        \ar[r]^-{s_{p,\pi^{\pr}}}_-{\cong}
        &\pi^{\pr}\otimes p
        \ar[r]^-{\gamma^{\pr}_p}
        &p
      }
    \end{aligned}
  \end{equation*}
  With resepect to the opposite bitorsor $\mathring{p}$,
  we can repeat what we did so far and obtain
  a strong $\kos{K}$-tensor functor
  $(\widebreve{\omega}\!^{\mathring{p}*},\what{\widebreve{\omega}}^{\mathring{p}*}):
  (\kos{R\!e\!p}(\pi^{\pr}),\kos{t}_{\pi^{\pr}}^*)
  \to (\kos{R\!e\!p}(\pi),\kos{t}_{\pi}^*)$.
  \begin{equation*}
    \vcenter{\hbox{
      \xymatrix@C=40pt{
        (\kos{R\!e\!p}(\pi^{\pr}),\kos{t}_{\pi^{\pr}}^*)
        \ar[r]^-{(\widebreve{\omega}\!^{\mathring{p}*},\what{\widebreve{\omega}}^{\mathring{p}*})}
        \ar@/_1pc/[dr]_-{(\omega^{\mathring{p}*},\what{\omega}^{\mathring{p}*})}
        &(\kos{R\!e\!p}(\pi),\kos{t}_{\pi}^*)
        \ar[d]^-{(\omega_{\pi},\what{\omega}_{\pi})}
        \\
        \text{ }
        &(\kos{K},\id_{\kos{K}})
      }
    }}
  \end{equation*}
  We are going to show that 
  $(\widebreve{\omega}\!^{p*}, \what{\widebreve{\omega}}^{p*})$
  and
  $(\widebreve{\omega}\!^{\mathring{p}*},\what{\widebreve{\omega}}^{\mathring{p}*})$
  are quasi-inverse to each other.
  We claim that there is a
  comonoidal $\kos{K}$-tensor natural isomorphism
  \begin{equation*}
    \vcenter{\hbox{
      \xymatrix@C=15pt{
        \widebreve{\vartheta}\!^p:
        (\id,\what{\id})
        \ar@2{->}[r]^-{\cong}
        &(\widebreve{\omega}\!^{\mathring{p}*},\what{\widebreve{\omega}}^{\mathring{p}*})(\widebreve{\omega}\!^{p*},\what{\widebreve{\omega}}^{p*})
        :(\kos{R\!e\!p}(\pi),\kos{t}_{\pi}^*)
        \to (\kos{R\!e\!p}(\pi),\kos{t}_{\pi}^*).
      }
    }}
  \end{equation*}
  We begin by describing the natural isomorphism
  $\xymatrix@C=15pt{
    \vartheta^p:
    \omega_{\pi}^*
    \ar@2{->}[r]^-{\cong}
    &\omega^{\mathring{p}*}\widebreve{\omega}\!^{p*}.
  }$
  Let $X=(x,\gamma_x)$ be an object in $\cat{Rep}(\pi)$.
  The object $\omega^{\mathring{p}*}\widebreve{\omega}\!^{p*}(X)$ in $\CK$
  is given by the following reflexive coequalizer.
  \begin{equation*}
    \vcenter{\hbox{
      \xymatrix@C=40pt{
        p\otimes \pi^{\pr}\otimes \omega^{p*}(X)
        \ar@<0.5ex>[r]^-{I_p\otimes \gamma^{\pr}_{\omega^{p*}(X)}}
        \ar@<-0.5ex>[r]_-{\lambda_p^{\pr}\otimes I_{\omega^{p*}(X)}}
        &p\otimes \omega^{p*}(X)
        \ar@{->>}[r]^-{\xi^{\mathring{p}}_{\widebreve{\omega}^{p*}(X)}}
        &\omega^{\mathring{p}*}\widebreve{\omega}\!^{p*}(X)
      }
    }}
  \end{equation*}
  The component of $\vartheta^p$ at $X=(x,\gamma_x)$
  is the unique isomorphism in $\CK$ which satisfies the relation below.
  \begin{equation*}
    \vcenter{\hbox{
      \xymatrix@C=70pt{
        ppx
        \ar@<0.5ex>[r]^-{e_pI_{px}}
        \ar@<-0.5ex>[r]_-{I_pe_pI_x}
        \ar[d]_-{I_p\kos{p}_!(\xi^p_X)}^-{\cong}
        &px
        \ar[r]^-{e_pI_x}
        \ar[d]^-{\kos{p}_!(\xi^p_X)}_-{\cong}
        &x
        \ar@{.>}[dd]^-{\vartheta^p_X}_-{\exists!\text{ }\cong}
        \\
        pp\omega^{p*}(X)
        \ar[d]_-{I_p\kos{p}_!(\xi^{\mathring{p}}_{\widebreve{\omega}^{p*}(X)})}^-{\cong}
        \ar@{}[r]|-{(\dagger)}
        &p\omega^{p*}(X)
        \ar[d]^-{\kos{p}_!(\xi^{\mathring{p}}_{\widebreve{\omega}^{p*}(X)})}_-{\cong}
        \ar@{->>}@/^0.5pc/[dr]|-{\xi^{\mathring{p}}_{\widebreve{\omega}^{p*}(X)}}
        &\text{ }
        \\
        pp\omega^{\mathring{p}*}\omega^{p*}(X)
        \ar@<0.5ex>[r]^-{e_pI_{p\omega^{\mathring{p}*}\omega^{p*}(X)}}
        \ar@<-0.5ex>[r]_-{I_pe_pI_{\omega^{\mathring{p}*}\omega^{p*}(X)}}
        &p\omega^{\mathring{p}*}\omega^{p*}(X)
        \ar[r]^-{e_pI_{\omega^{\mathring{p}*}\omega^{p*}(X)}}
        &\omega^{\mathring{p}*}\omega^{p*}(X)
      }
    }}
  \end{equation*}
  To conclude that the isomorphism $\vartheta^p_X$ is well-defined,
  we need to verify that the square diagram $(\dagger)$ commutes.
  The diagram $(\dagger)$ with upper horizontal morphisms is obviously commutative,
  and we verify the the diagram $(\dagger)$ with lower horizontal morphisms as follows.
  \begin{equation*}
    (\dagger):
    \!\!\!\!
    \!\!\!\!
    \!\!\!\!
    \vcenter{\hbox{
      \xymatrix@C=5pt{
        ppx
        \ar[dd]_-{I_pe_pI_x}
        \ar@{=}[r]
        &ppx
        \ar[d]^-{\cp_pI_{px}}
        \ar@{=}[rrr]
        &\text{ }
        &\text{ }
        &ppx
        \ar@{=}[dddd]
        \\
        \text{ }
        &pppx
        \ar[d]^-{I_ps_{p,p}I_{px}}_-{\cong}
        &\text{ }
        &\text{ }
        &\text{ }
        \\
        px
        \ar[ddd]_-{\kos{p}_!(\xi^p_X)}^-{\cong}
        \ar@/_1pc/[dddr]|-{\cp_p^{(3)}I_x}
        &pppx
        \ar[d]^-{I_{pp}\cp_pI_x}
        \ar@/^1.5pc/[ddrrr]|-{I_{pp}e_pI_x}
        \ar@{}[dddrrr]|-{(\dagger1)}
        &\text{ }
        &\text{ }
        &\text{ }
        \\
        \text{ }
        &ppppx
        \ar[dd]|-{I_pe_pI_{ppx}}
        \ar@/_0.5pc/[dr]^-{I_p(\tau^{\pr}_p)^{-1}I_{px}}_-{\cong}
        &\text{ }
        &\text{ }
        &\text{ }
        \\
        \text{ }
        &\text{ }
        &pp\pi^{\pr}px
        \ar@/_0.5pc/[dr]|-{I_{pp}\gamma^{\pr}_pI_x}
        &\text{ }
        &ppx
        \ar@/_0.5pc/[dl]|-{I_p\cp_pI_x}
        \ar[dd]^-{I_p\kos{p}_!(\xi^p_X)}_-{\cong}
        \\
        p\omega^{p*}(X)
        \ar[ddd]_-{\kos{p}_!(\xi^{\mathring{p}}_{\widebreve{\omega}^{p*}(X)})}^-{\cong}
        &pppx
        \ar@{->>}[d]^-{I_{pp}\xi^p_X}
        &\text{ }
        &pppx
        \ar@{->>}[d]^-{I_{pp}\xi^p_X}
        &\text{ }
        \\
        \text{ }
        &pp\omega^{p*}(X)
        \ar@{->>}[dd]^-{I_{p}\xi^{\mathring{p}}_{\widebreve{\omega}^{p*}(X)}}
        \ar@{}[rr]|-{(\dagger2)}
        &\text{ }
        &pp\omega^{p*}(X)
        \ar@{->>}[dd]^-{I_{p}\xi^{\mathring{p}}_{\widebreve{\omega}^{p*}(X)}}
        \ar@{=}[r]
        &pp\omega^{p*}(X)
        \ar[d]^-{I_p\kos{p}_!(\xi^{\mathring{p}}_{\widebreve{\omega}^{p*}(X)})}_-{\cong}
        \\
        \text{ }
        &\text{ }
        &\text{ }
        &\text{ }
        &pp\omega^{\mathring{p}*}\widebreve{\omega}^{p*}(X)
        \ar[d]^-{I_pe_pI_{\omega^{\mathring{p}*}\widebreve{\omega}^{p*}(X)}}
        \\
        p\omega^{\mathring{p}*}\widebreve{\omega}^{p*}(X)
        \ar@{=}[r]
        &p\omega^{\mathring{p}*}\widebreve{\omega}^{p*}(X)
        \ar@{=}[rr]
        &\text{ }
        &p\omega^{\mathring{p}*}\widebreve{\omega}^{p*}(X)
        \ar@{=}[r]
        &p\omega^{\mathring{p}*}\widebreve{\omega}^{p*}(X)
      }
    }}
  \end{equation*}
  where
  $\cp_p^{(3)}:p\to p\otimes p\otimes p$
  is the iterated coproduct.
  The diagram $(\dagger1)$
  is the relation (\ref{eq3 Torsors antipode formula})
  and the diagram $(\dagger2)$ is verified below.
  \begin{equation*}
    (\dagger2):
    \vcenter{\hbox{
      \xymatrix@C=60pt{
        pppx
        \ar[dd]_-{e_pI_{ppx}}
        \ar@{=}[r]
        &pppx
        \ar@{->>}[d]^-{I_{pp}\xi^p_X}
        \ar@{=}[r]
        &pppx
        \ar[d]^-{(\tau^{\pr}_p)^{-1}I_{px}}_-{\cong}
        \\
        \text{ }
        &pp\omega^{p*}(X)
        \ar@/_1.2pc/[ddl]|-{e_pI_{p\omega^{p*}(X)}}
        \ar[d]^-{(\tau^{\pr}_p)^{-1}I_{\omega^{p*}(X)}}_-{\cong}
        &p\pi^{\pr} px
        \ar[d]^-{I_p\gamma^{\pr}_pI_x}
        \ar@{->>}@/^0.5pc/[dl]|-{I_{p\pi^{\pr}}\xi^p_X}
        \\
        ppx
        \ar@{->>}[d]_-{I_p\xi^p_X}
        &p\pi^{\pr}\omega^{p*}(X)
        \ar@/_0.5pc/[dr]|-{I_p\gamma^{\pr}_{\omega^{p*}(X)}}
        \ar@/^0.5pc/[dl]|-{\lambda^{\pr}_pI_{\omega^{p*}(X)}}
        &ppx
        \ar@{->>}[d]^-{I_p\xi^p_X}
        \\
        p \omega^{p*}(X)
        \ar@{->>}[d]_-{\xi^{\mathring{p}}_{\widebreve{\omega}^{p*}(X)}}
        &\text{ }
        &p \omega^{p*}(X)
        \ar@{->>}[d]^-{\xi^{\mathring{p}}_{\widebreve{\omega}^{p*}(X)}}
        \\
        \omega^{\mathring{p}*}\omega^{p*}(X)
        \ar@{=}[rr]
        &\text{ }
        &\omega^{\mathring{p}*}\omega^{p*}(X)
      }
    }}
  \end{equation*}
  This shows that the natural isomorphism
  $\vartheta^p$ in $\CK$ is well-defined.
  Next we show that $\vartheta^p$
  is a comonoidal $\kos{K}$-tensor natural isomorphism
  \begin{equation*}
    \vcenter{\hbox{
      \xymatrix@C=15pt{
        \vartheta^p:
        (\omega_{\pi}^*,\what{\omega}_{\pi}^*)
        \ar@2{->}[r]^-{\cong}
        &(\omega^{\mathring{p}*},\what{\omega}^{\mathring{p}*})
        (\widebreve{\omega}\!^{p*},\what{\widebreve{\omega}}^{p*})
        :(\kos{R\!e\!p}(\pi),\kos{t}_{\pi}^*)\to (\kos{K},\id_{\kos{K}})
        .
      }
    }}
  \end{equation*}
  Let $X=(x,\gamma_x)$ and $Y=(y,\gamma_y)$
  be objects in $\cat{Rep}(\pi)$.
  We have
  \begin{equation*}
    \vcenter{\hbox{
      \xymatrix@C=20pt{
        x\otimes y
        \ar[d]_-{\vartheta^{p}_{X\tensor_{\pi} Y}}^-{\cong}
        \ar@{=}[rr]
        &\text{ }
        &x\otimes y
        \ar[d]^-{\vartheta^p_X\otimes \vartheta^p_Y}_-{\cong}
        \\
        \omega^{\mathring{p}*}\widebreve{\omega}\!^{p*}(X\tensor\!_{\pi} Y)
        \ar[r]_-{\omega^{\mathring{p}*}(\widebreve{\omega}^{p*}_{X,Y})}^-{\cong}
        &\omega^{\mathring{p}*}(\widebreve{\omega}\!^{p*}(X)\tensor\!_{\pi^{\pr}} \widebreve{\omega}\!^{p*}(Y))
        \ar[r]_(0.6){\omega^{\mathring{p}*}_{\widebreve{\omega}^{p*}(X),\widebreve{\omega}^{p*}(Y)}}^-{\cong}
        &\omega^{\mathring{p}*}\widebreve{\omega}\!^{p*}(X)
        \otimes \omega^{\mathring{p}*}\widebreve{\omega}\!^{p*}(Y)
      }
    }}
  \end{equation*}
  which we obtain by right-cancelling the epimorphism
  $e_p\otimes I_{x\otimes y}$ in the diagram below.
  \begin{equation*}
    \vcenter{\hbox{
      \xymatrix@C=7pt{
        pxy
        \ar[d]_-{e_pI_{xy}}
        \ar@{=}[r]
        &pxy
        \ar[d]^-{\cp_pI_{xy}}
        \ar@{=}[rr]
        &\text{ }
        &pxy
        \ar[dd]|-{(p\otimes)_{x,y}}
        \ar@/^2pc/@<2ex>[ddd]|-{e_pI_{xy}}
        \\
        xy
        \ar[dd]_-{\vartheta^p_{XY}}^-{\cong}
        &ppxy
        \ar@{->>}[d]_-{I_p\xi^p_{XY}}
        \ar@/^1pc/[dr]|-{I_p(p\otimes)_{x,y}}
        &\text{ }
        &\text{ }
        \\
        \text{ }
        &p\omega^{p*}(XY)
        \ar@/_0.5pc/@{->>}[dl]|-{\xi^{\mathring{p}}_{\widebreve{\omega}^{p*}(XY)}}
        \ar[d]_-{I_p\omega^{p*}_{X,Y}}^-{\cong}
        &ppxpy
        \ar@{->>}@/^0.5pc/[dl]|-{I_p\xi^p_X\xi^p_Y}
        \ar@<1ex>[d]|-{(p\otimes)_{px,py}}
        &pxpy
        \ar@/^0.5pc/[dl]|-{\cp_pI_x\cp_pI_y}
        \ar[d]|-{e_pI_xe_pI_y}
        \\
        \omega^{\mathring{p}*}\widebreve{\omega}\!^{p*}(XY)
        \ar[d]_-{\omega^{\mathring{p}*}(\widebreve{\omega}^{p*}_{X,Y})}^-{\cong}
        &p\omega^{p*}(X)\omega^{p*}(Y)
        \ar@/_0.5pc/@{->>}[dl]|-{\xi^{\mathring{p}}_{\widebreve{\omega}^{p*}(X)\widebreve{\omega}^{p*}(Y)}}
        \ar[d]|-{(p\otimes)_{\omega^{p*}(X),\omega^{p*}(Y)}}
        &ppxppy
        \ar@/^1pc/[dl]|-{I_p\xi^p_XI_p\xi^p_Y}
        &xy
        \ar[dd]^-{\vartheta^p_X\vartheta^p_Y}_-{\cong}
        \\
        \omega^{\mathring{p}*}(\widebreve{\omega}\!^{p*}(X)\widebreve{\omega}\!^{p*}(Y))
        \ar[d]|-{\omega^{\mathring{p}*}_{\widebreve{\omega}^{p*}(X),\widebreve{\omega}^{p*}(Y)}}
        &p\omega^{p*}(X)p\omega^{p*}(Y)
        \ar@/^0.5pc/@{->>}[dl]|-{\xi^{\mathring{p}}_{\widebreve{\omega}^{p*}(X)}\xi^{\mathring{p}}_{\widebreve{\omega}^{p*}(Y)}}
        &\text{ }
        &\text{ }
        \\
        \omega^{\mathring{p}*}\!\widebreve{\omega}\!^{p*}(X)\omega^{\mathring{p}*}\!\widebreve{\omega}\!^{p*}(Y)
        \ar@{=}[rrr]
        &\text{ }
        &\text{ }
        &\omega^{\mathring{p}*}\!\widebreve{\omega}\!^{p*}(X)\omega^{\mathring{p}*}\!\widebreve{\omega}\!^{p*}(Y)
      }
    }}
  \end{equation*}
  We also have
  \begin{equation*}
    \vcenter{\hbox{
      \xymatrix@C=40pt{
        \kappa
        \ar[d]_-{\vartheta^p_{\unit\!_{\pi}}}^-{\cong}
        \ar@{=}[rr]
        &\text{ }
        &\kappa
        \ar@{=}[d]
        \\
        \omega^{\mathring{p}*}\widebreve{\omega}\!^{p*}(\unit\!_{\pi})
        \ar[r]^-{\omega^{\mathring{p}*}(\omega^{p*}_{\unit\!_{\pi}})}_-{\cong}
        &\omega^{\mathring{p}*}(\unit\!_{\pi^{\pr}})
        \ar[r]^-{\omega^{\mathring{p}*}_{\unit\!_{\pi^{\pr}}}}_-{\cong}
        &\kappa
      }
    }}
  \end{equation*}
  which we obtain by right-cancelling the epimorphism 
  $e_p\otimes I_{\kappa}$ in the diagram below.
  \begin{equation*}
    \vcenter{\hbox{
      \xymatrix@C=50pt{
        p\kappa
        \ar[d]_-{e_pI_{\kappa}}
        \ar@{=}[r]
        &p\kappa
        \ar[d]^-{\cp_pI_{\kappa}}
        \ar@{=}[r]
        &p\kappa
        \ar@{=}[ddd]
        \\
        \kappa
        \ar[dd]_-{\vartheta^p_{\unit\!_{\pi}}}^-{\cong}
        &pp\kappa
        \ar@{->>}[d]^-{I_p\xi^p_{\unit\!_{\pi}}}
        \ar@/^1pc/[ddr]|-{I_pe_pI_{\kappa}}
        &\text{ }
        \\
        \text{ }
        &p\omega^{p*}(\unit\!_{\pi})
        \ar@/_0.5pc/@{->>}[dl]|-{\xi^{\mathring{p}}_{\widebreve{\omega}^{p*}(\unit\!_{\pi})}}
        \ar[d]^-{I_p\omega^{p*}_{\unit\!_{\pi}}}_-{\cong}
        &\text{ }
        \\
        \omega^{\mathring{p}*}\widebreve{\omega}\!^{p*}(\unit\!_{\pi})
        \ar[d]_-{\omega^{\mathring{p}*}(\widebreve{\omega}^{p*}_{\unit\!_{\pi}})}^-{\cong}
        &p\kappa
        \ar@/_0.5pc/@{->>}[dl]|-{\xi^{\mathring{p}}_{\unit\!_{\pi^{\pr}}}}
        \ar@{=}[r]
        &p\kappa
        \ar[dd]^-{e_pI_{\kappa}}
        \\
        \omega^{\mathring{p}*}(\unit\!_{\pi^{\pr}})
        \ar[d]_-{\omega^{\mathring{p}*}_{\unit\!_{\pi^{\pr}}}}^-{\cong}
        &\text{ }
        &\text{ }
        \\
        \kappa
        \ar@{=}[rr]
        &\text{ }
        &\kappa
      }
    }}
  \end{equation*}
  This shows that
  $\vartheta^p$
  is a comonoidal natural isomorphism.
  Let $z$ be an object in $\CK$.
  We have
  \begin{equation*}
    \what{\omega}_{\pi}^*:
    \xymatrix@C=25pt{
      \omega_{\pi}^*\kos{t}_{\pi}^*
      \ar@2{->}[r]^-{\vartheta^p\kos{t}_{\pi}^*}_-{\cong}
      &\omega^{\mathring{p}*}\widebreve{\omega}\!^{p*}\kos{t}_{\pi}^*
      \ar@2{->}[r]^-{\what{\omega^{\mathring{p}*}\widebreve{\omega}\!^{p*}}}_-{\cong}
      &\id_{\kos{K}}
    }
    \vcenter{\hbox{
      \xymatrix@C=25pt{
        z
        \ar[rr]^-{I_z=(\what{\omega}_{\pi}^*)_z}_-{\cong}
        \ar[d]_-{\vartheta^p_{\kos{t}_{\pi}^*(z)}}^-{\cong}
        &\text{ }
        &z
        \ar@{=}[d]
        \\
        \omega^{\mathring{p}*}\widebreve{\omega}\!^{p*}\kos{t}^*_{\pi}(z)
        \ar[r]^-{\omega^{\mathring{p}*}(\what{\widebreve{\omega}}^{p*}_z)}_-{\cong}
        &\omega^{\mathring{p}*}\kos{t}^*_{\pi^{\pr}}(z)
        \ar[r]^-{\what{\omega}^{\mathring{p}*}_z}_-{\cong}
        &z
      }
    }}
  \end{equation*}
  by right-cancelling the epimorphism
  $e_p\otimes I_z$ in the diagram below.
  \begin{equation*}
    \vcenter{\hbox{
      \xymatrix@C=50pt{
        pz
        \ar[d]_-{e_pI_z}
        \ar@{=}[r]
        &pz
        \ar[d]^-{\cp_pI_z}
        \ar@{=}[r]
        &pz
        \ar@{=}[ddd]
        \\
        z
        \ar[dd]_-{\vartheta^p_{\kos{t}_{\pi}^*(z)}}^-{\cong}
        &ppz
        \ar@{->>}[d]^-{I_p\xi^p_{\kos{t}_{\pi}^*(z)}}
        \ar@/^1pc/[ddr]|-{I_pe_pI_z}
        &\text{ }
        \\
        \text{ }
        &p\omega^{p*}\kos{t}^*_{\pi}(z)
        \ar@{->>}@/_0.5pc/[dl]|-{\xi^{\mathring{p}}_{\widebreve{\omega}\!^{p*}\kos{t}_{\pi}^*(z)}}
        \ar[d]^-{I_p\what{\omega}^{p*}_z}_-{\cong}
        &\text{ }
        \\
        \omega^{\mathring{p}*}\widebreve{\omega}\!^{p*}\kos{t}_{\pi}^*(z)
        \ar[d]_-{\omega^{\mathring{p}*}(\what{\widebreve{\omega}}^{p*}_z)}^-{\cong}
        &pz
        \ar@/_0.5pc/@{->>}[dl]|-{\xi^{\mathring{p}}_{\kos{t}_{\pi^{\pr}}^*(z)}}
        \ar@{=}[r]
        &pz
        \ar[dd]^-{e_pI_z}
        \\
        \omega^{\mathring{p}*}\kos{t}_{\pi^{\pr}}^*(z)
        \ar[d]_-{\what{\omega}^{\mathring{p}*}_z}^-{\cong}
        &\text{ }
        &\text{ }
        \\
        z
        \ar@{=}[rr]
        &\text{ }
        &z
      }
    }}
  \end{equation*}
  This shows that
  $\vartheta^p$ is a comonoidal $\kos{K}$-tensor natural isomorphism.
  Finally, we show that $\vartheta^p$ lifts to 
  a comonoidal $\kos{K}$-tensor natural isomorphism
  \begin{equation*}
    \vcenter{\hbox{
      \xymatrix@C=15pt{
        \widebreve{\vartheta}\!^p:
        (\id,\what{\id})
        \ar@2{->}[r]^-{\cong}
        &(\widebreve{\omega}\!^{\mathring{p}*},\what{\widebreve{\omega}}^{\mathring{p}*})
        (\widebreve{\omega}\!^{p*},\what{\widebreve{\omega}}^{p*})
        :(\kos{R\!e\!p}(\pi),\kos{t}_{\pi}^*)
        \to (\kos{R\!e\!p}(\pi),\kos{t}_{\pi}^*).
      }
    }}
  \end{equation*}
  Let $X=(x,\gamma_x)$ be a representation of $\pi$.
  We need to show that the isomorphism
  $\vartheta^p_X
  :x\xrightarrow{\cong}\omega^{\mathring{p}*}\widebreve{\omega}\!^{p*}(X)$
  in $\CK$ is compatible with left $\pi$-actions.
  Note that the left $\pi$-action
  $\gamma_{\omega^{\mathring{p}*}\widebreve{\omega}^{p*}(X)}$ is
  the unique morphism satisfying the relation below.
  \begin{equation*}
    \vcenter{\hbox{
      \xymatrix@C=40pt{
        \pi\otimes p\otimes \omega^{p*}(X)
        \ar@{->>}[r]^-{I_{\pi}\otimes \xi^{\mathring{p}}_{\widebreve{\omega}^{p*}(X)}}
        \ar[d]_-{\gamma_p\otimes I_{\omega^{p*}(X)}}
        &\pi\otimes \omega^{\mathring{p}*}\!\widebreve{\omega}^{p*}(X)
        \ar@{.>}[d]^-{\gamma_{\omega^{\mathring{p}*}\!\widebreve{\omega}^{p*}(X)}}_-{\exists!}
        \\
        p\otimes \omega^{p*}(X)
        \ar@{->>}[r]^-{\xi^{\mathring{p}}_{\widebreve{\omega}^{p*}(X)}}
        &\omega^{\mathring{p}*}\!\widebreve{\omega}^{p*}(X)
        .
      }
    }}
  \end{equation*}
  We begin our calculation as
  \begin{equation*}
    \vcenter{\hbox{
      \xymatrix@C=40pt{
        p\pi x
        \ar[ddd]_-{e_pI_{\pi x}}
        \ar@{=}[rr]
        &\text{ }
        &p\pi x
        \ar[d]^-{s_{p,\pi}I_x}_-{\cong}
        \ar@{=}[r]
        &p\pi x
        \ar[ddddddd]^-{\tau_pI_x}_-{\cong}
        \\
        \text{ }
        &\text{ }
        &\pi p x
        \ar[d]^-{I_{\pi}\cp_pI_x}
        \ar@/_1pc/[ddll]|-{I_{\pi}e_pI_x}
        &\text{ }
        \\
        \text{ }
        &\text{ }
        &\pi p p x
        \ar@/_0.5pc/@{->>}[dl]|-{I_{\pi p}\xi^p_X}
        \ar[d]^-{\varsigma_{\pi}I_{ppx}}_-{\cong}
        &\text{ }
        \\
        \pi x
        \ar[d]_-{I_{\pi}\vartheta^p_X}^-{\cong}
        &\pi p \omega^{p*}(X)
        \ar@/_0.5pc/[dl]|-{I_{\pi}\xi^{\mathring{p}}_{\widebreve{\omega}^{p*}(X)}}
        \ar@/_3.5pc/@<-3ex>[ddddd]|-{\gamma_pI_{\omega^{p*}(X)}}
        \ar[d]^-{\varsigma_{\pi}I_{p\omega^{p*}(X)}}_-{\cong}
        &\pi p px
        \ar@{->>}@/^0.5pc/[dl]|-{I_{\pi p}\xi^p_X}
        \ar[d]^-{s_{\pi,p}I_{px}}_-{\cong}
        &\text{ }
        \\
        \pi \omega^{\mathring{p}*}\widebreve{\omega}^{p*}(X)
        \ar[ddddd]|-{\gamma_{\omega^{\mathring{p}*}\widebreve{\omega}^{p*}(X)}}
        &\pi p \omega^{p*}(X)
        \ar[d]^-{s_{\pi,p}I_{\omega^{p*}(X)}}_-{\cong}
        &p \pi px
        \ar@{->>}@/^0.5pc/[dl]|-{I_{p\pi}\xi^p_X}
        \ar[d]^-{\tau_pI_{px}}_-{\cong}
        \ar@{}[ur]|-{(\ddagger)}
        &\text{ }
        \\
        \text{ }
        &p \pi \omega^{p*}(X)
        \ar@/_2pc/@<-1.5ex>[ddd]|-{\lambda_pI_{\omega^{p*}(X)}}
        \ar[d]^-{\tau_pI_{\omega^{p*}(X)}}_-{\cong}
        &p p px
        \ar@/^0.5pc/@{->>}[dl]|-{I_{pp}\xi^p_X}
        \ar[d]^-{(\tau^{\pr}_p)^{-1}I_{px}}_-{\cong}
        &\text{ }
        \\
        \text{ }
        &pp\omega^{p*}(X)
        \ar[dd]|-{e_pI_{p\omega^{p*}(X)}}
        \ar@/_0.5pc/[dr]^-{(\tau^{\pr}_p)^{-1}I_{\omega^{p*}(X)}}_-{\cong}
        &p\pi^{\pr}px
        \ar@{->>}[d]^-{I_{p\pi^{\pr}}\xi^p_X}
        \ar@/^0.5pc/[dr]|-{I_p\gamma^{\pr}_pI_x}
        &\text{ }
        \\
        \text{ }
        &\text{ }
        &p\pi^{\pr}\omega^{p*}(X)
        \ar@/^0.5pc/[dl]|-{\lambda^{\pr}_pI_{\omega^{p*}(X)}}
        \ar@/_0.5pc/[dr]|-{I_p\gamma^{\pr}_{\omega^{p*}(X)}}
        &ppx
        \ar@{->>}[d]^-{I_p\xi^p_X}
        \\
        \text{ }
        &p \omega^{p*}(X)
        \ar@{->>}[d]^-{\xi^{\mathring{p}}_{\widebreve{\omega}^{p*}(X)}}
        &\text{ }
        &p \omega^{p*}(X)
        \ar@{->>}[d]^-{\xi^{\mathring{p}}_{\widebreve{\omega}^{p*}(X)}}
        \\
        \omega^{\mathring{p}*}\widebreve{\omega}^{p*}(X)
        \ar@{=}[r]
        &\omega^{\mathring{p}*}\widebreve{\omega}^{p*}(X)
        \ar@{=}[rr]
        &\text{ }
        &\omega^{\mathring{p}*}\widebreve{\omega}^{p*}(X)
      }
    }}
  \end{equation*}
  where we postpone the proof of the diagram $(\ddagger)$.
  We finish our calculation as follows.
  \begin{equation*}
    \vcenter{\hbox{
      \xymatrix@C=40pt{
        p\pi x
        \ar[dd]_-{\tau_pI_x}^-{\cong}
        \ar@{=}[r]
        &p \pi x
        \ar[d]^-{\cp_pI_{\pi x}}
        \ar@{=}[r]
        &p \pi x
        \ar[d]^-{I_p\gamma_x}
        \ar@{=}[r]
        &p \pi x
        \ar[d]^-{e_pI_{\pi x}}
        \\
        \text{ }
        &pp\pi x
        \ar@/^0.5pc/[dl]|-{I_p\lambda_pI_x}
        \ar[d]^-{I_{pp}\gamma_x}
        &px
        \ar@/^0.5pc/[dl]|-{\cp_pI_x}
        \ar@/_0.5pc/[dr]|-{e_pI_x}
        &\pi x
        \ar[d]^-{\gamma_x}
        \\
        ppx
        \ar@{->>}[d]_-{I_p\xi^p_X}
        &ppx
        \ar@/^0.5pc/@{->>}[dl]|-{I_p\xi^p_X}
        &\text{ }
        &x
        \ar[dd]^-{\vartheta^p_X}_-{\cong}
        \\
        p\omega^{p*}(X)
        \ar@{->>}[d]_-{\xi^{\mathring{p}}_{\widebreve{\omega}^{p*}(X)}}
        &\text{ }
        &\text{ }
        &\text{ }
        \\
        \omega^{\mathring{p}*}\widebreve{\omega}^{p*}(X)
        \ar@{=}[rrr]
        &\text{ }
        &\text{ }
        &\omega^{\mathring{p}*}\widebreve{\omega}^{p*}(X)
      }
    }}
  \end{equation*}
  By right-cancelling the epimorphism
  $e_p\otimes I_{\pi\otimes x}$
  in the above diagrams,
  we obtain that $\vartheta^p$
  is compatible with left $\pi$-actions.
  We are left to verify the diagram
  $(\ddagger)$.
  We have
  \begin{equation*}
    (\ddagger1):
    \vcenter{\hbox{
      \xymatrix@C=50pt{
        p\pi
        \ar[d]_-{s_{p,\pi}}^-{\cong}
        \ar@{=}[r]
        &p\pi
        \ar[d]^-{\cp_pI_{\pi}}
        \ar@{=}[rr]
        &\text{ }
        &p\pi
        \ar[dd]^-{\cp_p\cp_{\pi}}
        \\
        \pi p
        \ar[d]_-{I_{\pi}\cp_p}
        &pp\pi
        \ar[d]^-{I_p\cp_pI_{\pi}}
        \ar@/^0.5pc/[dr]|-{I_p\cp_{p\pi}}
        &\text{ }
        &\text{ }
        \\
        \pi p p
        \ar[d]_-{\varsigma_{\pi}I_{pp}}^-{\cong}
        &ppp\pi
        \ar[d]^-{I_{pp}s_{p,\pi}}_-{\cong}
        &pp\pi p\pi
        \ar@/^0.5pc/[dl]|-{I_{pp\pi p}e_{\pi}}
        \ar[d]^-{I_{pp\pi p}\cp_{\pi}}
        &pp \pi\pi 
        \ar@/^0.5pc/[dl]|-{I_p\cp_{p\pi}I_{\pi}}
        \ar[d]^-{I_{pp}\varsigma_{\pi}I_{\pi}}_-{\cong}
        \\
        \pi p p
        \ar[d]_-{s_{p,\pi}I_p}^-{\cong}
        &pp\pi p
        \ar[ddd]^-{I_{pp}\varsigma_{\pi}I_p}_-{\cong}
        \ar@/_0.5pc/[ddr]|-{I_{pp}\varsigma_{\pi}I_p u_{\pi}}
        &pp\pi p\pi\pi
        \ar[d]^-{I_{pp}\varsigma_{\pi}I_p\varsigma_{\pi}I_{\pi}}_-{\cong}
        &pp\pi\pi
        \ar[d]^-{I_p\lambda_pI_{\pi}}
        \\
        p\pi p
        \ar[ddd]_-{\tau_pI_p}^-{\cong}
        \ar@/^0.5pc/[ddr]|-{\cp_p I_{\pi p}}
        &\text{ }
        &pp\pi p\pi\pi
        \ar[d]^-{I_{pp\pi p}\pc_{\pi}}
        \ar@/^0.5pc/[dr]|-{I_p\lambda_p\lambda_pI_{\pi}}
        &pp\pi
        \ar[d]^-{I_p\cp_pI_{\pi}}
        \\
        \text{ }
        &\text{ }
        &pp\pi p\pi
        \ar@/^0.5pc/[dl]|-{I_{pp\pi}\lambda_p}
        \ar[dd]^-{I_p\lambda_p\lambda_p}
        &ppp\pi
        \ar[dd]^-{I_{pp}\lambda_p}
        \\
        \text{ }
        &pp\pi p
        \ar[d]^-{I_p\lambda_pI_p}
        &\text{ }
        &\text{ }
        \\
        ppp
        \ar@{=}[r]
        &ppp
        \ar@{=}[r]
        &ppp
        \ar@{=}[r]
        &ppp
      }
    }}
  \end{equation*}
  and using the above diagram $(\ddagger1)$,
  we verify the diagram $(\ddagger)$ as follows.
  \begin{equation*}
    (\ddagger):
    \vcenter{\hbox{
      \xymatrix@C=60pt{
        p\pi
        \ar[d]_-{s_{p,\pi}}^-{\cong}
        \ar@{=}[r]
        &p\pi
        \ar[d]^-{\cp_p\cp_{\pi}}
        \ar@{=}[rr]
        &\text{ }
        &p\pi
        \ar@{=}[dddd]
        \\
        \pi p
        \ar[d]_-{I_{\pi}\cp_p}
        &pp\pi\pi
        \ar[d]^-{I_{pp}\varsigma_{\pi}I_{\pi}}_-{\cong}
        \ar@/^2.5pc/[dddrr]|-{I_pe_{p\pi}I_{\pi}}
        &\text{ }
        &\text{ }
        \\
        \pi p p
        \ar[d]_-{\varsigma_{\pi}I_{pp}}^-{\cong}
        \ar@{}[dr]|-{(\ddagger1)}
        &pp\pi\pi
        \ar[d]^-{I_p\lambda_pI_{\pi}}
        \ar@/^1.5pc/[ddrr]|-{I_pe_{p\pi}I_{\pi}}
        &\text{ }
        &\text{ }
        \\
        \pi p p
        \ar[d]_-{s_{\pi,p}I_p}^-{\cong}
        &pp\pi
        \ar[d]^-{I_p\cp_pI_{\pi}}
        \ar@/^0.5pc/[drr]|-{I_pe_pI_{\pi}}
        &\text{ }
        &\text{ }
        \\
        p\pi p
        \ar[d]_-{\tau_pI_p}^-{\cong}
        &ppp\pi
        \ar@/^0.5pc/[dl]|-{I_{pp}\lambda_p}
        \ar[d]^-{(\tau^{\pr}_p)^{-1}I_{p\pi}}
        \ar@{}[drr]|-{(\ddagger2)}
        &\text{ }
        &p\pi
        \ar@/_0.5pc/[ddl]|-{\cp_pI_{\pi}}
        \ar[ddd]^-{\tau_p}_-{\cong}
        \\
        ppp
        \ar[d]_-{(\tau^{\pr}_p)^{-1}I_p}^-{\cong}
        &p\pi^{\pr}p\pi
        \ar@/^0.5pc/[dl]|-{I_{p\pi^{\pr}}\lambda_p}
        \ar@/_0.5pc/[dr]|-{I_p\gamma^{\pr}_pI_{\pi}}
        &\text{ }
        &\text{ }
        \\
        p\pi^{\pr}p
        \ar[d]_-{I_p\gamma^{\pr}_p}
        &\text{ }
        &pp\pi
        \ar@/_0.5pc/[dr]|-{I_p\lambda_p}
        &\text{ }
        \\
        pp
        \ar@{=}[rrr]
        &\text{ }
        &\text{ }
        &pp
      }
    }}
  \end{equation*}
  We used the relation (\ref{eq3 Torsors antipode formula})
  in the diagram $(\ddagger2)$.
  This shows that we have
  a comonoidal $\kos{K}$-tensor natural isomorphism
  \begin{equation*}
    \vcenter{\hbox{
      \xymatrix@C=15pt{
        \widebreve{\vartheta}\!^p:
        (\id,\what{\id})
        \ar@2{->}[r]^-{\cong}
        &(\widebreve{\omega}\!^{\mathring{p}*},\what{\widebreve{\omega}}^{\mathring{p}*})
        (\widebreve{\omega}\!^{p*},\what{\widebreve{\omega}}^{p*})
        :(\kos{R\!e\!p}(\pi),\kos{t}_{\pi}^*)
        \to (\kos{R\!e\!p}(\pi),\kos{t}_{\pi}^*).
      }
    }}
  \end{equation*}
  If we repeat the construction of
  $\widebreve{\vartheta}\!^p$
  with respect to the opposite bitorsor $\mathring{p}$,
  then we obtain a comonoidal $\kos{K}$-tensor natural isomorphism
  \begin{equation*}
    \vcenter{\hbox{
      \xymatrix@C=15pt{
        \widebreve{\vartheta}\!^{\mathring{p}}:
        (\id,\what{\id})
        \ar@2{->}[r]^-{\cong}
        &(\widebreve{\omega}\!^{p*},\what{\widebreve{\omega}}^{p*})
        (\widebreve{\omega}\!^{\mathring{p}*},\what{\widebreve{\omega}}^{\mathring{p}*})
        :(\kos{R\!e\!p}(\pi^{\pr}),\kos{t}_{\pi^{\pr}}^*)
        \to (\kos{R\!e\!p}(\pi^{\pr}),\kos{t}_{\pi^{\pr}}^*).
      }
    }}
  \end{equation*}
  This shows that
  $(\widebreve{\omega}\!^{p*},\what{\widebreve{\omega}}^{p*})$
  and
  $(\widebreve{\omega}\!^{\mathring{p}*},\what{\widebreve{\omega}}^{\mathring{p}*})$
  are quasi-inverse to each other as strong $\kos{K}$-tensor functors.
  We conclude that
  $(\widebreve{\omega}\!^{p*},\what{\widebreve{\omega}}^{p*})$
  is an equivalence of strong $\kos{K}$-tensor categories.
  This completes the proof of Propoisition~\ref{prop twistbyatorsor}.
\qed\end{proof}

\begin{corollary}
  \label{cor twistbyatorsor}
  Let $\pi$ be a pre-Galois object in $\kos{K}$
  and let $p$ be a right $\pi$-torsor over $\kappa$.
  Recall Proposition~\ref{prop Torsors righttors become bitors}.
  We denote $\pi^p$ as the pre-Galois object in $\kos{K}$
  obtained by twisting $\pi$ by $p$,
  and $p$ becomes a $(\pi^p,\pi)$-bitorsor over $\kappa$.
  \begin{enumerate}
    \item 
    We have a pre-fiber functor
    $\varpi^p:\mathfrak{K}\to \mathfrak{Rep}(\pi)$
    which is obtained by twisting the pre-fiber functor
    $\varpi_{\pi}:\mathfrak{K}\to\mathfrak{Rep}(\pi)$.
    See Lemma~\ref{lem twistbyatorsor} and Proposition~\ref{prop twistbyatorsor}.

    \item
    The twisted pre-fiber functor
    $\varpi^p:\mathfrak{K}\to \mathfrak{Rep}(\pi)$
    factors through as an equivalence of Galois $\kos{K}$-prekosmoi
    $\widebreve{\varpi}\!^p:\mathfrak{Rep}(\pi^p)\xrightarrow[]{\simeq}\mathfrak{Rep}(\pi)$.
    \begin{equation*}
      \vcenter{\hbox{
        \xymatrix@C=40pt{
          \mathfrak{Rep}(\pi)
          &\mathfrak{Rep}(\pi^p)
          \ar[l]_-{\widebreve{\varpi}\!^p}^-{\simeq}
          \\
          \text{ }
          &\mathfrak{K}
          \ar@/^0.5pc/[ul]^-{\varpi^p}
          \ar[u]_-{\varpi_{\pi^{p}}}
        }
      }}
    \end{equation*}

    \item 
    The presheaf of groups of
    invertible Galois $\kos{K}$-transformations
    from $\varpi^p$ to $\varpi^p$
    is represented by the twisted pre-Galois object $\pi^p$ in $\kos{K}$.
    \begin{equation*}
      \xymatrix{
        \Hom_{\cat{Ens}(\kos{K})}(\slot,\pi^p)
        \ar@2{->}[r]^-{\cong}
        &\underline{\Aut}_{\mathbb{GAL}^{\cat{pre}}_{\kos{K}}}(\varpi^p)
        :\cat{Ens}(\kos{K})^{\op}\to\cat{Grp}
      }
    \end{equation*}
    \item 
    The presheaf of
    invertible Galois $\kos{K}$-transformations
    from $\varpi_{\pi}$ to $\varpi^p$
    is represented by $p$.
    \begin{equation*}
      \xymatrix{
        \Hom_{\cat{Ens}(\kos{K})}(\slot,p)
        \ar@2{->}[r]^-{\cong}
        &\underline{\Isom}_{\mathbb{GAL}^{\cat{pre}}_{\kos{K}}}(\varpi_{\pi},\varpi^p)
        :\cat{Ens}(\kos{K})^{\op}\to\cat{Set}
      }
    \end{equation*}
    The component of the universal element $\xi^p$
    of the representation 
    at each object $X=(x,\gamma_x)$ in $\cat{Rep}(\pi)$
    is the reflexive coequalizer
    $\xymatrix@C=15pt{
      \xi^p_X:p\otimes x\ar@{->>}[r]&\omega^{p*}(X)
    }$
    described in (\ref{eq twistbyatorsor}).
  \end{enumerate}
\end{corollary}

\begin{lemma} \label{lem twistbyatorsor functorial}
  Let $\pi$ be a pre-Galois object in $\kos{K}$
  and let $f:p\xrightarrow{\cong}q$
  be an isomorphism of right $\pi$-torsors over $\kappa$.
  We have an invertible Galois $\kos{K}$-transformation
  \begin{equation*}
    \vcenter{\hbox{
      \xymatrix{
        \mathfrak{Rep}(\pi)
        \\
        \mathfrak{K}
        \ar@/^1.2pc/[u]^-{\varpi^p}
        \ar@/_1.2pc/[u]_-{\varpi^q}
        \xtwocell[u]{}<>{<0>{\vartheta^f}}
      }
    }}
    \qquad\quad
    \vartheta^f:
    \xymatrix@C=15pt{
      \varpi^p
      \ar@2{->}[r]^-{\cong}
      &\varpi^q
      :\mathfrak{K}\to \mathfrak{Rep}(\pi)
    }
  \end{equation*}
  between twisted pre-fiber functors.
\end{lemma}
\begin{proof}
  Let $X=(x,\gamma_x)$ be an object in $\cat{Rep}(\pi)$.
  The component of
  $\vartheta^f:\omega^{p*}\cong \omega^{q*}$
  at each representation $X=(x,\gamma_x)$ of $\pi$ is
  the unique isomorphism which makes the following diagram commutative.
  \begin{equation*}
    \vcenter{\hbox{
      \xymatrix@C=40pt{
        p\otimes \pi\otimes x
        \ar[d]_-{f\otimes I_{\pi\otimes x}}^-{\cong}
        \ar@<0.5ex>[r]^-{I_p\otimes \gamma_x}
        \ar@<-0.5ex>[r]_-{\lambda_p\otimes I_x}
        &p\otimes x
        \ar@{->>}[r]^-{\xi^p_X}
        \ar[d]^-{f\otimes I_x}_-{\cong}
        &\omega^{p*}(X)
        \ar@{.>}[d]^-{\exists!\text{ }\vartheta^f_X}_-{\cong}
        \\
        q\otimes \pi\otimes x
        \ar@<0.5ex>[r]^-{I_q\otimes \gamma_x}
        \ar@<-0.5ex>[r]_-{\lambda_q\otimes I_x}
        &q\otimes x
        \ar@{->>}[r]^-{\xi^q_X}
        &\omega^{q*}(X)
      }
    }}
  \end{equation*}
  To conclude that $\vartheta^f:\varpi^p\cong \varpi^q$
  is an invertible Galois $\kos{K}$-transformation,
  we need to show that
  $\vartheta^f:(\omega^{p*},\what{\omega}^{p*})\cong (\omega^{q*},\what{\omega}^{q*})$
  is a comonoidal $\kos{K}$-tensor natural isomorphism.
  We obtain this by right-cancelling the epimorphisms 
  $\xi^p_{X\tensor\!_{\pi}Y}$,
  $\xi^p_{\unit\!_{\pi}}$
  and $\xi^p_{\kos{t}_{\pi}^*(z)}$
  in the diagrams below.
  Let $Y=(y,\gamma_y)$ be another object in $\cat{Rep}(\pi)$
  and let $z$ be an object in $\CK$.
  \begin{equation*}
    \vcenter{\hbox{
      \xymatrix@C=30pt{
        pxy
        \ar@{->>}[d]_-{\xi^p_{XY}}
        \ar@{=}[r]
        &pxy
        \ar[d]^-{fI_{xy}}_-{\cong}
        \ar@{=}[r]
        &pxy
        \ar[d]^-{(p\otimes)_{x,y}}
        \ar@{=}[r]
        &pxy
        \ar@{->>}[d]^-{\xi^p_{XY}}
        \\
        \omega^{p*}(XY)
        \ar[d]_-{\vartheta^p_{XY}}^-{\cong}
        &qxy
        \ar@{->>}@/^0.5pc/[dl]|-{\xi^q_{XY}}
        \ar@/_0.5pc/[dr]|-{(q\otimes)_{x,y}}
        &pxpy
        \ar[d]^-{f I_x f I_y}_-{\cong}
        \ar@{->>}@/^0.5pc/[dr]|-{\xi^p_X\xi^p_Y}
        &\omega^{p*}(XY)
        \ar[d]^-{\omega^{p*}_{X,Y}}_-{\cong}
        \\
        \omega^{q*}(XY)
        \ar[d]_-{\omega^{q*}_{X,Y}}^-{\cong}
        &\text{ }
        &qxqy
        \ar@{->>}@/_0.5pc/[dr]|-{\xi^q_X\xi^q_Y}
        &\omega^{p*}(X)\omega^{p*}(Y)
        \ar[d]^-{\vartheta^f_X\vartheta^f_Y}_-{\cong}
        \\
        \omega^{q*}(X)\omega^{q*}(Y)
        \ar@{=}[rrr]
        &\text{ }
        &\text{ }
        &\omega^{q*}(X)\omega^{q*}(Y)
      }
    }}
  \end{equation*}
  \begin{equation*}
    \vcenter{\hbox{
      \xymatrix@C=35pt{
        p\kappa
        \ar@{->>}[d]_-{\xi^p_{\unit}}
        \ar@{=}[r]
        &p\kappa
        \ar[d]_-{fI_{\kappa}}^-{\cong}
        \ar@{=}[r]
        \ar@/^0.5pc/[dddr]|-{e_pI_{\kappa}}
        &p\kappa
        \ar@{->>}[d]^-{\xi^p_{\unit}}
        \\
        \omega^{p*}(\unit)
        \ar[d]_-{\vartheta^f_{\unit}}^-{\cong}
        &q\kappa
        \ar@{->>}@/^0.5pc/[dl]|-{\xi^q_{\unit}}
        \ar[dd]^-{e_qI_{\kappa}}
        &\omega^{p*}(\unit)
        \ar[dd]^-{\omega^{p*}_{\unit}}_-{\cong}
        \\
        \omega^{q*}(\unit)
        \ar[d]_-{\omega^{q*}_{\unit}}^-{\cong}
        &\text{ }
        &\text{ }
        \\
        \kappa
        \ar@{=}[r]
        &\kappa
        \ar@{=}[r]
        &\kappa
      }
    }}
    \quad
    \vcenter{\hbox{
      \xymatrix@C=35pt{
        pz
        \ar@{->>}[d]_-{\xi^p_{\kos{t}^*_{\pi}(z)}}
        \ar@{=}[r]
        &pz
        \ar[d]_-{fI_z}^-{\cong}
        \ar@/^0.5pc/[dddr]|-{e_pI_z}
        \ar@{=}[r]
        &pz
        \ar@{->>}[d]^-{\xi^p_{\kos{t}^*_{\pi}(z)}}
        \\
        \omega^{p*}\kos{t}^*_{\pi}(z)
        \ar[d]_-{\vartheta^f_{\kos{t}^*_{\pi}(z)}}^-{\cong}
        &qz
        \ar@/^0.5pc/@{->>}[dl]|-{\xi^q_{\kos{t}^*_{\pi}(z)}}
        \ar[dd]^-{e_qI_z}
        &\omega^{p*}\kos{t}^*_{\pi}(z)
        \ar[dd]^-{\what{\omega}^{p*}_z}_-{\cong}
        \\
        \omega^{q*}\kos{t}^*_{\pi}(z)
        \ar[d]_-{\what{\omega}^{q*}_z}^-{\cong}
        &\text{ }
        &\text{ }
        \\
        \kappa
        \ar@{=}[r]
        &\kappa
        \ar@{=}[r]
        &\kappa
      }
    }}
  \end{equation*}
  This completes the proof of Lemma~\ref{lem twistbyatorsor functorial}.
\qed\end{proof}

\subsubsection{Torsors and pre-fiber functors}

\begin{theorem}
  \label{thm TorsFib adjequiv}
  Let $\mathfrak{T}=(\kos{T},\kos{t})$ be a pre-Galois $\kos{K}$-category
  and let $\varpi=(\omega,\what{\omega}):\mathfrak{K}\to\mathfrak{T}$
  be a pre-fiber functor for $\mathfrak{T}$.
  Recall Theorem~\ref{thm preGalKCat mainThm1}.
  We denote $\pi:=\omega^*\omega_!(\kappa)$
  as the pre-Galois object in $\kos{K}$
  which represents 
  $\underline{\Aut}_{\mathbb{GAL}^{\cat{pre}}_{\kos{K}}}(\varpi)$,
  and $\varpi$ factors through as equivalence of Galois $\kos{K}$-prekosmoi
  $\widebreve{\varpi}:\mathfrak{Rep}(\pi)\xrightarrow{\simeq}\mathfrak{T}$.
  \begin{equation*}
    \vcenter{\hbox{
      \xymatrix@C=40pt{
        \mathfrak{T}
        &\mathfrak{Rep}(\pi)
        \ar[l]_-{\widebreve{\varpi}}^-{\simeq}
        \\
        \text{ }
        &\mathfrak{K}
        \ar[u]_-{\varpi_{\pi}}
        \ar@/^0.5pc/[ul]^-{\varpi}
      }
    }}
  \end{equation*}
  Then we have an adjoint equivalence of groupoids
  \begin{equation*}
    \xymatrix{
      \cat{Fib}(\mathfrak{T})^{\cat{pre}}
      \ar@<0.5ex>[r]^-{\simeq}
      &(\cat{Tors-}\pi)_{\kappa}
      \ar@<0.5ex>[l]^-{\simeq}
    }
  \end{equation*}
  between the groupoid of pre-fiber functors for $\mathfrak{T}$
  and the groupoid of right $\pi$-torsors over $\kappa$.
  \begin{itemize}
    \item 
    Each pre-fiber functor $\varpi^{\pr}$
    for $\mathfrak{T}$
    is sent to the right $\pi$-torsor
    $\omega^{\pr*}\omega_!(\kappa)$ over $\kappa$,
    which represents the presheaf
    $\underline{\Isom}_{\mathbb{GAL}^{\cat{pre}}_{\kos{K}}}(\varpi,\varpi^{\pr})$
    of invertible Galois $\kos{K}$-transformations
    from $\varpi$ to $\varpi^{\pr}$.
    \begin{equation*}
      \vcenter{\hbox{
        \xymatrix{
          \Hom_{\cat{Ens}(\kos{K})}\bigl(\slot,\omega^{\pr*}\omega_*(\kappa)\bigr)
          \ar@2{->}[r]^-{\cong}
          &\underline{\Isom}_{\mathbb{GAL}^{\cat{pre}}_{\kos{K}}}(\varpi,\varpi^{\pr})
          :\cat{Ens}(\kos{K})^{\op}\to\cat{Set}
        }
      }}
    \end{equation*}
    The right $\pi$-action of $\omega^{\pr*}\omega_*(\kappa)$ is
    \begin{equation*}
      \xymatrix@C=35pt{
        \lambda_{\omega^{\pr*}\omega_!(\kappa)}:
        \omega^{\pr*}\omega_!(\kappa)\otimes \omega^*\omega_!(\kappa)
        \ar[r]^-{(\hatar{\omega^{\pr*}\omega_!}_{\omega^*\omega_!(\kappa)})^{-1}}_-{\cong}
        &\omega^{\pr*}\omega_!\omega^*\omega_!(\kappa)
        \ar[r]^-{\omega^{\pr*}(\epsilon_{\omega_!(\kappa)})}
        &\omega^{\pr*}\omega_!(\kappa)
        .
      }
    \end{equation*}

    \item 
    Each right $\pi$-torsor $p$ over $\kappa$
    is sent to the pre-fiber functor
    \begin{equation*}
      \widebreve{\varpi}\!\varpi^p:
      \xymatrix{
        \mathfrak{K}
        \ar[r]^-{\varpi^p}
        &\mathfrak{Rep}(\pi)
        \ar[r]^-{\widebreve{\varpi}}_-{\simeq}
        &\mathfrak{T}
      }
    \end{equation*}
    where $\varpi^p:\mathfrak{K}\to \mathfrak{Rep}(\pi)$
    is the pre-fiber functor
    twisted by the given torsor $p$
    which we introduced in Corollary~\ref{cor twistbyatorsor}.
  \end{itemize}
\end{theorem}
\begin{proof}
  We first describe the functor
  $\cat{Fib}(\mathfrak{T})^{\cat{pre}}
  \to (\cat{Tors}\text{-}\pi)_{\kappa}$.
  It sends each pre-fiber functor $\varpi^{\pr}:\mathfrak{K}\to\mathfrak{T}$
  to the right $\pi$-torsor $\omega^{\pr*}\omega_!(\kappa)$ over $\kappa$
  whose right $\pi$-action is
  \begin{equation*}
    \lambda_{\omega^{\pr*}\omega_!(\kappa)}:
    \xymatrix@C=45pt{
      \omega^{\pr*}\omega_!(\kappa)\otimes \omega^*\omega_!(\kappa)
      \ar[r]^-{\bigl(\hatar{\omega^{\pr*}\omega_!}_{\omega^*\omega_!(\kappa)}\bigr)^{-1}}_-{\cong}
      &\omega^{\pr*}\omega_!\omega^*\omega_!(\kappa)
      \ar[r]^-{\omega^{\pr*}(\epsilon_{\omega_!(\kappa)})}
      &\omega^{\pr*}\omega_!(\kappa)
      .
    }
  \end{equation*}
  We need to check that $(\omega^{\pr*}\omega_!(\kappa),\lambda_{\omega^{\pr*}\omega_!(\kappa)})$
  is indeed a right $\pi$-torsor over $\kappa$.
  The functor
  $\omega^{\pr*}\omega(\kappa)\otimes\slot:\CK\to \CK$
  is conservative and preserves reflexive coequalizers.
  This is because 
  we have
  $\hatar{\phi}^{\pr}:
  \xymatrix{
    \phi^{\pr}
    \ar@2{->}[r]^-{\cong}
    &\omega^{\pr*}\omega_!(\kappa)\otimes\slot
    :\CK\to \CK
  }$
  and the functor
  $\omega^{\pr*}\omega:\CK\to \CK$
  is conservative and preserves reflexive coequalizers.
  Recall that $\pi=\omega^*\omega_!(\kappa)$
  represents the presheaf of groups
  $\underline{\Aut}_{\mathbb{GAL}^{\cat{pre}}_{\kos{K}}}(\varpi)$
  and $\omega^{\pr*}\omega_!(\kappa)$
  represents the presheaf
  $\underline{\Isom}_{\mathbb{GAL}^{\cat{pre}}_{\kos{K}}}(\varpi,\varpi^{\pr})$
  By Lemma~\ref{lem RSKtensoradj composition},
  the morphism $\lambda_{\omega^{\pr*}\omega_!(\kappa)}$
  in $\cat{Ens}(\kos{K})$
  corresponds to the morphism of presheaves
  \begin{equation*}
    \underline{\Isom}_{\mathbb{GAL}^{\cat{pre}}_{\kos{K}}}(\varpi,\varpi^{\pr})
    \times
    \underline{\Aut}_{\mathbb{GAL}^{\cat{pre}}_{\kos{K}}}(\varpi)
    \Rightarrow
    \underline{\Isom}_{\mathbb{GAL}^{\cat{pre}}_{\kos{K}}}(\varpi,\varpi^{\pr})
  \end{equation*}
  obtained by the composition of Galois $\kos{K}$-transformations.
  This implies that $\lambda_{\omega^{\pr*}\omega_!(\kappa)}$ satisfies the right $\pi$-action relations.
  Moreover,
  the morphism of presheaves
  \begin{equation*}
    \vcenter{\hbox{
      \xymatrix@R=5pt@C=10pt{
        \underline{\Isom}_{\mathbb{GAL}^{\cat{pre}}_{\kos{K}}}(\varpi,\varpi^{\pr})
        \times
        \underline{\Aut}_{\mathbb{GAL}^{\cat{pre}}_{\kos{K}}}(\varpi)
        \ar@2{->}[r]
        &\underline{\Isom}_{\mathbb{GAL}^{\cat{pre}}_{\kos{K}}}(\varpi,\varpi^{\pr})
        \times
        \underline{\Isom}_{\mathbb{GAL}^{\cat{pre}}_{\kos{K}}}(\varpi,\varpi^{\pr})
        \\
        (\text{ }\vartheta^{\pr}
        \text{ },
        \text{ }\vartheta
        \text{ })
        \text{ }
        \ar@{|->}[r]
        &\text{ }
        (\text{ }
        \vartheta^{\pr}
        \text{ },\text{ }
        \vartheta^{\pr}\circ\vartheta
        \text{ })
      }
    }}
  \end{equation*}
  is an isomorphism of presheaves,
  whose inverse is given below.
  \begin{equation*}
    \vcenter{\hbox{
      \xymatrix@R=5pt@C=10pt{
        \underline{\Isom}_{\mathbb{GAL}^{\cat{pre}}_{\kos{K}}}(\varpi,\varpi^{\pr})
        \times
        \underline{\Aut}_{\mathbb{GAL}^{\cat{pre}}_{\kos{K}}}(\varpi)
        &\underline{\Isom}_{\mathbb{GAL}^{\cat{pre}}_{\kos{K}}}(\varpi,\varpi^{\pr})
        \times
        \underline{\Isom}_{\mathbb{GAL}^{\cat{pre}}_{\kos{K}}}(\varpi,\varpi^{\pr})
        \ar@2{->}[l]
        \\
        (\text{ }\vartheta^{\pr}\text{ }
        ,\text{ }
        \vartheta^{\pr-1}\circ \tilde{\vartheta}^{\pr}\text{ })
        \text{ }
        &\text{ }
        (\text{ }\vartheta^{\pr}\text{ }
        ,\text{ }
        \tilde{\vartheta}^{\pr}\text{ })
        \ar@{|->}[l]
      }
    }}
  \end{equation*}
  This implies that
  $\tau_{\omega^{\pr*}\omega_!(\kappa)}
  :\omega^{\pr*}\omega_!(\kappa)\otimes \pi
  \to \omega^{\pr*}\omega_!(\kappa)\otimes \omega^{\pr*}\omega_!(\kappa)$
  is an isomorphism.
  We conclude that
  $(\omega^{\pr*}\omega_!(\kappa),\lambda_{\omega^{\pr*}\omega_!(\kappa)})$
  is a right $\pi$-torsor over $\kappa$.
  Let $\varpi^{\ppr}:\mathfrak{K}\to\mathfrak{T}$
  be another pre-fiber functor and
  let
  $\vartheta:\!\!\xymatrix@C=15pt{\varpi^{\pr}\ar@2{->}[r]^-{\cong}&\varpi^{\ppr}}$
  be an invertible Galois $\kos{K}$-transformation.
  Then $\vartheta_{\omega_!(\kappa)}:\omega^{\pr*}\omega_!(\kappa)\to \omega^{\ppr*}\omega_!(\kappa)$
  is a morphism of right $\pi$-torsors over $\kappa$
  as we can see from the diagram below.
  \begin{equation*}
    \vcenter{\hbox{
      \xymatrix@C=40pt{
        \omega^{\pr*}\omega_!(\kappa)\otimes \omega^*\omega_!(\kappa)
        \ar[d]_-{\vartheta_{\omega_!(\kappa)}\otimes I_{\omega^*\omega_!(\kappa)}}
        \ar@/^0.5pc/[dr]|-{\bigl(\hatar{\omega^{\pr*}\omega_!}_{\omega^*\omega_!(\kappa)}\bigr)^{-1}}
        \ar@{=}[rr]
        &\text{ }
        &\omega^{\pr*}\omega_!(\kappa)\otimes \omega^*\omega_!(\kappa)
        \ar[dd]^-{\lambda_{\omega^{\pr*}\omega_!(\kappa)}}
        \\
        \omega^{\ppr*}\omega_!(\kappa)\otimes \omega^*\omega_!(\kappa)
        \ar[dd]_-{\lambda_{\omega^{\ppr*}\omega_!(\kappa)}}
        \ar@/_0.5pc/[dr]|-{\bigl(\hatar{\omega^{\ppr*}\omega_!}_{\omega^*\omega_!(\kappa)}\bigr)^{-1}}
        &\omega^{\pr*}\omega_!\omega^*\omega_!(\kappa)
        \ar[d]^-{\vartheta_{\omega_!\omega^*\omega_!(\kappa)}}
        \ar@/^0.5pc/[dr]|-{\omega^{\pr*}(\epsilon_{\omega_!(\kappa)})}
        &\text{ }
        \\
        \text{ }
        &\omega^{\ppr*}\omega_!\omega^*\omega_!(\kappa)
        \ar@/_0.5pc/[dr]|-{\omega^{\ppr*}(\epsilon_{\omega_!(\kappa)})}
        &\omega^{\pr*}\omega_!(\kappa)
        \ar[d]^-{\vartheta_{\omega_!(\kappa)}}
        \\
        \omega^{\ppr*}\omega_!(\kappa)
        \ar@{=}[rr]
        &\text{ }
        &\omega^{\ppr*}\omega_!(\kappa)
      }
    }}
  \end{equation*}
  This shows that the functor
  $\cat{Fib}(\mathfrak{T})^{\cat{pre}}
  \to (\cat{Tors}\text{-}\pi)_{\kappa}$
  is well-defined.
  Next we describe the functor
  $(\cat{Tors}\text{-}\pi)_{\kappa}\to
  \cat{Fib}(\mathfrak{T})^{\cat{pre}}$.
  It sends 
  each right $\pi$-torsor $p$ over $\kappa$
  to the pre-fiber functor
  \begin{equation*}
    \widebreve{\varpi}\!\varpi^p:
    \xymatrix{
      \mathfrak{K}
      \ar[r]^-{\varpi^p}
      &\mathfrak{Rep}(\pi)
      \ar[r]^-{\widebreve{\varpi}}_-{\simeq}
      &\mathfrak{T}
    }
  \end{equation*}
  where $\varpi^p:\mathfrak{K}\to \mathfrak{Rep}(\pi)$
  is the pre-fiber functor
  twisted by the given torsor $p$:
  see Corollary~\ref{cor twistbyatorsor}.
  Let $q$ be another right $\pi$-torsor over $\kappa$
  and let
  $f:p\xrightarrow{\cong}q$
  be an isomorphism of right $\pi$-torsors over $\kappa$.
  Then we have an isomorphism
  $\widebreve{\varpi}\vartheta^f:\widebreve{\varpi}\varpi^p\cong \widebreve{\varpi}\varpi^q$
  of pre-fiber functors for $\mathfrak{T}$,
  where $\vartheta^f:\varpi^p\cong \varpi^q$
  is described in Lemma~\ref{lem twistbyatorsor functorial}.
  This shows that the functor 
  $(\cat{Tors}\text{-}\pi)_{\kappa}\to
  \cat{Fib}(\mathfrak{T})^{\cat{pre}}$
  is well-defined.

  We explain the adjunction unit.
  Let $\varpi^{\pr}:\mathfrak{K}\to\mathfrak{Rep}(\pi)$
  be a pre-fiber functor
  and denote $p=\omega^{\pr*}\omega_!(\kappa)$.
  We claim that we have an isomorphism of pre-fiber functors
  \begin{equation*}
    \zeta^{\varpi^{\pr}}:
    \xymatrix{
      \varpi^{\pr}
      \ar@2{->}[r]^-{\cong}
      &\widebreve{\varpi}\varpi^p
      :\mathfrak{K}\to\mathfrak{T}
    }
  \end{equation*}
  whose component at each object $X$ in $\CT$ is
  the unique isomorphism
  $\zeta^{\varpi^{\pr}}_X:
  \omega^{\pr*}(X)\xrightarrow{\cong}
  \omega^{p*}\widebreve{\omega}\!^*(X)$
  in $\CK$ satisfying the relation below.
  Let us denote $\phi=\omega^*\omega_!:\kos{K}\to \kos{K}$
  as well as
  $\phi^{\pr}=\omega^{\pr*}\omega_!:\kos{K}\to \kos{K}$.
  \begin{equation} \label{eq TorsFib adjequiv}
    \vcenter{\hbox{
      \xymatrix@R=50pt@C=50pt{
        \omega^{\pr*}\omega_!\omega^*\omega_!\omega^*(X)
        \ar[d]_-{(\hatar{\phi}^{\pr}\hatar{\phi})_{\omega^*(X)}}^-{\cong}
        \ar@<0.5ex>[r]^-{\phi^{\pr}\omega^*(\epsilon_X)}
        \ar@<-0.5ex>[r]_-{\omega^{\pr*}(\epsilon_{\omega_!\omega^*(X)})}
        \ar@{}[dr]|-{(\dagger)}
        &\omega^{\pr*}\omega_!\omega^*(X)
        \ar[d]^-{\hatar{\phi}^{\pr}_{\omega^*(X)}}_-{\cong}
        \ar[r]^-{\omega^{\pr*}(\epsilon_X)}
        &\omega^{\pr*}(X)
        \ar@{.>}[d]^-{\exists!\text{ }\zeta^{\varpi^{\pr}}_X}_-{\cong}
        \\
        p\otimes \pi\otimes \omega^*(X)
        \ar@<0.5ex>[r]^-{I_p\otimes \gamma_{\omega^*(X)}}
        \ar@<-0.5ex>[r]_-{\lambda_p\otimes I_{\omega^*(X)}}
        &p\otimes \omega^*(X)
        \ar@{->>}[r]^-{\xi^p_{\widebreve{\omega}^*(X)}}
        &\omega^{p*}\widebreve{\omega}\!^*(X)
      }
    }}
  \end{equation}
  To conclude that the isomorphism
  $\zeta^{\varpi^{\pr}}_X$ is well-defined,
  we need to show that the diagram $(\dagger)$ commutes.
  We can check that the diagram $(\dagger)$ commutes with upper horizontal morphisms as follows.
  \begin{equation*}
    \vcenter{\hbox{
      \xymatrix{
        \phi^{\pr}\omega^*\omega_!\omega^*(X)
        \ar[d]_-{(\hatar{\phi}^{\pr}\hatar{\phi})_{\omega^*(X)}}^-{\cong}
        \ar@{=}[rr]
        &\text{ }
        &\phi^{\pr}\omega^*\omega_!\omega^*(X)
        \ar@/_1pc/[ddl]_-{\hatar{\phi}^{\pr}_{\omega^*\omega_!\omega^*(X)}}^-{\cong}
        \ar[d]^-{\phi^{\pr}\omega^*(\epsilon_X)}
        \\
        p\otimes \omega^*\omega_!(\kappa)\otimes \omega^*(X)
        \ar[dd]_-{I_p\otimes \gamma_{\omega^*(X)}}
        \ar@/_0.5pc/[dr]_-{I_p\otimes (\hatar{\phi}_{\omega^*(X)})^{-1}}^-{\cong}
        &\text{ }
        &\phi^{\pr}\omega^*(X)
        \ar[dd]^-{\hatar{\phi}^{\pr}_{\omega^*(X)}}_-{\cong}
        \\
        \text{ }
        &p\otimes \omega^*\omega_!\omega^*(X)
        \ar@/_0.5pc/[dr]|-{I_p\otimes \omega^*(\epsilon_X)}
        &\text{ }
        \\
        p\otimes \omega^*(X)
        \ar@{=}[rr]
        &\text{ }
        &p\otimes \omega^*(X)
      }
    }}
  \end{equation*}
  We can check that the diagram $(\dagger)$
  commutes with lower horizontal morphisms as follows.
  \begin{equation*}
    \vcenter{\hbox{
      \xymatrix@C=35pt{
        \omega^{\pr*}\omega_!\omega^*\omega_!\omega^*(X)
        \ar[d]_-{(\hatar{\phi}^{\pr}\hatar{\phi})_{\omega^*(X)}}^-{\cong}
        \ar@{}[dr]|-{(\dagger1)}
        \ar@{=}[rr]
        &\text{ }
        &\omega^{\pr*}\omega_!\omega^*\omega_!\omega^*(X)
        \ar@/^0.5pc/[dl]_-{\hatar{\phi^{\pr}\phi}_{\omega^*(X)}}^-{\cong}
        \ar[dd]^-{\omega^{\pr*}(\epsilon_{\omega_!\omega^*(X)})}
        \\
        p\otimes (\pi\otimes \omega^*(X))
        \ar[d]_-{a_{p,\pi,\omega^*(X)}}^-{\cong}
        &\omega^{\pr*}\omega_!\omega^*\omega_!(\kappa)\otimes \omega^*(X)
        \ar@/^0.5pc/[dl]_-{\hatar{\phi}^{\pr}_{\omega^*\omega_!(\kappa)}\otimes I_{\omega^*(X)}}^-{\cong}
        \ar@/_1pc/[ddr]|-{\omega^{\pr}(\epsilon_{\omega_!(\kappa)})\otimes I_{\omega^*(X)}}
        &\text{ }
        \\
        (p\otimes \pi)\otimes \omega^*(X)
        \ar[d]_-{\lambda_p\otimes I_{\omega^*(X)}}
        &\text{ }
        &\omega^{\pr*}\omega_!\omega^*(X)
        \ar[d]^-{\hatar{\phi}^{\pr}_{\omega^*(X)}}_-{\cong}
        \\
        p\otimes \omega^*(X)
        \ar@{=}[rr]
        &\text{ }
        &p\otimes \omega^*(X)
      }
    }}
  \end{equation*}
  We used the relation (\ref{eq Ens(T) KtensorEnd(K)=Ens(K)})
  in the diagram $(\dagger1)$.
  This shows that the isomorphism
  $\zeta^{\varpi^{\pr}}_X$ is well-defined.
  Next we show that
  $\zeta^{\varpi^{\pr}}$
  is a comonoidal $\kos{K}$-tensor natural isomorphism.
  Let $z\in\obj{\CK}$, $Y\in\obj{\CT}$
  and let us denote $\zeta=\zeta^{\varpi^{\pr}}$.
  We verify the relation
  \begin{equation*}
    \vcenter{\hbox{
      \xymatrix@C=30pt{
        \omega^{\pr*}(X\tensor Y)
        \ar[dd]_-{\omega^{\pr*}_{X,Y}}^-{\cong}
        \ar[r]^-{\zeta_{X\tensor Y}}_-{\cong}
        &\omega^{p*}\widebreve{\omega}\!^*(X\tensor Y)
        \ar[d]^-{\omega^{p*}(\widebreve{\omega}^*_{X,Y})}_-{\cong}
        \\
        \text{ }
        &\omega^{p*}(\widebreve{\omega}\!^*(X)\tensor\!_{\pi}\widebreve{\omega}\!^*(Y))
        \ar[d]^-{\omega^{p*}_{\widebreve{\omega}^*(X),\widebreve{\omega}^*(Y)}}_-{\cong}
        \\
        \omega^{\pr*}(X)\otimes \omega^{\pr*}(Y)
        \ar[r]^-{\zeta_X\otimes \zeta_Y}_-{\cong}
        &\omega^{p*}\widebreve{\omega}\!^*(X)\otimes \omega^{p*}\widebreve{\omega}\!^*(Y)
      }
    }}
  \end{equation*}
  by right-cancelling the epimorphism
  $\omega^{p*}(\epsilon_{X\tensor Y})$
  in the continued calculation of two diagrams presented below.
  \begin{equation*}
    \vcenter{\hbox{
      \xymatrix@C=40pt{
        \omega^{\pr*}\omega_!\omega^*(XY)
        \ar[d]_-{\omega^{\pr*}(\epsilon_{XY})}
        \ar@{=}[r]
        &\omega^{\pr*}\omega_!\omega^*(XY)
        \ar[d]_-{\hatar{\phi}^{\pr}_{\omega^*(XY)}}^-{\cong}
        \ar@{=}[r]
        &\omega^{\pr*}\omega_!\omega^*(XY)
        \ar[d]^-{\phi^{\pr}(\omega^*_{X,Y})}_-{\cong}
        \\
        \omega^{\pr*}(XY)
        \ar[d]_-{\zeta_{XY}}^-{\cong}
        &p\omega^*(XY)
        \ar@{->>}@/_0.5pc/[dl]|-{\xi^p_{\widebreve{\omega}^*(XY)}}
        \ar[d]_-{I_p\omega^*_{X,Y}}^-{\cong}
        &\phi^{\pr}(\omega^*(X)\omega^*(Y))
        \ar@/^0.5pc/[dl]_-{\hatar{\phi}^{\pr}_{\omega^*(X)\omega^*(Y)}}^-{\cong}
        \ar[d]^-{\phi^{\pr}_{\omega^*(X),\omega^*(Y)}}
        \\
        \omega^{p*}\widebreve{\omega}\!^*(XY)
        \ar[d]_-{\omega^{p*}(\widebreve{\omega}^*_{X,Y})}^-{\cong}
        &p\omega^*(X)\omega^*(Y)
        \ar@{->>}@/_0.5pc/[dl]|-{\xi^p_{\widebreve{\omega}^*(X)\widebreve{\omega}^*(Y)}}
        \ar[d]|-{(p\otimes)_{\omega^*(X),\omega^*(Y)}}
        &\phi^{\pr}\omega^*(X)\phi^{\pr}\omega^*(Y)
        \ar@/^0.5pc/[dl]_-{\hatar{\phi}^{\pr}_{\omega^*(X)}\hatar{\phi}^{\pr}_{\omega^*(Y)}}^-{\cong}
        \ar[d]^-{\omega^{\pr*}(\epsilon_X)\omega^{\pr*}(\epsilon_Y)}
        \\
        \omega^{p*}(\widebreve{\omega}\!^*(X)\widebreve{\omega}\!^*(Y))
        \ar[d]_-{\omega^{p*}_{\widebreve{\omega}\!^*(X),\widebreve{\omega}\!^*(Y)}}^-{\cong}
        &p\omega^*(X)p\omega^*(Y)
        \ar@/_0.5pc/@{->>}[dl]|-{\xi^p_{\widebreve{\omega}^*(X)}\xi^p_{\widebreve{\omega}^*(Y)}}
        &\omega^{\pr*}(X)\omega^{\pr*}(Y)
        \ar[d]^-{\zeta_X\zeta_Y}_-{\cong}
        \\
        \omega^{p*}\widebreve{\omega}\!^*(X)\omega^{p*}\widebreve{\omega}\!^*(Y)
        \ar@{=}[rr]
        &\text{ }
        &\omega^{p*}\widebreve{\omega}\!^*(X)\omega^{p*}\widebreve{\omega}\!^*(Y)
      }
    }}
  \end{equation*}
  \begin{equation*}
    \vcenter{\hbox{
      \xymatrix@C=15pt{
        \omega^{\pr*}\omega_!\omega^*(XY)
        \ar[d]_-{\phi^{\pr}(\omega^*_{X,Y})}^-{\cong}
        \ar@{=}[r]
        &\omega^{\pr*}\omega_!\omega^*(XY)
        \ar[d]^-{\omega^{\pr*}((\omega_!\omega^*)_{X,Y})}
        \ar@{=}[r]
        &\omega^{\pr*}\omega_!\omega^*(XY)
        \ar[dd]^-{\omega^{\pr*}(\epsilon_{XY})}
        \\
        \phi^{\pr}(\omega^*(X)\omega^*(Y))
        \ar[d]_-{\phi^{\pr}_{\omega^*(X),\omega^*(Y)}}
        &\omega^{\pr*}(\omega_!\omega^*(X)\omega_!\omega^*(Y))
        \ar@/^0.5pc/[dl]^-{\omega^{\pr*}_{\omega_!\omega^*(X),\omega_!\omega^*(Y)}}_-{\cong}
        \ar@/_0.5pc/[dr]|-{\omega^{\pr}(\epsilon_X\epsilon_Y)}
        &\text{ }
        \\
        \phi^{\pr}\omega^*(X)\phi^{\pr}\omega^*(Y)
        \ar[d]_-{\omega^{\pr*}(\epsilon_X)\omega^{\pr*}(\epsilon_Y)}
        &\text{ }
        &\omega^{\pr*}(XY)
        \ar[d]^-{\omega^{\pr*}_{X,Y}}_-{\cong}
        \\
        \omega^{\pr*}(X)\omega^{\pr*}(Y)
        \ar[d]_-{\zeta_X\zeta_Y}^-{\cong}
        \ar@{=}[rr]
        &\text{ }
        &\omega^{\pr*}(X)\omega^{\pr*}(Y)
        \ar[d]^-{\zeta_X\zeta_Y}_-{\cong}
        \\
        \omega^{p*}\widebreve{\omega}\!^*(X)\omega^{p*}\widebreve{\omega}\!^*(Y)
        \ar@{=}[rr]
        &\text{ }
        &\omega^{p*}\widebreve{\omega}\!^*(X)\omega^{p*}\widebreve{\omega}\!^*(Y)
      }
    }}
  \end{equation*}
  We also obtain the relations  
  \begin{equation*}
    \vcenter{\hbox{
      \xymatrix{
        \omega^{\pr*}(\unit)
        \ar[r]^-{\zeta_{\unit}}_-{\cong}
        \ar@/_1.5pc/[ddr]_-{\omega^{\pr*}_{\unit}}^-{\cong}
        &\omega^{p*}\widebreve{\omega}\!^*(\unit)
        \ar[d]^-{\omega^{p*}(\widebreve{\omega}\!^*_{\unit})}_-{\cong}
        \\
        \text{ }
        &\omega^{p*}(\unit\!_{\pi})
        \ar[d]^-{\omega^{p*}_{\unit_{\pi}}}_-{\cong}
        \\
        \text{ }
        &\kappa
      }
    }}
    \qquad\quad
    \vcenter{\hbox{
      \xymatrix{
        \omega^{\pr*}\kos{t}^*(z)
        \ar[r]^-{\zeta_{\kos{t}^*(z)}}_-{\cong}
        \ar@/_1.5pc/[ddr]_-{\what{\omega}^{\pr*}_z}^-{\cong}
        &\omega^{p*}\widebreve{\omega}\!^*\kos{t}^*(z)
        \ar[d]^-{\omega^{p*}(\what{\widebreve{\omega}}^*_z)}_-{\cong}
        \\
        \text{ }
        &\omega^{p*}\kos{t}^*_{\pi}(z)
        \ar[d]^-{\what{\omega}^{p*}_z}_-{\cong}
        \\
        \text{ }
        &z
      }
    }}
  \end{equation*}
  by right-cancelling the epimorphisms
  $\omega^{\pr*}(\epsilon_{\unit})$
  and 
  $\omega^{\pr*}(\epsilon_{\kos{t}^*(z)})$
  in the diagrams below.
  \begin{equation*}
    \vcenter
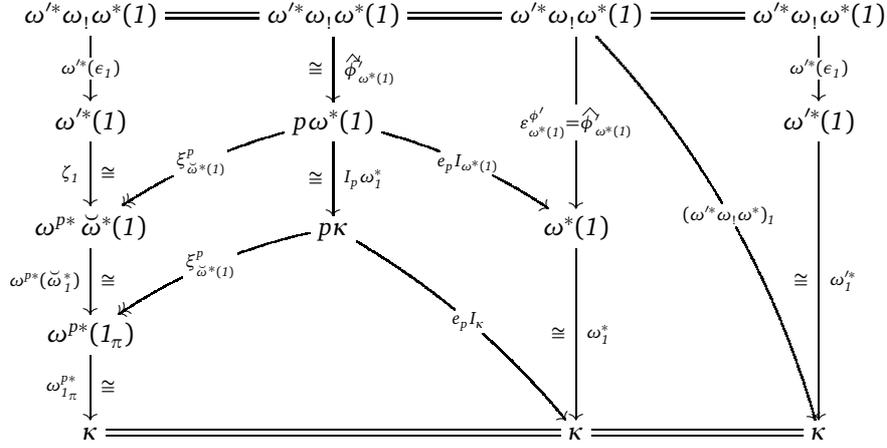
{\hbox{
      \xymatrix@C=35pt{
        \omega^{\pr*}\omega_!\omega^*(\unit)
        \ar[d]|-{\omega^{\pr*}(\epsilon_{\unit})}
        \ar@{=}[r]
        &\omega^{\pr*}\omega_!\omega^*(\unit)
        \ar[d]^-{\hatar{\phi}^{\pr}_{\omega^*(\unit)}}_-{\cong}
        \ar@{=}[r]
        &\omega^{\pr*}\omega_!\omega^*(\unit)
        \ar[dd]|-{\varepsilon^{\phi^{\pr}}_{\omega^*(\unit)}=\what{\phi}^{\pr}_{\omega^*(\unit)}}
        \ar@/^1pc/[ddddr]|-{(\omega^{\pr*}\omega_!\omega^*)_{\unit}}
        \ar@{=}[r]
        &\omega^{\pr*}\omega_!\omega^*(\unit)
        \ar[d]|-{\omega^{\pr*}(\epsilon_{\unit})}
        \\
        \omega^{\pr*}(\unit)
        \ar[d]_-{\zeta_{\unit}}^-{\cong}
        &p\omega^*(\unit)
        \ar@{->>}@/_0.5pc/[dl]|-{\xi^p_{\widebreve{\omega}^*(\unit)}}
        \ar[d]^-{I_p\omega^*_{\unit}}_-{\cong}
        \ar@/^0.5pc/[dr]|-{e_pI_{\omega^*(\unit)}}
        &\text{ }
        &\omega^{\pr*}(\unit)
        \ar[ddd]^-{\omega^{\pr*}_{\unit}}_-{\cong}
        \\
        \omega^{p*}\widebreve{\omega}\!^*(\unit)
        \ar[d]_-{\omega^{p*}(\widebreve{\omega}^*_{\unit})}^-{\cong}
        &p\kappa
        \ar@{->>}@/_0.5pc/[dl]|-{\xi^p_{\widebreve{\omega}^*(\unit)}}
        \ar@/^0.5pc/[ddr]|-{e_pI_{\kappa}}
        &\omega^*(\unit)
        \ar[dd]^-{\omega^*_{\unit}}_-{\cong}
        &\text{ }
        \\
        \omega^{p*}(\unit\!_{\pi})
        \ar[d]_-{\omega^{p*}_{\unit_{\pi}}}^-{\cong}
        &\text{ }
        &\text{ }
        &\text{ }
        \\
        \kappa
        \ar@{=}[rr]
        &\text{ }
        &\kappa
        \ar@{=}[r]
        &\kappa
      }
    }}
  \end{equation*}
  \begin{equation*}
    \vcenter
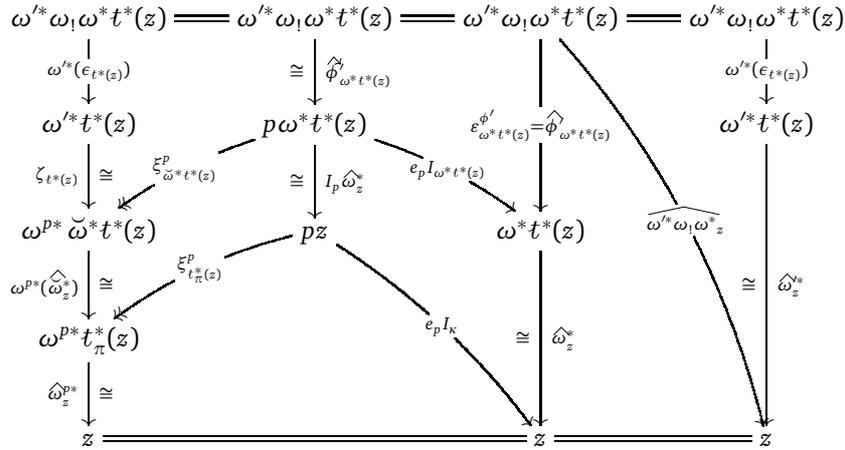
{\hbox{
      \xymatrix@C=20pt{
        \omega^{\pr*}\omega_!\omega^*\kos{t}^*(z)
        \ar[d]|-{\omega^{\pr*}(\epsilon_{\kos{t}^*(z)})}
        \ar@{=}[r]
        &\omega^{\pr*}\omega_!\omega^*\kos{t}^*(z)
        \ar[d]^-{\hatar{\phi}^{\pr}_{\omega^*\kos{t}^*(z)}}_-{\cong}
        \ar@{=}[r]
        &\omega^{\pr*}\omega_!\omega^*\kos{t}^*(z)
        \ar[dd]|-{\varepsilon^{\phi^{\pr}}_{\omega^*\kos{t}^*(z)}=\what{\phi}^{\pr}_{\omega^*\kos{t}^*(z)}}
        \ar@/^1pc/[ddddr]|-{\what{\omega^{\pr*}\omega_!\omega^*}_z}
        \ar@{=}[r]
        &\omega^{\pr*}\omega_!\omega^*\kos{t}^*(z)
        \ar[d]|-{\omega^{\pr*}(\epsilon_{\kos{t}^*(z)})}
        \\
        \omega^{\pr*}\kos{t}^*(z)
        \ar[d]_-{\zeta_{\kos{t}^*(z)}}^-{\cong}
        &p\omega^*\kos{t}^*(z)
        \ar@{->>}@/_0.5pc/[dl]|-{\xi^p_{\widebreve{\omega}^*\kos{t}^*(z)}}
        \ar[d]^-{I_p\what{\omega}^*_z}_-{\cong}
        \ar@/^0.5pc/[dr]|-{e_pI_{\omega^*\kos{t}^*(z)}}
        &\text{ }
        &\omega^{\pr*}\kos{t}^*(z)
        \ar[ddd]^-{\what{\omega}^{\pr*}_z}_-{\cong}
        \\
        \omega^{p*}\widebreve{\omega}\!^*\kos{t}^*(z)
        \ar[d]_-{\omega^{p*}(\what{\widebreve{\omega}}^*_z)}^-{\cong}
        &pz
        \ar@{->>}@/_0.5pc/[dl]|-{\xi^p_{\kos{t}^*_{\pi}(z)}}
        \ar@/^0.5pc/[ddr]|-{e_pI_{\kappa}}
        &\omega^*\kos{t}^*(z)
        \ar[dd]^-{\what{\omega}^*_z}_-{\cong}
        &\text{ }
        \\
        \omega^{p*}\kos{t}^*_{\pi}(z)
        \ar[d]_-{\what{\omega}^{p*}_z}^-{\cong}
        &\text{ }
        &\text{ }
        &\text{ }
        \\
        z
        \ar@{=}[rr]
        &\text{ }
        &z
        \ar@{=}[r]
        &z
      }
    }}
  \end{equation*}
  This shows that we have an isomorphism
  $\zeta^{\varpi^{\pr}}:\varpi^{\pr}\cong \widebreve{\varpi}\varpi^p$
  of pre-fiber functors as we claimed.
  To conclude that the adjunction unit is well-defined,
  we are left to show that $\zeta^{\varpi^{\pr}}$
  is natural in variable $\varpi^{\pr}$.
  Let $\varpi^{\ppr}:\mathfrak{K}\to\mathfrak{T}$ be another pre-fiber functor
  and let $\tilde{\vartheta}:\varpi^{\pr}\cong \varpi^{\ppr}$
  be an isomorphism of pre-fiber functors.
  Let us denote 
  $\phi^{\ppr}=\omega^{\ppr*}\omega_!:\kos{K}\to \kos{K}$
  and $q=\omega^{\ppr}\omega_!(\kappa)$
  as well as $f=\tilde{\vartheta}_{\omega_!(\kappa)}:p\xrightarrow{\cong}q$.
  We need to verify the relation below.
  \begin{equation*}
    \vcenter{\hbox{
      \xymatrix@C=40pt{
        \varpi^{\pr}
        \ar@2{->}[d]_-{\tilde{\vartheta}}^-{\cong}
        \ar@2{->}[r]^-{\zeta^{\varpi^{\pr}}}_-{\cong}
        &\widebreve{\varpi}\varpi^p
        \ar@2{->}[d]^-{\vartheta^f}_-{\cong}
        \\
        \varpi^{\ppr}
        \ar@2{->}[r]^-{\zeta^{\varpi^{\ppr}}}_-{\cong}
        &\widebreve{\varpi}\varpi^q
      }
    }}
  \end{equation*}
  We can check this by right-cancelling the epimorphism
  $\omega^{\pr*}(\epsilon_X)$ in the diagram below.
  Let $X$ be an object in $\CT$.
  \begin{equation*}
    \vcenter{\hbox{
      \xymatrix{
        \omega^{\pr*}\omega_!\omega^*(X)
        \ar[d]_-{\omega^{\pr*}(\epsilon_X)}
        \ar@{=}[r]
        &\omega^{\pr*}\omega_!\omega^*(X)
        \ar[d]^-{\hatar{\phi}^{\pr}_{\omega^*(X)}}_-{\cong}
        \ar@{=}[r]
        &\omega^{\pr*}\omega_!\omega^*(X)
        \ar[d]^-{\tilde{\vartheta}_{\omega_!\omega^*(X)}}_-{\cong}
        \ar@{=}[r]
        &\omega^{\pr*}\omega_!\omega^*(X)
        \ar[d]^-{\omega^{\pr*}(\epsilon_X)}
        \\
        \omega^{\pr*}(X)
        \ar[d]_-{\zeta^{\varpi^{\pr}}_X}^-{\cong}
        &p\otimes \omega^*(X)
        \ar@{->>}@/_0.5pc/[dl]|-{\xi^p_{\widebreve{\omega}^*(X)}}
        \ar[d]^-{f\otimes I_{\omega^*(X)}}
        &\omega^{\ppr*}\omega_!\omega^*(X)
        \ar@{->>}@/^0.5pc/[dl]^-{\hatar{\phi}^{\ppr}_{\omega^*(X)}}_-{\cong}
        \ar@/_0.5pc/[dr]|-{\omega^{\ppr*}(\epsilon_X)}
        &\omega^{\pr*}(X)
        \ar[d]^-{\tilde{\vartheta}_X}_-{\cong}
        \\
        \omega^{p*}\widebreve{\omega}\!^*(X)
        \ar[d]_-{\vartheta^f_{\widebreve{\omega}^*(X)}}^-{\cong}
        &q\otimes \omega^*(X)
        \ar@{->>}@/_0.5pc/@{->>}[dl]|-{\xi^q_{\widebreve{\omega}^*(X)}}
        &\text{ }
        &\omega^{\ppr*}(X)
        \ar[d]^-{\zeta^{\varpi^{\ppr}}_X}_-{\cong}
        \\
        \omega^{q*}\widebreve{\omega}\!^*(X)
        \ar@{=}[rrr]
        &\text{ }
        &\text{ }
        &\omega^{q*}\widebreve{\omega}\!^*(X)
      }
    }}
  \end{equation*}
  This shows that the adjunction unit is well-defined.
  Next we describe the adjunction counit.
  Let $p$ be a right $\pi$-torsor over $\kappa$.
  We have an isomorphism
  $f^p:\omega^{p*}\widebreve{\omega}\!^*\omega_!(\kappa)
  \xrightarrow{\cong} p$
  of right $\pi$-torsors over $\kappa$,
  which is the unique morphism satisfying the relation below.
  \begin{equation*}
    \vcenter{\hbox{
      \xymatrix@R=30pt@C=50pt{
        p\otimes \pi\otimes \omega^*\omega_!(\kappa)
        \ar@{=}[d]
        \ar@<0.5ex>[r]^-{I_p\otimes \gamma_{\omega^*\omega_!(\kappa)}}
        \ar@<-0.5ex>[r]_-{\lambda_p\otimes I_{\omega^*\omega_!(\kappa)}}
        &p\otimes \omega^*\omega_!(\kappa)
        \ar@{->>}[r]^-{\xi^p_{\widebreve{\omega}^*\omega_!(\kappa)}}
        \ar@{=}[d]
        &\omega^{p*}\widebreve{\omega}\!^*\omega_!(\kappa)
        \ar@{.>}[d]^-{\exists!\text{ }f^p}_-{\cong}
        \\
        p\otimes \pi\otimes \pi
        \ar@<0.5ex>[r]^-{I_p\otimes \pc_{\pi}}
        \ar@<-0.5ex>[r]_-{\lambda_p\otimes I_{\pi}}
        &p\otimes \pi
        \ar[r]^-{\lambda_p}
        &p
      }
    }}
  \end{equation*}
  Let $q$ be another right $\pi$-torsor over $\kappa$
  and $\tilde{f}:p\xrightarrow{\cong}q$
  be an isomorphism of right $\pi$-torsors over $\kappa$.
  Then we have the relation
  \begin{equation*}
    \vcenter{\hbox{
      \xymatrix@C=40pt{
        \omega^{p*}\widebreve{\omega}\!^*\omega_!(\kappa)
        \ar[r]^-{f^p}_-{\cong}
        \ar[d]_-{\vartheta^{\tilde{f}}_{\widebreve{\omega}^*\omega_!(\kappa)}}^-{\cong}
        &p
        \ar[d]^-{\tilde{f}}_-{\cong}
        \\
        \omega^{q*}\widebreve{\omega}\!^*\omega_!(\kappa)
        \ar[r]^-{f^q}_-{\cong}
        &q
      }
    }}
  \end{equation*}
  which we obtain by right-cancelling the epimorphism
  $\xi^p_{\widebreve{\omega}^*\omega_!(\kappa)}$ in the diagram below.
  \begin{equation*}
    \vcenter{\hbox{
      \xymatrix@C=50pt{
        p\otimes \omega^*\omega_!(\kappa)
        \ar@{->>}[d]_-{\xi^p_{\widebreve{\omega}^*\omega_!(\kappa)}}
        \ar@{=}[r]
        &p\otimes \omega^*\omega_!(\kappa)
        \ar[d]_-{\tilde{f}\otimes I_{\omega^*\omega_!(\kappa)}}
        \ar@/_0.5pc/[ddr]|-{\lambda_p}
        \ar@{=}[r]
        &p\otimes \omega^*\omega_!(\kappa)
        \ar@{->>}[d]^-{\xi^p_{\widebreve{\omega}^*\omega_!(\kappa)}}
        \\
        \omega^{p*}\widebreve{\omega}\!^*\omega_!(\kappa)
        \ar[d]_-{\vartheta^{\tilde{f}}_{\widebreve{\omega}^*\omega_!(\kappa)}}^-{\cong}
        &q\otimes \omega^*\omega_!(\kappa)
        \ar@/^0.5pc/[dl]|-{\xi^q_{\widebreve{\omega}^*\omega_!(\kappa)}}
        \ar[dd]^-{\lambda_q}
        &\omega^{p*}\widebreve{\omega}\!^*\omega_!(\kappa)
        \ar[d]^-{f^p}_-{\cong}
        \\
        \omega^{q*}\widebreve{\omega}\!^*\omega_!(\kappa)
        \ar[d]_-{f^q}^-{\cong}
        &\text{ }
        &p
        \ar[d]^-{\tilde{f}}_-{\cong}
        \\
        q
        \ar@{=}[r]
        &q
        \ar@{=}[r]
        &q
      }
    }}
  \end{equation*}
  This shows that the functors
  $\cat{Fib}(\mathfrak{T})^{\cat{pre}}
  \to (\cat{Tors-}\pi)_{\kappa}$
  and 
  $(\cat{Tors-}\pi)_{\kappa}
  \to \cat{Fib}(\mathfrak{T})^{\cat{pre}}$
  are quasi-inverse to each other.

  The final step is to show that the adjunction unit, counit
  satisfy the triangle relations.
  Let $\varpi^{\pr}:\mathfrak{K}\to\mathfrak{T}$
  be a pre-fiber functor
  and denote $\phi^{\pr}=\omega^{\pr*}\omega_!:\kos{K}\to \kos{K}$
  as well as $p=\omega^{\pr*}\omega_!(\kappa)$.
  We obtain one of the triangle relations
  \begin{equation*}
    \vcenter{\hbox{
      \xymatrix{
        \omega^{\pr*}\omega_!(\kappa)
        \ar[r]^-{\zeta^{\varpi^{\pr}}_{\omega_!(\kappa)}}_-{\cong}
        \ar@{=}@/_1pc/[dr]
        &\omega^{p*}\widebreve{\omega}\!^*\omega_!(\kappa)
        \ar[d]^-{f^p}_-{\cong}
        \\
        \text{ }
        &p
      }
    }}
  \end{equation*}
  by right-cancelling the epimorphism
  $\omega^{\pr*}(\epsilon_{\omega_!(\kappa)})$
  in the diagram below.
  \begin{equation*}
    \vcenter{\hbox{
      \xymatrix{
        \omega^{\pr*}\omega_!\omega^*\omega_!(\kappa)
        \ar[d]_-{\omega^{\pr*}(\epsilon_{\omega_!(\kappa)})}
        \ar@{=}[rr]
        &\text{ }
        &\omega^{\pr*}\omega_!\omega^*\omega_!(\kappa)
        \ar@/^0.5pc/[dl]_-{\hatar{\phi}^{\pr}_{\omega^*\omega_!(\kappa)}}^-{\cong}
        \ar[ddd]^-{\omega^{\pr*}(\epsilon_{\omega_!(\kappa)})}
        \\
        \omega^{\pr*}\omega_!(\kappa)
        \ar[d]_-{\zeta^{\varpi^{\pr}}_{\omega_!(\kappa)}}^-{\cong}
        &p\otimes \omega^*\omega_!(\kappa)
        \ar@{->>}@/^0.5pc/[dl]|-{\xi^p_{\widebreve{\omega}^*\omega_!(\kappa)}}
        \ar@/_1pc/[ddr]|-{\lambda_p}
        &\text{ }
        \\
        \omega^{p*}\widebreve{\omega}\!^*\omega_!(\kappa)
        \ar[d]_-{f^p}^-{\cong}
        &\text{ }
        &\text{ }
        \\
        p
        \ar@{=}[rr]
        &\text{ }
        &p
      }
    }}
  \end{equation*}
  Let $p$ be a right $\pi$-torsor over $\kappa$.
  Let us denote 
  $\varpi^{\pr}=\widebreve{\varpi}\omega^p$
  and $\phi=\omega^*\omega_!:\kos{K}\to \kos{K}$,
  $\phi^{\pr}=\omega^{\pr*}\omega_!:\kos{K}\to \kos{K}$
  as well as
  $\tilde{p}=\omega^{\pr*}\omega_!(\kappa)
  =\omega^{\pr*}\widebreve{\omega}\!^*\omega_!(\kappa)$.
  Recall that we have
  \begin{equation*}
    \xymatrix@C=40pt{
      \tilde{p}=
      \omega^{\pr*}\widebreve{\omega}\!^*\omega_!(\kappa)
      \ar[r]^-{f^p}_-{\cong}
      &p.
    }
  \end{equation*}
  We verify the other triangle relation
  \begin{equation*}
    \vcenter{\hbox{
      \xymatrix@C=50pt{
        \varpi^{\pr}
        \ar@/_1pc/@{=}[dr]
        \ar@2{->}[r]^-{\zeta^{\varpi^{\pr}}}_-{\cong}
        &\widebreve{\varpi}\omega^{\tilde{p}}
        \ar@2{->}[d]^-{\widebreve{\varpi}\vartheta^{f^p}}_-{\cong}
        \\
        \text{ }
        &\widebreve{\varpi}\omega^p
      }
    }}
  \end{equation*}
  by right-cancelling the epimorphisms
  $\xi^p_{\widebreve{\omega}^*\omega_!\omega^*(X)}$
  and
  $\omega^{\pr*}(\epsilon_X)$
  in the diagram below.
  Let $X$ be an object in $\CT$.
  \begin{equation*}
    \vcenter{\hbox{
      \xymatrix@C=2pt{
        p\!\otimes\! \omega^*\omega_!\omega^*(X)
        \ar@{->>}[d]|-{\xi^p_{\widebreve{\omega}^*\omega_!\omega^*(X)}}
        \ar@{=}[rrr]
        \ar@{}[drr]|-{(\ddagger)}
        &\text{ }
        &\text{ }
        \ar@{=}[r]
        &p\!\otimes\! \omega^*\omega_!\omega^*(X)
        \ar@/_0.3pc/[dl]_-{I_p\otimes \hatar{\phi}_{\omega^*(X)}}^-{\cong}
        \ar@/^1.5pc/[dddl]|(0.45){I_p\otimes \omega^*(\epsilon_X)}
        \ar@{->>}[ddd]|(0.7){\xi^p_{\widebreve{\omega}^*\omega_!\omega^*(X)}}
        \\
        \omega^{p*}\widebreve{\omega}\!^*\omega_!\omega^*(X)
        \ar[d]_-{\omega^{\pr*}(\epsilon_X)}
        \ar@/^0.5pc/[dr]^-{\hatar{\phi}^{\pr}_{\omega^*(X)}}_-{\cong}
        &\text{ }
        &p\!\otimes\! \omega^*\omega_!(\kappa)\!\otimes\! \omega^*(X)
        \ar@{->>}@/_0.5pc/[dl]|-{\xi^p_{\widebreve{\omega}^*\omega_!(\kappa)}}
        \ar@/^1.5pc/[ddl]|-{\lambda_p\otimes I_{\omega^*(X)}}
        \ar@<2ex>[dd]|-{I_p\otimes \gamma_{\omega^*\omega_!(\kappa)}}
        &\text{ }
        \\
        \omega^{p*}\widebreve{\omega}\!^*(X)
        \ar[d]_-{\zeta^{\varpi^{\pr}}_X}^-{\cong}
        &\omega^{p*}\widebreve{\omega}\!^*\omega_!(\kappa)\otimes \omega^*(X)
        \ar@{->>}@/^0.5pc/[dl]|-{\xi^{\tilde{p}}_{\widebreve{\omega}^*(X)}}
        \ar[d]^-{f^p\otimes I_{\omega^*(X)}}
        &\text{ }
        &\text{ }
        \\
        \omega^{\tilde{p}*}\widebreve{\omega}\!^*(X)
        \ar[d]_-{\vartheta^{f^p}_{\widebreve{\omega}^*(X)}}^-{\cong}
        &p\otimes \omega^*(X)
        \ar@/^0.5pc/@{->>}[dl]|-{\xi^p_{\widebreve{\omega}^*(X)}}
        &p\otimes \omega^*(X)
        \ar@/_0.5pc/@{->>}[dr]|-{\xi^p_{\widebreve{\omega}^*(X)}}
        &\omega^{p*}\widebreve{\omega}\!^*\omega_!\omega^*(X)
        \ar[d]^-{\omega^{\pr}(\epsilon_X)}
        \\
        \omega^{p*}\widebreve{\omega}\!^*(X)
        \ar@{=}[rrr]
        &\text{ }
        &\text{ }
        &\omega^{p*}\widebreve{\omega}\!^*(X)
      }
    }}
  \end{equation*}
  Note that the diagram $(\ddagger)$ is verified below.
  \begin{equation*}
    (\ddagger):
    \vcenter{\hbox{
      \xymatrix@C=5pt{
        p\otimes \omega^*\omega_!\omega^*(X)
        \ar[d]_-{I_p\otimes \hatar{\phi}_{\omega^*(X)}}^-{\cong}
        \ar@{=}[rr]
        &\text{ }
        &p\otimes \omega^*\omega_!\omega^*(X)
        \ar@{->>}[d]^-{\xi^p_{\widebreve{\omega}^*\omega_!\omega^*(X)}}
        \\
        p\otimes (\omega^*\omega_!(\kappa)\otimes \omega^*(X))
        \ar[d]_-{a_{p,\omega^*\omega_!(\kappa),\omega^*(X)}}^-{\cong}
        \ar@{->>}@/^0.5pc/[dr]|-{\xi^p_{\widebreve{\omega}^*\omega_!(\kappa)\otimes \omega^*(X)}}
        &\text{ }
        &\omega^{p*}\widebreve{\omega}\!^*\omega_!\omega^*(X)
        \ar[dd]^-{\hatar{\phi}^{\pr}_{\omega^*(X)}}_-{\cong}
        \ar@/_0.5pc/[dl]_-{\omega^{p*}\bigl(\hatar{\widebreve{\omega}\!^*\omega_!}_{\omega^*(X)}\bigr)}^-{\cong}
        \\
        (p\otimes \omega^*\omega_!(\kappa))\otimes \omega^*(X)
        \ar[d]_-{\xi^p_{\widebreve{\omega}^*\omega_!(\kappa)}\otimes I_{\omega^*(X)}}
        &\omega^{p*}(\widebreve{\omega}\!^*\omega_!(\kappa)\otimes \omega^*(X))
        \ar@/^0.5pc/[dr]_(0.4){\vecar{\omega}^{p*}_{\omega^*(X),\widebreve{\omega}^*\omega_!(\kappa)}}^-{\cong}
        &\text{ }
        \\
        \omega^{p*}\!\widebreve{\omega}\!^*\omega_!(\kappa)\otimes \omega^*(X)
        \ar@{=}[rr]
        &\text{ }
        &\omega^{p*}\!\widebreve{\omega}\!^*\omega_!(\kappa)\otimes \omega^*(X)
      }
    }}
  \end{equation*}
  This completes the proof of Theorem~\ref{thm TorsFib adjequiv}.
\qed\end{proof}

\newpage

\newpage
\section{Lax $\kos{K}$-tensor categories}
\label{sec LaxKtensorCat}
Throughout \textsection~\ref{sec LaxKtensorCat},
we fix a symmetric monoidal category
$\kos{K}=(\CK,\otimes,\kappa)$.
We study about
symmetric monoidal categories
$\kos{T}=(\CT,\tensor,\unit)$
equipped with lax symmetric monoidal functors
$\mon{t}:\kos{K}\to\kos{T}$,
and about the induced actions of $\kos{K}$ on $\kos{T}$.

\begin{definition} \label{def LaxKtensorCat}
  We define the $2$-category
  \begin{equation*}
    \mathbb{SMC}_{\cat{lax}}^{\kos{K}\backslash\!\!\backslash}
  \end{equation*}
  of
  lax $\kos{K}$-tensor categories,
  lax $\kos{K}$-tensor functors
  and
  monoidal $\kos{K}$-tensor natural transformations
  as the lax coslice $2$-category 
  of $\mathbb{SMC}_{\cat{lax}}$ under $\kos{K}$.
\end{definition}

Let us explain Definition~\ref{def LaxKtensorCat}
in detail.
A \emph{lax $\kos{K}$-tensor category}
$(\kos{T},\mon{t})$ is a pair
of a symmetric monoidal category
$\kos{T}=(\CT,\tensor,\unit)$
and a lax symmetric monoidal functor
$\mon{t}:\kos{K}\to\kos{T}$.
We denote the monoidal coherence morphisms of $\mon{t}$ as
\begin{equation}
  \label{eq LaxKtensorCat coherence}
  \mon{t}_{x,y}:
  \mon{t}(x)\tensor \mon{t}(y)
  \to
  \mon{t}(x\otimes y)
  ,
  \quad
  x,y\in\obj{\CK}
  ,
  \qquad
  \mon{t}_{\kappa}:
  \unit\to \mon{t}(\kappa)
  . 
\end{equation}
Let $(\kos{S},\mon{s})$ be another lax $\kos{K}$-tensor category
where $\text{$\kos{S}$}=(\CS,\ctimes,\pzc{1})$ is the underlying symmetric monoidal category.
A \emph{lax $\kos{K}$-tensor functor}
$(\phi,\what{\phi}):(\kos{S},\mon{s})\to (\kos{T},\mon{t})$
is a pair of a lax symmetric monoidal functor
$\phi:\kos{S}\to\kos{T}$
and a monoidal natural transformation
\begin{equation*}
  \vcenter{\hbox{
    \xymatrix@C=15pt{
      \text{ }
      &\text{ }
      &\kos{S}
      \ar[d]^-{\phi}
      \\
      \kos{K}
      \ar@/^0.7pc/[urr]^-{\mon{s}}
      \ar[rr]_-{\mon{t}}
      &\text{ }
      &\kos{T}
      \xtwocell[l]{}<>{<2.5>{\what{\phi}\text{ }\text{ }}}
    }
  }}
  \qquad\quad
  \what{\phi}:\mon{t}\Rightarrow\phi\mon{s}:\kos{K}\to \kos{S}
  .
\end{equation*}
Let $(\psi,\what{\psi}):(\kos{S},\mon{s})\to (\kos{T},\mon{t})$
be another lax $\kos{K}$-tensor functor.
A \emph{monoidal $\kos{K}$-tensor natural transformation}
$\vartheta:(\phi,\what{\phi})\Rightarrow (\psi,\what{\psi})
:(\kos{S},\mon{s})\to (\kos{T},\mon{t})$
is a monoidal natural transformation
$\vartheta:\phi\Rightarrow\psi:\kos{S}\to\kos{T}$
which satisfies the relation
\begin{equation}
  \label{eq LaxKtensorCat KtensorNat}
  \vcenter{\hbox{
    \xymatrix@C=20pt{
      \phi\mon{s}
      \ar@2{->}[rr]^-{\vartheta\mon{s}}
      &\text{ }
      &\psi\mon{s}
      \\
      \text{ }
      &\mon{t}
      \ar@2{->}[ul]^-{\what{\phi}}
      \ar@2{->}[ur]_-{\what{\psi}}
      &\text{ }
    }
  }}
  \qquad\quad
  \what{\psi}=(\vartheta\mon{s})\circ \what{\phi}:\mon{t}\Rightarrow\phi\mon{s}
  .
\end{equation}
We may also describe the relation
(\ref{eq LaxKtensorCat KtensorNat})
as the following diagram.
\begin{equation*}
  \vcenter{\hbox{
    \xymatrix@C=40pt{
      \text{ }
      &\kos{S}
      \ar[d]^-{\psi}
      \\
      \kos{K}
      \ar@/^0.7pc/[ur]^-{\mon{s}}
      \ar[r]_-{\mon{t}}
      &\kos{T}
      \xtwocell[l]{}<>{<2.5>{\what{\psi}\text{ }\text{ }}}
    }
  }}
  \quad=
  \vcenter{\hbox{
    \xymatrix@C=15pt{
      \text{ }
      &\text{ }
      &\text{ }
      &\kos{S}
      \ar@/^1.1pc/[d]^-{\psi}
      \ar@/_1.1pc/[d]_-{\phi}
      \\
      \kos{K}
      \ar@/^0.8pc/[urrr]^-{\mon{s}}
      \ar[rrr]_-{\mon{t}}
      &\text{ }
      &\text{ }
      &\kos{T}
      \xtwocell[u]{}<>{<0>{\vartheta}}
      \xtwocell[lll]{}<>{<2>{\what{\phi}\text{ }\text{ }}}
    }
  }}
\end{equation*}

We introduce a special kind of lax $\kos{K}$-tensor categories.

\begin{definition}
  A lax $\kos{K}$-tensor category
  $(\kos{T},\mon{t})$
  is called a \emph{strong $\kos{K}$-tensor category}
  if $\mon{t}:\kos{K}\to \kos{T}$
  is a strong $\kos{K}$-tensor functor,
  i.e., the monoidal coherence morphisms
  in (\ref{eq ColaxKtensorCat coherence})
  are isomorphisms.
\end{definition}

The pair $(\kos{K},\id_{\kos{K}})$
of $\kos{K}$ and the identity functor $\id_{\kos{K}}:\kos{K}\to \kos{K}$
is a strong $\kos{K}$-tensor category.

We also introduce a special kind of lax $\kos{K}$-tensor functors.

\begin{definition}
  Let 
  $(\kos{T},\mon{t})$, $(\kos{S},\mon{s})$
  be lax $\kos{K}$-tensor categories.
  A lax $\kos{K}$-tensor functor
  $(\phi,\what{\phi}):(\kos{S},\mon{s})\to (\kos{T},\mon{t})$
  is called a \emph{strong $\kos{K}$-tensor functor}
  if
  $\phi:\kos{S}\to\kos{T}$ is a strong symmetric monoidal functor
  and 
  $\what{\phi}:\mon{t}\Rightarrow\phi\mon{s}$
  is a monoidal natural isomorphism.
\end{definition}

Let $(\kos{T},\mon{t})$ be a lax $\kos{K}$-tensor category.
The pair of the lax symmetric monoidal functor
$\mon{t}:\kos{K}\to \kos{T}$
and the identity natural transformation
$I_{\mon{t}}$ of $\mon{t}$
defines a lax $\kos{K}$-tensor functor
$(\text{$\mon{t}$},I_{\text{$\mon{t}$}}):(\text{$\kos{K}$},\id_{\kos{K}})\to (\kos{T},\mon{t})$,
which is a strong $\kos{K}$-tensor functor
if and only if $(\kos{T},\mon{t})$ is a strong $\kos{K}$-tensor category.

We introduce the action of $\kos{K}$
associated to each lax $\kos{K}$-tensor category.

\begin{definition}
  Let $(\kos{T},\mon{t})$
  be a lax $\kos{K}$-tensor category
  where $\kos{T}=(\CT,\tensor,\unit)$
  is the underlying symmetric monoidal category.
  We define the functor
  \begin{equation*}
    \acts\!_{\mon{t}}:\CK\times \CT\to \CT,
    \qquad
    z\acts\!_{\mon{t}} X
    :=
    \mon{t}(z)\tensor X
    ,
    \qquad
    z\in \obj{\CK}
    ,
    \qquad
    X\in \obj{\CT}
  \end{equation*}
  and call it as the \emph{associated (lax) $\kos{K}$-action on $\kos{T}$}.
\end{definition}

We introduce how a lax $\kos{K}$-tensor functor
between lax $\kos{K}$-tensor categories
interacts with the associated $\kos{K}$-actions.

\begin{definition}
  Let $(\kos{T},\mon{t})$, $(\kos{S},\mon{s})$
  be lax $\kos{K}$-tensor categories
  where $\kos{T}=(\CT,\tensor,\unit)$,
  $\text{$\kos{S}$}=(\CS,\ctimes,\pzc{1})$
  are underlying symmetric monoidal categories.
  For each lax $\kos{K}$-tensor functor
  $(\phi,\what{\phi}):(\kos{S},\mon{s})\to (\kos{T},\mon{t})$,
  we define the natural transformation
  \begin{equation*}
    \cevar{\phi}_{\slot,\slot}:
    \slot\acts\!_{\mon{t}}\phi(\slot)
    \Rightarrow
    \phi(\slot\acts\!_{\mon{s}}\slot)
    :\CK\times \CS\to \CT
  \end{equation*}
  which we call as the \emph{associated (lax) $\kos{K}$-equivariance of $(\phi,\what{\phi})$},
  whose component at $z\in \obj{\CK}$, $\pzc{X}\in \obj{\CS}$ is
  described below.
  \begin{equation*}
    \xymatrix@C=40pt{
      \phi(z\acts\!_{\mon{s}} \pzc{X})
      \ar@{=}[d]
      &\text{ }
      &z\acts\!_{\mon{t}} \phi(\pzc{X})
      \ar@{=}[d]
      \ar[ll]_-{\cevar{\phi}_{z,\pzc{X}}}
      \\
      \phi(\text{$\mon{s}(z)$}\ctimes \text{$\pzc{X}$})
      &\text{$\phi\mon{s}(z)$}\tensor \phi(\text{$\pzc{X}$})
      \ar[l]_-{\phi_{\mon{s}(z),\pzc{X}}}
      &\text{$\mon{t}$}(z)\tensor \phi(\text{$\pzc{X}$})
      \ar[l]_-{\what{\phi}_z\tensor I_{\phi(\pzc{X})}}
    }
  \end{equation*}
\end{definition}

In particular, the $\kos{K}$-equivariance
$\cevar{\phi}$
associated to a strong $\kos{K}$-tensor functor
$(\phi,\what{\phi})$
between lax $\kos{K}$-tensor categories
is a natural isomorphism.

\begin{lemma}
  \label{lem LaxKtensorCat whatphi via cevarphi}
  Let 
  $(\phi,\what{\phi}):(\kos{S},\mon{s})\to (\kos{T},\mon{t})$
  be a lax $\kos{K}$-tensor functor
  between lax $\kos{K}$-tensor categories.
  The given monoidal natural transformation
  $\what{\phi}:\mon{t}\Rightarrow\phi\mon{s}$
  is described in terms of the associated $\kos{K}$-equivariance
  $\cevar{\phi}$ as follows.
  Let $z\in\obj{\CK}$.
  \begin{equation}\label{eq LaxKtensorCat whatphi via cevarphi}
    \what{\phi}_z:
    \xymatrix{
      \mon{t}(z)
      \ar[r]^-{\jmath_{\mon{t}(z)}}_-{\cong}
      &\mon{t}(z)\tensor \unit
      \ar[r]^-{I_{\mon{t}(z)}\tensor \phi_{\text{$\pzc{1}$}}}
      &\mon{t}(z)\tensor \phi(\text{$\pzc{1}$})
      \ar[r]^-{\cevar{\phi}_{z,\pzc{1}}}
      &\phi(\mon{s}(z)\ctimes \text{$\pzc{1}$})
      \ar[r]^-{\phi(\jmath_{\mon{s}(z)}^{-1})}_-{\cong}
      &\phi\mon{s}(z)
    }
  \end{equation}
\end{lemma}
\begin{proof}
  We can check the relation (\ref{eq LaxKtensorCat whatphi via cevarphi}) as follows.
  \begin{equation*}
    \vcenter{\hbox{
      \xymatrix@C=60pt{
        \mon{t}(z)
        \ar[d]_-{\jmath_{\mon{t}(z)}}^-{\cong}
        \ar@{=}[rr]
        &\text{ }
        &\mon{t}(z)
        \ar[d]^-{\what{\phi}_z}
        \\
        \mon{t}(z)\tensor \unit
        \ar[d]_-{I_{\mon{t}(z)}\tensor \phi_{\text{$\pzc{1}$}}}
        \ar@/^0.5pc/[dr]|-{\what{\phi}_z\tensor I_{\unit}}
        &\text{ }
        &\phi\mon{s}(z)
        \ar@/^0.5pc/[dl]^-{\jmath_{\phi\mon{s}(z)}}_-{\cong}
        \ar@{=}[dddd]
        \\
        \text{$\mon{t}$}(z)\tensor \phi(\text{$\pzc{1}$})
        \ar[dd]_-{\cevar{\phi}_{z,\text{$\pzc{1}$}}}
        \ar@/^0.5pc/[dr]|-{\what{\phi}_z\tensor I_{\phi(\text{$\pzc{1}$})}}
        &\phi\text{$\mon{s}$}(z)\tensor \unit
        \ar[d]^-{I_{\text{$\phi\mon{s}(z)$}}\tensor \phi_{\text{$\pzc{1}$}}}
        &\text{ }
        \\
        \text{ }
        &\phi\text{$\mon{s}$}(z)\tensor \phi(\text{$\pzc{1}$})
        \ar@/^0.5pc/[dl]|-{\phi_{\text{$\mon{s}(z)$},\text{$\pzc{1}$}}}
        &\text{ }
        \\
        \phi(\text{$\mon{s}$}(z)\ctimes \text{$\pzc{1}$})
        \ar[d]_-{\phi(\jmath^{-1}_{\mon{s}(z)})}^-{\cong}
        &\text{ }
        &\text{ }
        \\
        \phi\mon{s}(z)
        \ar@{=}[rr]
        &\text{ }
        &\phi\mon{s}(z)
      }
    }}
  \end{equation*}
  This completes the proof of Lemma~\ref{lem LaxKtensorCat whatphi via cevarphi}.
\qed\end{proof}

The following lemma describes the $\kos{K}$-equivariances
associated to compositions of lax $\kos{K}$-tensor functors.

\begin{lemma} \label{lem LaxKtensorCat compositionKequiv}
  Let $(\kos{T},\mon{t})$, $(\kos{T}^{\pr},\mon{t}^{\pr})$, $(\kos{S},\mon{s})$
  be lax $\kos{K}$-tensor categories
  and let
  \begin{equation*}
    \vcenter{\hbox{
      \xymatrix@C=30pt{
        (\kos{S},\mon{s})
        \ar[r]^-{(\phi,\what{\phi})}
        &(\kos{T},\mon{t})
        \ar[r]^-{(\psi,\what{\psi})}
        &(\kos{T}^{\pr},\mon{t}^{\pr})
      }
    }}
  \end{equation*}
  be lax $\kos{K}$-tensor functors.
  The $\kos{K}$-equivariance
  associated to the composition lax $\kos{K}$-tensor functor
  $(\psi\phi,\what{\psi\phi}):(\kos{S},\mon{s})\to (\kos{T}^{\pr},\mon{t}^{\pr})$
  is given as follows.
  Let $z\in \obj{\CK}$ and $\pzc{X}\in \obj{\CS}$.
  \begin{equation}\label{eq LaxKtensorCat compositionKequiv}
    \cevar{\psi\phi}_{z,\text{$\pzc{X}$}}:
    \xymatrix@C=30pt{
      z\acts\!_{\text{$\mon{t}^{\pr}$}}\psi\phi(\text{$\pzc{X}$})
      \ar[r]^-{\cevar{\psi}_{z,\text{$\phi(\pzc{X})$}}}
      &\psi(z\acts\!_{\text{$\mon{t}$}}\phi(\text{$\pzc{X}$}))
      \ar[r]^-{\psi(\cevar{\phi}_{z,\text{$\pzc{X}$}})}
      &\psi\phi(z\acts\!_{\text{$\mon{s}$}}\text{$\pzc{X}$})
    }
  \end{equation}
\end{lemma}
\begin{proof}
  The monoidal natural transformation
  $\what{\psi\phi}$
  of the composition lax $\kos{K}$-tensor functor is given by
  $\what{\psi\phi}:\!\!
  \xymatrix@C=17pt{
    \mon{t}^{\pr}
    \ar@2{->}[r]^-{\what{\psi}}
    &\psi\mon{t}
    \ar@2{->}[r]^-{\psi\what{\phi}}
    &\psi\phi\mon{s}
    .
  }$
  We can check the relation (\ref{eq LaxKtensorCat compositionKequiv}) as follows.
  \begin{equation*}
    \vcenter{\hbox{
      \xymatrix@C=60pt{
        \text{$\mon{t}$}^{\pr}(z)\tensor\!^{\pr}\psi\phi(\text{$\pzc{X}$})
        \ar[dd]_-{\cevar{\psi}_{z,\phi(\text{$\pzc{X}$})}}
        \ar@{=}[r]
        &\text{$\mon{t}$}^{\pr}(z)\tensor\!^{\pr}\psi\phi(\text{$\pzc{X}$})
        \ar[d]^-{\what{\psi}_z\tensor\!^{\pr}I_{\psi\phi(\text{$\pzc{X}$})}}
        \ar@{=}[r]
        &\text{$\mon{t}$}^{\pr}(z)\tensor\!^{\pr}\psi\phi(\text{$\pzc{X}$})
        \ar[dddd]^-{\cevar{\psi\phi}_{z,\text{$\pzc{X}$}}}
        \ar@/^1.8pc/[ddl]|-{\what{\psi\phi}_z\tensor\!^{\pr}I_{\psi\phi(\text{$\pzc{X}$})}}
        \\
        \text{ }
        &\psi\text{$\mon{t}$}(z)\tensor\!^{\pr}\psi\phi(\text{$\pzc{X}$})
        \ar@/_0.5pc/[dl]|-{\psi_{\text{$\mon{t}$}(z),\phi(\text{$\pzc{X}$})}}
        \ar[d]^-{\psi(\what{\phi}_z)\tensor\!^{\pr}I_{\psi\phi(\text{$\pzc{X}$})}}
        &\text{ }
        \\
        \psi(\text{$\mon{t}$}(z)\tensor \phi(\text{$\pzc{X}$}))
        \ar[dd]_-{\psi(\cevar{\phi}_{z,\text{$\pzc{X}$}})}
        \ar@/_0.5pc/[dr]|-{\psi(\what{\phi}_z\tensor I_{\phi(\text{$\pzc{X}$})})}
        &\psi\phi\text{$\mon{s}$}(z)\tensor\!^{\pr} \psi\phi(\text{$\pzc{X}$})
        \ar[d]^-{\psi_{\phi\text{$\mon{s}$}(z),\phi(\text{$\pzc{X}$})}}
        \ar@/^1.8pc/[ddr]|-{(\psi\phi)_{\text{$\mon{s}$}(z),\text{$\pzc{X}$}}}
        &\text{ }
        \\
        \text{ }
        &\psi(\phi\text{$\mon{s}$}(z)\tensor \phi(\text{$\pzc{X}$}))
        \ar[d]^-{\psi(\phi_{\text{$\mon{s}$}(z),\text{$\pzc{X}$}})}
        &\text{ }
        \\
        \psi\phi(\text{$\mon{s}$}(z)\ctimes \text{$\pzc{X}$})
        \ar@{=}[r]
        &\psi\phi(\text{$\mon{s}$}(z)\ctimes \text{$\pzc{X}$})
        \ar@{=}[r]
        &\psi\phi(\text{$\mon{s}$}(z)\ctimes \text{$\pzc{X}$})
      }
    }}
  \end{equation*}
  This completes the proof of Lemma~\ref{lem LaxKtensorCat compositionKequiv}.
\qed\end{proof}

The following lemma states that
a monoidal $\kos{K}$-tensor natural transformation
between lax $\kos{K}$-tensor functors
is precisely
a monoidal natural transformation
that is compatible with associated $\kos{K}$-equivariances.

\begin{lemma}
  \label{lem LaxKtensorCat KtensorNat preserves Kaction}
  Let $(\phi,\what{\phi})$,
  $(\psi,\what{\psi}):(\kos{S},\mon{s})\to (\kos{T},\mon{t})$
  be lax $\kos{K}$-tensor functors
  between lax $\kos{K}$-tensor categories
  and let
  $\vartheta:\phi\Rightarrow\psi:\kos{S}\to \kos{T}$
  be a monoidal natural transformation.
  Then
  $\vartheta:(\phi,\what{\phi})\Rightarrow (\psi,\what{\psi})
  :(\kos{S},\mon{s})\to (\kos{T},\mon{t})$
  is a monoidal $\kos{K}$-tensor natural transformation
  if and only if
  the following relation
  holds for all
  $z\in \obj{\CK}$, $\pzc{X}\in\obj{\CS}$.
  \begin{equation}\label{eq LaxKtensorCat KtensorNat preserves Kaction}
    \vcenter{\hbox{
      \xymatrix@C=50pt{
        \phi(z\acts\!_{\text{$\mon{s}$}}\text{$\pzc{X}$})
        \ar[d]_-{\vartheta_{z\acts\!_{\text{$\mon{s}$}}\text{$\pzc{X}$}}}
        &z\acts\!_{\text{$\mon{t}$}}\phi(\text{$\pzc{X}$})
        \ar[d]^-{I_z\acts\!_{\text{$\mon{t}$}}\vartheta_{\text{$\pzc{X}$}}}
        \ar[l]_-{\cevar{\phi}_{z,\pzc{X}}}
        \\
        \psi(z\acts\!_{\text{$\mon{s}$}}\text{$\pzc{X}$})
        &z\acts\!_{\text{$\mon{t}$}}\psi(\text{$\pzc{X}$})
        \ar[l]_-{\cevar{\psi}_{z,\pzc{X}}}
      }
    }}
  \end{equation}
\end{lemma}
\begin{proof}
  We first prove the only if part.
  Suppose that $\vartheta:(\phi,\what{\phi})\Rightarrow (\psi,\what{\psi})$
  is a monoidal $\kos{K}$-tensor natural transformation.
  We obtain the relation (\ref{eq LaxKtensorCat KtensorNat preserves Kaction})
  as follows.
  \begin{equation*}
    \vcenter
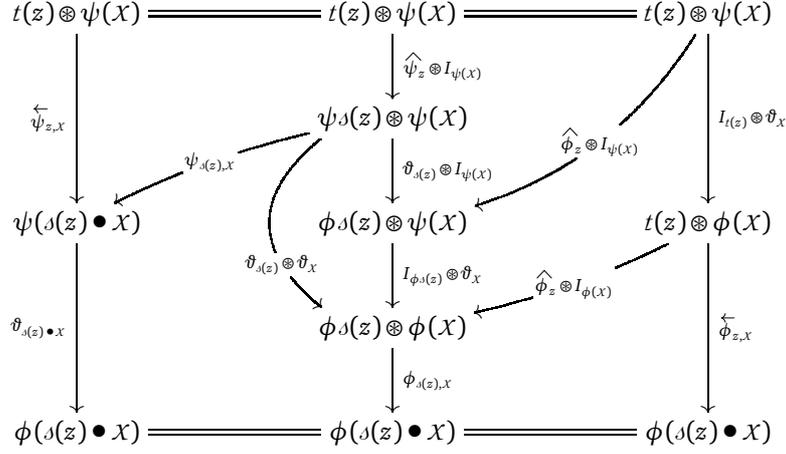
{\hbox{
      \xymatrix@C=60pt{
        \text{$\mon{t}$}(z)\tensor \psi(\text{$\pzc{X}$})
        \ar[dd]_-{\cevar{\psi}_{z,\text{$\pzc{X}$}}}
        \ar@{=}[r]
        &\text{$\mon{t}$}(z)\tensor \psi(\text{$\pzc{X}$})
        \ar[d]^-{\what{\psi}_z\tensor I_{\psi(\text{$\pzc{X}$})}}
        \ar@{=}[r]
        &\text{$\mon{t}$}(z)\tensor \psi(\text{$\pzc{X}$})
        \ar[dd]^-{I_{\text{$\mon{t}$}(z)}\tensor \vartheta_{\text{$\pzc{X}$}}}
        \ar@/^1.5pc/[ddl]|-{\what{\phi}_z\tensor I_{\psi(\text{$\pzc{X}$})}}
        \\
        \text{ }
        &\psi\text{$\mon{s}$}(z)\tensor \psi(\text{$\pzc{X}$})
        \ar[d]^-{\vartheta_{\text{$\mon{s}$}(z)}\tensor I_{\psi(\text{$\pzc{X}$})}}
        \ar@/_0.5pc/[dl]|-{\psi_{\text{$\mon{s}$}(z),\text{$\pzc{X}$}}}
        \ar@<-3.5ex>@/_2.5pc/[dd]|(0.7){\vartheta_{\text{$\mon{s}$}(z)}\tensor \vartheta_{\text{$\pzc{X}$}}}
        &\text{ }
        \\
        \psi(\text{$\mon{s}$}(z)\ctimes \text{$\pzc{X}$})
        \ar[dd]_-{\vartheta_{\text{$\mon{s}$}(z)\ctimes \text{$\pzc{X}$}}}
        &\phi\text{$\mon{s}$}(z)\tensor \psi(\text{$\pzc{X}$})
        \ar[d]^-{I_{\phi\text{$\mon{s}$}(z)}\tensor \vartheta_{\text{$\pzc{X}$}}}
        &\text{$\mon{t}$}(z)\tensor \phi(\text{$\pzc{X}$})
        \ar[dd]^-{\cevar{\phi}_{z,\text{$\pzc{X}$}}}
        \ar@/^0.5pc/[dl]|-{\what{\phi}_z\tensor I_{\phi(\text{$\pzc{X}$})}}
        \\
        \text{ }
        &\phi\text{$\mon{s}$}(z)\tensor \phi(\text{$\pzc{X}$})
        \ar[d]^-{\phi_{\text{$\mon{s}$}(z),\text{$\pzc{X}$}}}
        &\text{ }
        \\
        \phi(\text{$\mon{s}$}(z)\ctimes \text{$\pzc{X}$})
        \ar@{=}[r]
        &\phi(\text{$\mon{s}$}(z)\ctimes \text{$\pzc{X}$})
        \ar@{=}[r]
        &\phi(\text{$\mon{s}$}(z)\ctimes \text{$\pzc{X}$})
      }
    }}
  \end{equation*}
  Next we prove the if part.
  Assume that $\vartheta:\phi\Rightarrow\psi$
  satisfies the relation (\ref{eq LaxKtensorCat KtensorNat preserves Kaction})
  for all $z$, $\pzc{X}$.
  Using the descriptions
  of $\what{\phi}$, $\what{\psi}$
  given in (\ref{eq LaxKtensorCat whatphi via cevarphi}),
  we can check that
  $\what{\psi}=(\vartheta\mon{s})\circ \what{\phi}$
  as follows.
  \begin{equation*}
    \vcenter
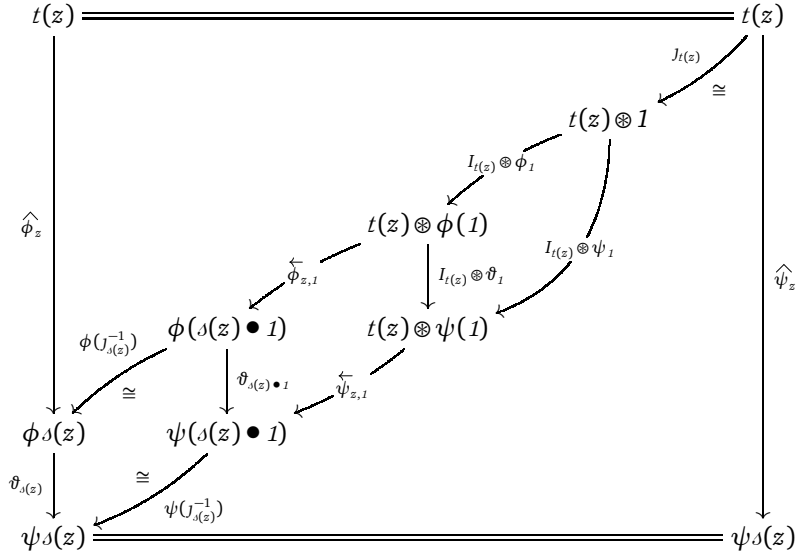
{\hbox{
      \xymatrix{
        \mon{t}(z)
        \ar[dddd]_-{\what{\phi}_z}
        \ar@{=}[rrrr]
        &\text{ }
        &\text{ }
        &\text{ }
        &\mon{t}(z)
        \ar[ddddd]^-{\what{\psi}_z}
        \ar@/^0.5pc/[dl]_-{\jmath_{\mon{t}(z)}}^-{\cong}
        \\
        \text{ }
        &\text{ }
        &\text{ }
        &\mon{t}(z)\tensor \unit
        \ar@/_0.5pc/[dl]|-{I_{\text{$\mon{t}$}(z)}\tensor \phi_{\text{$\pzc{1}$}}}
        \ar@/^2pc/[ddl]|-{I_{\text{$\mon{t}$}(z)}\tensor \psi_{\text{$\pzc{1}$}}}
        &\text{ }
        \\
        \text{ }
        &\text{ }
        &\text{$\mon{t}$}(z)\tensor \phi(\text{$\pzc{1}$})
        \ar[d]^-{I_{\text{$\mon{t}$}(z)}\tensor \vartheta_{\text{$\pzc{1}$}}}
        \ar@/_0.5pc/[dl]|-{\cevar{\phi}_{z,\text{$\pzc{1}$}}}
        &\text{ }
        &\text{ }
        \\
        \text{ }
        &\phi(\text{$\mon{s}$}(z)\ctimes \text{$\pzc{1}$})
        \ar[d]^-{\vartheta_{\text{$\mon{s}$}(z)\ctimes \text{$\pzc{1}$}}}
        \ar@/_0.5pc/[dl]_(0.5){\phi(\jmath_{\mon{s}(z)}^{-1})}^-{\cong}
        &\text{$\mon{t}$}(z)\tensor \psi(\text{$\pzc{1}$})
        \ar@/^0.5pc/[dl]|-{\cevar{\psi}_{z,\text{$\pzc{1}$}}}
        &\text{ }
        &\text{ }
        \\
        \phi\mon{s}(z)
        \ar[d]_-{\vartheta_{\mon{s}(z)}}
        &\psi(\text{$\mon{s}$}(z)\ctimes \text{$\pzc{1}$})
        \ar@/^0.5pc/[dl]^(0.4){\psi(\jmath_{\mon{s}(z)}^{-1})}_-{\cong}
        &\text{ }
        &\text{ }
        &\text{ }
        \\
        \psi\mon{s}(z)
        \ar@{=}[rrrr]
        &\text{ }
        &\text{ }
        &\text{ }
        &\psi\mon{s}(z)
      }
    }}
  \end{equation*}
  Therefore
  $\vartheta:(\phi,\what{\phi})\Rightarrow (\psi,\what{\psi})$
  is a monoidal $\kos{K}$-tensor natural transformation.
  This completes the proof of Lemma~\ref{lem LaxKtensorCat KtensorNat preserves Kaction}.
\qed\end{proof}

\subsection{Commutative monoids as lax $\kos{K}$-tensor functors}
\label{subsec Comm(T)}
Let $\kos{T}=(\CT,\tensor,\unit)$
be a symmetric monoidal category.
We denote
\begin{equation*}
  \kos{C\!o\!m\!m}(\kos{T})
  =(\cat{Comm}(\kos{T}),\tensor,\unit)
\end{equation*}
as the cocartesian monoidal category
of commutative monoids in $\kos{T}$.
An object in $\cat{Comm}(\kos{T})$
is a triple
$(B,\pc_B,u_B)$ 
where $B$ is an object in $\CT$ and 
$\pc_B:B\to B\tensor B$,
$u_B:\unit\to B$
are product, unit morphisms in $\CT$
satisfying the associativity, unital, commutativity relations.
We often omit product, unit morphisms
and simply denote an object
$(B,\pc_B,u_B)$ in $\cat{Comm}(\kos{T})$ as $B$.
A morphism $B\to B^{\pr}$ in $\cat{Comm}(\kos{T})$
is a morphism in $\CT$
which is compatible with products, units morphisms of $B$, $B^{\pr}$.
The category $\cat{Comm}(\kos{T})$ has finite coproducts.
The object $\unit$ equipped with
$\pc_{\unit}=\imath^{-1}_{\unit}=\jmath^{-1}_{\unit}:\unit\tensor \unit\xrightarrow{\cong}\unit$
and
$u_{\unit}=I_{\unit}:\unit\xrightarrow{\cong}\unit$
is an initial object in $\cat{Comm}(\kos{T})$.
The product of objects $B$, $B^{\pr}$ in $\cat{Comm}(\kos{T})$
is $B\tensor B^{\pr}$ whose product, unit morphisms are described below.
\begin{equation} \label{eq Comm(T) BB' definition}
  \begin{aligned}
    \pc_{B\tensor B^{\pr}}
    &:
    \xymatrix@C=45pt{
      B\tensor B^{\pr}\tensor B\tensor B^{\pr}
      \ar[r]^-{I_B\tensor s_{B^{\pr},B}\tensor I_{B^{\pr}}}_-{\cong}
      &B\tensor B\tensor B^{\pr}\tensor B^{\pr}
      \ar[r]^-{\pc_B\tensor \pc_{B^{\pr}}}
      &B\tensor B^{\pr}
    }
    \\
    u_{B\tensor B^{\pr}}
    &:
    \xymatrix@C=30pt{
      \unit
      \ar[r]^-{\pc_{\unit}^{-1}}_-{\cong}
      &\unit\tensor \unit
      \ar[r]^-{u_B\tensor u_{B^{\pr}}}
      &B\tensor B^{\pr}
    }
  \end{aligned}
\end{equation}
The symmetric monoidal coherence isomorphisms
$a$, $\imath$, $\jmath$, $s$
of
$\kos{C\!o\!m\!m}(\kos{T})$
are given by those of $\kos{T}$.
We also denote 
\begin{equation*}
  \kos{A\!f\!f}(\kos{T})
  =
  (\cat{Aff}(\kos{T}),\times,\cat{Spec}(\unit))
\end{equation*}
as the opposite cartesian category of
$\kos{C\!o\!m\!m}(\kos{T})
=(\cat{Comm}(\kos{T}),\tensor,\unit)$.
For each object $B$ in $\cat{Comm}(\kos{T})$,
the corresponding object in $\cat{Aff}(\kos{T})$
is denoted as
$\cat{Spec}(B)$.
The cartesian product of objects
$\cat{Spec}(B)$, $\cat{Spec}(B^{\pr})$
in $\cat{Aff}(\kos{T})$
is
$\cat{Spec}(B)\times \cat{Spec}(B^{\pr})
=\cat{Spec}(B\tensor B^{\pr})$.

For each object $B$ in $\cat{Comm}(\kos{T})$,
we have a lax symmetric monoidal monad 
\begin{equation*}
  \langle \tensor B\rangle
  =(\tensor B,\mu^{\tensor B},\upsilon^{\tensor B})
\end{equation*}
on $\kos{T}$.  
The underlying lax symmetric monoidal endofunctor
$\tensor B:\kos{T}\to \kos{T}$
is the endofunctor
$\tensor B=\slot\tensor B:\CT\to \CT$
equipped with the following monoidal coherence morphisms
\begin{equation}\label{eq Comm(T) tensorB definition}
  \begin{aligned}
    (\tensor B)_{X,Y}
    &:
    \xymatrix@C=40pt{
      X\tensor B\tensor Y\tensor B
      \ar[r]^-{I_X\tensor s_{B,Y}\tensor I_B}_-{\cong}
      &X\tensor Y\tensor B\tensor B
      \ar[r]^-{I_{X\tensor Y}\tensor \pc_B}
      &X\tensor Y\tensor B
    }
    \\
    (\tensor B)_{\unit}
    &:
    \xymatrix@C=30pt{
      \unit
      \ar[r]^-{\pc_{\unit}^{-1}}_-{\cong}
      &\unit\tensor \unit
      \ar[r]^-{I_{\unit}\tensor u_B}
      &\unit\tensor B
    }
    \qquad
    X,Y\in\obj{\CT}
  \end{aligned}
\end{equation}
and the monoidal natural transformations
$\mu^{\tensor B}$, $\upsilon^{\tensor B}$
are described below.
\begin{equation}\label{eq Comm(T) muBupsilonB definition}
  \begin{aligned}
    \mu^{\tensor B}_X
    &:
    \xymatrix@C=30pt{
      (X\tensor B)\tensor B
      \ar[r]^-{a^{-1}_{X,B,B}}_-{\cong}
      &X\tensor (B\tensor B)
      \ar[r]^-{I_X\tensor \pc_B}
      &X\tensor B
    }
    \\
    \upsilon^{\tensor B}_X
    &:
    \xymatrix@C=30pt{
      X
      \ar[r]^-{\jmath_X}_-{\cong}
      &X\tensor \unit
      \ar[r]^-{I_X\tensor u_B}
      &X\tensor B
    }
    \qquad\quad
    X\in\obj{\CT}
  \end{aligned}
\end{equation}

\begin{proposition} \label{prop Comm(T) Comm(T)laxK(K,T) adjunction}
  Let $(\kos{T},\mon{t})$ be a strong $\kos{K}$-tensor category.
  We have a coreflective adjunction
  \begin{equation*}
    \iota:
    \cat{Comm}(\kos{T})
    \rightleftarrows
    \mathbb{SMC}_{\cat{lax}}^{\kos{K}\backslash\!\!\backslash}\bigl((\kos{K},\id_{\kos{K}}),(\kos{T},\mon{t})\bigr)
    :\CR
  \end{equation*}
  between 
  the category of commutative monoids in $\kos{T}$
  and 
  the category of lax $\kos{K}$-tensor functors
  $(\kos{K},\id_{\kos{K}})\to (\kos{T},\mon{t})$.
  \begin{itemize}
    \item 
    The left adjoint $\iota$
    sends each object
    $B$ in $\cat{Comm}(\kos{T})$
    to the lax $\kos{K}$-tensor functor
    \begin{equation*}
      \iota(B)
      :=
      (\text{$\mon{t}$}\tensor B,\what{\mon{t}\tensor B}):
      (\kos{K},\id_{\kos{K}})
      \to
      (\kos{T},\mon{t})
      .
    \end{equation*}
    The underlying lax symmetric monoidal functor is
    $\text{$\mon{t}$}\tensor B:
    \kos{K}
    \xrightarrow[]{\text{$\mon{t}$}}
    \kos{T}
    \xrightarrow[]{\tensor B}
    \kos{T}$
    and the component of the monoidal natural transformation
    $\what{\mon{t}\tensor B}$
    at $x\in \obj{\CK}$ is
    \begin{equation*}
      \vcenter{\hbox{
        \xymatrix@C=15pt{
          \text{ }
          &\text{ }
          &\kos{K}
          \ar[d]^-{\mon{t}\tensor B}
          \\
          \kos{K}
          \ar@/^0.7pc/[urr]^-{\id_{\kos{K}}}
          \ar[rr]_-{\mon{t}}
          &\text{ }
          &\kos{T}
          \xtwocell[l]{}<>{<2.5>{\what{\mon{t}\tensor B}\quad\text{ }\text{ }}}
        }
      }}
      \qquad
      \what{\mon{t}\tensor B}_x
      =
      \upsilon^{\tensor B}_{\mon{t}(x)}
      :\!\!
      \xymatrix@C=25pt{
        \mon{t}(x)
        \ar[r]^-{\jmath_{\mon{t}(x)}}_-{\cong}
        &\mon{t}(x)\tensor \unit
        \ar[r]^-{I_{\mon{t}(x)}\tensor u_B}
        &\mon{t}(x)\tensor B
        .
      }
    \end{equation*}

    \item 
    The right adjoint $\CR$ sends each lax $\kos{K}$-tensor functor
    $(\phi,\what{\phi}):(\kos{K},\id_{\kos{K}})\to (\kos{T},\mon{t})$
    to the object
    $\CR(\phi,\what{\phi}):=\phi(\kappa)$
    in $\cat{Comm}(\kos{T})$,
    where
    \begin{equation*}
      \pc_{\phi(\kappa)}:\!
      \xymatrix{
        \phi(\kappa)\tensor \phi(\kappa)
        \ar[r]^-{\phi_{\kappa,\kappa}}
        &\phi(\kappa\otimes\kappa)
        \ar[r]^-{\phi(\pc_{\kappa})}_-{\cong}
        &\phi(\kappa)
        ,
      }
      \quad
      u_{\phi(\kappa)}:\!
      \xymatrix{
        \unit
        \ar[r]^-{\phi_{\kappa}}
        &\phi(\kappa)
        .
      }
    \end{equation*}

    \item 
    The component of the adjunction unit at each object
    $B$ in $\cat{Comm}(\kos{T})$ is
    \begin{equation*}
      \xymatrix@C=30pt{
        B
        \ar[r]^-{\imath_B}_-{\cong}
        &\unit\tensor B
        \ar[r]^-{\mon{t}_{\kappa}\tensor I_B}_-{\cong}
        &\mon{t}(\kappa)\tensor B
        .
      }
    \end{equation*}

    \item
    The component of the adjunction counit at each
    lax $\kos{K}$-tensor functor
    $(\phi,\what{\phi}):(\kos{K},\id_{\kos{K}})\to (\kos{T},\mon{t})$
    is the monoidal $\kos{K}$-tensor natural transformation
    \begin{equation*}
      \tahar{\phi}
      :
      (\text{$\mon{t}$}\tensor\phi(\kappa),\what{\mon{t}\tensor\phi(\kappa)})
      \Rightarrow
      (\phi,\what{\phi})
      :(\kos{K},\id_{\kos{K}})\to (\kos{T},\mon{t})
    \end{equation*}
    whose component at $x\in\obj{\CK}$ is
    \begin{equation}
      \label{eq Comm(T) Comm(T)laxK(K,T) adjunction}
      \tahar{\phi}_x:
      \xymatrix@C=25pt{
        \mon{t}(x)\tensor \phi(\kappa)
        \ar@/_1pc/@<-1ex>[rr]|-{\cevar{\phi}_{x,\kappa}}
        \ar[r]^-{\what{\phi}_x\tensor I_{\phi(\kappa)}}
        &\phi(x)\tensor \phi(\kappa)
        \ar[r]^-{\phi_{x,\kappa}}
        &\phi(x\otimes \kappa)
        \ar[r]^-{\phi(\jmath_x^{-1})}_-{\cong}
        &\phi(x)
        .
      }
    \end{equation}
    In particular, the component of $\tahar{\phi}$ at $\kappa$ is
    \begin{equation}
      \label{eq2 Comm(T) Comm(T)laxK(K,T) adjunction}
      \tahar{\phi}_{\kappa}:
      \xymatrix@C=40pt{
        \mon{t}(\kappa)\tensor \phi(\kappa)
        \ar[r]^-{\mon{t}_{\kappa}^{-1}\tensor I_{\phi(\kappa)}}_-{\cong}
        &\unit\tensor \phi(\kappa)
        \ar[r]^-{\imath_{\phi(\kappa)}^{-1}}_-{\cong}
        &\phi(\kappa).
      }
    \end{equation}
  \end{itemize}
\end{proposition}
\begin{proof}
  For each object $B$ in $\cat{Comm}(\kos{T})$,
  the lax $\kos{K}$-tensor functor
  $\iota(B)=(\text{$\mon{t}$}\tensor B,\what{\mon{t}\tensor B}):
  (\kos{K},\id_{\kos{K}})\to(\kos{T},\mon{t})$
  is well-defined.
  For each morphism $f:B\to B^{\pr}$ in $\cat{Comm}(\kos{T})$,
  one can check that the monoidal natural transformation
  $\mon{t}\tensor f:\mon{t}\tensor B\Rightarrow \mon{t}\tensor B^{\pr}$
  becomes a monoidal $\kos{K}$-tensor natural transformation
  $\iota(f):\iota(B)\Rightarrow\iota(B^{\pr})$,
  by verifying the relation
  $\what{\mon{t}\tensor B}:
  \xymatrix@C=15pt{
  \mon{t}\tensor B
  \ar@2{->}[r]^-{\mon{t}\tensor f}
  &\mon{t}\tensor B^{\pr}
  \ar@2{->}[r]^-{\what{\mon{t}\tensor B^{\pr}}}
  &\mon{t}
  }$.
  Thus the left adjoint functor $\iota$ is well-defined.
  It is immediate that the right adjoint functor $\CR$ is also well-defined.
  Furthermore,
  the isomorphism
  $\xymatrix@C=20pt{
    B
    \ar[r]^-{\imath_B}_-{\cong}
    &\unit\tensor B
    \ar[r]^-{\mon{t}_{\kappa}\tensor I_B}_-{\cong}
    &\mon{t}(\kappa)\tensor B
  }$
  in $\cat{Comm}(\kos{T})$
  is natural in variable $B\in\obj{\cat{Comm}(\kos{T})}$.
  Thus the adjunction unit is also well-defined.

  We show that the adjunction counit is well-defined in two steps.
  First step is to show that 
  for each lax $\kos{K}$-tensor functor
  $(\phi,\what{\phi}):(\kos{K},\id_{\kos{K}})\to (\kos{T},\mon{t})$,
  we have the monoidal $\kos{K}$-tensor natural transformation
  $\tahar{\phi}:
  (\text{$\mon{t}$}\tensor\phi(\kappa),\what{\mon{t}\tensor\phi(\kappa)})
  \Rightarrow(\phi,\what{\phi})$
  as we claimed in (\ref{eq Comm(T) Comm(T)laxK(K,T) adjunction}).
  We leave for the readers to check that
  $\tahar{\phi}:
  \text{$\mon{t}$}\tensor\phi(\kappa)
  \Rightarrow \phi:\kos{K}\to\kos{T}$
  is a monoidal natural transformation.
  We can check the relation
  \begin{equation*}
    \vcenter{\hbox{
      \xymatrix{
        \mon{t}\tensor \phi(\kappa)
        \ar@2{->}[rr]^-{\tahar{\phi}}
        &\text{ }
        &\phi
        \\
        \text{ }
        &\mon{t}
        \ar@2{->}[ul]^-{\what{\mon{t}\tensor\phi(\kappa)}}
        \ar@2{->}[ur]_-{\what{\phi}}
        &\text{ }
      }
    }}
  \end{equation*}
  as follows.
  Let $x$ be an object in $\CK$.
  \begin{equation*}
    \xymatrix@C=70pt{
      \mon{t}(x)
      \ar[dd]_-{\what{\mon{t}\tensor \phi(\kappa)}_x}
      \ar@{=}[r]
      &\mon{t}(x)
      \ar[d]^-{\jmath_{\mon{t}(x)}}_-{\cong}
      \ar@{=}[r]
      &\mon{t}(x)
      \ar[d]^-{\what{\phi}_x}
      \\
      \text{ }
      &\mon{t}(x)\tensor \unit
      \ar@/_0.5pc/[dl]|-{I_{\mon{t}(x)}\tensor \phi_{\kappa}}
      \ar[d]^-{\what{\phi}_x\tensor I_{\unit}}
      &\phi(x)
      \ar@{=}[dddd]
      \ar@/^0.5pc/[dl]^-{\jmath_{\phi(x)}}_-{\cong}
      \\
      \mon{t}(x)\tensor \phi(\kappa)
      \ar[ddd]_-{\tahar{\phi}_x}
      \ar@/_0.5pc/[dr]|-{\what{\phi}_x\tensor I_{\phi(\kappa)}}
      &\phi(x)\tensor \unit
      \ar[d]^-{I_{\phi(x)}\tensor \phi_{\kappa}}
      &\text{ }
      \\
      \text{ }
      &\phi(x)\tensor \phi(\kappa)
      \ar[d]^-{\phi_{x,\kappa}}
      &\text{ }
      \\
      \text{ }
      &\phi(x\otimes \kappa)
      \ar[d]^-{\phi(\jmath_x^{-1})}_-{\cong}
      &\text{ }
      \\
      \phi(x)
      \ar@{=}[r]
      &\phi(x)
      \ar@{=}[r]
      &\phi(x)
    }
  \end{equation*}
  Thus the monoidal $\kos{K}$-tensor natural transformation
  $\tahar{\phi}:
  (\text{$\mon{t}$}\tensor\phi(\kappa),\what{\mon{t}\tensor\phi(\kappa)})
  \Rightarrow(\phi,\what{\phi})$
  is well-defined.
  Next step is to show that
  $\tahar{\phi}$
  is natural in variable $(\phi,\what{\phi})$.
  Let
  $(\psi,\what{\psi}):(\kos{K},\id_{\kos{K}})\to (\kos{T},\mon{t})$
  be another lax $\kos{K}$-tensor functor
  and let
  $\vartheta:(\phi,\what{\phi})\Rightarrow (\psi,\what{\psi})$
  be a monoidal $\kos{K}$-tensor natural transformation.
  Recall that we have $\what{\psi}=\vartheta\circ \what{\phi}
  :\mon{t}\Rightarrow \psi$.
  We need to verify the relation
  \begin{equation*}
    \vcenter{\hbox{
      \xymatrix@R=30pt@C=40pt{
        \mon{t}\tensor \phi(\kappa)
        \ar@2{->}[d]_-{\mon{t}\tensor \vartheta_{\kappa}}
        \ar@2{->}[r]^-{\tahar{\phi}}
        &\phi
        \ar@2{->}[d]^-{\vartheta}
        \\
        \mon{t}\tensor \psi(\kappa)
        \ar@2{->}[r]^-{\tahar{\psi}}
        &\psi
      }
    }}
    \qquad\quad
    \tahar{\psi}\circ (\mon{t}\tensor \vartheta_{\kappa})
    =
    \vartheta\circ \tahar{\phi}:
    \mon{t}\tensor \phi(\kappa)\Rightarrow\psi
    .
  \end{equation*}
  We can check this as follows.
  Let $x$ be an object in $\CK$.
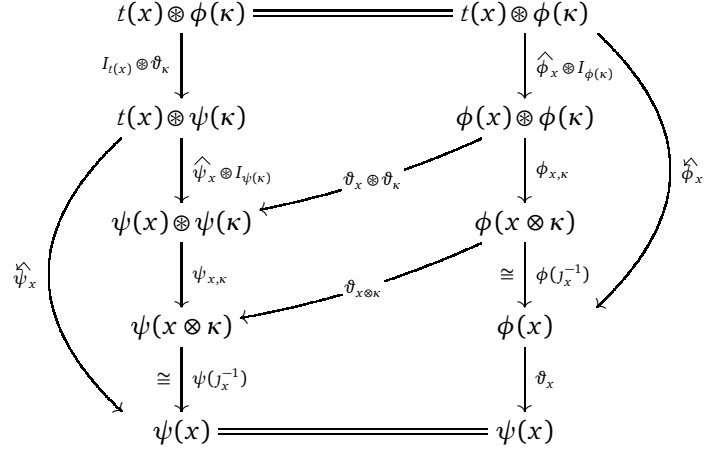
  \begin{equation*}
    \vcenter{\hbox{
      \xymatrix@C=70pt{
        \mon{t}(x)\tensor \phi(\kappa)
        \ar[d]_-{I_{\mon{t}(x)}\tensor \vartheta_{\kappa}}
        \ar@{=}[r]
        &\mon{t}(x)\tensor \phi(\kappa)
        \ar[d]^-{\what{\phi}_x\tensor I_{\phi(\kappa)}}
        \ar@/^3pc/@<4ex>[ddd]^-{\tahar{\phi}_x}
        \\
        \mon{t}(x)\tensor \psi(\kappa)
        \ar[d]^-{\what{\psi}_x\tensor I_{\psi(\kappa)}}
        \ar@/_3pc/@<-3ex>[ddd]_-{\tahar{\psi}_x}
        &\phi(x)\tensor \phi(\kappa)
        \ar[d]^-{\phi_{x,\kappa}}
        \ar@/^0.5pc/[dl]|-{\vartheta_x\tensor \vartheta_{\kappa}}
        \\
        \psi(x)\tensor \psi(\kappa)
        \ar[d]^-{\psi_{x,\kappa}}
        &\phi(x\otimes \kappa)
        \ar[d]^-{\phi(\jmath_x^{-1})}_-{\cong}
        \ar@/^0.5pc/[dl]|-{\vartheta_{x\otimes \kappa}}
        \\
        \psi(x\otimes \kappa)
        \ar[d]^-{\psi(\jmath_x^{-1})}_-{\cong}
        &\phi(x)
        \ar[d]^-{\vartheta_x}
        \\
        \psi(x)
        \ar@{=}[r]
        &\psi(x)
      }
    }}
  \end{equation*}
  This shows that the adjunction counit is well-defined.

  We are left to show that the adjunction unit, counit
  satisfy the triangle identities.
  Let 
  $(\phi,\what{\phi}):(\kos{K},\id_{\kos{K}})\to (\kos{T},\mon{t})$
  be a lax $\kos{K}$-tensor functor.
  we can easily check one of the traingle identities 
  \begin{equation*}
    I_{\phi(\kappa)}:
    \xymatrix@C=25pt{
      \phi(\kappa)
      \ar[r]^-{\imath_{\phi(\kappa)}}_-{\cong}
      &\unit\tensor \phi(\kappa)
      \ar[r]^-{\mon{t}_{\kappa}\tensor I_{\phi(\kappa)}}_-{\cong}
      &\mon{t}(\kappa)\otimes \phi(\kappa)
      \ar[r]^-{\tahar{\phi}_{\kappa}}
      &\phi(\kappa)
    }
  \end{equation*}
  which is the description
  (\ref{eq2 Comm(T) Comm(T)laxK(K,T) adjunction})
  of $\hatar{\phi}_{\kappa}$ that we claimed.
  Let $B$ be an object in $\cat{Comm}(\kos{T})$.
  For each object $x$ in $\CK$,
  we can rewrite $\tahar{\mon{t}\tensor B}_x$ as
  \begin{equation*}
    \tahar{\mon{t}\tensor B}_x:
    \xymatrix{
      \mon{t}(x)\tensor (\mon{t}(\kappa)\tensor B)
      \ar[r]^-{a_{\mon{t}(x),\mon{t}(\kappa),B}}_-{\cong}
      &(\mon{t}(x)\tensor \mon{t}(\kappa))\tensor B
      \ar[r]^-{\mon{t}_{x,\kappa}\tensor I_B}_-{\cong}
      &\mon{t}(x\otimes \kappa)\tensor B
      \ar[r]^-{\mon{t}(\jmath_x^{-1})\tensor I_B}_-{\cong}
      &\mon{t}(x)\tensor B
    }
  \end{equation*}
  and using the above relation,
  we verify the other triangle identity
  \begin{equation*}
      I_{\mon{t}\tensor B}:
      \xymatrix@C=30pt{
        \mon{t}\tensor B
        \ar@2{->}[r]^-{\mon{t}\tensor \imath_B}_-{\cong}
        &\mon{t}\tensor (\unit\tensor B)
        \ar@2{->}[r]^-{\mon{t}\tensor (\mon{t}_{\kappa}\tensor I_B)}_-{\cong}
        &\mon{t}\tensor (\mon{t}(\kappa)\tensor B)
        \ar[r]^-{\tahar{\mon{t}\tensor B}}
        &\mon{t}\tensor B
      }
  \end{equation*}
  as follows. Let $x$ be an object in $\CK$.
  \begin{equation*}
    \vcenter{\hbox{
      \xymatrix@C=60pt{
        \mon{t}(x)\tensor B
        \ar[d]_-{I_{\mon{t}(x)}\tensor \imath_B}^-{\cong}
        \ar@{=}[r]
        &\mon{t}(x)\tensor B
        \ar[dd]^-{\jmath_{\mon{t}(x)}\tensor I_B}_-{\cong}
        \ar@/^3pc/@<5ex>@{=}[ddddd]
        \\
        \mon{t}(x)\tensor (\unit\tensor B)
        \ar[d]_-{I_{\mon{t}(x)}\tensor (\mon{t}_{\kappa}\tensor I_B)}
        \ar@/^0.5pc/[dr]^-{a_{\mon{t}(x),\unit,B}}_-{\cong}
        &\text{ }
        \\
        \mon{t}(x)\tensor (\mon{t}(\kappa)\tensor B)
        \ar@/_3pc/@<-4ex>[ddd]_-{\tahar{\mon{t}\tensor B}_x}
        \ar[d]^-{a_{\mon{t}(x),\mon{t}(\kappa),B}}_-{\cong}
        &(\mon{t}(x)\tensor \unit)\tensor B
        \ar@/^0.5pc/[dl]|-{(I_{\mon{t}(x)}\tensor \mon{t}_{\kappa})\tensor I_B}
        \ar[ddd]^-{\jmath^{-1}_{\mon{t}(x)}\tensor I_B}
        \\
        (\mon{t}(x)\tensor \mon{t}(\kappa))\tensor B
        \ar[d]^-{\mon{t}_{x,\kappa}\tensor I_B}_-{\cong}
        &\text{ }
        \\
        \mon{t}(x\otimes \kappa)\tensor B
        \ar[d]^-{\mon{t}(\jmath_x^{-1})\tensor I_B}_-{\cong}
        &\text{ }
        \\
        \mon{t}(x)\tensor B
        \ar@{=}[r]
        &\mon{t}(x)\tensor B
      }
    }}
  \end{equation*}
  This completes the proof of Proposition~\ref{prop Comm(T) Comm(T)laxK(K,T) adjunction}.
\qed\end{proof}

\begin{definition}
  Let $(\kos{T},\mon{t})$ be a strong $\kos{K}$-tensor category
  and let 
  $(\phi,\what{\phi}):(\kos{K},\id_{\kos{K}})\to (\kos{T},\mon{t})$
  be a lax $\kos{K}$-tensor functor.
  We define the \emph{coreflection} of $(\phi,\what{\phi})$
  as the monoidal $\kos{K}$-tensor natural transformation
  \begin{equation*}
    \tahar{\phi}
    :
    (\text{$\mon{t}$}\tensor\phi(\kappa),\what{\mon{t}\tensor\phi(\kappa)})
    \Rightarrow
    (\phi,\what{\phi})
    :(\kos{K},\id_{\kos{K}})\to (\kos{T},\mon{t})
  \end{equation*}
  introduced in Proposition~\ref{prop Comm(T) Comm(T)laxK(K,T) adjunction}.
  We say $(\phi,\what{\phi})$ is \emph{coreflective}
  if its coreflection $\tahar{\phi}$
  is a monoidal $\kos{K}$-tensor natural isomorphism.
\end{definition}

\begin{corollary} \label{cor Comm(T) strongKcat coreflective adjunction}
  Let $(\kos{T},\mon{t})$ be a strong $\kos{K}$-tensor category
  and recall the coreflective adjunction
  $\iota\dashv \CR$ in Proposition~\ref{prop Comm(T) Comm(T)laxK(K,T) adjunction}.
  The essential image of the left adjoint $\iota$
  is the coreflective full subcategory
  \begin{equation*}
    \mathbb{SMC}_{\cat{lax}}^{\kos{K}\backslash\!\!\backslash}\bigl((\kos{K},\id_{\kos{K}}),(\kos{T},\mon{t})\bigr)_{\cat{crfl}}
  \end{equation*}
  of coreflective lax $\kos{K}$-tensor functors
  from $(\kos{K},\id_{\kos{K}})$ to $(\kos{T},\mon{t})$,
  and the adjunction $\iota\dashv \CR$ restricts to an adjoint equivalence of categories
  \begin{equation*}
    \iota:
    \cat{Comm}(\kos{T})
    \simeq
    \mathbb{SMC}_{\cat{lax}}^{\kos{K}\backslash\!\!\backslash}\bigl((\kos{K},\id_{\kos{K}}),(\kos{T},\mon{t})\bigr)_{\cat{crfl}}
    :\CR
    .
  \end{equation*}
\end{corollary}

Let
$(\phi,\what{\phi}):(\kos{K},\id_{\kos{K}})\to (\kos{T},\mon{t})$
be a lax $\kos{K}$-tensor functor
to a strong $\kos{K}$-tensor category $(\kos{T},\mon{t})$.
From the definition 
(\ref{eq Comm(T) Comm(T)laxK(K,T) adjunction}) of the coreflection $\tahar{\phi}$,
we see that 
$(\phi,\what{\phi})$
is coreflective
if the associated $\kos{K}$-equivariance
$\cevar{\phi}$ is a natural isomorphism.
The following lemma states that the converse is also true.

\begin{lemma} \label{lem Comm(T) equivariance and coreflection}
  Let $(\kos{T},\mon{t})$ be a strong $\kos{K}$-tensor category
  and let $(\phi,\what{\phi}):(\kos{K},\id_{\kos{K}})\to (\kos{T},\mon{t})$
  be a lax $\kos{K}$-tensor functor.
  Then
  $(\phi,\what{\phi})$ is coreflective
  if and only if
  the associated $\kos{K}$-equivariance $\cevar{\phi}$
  is a natural isomorphism.
\end{lemma}
\begin{proof}
  We explained the if part.
  To verify the only if part,
  it suffices to show that the following diagram
  commutes for all $x$, $y\in\obj{\CK}$.
  \begin{equation*}
    \vcenter{\hbox{
      \xymatrix@C=30pt{
        \mon{t}(x)\tensor (\mon{t}(y)\tensor \phi(\kappa))
        \ar[d]_-{a_{\mon{t}(x),\mon{t}(y),\phi(\kappa)}}^-{\cong}
        \ar[rr]^-{I_{\mon{t}(x)}\tensor \tahar{\phi}_y}
        &\text{ }
        &\mon{t}(x)\tensor \phi(y)
        \ar[d]^-{\cevar{\phi}_{x,y}}
        \\
        (\mon{t}(x)\tensor \mon{t}(y))\tensor \phi(\kappa)
        \ar[r]^-{\mon{t}_{x,y}\tensor I_{\phi(\kappa)}}_-{\cong}
        &\mon{t}(x\otimes y)\tensor \phi(\kappa)
        \ar[r]^-{\tahar{\phi}_{x\otimes y}}
        &\phi(x\otimes y)
      }
    }}
  \end{equation*}
  We can check this as follows.
  \begin{equation*}
    \vcenter{\hbox{
      \xymatrix@C=3pt{
        \mon{t}(x)\tensor (\mon{t}(y)\tensor \phi(\kappa))
        \ar[d]_-{a_{\mon{t}(x),\mon{t}(y),\phi(\kappa)}}^-{\cong}
        \ar@{=}[rr]
        &\text{ }
        &\mon{t}(x)\tensor (\mon{t}(y)\tensor \phi(\kappa))
        \ar[d]_-{\what{\phi}_x\tensor (\what{\phi}_y\tensor I_{\phi(\kappa)})}
        \ar@/^2pc/[dddr]|(0.7){I_{\mon{t}(x)}\tensor \tahar{\phi}_y}
        &\text{ }
        \\
        (\mon{t}(x)\tensor \mon{t}(y))\tensor \phi(\kappa)
        \ar[d]_-{\mon{t}_{x,y}\tensor I_{\phi(\kappa)}}^-{\cong}
        \ar@/^0.5pc/[dr]|-{(\what{\phi}_x\tensor \what{\phi}_y)\tensor I_{\phi(\kappa)}}
        &\text{ }
        &\phi(x)\tensor (\phi(y)\tensor \phi(\kappa))
        \ar@/_0.5pc/[dl]_-{a_{\phi(x),\phi(y),\phi(\kappa)}}^-{\cong}
        \ar[d]^-{I_{\phi(\kappa)}\tensor \phi_{y,\kappa}}
        &\text{ }
        \\
        \mon{t}(x\otimes y)\tensor \phi(\kappa)
        \ar@/_2pc/@<-2ex>@{.>}[ddd]|-{\tahar{\phi}_{x\otimes y}}
        \ar[d]|-{\what{\phi}_{x\otimes y}\tensor I_{\phi(\kappa)}}
        &(\phi(x)\tensor \phi(y))\tensor \phi(\kappa)
        \ar@/^0.5pc/[dl]|-{\phi_{x,y}\tensor I_{\phi(\kappa)}}
        &\phi(x)\tensor \phi(y\otimes \kappa)
        \ar@/^0.5pc/[dl]|-{\phi_{x,y\otimes \kappa}}
        \ar[dd]_(0.6){I_{\phi(x)}\tensor \phi(\jmath_y^{-1})}^(0.6){\cong}
        &\text{ }
        \\
        \phi(x\otimes y)\tensor \phi(\kappa)
        \ar[d]^-{\phi_{x\otimes y,\kappa}}
        &\phi(x\otimes (y\otimes \kappa))
        \ar@/^0.5pc/[dl]^-{a_{x,y,\kappa}}_-{\cong}
        \ar[dd]^-{\phi(I_x\otimes \jmath_y^{-1})}_-{\cong}
        &\text{ }
        &\mon{t}(x)\!\tensor\! \phi(y)
        \ar[dd]^-{\cevar{\phi}_{x,y}}
        \ar@/^0.5pc/[dl]|-{\what{\phi}_x\tensor I_{\phi(y)}}
        \\
        \phi((x\otimes y)\otimes \kappa)
        \ar[d]^-{\phi(\jmath_{x\otimes y}^{-1})}_-{\cong}
        &\text{ }
        &\phi(x)\tensor \phi(y)
        \ar[d]^-{\phi_{x,y}}
        &\text{ }
        \\
        \phi(x\otimes y)
        \ar@{=}[r]
        &\phi(x\otimes y)
        \ar@{=}[r]
        &\phi(x\otimes y)
        \ar@{=}[r]
        &\phi(x\otimes y)
      }
    }}
  \end{equation*}
  This completes the proof of Lemma~\ref{lem Comm(T) equivariance and coreflection}.
\qed\end{proof}

For the rest of this subsection,
we focus on the case
$(\kos{T},\mon{t})=(\kos{K},\id_{\kos{K}})$.
Denote
\begin{equation}\label{eq Comm(T) Comm(T)def}
  \kos{C\!o\!m\!m}(\kos{K})
  =(\cat{Comm}(\kos{K}),\otimes,\kappa)
\end{equation}
as the cocartesian monoidal category
of commutative monoids in $\kos{K}$.
An object in $\cat{Comm}(\kos{K})$
is a triple
$(b,\pc_b,u_b)$ 
which we often abbreviate as $b$.
The object $\kappa$ equipped with
$\pc_{\kappa}=\imath^{-1}_{\kappa}=\jmath^{-1}_{\kappa}
:\kappa\otimes \kappa\xrightarrow{\cong}\kappa$
and 
$u_{\kappa}=I_{\kappa}:\kappa\xrightarrow{\cong}\kappa$
is an initial object in $\cat{Comm}(\kos{K})$.
The coproduct of objects $b$, $b^{\pr}$ in $\cat{Comm}(\kos{K})$
is $b\otimes b^{\pr}$ whose product, unit morphisms
$\pc_{b\otimes b^{\pr}}$, $u_{b\otimes b^{\pr}}$ are
defined as in (\ref{eq Comm(T) BB' definition}).
We also denote 
\begin{equation}\label{eq Comm(T) Aff(T)def}
  \kos{A\!f\!f}(\kos{K})
  =
  (\cat{Aff}(\kos{K}),\times,\cat{Spec}(\kappa))
\end{equation}
as the opposite cartesian category of
$\kos{C\!o\!m\!m}(\kos{K})$.
For each object $b$ in $\cat{Comm}(\kos{K})$,
we denote $\cat{Spec}(b)$ as
the corresponding object in $\cat{Aff}(\kos{K})$.

Let $b$ be an object in $\cat{Comm}(\kos{K})$.
We have a lax symmetric monoidal monad 
\begin{equation*}
  \langle \otimes b\rangle=(\otimes b,\mu^{\otimes b},\upsilon^{\otimes b})
\end{equation*}
on $\kos{K}$.  
The underlying lax symmetric monoidal endofunctor
$\otimes b:\kos{K}\to \kos{K}$
is the endofunctor
$\otimes b= \slot\otimes b:\CK\to \CK$
equipped with monoidal coherence morphisms
\begin{equation*}
  \xymatrix{
    (x\otimes b)\otimes (y\otimes b)
    \ar[r]^-{(\otimes b)_{x,y}}
    &(x\otimes y)\otimes b
  }
  \qquad
  x,y\in\obj{\CK}
  \qquad
  \xymatrix{
    \kappa
    \ar[r]^-{(\otimes b)_{\kappa}}
    &\kappa\otimes b
  }
\end{equation*}
defined as in (\ref{eq Comm(T) tensorB definition}),
and the monoidal natural transformations
\begin{equation*}
  \xymatrix{
    (x\otimes b)\otimes b
    \ar[r]^-{\mu^{\otimes b}_x}
    &x\otimes b
  }
  \qquad
  \xymatrix{
    x
    \ar[r]^-{\upsilon^{\otimes b}_x}
    &x\otimes b
  }
  \qquad
  x\in\obj{\CK}
\end{equation*}
are defined as in (\ref{eq Comm(T) muBupsilonB definition}).

Let $(\phi,\what{\phi})$,
$(\psi,\what{\psi}):(\kos{K},\id_{\kos{K}})\to (\kos{K},\id_{\kos{K}})$
be coreflective lax $\kos{K}$-tensor endofunctors
on $(\kos{K},\id_{\kos{K}})$.
By Lemma~\ref{lem Comm(T) equivariance and coreflection},
the associated $\kos{K}$-equivariances $\cevar{\phi}$, $\cevar{\psi}$
are natural isomorphisms.
By Lemma~\ref{lem LaxKtensorCat compositionKequiv},
the $\kos{K}$-equivariance
$\cevar{\psi\phi}$
associated to the composition $(\psi\phi,\what{\psi\phi})$
is also a natural isomorphism.
In particular,
the composition $(\psi\phi,\what{\psi\phi})$ is also coreflective.
Thus we denote
\begin{equation}
  \label{eq Comm(T) laxKTend(K,idK) def}
  \kos{E\!n\!\!d}_{\cat{lax}}^{\kos{K}\backslash\!\!\backslash}(\kos{K},\id_{\kos{K}})_{\cat{crfl}}
  =
  \bigl(
    \mathbb{SMC}_{\cat{lax}}^{\kos{K}\backslash\!\!\backslash}\bigl((\kos{K},\id_{\kos{K}}),(\kos{K},\id_{\kos{K}})\bigr)_{\cat{crfl}}
    ,
    \circ
    ,
    (\id_{\kos{K}},I)
  \bigr)
\end{equation}
as the monoidal category of
coreflective lax $\kos{K}$-tensor endofunctors
on $(\kos{K},\id_{\kos{K}})$.
The monoidal product is given by composition $\circ$
and the unit object
$(\id_{\kos{K}},I)$
is the identity strong $\kos{K}$-tensor endofunctor
of $(\kos{K},\id_{\kos{K}})$.
We also denote
\begin{equation*}
  \kos{E\!n\!\!d}_{\cat{lax}}^{\kos{K}\backslash\!\!\backslash}(\kos{K},\id_{\kos{K}})_{\cat{crfl}}^{\cat{rev}}
\end{equation*}
as the reversed monoidal category of
$\kos{E\!n\!\!d}_{\cat{lax}}^{\kos{K}\backslash\!\!\backslash}(\kos{K},\id_{\kos{K}})_{\cat{crfl}}$,
whose reversed monoidal product $\circ^{\cat{rev}}$
is given by $(\phi,\what{\phi})\circ^{\cat{rev}}(\psi,\what{\psi})
=(\psi\phi,\what{\psi\phi})$.

\begin{proposition} \label{prop Comm(T) Comm(K)=KtensorEnd(K)}
  The adjoint equivalence of categories
  \begin{equation*}
    \iota:
    \cat{Comm}(\kos{K})
    \simeq
    \mathbb{SMC}_{\cat{lax}}^{\kos{K}\backslash\!\!\backslash}\bigl((\kos{K},\id_{\kos{K}}),(\kos{K},\id_{\kos{K}})\bigr)_{\cat{crfl}}
    :\CR
  \end{equation*}
  described in Corollary~\ref{cor Comm(T) strongKcat coreflective adjunction}
  becomes an adjoint equivalence of monoidal categories
  \begin{equation*}
    \iota:
    \kos{C\!o\!m\!m}(\kos{K})
    \simeq
    \kos{E\!n\!\!d}_{\cat{lax}}^{\kos{K}\backslash\!\!\backslash}(\kos{K},\id_{\kos{K}})_{\cat{crfl}}^{\cat{rev}}
    :\CR
    .
  \end{equation*}
  Let $b$, $b^{\pr}$ be objects in $\cat{Comm}(\kos{K})$
  and let
  $(\phi,\what{\phi})$,
  $(\psi,\what{\psi})$ 
  be coreflective lax $\kos{K}$-tensor endofunctors on $(\kos{K},\id_{\kos{K}})$.
  \begin{itemize}
    \item 
    The monoidal coherence isomorphisms of the left adjoint $\iota$
    are the following monoidal $\kos{K}$-tensor natural isomorphisms.
    Let $x\in\obj{\CK}$.
    \begin{equation*}
      \begin{aligned}
        a_{\slot,b,b^{\pr}}
        &:\!
        \xymatrix@C=15pt{
          \iota(b\otimes b^{\pr})
          \ar@2{->}[r]^-{\cong}
          &\iota(b)\circ^{\cat{rev}}\iota(b^{\pr})
          ,
        }
        &\qquad
        \jmath^{-1}_{\slot}
        &:\!
        \xymatrix@C=15pt{
          \iota(\kappa)
          \ar@2{->}[r]^-{\cong}
          &(\id_{\kos{K}},I)
          ,
        }
        \\
        a_{x,b,b^{\pr}}
        &:\!
        \xymatrix@C=15pt{
          x\otimes (b\otimes b^{\pr})
          \ar[r]^-{\cong}
          &(x\otimes b)\otimes b^{\pr}
          ,
        }
        &\qquad
        \jmath^{-1}_x
        &:\!
        \xymatrix@C=15pt{
          x\otimes \kappa
          \ar[r]^-{\cong}
          &x
          .
        }
      \end{aligned}
    \end{equation*}

    \item 
    The monoidal coherence isomorphisms of the right adjoint $\CR$
    are the following isomorphisms in $\cat{Comm}(\kos{K})$.
    \begin{equation*}
      \vcenter{\hbox{
        \xymatrix{
          \CR(\phi,\what{\phi})\otimes \CR(\psi,\what{\psi})
          \ar@{=}[d]
          \ar[r]^-{\cong}
          &\CR(\psi\phi,\what{\psi\phi})
          \ar@{=}[d]
          \\
          \phi(\kappa)\otimes \psi(\kappa)
          \ar[r]^-{\tahar{\psi}_{\phi(\kappa)}}_-{\cong}
          &\psi\phi(\kappa)
        }
      }}
      \qquad\qquad
      \kappa=\CR(\id_{\kos{K}},I)
    \end{equation*}
  \end{itemize}
  In particular, we have the following relation
  for each object $x$ in $\CK$.
  \begin{equation}\label{eq Comm(T) Comm(K)=KtensorEnd(K)}
    \vcenter{\hbox{
      \xymatrix@C=60pt{
        x\otimes (\phi(\kappa)\otimes \psi(\kappa))
        \ar[r]^-{a_{x,\phi(\kappa),\psi(\kappa)}}_-{\cong}
        \ar[d]_-{I_x\otimes \tahar{\psi}_{\phi(\kappa)}}^-{\cong}
        &(x\otimes \phi(\kappa))\otimes \psi(\kappa)
        \ar[d]^-{(\tahar{\psi}\tahar{\phi})_x}_-{\cong}
        \\
        x\otimes \psi\phi(\kappa)
        \ar[r]^-{\tahar{\psi\phi}_x}_-{\cong}
        &\psi\phi(x)
      }
    }}
  \end{equation}
\end{proposition}
\begin{proof}
  We leave for the readers to check that
  $a_{\slot,b,b^{\pr}}$
  and
  $\jmath^{-1}_{\slot}$
  are monoidal $\kos{K}$-tensor natural transformations,
  and that $a_{\slot,b,b^{\pr}}$
  is natural in variables $b$, $b^{\pr}$.
  To conclude that the left adjoint $\iota$
  is a strong monoidal functor
  with monoidal coherence isomorphisms
  $a_{\slot,b,b^{\pr}}$,
  $\jmath^{-1}_{\slot}$
  we need to check the monoidal coherence relations.
  These relations are obviously true
  due to the Mac Lane coherence theorem.
  Thue the left adjoint $\iota$
  is a strong monoidal functor as we claimed.
  The right adjoint $\CR$ has a unique strong monoidal functor structure
  such that the adjoint equivalence of categories $\iota\dashv\CR$
  becomes an equivalence of monoidal categories.
  The monoidal coherence isomorphisms of $\CR$ are given as follows.
  \begin{equation*}
    \vcenter{\hbox{
      \xymatrix@C=5pt{
        \phi(\kappa)\otimes \psi(\kappa)
        \ar[d]_-{\imath_{\phi(\kappa)\otimes \psi(\kappa)}}^-{\cong}
        \ar@{=}[rr]
        &\text{ }
        &\phi(\kappa)\otimes \psi(\kappa)
        \ar@/^1pc/[ddll]_-{\imath_{\phi(\kappa)}\otimes I_{\psi(\kappa)}}^-{\cong}
        \ar[dd]^-{\tahar{\psi}_{\phi(\kappa)}}_-{\cong}
        \\
        \kappa\otimes (\phi(\kappa)\otimes \psi(\kappa))
        \ar[d]_-{a_{\kappa,\phi(\kappa),\psi(\kappa)}}^-{\cong}
        &\text{ }
        &\text{ }
        \\
        (\kappa\otimes \phi(\kappa))\otimes \psi(\kappa)
        \ar[dd]_-{(\tahar{\psi}\tahar{\phi})_{\kappa}}^-{\cong}
        \ar@/^0.5pc/[dr]_-{\tahar{\psi}_{\kappa\otimes\phi(\kappa)}}^{\cong}
        &\text{ }
        &\psi\phi(\kappa)
        \ar@{=}[dd]
        \ar@/^0.5pc/[dl]_-{\psi(\imath_{\phi(\kappa)})}^-{\cong}
        \\
        \text{ }
        &\psi(\kappa\otimes\phi(\kappa))
        \ar@/^0.5pc/[dr]_-{\psi(\tahar{\phi}_{\kappa})}^-{\cong}
        \\
        \psi\phi(\kappa)
        \ar@{=}[rr]
        &\text{ }
        &\psi\phi(\kappa)
      }
    }}
    \vcenter{\hbox{
      \xymatrix@C=5pt{
        \kappa
        \ar[d]_-{\imath_{\kappa}}^-{\cong}
        \ar@{=}[r]
        &\kappa
        \ar[d]^-{\imath_{\kappa}}_-{\cong}
        \ar@<2ex>@/^2pc/@{=}[dd]
        \\
        \CR\iota(\kappa)
        \ar[d]_-{\CR(\jmath^{-1}_{\slot})}^-{\cong}
        \ar@{=}[r]
        &\kappa\otimes\kappa
        \ar[d]^-{\jmath_{\kappa}^{-1}}_-{\cong}
        \\
        \CR(\id_{\kos{K}},I)
        \ar@{=}[r]
        &\kappa
      }
    }}
  \end{equation*}
  This shows that we have the adjoint equivalence $\iota\dashv \CR$
  of monoidal categories as we claimed.
  In particular, we obtain the relation (\ref{eq Comm(T) Comm(K)=KtensorEnd(K)})
  which is a consequence of the fact that 
  $\iota\dashv \CR$
  is an adjoint equivalence of monoidal categories.
  \begin{equation*}
    \vcenter{\hbox{
      \xymatrix@C=30pt{
        \iota(\CR(\phi,\what{\phi})\otimes \CR(\psi,\what{\psi}))
        \ar@2{->}[d]_-{\iota(\tahar{\psi}_{\phi(\kappa)})}^-{\cong}
        \ar@2{->}[rr]^-{a_{\slot,\CR(\phi,\what{\phi}),\CR(\psi,\what{\psi})}}_-{\cong}
        &\text{ }
        &\iota\CR(\phi,\what{\phi})\circ^{\cat{rev}}\iota\CR(\psi,\what{\psi})
        \ar@2{->}[d]^-{\tahar{\psi}\tahar{\phi}}
        \\
        \iota\CR(\psi\phi,\what{\psi\phi})
        \ar@2{->}[r]^-{\tahar{\psi\phi}}_-{\cong}
        &(\psi\phi,\what{\psi,\phi})
        \ar@{=}[r]
        &(\phi,\what{\phi})\circ^{\cat{rev}}(\psi,\what{\psi})
      }
    }}
  \end{equation*}
  This completes the proof of Proposition~\ref{prop Comm(T) Comm(K)=KtensorEnd(K)}.
\qed\end{proof}

Since $\kos{C\!o\!m\!m}(\kos{K})$ is a cocartesian monoidal category,
each object $b$ in $\cat{Comm}(\kos{K})$
has a unique structure of a monoid in $\kos{C\!o\!m\!m}(\kos{K})$
with product $\pc_b:b\otimes b\to b$
and unit $u_b:\kappa\to b$.
Recall the adjoint equivalence $\iota\dashv \CR$ of monoidal categories
in Proposition~\ref{prop Comm(T) Comm(K)=KtensorEnd(K)}.
The left adjoint $\iota$
sends the monoid $(b,\pc_b,u_b)$ in $\kos{C\!o\!m\!m}(\kos{K})$
to the lax $\kos{K}$-tensor monad
$\langle \otimes b, \what{\otimes b}\rangle
=
(\otimes b,\what{\otimes b},\mu^{\otimes b},\upsilon^{\otimes b})$
on $(\kos{K},\id_{\kos{K}})$
where
\begin{equation*}
  \vcenter{\hbox{
    \xymatrix@C=15pt{
      \text{ }
      &\text{ }
      &\kos{K}
      \ar[d]^-{\otimes b}
      \\
      \kos{K}
      \ar@/^0.7pc/[urr]^-{\id_{\kos{K}}}
      \ar[rr]_-{\id_{\kos{K}}}
      &\text{ }
      &\kos{K}
      \xtwocell[l]{}<>{<2.5>{\what{\otimes b}\quad}}
    }
  }}
  \qquad
  \what{\otimes b}_x=\upsilon^{\otimes b}_x:
  \xymatrix@C=18pt{
    x
    \ar[r]^-{\jmath_x}_-{\cong}
    &x\otimes \kappa
    \ar[r]^-{I_x\otimes u_b}
    &x\otimes b
    ,
  }
  \quad
  x\in \obj{\CK}.
\end{equation*}

\begin{proposition} \label{prop Comm(T) KtensorEnd(K) is cocartesian}
  The monoidal category
  $\kos{E\!n\!\!d}_{\cat{lax}}^{\kos{K}\backslash\!\!\backslash}(\kos{K},\id_{\kos{K}})_{\cat{crfl}}$
  is cocartesian.
  \begin{enumerate}
    \item 
    Every coreflective lax $\kos{K}$-tensor endofunctor
    $(\phi,\what{\phi})$ on $(\kos{K},\id_{\kos{K}})$
    has a unique lax $\kos{K}$-tensor monad structure
    \begin{equation*}
      \langle\phi,\what{\phi}\rangle
      =
      (\phi,\what{\phi},\mu^{\phi},\upsilon^{\phi})
      \qquad\quad
      \begin{aligned}
        \mu^{\phi}
        &:
        (\phi\phi,\what{\phi\phi})\Rightarrow(\phi,\what{\phi})
        \\
        \upsilon^{\phi}
        &:
        (\id_{\kos{K}},I)\Rightarrow (\phi,\what{\phi})
      \end{aligned}
    \end{equation*}
    where the components of $\mu^{\phi}$, $\upsilon^{\phi}$ at
    $x\in\obj{\CK}$ are given below.
    \begin{equation}
      \label{eq Comm(T) KtensorEnd(K) is cocartesian}
      \begin{aligned}
       \mu^{\phi}_x
        &:\!
        \xymatrix@C=25pt{
          \phi\phi(x)
          \ar[r]^-{(\tahar{\phi}_{\phi(x)})^{-1}}_-{\cong}
          &\phi(x)\otimes \phi(\kappa)
          \ar[r]^-{\phi_{x,\kappa}}
          &\phi(x\otimes \kappa)
          \ar[r]^-{\phi(\jmath^{-1}_x)}_-{\cong}
          &\phi(x)
        }
        \\
        \upsilon^{\phi}_x
        &:\!
        \xymatrix@C=20pt{
          x
          \ar[r]^-{\what{\phi}_x}
          &\phi(x)
        }
      \end{aligned}
    \end{equation}

    \item 
    Every monoidal $\kos{K}$-tensor natural transformation
    $\vartheta:(\phi,\what{\phi})\Rightarrow (\psi,\what{\psi})$
    between coreflective lax $\kos{K}$-tensor endofunctors
    $(\phi,\what{\phi})$, $(\psi,\what{\psi})$
    on $(\kos{K},\id_{\kos{K}})$
    becomes a morphism 
    $\vartheta:\langle\phi,\what{\phi}\rangle\Rightarrow \langle\psi,\what{\psi}\rangle$
    of lax $\kos{K}$-tensor monads.
  \end{enumerate}
\end{proposition}
\begin{proof}
  As we have the adjoint equivalence $\iota\dashv\CR$
  of monoidal categories described in Proposition~\ref{prop Comm(T) Comm(K)=KtensorEnd(K)},
  we conclude that 
  $\kos{E\!n\!\!d}_{\cat{lax}}^{\kos{K}\backslash\!\!\backslash}(\kos{K},\id_{\kos{K}})_{\cat{crfl}}$
  is cocartesian monoidal.
  Therefore statements 1 and 2 are true,
  but we still need to verify
  the explicit descriptions 
  (\ref{eq Comm(T) KtensorEnd(K) is cocartesian})
  of $\mu^{\phi}$ and $\upsilon^{\phi}$.
  For each coreflective lax $\kos{K}$-tensor endofunctor $(\phi,\what{\phi})$
  on $(\kos{K},\id_{\kos{K}})$,
  the coreflection $\tahar{\phi}$
  becomes an isomorphism
  \begin{equation*}
    \xymatrix@C=15pt{
      \hatar{\phi}:
      \langle \otimes\phi(\kappa),\what{\otimes\phi(\kappa)}\rangle
      \ar@2{->}[r]^-{\cong}
      &\langle \phi,\what{\phi}\rangle
    }
  \end{equation*}
  of coreflective lax $\kos{K}$-tensor monads on $(\kos{K},\id_{\kos{K}})$.
  Thus we have the following relations
  for each object $x$ in $\CK$.
  \begin{equation*}
    \vcenter{\hbox{
      \xymatrix@C=20pt{
        (x\otimes \phi(\kappa))\otimes \phi(\kappa)
        \ar[d]_-{\mu^{\otimes\phi(\kappa)}_x}
        \ar[rr]^-{(\tahar{\phi}\tahar{\phi})_x}_-{\cong}
        &\text{ }
        &\phi\phi(x)
        \ar[d]^-{\mu^{\phi}_x}
        &\text{ }
        &x
        \ar[dl]_-{\upsilon^{\otimes\phi(\kappa)}_x}^-{\cong}
        \ar[dr]^-{\upsilon^{\phi}_x}_-{\cong}
        &\text{ }
        \\
        x\otimes \phi(\kappa)
        \ar[rr]^-{\tahar{\phi}_x}_-{\cong}
        &\text{ }
        &\phi(x)
        &x\otimes \phi(\kappa)
        \ar[rr]^-{\tahar{\phi}_x}_-{\cong}
        &\text{ }
        &\phi(x)
      }
    }}
  \end{equation*}
  These relations determine $\mu^{\phi}$ and $\upsilon^{\phi}$.
  One can easily check that $\upsilon^{\phi}=\what{\phi}$.
  We are left to show that $\mu^{\phi}_x$ is as claimed in (\ref{eq Comm(T) KtensorEnd(K) is cocartesian}).
  We have the relation
  \begin{equation*}
    (\dagger):
    \vcenter{\hbox{
      \xymatrix@C=20pt{
        (x\otimes\phi(\kappa))\otimes \phi(\kappa)
        \ar[dd]_-{\mu^{\otimes\phi(\kappa)}_x}
        \ar@{=}[rr]
        &\text{ }
        &(x\otimes\phi(\kappa))\otimes \phi(\kappa)
        \ar[ddd]^-{\tahar{\phi}_x\otimes I_{\phi(\kappa)}}_-{\cong}
        \ar@/^0.5pc/[dl]|-{(\what{\phi}_x\otimes I_{\phi(\kappa)})\otimes I_{\phi(\kappa)}}
        \\
        \text{ }
        &(\phi(x)\otimes \phi(\kappa))\otimes \phi(\kappa)
        \ar[d]^-{\phi_{x,\kappa}\otimes I_{\phi(\kappa)}}
        \ar@/_1pc/[ddl]|-{\mu^{\otimes\phi(\kappa)}_{\phi(x)}}
        &\text{ }
        \\
        x\otimes \phi(\kappa)
        \ar[d]_-{\what{\phi}_x\otimes I_{\phi(\kappa)}}
        &\phi(x\otimes \kappa)\otimes \phi(\kappa)
        \ar@/^0.5pc/[dr]^-{\phi(\jmath_x^{-1})\otimes I_{\phi(\kappa)}}_-{\cong}
        \ar@{}[dd]|-{(\dagger1)}
        &\text{ }
        \\
        \phi(x)\otimes \phi(\kappa)
        \ar[d]_-{\phi_{x,\kappa}}
        &\text{ }
        &\phi(x)\otimes \phi(\kappa)
        \ar[d]^-{\phi_{x,\kappa}}
        \\
        \phi(x\otimes \kappa)
        \ar@{=}[rr]
        &\text{ }
        &\phi(x\otimes \kappa)
      }
    }}
  \end{equation*}
  and we leave for the readers to check the diagram $(\dagger1)$.
  Using the relation $(\dagger)$,
  we obtain the description of $\mu^{\phi}_x$
  as we claimed.
  \begin{equation*}
    \vcenter{\hbox{
      \xymatrix@C=20pt{
        \phi\phi(x)
        \ar[ddddd]_-{\mu^{\phi}_x}
        \ar@{=}[r]
        &\phi\phi(x)
        \ar[d]^-{(\tahar{\phi}\tahar{\phi})^{-1}_x}_-{\cong}
        \ar@{=}[rr]
        &\text{ }
        &\phi\phi(x)
        \ar[ddd]^-{(\tahar{\phi}_{\phi(x)})^{-1}}_-{\cong}
        \\
        \text{ }
        &(x\otimes \phi(\kappa))\otimes \phi(\kappa)
        \ar[d]^-{\mu^{\otimes\phi(\kappa)}_x}
        \ar@/^1pc/[ddrr]^-{\tahar{\phi}_x\otimes I_{\phi(\kappa)}}_-{\cong}
        &\text{ }
        &\text{ }
        \\
        \text{ }
        &x\otimes \phi(\kappa)
        \ar[ddd]^-{\tahar{\phi}_x}_-{\cong}
        \ar@/_0.5pc/[dr]|-{\what{\phi}_x\otimes I_{\phi(\kappa)}}
        \ar@{}[drr]|-{(\dagger)}
        &\text{ }
        &\text{ }
        \\
        \text{ }
        &\text{ }
        &\phi(x)\otimes \phi(\kappa)
        \ar@/_0.5pc/[dr]|-{\phi_{x,\kappa}}
        &\phi(x)\otimes \phi(\kappa)
        \ar[d]^-{\phi_{x,\kappa}}
        \\
        \text{ }
        &\text{ }
        &\text{ }
        &\phi(x\otimes \kappa)
        \ar[d]^-{\phi(\jmath_x^{-1})}_-{\cong}
        \\
        \phi(x)
        \ar@{=}[r]
        &\phi(x)
        \ar@{=}[rr]
        &\text{ }
        &\phi(x)
      }
    }}
  \end{equation*}
  This completes the proof of Proposition~\ref{prop Comm(T) KtensorEnd(K) is cocartesian}.
\qed\end{proof}

\subsection{Left-strong $\kos{K}$-tensor adjunctions}
\label{subsec LSKtensoradj}

In this subsection we study about adjunctions
between lax $\kos{K}$-tensor categories.
We begin by introducing adjunctions
between symmetric monoidal categories
in lax context.

\begin{definition} \label{def LSKtensoradj LSSMadj}
  A \emph{left-strong symmetric monoidal (LSSM) adjunction}
  is an adjunction
  internal to the $2$-category
  $\mathbb{SMC}_{\cat{lax}}$.
\end{definition}

Let us explain Definition~\ref{def LSKtensoradj LSSMadj}
in detail.
Let $\kos{T}=(\CT,\tensor,\unit)$,
$\text{$\kos{S}$}=(\CS,\ctimes,\pzc{1})$
be symmetric monoidal categories.
A LSSM adjunction
\begin{equation*}
  \kos{f}:\kos{S}\to \kos{T}
\end{equation*}
is a data of a pair of lax symmetric monoidal functors
$\kos{f}_*:\kos{S}\to \kos{T}$,
$\kos{f}^*:\kos{T}\to \kos{S}$
and an adjunction
$\kos{f}^*:\CT\rightleftarrows \CS:\kos{f}_*$
such that the adjunction unit
$\eta:\id_{\kos{T}}\Rightarrow \kos{f}_*\kos{f}^*$
and the adjunction counit
$\epsilon:\kos{f}^*\kos{f}_*\Rightarrow \id_{\kos{S}}$
are monoidal natural transformations.
We often denote a LSSM adjunction
$\kos{f}:\kos{S}\to\kos{T}$
as follows.
\begin{equation*}
  \vcenter{\hbox{
    \xymatrix{
      \kos{T}
      \ar@/_1pc/[d]_-{\kos{f}^*}
      \\
      \kos{S}
      \ar@/_1pc/[u]_-{\kos{f}_*}
    }
  }}
\end{equation*}
In this case,
the left adjoint $\kos{f}^*$
is a strong symmetric monoidal functor.
\begin{equation} \label{eq LSKtensoradj f*coherence}
  \text{$\kos{f}$}^*_{X,Y}:
  \text{$\kos{f}$}^*(X)\ctimes \text{$\kos{f}$}^*(Y)
  \xrightarrow{\cong}
  \text{$\kos{f}$}^*(X\tensor Y)
  \qquad
  X,Y\in\obj{\CT}
  \qquad
  \text{$\kos{f}$}^*_{\unit}:
  \text{$\pzc{1}$}
  \xrightarrow{\cong}
  \text{$\kos{f}$}^*(\unit)
\end{equation}
The inverses
$(\kos{f}^*_{X,Y})^{-1}:
\text{$\kos{f}$}^*(X\tensor Y)
\xrightarrow{\cong}
\text{$\kos{f}$}^*(X)\ctimes \text{$\kos{f}$}^*(Y)$
and
$(\text{$\kos{f}$}^*_{\unit})^{-1}:
\text{$\kos{f}$}^*(\unit)
\xrightarrow{\cong}
\text{$\pzc{1}$}$
of the coherence morphisms
of $\kos{f}^*$ are morphisms in $\CS$
that correspond to the following morphisms
in $\CT$ under the adjunction $\kos{f}^*\dashv\kos{f}_*$.
\begin{equation*}
  \xymatrix@C=35pt{
    X\tensor Y
    \ar[r]^-{\eta_X\tensor \eta_Y}
    &\kos{f}_*\kos{f}^*(X)\tensor \kos{f}_*\kos{f}^*(Y)
    \ar[r]^-{\kos{f}_{*\kos{f}^*(X),\kos{f}^*(Y)}}
    &\kos{f}_*(\kos{f}^*(X)\ctimes \kos{f}^*(Y))
  }
  \qquad
  \xymatrix{
    \unit
    \ar[r]^-{\kos{f}_{*\pzc{1}}}
    &\kos{f}_*(\pzc{1})
  }
\end{equation*}

Conversely,
suppose we are given an adjunction
$\kos{f}^*:\CT\rightleftarrows \CS:\kos{f}_*$
between the underlying categories $\CT$, $\CS$
and the left adjoint $\kos{f}^*$
has a strong symmetric monoidal functor structure
$\kos{f}^*:\kos{T}\to \kos{S}$.
Then the right adjoint $\kos{f}_*$ has a unique lax symmetric monoidal functor structure
$\kos{f}_*:\kos{S}\to \kos{T}$
such that the given adjunction
$\kos{f}^*\dashv\kos{f}_*$
becomes a LSSM adjunction
$\kos{f}:\kos{S}\to \kos{T}$.
The monoidal coherence morphisms
\begin{equation}\label{eq LSKtensoradj f_*coherence}
  \text{$\kos{f}$}_{*\pzc{X},\pzc{Y}}
  :
  \text{$\kos{f}$}_*(\pzc{X})\tensor \text{$\kos{f}$}_*(\pzc{Y})
  \to
  \text{$\kos{f}$}_*(\pzc{X}\ctimes \pzc{Y})
  \qquad
  \pzc{X},\pzc{Y}\in \obj{\CS}
  \qquad
  \text{$\kos{f}$}_{*\pzc{1}}:
  \unit
  \to
  \text{$\kos{f}$}_*(\pzc{1})
\end{equation}
of $\kos{f}_*$
are morphisms in $\CT$
that correspond to the following morphisms in $\CS$
under the adjunction $\kos{f}^*\dashv\kos{f}_*$.
\begin{equation*}
  \xymatrix@C=35pt{
    \text{$\kos{f}$}^*(\text{$\kos{f}$}_*(\text{$\pzc{X}$})\tensor \text{$\kos{f}$}_*(\text{$\pzc{Y}$}))
    \ar[r]^-{(\text{$\kos{f}$}^*_{\text{$\kos{f}$}_*(\text{$\pzc{X}$}),\text{$\kos{f}$}_*(\text{$\pzc{Y}$})})^{-1}}_-{\cong}
    &\text{$\kos{f}$}^*\text{$\kos{f}$}_*(\text{$\pzc{X}$})\ctimes \text{$\kos{f}$}^*\text{$\kos{f}$}_*(\text{$\pzc{Y}$})
    \ar[r]^-{\epsilon_{\text{$\pzc{X}$}}\ctimes \epsilon_{\text{$\pzc{Y}$}}}
    &\text{$\pzc{X}$}\ctimes \text{$\pzc{Y}$}
  }
  \qquad
  \xymatrix{
    \kos{f}^*(\unit)
    \ar[r]^-{(\kos{f}^*_{\unit})^{-1}}_-{\cong}
    &\text{$\pzc{1}$}
  }
\end{equation*}

\begin{definition}
  Let $\kos{T}$, $\kos{S}$ be symmetric monoidal categories.
  A \emph{morphism} 
  $\vartheta:\kos{f}\Rightarrow\kos{g}$
  of LSSM adjunctions
  $\kos{f}$, $\kos{g}:\kos{S}\to\kos{T}$
  is a monoidal natural transformation
  $\vartheta:\kos{f}^*\Rightarrow \kos{g}^*:\kos{T}\to \kos{S}$
  between left adjoints.
\end{definition}

We denote the $2$-category of
symmetric monoidal categories,
LSSM adjunctions
and morphisms of LSSM adjunctions
as
\begin{equation} \label{eq LSKtensoradj LSSM 2-cat}
  \mathbb{ADJ}_{\cat{left}}(\mathbb{SMC}_{\cat{lax}}).
\end{equation}
We use the subscript $_{\cat{left}}$
to indicate that $2$-cells in
$\mathbb{ADJ}_{\cat{left}}(\mathbb{SMC}_{\cat{lax}})$
are defined as $2$-cells 
in $\mathbb{SMC}_{\cat{lax}}$
between left adjoints.
The $2$-category structure of
$\mathbb{ADJ}_{\cat{left}}(\mathbb{SMC}_{\cat{lax}})$
is defined as one expects.
For instance,
the vertical composition of $2$-cells
$\kos{f}\Rightarrow\kos{g}\Rightarrow\kos{h}$
in
$\mathbb{ADJ}_{\cat{left}}(\mathbb{SMC}_{\cat{lax}})$
is given by
the vertical composition of $2$-cells
$\kos{f}^*\Rightarrow\kos{g}^*\Rightarrow\kos{h}^*$
between left adjoints.

We introduce an important property of LSSM adjunctions.

\begin{definition} \label{def LSKtensoradj projformula}
  Let $\kos{f}:\kos{S}\to \kos{T}$
  be a LSSM adjunction between 
  symmetric monoidal categories $\kos{T}$, $\kos{S}$.
  We define a canonical natural transformation
  \begin{equation*}
    \varphi_{\slot,\slot}
    :
    \slot\tensor \kos{f}_*(\slot)
    \Rightarrow
    \kos{f}_*(\kos{f}^*(\slot)\ctimes \slot)
    :\CT\times \CS\to \CT
  \end{equation*}
  whose component at 
  $X\in\obj{\CT}$,
  $\text{$\pzc{Z}$}\in\obj{\CS}$
  is given by
  \begin{equation*}
    \varphi_{X,\text{$\pzc{Z}$}}
    :
    \xymatrix@C=30pt{
      X\tensor \kos{f}_*(\pzc{Z})
      \ar[r]^-{\eta_X\tensor I_{\kos{f}_*(\pzc{Z})}}
      &\kos{f}_*\kos{f}^*(X)\tensor \kos{f}_*(\pzc{Z})
      \ar[r]^-{\text{$\kos{f}$}_{*\text{$\kos{f}$}^*(X),\text{$\pzc{Z}$}}}
      &\text{$\kos{f}$}_*(\text{$\kos{f}$}^*(X)\ctimes \text{$\pzc{Z}$})
      .
    }
  \end{equation*}
  We say a LSSM adjunction $\kos{f}:\kos{S}\to \kos{T}$
  \emph{satisfies the projection formula}
  if the natural transformation
  $\varphi_{\slot,\slot}$
  associated to $\kos{f}$ is a natural isomorphism.
\end{definition}

Next we study about adjunctions
between lax $\kos{K}$-tensor categories.

\begin{definition} \label{def LSKtensoradj RSKtensoradjunction}
  A \emph{left-strong $\kos{K}$-tensor adjunction}
  is an adjunction
  internal to the $2$-category
  $\mathbb{SMC}^{\kos{K}\backslash\!\!\backslash}_{\cat{lax}}$.
\end{definition}

Let us explain Definition~\ref{def LSKtensoradj RSKtensoradjunction} in detail.
Let $(\kos{T},\mon{t})$, $(\kos{S},\mon{s})$
be lax $\kos{K}$-tensor categories.
A left-strong $\kos{K}$-tensor adjunction
\begin{equation*}
  (\kos{f},\what{\kos{f}}):(\kos{S},\mon{s})\to (\kos{T},\mon{t})
\end{equation*}
is a data of a pair of lax $\kos{K}$-tensor functors
$(\kos{f}_*,\what{\kos{f}}_*):(\kos{S},\mon{s})\to (\kos{T},\mon{t})$,
$(\kos{f}^*,\what{\kos{f}}^*):(\kos{T},\mon{t})\to (\kos{S},\mon{s})$
and an adjunction
$\kos{f}^*:\CT\rightleftarrows \CS:\kos{f}_*$
such that the adjunction unit
$\eta:(\id_{\kos{T}},I)\Rightarrow (\kos{f}_*\kos{f}^*,\what{\kos{f}_*\kos{f}^*})$
and the adjunction counit
$\epsilon:(\kos{f}^*\kos{f}_*,\what{\kos{f}^*\kos{f}_*})\Rightarrow (\id_{\kos{S}},I)$
are monoidal $\kos{K}$-tensor natural transformations.
We often denote a left-strong $\kos{K}$-tensor adjunction
$(\kos{f},\what{\kos{f}}):(\kos{S},\mon{s})\to (\kos{T},\mon{t})$
as follows.
\begin{equation*}
  \vcenter{\hbox{
    \xymatrix{
      (\kos{T},\mon{t})
      \ar@/_1pc/[d]_-{(\kos{f}^*,\what{\kos{f}}^*)}
      \\
      (\kos{S},\mon{s})
      \ar@/_1pc/[u]_-{(\kos{f}_*,\what{\kos{f}}_*)}
    }
  }}
\end{equation*}
In this case, the left adjoint 
$(\kos{f}^*,\what{\kos{f}}^*):(\kos{T},\mon{t})\to (\kos{S},\mon{s})$
is a strong $\kos{K}$-tensor functor.
We explained the coherence isomorphisms
of the underlying strong symmetric monoidal functor
$\kos{f}^*:\kos{T}\to \kos{S}$
in (\ref{eq LSKtensoradj f*coherence}).
The inverse of the monoidal natural isomorphism
$\xymatrix@C=15pt{\what{\kos{f}}^*:\mon{s}\ar@2{->}[r]^-{\cong}&\kos{f}^*\mon{t}}$
is explicitly given below.
\begin{equation*}
  \xymatrix@C=30pt{
    (\what{\kos{f}}^*)^{-1}:
    \kos{f}^*\mon{t}
    \ar@2{->}[r]^-{\kos{f}^*\what{\kos{f}}_*}
    &\kos{f}^*\kos{f}_*\mon{s}
    \ar@2{->}[r]^-{\epsilon\mon{s}}
    &\mon{s}
  }
\end{equation*}

Conversely,
suppose we are given an adjunction
$\kos{f}^*:\CT\rightleftarrows \CS:\kos{f}_*$
between the underlying categories $\CT$, $\CS$
and the left adjoint has a strong $\kos{K}$-tensor functor structure
$(\kos{f}^*,\what{\kos{f}}^*):(\kos{T},\mon{t})\to (\kos{S},\mon{s})$.
Then the right adjoint has a unique lax $\kos{K}$-tensor functor structure
$(\kos{f}_*,\what{\kos{f}}_*):
(\kos{S},\mon{s})\to (\kos{T},\mon{t})$
such that the given adjunction
$\kos{f}^*\dashv\kos{f}_*$
becomes a left-strong $\kos{K}$-tensor adjunction
$(\kos{f},\what{\kos{f}}):(\kos{S},\mon{s})\to (\kos{T},\mon{t})$.
We explained the coherence morphisms of $\kos{f}_*$
in (\ref{eq LSKtensoradj f_*coherence}).
The monoidal natural transformation
$\what{\kos{f}}_*:\mon{t}\Rightarrow\kos{f}_*\mon{s}$
is described below.
\begin{equation}\label{eq LSKtensoradj rightadj uniqueKTstr}
  \xymatrix@C=30pt{
    \what{\kos{f}}_*:
    \mon{t}
    \ar@2{->}[r]^-{\eta\mon{t}}
    &\kos{f}_*\kos{f}^*\mon{t}
    \ar@2{->}[r]^-{\kos{f}_*(\what{\kos{f}}^*)^{-1}}_-{\cong}
    &\kos{f}_*\mon{s}
  }
\end{equation}

\begin{definition}
  Let $(\kos{T},\mon{t})$,
  $(\kos{S},\mon{s})$
  be lax $\kos{K}$-tensor categories.
  A \emph{morphism} 
  $\vartheta:(\kos{f},\what{\kos{f}})\Rightarrow(\kos{g},\what{\kos{g}})$
  of left-strong $\kos{K}$-tensor adjunctions
  $(\kos{f},\what{\kos{f}})$,
  $(\kos{g},\what{\kos{g}}):(\kos{S},\mon{s})\to(\kos{T},\mon{t})$
  is a monoidal $\kos{K}$-tensor natural transformation
  $\vartheta:(\kos{f}^*,\what{\kos{f}}^*)\Rightarrow (\kos{g}^*,\what{\kos{g}}^*)
  :(\kos{T},\mon{t})\to (\kos{S},\mon{s})$
  between left adjoints.
\end{definition}

We denote the $2$-category of
lax $\kos{K}$-tensor categories,
left-strong $\kos{K}$-tensor adjunctions
and morphisms of left-strong $\kos{K}$-tensor adjunctions
as
\begin{equation} \label{eq LSKtensoradj LSKTadj 2-cat}
  \mathbb{ADJ}_{\cat{left}}\bigl(\mathbb{SMC}^{\kos{K}\backslash\!\!\backslash}_{\cat{lax}}\bigr).
\end{equation}
Again, we use the subscript $_{\cat{left}}$
to indicate that $2$-cells in
$\mathbb{ADJ}_{\cat{left}}\bigl(\mathbb{SMC}^{\kos{K}\backslash\!\!\backslash}_{\cat{lax}}\bigr)$
are defined as $2$-cells in
$\mathbb{SMC}^{\kos{K}\backslash\!\!\backslash}_{\cat{lax}}$
between left adjoints.
The $2$-category structure of
$\mathbb{ADJ}_{\cat{left}}\bigl(\mathbb{SMC}^{\kos{K}\backslash\!\!\backslash}_{\cat{lax}}\bigr)$
is defined as one expects.
For instance,
the vertical composition of $2$-cells in
$\mathbb{ADJ}_{\cat{left}}\bigl(\mathbb{SMC}^{\kos{K}\backslash\!\!\backslash}_{\cat{lax}}\bigr)$
is given by
the vertical composition of $2$-cells
between left adjoints.

\begin{definition}
  Let $(\kos{T},\mon{t})$, $(\kos{S},\mon{s})$
  be lax $\kos{K}$-tensor categories.
  We say a left-strong $\kos{K}$-tensor adjunction
  $(\kos{f},\what{\kos{f}}):(\kos{S},\mon{s})\to (\kos{T},\mon{t})$
  \emph{satisfies the projection formula}
  if the underlying LSSM adjunction
  $\kos{f}:\kos{S}\to \kos{T}$
  satisfies the projection formula.
\end{definition}

\begin{remark} \label{rem LSKtensoradj projformula as Tequiv}
  Let $\kos{T}$, $\kos{S}$ be symmetric monoidal categories
  and let $\kos{f}:\kos{S}\to \kos{T}$ be a LSSM adjunction.
  We explain how we can enhence the given LSSM adjunction
  $\kos{f}$
  to a left-strong $\kos{T}$-tensor adjunction.
  \begin{equation*}
    \vcenter{\hbox{
      \xymatrix{
        \kos{T}
        \ar@/_1pc/[d]_-{\kos{f}^*}
        \\
        \kos{S}
        \ar@/_1pc/[u]_-{\kos{f}_*}
      }
    }}
  \end{equation*}
  We can see $\kos{T}$ as a strong $\kos{T}$-tensor category
  which we denote as $(\kos{T},\id_{\kos{T}})$.
  We have a strong $\kos{T}$-tensor category
  $(\kos{S},\kos{f}^*)$
  as the left adjoint $\kos{f}^*:\kos{T}\to \kos{S}$
  is a strong symmetric monoidal functor.
  Then we obtain a left-strong $\kos{T}$-tensor adjunction
  \begin{equation*}
    \vcenter{\hbox{
      \xymatrix{
        (\kos{T},\id_{\kos{T}})
        \ar@/_1pc/[d]_-{(\kos{f}^*,I)}
        \\
        (\kos{S},\kos{f}^*)
        \ar@/_1pc/[u]_-{(\kos{f}_*,\eta)}
      }
    }}
  \end{equation*}
  as the left adjoint
  has a strong $\kos{T}$-tensor functor structure
  $(\kos{f}^*,I):(\kos{T},\id_{\kos{T}})\to (\kos{S},\kos{f}^*)$.
  The right adjoint lax $\kos{T}$-tensor functor
  is $\kos{f}_*$ equipped with
  $\eta:\id_{\kos{T}}\Rightarrow \kos{f}_*\kos{f}^*$.
  \begin{equation*}
    \vcenter{\hbox{
      \xymatrix@R=30pt@C=40pt{
        \text{ }
        &\kos{S}
        \ar[d]^-{\kos{f}_*}
        \\
        \kos{T}
        \ar@/^0.7pc/[ur]^-{\kos{f}^*}
        \ar[r]_-{\id_{\kos{T}}}
        &\kos{T}
        \xtwocell[l]{}<>{<2.5>{\eta\text{ }}}
      }
    }}
  \end{equation*}
  One can check that the associated $\kos{T}$-equivariance
  $\cevar{\kos{f}}_*^{\kos{T}}$
  of the right adjoint 
  is equal to the natural transformation
  $\varphi_{\slot,\slot}$
  defined in Definition~\ref{def LSKtensoradj projformula}.
  Let $X\in\obj{\CT}$, $\pzc{Z}\in\obj{\CS}$.
  \begin{equation*}
    \vcenter{\hbox{
      \xymatrix@C=50pt{
        \kos{f}_*(X\acts_{\kos{f}^*}\pzc{Z})
        \ar@{=}[d]
        &X\acts_{\id_{\kos{T}}}\kos{f}_*(\pzc{Z})
        \ar[l]_-{(\cevar{\kos{f}}_!^{\kos{T}})_{X,\pzc{Z}}}
        \ar@{=}[d]
        \\
        \text{$\kos{f}$}_*(\text{$\kos{f}$}^*(X)\ctimes \pzc{Z})
        &X\tensor \kos{f}_*(\pzc{Z})
        \ar[l]_-{\varphi_{X,\pzc{Z}}}
      }
    }}
  \end{equation*}
\end{remark}

Recall the category
$\cat{Comm}(\kos{K})$
of commutative monoids in $\kos{K}$
as well as the opposite category
$\cat{Aff}(\kos{K})$
introduced in (\ref{eq Comm(T) Comm(T)def})
and (\ref{eq Comm(T) Aff(T)def}).
For each object $b$ in $\cat{Comm}(\kos{K})$,
we have a symmetric monoidal category
\begin{equation*}
  \kos{K}_b=(\CK_b,\otimes_b,\kos{b}^*(\kappa))
\end{equation*}
and a LSSM adjunction
$\kos{b}:\kos{K}_b\to \kos{K}$
satisfying the projection formula.
\begin{equation*}
  \vcenter{\hbox{
    \xymatrix{
      \kos{K}
      \ar@/_1pc/[d]_-{\kos{b}^*}
      \\
      \kos{K}_b
      \ar@/_1pc/[u]_-{\kos{b}_*}
    }
  }}
\end{equation*}
\begin{itemize}
  \item 
  The underlying category $\CK_b$ of $\kos{K}_b$
  is the Kleisli category associated to the monad
  $\langle \otimes b\rangle=(\otimes b,\mu^{\otimes b},\upsilon^{\otimes b})$
  on $\CK$.
  An object in $\CK_b$ is denoted as $\kos{b}^*(x)$ where $x$ is an object in $\CK$.
  A morphism $\kos{b}^*(x)\xrightarrow{l} \kos{b}^*(y)$ in $\CK_b$
  is a morphism $x\xrightarrow{l} y\otimes b$ in $\CK$.
  The composition of morphisms
  $\kos{b}^*(x)\xrightarrow{l}\kos{b}^*(y)\xrightarrow{l^{\pr}}\kos{b}^*(z)$
  in $\CK_b$ is given by
  $x
  \xrightarrow{l}
  y\otimes b
  \xrightarrow{l^{\pr}\otimes I_b}
  (z\otimes b)\otimes b
  \xrightarrow{\mu^{\otimes b}_z}
  z\otimes b$
  and the identity morphism of $\kos{b}^*(x)$ in $\CK_b$ is
  $I_{\kos{b}^*(x)}=\upsilon^{\otimes b}_x:x\to x\otimes b$.

  \item
  We have an adjunction
  $\kos{b}^*:\CK\rightleftarrows \CK_b:\kos{b}_*$.
  The left adjoint $\kos{b}^*$ sends each object $x$ in $\CK$ to $\kos{b}^*(x)\in\obj{\CK_b}$.
  The right adjoint $\kos{b}_*$ sends each object $\kos{b}^*(z)$ in $\CK_b$
  to $\kos{b}_*\kos{b}^*(z)=z\otimes b\in\obj{\CK}$.
  The component of the adjunction unit at $x\in\obj{\CK}$
  is $\upsilon^{\otimes b}_x:x\to x\otimes b=\kos{b}_*\kos{b}^*(x)$.
  The component $\kos{b}^*\kos{b}_*\kos{b}^*(z)\to \kos{b}^*(z)$
  of the adjunction counit at $\kos{b}^*(z)\in\obj{\CK_b}$
  is $I_{z\otimes b}:z\otimes b\xrightarrow{\cong}z\otimes b$.
  
  \item
  The symmetric monoidal category structure of
  $\kos{K}_b$
  is induced from that of $\kos{K}$
  via the left adjoint $\kos{b}^*$.
  For instance, the monoidal product in $\kos{K}_b$ is given by
  $\kos{b}^*(x)\otimes_b\kos{b}^*(y)=\kos{b}^*(x\otimes y)$
  where $x$, $y\in\obj{\CK}$.
  The left adjoint is then a strong symmetric monoidal functor
  $\kos{b}^*:\kos{K}\to \kos{K}_b$
  whose coherence isomorphisms are identity morphisms.
  The adjunction $\kos{b}^*\dashv \kos{b}_*$
  becomes a LSSM adjunction $\kos{b}:\kos{K}_b\to \kos{K}$,
  and the lax symmetric monoidal coherences of
  $\kos{b}_*:\kos{K}_b\to \kos{K}$ are described below.
  \begin{equation*}
    \vcenter{\hbox{
      \xymatrix@C=33pt{
        \kos{b}_*\kos{b}^*(b)\otimes \kos{b}_*\kos{b}^*(y)
        \ar@{=}[d]
        \ar[r]^-{(\kos{b}_*)_{\kos{b}^*(x),\kos{b}^*(y)}}
        &\kos{b}_*(\kos{b}^*(x)\otimes_b\kos{b}^*(y))
        \ar@{=}[d]
        &\kappa
        \ar@{=}[d]
        \ar[r]^-{(\kos{b}_*)_{\kos{b}^*(\kappa)}}
        &\kos{b}_*\kos{b}^*(\kappa)
        \ar@{=}[d]
        \\
        (x\otimes b)\otimes (y\otimes b)
        \ar[r]^-{(\otimes b)_{x,y}}
        &(x\otimes y)\otimes b
        &\kappa
        \ar[r]^-{(\otimes b)_{\kappa}}
        &\kappa\otimes b
      }
    }}
  \end{equation*}

  \item
  The LSSM adjunction $\kos{b}:\kos{K}_b\to \kos{K}$
  satisfies the projection formula,
  as we can see from the diagram below.
  Let $x\in\obj{\CK}$ and $\kos{b}^*(z)\in \obj{\CK_b}$.
  \begin{equation*}
    \vcenter{\hbox{
      \xymatrix@C=40pt{
        x\otimes \kos{b}_*\kos{b}^*(z)
        \ar[d]^-{\upsilon^{\otimes b}_x\otimes I_{\kos{b}_*\kos{b}^*(z)}}
        \ar@/_2pc/@<-5ex>[dd]_-{\varphi_{x,\kos{b}^*(z)}}
        \ar@{=}[r]
        &x\otimes (z\otimes b)
        \ar[d]^-{\upsilon^{\otimes b}_x\otimes I_{z\otimes b}}
        \ar@<6ex>@/^2pc/[dd]^-{a_{x,z,b}}_-{\cong}
        \\
        \kos{b}_*\kos{b}^*(x)\otimes \kos{b}_*\kos{b}^*(z)
        \ar[d]^-{(\kos{b}_*)_{\kos{b}^*(x),\kos{b}^*(z)}}
        \ar@{=}[r]
        &(x\otimes b)\otimes (z\otimes b)
        \ar[d]^-{(\otimes b)_{x,z}}
        \\
        \kos{b}_*(\kos{b}^*(x)\otimes_b\kos{b}^*(z))
        \ar@{=}[r]
        &(x\otimes z)\otimes b
      }
    }}
  \end{equation*}
\end{itemize}
Thus we obtain a left-strong $\kos{K}$-tensor adjunction
\begin{equation*}
  \vcenter{\hbox{
    \xymatrix{
      (\kos{K},\id_{\kos{K}})
      \ar@/_1pc/[d]_-{(\kos{b}^*,\what{\kos{b}}^*)}
      \\
      (\kos{K}_b,\kos{b}^*)
      \ar@/_1pc/[u]_-{(\kos{b}_*,\what{\kos{b}}_*)}
    }
  }}
\end{equation*}
as explained in Remark~\ref{rem LSKtensoradj projformula as Tequiv}.
\begin{itemize}
  \item 
  The monoidal natural isomorphism
  $\what{\kos{b}}^*:\kos{b}^*=\kos{b}^*$
  is the identity natural transformation of $\kos{b}^*$,
  and the monoidal natural transformation
  $\what{\kos{b}}_*:\id_{\kos{K}}\Rightarrow\kos{b}_*\kos{b}^*$
  is given by $\upsilon^{\otimes b}:\id_{\kos{K}}\Rightarrow \otimes b$.
  One can check that
  \begin{equation*}
    (\kos{b}_*\kos{b}^*,\what{\kos{b}_*\kos{b}^*})
    =(\otimes b,\what{\otimes b})
    :(\kos{K},\id_{\kos{K}})\to (\kos{K},\id_{\kos{K}})
    .
  \end{equation*}

  \item
  The associated $\kos{K}$-equivariance $\cevar{\kos{b}}_*$
  of the right adjoint
  $(\kos{b}_*,\what{\kos{b}}_*)$
  is a natural isomorphism:
  \begin{equation*}
    \cevar{\kos{b}}_{*x,\kos{b}^*(z)}
    =
    \varphi_{x,\kos{b}^*(z)}
    =
    a_{x,z,b}:
    \xymatrix{
      x\otimes (z\otimes b)
      \ar[r]^-{\cong}
      &(x\otimes z)\otimes b
      .
    }
  \end{equation*}
\end{itemize}
Let $f:b\to b^{\pr}$ be a morphism in $\cat{Comm}(\kos{K})$
and consider the corresponding morphism
$f:\cat{Spec}(b^{\pr})\to \cat{Spec}(b)$
in $\cat{Aff}(\kos{K})$.
We have a strong $\kos{K}$-tensor functor
\begin{equation} \label{eq LSKtensoradj functorbetweenKelisliCat}
  (\kos{f}^*,\what{\kos{f}}^*):
  (\kos{K}_b,\kos{b}^*)\to (\kos{K}_{b^{\pr}},\kos{b}^{\pr*})
\end{equation}
which satisfies the relation
\begin{equation*}
  \vcenter{\hbox{
    \xymatrix@C=0pt{
      \text{ }
      &(\kos{K},\id_{\kos{K}})
      \ar[dl]_-{(\kos{b}^*,\what{\kos{b}}^*)}
      \ar[dr]^-{(\kos{b}^{\pr*},\what{\kos{b}}^{\pr*})}
      \\
      (\kos{K}_b,\kos{b}^*)
      \ar[rr]^-{(\kos{f}^*,\what{\kos{f}}^*)}
      &\text{ }
      &(\kos{K}_{b^{\pr}},\kos{b}^{\pr*})
    }
  }}
  \qquad\quad
  (\kos{f}^*,\what{\kos{f}}^*)(\kos{b}^*,\what{\kos{b}}^*)
  =
  (\kos{b}^{\pr*},\what{\kos{b}}^{\pr*}).
\end{equation*}
\begin{itemize}
  \item 
  The functor $\kos{f}^*$ sends each object $\kos{b}^*(x)$ in $\CK_b$
  to the object
  $\kos{f}^*\kos{b}^*(x)=\kos{b}^{\pr*}(x)$ in $\CK_{b^{\pr}}$,
  and each morphism $\kos{b}^*(x)\xrightarrow{l} \kos{b}^*(y)$ in $\CK_b$
  to the morphism $\kos{b}^{\pr*}(x)\xrightarrow{\kos{f}^*(l)} \kos{b}^{\pr*}(y)$ in $\CK_{b^{\pr}}$
  where $\kos{f}^*(l):x\xrightarrow{l}y\otimes b\xrightarrow{I_y\otimes f}y\otimes b^{\pr}$.

  \item
  The symmetric monoidal coherence isomorphisms of
  $\kos{f}^*:\kos{K}_b\to \kos{K}_{b^{\pr}}$
  are identity morphisms,
  and the monoidal natural isomorphism
  $\what{\kos{f}}^*:\kos{b}^{\pr*}=\kos{f}^*\kos{b}^*$
  is the identity natural transformation.
\end{itemize}

\begin{definition} \label{def LSKtensoradj Hom,Isom}
  Let $(\kos{T},\mon{t})$ be a lax $\kos{K}$-tensor category
  and let $(\phi,\what{\phi})$,
  $(\psi,\what{\psi}):(\kos{T},\mon{t})\to (\kos{K},\id_{\kos{K}})$
  be lax $\kos{K}$-tensor functors.
  \begin{enumerate}
    \item 
    We define the \emph{presheaf of monoidal $\kos{K}$-tensor natural transformations
    from $(\phi,\what{\phi})$ to $(\psi,\what{\psi})$}
    as the presheaf on $\cat{Aff}(\kos{K})$
    \begin{equation*}
      \underline{\Hom}_{\mathbb{SMC}_{\cat{lax}}^{\kos{K}\backslash\!\!\backslash}}
      \bigl(
        (\phi,\what{\phi})
        ,
        (\psi,\what{\psi})
      \bigr)
      :
      \cat{Aff}(\kos{K})^{\op}
      \to
      \cat{Set}
    \end{equation*}
    which sends each object $\cat{Spec}(b)$ in $\cat{Aff}(\kos{K})$
    to the set of monoidal $\kos{K}$-tensor natural transformations
    $(\kos{b}^*\phi,\what{\kos{b}^*\phi})
    \Rightarrow
    (\kos{b}^*\psi,\what{\kos{b}^*\psi})
    :(\kos{T},\mon{t})\to (\kos{K}_b,\kos{b}^*)$.
    
    \item 
    We define the \emph{presheaf of monoidal $\kos{K}$-tensor natural isomorphisms
    from $(\phi,\what{\phi})$ to $(\psi,\what{\psi})$}
    as the presheaf on $\cat{Aff}(\kos{K})$
    \begin{equation*}
      \underline{\Isom}_{\mathbb{SMC}_{\cat{lax}}^{\kos{K}\backslash\!\!\backslash}}
      \bigl(
        (\phi,\what{\phi})
        ,
        (\psi,\what{\psi})
      \bigr)
      :
      \cat{Aff}(\kos{K})^{\op}
      \to
      \cat{Set}
    \end{equation*}
    which sends each object $\cat{Spec}(b)$ in $\cat{Aff}(\kos{K})$
    to the set of monoidal $\kos{K}$-tensor natural isomorphisms
    $\xymatrix@C=12pt{(\kos{b}^*\phi,\what{\kos{b}^*\phi})
    \ar@2{->}[r]^-{\cong}
    &(\kos{b}^*\psi,\what{\kos{b}^*\psi})
    :(\kos{T},\mon{t})\to (\kos{K}_b,\kos{b}^*).}$
  \end{enumerate}
\end{definition}





For the rest of this subsection,
we study about left-strong $\kos{K}$-tensor adjunctions
from $(\kos{K},\id_{\kos{K}})$
to a strong $\kos{K}$-tensor category $(\kos{T},\mon{t})$.

\begin{definition} \label{def LSKtensoradj coreflective}
  Let $(\kos{T},\mon{t})$ be a strong $\kos{K}$-tensor category
  and let $(\omega,\what{\omega}):(\kos{K},\id_{\kos{K}})\to (\kos{T},\mon{t})$
  be a left-strong $\kos{K}$-tensor adjunction.
  \begin{equation*}
    \vcenter{\hbox{
      \xymatrix{
        (\kos{T},\mon{t})
        \ar@/_1pc/[d]_-{(\omega^*,\what{\omega}^*)}
        \\
        (\kos{K},\id_{\kos{K}})
        \ar@/_1pc/[u]_-{(\omega_*,\what{\omega}_*)}
      }
    }}
  \end{equation*}
  \begin{enumerate}
    \item
    We say $(\omega,\what{\omega})$
    is \emph{coreflective} if the right adjoint
    $(\omega_*,\what{\omega}_*):(\kos{K},\id_{\kos{K}})\to (\kos{T},\mon{t})$
    is a coreflective lax $\kos{K}$-tensor functor.
    This amounts to saying that
    the associated $\kos{K}$-equivariance
    $\cevar{\omega}_*$
    of the right adjoint
    is a natural isomorphism:
    see Lemma~\ref{lem Comm(T) equivariance and coreflection}.

    \item 
    We define a natural transformation
    \begin{equation*}
      \varphi_{\slot}:
      \slot\tensor \omega_*(\kappa)
      \Rightarrow
      \omega_*\omega^*:\CT\to \CT
    \end{equation*}
    whose component at $X\in\obj{\CT}$ is
    \begin{equation} \label{eq LSKtensoradj modifiedvarphi def}
      \varphi_X:
      \xymatrix@C=30pt{
        X\tensor \omega_*(\kappa)
        \ar[r]^-{\varphi_{X,\kappa}}
        &\omega_*(\omega^*(X)\otimes \kappa)
        \ar[r]^-{\omega^*(\jmath^{-1}_{\omega^*(X)})}_-{\cong}
        &\omega_*\omega^*(X)
        .
      }
    \end{equation}
  \end{enumerate}
\end{definition}

\begin{lemma} \label{lem LSKtensoradj phi-modifiedphi}
  Let $(\kos{T},\mon{t})$ be a strong $\kos{K}$-tensor category
  and let $(\omega,\what{\omega}):(\kos{K},\id_{\kos{K}})\to (\kos{T},\mon{t})$
  be a coreflective left-strong $\kos{K}$-tensor adjunction.
  \begin{equation*}
    \vcenter{\hbox{
      \xymatrix{
        (\kos{T},\mon{t})
        \ar@/_1pc/[d]_-{(\omega^*,\what{\omega}^*)}
        \\
        (\kos{K},\id_{\kos{K}})
        \ar@/_1pc/[u]_-{(\omega_*,\what{\omega}_*)\text{ coreflective}}
      }
    }}
  \end{equation*}
  For each object $X$ in $\CT$, the following are equivalent:
  \begin{enumerate}
    \item 
    $\varphi_{X,\slot}
    :X\tensor \omega_*(\slot)
    \Rightarrow\omega_*(\omega^*(X)\otimes \slot)
    :\CK\to \CT$
    is a natural isomorphism;

    \item
    $\varphi_X:
    X\tensor \omega_*(\kappa)
    \to \omega_*\omega^*(X)$
    is an isomorphism in $\CT$.
  \end{enumerate}
  In particular,
  $(\omega,\what{\omega})$ satisfies the projection formula
  if and only if
  $\varphi_{\slot}:
  \slot\tensor \omega_*(\kappa)
  \Rightarrow
  \omega_*\omega^*:\CT\to \CT$
  is a natural isomorphism.
\end{lemma}
\begin{proof}
  Statement 1 implies statemet 2
  as we can see from the definition (\ref{eq LSKtensoradj modifiedvarphi def})
  of $\varphi_X$.
  To prove the converse,
  it suffices to show that
  for each object $z$ in $\CK$,
  the following diagram commutes.
  \begin{equation} \label{eq LSKtensoradj phi}
    \vcenter{\hbox{
      \xymatrix@C=40pt{
        \mon{t}(z)\tensor X\tensor \omega_*(\kappa)
        \ar[d]_-{I_{\mon{t}(z)}\tensor \varphi_X}
        \ar[r]^-{s_{\mon{t}(z),X}\tensor I_{\omega_*(\kappa)}}_-{\cong}
        &X\tensor \mon{t}(z)\tensor \omega_*(\kappa)
        \ar[r]^-{I_X\tensor (\tahar{\omega}_*)_z}_-{\cong}
        &X\tensor \omega_*(z)
        \ar[d]^-{\varphi_{X,z}}
        \\
        \mon{t}(z)\tensor \omega_*\omega^*(X)
        \ar[r]^-{(\cevar{\omega}_*)_{z,\omega^*(X)}}_-{\cong}
        &\omega_*(z\otimes \omega^*(X))
        \ar[r]^-{\omega^*(s_{z,\omega^*(X)})}_-{\cong}
        &\omega_*(\omega^*(X)\otimes z)
      }
    }}
  \end{equation}
  We can check this as follows.
  We leave for the readers to check the relation
  \begin{equation*}
    (\dagger):
    \vcenter{\hbox{
      \xymatrix{
        \mon{t}(z)\tensor X\tensor \omega_*(\kappa)
        \ar[d]_-{I_{\mon{t}(z)}\tensor \varphi_X}
        \ar@{=}[r]
        &\mon{t}(z)\tensor X\tensor \omega_*(\kappa)
        \ar[d]^-{(\what{\omega}_*)_z\tensor \eta_X\tensor I_{\omega_*(\kappa)}}
        \\
        \mon{t}(z)\tensor \omega_*\omega^*(X)
        \ar[dd]_-{(\cevar{\omega}_*)_{z,\omega^*(X)}}^-{\cong}
        &\omega_*(z)\tensor \omega_*\omega^*(X)\tensor \omega_*(\kappa)
        \ar[d]^-{\omega_{*z,\omega^*(X),\kappa}}
        \\
        \text{ }
        &\omega_*(z\otimes \omega^*(X)\otimes \kappa)
        \ar[d]^-{\omega_*(I_z\otimes \jmath_{\omega^*(X)}^{-1})}_-{\cong}
        \\
        \omega_*(z\otimes \omega^*(X))
        \ar@{=}[r]
        &\omega_*(z\otimes \omega^*(X))
      }
    }}
  \end{equation*}
  where the notation $\omega_{*z,\omega^*(X),\kappa}$
  means the lax symmetric monoidal coherence of $\omega_*:\kos{K}\to\kos{T}$
  with respect to three input objects
  $z$, $\omega^*(X)$, $\kappa$ in $\CK$.
  One can also check the following relation
  \begin{equation*}
    (\dagger\!\dagger):
    \vcenter{\hbox{
      \xymatrix{
        X\tensor \mon{t}(z)\tensor \omega_*(\kappa)
        \ar[dd]_-{I_X\tensor (\tahar{\omega}_*)_z}^-{\cong}
        \ar@{=}[r]
        &X\tensor \mon{t}(z)\tensor \omega_*(\kappa)
        \ar[d]^-{\eta_X\tensor (\what{\omega}_*)_z\tensor I_{\omega_*(\kappa)}}
        \\
        \text{ }
        &\omega_*\omega*(X)\tensor \omega_*(z)\tensor \omega_*(\kappa)
        \ar[d]^-{\omega_{*\omega^*(X),z,\kappa}}
        \\
        X\tensor \omega_*(z)
        \ar[d]_-{\varphi_{X,z}}
        &\omega_*(\omega^*(X)\otimes z\otimes \kappa)
        \ar[d]^-{\omega_*(I_{\omega^*(X)}\otimes \jmath_z^{-1})}_-{\cong}
        \\
        \omega_*(\omega^*(X)\otimes z)
        \ar@{=}[r]
        &\omega_*(\omega^*(X)\otimes z)
      }
    }}
  \end{equation*}
  where the notation $\omega_{*\omega^*(X),z,\kappa}$
  means the lax symmetric monoidal coherence of $\omega_*:\kos{K}\to\kos{T}$
  with respect to three input objects
  $\omega^*(X)$, $z$, $\kappa$ in $\CK$.
  Using these relations,
  we verify the diagram (\ref{eq LSKtensoradj phi}) as follows.
  \begin{equation*}
    \vcenter{\hbox{
      \xymatrix@C=0pt{
        \mon{t}(z)\tensor X\tensor \omega_*(\kappa)
        \ar[d]^-{I_{\mon{t}(z)}\tensor \varphi_X}
        \ar@{=}[rr]
        &\text{ }
        &\mon{t}(z)\tensor X\tensor \omega_*(\kappa)
        \ar@<-1ex>[d]|-{(\what{\omega}_*)_z\tensor \eta_X\tensor I_{\omega_*(\kappa)}}
        \ar@{=}[r]
        &\mon{t}(z)\tensor X\tensor \omega_*(\kappa)
        \ar@<1.5ex>[d]_-{s_{\mon{t}(z),X}\tensor I_{\omega_*(\kappa)}}^-{\cong}
        \\
        \mon{t}(z)\tensor \omega_*\omega^*(X)
        \ar[dd]^(0.35){(\cevar{\omega}_*)_{z,\omega^*(X)}}_(0.35){\cong}
        \ar@{}[rr]|-{(\dagger)}
        &\text{ }
        &\omega_*(z)\!\tensor\! \omega_*\omega^*(X)\!\tensor\! \omega_*(\kappa)
        \ar@/_1.5pc/[dl]|-{\omega_{*z,\omega^*(X),\kappa}}
        \ar@<-1ex>[d]|-{s_{\omega_*(z),\omega_*\omega^*(X)}\tensor I_{\omega_*(\kappa)}}
        &X\tensor \mon{t}(z)\tensor \omega_*(\kappa)
        \ar@/^1pc/[dl]|-{\eta_X\tensor (\what{\omega}_*)_z\tensor I_{\omega_*(\kappa)}}
        \ar@<1.5ex>[dd]_(0.65){I_X\tensor (\tahar{\omega}_*)_z}^(0.65){\cong}
        \\
        \text{ }
        &\omega_*(z\otimes \omega^*(X)\otimes \kappa)
        \ar@/_0.5pc/[dl]^-{\omega_*(I_z\otimes \jmath_{\omega^*(X)}^{-1})}_-{\cong}
        \ar@/_0.5pc/[dr]_-{\omega_*(s_{z,\omega^*(X)}\otimes I_{\kappa})}^-{\cong}
        &\omega_*\omega^*(X)\!\tensor\! \omega_*(z)\!\tensor\! \omega_*(\kappa)
        \ar@<-1ex>[d]|-{\omega_{*\omega^*(X),z,\kappa}}
        &\text{ }
        \\
        \omega_*(z\otimes \omega^*(X))
        \ar[d]^-{\omega^*(s_{z,\omega^*(X)})}_-{\cong}
        &\text{ }
        &\omega_*(\omega^*(X)\otimes z\otimes \kappa)
        \ar@<-1ex>[d]_-{\omega_*(I_{\omega^*(X)}\otimes \jmath_z^{-1})}^-{\cong}
        \ar@{}[dr]|-{(\dagger\!\dagger)}
        &X\tensor \omega_*(z)
        \ar@<1.5ex>[d]_-{\varphi_{X,z}}
        \\
        \omega_*(\omega^*(X)\otimes z)
        \ar@{=}[rr]
        &\text{ }
        &\omega_*(\omega^*(X)\otimes z)
        \ar@{=}[r]
        &\omega_*(\omega^*(X)\otimes z)
      }
    }}
  \end{equation*}
  This completes the proof of Lemma~\ref{lem LSKtensoradj phi-modifiedphi}.
\qed\end{proof}

The following proposition is the main result of this subsection.

\begin{proposition} \label{prop LSKtensoradj HomRep,Hom=Isom}
  Let $(\kos{T},\mon{t})$
  be a strong $\kos{K}$-tensor category.
  Suppose we are given a coreflective left-strong $\kos{K}$-tensor adjunction
  $(\omega,\what{\omega})
  :(\kos{K},\id_{\kos{K}})\to (\kos{T},\mon{t})$
  and a strong $\kos{K}$-tensor functor
  $(\omega^{\pr*},\what{\omega}^{\pr*})
  :(\kos{T},\mon{t})\to (\kos{K},\id_{\kos{K}})$.
  \begin{equation*}
    \vcenter{\hbox{
      \xymatrix{
        (\kos{T},\mon{t})
        \ar[d]_-{(\omega^{\pr*},\what{\omega}^{\pr*})}
        \\
        (\kos{K},\id_{\kos{K}})
      }
    }}
    \qquad\qquad
    \vcenter{\hbox{
      \xymatrix{
        (\kos{T},\mon{t})
        \ar@/_1pc/[d]_-{(\omega^*,\what{\omega}^*)}
        \\
        (\kos{K},\id_{\kos{K}})
        \ar@/_1pc/[u]_-{(\omega_*,\what{\omega}_*)\text{ coreflective}}
      }
    }}
  \end{equation*}
  \begin{enumerate}
    \item 
    The presheaf of monoidal $\kos{K}$-tensor natural transformations
    from
    $(\omega^{\pr*},\what{\omega}^{\pr*})$
    to
    $(\omega^*,\what{\omega}^*)$
    is represented by the object
    $\cat{Spec}(\omega^{\pr*}\omega_*(\kappa))$ in $\cat{Aff}(\kos{K})$.
    \begin{equation*}
      \hspace{-1cm}
      \xymatrix@C=12pt{
        \Hom_{\cat{Aff}(\kos{K})}\bigl(\slot, \cat{Spec}(\omega^{\pr*}\omega_*(\kappa))\bigr)
        \ar@2{->}[r]^-{\cong}
        &\underline{\Hom}_{\mathbb{SMC}_{\cat{lax}}^{\kos{K}\backslash\!\!\backslash}}
        \!\bigl(
          (\omega^{\pr*}\!,\what{\omega}^{\pr*})
          ,
          (\omega^*\!,\what{\omega}^*)
        \bigr)
        \!:
        \cat{Aff}(\kos{K})^{\op}
        \!\to
        \cat{Set}
      }
    \end{equation*}
    The universal element is the monoidal $\kos{K}$-tensor natural transformation
    \begin{equation*}
      \begin{aligned}
        \xi
        &:
        \xymatrix{
          (\kos{p}^*\omega^{\pr*},\what{\kos{p}^*\omega^{\pr*}})
          \Rightarrow
          (\kos{p}^*\omega^*,\what{\kos{p}^*\omega^*})
          ,
        }
        \qquad
        p:=\omega^{\pr*}\omega_*(\kappa)
      \end{aligned}
    \end{equation*}
    whose component at each object $X$ in $\CT$ is
    \begin{equation*}
      \begin{aligned}
        \xi_X
        &:
        \kos{p}^*\omega^{\pr*}(X)
        \to
        \kos{p}^*\omega^*(X)
        ,
        \\
        \xi_X
        &:
        \xymatrix@C=40pt{
          \omega^{\pr*}(X)
          \ar[r]^-{\omega^{\pr*}(\eta_X)}
          &\omega^{\pr*}\omega_*\omega^*(X)
          \ar[r]^-{(\tahar{\omega^{\pr*}\omega_*}_{\omega^*(X)})^{-1}}_-{\cong}
          &\omega^*(X)\otimes \omega^{\pr*}\omega_*(\kappa).
        }
      \end{aligned}
    \end{equation*}

    \item
    Let $X$ be an object in $\CT$.
    The following are equivalent:
    \begin{enumerate}
      \item 
      $\omega^{\pr*}(\varphi_X):
      \omega^{\pr*}(X\tensor \omega_*(\kappa))
      \to
      \omega^{\pr*}\omega_*\omega^*(X)$
      is an isomorphism in $\CT$.
      
      \item 
      The component
      $\xi_X:\kos{p}^*\omega^{\pr*}(X)\to\kos{p}^*\omega^*(X)$
      of the universal element $\xi$ at $X$
      is an isomorphism in the Kleisli category $\CK_p$.
      
      \item
      For each object $\cat{Spec}(b)$ in $\cat{Aff}(\kos{K})$
      and for each monoidal $\kos{K}$-tensor natural transformation
      $\vartheta:
      (\kos{b}^*\omega^{\pr*},\what{\kos{b}^*\omega^{\pr*}})
      \Rightarrow
      (\kos{b}^*\omega^*,\what{\kos{b}^*\omega^*})
      :(\kos{T},\mon{t})\to (\kos{K}_b,\kos{b}^*)$,
      the component
      $\vartheta_X:\kos{b}^*\omega^{\pr*}(X)\to \kos{b}^*\omega^*(X)$
      of $\vartheta$ at $X$ is an isomorphism
      in the Kleisli category $\CK_b$.
    \end{enumerate}
    In this case, if we denote
    \begin{equation*}
      \check{\xi}_X:
      \xymatrix@C=30pt{
        \omega^*(X)
        \ar[r]^-{\xi_X^{-1}}
        &\omega^{\pr*}(X)\otimes \omega^{\pr*}\omega_*(\kappa)
        \ar[r]^-{\tahar{\omega^{\pr*}\omega_*}_{\omega^{\pr*}(X)}}_-{\cong}
        &\omega^{\pr*}\omega_*\omega^{\pr*}(X)
      }
    \end{equation*}
    then the inverse $\vartheta_X^{-1}$ is explicitly described as
    \begin{equation*}
      \vartheta_X^{-1}:\!\!
      \xymatrix@C=28pt{
        \omega^*(X)
        \ar[r]^-{\check{\xi}_X}
        &\omega^{\pr*}\omega_*\omega^{\pr*}(X)
        \ar[r]^-{\vartheta_{\omega_*\omega^{\pr*}(X)}}
        &\omega^*\omega_!\omega^{\pr*}(X)\otimes b
        \ar[r]^-{\epsilon_{\omega^{\pr*}(X)}\otimes I_b}
        &\omega^{\pr*}(X)\otimes b
        .
      }
    \end{equation*}

    \item
    We have
    \begin{equation*}
      \underline{\Hom}_{\mathbb{SMC}_{\cat{lax}}^{\kos{K}\backslash\!\!\backslash}}
      \bigl(
        (\omega^{\pr*},\what{\omega}^{\pr*})
        ,
        (\omega^*,\what{\omega}^*)
      \bigr)
      =
      \underline{\Isom}_{\mathbb{SMC}_{\cat{lax}}^{\kos{K}\backslash\!\!\backslash}}
      \bigl(
        (\omega^{\pr*},\what{\omega}^{\pr*})
        ,
        (\omega^*,\what{\omega}^*)
      \bigr)
    \end{equation*}
    if and only if
    $\omega^{\pr*}(\varphi_{\slot}):
    \omega^{\pr*}(\slot\tensor \omega_*(\kappa))
    \Rightarrow
    \omega^{\pr*}\omega_*\omega^*
    :\CT\to \CK$
    is a natural isomorphism.
  \end{enumerate}
\end{proposition}
\begin{proof}
  We first prove statement 1.
  For each object $\cat{Spec}(b)$ in $\cat{Aff}(\kos{K})$,
  we have bijections
  \begin{equation} \label{eq1 LSKtensoradj HomRep,Hom=Isom}
    \begin{aligned}
      \Hom_{\cat{Aff}(\kos{K})}\bigl(\cat{Spec}(b),\text{ }\cat{Spec}(\omega^{\pr*}\omega_*(\kappa))\bigr)
      &\cong
      \Hom_{\cat{Comm}(\kos{K})}(\omega^{\pr*}\omega_*(\kappa),b)
      \\
      &\cong
      \Hom_{\cat{Comm}(\kos{K})}(\omega^{\pr*}\omega_*(\kappa),\kos{b}_*\kos{b}^*(\kappa))
      \\
      &\cong
      \mathbb{SMC}_{\cat{lax}}^{\kos{K}\backslash\!\!\backslash}\bigl((\omega^{\pr*}\omega_*,\what{\omega^{\pr*}\omega_*}),\text{ }(\kos{b}_*\kos{b}^*,\what{\kos{b}_*\kos{b}^*})\bigr)
      \\
      &\cong
      \mathbb{SMC}_{\cat{lax}}^{\kos{K}\backslash\!\!\backslash}\bigl((\kos{b}^*\omega^{\pr*},\what{\kos{b}^*\omega^{\pr*}}),\text{ }(\kos{b}^*\omega^*,\what{\kos{b}^*\omega^*})\bigr)
      .
    \end{aligned}
  \end{equation}
  The second bijection in (\ref{eq1 LSKtensoradj HomRep,Hom=Isom})
  is obtained by post-composing the isomorphism
  $\imath_b:b\xrightarrow{\cong}\kappa\otimes b=\kos{b}_*\kos{b}^*(\kappa)$
  in $\cat{Comm}(\kos{K})$.
  Recall the adjoint equivalence of categories
  \begin{equation*}
    \iota:
    \cat{Comm}(\kos{T})
    \simeq
    \mathbb{SMC}_{\cat{lax}}^{\kos{K}\backslash\!\!\backslash}\bigl((\kos{K},\id_{\kos{K}}),(\kos{T},\mon{t})\bigr)_{\cat{crfl}}
    :\CR
  \end{equation*}
  described in Corollary~\ref{cor Comm(T) strongKcat coreflective adjunction}.
  As the lax $\kos{K}$-tensor functors
  \begin{equation*}
    (\omega^{\pr*}\omega_*,\what{\omega^{\pr*}\omega_*})
    ,
    (\kos{b}_*\kos{b}^*,\what{\kos{b}_*\kos{b}^*})
    :(\kos{K},\id_{\kos{K}})\to (\kos{K},\id_{\kos{K}})
  \end{equation*}
  are both coreflective,
  we obtain the third bijection in (\ref{eq1 LSKtensoradj HomRep,Hom=Isom}).
  The last bijection in (\ref{eq1 LSKtensoradj HomRep,Hom=Isom})
  is the bijective correspondence of mates
  in the $2$-category $\mathbb{SMC}_{\cat{lax}}^{\kos{K}\backslash\!\!\backslash}$.
  \begin{equation*}
    \vcenter{\hbox{
      \xymatrix@R=30pt@C=40pt{
        (\kos{T},\mon{t})
        \ar[r]^-{(\omega^*,\what{\omega}^*)}
        \ar[d]_-{(\omega^{\pr*},\what{\omega}^{\pr*})}
        &(\kos{K},\id_{\kos{K}})
        \ar[d]^-{(\kos{b}^*,\what{\kos{b}}^*)}
        \\
        (\kos{K},\id_{\kos{K}})
        \ar[r]_-{(\kos{b}^*,\what{\kos{b}}^*)}
        &(\kos{K}_b,\kos{b}^*)
        \xtwocell[ul]{}<>{<0>{\vartheta\text{ }\text{ }}}
      }
    }}
    \quad\leftrightsquigarrow\quad
    \vcenter{\hbox{
      \xymatrix@R=30pt@C=40pt{
        (\kos{T},\kos{t})
        \ar[d]_-{(\omega^{\pr*},\what{\omega}^{\pr*})}
        &(\kos{K},\id_{\kos{K}})
        \ar[l]_-{(\omega_*,\what{\omega}_*)}
        \ar[d]^-{(\kos{b}^*,\what{\kos{b}}^*)}
        \\
        (\kos{K},\id_{\kos{K}})
        \xtwocell[ur]{}<>{<0>{\text{ }\vartheta_*}}
        &(\kos{K}_b,\kos{b}^*)
        \ar[l]^-{(\kos{b}_*,\what{\kos{b}}_*)}
      }
    }}
  \end{equation*}
  The mate pair $\vartheta_*$ and $\vartheta$
  are related as follows.
  \begin{equation*}
    \begin{aligned}
      \vartheta_X
      &:
      \xymatrix@C=30pt{
        \omega^{\pr*}(X)
        \ar[r]^-{\omega^{\pr*}(\eta_X)}
        &\omega^{\pr*}\omega_*\omega^*(X)
        \ar[r]^-{(\vartheta_*)_{\omega^*(X)}}
        &\omega^*(X)\otimes b
      }
      ,
      \qquad
      X\in\obj{\CT}
      ,
      \\
      (\vartheta_*)_x
      &:
      \xymatrix@C=30pt{
        \omega^{\pr*}\omega_*(x)
        \ar[r]^-{\vartheta_{\omega_*(x)}}
        &\omega^*\omega_*(x)\otimes b
        \ar[r]^-{\epsilon_x\otimes I_b}
        &x\otimes b
      }   
      ,
      \qquad
      x\in\obj{\CK}
      .
    \end{aligned}
  \end{equation*}
  To complete the proof of statement 1,
  it suffices to show that the bijection (\ref{eq1 LSKtensoradj HomRep,Hom=Isom})
  sends each morphism $g:\omega^{\pr*}\omega_*(\kappa)\to b$ in $\cat{Comm}(\kos{K})$
  to the monoidal $\kos{K}$-tensor natural transformation
  $\vartheta:(\kos{b}^*\omega^{\pr*},\what{\kos{b}^*\omega^{\pr*}})\Rightarrow (\kos{b}^*\omega^*,\what{\kos{b}^*\omega^*})$
  whose component at each object $X$ in $\CT$ is
  \begin{equation*}
    \vartheta_X:
    \xymatrix@C=30pt{
      \omega^{\pr*}(X)
      \ar[r]^-{\xi_X}
      &\omega^*(X)\otimes \omega^{\pr*}\omega_*(\kappa)
      \ar[r]^-{I_{\omega^*(X)}\otimes g}
      &\omega^*(X)\otimes b
      .
    }
  \end{equation*}
  We can check this as follows.
  The component of the corresponding lax $\kos{K}$-tensor natural transformation
  $\vartheta_*:(\omega^{\pr*}\omega_*,\what{\omega^{\pr*}\omega_*})
  \Rightarrow 
  (\kos{b}_*\kos{b}^*,\what{\kos{b}_*\kos{b}^*})$
  at each object $x$ in $\CK$ is
  \begin{equation*}
    (\vartheta_*)_x:
      \xymatrix@C=35pt{ 
      \omega^{\pr*}\omega_*(x)
      \ar[r]^-{(\tahar{\omega^{\pr*}\omega_*}_x)^{-1}}_-{\cong}
      &x\otimes \omega^{\pr*}\omega_*(\kappa)
      \ar[r]^-{I_x\otimes g}
      &x\otimes b
      .
    }
  \end{equation*}
  Therefore the component of $\vartheta$ at $X$ is given as follows.
  \begin{equation*}
    \vcenter{\hbox{
      \xymatrix@C=35pt{
        \omega^{\pr*}(X)
        \ar[d]^-{\omega^{\pr*}(\eta_X)}
        \ar@/_3pc/@<-3ex>[ddd]_-{\vartheta_X}
        \ar@{=}[r]
        &\omega^{\pr*}(X)
        \ar[dd]^-{\xi_X}
        \\
        \omega^{\pr*}\omega_*\omega^*(X)
        \ar[dd]^-{(\vartheta_*)_{\omega^*(X)}}
        \ar@/^0.5pc/[dr]^-{(\tahar{\omega^{\pr*}\omega_*}_{\omega^*(X)})^{-1}}_-{\cong}
        &\text{ }
        \\
        \text{ }
        &\omega^*(X)\otimes \omega^{\pr*}\omega_*(\kappa)
        \ar[d]^-{I_{\omega^*(X)}\otimes g}
        \\
        \omega^*(X)\otimes b
        \ar@{=}[r]
        &\omega^*(X)\otimes b
      }
    }}
  \end{equation*}
  This completes the proof of statement 1.
  We mention here that under the bijection (\ref{eq1 LSKtensoradj HomRep,Hom=Isom}),
  each monoidal $\kos{K}$-tensor natural transformation
  $\vartheta:
  (\kos{b}^*\omega^{\pr*},\what{\kos{b}^*\omega^{\pr*}})
  \Rightarrow
  (\kos{b}^*\omega^*,\what{\kos{b}^*\omega^*})$
  corresponds to the following morphism in $\cat{Comm}(\kos{K})$.
  \begin{equation*}
    g:
    \xymatrix{
      \omega^{\pr*}\omega_*(\kappa)
      \ar[r]^-{\vartheta_{\omega_*(\kappa)}}
      \ar@/_1pc/@<-1ex>[rr]|-{(\vartheta_*)_{\kappa}}
      &\omega^*\omega_*(\kappa)\otimes b
      \ar[r]^-{\epsilon_{\kappa}\otimes I_b}
      &\kappa\otimes b
      \ar[r]^-{\imath_b^{-1}}_-{\cong}
      &b
    }
  \end{equation*}

  We are left to prove statement 2,
  as it implies statment 3.
  One can check that the morphism
  $\kos{p}_*(\xi_X)$ in $\CK$ 
  is explicitly given as follows.
  \begin{equation*}\label{eq2 LSKtensoradj HomRep,Hom=Isom}
    \vcenter{\hbox{
      \xymatrix@C=20pt{
        \omega^{\pr*}(X)\otimes \omega^{\pr*}\omega_*(\kappa)
        \ar[r]^-{\omega^{\pr*}_{X,\omega_*(\kappa)}}_-{\cong}
        \ar@/_1.2pc/[drr]_-{\kos{p}_*(\xi_X)}
        &\omega^{\pr*}(X\tensor \omega_*(\kappa))
        \ar[r]^-{\omega^{\pr*}(\varphi_X)}
        &\omega^{\pr*}\omega_*\omega^*(X)
        \ar[d]^-{(\tahar{\omega^{\pr*}\omega_*}_{\omega^*(X)})^{-1}}_-{\cong}
        \\
        \text{ }
        &\text{ }
        &\omega^*(X)\otimes \omega^{\pr*}\omega_*(\kappa)
      }
    }}
  \end{equation*}
  From the explicit description
  (\ref{eq2 LSKtensoradj HomRep,Hom=Isom}) of
  $\kos{p}_*(\xi_X)$ 
  and the fact that the functor $\kos{p}_*$ is conservative,
  we deduce that
  $\xi_X$ is an isomorphism in $\CK_p$
  if and only if
  $\kos{p}_*(\xi_X)$ is an isomorphism in $\CK$
  if and only if
  $\omega^{\pr*}(\varphi_X)$ is an isomorphism in $\CK$.
  This shows that statements \emph{(a)} and \emph{(b)} are equivalent.

  As statement \emph{(c)} impiles statement \emph{(b)},
  we are left to show that \emph{(a)}$\Leftrightarrow$\emph{(b)} implies \emph{(c)}.
  We leave for the readers to check that 
  the morphism
  $\check{\xi}_X:
  \omega^*(X)\to \omega^{\pr*}\omega_*\omega^{\pr*}(X)$
  in $\CK$
  satisfies the following relations.
  \begin{equation} \label{eq3 LSKtensoradj HomRep,Hom=Isom}
    \vcenter{\hbox{
      \xymatrix@C=20pt{
        \omega^*(X)
        \ar[r]^-{\check{\xi}_X}
        \ar@/_1.5pc/[ddr]_-{\upsilon^{\omega^{\pr*}\omega_*}_{\omega^*(X)}}
        &\omega^{\pr*}\omega_*\omega^{\pr*}(X)
        \ar[d]^-{\omega^{\pr*}\omega_*\omega^{\pr*}(\eta_X)}
        \\
        \text{ }
        &\omega^{\pr*}\omega_*\omega^{\pr*}\omega_*\omega^*(X)
        \ar[d]^-{\mu^{\omega^{\pr*}\omega_*}_{\omega^*(X)}}
        \\
        \text{ }
        &\omega^{\pr*}\omega_*\omega^*(X)
      }
    }}
    \quad
    \vcenter{\hbox{
      \xymatrix@C=20pt{
        \omega^{\pr*}(X)
        \ar[r]^-{\omega^{\pr*}(\eta_X)}
        \ar@/_1.5pc/[ddr]_-{\upsilon^{\omega^{\pr*}\omega_*}_{\omega^{\pr*}(X)}}
        &\omega^{\pr*}\omega_*\omega^*(X)
        \ar[d]^-{\omega^{\pr*}\omega_*(\check{\xi}_X)}
        \\
        \text{ }
        &\omega^{\pr*}\omega_*\omega^{\pr*}\omega_*\omega^{\pr*}(X)
        \ar[d]^-{\mu^{\omega^{\pr*}\omega_*}_{\omega^{\pr*}(X)}}
        \\
        \text{ }
        &\omega^{\pr*}\omega_*\omega^{\pr*}(X)
      }
    }}
  \end{equation}
  Using these relations,
  we are going to show that $\vartheta_X$ and $\vartheta_X^{-1}$ are inverse to each other.
  We can describe both $\vartheta_X$ and $\vartheta_X^{-1}$
  using $\vartheta_*$ as follows.
  \begin{equation*}
    \begin{aligned}
      \vartheta_X
      &:
      \xymatrix@C=30pt{
        \omega^{\pr*}(X)
        \ar[r]^-{\omega^{\pr*}(\eta_X)}
        &\omega^{\pr*}\omega_*\omega^*(X)
        \ar[r]^-{(\vartheta_*)_{\omega^*(X)}}
        &\omega^*(X)\otimes b
      }
      \\
      \vartheta^{-1}_X
      &:
      \xymatrix@C=30pt{
        \omega^*(X)
        \ar[r]^-{\check{\xi}_X}
        &\omega^{\pr*}\omega_*\omega^{\pr*}(X)
        \ar[r]^-{(\vartheta_*)_{\omega^{\pr*}(X)}}
        &\omega^{\pr*}(X)\otimes b
      }   
    \end{aligned}
  \end{equation*}
  We have
  \begin{equation*}
    \vcenter{\hbox{
      \xymatrix@C=25pt{
        \omega^{\pr*}(X)
        \ar[dd]_-{\vartheta_X}
        \ar@{=}[r]
        &\omega^{\pr*}(X)
        \ar[d]^-{\omega^{\pr*}(\eta_X)}
        \ar@{=}[r]
        &\omega^{\pr*}(X)
        \ar[ddd]^-{\upsilon^{\omega^{\pr*}\omega_*}_{\omega^{\pr*}(X)}}
        \ar@/^3pc/@<4ex>[ddddd]^-{\upsilon^{\otimes b}_{\omega^{\pr*}(X)}}
        \\
        \text{ }
        &\omega^{\pr*}\omega_*\omega^*(X)
        \ar@/_0.5pc/[dl]|-{(\vartheta_*)_{\omega^*(X)}}
        \ar[d]^-{\omega^{\pr*}\omega_*(\check{\xi}_X)}
        &\text{ }
        \\
        \omega^*(X)\otimes b
        \ar[dd]_-{\vartheta_X^{-1}\otimes I_b}
        \ar@/_0.5pc/[dr]|-{\check{\xi}_X\otimes I_b}
        &\omega^{\pr*}\omega_*\omega^{\pr*}\omega_*\omega^{\pr*}(X)
        \ar[d]^-{(\vartheta_*)_{\omega^{\pr*}\omega_*\omega^{\pr*}(X)}}
        \ar@/^1pc/[dr]^-{\mu^{\omega^{\pr*}\omega_*}_{\omega^{\pr*}(X)}}
        &\text{ }
        \\
        \text{ }
        &\omega^{\pr*}\omega_*\omega^{\pr*}(X)\otimes b
        \ar@/^0.5pc/[dl]|-{(\vartheta_*)_{\omega^{\pr*}(X)}\otimes I_b}
        &\omega^{\pr*}\omega_*\omega^{\pr*}(X)
        \ar[dd]^-{(\vartheta_*)_{\omega^{\pr*}(X)}}
        \\
        (\omega^{\pr*}(X)\otimes b)\otimes b
        \ar[d]_-{\mu^{\otimes b}_{\omega^{\pr*}(X)}}
        &\text{ }
        &\text{ }
        \\
        \omega^{\pr*}(X)\otimes b
        \ar@{=}[rr]
        &\text{ }
        &\omega^{\pr*}(X)\otimes b
      }
    }}
  \end{equation*}
  as well as
  \begin{equation*}
    \vcenter
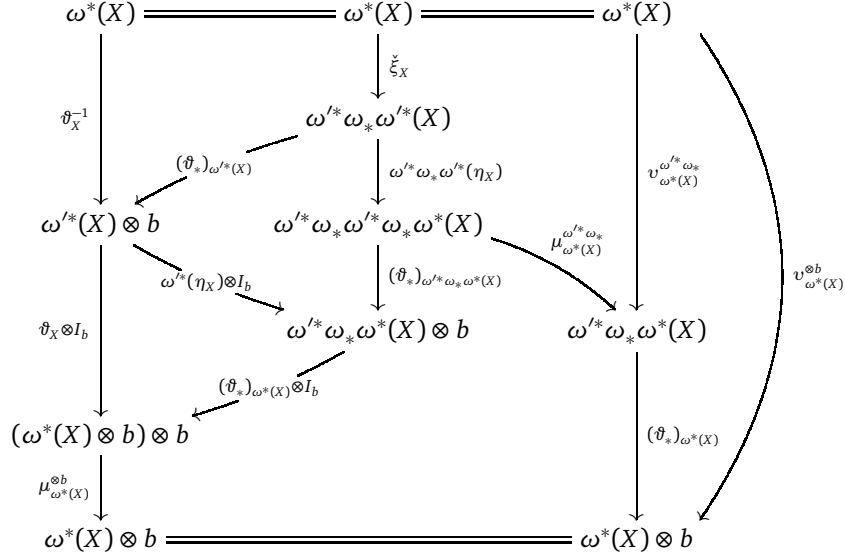
{\hbox{
      \xymatrix@C=25pt{
        \omega^*(X)
        \ar[dd]_-{\vartheta^{-1}_X}
        \ar@{=}[r]
        &\omega^*(X)
        \ar[d]^-{\check{\xi}_X}
        \ar@{=}[r]
        &\omega^*(X)
        \ar[ddd]^-{\upsilon^{\omega^{\pr*}\omega_*}_{\omega^*(X)}}
        \ar@/^3pc/@<4ex>[ddddd]^-{\upsilon^{\otimes b}_{\omega^*(X)}}
        \\
        \text{ }
        &\omega^{\pr*}\omega_*\omega^{\pr*}(X)
        \ar@/_0.5pc/[dl]|-{(\vartheta_*)_{\omega^{\pr*}(X)}}
        \ar[d]^-{\omega^{\pr*}\omega_*\omega^{\pr*}(\eta_X)}
        &\text{ }
        \\
        \omega^{\pr*}(X)\otimes b
        \ar[dd]_-{\vartheta_X\otimes I_b}
        \ar@/_0.5pc/[dr]|-{\omega^{\pr*}(\eta_X)\otimes I_b}
        &\omega^{\pr*}\omega_*\omega^{\pr*}\omega_*\omega^*(X)
        \ar[d]^-{(\vartheta_*)_{\omega^{\pr*}\omega_*\omega^*(X)}}
        \ar@/^1pc/[dr]^-{\mu^{\omega^{\pr*}\omega_*}_{\omega^*(X)}}
        &\text{ }
        \\
        \text{ }
        &\omega^{\pr*}\omega_*\omega^*(X)\otimes b
        \ar@/^0.5pc/[dl]|-{(\vartheta_*)_{\omega^*(X)}\otimes I_b}
        &\omega^{\pr*}\omega_*\omega^*(X)
        \ar[dd]^-{(\vartheta_*)_{\omega^*(X)}}
        \\
        (\omega^*(X)\otimes b)\otimes b
        \ar[d]_-{\mu^{\otimes b}_{\omega^*(X)}}
        &\text{ }
        &\text{ }
        \\
        \omega^*(X)\otimes b
        \ar@{=}[rr]
        &\text{ }
        &\omega^*(X)\otimes b
      }
    }}
  \end{equation*}
  where we used both relations
  in (\ref{eq3 LSKtensoradj HomRep,Hom=Isom}).
  This shows that the morphisms
  $\vartheta_X$ and $\vartheta_X^{-1}$
  are inverse to each other in the Kleisli category
  $\CK_b$.
  This completes the proof of Proposition~\ref{prop RSKtensoradj HomRep,Hom=Isom}.
\qed\end{proof}

\begin{lemma}\label{lem LSKtensoradj antipode}
  Let $(\kos{T},\mon{t})$
  be a strong $\kos{K}$-tensor category
  and let 
  $(\omega,\what{\omega})$,
  $(\omega^{\pr},\what{\omega}^{\pr})
  :(\kos{K},\id_{\kos{K}})\to (\kos{T},\mon{t})$
  be coreflective left-strong $\kos{K}$-tensor adjunctions.
  \begin{equation*}
    \vcenter{\hbox{
      \xymatrix{
        (\kos{T},\mon{t})
        \ar@/_1pc/[d]_-{(\omega^{\pr*},\what{\omega}^{\pr*})}
        \\
        (\kos{K},\id_{\kos{K}})
        \ar@/_1pc/[u]_-{(\omega^{\pr}_*,\what{\omega}^{\pr}_*)\text{ coreflective}}
      }
    }}
    \qquad
    \vcenter{\hbox{
      \xymatrix{
        (\kos{T},\mon{t})
        \ar@/_1pc/[d]_-{(\omega^*,\what{\omega}^*)}
        \\
        (\kos{K},\id_{\kos{K}})
        \ar@/_1pc/[u]_-{(\omega_*,\what{\omega}_*)\text{ coreflective}}
      }
    }}
  \end{equation*}
  Assume that the natural transformations
  \begin{equation*}
    \omega^{\pr*}(\varphi_{\slot}):
    \omega^{\pr*}(\slot\tensor \omega_*(\kappa))
    \Rightarrow
    \omega^{\pr*}\omega_*\omega^*
    ,
    \quad
    \omega^*(\varphi^{\pr}_{\slot}):
    \omega^*(\slot\tensor \omega^{\pr}_*(\kappa))
    \Rightarrow
    \omega^*\omega^{\pr}_*\omega^{\pr*}
  \end{equation*}
  are natural isomorphisms.
  \begin{enumerate}
    \item 
    We have a monoidal $\kos{K}$-tensor natural isomorphism
    \begin{equation*}
      \varsigma^{\omega^{\pr},\omega}:
      \xymatrix@C=18pt{
        (\omega^*\omega^{\pr}_*,\what{\omega^*\omega^{\pr}_*})
        \ar@2{->}[r]^-{\cong}
        &(\omega^{\pr*}\omega_*,\what{\omega^{\pr*}\omega_*})
        :(\kos{K},\id_{\kos{K}})\to (\kos{K},\id_{\kos{K}})
      }
    \end{equation*}
    whose component at each object $x$ in $\CK$ is
    \begin{equation*}
      \varsigma^{\omega^{\pr},\omega}_x:
      \xymatrix@C=30pt{
        \omega^*\omega^{\pr}_*(x)
        \ar[r]^-{\check{\xi}_{\omega^{\pr}_*(x)}}
        &\omega^{\pr*}\omega_*\omega^{\pr*}\omega^{\pr}_*(x)
        \ar[r]^-{\omega^{\pr*}\omega_*(\epsilon^{\pr}_x)}
        &\omega^{\pr*}\omega_*(x)
        .
      }
    \end{equation*}
    
    \item
    The following diagram of presheaves on $\cat{Aff}(\kos{K})$ strictly commutes.
    \begin{equation*}
      \hspace{-0.5cm}
      \xymatrix@C=18pt{
        \Hom_{\cat{Aff}(\kos{K})}\bigl(\slot, \cat{Spec}(\omega^{\pr*}\omega_*(\kappa))\bigr)
        \ar@2{->}[r]^-{\cong}
        \ar@2{->}[d]_-{(\varsigma^{\omega^{\pr},\omega}_{\kappa})\circ \text{ }\slot}^-{\cong}
        &\underline{\Isom}_{\mathbb{SMC}_{\cat{lax}}^{\kos{K}\backslash\!\!\backslash}}
        \bigl(
          (\omega^{\pr*},\what{\omega}^{\pr*})
          ,
          (\omega^*,\what{\omega}^*)
        \bigr)
        \ar@2{->}[d]^-{\text{taking inverse}}_-{\cong}
        &\vartheta
        \ar@{|->}[d]^-{\cong}
        \\
        \Hom_{\cat{Aff}(\kos{K})}\bigl(\slot, \cat{Spec}(\omega^*\omega^{\pr}_*(\kappa))\bigr)
        \ar@2{->}[r]^-{\cong}
        &\underline{\Isom}_{\mathbb{SMC}_{\cat{lax}}^{\kos{K}\backslash\!\!\backslash}}
        \bigl(
          (\omega^*,\what{\omega}^*)
          ,
          (\omega^{\pr*},\what{\omega}^{\pr*})
        \bigr)
        &\vartheta^{-1}
      }
    \end{equation*}
  \end{enumerate}
\end{lemma}
\begin{proof}
  Let us denote 
  $(\phi,\what{\phi})=(\omega^{\pr*}\omega_*,\what{\omega^{\pr*}\omega_*})
  :(\kos{K},\id_{\kos{K}})\to (\kos{K},\id_{\kos{K}})$
  and $p=\omega^{\pr*}\omega_*(\kappa)$.
  By Proposition~\ref{prop LSKtensoradj HomRep,Hom=Isom},
  the universal element
  \begin{equation*}
    \xymatrix{
      \xi:
      (\kos{p}^*\omega^{\pr*},\what{\kos{p}^*\omega^{\pr*}})
      \ar@2{->}[r]^-{\cong}
      &(\kos{p}^*\omega^*,\what{\kos{p}^*\omega^*})
      :(\kos{T},\mon{t})\to (\kos{K}_p,\kos{p}^*)
    }
  \end{equation*}
  is a monoidal $\kos{K}$-tensor natural isomorphism.
  Therefore
  the following are also monoidal $\kos{K}$-tensor natural transformations.
  \begin{equation*}
    \begin{aligned}
      \xi
      &:\!
      \xymatrix@C=20pt{
        (\omega^{\pr*},\what{\omega}^{\pr*})
        \ar@2{->}[r]^-{\upsilon^{\otimes p}\omega^{\pr*}}
        &(\omega^{\pr*}\otimes p,\what{\omega^{\pr*}\otimes p})
        \ar@2{->}[r]^-{\kos{p}_*\xi}_-{\cong}
        &(\omega^*\otimes p,\what{\omega^*\otimes p})
        :(\kos{T},\mon{t})\to (\kos{K},\id_{\kos{K}})
      }
      \\
      \xi^{-1}
      &:\!
      \xymatrix@C=20pt{
        (\omega^*,\what{\omega}^*)
        \ar@2{->}[r]^-{\upsilon^{\otimes p}\omega^*}
        &(\omega^*\otimes p,\what{\omega^*\otimes p})
        \ar@2{->}[r]^-{(\kos{p}_*\xi)^{-1}}_-{\cong}
        &(\omega^{\pr*}\otimes p,\what{\omega^{\pr*}\otimes p})
        :(\kos{T},\mon{t})\to (\kos{K},\id_{\kos{K}})
      }
      \\
      \check{\xi}
      &:\!
      \xymatrix@C=14pt{
        (\omega^*,\what{\omega}^*)
        \ar@2{->}[r]^-{\xi^{-1}}
        &(\omega^{\pr*}\otimes p,\what{\omega^{\pr*}\otimes p})
        \ar@2{->}[r]^-{\tahar{\phi}\omega^{\pr*}}_-{\cong}
        &(\omega^{\pr*}\omega_*\omega^{\pr*},\what{\omega^{\pr*}\omega_*\omega^{\pr*}})
        :(\kos{T},\mon{t})\to (\kos{K},\id_{\kos{K}})
      }
    \end{aligned}
  \end{equation*}
  We conclude that
  \begin{equation*}
    \xymatrix@C=50pt{
      \varsigma^{\omega^{\pr},\omega}:
      (\omega^*\omega^{\pr}_*,\what{\omega^*\omega^{\pr}_*})
      \ar@2{->}[r]^-{(\omega^{\pr*}\omega_*\epsilon^{\pr})\circ (\check{\xi}\omega^{\pr}_*)}
      &(\omega^{\pr*}\omega_*,\what{\omega^{\pr*}\omega_*})
      :(\kos{K},\id_{\kos{K}})\to (\kos{K},\id_{\kos{K}})
    }
  \end{equation*}
  is a monoidal $\kos{K}$-tensor natural transformation.
  To check that the diagram of presheaves commutes,
  it suffices to show that
  $\cat{Spec}(\omega^{\pr*}\omega_*(\kappa))
  \xrightarrow{\varsigma^{\omega^{\pr},\omega}_{\kappa}}
  \cat{Spec}(\omega^*\omega^{\pr}_*(\kappa))$
  is the unique morphism 
  in $\cat{Aff}(\kos{K})$ corresponding to 
  the inverse $\xi^{-1}$ of the universal element.
  We can check this as follows.
  \begin{equation*}
    \vcenter{\hbox{
      \xymatrix@C=40pt{
        \omega^*\omega^{\pr}_*(\kappa)
        \ar[dd]_-{\xi^{-1}_{\omega^{\pr}_*(\kappa)}}
        \ar@{=}[r]
        &\omega^*\omega^{\pr}_*(\kappa)
        \ar[d]^-{\check{\xi}_{\omega^{\pr}_*(\kappa)}}
        \ar@/^2pc/@<5ex>[dddd]^-{\varsigma^{\omega^{\pr},\omega}_{\kappa}}
        \\
        \text{ }
        &\omega^{\pr*}\omega_*\omega^{\pr*}\omega^{\pr}_*(\kappa)
        \ar[ddd]^-{\omega^{\pr*}\omega_*(\epsilon^{\pr}_{\kappa})}
        \ar@/^0.5pc/[dl]_-{(\tahar{\omega^{\pr*}\omega_*}_{\omega^{\pr*}\omega^{\pr}_*(\kappa)})^{-1}}^-{\cong}
        \\
        \omega^{\pr*}\omega^{\pr}_*(\kappa)\otimes \omega^{\pr*}\omega_*(\kappa)
        \ar[d]_-{\epsilon^{\pr}_{\kappa}\otimes I_{\omega^{\pr*}\omega_*(\kappa)}}
        &\text{ }
        \\
        \kappa\otimes \omega^{\pr*}\omega_*(\kappa)
        \ar[d]_-{\imath_{\omega^{\pr*}\omega_*(\kappa)}^{-1}}^-{\cong}
        &\text{ }
        \\
        \omega^{\pr*}\omega_*(\kappa)
        \ar@{=}[r]
        &\omega^{\pr*}\omega_*(\kappa)
      }
    }}
  \end{equation*}
  This shows that the diagram of presheaves commutes.
  As a consequence,
  we obtain that
  $\varsigma^{\omega^{\pr},\omega}:
  (\omega^*\omega^{\pr}_*,\what{\omega^*\omega^{\pr}_*})
  \Rightarrow
  (\omega^{\pr*}\omega_*,\what{\omega^{\pr*}\omega_*})$
  is a monoidal $\kos{K}$-tensor natural isomorphism,
  whose inverse is
  $\varsigma^{\omega,\omega^{\pr}}$.
  This completes the proof of Lemma~\ref{lem LSKtensoradj antipode}.
\qed\end{proof}

\begin{lemma}\label{lem LSKtensoradj composition}
  Let $(\kos{T},\mon{t})$
  be a strong $\kos{K}$-tensor category.
  Suppose we are given coreflective left-strong $\kos{K}$-tensor adjunctions
  $(\omega,\what{\omega})$,
  $(\omega^{\pr},\what{\omega}^{\pr})
  :(\kos{K},\id_{\kos{K}})\to (\kos{T},\mon{t})$
  and a strong $\kos{K}$-tensor functor
  $(\omega^{\ppr*},\what{\omega}^{\ppr*})
  :(\kos{T},\mon{t})\to (\kos{K},\id_{\kos{K}})$.
  \begin{equation*}
    \vcenter{\hbox{
      \xymatrix{
        (\kos{T},\mon{t})
        \ar[d]_-{(\omega^{\ppr*},\what{\omega}^{\ppr*})}
        \\
        (\kos{K},\id_{\kos{K}})
      }
    }}
    \quad
    \vcenter{\hbox{
      \xymatrix{
        (\kos{T},\mon{t})
        \ar@/_1pc/[d]_-{(\omega^{\pr*},\what{\omega}^{\pr*})}
        \\
        (\kos{K},\id_{\kos{K}})
        \ar@/_1pc/[u]_-{(\omega^{\pr}_*,\what{\omega}^{\pr}_*)\text{ coreflective}}
      }
    }}
    \quad
    \vcenter{\hbox{
      \xymatrix{
        (\kos{T},\mon{t})
        \ar@/_1pc/[d]_-{(\omega^*,\what{\omega}^*)}
        \\
        (\kos{K},\id_{\kos{K}})
        \ar@/_1pc/[u]_-{(\omega_*,\what{\omega}_*)\text{ coreflective}}
      }
    }}
  \end{equation*}
  Then the morphism of presheaves
  \begin{equation*}
    \begin{aligned}
      \underline{\Hom}_{\mathbb{SMC}_{\cat{lax}}^{\kos{K}\backslash\!\!\backslash}}
      \bigl(
        (\omega^{\pr*},\what{\omega}^{\pr*})
        ,&
        (\omega^*,\what{\omega}^*)
      \bigr)
      \text{ }\!\times\text{ }
      \underline{\Hom}_{\mathbb{SMC}_{\cat{lax}}^{\kos{K}\backslash\!\!\backslash}}
      \bigl(
        (\omega^{\ppr*},\what{\omega}^{\ppr*})
        ,
        (\omega^{\pr*},\what{\omega}^{\pr*})
      \bigr)
      \\
      &\Rightarrow\text{ }\!
      \underline{\Hom}_{\mathbb{SMC}_{\cat{lax}}^{\kos{K}\backslash\!\!\backslash}}
      \bigl(
        (\omega^{\ppr*},\what{\omega}^{\ppr*})
        ,
        (\omega^*,\what{\omega}^*)
      \bigr)
      :\cat{Aff}(\kos{K})^{\op}\to\cat{Set}
    \end{aligned}
  \end{equation*}
  obtained by composing monoidal $\kos{K}$-tensor natural transformations
  corresponds to the morphism 
  \begin{equation*}
    \cat{Spec}(\omega^{\pr*}\omega_*(\kappa))
    \times
    \cat{Spec}(\omega^{\ppr*}\omega^{\pr}_*(\kappa))
    \to
    \cat{Spec}(\omega^{\ppr*}\omega_*(\kappa))
  \end{equation*}
  between the representing objects in $\cat{Aff}(\kos{K})$,
  whose opposite morphism in $\cat{Comm}(\kos{K})$ is
  \begin{equation*}
    \xymatrix@C=50pt{
      \omega^{\ppr*}\omega_*(\kappa)
      \ar[r]^-{\omega^{\ppr*}(\eta^{\pr}_{\omega_*(\kappa)})}
      &\omega^{\ppr*}\omega^{\pr}_*\omega^{\pr*}\omega_*(\kappa)
      \ar[r]^-{(\tahar{\omega^{\ppr*}\omega^{\pr}_*}_{\omega^{\pr*}\omega_*(\kappa)})^{-1}}_-{\cong}
      &\omega^{\pr*}\omega_*(\kappa)
      \otimes
      \omega^{\ppr*}\omega^{\pr}_*(\kappa)
      .
    }
  \end{equation*}
\end{lemma}
\begin{proof}
  Let $\cat{Spec}(b)$ be an object in $\cat{Aff}(\kos{K})$.
  Suppose we are given monoidal $\kos{K}$-tensor natural transformations
  \begin{equation*}
    \vcenter{\hbox{
      \xymatrix{
        (\kos{b}^*\omega^{\ppr*},\what{\kos{b}^*\omega^{\ppr*}})
        \ar@2{->}[r]^-{\vartheta^{\pr}}
        &(\kos{b}^*\omega^{\pr*},\what{\kos{b}^*\omega^{\pr*}})
        \ar@2{->}[r]^-{\vartheta}
        &(\kos{b}^*\omega^*,\what{\kos{b}^*\omega^*})
      }
    }}
  \end{equation*}
  and consider their corresponding mates
  $\vartheta_*:\omega^{\pr*}\omega_*\Rightarrow \otimes b$
  and
  $\vartheta^{\pr}_*:\omega^{\ppr*}\omega^{\pr}_*\Rightarrow\otimes b$.
  Then the corresponding mate of the composition
  $\tilde{\vartheta}=\vartheta\circ\vartheta^{\pr}:
  (\kos{b}^*\omega^{\ppr*},\what{\kos{b}^*\omega^{\ppr*}})
  \Rightarrow
  (\kos{b}^*\omega^*,\what{\kos{b}^*\omega^*})$
  is given by
  \begin{equation*}
    \vcenter{\hbox{
      \xymatrix{
        \tilde{\vartheta}_*:
        \omega^{\ppr*}\omega_*
        \ar@2{->}[r]^-{\omega^{\ppr*}\eta^{\pr}\omega_*}
        &\omega^{\ppr*}\omega^{\pr}_*\omega^{\pr*}\omega_*
        \ar@2{->}[r]^-{\vartheta^{\pr}_*\vartheta_*}
        &(\otimes b)\otimes b
        \ar@2{->}[r]^-{\mu^{\otimes b}}
        &\otimes b
      }
    }}
  \end{equation*}
  as we can see from the diagram below.
  Let $x$ be an object in $\CK$.
  \begin{equation*}
    \vcenter{\hbox{
      \xymatrix@C=60pt{
        \omega^{\ppr*}\omega_*(x)
        \ar[d]_-{\omega^{\ppr*}(\eta^{\pr}_{\omega_*(x)})}
        \ar@{=}[r]
        &\omega^{\ppr*}\omega_*(x)
        \ar[dd]^-{\vartheta^{\pr}_{\omega_*(x)}}
        \ar@{=}[r]
        &\omega^{\ppr*}\omega_*(x)
        \ar@/^3pc/[ddddl]_-{\tilde{\vartheta}_{\omega^*(x)}}
        \ar[ddddd]^-{(\tilde{\vartheta}_*)_x}
        \\
        \omega^{\ppr*}\omega^{\pr}_*\omega^{\pr*}\omega_*(x)
        \ar[d]_-{(\vartheta^{\pr}_*)_{\omega^{\pr*}\omega_*(x)}}
        &\text{ }
        &\text{ }
        \\
        \omega^{\pr*}\omega_*(x)\otimes b
        \ar[dd]_-{(\vartheta_*)_x\otimes I_b}
        \ar@{=}[r]
        &\omega^{\pr*}\omega_*(x)\otimes b
        \ar[d]^-{\vartheta_{\omega_*(X)}\otimes I_b}
        &\text{ }
        \\
        \text{ }
        &(\omega^*\omega_*(x)\otimes b)\otimes b
        \ar[d]^-{\mu^{\otimes b}_{\omega^*\omega_*(x)}}
        \ar@/^0.5pc/[dl]|-{(\epsilon_x\otimes I_b)\otimes I_b}
        &\text{ }
        \\
        (x\otimes b)\otimes b
        \ar[d]_-{\mu^{\otimes b}_x}
        &\omega^*\omega_*(x)\otimes b
        \ar[d]^-{\epsilon_x\otimes I_b}
        &\text{ }
        \\
        x\otimes b
        \ar@{=}[r]
        &x\otimes b
        \ar@{=}[r]
        &x\otimes b
      }
    }}
  \end{equation*}
  Let us denote 
  $\omega^{\pr*}\omega_*(\kappa)\xrightarrow{g}b$,
  $\omega^{\ppr*}\omega^{\pr}_*(\kappa)\xrightarrow{g^{\pr}}b$,
  $\omega^{\ppr*}\omega_*(\kappa)\xrightarrow{\tilde{g}}b$
  as the morphisms in $\cat{Comm}(\kos{K})$
  corresponding to
  $\vartheta$, $\vartheta^{\pr}$, $\tilde{\vartheta}$.
  Then we have the following relation.
  \begin{equation*}
    \vcenter{\hbox{
      \xymatrix@C=20pt{
        \omega^{\ppr*}\omega_*(\kappa)
        \ar[d]_-{\omega^{\ppr*}(\eta^{\pr}_{\omega_*(\kappa)})}
        \ar@{=}[rr]
        &\text{ }
        &\omega^{\ppr*}\omega_*(\kappa)
        \ar[ddd]^-{(\tilde{\vartheta}_*)_{\kappa}}
        \ar@{=}[r]
        &\omega^{\ppr*}\omega_*(\kappa)
        \ar[ddddd]^-{\tilde{g}}
        \\
        \omega^{\ppr*}\omega^{\pr}_*\omega^{\pr*}\omega_*(\kappa)
        \ar[d]_-{(\tahar{\omega^{\ppr*}\omega^{\pr}_*}_{\omega^{\pr*}\omega_*(\kappa)})}^-{\cong}
        \ar@/^0.5pc/[dr]^-{(\vartheta^{\pr}_*\vartheta_*)_{\kappa}}
        &\text{ }
        &\text{ }
        \\
        \omega^{\pr*}\omega_*(\kappa)\otimes \omega^{\ppr*}\omega^{\pr}_*(\kappa)
        \ar[dd]_-{g\otimes g^{\pr}}
        \ar@/^0.5pc/[dr]|-{(\vartheta_*)_{\kappa}\otimes (\vartheta^{\pr}_*)_{\kappa}}
        &(\kappa\otimes b)\otimes b
        \ar[d]^-{(\tahar{\otimes b}_{\kappa\otimes b})^{-1}}_-{\cong}
        \ar@/^0.5pc/[dr]^-{\mu^{\otimes b}_{\kappa}}
        &\text{ }
        \\
        \text{ }
        &(\kappa\otimes b)\otimes (\kappa\otimes b)
        \ar@/^0.5pc/[dl]^-{\imath_b^{-1}\otimes \imath_b^{-1}}_-{\cong}
        &\kappa\otimes b 
        \ar[dd]^-{\imath_b^{-1}}_-{\cong}
        \\
        b\otimes b
        \ar[d]_-{\pc_b}
        &\text{ }
        &\text{ }
        \\
        b
        \ar@{=}[rr]
        &\text{ }
        &b
        \ar@{=}[r]
        &b
      }
    }}
  \end{equation*}
  This completes the proof of Lemma~\ref{lem LSKtensoradj composition}.
\qed\end{proof}

\subsection{Lax $\kos{K}$-tensor comonad}

\begin{definition}\label{def LaxKTcomonad}
  A \emph{lax $\kos{K}$-tensor comonad}
  is a comonad internal to the $2$-category $\mathbb{SMC}_{\cat{lax}}^{\kos{K}\backslash\!\!\backslash}$.
\end{definition}

Let us explain Definition~\ref{def LaxKTcomonad} in detail.
Let $(\kos{S},\mon{s})$ be a lax $\kos{K}$-tensor category
where
$\text{$\kos{S}$}=(\CS,\ctimes,\pzc{1})$
is the underlying symmetric monoidal category.
A lax $\kos{K}$-tensor comonad on $(\kos{S},\mon{s})$
is a tuple
\begin{equation*}
  \langle \phi,\what{\phi}\rangle
  =(\phi,\what{\phi},\delta,\epsilon)  
\end{equation*}
of a lax $\kos{K}$-tensor endofunctor
$(\phi,\what{\phi}):(\kos{S},\mon{s})\to (\kos{S},\mon{s})$
and monoidal $\kos{K}$-tensor natural transformations
$\delta:(\phi,\what{\phi})
\Rightarrow (\phi\phi,\what{\phi\phi})$,
$\epsilon:(\phi,\what{\phi})\Rightarrow (\id_{\kos{S}},I)$
which satisfy the coassociativity, counital relations.
Equivalently,
a lax $\kos{K}$-tensor comonad on $(\kos{S},\mon{s})$
is a comonoid object in the monoidal category
\begin{equation*}
  \kos{E\!n\!\!d}_{\cat{lax}}^{\kos{K}\backslash\!\!\backslash}(\kos{S},\mon{s})
  =
  \Bigl(
    \mathbb{SMC}_{\cat{lax}}^{\kos{K}\backslash\!\!\backslash}\bigl((\kos{S},\mon{s}),(\kos{S},\mon{s})\bigr)
    ,
    \circ
    ,
    (\id_{\kos{S}},I)
  \Bigr)
\end{equation*}
of lax $\kos{K}$-tensor endofunctors on $(\kos{S},\mon{s})$.

\begin{definition}
  \label{def LaxKTcomonad HopfKTcomonad}
  Let $(\kos{S},\mon{s})$ be a lax $\kos{K}$-tensor category.
  \begin{enumerate}
    \item 
    The \emph{fusion operator} associated to
    a lax $\kos{K}$-tensor comonad
    $\langle\phi,\what{\phi}\rangle$ on $(\kos{S},\mon{s})$
    is the natural transformation
    $\chi_{\slot,\slot}:\phi(\slot)\ctimes\phi(\slot)\Rightarrow\phi(\phi(\slot)\ctimes \slot):\CS\times\CS\to \CS$
    whose component at $\pzc{X}$, $\pzc{Y}\in\obj{\CS}$ is
    \begin{equation*}
      \chi_{\pzc{X},\pzc{Y}}:
      \xymatrix@C=35pt{
        \phi(\pzc{X})\ctimes \phi(\pzc{Y})
        \ar[r]^-{\delta_{\pzc{X}}\ctimes I_{\phi(\pzc{Y})}}
        &\phi\phi(\pzc{X})\ctimes \phi(\pzc{Y})
        \ar[r]^-{\phi_{\phi(\pzc{X}),\pzc{Y}}}
        &\phi(\phi(\pzc{X})\ctimes \pzc{Y})
        .
      }
    \end{equation*}

    \item 
    A \emph{Hopf $\kos{K}$-tensor comonad} on $(\kos{S},\mon{s})$
    is a lax $\kos{K}$-tensor comonad $\langle\phi,\what{\phi}\rangle$ on $(\kos{S},\mon{s})$
    whose associated fusion operator $\chi_{\slot,\slot}$ is a natural isomorphism.
  \end{enumerate}
\end{definition}

Let $\langle \phi,\what{\phi}\rangle=(\phi,\what{\phi},\delta,\epsilon)$
be a lax $\kos{K}$-tensor comonad on
a lax $\kos{K}$-tensor category
$(\kos{S},\mon{s})$.
We have the \emph{Eilenberg-Moore lax $\kos{K}$-tensor category}
\begin{equation}\label{eq LaxKTcomonad defE-Mcat}
  (\text{$\kos{S}$},\mon{s})^{\langle\phi,\what{\phi}\rangle}
  =(\text{$\kos{S}$}^{\langle\phi,\what{\phi}\rangle},\mon{s}^{\langle\phi,\what{\phi}\rangle})
\end{equation}
where we denote
$\text{$\kos{S}$}^{\langle\phi,\what{\phi}\rangle}
=(\CS^{\langle\phi,\what{\phi}\rangle},\ctimes\!^{\langle\phi,\what{\phi}\rangle},\pzc{1}^{\langle\phi,\what{\phi}\rangle})$
as the underlying symmetric monoidal category.
\begin{itemize}
  \item 
  The category $\CS^{\langle\phi,\what{\phi}\rangle}$
  is the Eilenberg-Moore category
  associated to the comonad $\langle\phi\rangle=(\phi,\delta,\epsilon)$ on $\CS$.
  An object in
  $\CS^{\langle\phi,\what{\phi}\rangle}$
  is a pair $X=(\pzc{X},\rho_{\pzc{X}})$
  of an object $\pzc{X}$ in $\CS$
  and a morphism $\rho_{\pzc{X}}:\pzc{X}\to\phi(\pzc{X})$ in $\CS$
  which satisfies the following relations.
  \begin{equation*}
    \vcenter
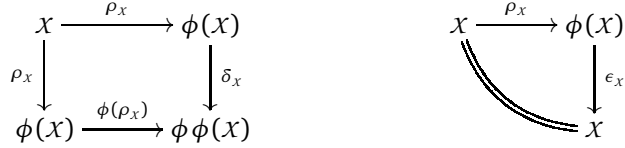
{\hbox{
      \xymatrix@C=30pt{
        \pzc{X}
        \ar[r]^-{\rho_{\pzc{X}}}
        \ar[d]_-{\rho_{\pzc{X}}}
        &\phi(\pzc{X})
        \ar[d]^-{\delta_{\pzc{X}}}
        &\text{ }
        &\pzc{X}
        \ar[r]^-{\rho_{\pzc{X}}}
        \ar@/_1pc/@{=}[dr]
        &\phi(\pzc{X})
        \ar[d]^-{\epsilon_{\pzc{X}}}
        \\
        \phi(\pzc{X})
        \ar[r]^-{\phi(\rho_{\pzc{X}})}
        &\phi\phi(\pzc{X})
        &\text{ }
        &\text{ }
        &\pzc{X}
      }
    }}
  \end{equation*} 
  Let $Y=(\pzc{Y},\gamma_{\pzc{Y}})$
  be another object in $\CS^{\langle\phi,\what{\phi}\rangle}$.
  A morphism $X\to Y$ in $\CS^{\langle\phi,\what{\phi}\rangle}$
  is a morphism $\pzc{X}\to \pzc{Y}$ in $\CS$ which is compatible with
  $\rho_{\pzc{X}}$, $\rho_{\pzc{Y}}$.

  \item
  The monoidal product 
  of objects
  $X$, $Y$
  in $\CS^{\langle\phi,\what{\phi}\rangle}$
  is the pair
  $X\ctimes\!^{\langle\phi,\what{\phi}\rangle} Y
  =(\pzc{X}\ctimes \pzc{Y},\rho_{\pzc{X}\ctimes \pzc{Y}})$
  where
  $\rho_{\pzc{X}\ctimes \pzc{Y}}:
  \pzc{X}\ctimes \pzc{Y}
  \xrightarrow{\rho_{\pzc{X}}\ctimes \rho_{\pzc{Y}}}
  \phi(\pzc{X})\ctimes \phi(\pzc{Y})
  \xrightarrow{\phi_{\pzc{X},\pzc{Y}}}
  \phi(\pzc{X}\ctimes \pzc{Y})$
  and the unit object 
  in $\CS^{\langle\phi,\what{\phi}\rangle}$
  is the pair
  $\pzc{1}^{\langle\phi,\what{\phi}\rangle}
  =(\pzc{1},\rho_{\pzc{1}})$ where
  $\rho_{\pzc{1}}:
  \pzc{1}\xrightarrow{\phi_{\pzc{1}}}\phi(\pzc{1})$.
  The symmetric monoidal coherence morphisms
  $a$, $\imath$, $\jmath$, $s$ 
  of $\kos{S}^{\langle\phi,\what{\phi}\rangle}$
  are given by those of $\kos{S}$.

  \item
  The lax symmetric monoidal functor
  $\text{$\mon{s}$}^{\langle\phi,\what{\phi}\rangle}:
  \kos{K}\to \kos{S}^{\langle\phi,\what{\phi}\rangle}$
  sends each object $x$ in $\CK$ to the pair
  $\mon{s}^{\langle\phi,\what{\phi}\rangle}(x)
  =(\mon{s}(x),\what{\phi}_x)$
  where
  $\what{\phi}_x:\mon{s}(x)\to \phi\mon{s}(x)$.
  The lax symmetric monoidal coherence morphisms
  of $\mon{s}^{\langle\phi,\what{\phi}\rangle}$ are given by those of $\mon{s}$.
  In particular,
  $(\text{$\kos{S}$},\mon{s})^{\langle\phi,\what{\phi}\rangle}
  =(\text{$\kos{S}$}^{\langle\phi,\what{\phi}\rangle},\mon{s}^{\langle\phi,\what{\phi}\rangle})$
  is a strong $\kos{K}$-tensor category
  if and only if
  $(\kos{S},\mon{s})$
  is a strong $\kos{K}$-tensor category.
\end{itemize}
We also have a left-strong $\kos{K}$-tensor adjunction
$(\omega,\what{\omega}):(\kos{S},\mon{s})\to 
(\kos{S},\mon{s})^{\langle\phi,\what{\phi}\rangle}$.
\begin{equation}\label{eq LaxKTcomonad E-McatLSKTadj}
  \vcenter{\hbox{
    \xymatrix{
      (\kos{S},\mon{s})^{\langle\phi,\what{\phi}\rangle}
      \ar@/_1pc/[d]_-{(\omega^*,\what{\omega}^*)}
      \\
      (\kos{S},\mon{s})
      \ar@/_1pc/[u]_-{(\omega_*,\what{\omega}_*)}
    }
  }}
\end{equation}
\begin{itemize}
  \item 
  The left adjoint $\omega^*$ is
  the forgetful functor
  which sends each object
  $X=(\pzc{X},\rho_{\pzc{X}})$ in $\CS^{\langle\phi,\what{\phi}\rangle}$
  to the underlying object $\omega^*(X)=\pzc{X}$ in $\CS$.
  The symmetric monoidal coherence isomorphisms of
  $\omega^*:\text{$\kos{S}$}^{\langle\phi,\what{\phi}\rangle}\to \text{$\kos{S}$}$
  are identity morphisms,
  and the monoidal natural isomorphism
  $\what{\omega}^*:\mon{s}=\text{$\omega$}^*\mon{s}^{\langle\phi,\what{\phi}\rangle}$
  is the identity natural transformation.
  
  \item
  The right adjoint $\omega_*$ sends each
  $\pzc{X}\in\obj{\CS}$
  to the object
  $\omega_*(\pzc{X})=(\phi(\pzc{X}),\delta_{\pzc{X}})$
  in $\CS^{\langle\phi,\what{\phi}\rangle}$
  where $\delta_{\pzc{X}}:\phi(\pzc{X})\to\phi\phi(\pzc{X})$.
  The lax symmetric monoidal coherence morphisms of
  $\omega_*:\kos{S}\to \kos{S}^{\langle\phi,\what{\phi}\rangle}$ are
  given by those of $\phi:\kos{S}\to\kos{S}$,
  and the monoidal natural transformation
  $\what{\omega}_*:\text{$\mon{s}$}^{\langle\phi,\what{\phi}\rangle}\Rightarrow \omega_*\mon{s}$
  is given by
  $\what{\phi}:\mon{s}\Rightarrow\phi\mon{s}$.

  \item
  The component of the adjunction unit at each object
  $X=(\pzc{X},\gamma_{\pzc{X}})$ in $\CS^{\langle\phi,\what{\phi}\rangle}$
  is $\rho_{\pzc{X}}:X\to \omega_*\omega^*(X)$.
  The component of the adjunction counit at each object
  $\pzc{X}$ in $\CS$ is
  $\epsilon_{\pzc{X}}:\omega^*\omega_*(\pzc{X})=\phi(\pzc{X})\to \pzc{X}$.
  The adjunction lax $\kos{K}$-tensor comonad on $(\kos{S},\mon{s})$
  induced from the left-strong $\kos{K}$-tensor adjunction
  $(\omega,\what{\omega}):(\kos{S},\mon{s})
  \to (\kos{S},\mon{s})^{\langle\phi,\what{\phi}\rangle}$
  is
  $\langle\omega^*\omega_*,\what{\omega^*\omega_*}\rangle
  =\langle\phi,\what{\phi}\rangle$.
\end{itemize}

Moreover, one can check that the following are true.
\begin{enumerate}
  \item 
  The $\kos{K}$-equivariances
  $\cevar{\phi}$, $\cevar{\omega}_*$
  associated to
  $(\phi,\what{\phi}):(\kos{S},\mon{s})\to (\kos{S},\mon{s})$
  and
  $(\omega_*,\what{\omega}_*):
  (\kos{S},\mon{s})\to
  (\kos{S},\mon{s})^{\langle\phi,\what{\phi}\rangle}$
  satisfy the relation
  $\cevar{\phi}=\omega^*\cevar{\omega}_*$.
  In particular,
  $\cevar{\phi}$ is a natural isomorphism
  if and only if $\cevar{\omega}_*$ is a natural isomorphism.
  
  \item
  The left-strong $\kos{K}$-tensor adjunction
  $(\omega,\what{\omega}):
  (\kos{S},\mon{s})\to
  (\kos{S},\mon{s})^{\langle\phi,\what{\phi}\rangle}$
  satisfies the projection formula
  if and only if
  for every object
  $X=(\pzc{X},\gamma_{\pzc{X}})$ in $\CS^{\langle\phi,\what{\phi}\rangle}$
  and every object $\pzc{Z}$ in $\CS$,
  \begin{equation*}
    \varphi_{X,\pzc{Z}}:
    \xymatrix@C=30pt{
      \pzc{X}\ctimes \phi(\pzc{Z})
      \ar[r]^-{\rho_{\pzc{X}}\ctimes I_{\phi(\pzc{Z})}}
      &\phi(\pzc{X})\ctimes \phi(\pzc{Z})
      \ar[r]^-{\phi_{\pzc{X},\pzc{Z}}}
      &\phi(\pzc{X}\ctimes \pzc{Z})
    }
  \end{equation*}
  is an isomorphism in $\CS$.

  \item
  The given lax $\kos{K}$-tensor comonad
  $\langle\phi,\what{\phi}\rangle=(\phi,\what{\phi},\delta,\epsilon)$
  on $(\kos{S},\mon{s})$
  is a Hopf $\kos{K}$-tensor comonad
  if and only if 
  for every pair of objects $\pzc{X}$, $\pzc{Z}$ in $\CS$,
  \begin{equation*}
    \chi_{\pzc{X},\pzc{Z}}=
    \varphi_{\omega_*(\pzc{X}),\pzc{Z}}:
    \xymatrix@C=30pt{
      \phi(\pzc{X})\ctimes \phi(\pzc{Z})
      \ar[r]^-{\delta_{\pzc{X}}\ctimes I_{\phi(\pzc{Z})}}
      &\phi\phi(\pzc{X})\ctimes \phi(\pzc{Z})
      \ar[r]^-{\phi_{\phi(\pzc{X}),\pzc{Z}}}
      &\phi(\phi(\pzc{X})\ctimes \pzc{Z})
    }
  \end{equation*}
  is an isomorphism in $\CS$.
\end{enumerate}

\subsection{Representations of a group object in $\kos{A\!f\!f}(\kos{K})$}
\label{subsec RepGrpAff(K)}

Recall the cocartesian monoidal category
$\kos{C\!o\!m\!m}(\kos{K})=(\cat{Comm}(\kos{K}),\otimes,\kappa)$
of commutative monoids in $\kos{K}$
and the opposite monoidal category
\begin{equation*}
  \kos{A\!f\!f}(\kos{K})=(\cat{Aff}(\kos{K}),\times,\cat{Spec}(\kappa))
\end{equation*}
introduced in (\ref{eq Comm(T) Comm(T)def}) and (\ref{eq Comm(T) Aff(T)def}).
We denote
\begin{equation*}
  \cat{Mon}(\kos{A\!f\!f}(\kos{K}))
\end{equation*}
as the category of monoid objects in $\kos{A\!f\!f}(\kos{K})$.
A monoid object in $\kos{A\!f\!f}(\kos{K})$
is a triple 
\begin{equation*}
  (\cat{Spec}(\pi),\cp_{\pi},e_{\pi})
\end{equation*}
where $\cat{Spec}(\pi)$ is an object in $\cat{Aff}(\kos{K})$
and 
$\cat{Spec}(\pi)\times \cat{Spec}(\pi)\xrightarrow{\cp_{\pi}} \cat{Spec}(\pi)$,
$\cat{Spec}(\kappa)\xrightarrow{e_{\pi}} \cat{Spec}(\pi)$
are product, unit morphisms
in $\cat{Aff}(\kos{K})$
satisfying the associativity, unital relations.
Equivalently,
a monoid object $(\cat{Spec}(\pi),\cp_{\pi},e_{\pi})$
in $\kos{A\!f\!f}(\kos{K})$ 
corresponds to a commutative bimonoid $(\pi,\pc_{\pi},u_{\pi},\cp_{\pi},e_{\pi})$ in $\kos{K}$
where $\pi\xrightarrow{\cp_{\pi}} \pi\otimes \pi$, $\pi\xrightarrow{e_{\pi}}\kappa$
are coproduct, counit morphisms of $\pi$.
We often denote a monoid object
$(\cat{Spec}(\pi),\cp_{\pi},e_{\pi})$
in $\kos{A\!f\!f}(\kos{K})$ as $\cat{Spec}(\pi)$.
A morphism $\cat{Spec}(\pi^{\pr})\to \cat{Spec}(\pi)$
of monoid objects in $\kos{A\!f\!f}(\kos{K})$
is a morphism in $\cat{Aff}(\kos{K})$
which is compatible with product, unit morphisms of $\cat{Spec}(\pi^{\pr})$, $\cat{Spec}(\pi)$.

Let $\cat{Spec}(\pi)$ be a monoid object in $\kos{A\!f\!f}(\kos{K})$.
We have the following representable presheaf
on $\cat{Aff}(\kos{K})$
which factors through the category $\cat{Mon}$ of monoids.
\begin{equation*}
  \Hom_{\cat{Aff}(\kos{K})}\bigl(\slot,\cat{Spec}(\pi)\bigr):
  \cat{Aff}(\kos{K})^{\op}\to\cat{Mon}
\end{equation*}
For each object $\cat{Spec}(b)$ in $\cat{Aff}(\kos{K})$,
the set
\begin{equation*}
  \Hom_{\cat{Aff}(\kos{K})}\bigl(\cat{Spec}(b),\cat{Spec}(\pi)\bigr)
  =\Hom_{\cat{Comm}(\kos{K})}(\pi,b)  
\end{equation*}
is a monoid equipped with convolution product
$(\pi\xrightarrow{f} b,\pi\xrightarrow{g}b)
\mapsto
f\star g:\pi\xrightarrow{\cp_{\pi}}\pi\otimes \pi\xrightarrow{f\otimes g}b\otimes b\xrightarrow{\pc_b}b$
and identity object
$\pi\xrightarrow{e_{\pi}}\kappa\xrightarrow{u_b}b$.
Using the Yoneda lemma,
one can check that the correspondence
$\cat{Spec}(\pi)\mapsto \Hom_{\cat{Aff}(\kos{K})}\bigl(\slot,\cat{Spec}(\pi)\bigr)$
defines an equivalence of categories
from
$\cat{Mon}(\kos{A\!f\!f}(\kos{K}))$
to the category of
representable presheaves
on $\cat{Aff}(\kos{K})$
which factors through $\cat{Mon}$.

\begin{remark} \label{rem RepGrpAff(K) tensorpi conservative}
  Let $\cat{Spec}(\pi)$ be a monoid object in $\kos{A\!f\!f}(\kos{K})$.
  Then the functor $\slot\otimes\pi:\CK\to \CK$ is conservative.
  Indeed, if $x\xrightarrow{f} y$ is a morphism in $\CK$
  such that $x\otimes \pi\xrightarrow{f\otimes I_{\pi}} y\otimes \pi$ is an isomorphism,
  then $f$ has an inverse
  $f^{-1}:
  y
  \xrightarrow[\cong]{\jmath_y}
  y\otimes \kappa
  \xrightarrow{I_y\otimes u_{\pi}}
  y\otimes \pi
  \xrightarrow[\cong]{(f\otimes I_{\pi})^{-1}}
  x\otimes \pi
  \xrightarrow{I_x\otimes e_{\pi}}
  x\otimes \kappa
  \xrightarrow[\cong]{\jmath_x^{-1}}
  x$.
\end{remark}

In (\ref{eq Comm(T) laxKTend(K,idK) def})
we introduced the monoidal category
$\kos{E\!n\!\!d}_{\cat{lax}}^{\kos{K}\backslash\!\!\backslash}(\kos{K},\id_{\kos{K}})_{\cat{crfl}}$
of coreflective lax $\kos{K}$-tensor endofunctors on $(\kos{K},\id_{\kos{K}})$.
Recall the definition of lax $\kos{K}$-tensor comonads
in Definition~\ref{def LaxKTcomonad}.

\begin{definition}
  A lax $\kos{K}$-tensor comonad
  $\langle\phi,\what{\phi}\rangle=(\phi,\what{\phi},\delta,\epsilon)$ on
  $(\kos{K},\id_{\kos{K}})$
  is called \emph{coreflective}
  if the underlying lax $\kos{K}$-tensor endofunctor
  $(\phi,\what{\phi})$ on $(\kos{K},\id_{\kos{K}})$ is coreflective.
  We denote
  \begin{equation*}
    \cat{Comonad}_{\cat{lax}}^{\kos{K}\backslash\!\!\backslash}(\kos{K},\id_{\kos{K}})_{\cat{crfl}}
  \end{equation*}
  as the category of coreflective lax $\kos{K}$-tensor comonads on $(\kos{K},\id_{\kos{K}})$,
  which is the category of comonoid objects in
  $\kos{E\!n\!\!d}_{\cat{lax}}^{\kos{K}\backslash\!\!\backslash}(\kos{K},\id_{\kos{K}})_{\cat{crfl}}$.
\end{definition}

Recall the adjoint equivalence of monoidal categories
\begin{equation*}
  \iota:
  \kos{C\!o\!m\!m}(\kos{K})
  \simeq
  \kos{E\!n\!\!d}_{\cat{lax}}^{\kos{K}\backslash\!\!\backslash}(\kos{K},\id_{\kos{K}})_{\cat{crfl}}^{\cat{rev}}
  :\CR
\end{equation*}
in Proposition~\ref{prop Comm(T) Comm(K)=KtensorEnd(K)}.
By considering the category of comonoid objects on both sides,
we obtain the following corollary.
Note that
\begin{itemize}
  \item 
  a comonoid object in $\kos{C\!o\!m\!m}(\kos{K})$
  is the same as a monoid object in $\kos{A\!f\!f}(\kos{K})$.

  \item
  a comonoid object in the reversed monoidal category
  $\kos{E\!n\!\!d}_{\cat{lax}}^{\kos{K}\backslash\!\!\backslash}(\kos{K},\id_{\kos{K}})_{\cat{crfl}}^{\cat{rev}}$
  is the same as a comonoid object in
  $\kos{E\!n\!\!d}_{\cat{lax}}^{\kos{K}\backslash\!\!\backslash}(\kos{K},\id_{\kos{K}})_{\cat{crfl}}$.
\end{itemize}

\begin{corollary}
  \label{cor RepGrpAff(K) Mon(Aff(K))=crflLKTcomonad}
  We have an adjoint equivalence of categories
  \begin{equation*}
    \iota:
    \cat{Mon}(\kos{A\!f\!f}(\kos{K}))
    \simeq
    \cat{Comonad}_{\cat{lax}}^{\kos{K}\backslash\!\!\backslash}(\kos{K},\id_{\kos{K}})_{\cat{crfl}}
    :\CR
  \end{equation*}
  \begin{itemize}
    \item 
    The left adjoint
    $\iota$
    sends each monoid object $\cat{Spec}(\pi)$ in $\kos{A\!f\!f}(\kos{K})$
    to the coreflective lax $\kos{K}$-tensor comonad
    $\iota(\cat{Spec}(\pi))
    =\langle\otimes\pi,\what{\otimes\pi}\rangle
    =(\otimes\pi,\what{\otimes\pi},\delta^{\otimes\pi},\epsilon^{\otimes\pi})$
    on $(\kos{K},\id_{\kos{K}})$.
    The components of $\delta^{\otimes\pi}$, $\epsilon^{\otimes\pi}$
    at each object $x$ in $\CK$ are described below.
    \begin{equation*}
      \hspace*{-0.5cm}
      \delta^{\otimes \pi}_x
      :\!\!
      \xymatrix@C=18pt{
        x\otimes \pi
        \ar[r]^-{I_x\otimes \cp_{\pi}}
        &x\otimes (\pi\otimes \pi)
        \ar[r]^-{a_{x,\pi,\pi}}_-{\cong}
        &(x\otimes \pi)\otimes \pi
      }
      \quad
      \epsilon^{\otimes\pi}_x
      :\!\!
      \xymatrix@C=18pt{
        x\otimes \pi
        \ar[r]^-{I_x\otimes e_{\pi}}
        &x\otimes \kappa
        \ar[r]^-{\jmath_x^{-1}}_-{\cong}
        &x
      }
    \end{equation*}
   
    \item 
    The right adjoint
    $\CR$
    sends each coreflective lax $\kos{K}$-tensor comonad
    $\langle\phi,\what{\phi}\rangle
    =(\phi,\what{\phi},\mu,\epsilon)$
    on $(\kos{K},\id_{\kos{K}})$
    to the monoid object
    $\CL\langle \phi,\what{\phi}\rangle
    =\cat{Spec}(\phi(\kappa))$ in $\kos{A\!f\!f}(\kos{K})$
    where
    \begin{equation}\label{eq RepGrpAff(K) Mon(Aff(K))=crflLKTcomonad}
      \cp_{\phi(\kappa)}
      :
      \xymatrix{
        \phi(\kappa)
        \ar[r]^-{\delta_{\kappa}}
        &\phi\phi(\kappa)
        \ar[r]^-{(\tahar{\phi}_{\phi(\kappa)})^{-1}}_-{\cong}
        &\phi(\kappa)\otimes \phi(\kappa)
        ,
      }
      \qquad
      e_{\phi(\kappa)}
      :
      \xymatrix{
        \phi(\kappa)
        \ar[r]^-{\epsilon_{\kappa}}
        &\kappa
        .
      }
    \end{equation}
  \end{itemize}
\end{corollary}

A group object in $\kos{A\!f\!f}(\kos{K})$
is a tuple $(\cat{Spec}(\pi),\cp_{\pi},e_{\pi},\varsigma_{\pi})$
where $(\cat{Spec}(\pi),\cp_{\pi},e_{\pi})$ is a monoid object in $\kos{A\!f\!f}(\kos{K})$
and
$\varsigma_{\pi}:\cat{Spec}(\pi)\to \cat{Spec}(\pi)$
is the antipode morphism,
which is a morphism in $\cat{Aff}(\kos{K})$
that satisfies the relation
$\varsigma_{\pi}\star I_{\pi}=I_{\pi}\star \varsigma_{\pi}=u_{\pi}\circ e_{\pi}:\pi\to \pi$.
Equivalently,
a group object $(\cat{Spec}(\pi),\cp_{\pi},e_{\pi},\varsigma_{\pi})$
in $\kos{A\!f\!f}(\kos{K})$ is
a commutative Hopf monoid $(\pi,\pc_{\pi},u_{\pi},\cp_{\pi},e_{\pi},\varsigma_{\pi})$ in $\kos{K}$.
One can check that the antipode morphism
$\varsigma_{\pi}$
is involutive (hence is an isomorphism),
and is an anti-morphism of monoid objects in $\kos{A\!f\!f}(\kos{K})$.
\begin{equation*}
  \vcenter{\hbox{
    \xymatrix@C=30pt{
      \pi
      \ar[r]^-{\varsigma_{\pi}}_-{\cong}
      \ar@/_1pc/@{=}[dr]
      &\pi
      \ar[d]^-{\varsigma_{\pi}}_-{\cong}
      &\pi
      \ar[d]_-{\cp_{\pi}}
      \ar[rr]^-{\varsigma_{\pi}}_-{\cong}
      &\text{ }
      &\pi
      \ar[d]^-{\cp_{\pi}}
      &\pi
      \ar[r]^-{\varsigma_{\pi}}_-{\cong}
      \ar@/_1pc/[dr]_-{e_{\pi}}
      &\pi
      \ar[d]^-{e_{\pi}}
      \\
      \text{ }
      &\pi
      &\pi\otimes \pi
      \ar[r]^-{s_{\pi,\pi}}_-{\cong}
      &\pi\otimes \pi
      \ar[r]^-{\varsigma_{\pi}\otimes \varsigma_{\pi}}_-{\cong}
      &\pi\otimes \pi
      &\text{ }
      &\kappa
    }
  }}
\end{equation*}
We also denote a group object
$(\cat{Spec}(\pi),\pc_{\pi},u_{\pi},\varsigma_{\pi})$
in $\kos{A\!f\!f}(\kos{K})$ as $\cat{Spec}(\pi)$.
A morphism
of group objects
$\cat{Spec}(\pi^{\pr})\to \cat{Spec}(\pi)$
 in $\kos{A\!f\!f}(\kos{K})$
is a morphism as monoid objects in $\kos{A\!f\!f}(\kos{K})$.

Let $\cat{Spec}(\pi)$ be a group object in $\kos{A\!f\!f}(\kos{K})$.
We have the following representable presheaf
on $\cat{Aff}(\kos{K})$
which factors through the category $\cat{Grp}$ of groups.
\begin{equation}\label{eq RepGrpAff(K) Hom(-,Spec(pi))}
  \Hom_{\cat{Aff}(\kos{K})}\bigl(\slot,\cat{Spec}(\pi)\bigr):
  \cat{Aff}(\kos{K})^{\op}\to\cat{Grp}
\end{equation}
For each object $\cat{Spec}(b)$ in $\cat{Aff}(\kos{K})$,
the set
\begin{equation*}
  \Hom_{\cat{Aff}(\kos{K})}\bigl(\cat{Spec}(b),\cat{Spec}(\pi)\bigr)
  =\Hom_{\cat{Comm}(\kos{K})}(\pi,b)
\end{equation*}
is a group
where the inverse of 
$\pi\xrightarrow{f}b$
with respect to the convolution product $\star$
is $\pi\xrightarrow[\cong]{\varsigma_{\pi}}\pi\xrightarrow{f}b$.
Using the Yoneda lemma,
one can check that
$\cat{Spec}(\pi)\mapsto \Hom_{\cat{Aff}(\kos{K})}\bigl(\slot,\cat{Spec}(\pi)\bigr)$
defines an equivalence of categories
from
the category of group objects in $\kos{A\!f\!f}(\kos{K})$
to the category of representable presheaves
on $\cat{Aff}(\kos{K})$
which factors through $\cat{Grp}$.

Recall the definition of Hopf $\kos{K}$-tensor comonads
in Definition~\ref{def LaxKTcomonad HopfKTcomonad}.

\begin{definition}
  A Hopf $\kos{K}$-tensor comonad
  $\langle\phi,\what{\phi}\rangle=(\phi,\what{\phi},\delta,\epsilon)$ on
  $(\kos{K},\id_{\kos{K}})$
  is called \emph{coreflective}
  if the underlying lax $\kos{K}$-tensor endofunctor
  $(\phi,\what{\phi})$ on $(\kos{K},\id_{\kos{K}})$ is coreflective.
\end{definition}

Each group object $\cat{Spec}(\pi)$ in $\kos{A\!f\!f}(\kos{K})$
determines a coreflective Hopf $\kos{K}$-tensor comonad
$\iota(\pi)=\langle\otimes\pi,\what{\otimes\pi}\rangle$
on $(\kos{K},\id_{\kos{K}})$.
The component of the associated fusion operator
at $x$, $y\in\obj{\CK}$ is
\begin{equation*}
  \chi_{x,y}:
  \xymatrix@C=35pt{
    x\otimes \pi\otimes y\otimes \pi
    \ar[r]^-{I_x\otimes \cp_{\pi}\otimes I_{y\otimes \pi}}
    &x\otimes \pi\otimes \pi\otimes y\otimes \pi
    \ar[r]^-{(\otimes\pi)_{x\otimes \pi,y}}
    &x\otimes \pi\otimes y\otimes \pi
  }
\end{equation*}
which is an isomorphism in $\CK$
whose inverse is described below.
\begin{equation*}
  \xymatrix@C=40pt{
    x\otimes \pi\otimes y\otimes \pi
    \ar[d]_-{\chi_{x,y}^{-1}}
    \ar[r]^-{I_x\otimes \cp_{\pi}\otimes I_{y\otimes \pi}}
    &x\otimes \pi\otimes \pi\otimes y\otimes \pi
    \ar[d]^-{I_{x\otimes \pi}\otimes \varsigma_{\pi}\otimes I_{y\otimes \pi}}_-{\cong}
    \\
    x\otimes \pi\otimes y\otimes \pi
    &x\otimes \pi\otimes \pi\otimes y\otimes \pi
    \ar[l]_-{(\otimes\pi)_{x\otimes \pi,y}}
  }
\end{equation*}

Let $\cat{Spec}(\pi)$ be a group object in $\kos{A\!f\!f}(\kos{K})$.
We define the strong $\kos{K}$-tensor category
of representations of $\cat{Spec}(\pi)$ in $\kos{K}$ as follows.
Recall the Eilenberg-Moore category 
associated to a lax $\kos{K}$-tensor comonad
introduced in (\ref{eq ColaxKTmonad defE-Mcat}).
  
\begin{definition}\label{def RepGrpAff(K) Rep(pi)Ktensorcat}
  Let $\cat{Spec}(\pi)$ be a group object in $\kos{A\!f\!f}(\kos{K})$.
  We define the strong $\kos{K}$-tensor category
  \begin{equation*}
    (\kos{R\!e\!p}(\pi),\kos{t}^*_{\pi})    
  \end{equation*}
  of representations of $\cat{Spec}(\pi)$ in $\kos{K}$
  as the Eilenberg-Moore strong $\kos{K}$-tensor category
  associated to the coreflective Hopf $\kos{K}$-tensor comonad
  $\langle \otimes\pi,\what{\otimes\pi}\rangle$ on $(\kos{K},\id_{\kos{K}})$.
  The underlying symmetric monoidal category is denoted as
  \begin{equation*}
    \kos{R\!e\!p}(\pi)=(\cat{Rep}(\pi),\tensor\!_{\pi},\unit\!_{\pi})
  \end{equation*}
  and the strong symmetric monoidal functor
  $\kos{t}_{\pi}^*:\kos{K}\to \kos{R\!e\!p}(\pi)$
  is the functor of constructing trivial representations.
\end{definition}

Let us explain Definition~\ref{def RepGrpAff(K) Rep(pi)Ktensorcat} in detail.
Let $\cat{Spec}(\pi)$ be a group object in $\kos{A\!f\!f}(\kos{K})$.
An object in
$\cat{Rep}(\pi)$
is a pair $X=(x,\rho_x)$
of an object $x$ in $\CK$
and a morphism $\rho_x:x\to x\otimes \pi$ in $\CK$
satisfying the right $\pi$-coaction relations.
Let $Y=(y,\rho_y)$ be another object in $\cat{Rep}(\pi)$.
A morphism $X\to Y$ in $\cat{Rep}(\pi)$
is a morphism $x\to y$ in $\CK$ which is compatible with
$\rho_x$, $\rho_y$.
The monoidal product 
of objects $X$, $Y$ in $\cat{Rep}(\pi)$
is the pair
$X\tensor\!_{\pi}Y=(x\otimes y,\rho_{x\otimes y})$
where 
$\rho_{x\otimes y}:
x\otimes y
\xrightarrow{\rho_x\otimes \rho_y}
(x\otimes \pi)\otimes (y\otimes \pi)
\xrightarrow{(\otimes\pi)_{x,y}}
(x\otimes y)\otimes \pi$
and the unit object in $\cat{Rep}(\pi)$
is the pair
$\unit\!_{\pi}=(\kappa,\rho_{\kappa})$ where
$\rho_{\kappa}=\upsilon^{\otimes\pi}_{\kappa}:\kappa\xrightarrow[\cong]{\pc_{\kappa}^{-1}}\kappa\otimes \kappa\xrightarrow{I_{\kappa}\otimes u_{\pi}}\kappa\otimes \pi$.
The symmetric monoidal coherence morphisms
$a$, $\imath$, $\jmath$, $s$ 
of $\kos{R\!e\!p}(\pi)$ are given by those of $\kos{K}$.
The functor
$\kos{t}_{\pi}^*$
sends each object $z$ in $\CK$ to the pair
\begin{equation*}
  \kos{t}_{\pi}^*(z)
  =
  (z,\what{\otimes\pi}_z:
  z\to z\otimes \pi),
  \qquad
  \what{\otimes\pi}_z
  =
  \upsilon^{\otimes\pi}_z:
  \xymatrix{
    z
    \ar[r]^-{\jmath_z}_-{\cong}
    &z\otimes \kappa
    \ar[r]^-{I_z\otimes u_{\pi}}
    &z\otimes \pi
    .
  }
\end{equation*}
The symmetric monoidal coherence isomorphisms of
$\kos{t}_{\pi}^*:\kos{K}\to\kos{R\!e\!p}(\pi)$
are given by identity morphisms.

\begin{remark}
  Let $\cat{Spec}(\pi)$ be a group object in $\kos{A\!f\!f}(\kos{K})$.
  We explain why we define representations of $\cat{Spec}(\pi)$ in $\kos{K}$
  as in Definition~\ref{def RepGrpAff(K) Rep(pi)Ktensorcat}.
  For each object $x$ in $\CK$,
  we denote 
  \begin{equation*}
    \underline{\Aut}(x):\cat{Aff}(\kos{K})^{\op}\to\cat{Grp}
  \end{equation*}
  as the presheaf of groups of automorphisms of $x$.
  It sends each object $\cat{Spec}(b)$ in $\cat{Aff}(\kos{K})$
  to the group $\Aut_{\CK_b}(\kos{b}^*(x))$
  of automorphisms of $\kos{b}^*(x)$ in the Kleisli category $\CK_b$,
  and 
  each morphism $\cat{Spec}(b^{\pr})\xrightarrow{f}\cat{Spec}(b)$ in $\cat{Aff}(\kos{K})$
  to the morphism
  $\Aut_{\CK_b}(\kos{b}^*(x))
  \to
  \Aut_{\CK_{b^{\pr}}}(\kos{b}^{\pr*}(x))$
  of groups induced by the functor
  $\kos{f}^*:\kos{K}_b\to\kos{K}_{b^{\pr}}$
  introduced in (\ref{eq LSKtensoradj functorbetweenKelisliCat}).
  \begin{itemize}
    \item 
    Let $X=(x,\rho_x:x\to x\otimes \pi)$
    be a representation of $\cat{Spec}(\pi)$ in $\kos{K}$.
    Then we have a morphism of presheaves of groups
    \begin{equation}\label{eq RepGrpAff(K) mor of presheaves}
      r^x:
      \Hom_{\cat{Aff}(\kos{K})}\bigl(\slot,\cat{Spec}(\pi)\bigr)
      \Rightarrow
      \underline{\Aut}(x):
      \cat{Aff}(\kos{K})^{\op}\to\cat{Grp}
    \end{equation}
    which sends each morphism $\cat{Spec}(b)\xrightarrow{g}\cat{Spec}(\pi)$ in $\cat{Aff}(\kos{K})$
    to the automorphism $r^x_b(g):\kos{b}^*(x)\xrightarrow{\cong}\kos{b}^*(x)$
    in $\CK_b$ where
    $r^x_b(g):x\xrightarrow{\rho_x}x\otimes \pi\xrightarrow{I_x\otimes g}x\otimes b$.

    \item
    Conversely, suppose we are given a pair
    $(x,r^x)$ of an object $x$ in $\CK$
    and a morphism of presheaves $r^x$ as in (\ref{eq RepGrpAff(K) mor of presheaves}).
    Then we have a representation $X=(x,\rho_x)$ of $\cat{Spec}(\pi)$ in $\kos{K}$
    where $\rho_x=r^x_{\pi}(I_{\pi}):x\to x\otimes \pi$.
  \end{itemize}
  Using the Yoneda Lemma,
  one can show that the correspondences
  $X=(x,\rho_x)\mapsto (x,r^x)$
  and 
  $(x,r^x)\mapsto X=(x,\rho_x)$
  described above are inverse to each other.
\end{remark}

\begin{lemma} \label{lem RepGrpAff(K) Rep(pi)}
  Let $\cat{Spec}(\pi)$ be a group object in $\kos{A\!f\!f}(\kos{K})$
  and recall the strong $\kos{K}$-tensor category
  $(\kos{R\!e\!p}(\pi),\kos{t}^*_{\pi})$
  of representations of $\cat{Spec}(\pi)$ in $\kos{K}$
  defined in Definition~\ref{def RepGrpAff(K) Rep(pi)Ktensorcat}.
  \begin{enumerate}
    \item 
    The functor $\kos{t}^*_{\pi}$ of constructing trivial representations
    is fully faithful.
    
    \item 
    We have a left-strong $\kos{K}$-tensor adjunction
    \begin{equation*}
      \vcenter{\hbox{
        \xymatrix{
          (\kos{R\!e\!p}(\pi),\kos{t}_{\pi}^*)
          \ar@/_1pc/[d]_-{(\omega_{\pi}^*,\what{\omega}_{\pi}^*)}
          \\
          (\kos{K},\id_{\kos{K}})
          \ar@/_1pc/[u]_-{(\omega_{\pi*},\what{\omega}_{\pi*})}
        }
      }}
      \qquad
      (\omega_{\pi},\what{\omega}_{\pi})
      :(\kos{K},\id_{\kos{K}})\to 
      (\kos{R\!e\!p}(\pi),\kos{t}_{\pi}^*)
    \end{equation*}
    satisfying the projection formula,
    such that the right adjoint
    $(\omega_{\pi*},\what{\omega}_{\pi*})$
    is a coreflective lax $\kos{K}$-tensor functor
    and the underlying functor $\omega_{\pi*}$ is conservative.
  \end{enumerate}
\end{lemma}
\begin{proof}
  One can easily check that the functor
  $\kos{t}^*_{\pi}:\CK\to \cat{Rep}(\pi)$
  of constructing trivial representations
  is fully faithful.
  The left-strong $\kos{K}$-tensor adjunction
  $(\omega_{\pi},\what{\omega}_{\pi})$
  is defined as in (\ref{eq LaxKTcomonad E-McatLSKTadj}).
  \begin{itemize}
    \item 
    The left adjoint $\omega^*_{\pi}$ is
    the forgetful functor
    which sends each object
    $X=(x,\rho_x)$ in $\cat{Rep}(\pi)$
    to the underlying object $\omega^*_{\pi}(X)=x$ in $\CK$.
    The symmetric monoidal coherence isomorphisms of
    $\omega^*_{\pi}:\kos{R\!e\!p}(\pi)\to \kos{K}$
    are identity morphisms,
    and the comonoidal natural isomorphism
    $\what{\omega}_{\pi}^*:\id_{\kos{K}}=\omega_{\pi}^*\kos{t}_{\pi}^*$
    is the identity natural transformation.
    
    \item
    The right adjoint $\omega_{\pi*}$ sends each object
    $x$ in $\CK$
    to the object
    $\omega_{\pi*}(x)=(x\otimes \pi,\delta^{\otimes\pi}_x)$
    in $\cat{Rep}(\pi)$
    where
    $\delta^{\otimes \pi}_x:x\otimes \pi
    \xrightarrow{I_x\otimes \cp_{\pi}}
    x\otimes (\pi\otimes \pi)
    \xrightarrow[\cong]{a_{x,\pi,\pi}}
    (x\otimes \pi)\otimes \pi$.
    The lax symmetric monoidal coherence morphisms of
    $\omega_{\pi*}:\kos{K}\to \kos{R\!e\!p}(\pi)$ are
    given by those of $\otimes\pi:\kos{K}\to\kos{K}$,
    and the monoidal natural transformation
    $\what{\omega}_{\pi*}:\kos{t}_{\pi}^*\Rightarrow\omega_{\pi*}$
    is given by
    $\what{\otimes\pi}:\id_{\kos{K}}\Rightarrow\otimes\pi$.
  
    \item
    The component of the adjunction unit at each object
    $X=(x,\rho_x)$ in $\cat{Rep}(\pi)$
    is $\rho_x:X\to \omega_{\pi*}\omega_{\pi}^*(X)$.
    The component of the adjunction counit at each object
    $x$ in $\CK$ is
    $\epsilon^{\otimes\pi}_x:\omega_{\pi}^*\omega_{\pi*}(x)=x\otimes \pi\xrightarrow{I_x\otimes e_{\pi}}x\otimes \kappa \xrightarrow[\cong]{\jmath_x^{-1}}x$.
  \end{itemize}
  The left-strong $\kos{K}$-tensor adjunction
  $(\omega_{\pi},\what{\omega}_{\pi})$
  satisfies the projection formula.
  For each object
  $X=(x,\rho_x)$ in $\cat{Rep}(\pi)$
  and each object $z$ in $\CK$,
  \begin{equation*}
    \varphi_{X,z}:
    \xymatrix@C=30pt{
      x\otimes z\otimes \pi
      \ar[r]^-{\rho_x\otimes I_{z\otimes \pi}}
      &x\otimes \pi\otimes z\otimes \pi
      \ar[r]^-{(\otimes \pi)_{x,z}}
      &x\otimes z\otimes \pi
    }
  \end{equation*}
  is an isomorphism in $\CK$,
  whose inverse is
  \begin{equation*}
    \varphi_{X,z}^{-1}:
    \xymatrix@C=35pt{
      x\otimes z\otimes \pi
      \ar[r]^-{\rho_x\otimes I_{z\otimes \pi}}
      &x\otimes \pi\otimes z\otimes \pi
      \ar[r]^-{I_x\otimes \varsigma_{\pi}\otimes I_{z\otimes \pi}}
      &x\otimes \pi\otimes z\otimes \pi
      \ar[r]^-{(\otimes \pi)_{x,z}}
      &x\otimes z\otimes \pi
      .
    }
  \end{equation*}
  The lax $\kos{K}$-tensor endofunctor
  $(\otimes\pi,\what{\otimes\pi})$ on $(\kos{K},\id_{\kos{K}})$ is coreflective
  and the underlying functor
  $\otimes\pi:\CK\to \CK$ 
  is conservative as explained in Remark~\ref{rem RepGrpAff(K) tensorpi conservative}.
  Since we have
  $(\omega_{\pi}^*\omega_{\pi*},\what{\omega_{\pi}^*\omega_{\pi*}})
  =(\otimes\pi,\what{\otimes\pi})$
  as lax $\kos{K}$-tensor endofunctor on $(\kos{K},\id_{\kos{K}})$,
  we obtain that the right adjoint
  $(\omega_{\pi*},\what{\omega}_{\pi*})
  :(\kos{K},\id_{\kos{K}})\to 
  (\kos{R\!e\!p}(\pi),\kos{t}_{\pi}^*)$
  is also coreflective,
  and the underlying functor
  $\omega_{\pi*}:\CK\to \cat{Rep}(\pi)$
  is also conservative.
  This completes the proof of Lemma~\ref{lem RepGrpAff(K) Rep(pi)}.
\qed\end{proof}

We define as follows,
which is analogous to
Definition~\ref{def LSKtensoradj Hom,Isom}.

\begin{definition} \label{def RepGrpAff(K) End,Aut}
  Let $(\kos{T},\mon{t})$ be a lax $\kos{K}$-tensor category
  and let $(\phi,\what{\phi}):(\kos{T},\mon{t})\to (\kos{K},\id_{\kos{K}})$
  be a lax $\kos{K}$-tensor functor.
  \begin{enumerate}
    \item 
    We define the \emph{presheaf of monoids of monoidal $\kos{K}$-tensor natural endomorphisms
    of $(\phi,\what{\phi})$}
    as the presheaf of monoids on $\cat{Aff}(\kos{K})$
    \begin{equation*}
      \underline{\End}_{\mathbb{SMC}_{\cat{lax}}^{\kos{K}\backslash\!\!\backslash}}
      (\phi,\what{\phi})
      :
      \cat{Aff}(\kos{K})^{\op}
      \to
      \cat{Mon}
    \end{equation*}
    which sends each object $\cat{Spec}(b)$ in $\cat{Aff}(\kos{K})$
    to the set of monoidal $\kos{K}$-tensor natural endomorphisms
    $(\kos{b}^*\phi,\what{\kos{b}^*\phi})
    \Rightarrow
    (\kos{b}^*\phi,\what{\kos{b}^*\phi})
    :(\kos{T},\mon{t})\to (\kos{K}_b,\kos{b}^*)$.
    
    \item 
    We define the \emph{presheaf of groups of monoidal $\kos{K}$-tensor natural automorphisms
    of $(\phi,\what{\phi})$}
    as the presheaf of groups on $\cat{Aff}(\kos{K})$
    \begin{equation*}
      \underline{\Aut}_{\mathbb{SMC}_{\cat{lax}}^{\kos{K}\backslash\!\!\backslash}}
      (\phi,\what{\phi})
      :
      \cat{Aff}(\kos{K})^{\op}
      \to
      \cat{Grp}
    \end{equation*}
    which sends each object $\cat{Spec}(b)$ in $\cat{Aff}(\kos{K})$
    to the group of monoidal $\kos{K}$-tensor natural automorphisms
    $\xymatrix@C=12pt{(\kos{b}^*\phi,\what{\kos{b}^*\phi})
    \ar@2{->}[r]^-{\cong}
    &(\kos{b}^*\phi,\what{\kos{b}^*\phi})
    :(\kos{T},\mon{t})\to (\kos{K}_b,\kos{b}^*).}$
  \end{enumerate}
\end{definition}

The following is the main result of this subsection,
which we obtain by applying the results in \textsection~\ref{subsec LSKtensoradj}.

\begin{proposition} \label{prop RepGrpAff(K) comparison functor}
  Let $(\kos{T},\mon{t})$
  be a strong $\kos{K}$-tensor category
  and let 
  $(\omega,\what{\omega})
  :(\kos{K},\id_{\kos{K}})\to (\kos{T},\mon{t})$
  be a coreflective left-strong $\kos{K}$-tensor adjunction.
  \begin{equation*}
    \vcenter{\hbox{
      \xymatrix{
        (\kos{T},\mon{t})
        \ar@/_1pc/[d]_-{(\omega^*,\what{\omega}^*)}
        \\
        (\kos{K},\id_{\kos{K}})
        \ar@/_1pc/[u]_-{(\omega_*,\what{\omega}_*)\text{ coreflective}}
      }
    }}
  \end{equation*}
  Assume further that
  $\omega^*(\varphi_{\slot}):
  \omega^*(\slot\tensor \omega_*(\kappa))
  \Rightarrow
  \omega^*\omega_*\omega^*
  :\CT\to \CK$
  is a natural isomorphism.
  \begin{enumerate}
    \item 
    We have a group object $\cat{Spec}(\pi)$, $\pi:=\omega^*\omega_*(\kappa)$ in $\kos{A\!f\!f}(\kos{K})$
    which represents the presheaf of groups
    of monoidal $\kos{K}$-tensor natural
    automorphisms of $(\omega^*,\what{\omega}^*)$.
    \begin{equation*}
      \xymatrix@C=18pt{
        \Hom_{\cat{Aff}(\kos{K})}\bigl(\slot, \cat{Spec}(\omega^*\omega_*(\kappa))\bigr)
        \ar@2{->}[r]^-{\cong}
        &\underline{\Aut}_{\mathbb{SMC}_{\cat{lax}}^{\kos{K}\backslash\!\!\backslash}}
        (\omega^*,\what{\omega}^*)
        :
        \cat{Aff}(\kos{K})^{\op}
        \to
        \cat{Grp}
      }
    \end{equation*}
    The product, unit, antipode morphisms of $\cat{Spec}(\omega^*\omega_*(\kappa))$
    are given as follows.
    Recall the monoidal $\kos{K}$-tensor natural automorphism
    $\varsigma^{\omega,\omega}$
    of $(\omega^*\omega_*,\what{\omega^*\omega_*})$
    defined in Lemma~\ref{lem LSKtensoradj antipode}.
    \begin{equation}\label{eq RepGrpAff(K) comparison functor}
      \begin{aligned}
        \cp_{\omega^*\omega_*(\kappa)}
        &:
        \xymatrix@C=40pt{
          \omega^*\omega_*(\kappa)
          \ar[r]^-{\omega^*(\eta_{\omega_*(\kappa)})}
          &\omega^*\omega_*\omega^*\omega_*(\kappa)
          \ar[r]^-{(\tahar{\omega^*\omega_*}_{\omega^*\omega_*(\kappa)})^{-1}}_-{\cong}
          &\omega^*\omega_*(\kappa)
          \otimes
          \omega^*\omega_*(\kappa)
        }
        \\
        e_{\omega^*\omega_*(\kappa)}
        &:
        \xymatrix{
          \omega^*\omega_*(\kappa)
          \ar[r]^-{\epsilon_{\kappa}}
          &\kappa
        }
        \\
        \varsigma_{\omega^*\omega_*(\kappa)}
        &:
        \xymatrix{
          \omega^*\omega_*(\kappa)
          \ar[r]^-{\varsigma^{\omega,\omega}_{\kappa}}_-{\cong}
          &\omega^*\omega_*(\kappa)
        }
      \end{aligned}
    \end{equation}
    
    \item 
    The strong $\kos{K}$-tensor functor
    $(\omega^*,\what{\omega}^*):
    (\kos{T},\mon{t})\to (\kos{K},\id_{\kos{K}})$
    factors through as a strong $\kos{K}$-tensor functor
    $(\widebreve{\omega}\!^*,\what{\widebreve{\omega}}^*)
    :(\kos{T},\mon{t})\to (\kos{R\!e\!p}(\pi),\kos{t}^*_{\pi})$.
    \begin{equation*}
      \vcenter{\hbox{
        \xymatrix@C=40pt{
          (\kos{T},\mon{t})
          \ar[r]^-{(\widebreve{\omega}\!^*,\what{\widebreve{\omega}}^*)}
          \ar@/_1pc/[dr]_-{(\omega^*,\what{\omega}^*)}
          &(\kos{R\!e\!p}(\pi),\kos{t}^*_{\pi})
          \ar[d]^-{(\omega_{\pi}^*,\what{\omega}_{\pi}^*)}
          \\
          \text{ }
          &(\kos{K},\id_{\kos{K}})
        }
      }}
      \qquad\quad
      \pi=\omega^*\omega_*(\kappa)
    \end{equation*}
  \end{enumerate}
\end{proposition}
\begin{proof}
  Consider the case
  $(\omega^{\pr*},\what{\omega}^{\pr*})=(\omega^*,\what{\omega}^*)$
  in Proposition~\ref{prop LSKtensoradj HomRep,Hom=Isom}.
  Statements 1,3 of Proposition~\ref{prop LSKtensoradj HomRep,Hom=Isom}
  imply that we have an equality
  of presheaves on $\cat{Aff}(\kos{K})$
  \begin{equation*}
    \underline{\End}_{\mathbb{SMC}_{\cat{lax}}^{\kos{K}\backslash\!\!\backslash}}
    (\omega^*,\what{\omega}^*)
    =
    \underline{\Aut}_{\mathbb{SMC}_{\cat{lax}}^{\kos{K}\backslash\!\!\backslash}}
    (\omega^*,\what{\omega}^*)
    :\cat{Aff}(\kos{K})^{\op}\to\cat{Set}
  \end{equation*}
  which is represented by the object
  $\cat{Spec}(\pi)$,
  $\pi:=\omega^*\omega_*(\kappa)$ in $\cat{Aff}(\kos{K})$.
  The universal element is
  the monoidal $\kos{K}$-tensor natural automorphism
  \begin{equation*}
    \xi:
    \xymatrix@C=15pt{
      (\pi^*\omega^*,\what{\pi^*\omega^*})
      \ar@2{->}[r]^-{\cong}
      &(\pi^*\omega^*,\what{\pi^*\omega^*})
      :(\kos{T},\mon{t})\to (\kos{K}_{\pi},\pi^*)
    }
  \end{equation*}
  whose component
  $\xi_X:\pi^*\omega^*(X)\to \pi^*\omega^*(X)$
  at each object $X$ in $\CT$ is
  \begin{equation*}
    \xi_X:
    \xymatrix@C=40pt{
      \omega^*(X)
      \ar[r]^-{\omega^*(\eta_X)}
      &\omega^*\omega_*\omega^*(X)
      \ar[r]^-{(\tahar{\omega^*\omega_*}_{\omega^*(X)})^{-1}}_-{\cong}
      &\omega^*(X)\otimes \omega^*\omega_*(\kappa)
      .
    }
  \end{equation*}
  Thus the representing object
  $\cat{Spec}(\omega^*\omega_*(\kappa))$
  has a unique structure 
  of a group object in $\kos{A\!f\!f}(\kos{K})$
  so that the representation
  \begin{equation*}
    \xymatrix@C=15pt{
      \Hom_{\cat{Aff}(\kos{K})}\bigl(\slot,\cat{Spec}(\omega^*\omega_*(\kappa))\bigr)
      \ar@2{->}[r]^-{\cong}
      &\underline{\Aut}_{\mathbb{SMC}_{\cat{lax}}^{\kos{K}\backslash\!\!\backslash}}(\omega^*,\what{\omega}^*)
      :\cat{Aff}(\kos{K})^{\op}\to \cat{Grp}
    }
  \end{equation*}
  becomes an isomorphism between presheaves of groups
  on $\cat{Aff}(\kos{K})$.
  One can check that the identity natural transformations of 
  $(\omega^*,\what{\omega^*})$
  corresponds to
  $e_{\omega^*\omega_*(\kappa)}=\epsilon_{\kappa}:
  \cat{Spec}(\kappa)\to \cat{Spec}(\omega^*\omega_*(\kappa))$
  as described in (\ref{eq RepGrpAff(K) comparison functor}).
  From the case
  $(\omega^{\pr},\what{\omega}^{\pr})=(\omega,\what{\omega})$
  in Lemma~\ref{lem LSKtensoradj antipode}
  and the case
  $(\omega^{\pr},\what{\omega}^{\pr})=(\omega,\what{\omega})$,
  $(\omega^{\ppr*},\what{\omega}^{\ppr*})=(\omega^*,\what{\omega}^*)$
  in Lemma~\ref{lem LSKtensoradj composition},
  we see that the product, antipode morphisms
  $\cp_{\omega^*\omega_*(\kappa)}$, $\varsigma_{\omega^*\omega_*(\kappa)}$
  of
  $\cat{Spec}(\omega^*\omega_*(\kappa))$
  are also given as described in (\ref{eq RepGrpAff(K) comparison functor}).
  
  The given left-strong $\kos{K}$-tensor adjunction
  $(\omega,\what{\omega}):(\kos{K},\id_{\kos{K}})\to(\kos{T},\mon{t})$
  induces a lax $\kos{K}$-tensor comonad
  $\langle\omega^*\omega_*,\what{\omega^*\omega_*}\rangle
  =(\omega^*\omega_*,\what{\omega^*\omega_*},\delta,\epsilon)$
  on $(\kos{K},\id_{\kos{K}})$
  where
  \begin{equation*}
    \delta=\omega^*\eta\omega_*:
    (\omega^*\omega_*,\what{\omega^*\omega_*})
    \Rightarrow
    (\omega^*\omega_*\omega^*\omega_*,\what{\omega^*\omega_*\omega^*\omega_*}).
  \end{equation*}
  Recall the adjoint equivalence of categories
  $\iota\dashv\CR$
  in Corollary~\ref{cor RepGrpAff(K) Mon(Aff(K))=crflLKTcomonad}.
  We have an equality
  $\CR\langle\omega^*\omega_*,\what{\omega^*\omega_*}\rangle
  =\cat{Spec}(\omega^*\omega_*(\kappa))$
  of monoid objects in $\kos{A\!f\!f}(\kos{K})$
  as the product, unit morphisms of
  $\cat{Spec}(\omega^*\omega_*(\kappa))$
  in (\ref{eq RepGrpAff(K) comparison functor})
  is the same as those described in (\ref{eq RepGrpAff(K) Mon(Aff(K))=crflLKTcomonad}).
  The component of the adjunction counit
  of $\iota\dashv\CR$ in Corollary~\ref{cor RepGrpAff(K) Mon(Aff(K))=crflLKTcomonad}
  at $\langle\omega^*\omega_*,\what{\omega^*\omega_*}\rangle$
  is the following isomorphism of lax $\kos{K}$-tensor comonads on $(\kos{K},\id_{\kos{K}})$.
  \begin{equation}\label{eq2 RepGrpAff(K) comparison functor}
    \tahar{\omega^*\omega_*}:
    \xymatrix{
      \langle\otimes \omega^*\omega_*(\kappa),\what{\otimes \omega^*\omega_*(\kappa)}\rangle
      \ar@2{->}[r]^-{\cong}
      &\langle\omega^*\omega_*,\what{\omega^*\omega_*}\rangle
    }
  \end{equation}
  Using the above isomorphism of lax $\kos{K}$-tensor comonads,
  one can show that the strong $\kos{K}$-tensor functor
  $(\omega^*,\what{\omega}^*)$
  factors through as a strong $\kos{K}$-tensor functor
  $(\widebreve{\omega}\!^*,\what{\widebreve{\omega}}^*)
  :(\kos{T},\mon{t})\to (\kos{R\!e\!p}(\pi),\kos{t}^*_{\pi})$.
  The underlying functor $\widebreve{\omega}\!^*$ sends each object $X$ in $\CT$
  to the object
  $\widebreve{\omega}\!^*(X)=
  (\omega^*(X),\rho_{\omega^*(X)})$
  in $\cat{Rep}(\pi)$
  where
  \begin{equation*}
    \rho_{\omega^*(X)}=\xi_X:
    \!\!
    \xymatrix@C=40pt{
      \omega^*(X)
      \ar[r]^-{\omega^*(\eta_X)}
      &\omega^*\omega_*\omega^*(X)
      \ar[r]^-{(\tahar{\omega^*\omega_*}_{\omega^*(X)})^{-1}}_-{\cong}
      &\omega^*(X)\otimes \omega^*\omega_*(\kappa)
      .
    }
  \end{equation*}
  Let $Y$ be another object in $\CT$.
  The symmetric monoidal coherernce isomorphisms
  $\omega^*_{X,Y}:
  \omega^*(X)\otimes \omega^*(Y)
  \xrightarrow{\cong}
  \omega^*(X\tensor Y)$,
  $\omega^*_{\unit}:
  \kappa
  \xrightarrow{\cong}
  \omega^*(\unit)$
  of $\omega^*:\kos{T}\to\kos{K}$
  become isomorphisms
  $\widebreve{\omega}\!^*_{X,Y}:
  \widebreve{\omega}\!^*(X)\tensor\!_{\pi} \widebreve{\omega}^*(Y)
  \xrightarrow{\cong}
  \widebreve{\omega}\!^*(X\tensor Y)$,
  $\widebreve{\omega}\!^*_{\unit}:
  \unit\!_{\pi}
  \xrightarrow{\cong}
  \widebreve{\omega}\!^*(\unit)$
  in $\cat{Rep}(\pi)$.
  Moreover, the monoidal natural isomorphism
  $\what{\omega}^*:\id_{\kos{K}}\cong \omega^*\text{$\mon{t}$}:\kos{K}\to\kos{K}$
  becomes a monoidal natural isomorphism
  $\what{\widebreve{\omega}}^*:
  \kos{t}^*_{\pi}\cong \widebreve{\omega}\!^*\text{$\mon{t}$}:\kos{K}\to\kos{R\!e\!p}(\pi)$.
  This completes the proof of Proposition~\ref{prop RepGrpAff(K) comparison functor}.
\qed\end{proof}

\newpage
\section{Prekosmic Grothendieck categories}
\label{sec PreGroCat}

Recall the definition of
left-strong symmetric monoidal (LSSM) adjunctions
in Definition~\ref{def LSKtensoradj LSSMadj}
as well as the $2$-category
$\mathbb{ADJ}_{\cat{left}}(\mathbb{SMC}_{\cat{lax}})$
of 
symmetric monoidal categories,
LSSM adjunctions
and morphisms of LSSM adjunctions
introduced in (\ref{eq LSKtensoradj LSSM 2-cat}).

\subsection{Grothendieck prekosmoi}
\label{subsec Groprekosmoi}

\begin{definition} \label{def Groprekosmoi GROpre}
  A \emph{Grothendieck prekosmos}
  is a symmetric monoidal category
  \begin{equation*}
    \kos{K}=(\CK,\otimes,\kappa)
  \end{equation*}
  whose underlying category $\CK$
  has coreflexive equalizers.
  Given another Grothendieck prekosmos
  $\kos{T}=(\CT,\tensor,\unit)$,
  a \emph{Grothendieck morphism}
  from $\kos{T}$ to $\kos{K}$
  is a left-strong symmetric monoidal (LSSM) adjunction
  $\kos{f}:\kos{T}\to \kos{K}$.
  \begin{equation*}
    \vcenter{\hbox{
      \xymatrix{
        \kos{K}
        \ar@/_1pc/[d]_-{\kos{f}^*}
        \\
        \kos{T}
        \ar@/_1pc/[u]_-{\kos{f}_*}
      }
    }}
  \end{equation*}
  Let
  $\kos{g}:\kos{T}\to \kos{K}$
  be another Grothendieck morphism.
  A \emph{Grothendieck transformation}
  $\vartheta:\kos{f}\Rightarrow\kos{g}:\kos{T}\to\kos{K}$
  is a monoidal natural transformation
  $\vartheta:\kos{f}^*\Rightarrow\kos{g}^*:\kos{K}\to\kos{T}$
  between left adjoints.
  We define the $2$-category
  $\mathbb{GRO}^{\cat{pre}}$
  of Grothendieck prekosmoi,
  Grothendieck morphisms and Grothendieck transformations
  as the full sub-$2$-category of
  $\mathbb{ADJ}_{\cat{left}}(\mathbb{SMC}_{\cat{lax}})$.
  \begin{equation*}
    \mathbb{GRO}^{\cat{pre}}
    \hookrightarrow
    \mathbb{ADJ}_{\cat{left}}(\mathbb{SMC}_{\cat{lax}})
  \end{equation*}
\end{definition}

We introduced the $2$-category structure of
$\mathbb{GRO}^{\cat{pre}}$
in Definition~\ref{def Groprekosmoi GROpre}.
Let $\kos{T}$, $\kos{K}$ be Grothendieck prekosmoi
and let
$\kos{f}$, $\kos{g}$, $\kos{h}:\kos{T}\to\kos{K}$
be Grothendieck morphisms.
The vertical composition of Grothendieck transformations
$\kos{f}\Rightarrow\kos{g}\Rightarrow\kos{h}$
is defined as the vertical composition
$\kos{f}^*\Rightarrow\kos{g}^*\Rightarrow\kos{h}^*$
between left adjoints.

We define as follows,
which is dual to the terminologies in Definition~\ref{def Galprekosmoi Galmor property}.

\begin{definition}\label{def Groprekosmoi Gromor property}
  We say a Grothendieck morphism
  $\kos{f}$
  between Grothendieck prekosmoi is
  \begin{itemize}
    \item 
    \emph{connected}
    if the left adjoint $\kos{f}^*$
    is fully faithful.
  
    \item 
    \emph{locally connected}
    if the LSSM adjunction $\kos{f}$
    satisfies the projection formula:
    see Definition~\ref{def LSKtensoradj projformula}.
  
    \item
    \emph{\'{e}tale}
    if the LSSM adjunction $\kos{f}$ satisfies the projection formula
    and the right adjoint $\kos{f}_*$ is conservative;
  
    \item 
    \emph{surjective}
    if the left adjoint
    $\kos{f}^*$
    is conservative and preserves coreflexive equalizers.
  \end{itemize}  
\end{definition}

\begin{lemma}
  Let $\kos{K}=(\CK,\otimes,\kappa)$ be a Grothendieck prekosmos.
  Then the cocartesian monoidal category
  $\kos{C\!o\!m\!m}(\kos{K})=(\cat{Comm}(\kos{K}),\otimes,\kappa)$
  of commutative monoids in $\kos{K}$
  is also a Grothendieck prekosmos.
\end{lemma}


For the rest of this section
\textsection~\ref{sec PreGroCat},
we fix a Grothendieck prekosmos
$\kos{K}=(\CK,\otimes,\kappa)$.

\begin{definition} \label{def Groprekosmoi GROpreK}
  We define the $2$-category
  \begin{equation*}
    \mathbb{GRO}^{\cat{pre}}_{\kos{K}}
  \end{equation*}
  of Grothendieck $\kos{K}$-prekosmoi,
  Grothendieck $\kos{K}$-morphisms
  and Grothendieck $\kos{K}$-transformations
  as the slice $2$-category of $\mathbb{GRO}^{\cat{pre}}$ over $\kos{K}$.
\end{definition}

Let us explain Definition~\ref{def Groprekosmoi GROpreK}
in detail.
A \emph{Grothendieck $\kos{K}$-prekosmos}
$\mathfrak{T}=(\kos{T},\kos{t})$ is a pair
of a Grothendieck prekosmos
$\kos{T}=(\CT,\tensor,\unit)$
and a Grothendieck morphism 
\begin{equation*}
  \vcenter{\hbox{
    \xymatrix{
      \kos{K}
      \ar@/_1pc/[d]_-{\kos{t}^*}
      \\
      \kos{T}
      \ar@/_1pc/[u]_-{\kos{t}_*}
    }
  }}
  \qquad\quad
  \kos{t}:\kos{T}\to\kos{K}.
\end{equation*}
We can see $\kos{K}$ itself as a Grothendieck $\kos{K}$-prekosmos
which we denote as
$\text{$\mathfrak{K}$}=(\kos{K},\id_{\kos{K}})$.
Let $\mathfrak{S}=(\kos{S},\kos{s})$
be another Grothendieck $\kos{K}$-prekosmos
and denote
$\text{$\kos{S}$}=(\CS,\ctimes,\pzc{1})$
as the underlying Grothendieck prekosmos.
A \emph{Grothendieck $\kos{K}$-morphism}
$\mathfrak{f}=(\kos{f},\what{\kos{f}}^*):
\mathfrak{S}\to \mathfrak{T}$
is a pair of a Grothendieck morphism
$\kos{f}:\kos{S}\to\kos{T}$
and an invertible Grothendieck transformation
\begin{equation*}
  \xymatrix@C=15pt{
    \what{\kos{f}}^*:
    \kos{s}
    \ar@2{->}[r]^-{\cong}
    &\kos{t}\kos{f}
    :\kos{S}\to \kos{K}
  }  
\end{equation*}
which is a monoidal natural isomorphism between left adjoints
\begin{equation*}
  \vcenter{\hbox{
    \xymatrix@R=30pt@C=40pt{
      \text{ }
      &\kos{T}
      \ar[d]^-{\kos{f}^*}
      \\
      \kos{K}
      \ar@/^0.7pc/[ur]^-{\kos{t}^*}
      \ar[r]_-{\kos{s}^*}
      &\kos{S}
      \xtwocell[l]{}<>{<3>{\what{\kos{f}}^*\text{ }}}
    }
  }}
  \qquad\quad
  \xymatrix@C=15pt{
    \what{\kos{f}}^*:
    \kos{s}^*
    \ar@2{->}[r]^-{\cong}
    &\kos{f}^*\kos{t}^*
    :\kos{K}\to \kos{S}
    .
  }
\end{equation*}
Let $\mathfrak{g}=(\kos{g},\what{\kos{g}}^*):
\mathfrak{S}\to\mathfrak{T}$
be another Grothendieck $\kos{K}$-morphism.
A \emph{Grothendieck $\kos{K}$-transformation}
$\vartheta:\mathfrak{f}\Rightarrow\mathfrak{g}:\mathfrak{S}\to\mathfrak{T}$
is a Grothendieck transformation
$\vartheta:\kos{f}\Rightarrow\kos{g}:\kos{S}\to\kos{T}$
which satisfies the relation
\begin{equation} \label{eq Groprekosmoi GroKtransformation}
  \vcenter{\hbox{
    \xymatrix@C=15pt{
      \kos{f}^*\kos{t}^*
      \ar@2{->}[rr]^-{\vartheta\kos{t}^*}
      &\text{ }
      &\kos{g}^*\kos{t}^*
      \\
      \text{ }
      &\kos{s}^*
      \ar@2{->}[ul]^-{\what{\kos{f}}^*}_-{\cong}
      \ar@2{->}[ur]_-{\what{\kos{g}}^*}^-{\cong}
    }
  }}
  \qquad\quad
  \what{\kos{g}}^*
  =
  (\vartheta^*\kos{t}^*)
  \circ \what{\kos{f}}^*
  :\kos{K}\to \kos{S}
  .
\end{equation}
We can also describe relation (\ref{eq Groprekosmoi GroKtransformation})
as follows.
\begin{equation*}
  \vcenter{\hbox{
    \xymatrix@R=30pt@C=40pt{
      \text{ }
      &\kos{T}
      \ar[d]^-{\kos{g}^*}
      \\
      \kos{K}
      \ar@/^0.7pc/[ur]^-{\kos{t}^*}
      \ar[r]_-{\kos{s}^*}
      &\kos{S}
      \xtwocell[l]{}<>{<3>{\what{\kos{g}}^*\text{ }}}
    }
  }}
  \quad=
  \vcenter{\hbox{
    \xymatrix@R=30pt@C=25pt{
      \text{ }
      &\text{ }
      &\kos{T}
      \ar@/^1pc/[d]^-{\kos{g}^*}
      \ar@/_1pc/[d]_-{\kos{f}^*}
      \\
      \kos{K}
      \ar@/^1pc/[urr]^-{\kos{t}^*}
      \ar[rr]_-{\kos{s}^*}
      &\text{ }
      &\kos{S}
      \xtwocell[u]{}<>{<0>{\vartheta}}
      \xtwocell[ll]{}<>{<2>{\what{\kos{f}}^*\text{ }}}
    }
  }}
\end{equation*}

Recall the $2$-category
$\mathbb{ADJ}_{\cat{left}}\bigl(\mathbb{SMC}_{\cat{lax}}^{\kos{K}\backslash\!\!\backslash}\bigr)$
of lax $\kos{K}$-tensor categories,
left-strong $\kos{K}$-tensor adjunctions,
and morphisms of left-strong $\kos{K}$-tensor adjunctions
introduced in (\ref{eq LSKtensoradj LSKTadj 2-cat}).
The $2$-category
$\mathbb{GRO}^{\cat{pre}}_{\kos{K}}$
of Grothendieck $\kos{K}$-prekosmoi
is a full sub-$2$-category of
$\mathbb{ADJ}_{\cat{left}}\bigl(\mathbb{SMC}_{\cat{lax}}^{\kos{K}\backslash\!\!\backslash}\bigr)$.
\begin{equation} \label{eq Groprekosmoi GropreK LSKTadj}
  \mathbb{GRO}^{\cat{pre}}_{\kos{K}}
  \hookrightarrow
  \mathbb{ADJ}_{\cat{left}}\bigl(\mathbb{SMC}_{\cat{lax}}^{\kos{K}\backslash\!\!\backslash}\bigr)
\end{equation}
\begin{itemize}
  \item 
  Let $\mathfrak{T}=(\kos{T},\kos{t})$
  be a Grothendieck $\kos{K}$-prekosmos.
  Then the pair $(\kos{T},\kos{t}^*)$ is a
  strong $\kos{K}$-tensor category.

  \item 
  Let 
  $\mathfrak{S}=(\kos{S},\kos{s})$
  be another Grothendieck $\kos{K}$-prekosmos
  and let
  $\mathfrak{f}=(\kos{f},\what{\kos{f}}^*):
  \mathfrak{S}\to \mathfrak{T}$
  be a Grothendieck $\kos{K}$-morphism.
  Then the pair
  $(\kos{f}^*,\what{\kos{f}}^*):
  (\kos{T},\kos{t}^*)\to (\kos{S},\kos{s}^*)$
  is a strong $\kos{K}$-tensor functor.
  Thus the right adjoint has a unique lax $\kos{K}$-tensor functor structure
  $(\kos{f}_*,\what{\kos{f}}_*):
  (\kos{S},\kos{s}^*)\to (\kos{T},\kos{t}^*)$
  such that the given Grothendieck morphism
  $\kos{f}:\kos{S}\to \kos{T}$
  becomes a left-strong $\kos{K}$-tensor adjunction
  \begin{equation} \label{eq Groprekosmoi GroKmor associated LSKTadj}
    \vcenter{\hbox{
      \xymatrix{
        (\kos{T},\kos{t}^*)
        \ar@/_1pc/[d]_-{(\kos{f}^*,\what{\kos{f}}^*)}
        \\
        (\kos{S},\kos{s}^*)
        \ar@/_1pc/[u]_-{(\kos{f}_*,\what{\kos{f}}_*)}
      }
    }}
    \qquad\qquad
    (\kos{f},\what{\kos{f}}):
    (\kos{S},\kos{s}^*)\to (\kos{T},\kos{t}^*)
    .
  \end{equation}
  The monoidal natural transformation
  $\what{\kos{f}}_*$ is 
  given as in (\ref{eq LSKtensoradj rightadj uniqueKTstr}).
  \begin{equation*}
    \vcenter{\hbox{
      \xymatrix@R=30pt@C=40pt{
        \text{ }
        &\kos{S}
        \ar[d]^-{\kos{f}_*}
        \\
        \kos{K}
        \ar@/^0.7pc/[ur]^-{\kos{s}^*}
        \ar[r]_-{\kos{t}^*}
        &\kos{T}
        \xtwocell[l]{}<>{<3>{\what{\kos{f}}_*\text{ }}}
      }
    }}
    \qquad\quad
    \xymatrix@C=30pt{
      \what{\kos{f}}_*:
      \kos{t}^*
      \ar@2{->}[r]^-{\eta\kos{t}^*}
      &\kos{f}_*\kos{f}^*\kos{t}^*
      \ar@2{->}[r]^-{\kos{f}_*(\what{\kos{f}}^*)^{-1}}_-{\cong}
      &\kos{f}_*\kos{s}^*
    }
  \end{equation*}

  \item 
  Let $\mathfrak{g}=(\kos{g},\what{\kos{g}}^*):\mathfrak{S}\to\mathfrak{T}$
  be another Grothendieck $\kos{K}$-morphism.
  Then a Grothendieck $\kos{K}$-transformation
  $\vartheta:\mathfrak{f}\Rightarrow\mathfrak{g}:\mathfrak{S}\to\mathfrak{T}$
  is precisely 
  a monoidal $\kos{K}$-tensor natural transformation
  $\vartheta:(\kos{f}^*,\what{\kos{f}}^*)\Rightarrow (\kos{g}^*,\what{\kos{g}}^*)
  :(\kos{T},\kos{t}^*)\to (\kos{S},\kos{s}^*)$.
\end{itemize}
We denote a Grothendieck $\kos{K}$-morphism
$\mathfrak{f}:\mathfrak{S}\to \mathfrak{T}$
between Grothendieck $\kos{K}$-prekosmoi
$\mathfrak{T}=(\kos{T},\kos{t})$,
$\mathfrak{S}=(\kos{S},\kos{s})$
as
\begin{equation*}
  \mathfrak{f}=(\kos{f},\what{\kos{f}}):\mathfrak{S}\to\mathfrak{T}
\end{equation*}
where
$(\kos{f},\what{\kos{f}}):
(\kos{S},\kos{s}^*)\to (\kos{T},\kos{t}^*)$
is the associated left-strong $\kos{K}$-tensor adjunction
described in (\ref{eq Groprekosmoi GroKmor associated LSKTadj}).

\begin{definition}
  A Grothendieck $\kos{K}$-prekosmoi 
  $\mathfrak{T}=(\kos{T},\kos{t})$
  is called 
  \begin{equation*}
    \text{\emph{connected|locally connected|\'{e}tale|surjective}}
  \end{equation*}
  if the Grothendieck morphism
  $\kos{t}:\kos{T}\to \kos{K}$ has such property.
\end{definition}

\begin{definition}
  Let
  $\mathfrak{T}=(\kos{T},\kos{t})$,
  $\mathfrak{S}=(\kos{S},\kos{s})$
  be Grothendieck $\kos{K}$-prekosmoi.
  A Grothendieck $\kos{K}$-morphism
  $\mathfrak{f}=(\kos{f},\what{\kos{f}}):\mathfrak{S}\to\mathfrak{T}$
  is called \emph{connected|locally connected|\'{e}tale|surjective}
  if the underlying Grothendieck morphism $\kos{f}:\kos{S}\to\kos{T}$ has such property.
\end{definition}

\begin{definition}
  Let $\mathfrak{T}=(\kos{T},\kos{t})$
  be a Grothendieck $\kos{K}$-prekosmos
  and recall that we denote
  $\text{$\mathfrak{K}$}=(\kos{K},\id_{\kos{K}})$.
  We say a Grothendieck $\kos{K}$-morphism
  \begin{equation*}
    \varpi=(\omega,\what{\omega})
    :\mathfrak{K}\to\mathfrak{T}
  \end{equation*}
  is \emph{coreflective}
  if the associated left-strong $\kos{K}$-tensor adjunction
  $(\omega,\what{\omega}):(\kos{K},\id_{\kos{K}})\to (\kos{T},\kos{t}^*)$
  is coreflective:
  see Definition~\ref{def LSKtensoradj coreflective}.
  This amounts to saying that
  the right adjoint
  $(\omega_*,\what{\omega}_*):
  (\kos{K},\id_{\kos{K}})\to (\kos{T},\kos{t}^*)$
  is a coreflective lax $\kos{K}$-tensor functor,
  or equivalently,
  the associated $\kos{K}$-equivariance
  $\cevar{\omega}_*$
  is a natural isomorphism:
  see Lemma~\ref{lem Comm(T) equivariance and coreflection}.
\end{definition}

\begin{definition}
  Let $\mathfrak{T}=(\kos{T},\kos{t})$
  be a Grothendieck $\kos{K}$-prekosmos.
  \begin{enumerate}
    \item 
    A \emph{pre-fiber functor} for $\mathfrak{T}$
    is a coreflective \'{e}tale Grothendieck $\kos{K}$-morphism
    \begin{equation*}
      \varpi=(\omega,\what{\omega}):\mathfrak{K}\to \mathfrak{T}.
    \end{equation*}
    Equivalently, it is a Grothendieck $\kos{K}$-morphism
    $\varpi=(\omega,\what{\omega}):\mathfrak{K}\to\mathfrak{T}$
    which has the following properties:
    \begin{itemize}
      \item 
      the Grothendieck morphism $\omega:\kos{K}\to \kos{T}$
      satisfies the projection formula;
      
      \item 
      the right adjoint $(\omega_*,\what{\omega}_*):(\kos{K},\id_{\kos{K}})\to (\kos{T},\kos{t}^*)$
      is a coreflective lax $\kos{K}$-tensor functor,
      and the underlying functor
      $\omega_*$ is conservative.
    \end{itemize}
    
    \item
    A pre-fiber functor $\varpi=(\omega,\what{\omega}):\mathfrak{K}\to\mathfrak{T}$
    is called \emph{surjective}
    if it is surjective as a Grothendieck $\kos{K}$-morphism.
    This amounts to saying that
    the left adjoint $\omega^*$ is conservative and preserves coreflexive equalizers.
  \end{enumerate}
  A \emph{morphism} of pre-fiber functors $\varpi$, $\varpi^{\pr}$ for $\mathfrak{T}$
  is a Grothendieck $\kos{K}$-transformation 
  $\vartheta:\varpi\Rightarrow\varpi^{\pr}:\mathfrak{K}\to\mathfrak{T}$.
  We denote $\cat{Fib}(\mathfrak{T})^{\cat{pre}}$
  as the category of pre-fiber functors for $\mathfrak{T}$.
\end{definition}

\begin{definition}
  Let $\mathfrak{T}=(\kos{T},\kos{t})$ be a Grothendieck $\kos{K}$-prekosmos
  and let $\varpi=(\omega,\what{\omega})$,
  $\varpi^{\pr}=(\omega^{\pr},\what{\omega}^{\pr}):\mathfrak{K}\to\mathfrak{T}$
  be pre-fiber functors.
  Consider the strong $\kos{K}$-tensor functors
  $(\omega^*,\what{\omega}^*)$,
  $(\omega^{\pr*},\what{\omega}^{\pr*})
  :(\kos{T},\kos{t}^*)\to (\kos{K},\id_{\kos{K}})$
  and recall Definition~\ref{def LSKtensoradj Hom,Isom}
  as well as Definition~\ref{def RepGrpAff(K) End,Aut}.
  \begin{enumerate}
    \item 
    We define the \emph{presheaf of Grothendieck $\kos{K}$-transformations
    from $\varpi^{\pr}$ to $\varpi$}
    as
    \begin{equation*}
      \underline{\Hom}_{\mathbb{GRO}^{\cat{pre}}_{\kos{K}}}(\varpi^{\pr},\varpi)
      :=
      \underline{\Hom}_{\mathbb{SMC}_{\cat{lax}}^{\kos{K}\backslash\!\!\backslash}}
      \bigl((\omega^{\pr*},\what{\omega}^{\pr*}),(\omega^*,\what{\omega}^*)\bigr)
      :\cat{Aff}(\kos{K})^{\op}\to\cat{Set}
    \end{equation*}
    and the \emph{presheaf of monoids of Grothendieck $\kos{K}$-transformations
    from $\varpi$ to $\varpi$}
    as
    \begin{equation*}
      \underline{\End}_{\mathbb{GAL}^{\cat{pre}}_{\kos{K}}}(\varpi)
      :=
      \underline{\End}_{\mathbb{SMC}_{\cat{lax}}^{\kos{K}\backslash\!\!\backslash}}
      (\omega^*,\what{\omega}^*)
      :\cat{Aff}(\kos{K})^{\op}\to\cat{Mon}
      . 
    \end{equation*}

    \item 
    We define the \emph{presheaf of invertible Grothendieck $\kos{K}$-transformations
    from $\varpi^{\pr}$ to $\varpi$}
    as
    \begin{equation*}
      \underline{\Isom}_{\mathbb{GRO}^{\cat{pre}}_{\kos{K}}}(\varpi^{\pr},\varpi)
      :=
      \underline{\Isom}_{\mathbb{SMC}_{\cat{lax}}^{\kos{K}\backslash\!\!\backslash}}
      \bigl((\omega^{\pr*},\what{\omega}^{\pr*}),(\omega^*,\what{\omega}^*)\bigr)
      :\cat{Aff}(\kos{K})^{\op}\to\cat{Set}
    \end{equation*}
    and the \emph{presheaf of groups of invertible Grothendieck $\kos{K}$-transformations
    from $\varpi$ to $\varpi$} as
    \begin{equation*}
      \underline{\Aut}_{\mathbb{GRO}^{\cat{pre}}_{\kos{K}}}(\varpi)
      :=
      \underline{\Aut}_{\mathbb{SMC}_{\cat{lax}}^{\kos{K}\backslash\!\!\backslash}}
      (\omega^*,\what{\omega}^*)
      :\cat{Aff}(\kos{K})^{\op}\to\cat{Grp}
      . 
    \end{equation*}
  \end{enumerate}
\end{definition}

\begin{remark}
  Let $\mathfrak{T}=(\kos{T},\kos{t})$ be a Grothendieck $\kos{K}$-prekosmos
  and let $\varpi=(\omega,\what{\omega})$,
  $\varpi^{\pr}=(\omega^{\pr},\what{\omega}^{\pr}):\mathfrak{K}\to\mathfrak{T}$
  be pre-fiber functors.
  For each object $\cat{Spec}(b)$ in $\cat{Aff}(\kos{K})$,
  an element $\vartheta$ in 
  $\underline{\Hom}_{\mathbb{GRO}^{\cat{pre}}_{\kos{K}}}(\varpi^{\pr},\varpi)(\cat{Spec}(b))$
  is a monoidal $\kos{K}$-tensor natural transformation
  \begin{equation*}
    \vartheta:
    (\kos{b}^*\omega^{\pr*},\what{\kos{b}^*\omega^{\pr*}})
    \Rightarrow
    (\kos{b}^*\omega^*,\what{\kos{b}^*\omega^*})
    :(\kos{T},\kos{t}^*)\to (\kos{K}_b,\kos{b}^*).
  \end{equation*}
  Note that $\vartheta$ is not a Grothendieck $\kos{K}$-transformation,
  since the symmetric monoidal category $\kos{K}_b$ is in general not 
  a Grothendieck prekosmos.
  Nontheless, 
  we are going to call 
  $\underline{\Hom}_{\mathbb{GRO}^{\cat{pre}}_{\kos{K}}}(\varpi^{\pr},\varpi)$
  as the presheaf of Grothendieck $\kos{K}$-transformations
  from $\varpi^{\pr}$ to $\varpi$.
\end{remark}

The following is a consequence of results in \textsection~\ref{subsec LSKtensoradj}.

\begin{corollary} \label{cor Groprekosmoi summary}
  Let $\mathfrak{T}=(\kos{T},\kos{t})$ be a Grothendieck $\kos{K}$-prekosmos
  and let $\varpi=(\omega,\what{\omega})$,
  $\varpi^{\pr}=(\omega^{\pr},\what{\omega}^{\pr}):\mathfrak{K}\to\mathfrak{T}$
  be pre-fiber functors.
  \begin{enumerate}
    \item 
    The presheaf of Grothendieck $\kos{K}$-transformations from $\varpi^{\pr}$ to $\varpi$
    is represented by the object
    $\cat{Spec}(\omega^{\pr*}\omega_*(\kappa))$ in $\cat{Aff}(\kos{K})$.
    \begin{equation*}
      \vcenter{\hbox{
        \xymatrix{
          \Hom_{\cat{Aff}(\kos{K})}\bigl(\slot,\cat{Spec}(\omega^{\pr*}\omega_*(\kappa))\bigr)
          \ar@2{->}[r]^-{\cong}
          &\underline{\Hom}_{\mathbb{GRO}^{\cat{pre}}_{\kos{K}}}(\varpi^{\pr},\varpi)
          :\cat{Aff}(\kos{K})^{\op}\to\cat{Set}
        }
      }}
    \end{equation*}

    \item
    We have
    $\underline{\Hom}_{\mathbb{GRO}^{\cat{pre}}_{\kos{K}}}(\varpi^{\pr},\varpi)
    =
    \underline{\Isom}_{\mathbb{GRO}^{\cat{pre}}_{\kos{K}}}(\varpi^{\pr},\varpi)
    :\cat{Aff}(\kos{K})^{\op}\to\cat{Set}$.

    \item
    We have a monoidal $\kos{K}$-tensor natural isomorphism
    \begin{equation*}
      \varsigma^{\omega^{\pr},\omega}:
      \xymatrix@C=18pt{
        (\omega^*\omega^{\pr}_*,\what{\omega^*\omega^{\pr}_*})
        \ar@2{->}[r]^-{\cong}
        &(\omega^{\pr*}\omega_*,\what{\omega^{\pr*}\omega_*})
        :(\kos{K},\id_{\kos{K}})\to (\kos{K},\id_{\kos{K}})
      }
    \end{equation*}
    whose component at $\kappa$ is an isomorphism
    $\varsigma^{\omega^{\pr},\omega}_{\kappa}:
    \omega^*\omega^{\pr}_*(\kappa)
    \xrightarrow[]{\cong}
    \omega^{\pr*}\omega_*(\kappa)$
    in $\cat{Comm}(\kos{K})$,
    and the following diagram of presheaves strictly commutes.
    \begin{equation*}
      \vcenter{\hbox{
        \xymatrix{
          \Hom_{\cat{Aff}(\kos{K})}\bigl(\slot,\cat{Spec}(\omega^{\pr*}\omega_*(\kappa))\bigr)
          \ar@2{->}[d]_-{\varsigma^{\omega^{\pr},\omega}_{\kappa}\circ(\slot)}^-{\cong}
          \ar@2{->}[r]^-{\cong}
          &\underline{\Isom}_{\mathbb{GRO}^{\cat{pre}}_{\kos{K}}}(\varpi^{\pr},\varpi)
          \ar@2{->}[d]^-{\cat{inverse}}_-{\cong}
          \\
          \Hom_{\cat{Aff}(\kos{K})}\bigl(\slot,\cat{Spec}(\omega^*\omega^{\pr}_*(\kappa))\bigr)
          \ar@2{->}[r]^-{\cong}
          &\underline{\Isom}_{\mathbb{GRO}^{\cat{pre}}_{\kos{K}}}(\varpi,\varpi^{\pr})
        }
      }}
    \end{equation*}
  \end{enumerate}
  In particular, the category
  $\cat{Fib}(\mathfrak{T})^{\cat{pre}}$
  of pre-fiber functors for $\mathfrak{T}$ is a groupoid.
\end{corollary}
\begin{proof}
  The left-strong $\kos{K}$-tensor adjunctions
  $(\omega,\what{\omega})$,
  $(\omega^{\pr},\what{\omega}^{\pr})$
  associated to $\varpi$, $\varpi^{\pr}$
  are coreflective and satisfy the projection formula.
  In particular, the associated natural transformations
  $\varphi_{\slot}$,
  $\varphi^{\pr}_{\slot}$
  defined as in Definition~\ref{def LSKtensoradj coreflective}
  are natural isomorphisms.
  By 
  Proposition~\ref{prop LSKtensoradj HomRep,Hom=Isom}
  and 
  Lemma~\ref{lem LSKtensoradj antipode},
  we immediately obtain statements 1-3.
  Moreover, statement 2 implies that
  $\cat{Fib}(\mathfrak{T})_{\kappa}$
  is a groupoid.
\qed\end{proof}

\subsection{Pre-Grothendieck objects in a Grothendieck prekosmos $\kos{K}$}
\label{subsec preGroobj}

Let $\kos{K}=(\CK,\otimes,\kappa)$ be a Grothendieck prekosmos.
Recall the definition of a group object $\cat{Spec}(\pi)$ in $\kos{A\!f\!f}(\kos{K})$
introduced in \textsection\!~\ref{subsec RepGrpAff(K)}

\begin{definition}
  A \emph{pre-Grothendieck object} in $\kos{K}$
  is a group object $\cat{Spec}(\pi)$ in $\kos{A\!f\!f}(\kos{K})$
  such that the functor $\slot\otimes\pi:\CK\to\CK$
  preserves coreflexive equalizers.
\end{definition}

Let $\cat{Spec}(\pi)$ be a pre-Grothendieck object in $\kos{K}$.
We denote the product, unit, antipode morphisms of $\cat{Spec}(\pi)$ as
$\pi\xrightarrow{\cp_{\pi}} \pi\otimes \pi$,
$\pi\xrightarrow{e_{\pi}} \kappa$,
$\pi\xrightarrow[\cong]{\varsigma_{\pi}}\pi$.
The functor $\slot\otimes\pi:\CK\to \CK$
is conservative
as explained in Remark~\ref{rem RepGrpAff(K) tensorpi conservative}.

A \emph{morphism} $\cat{Spec}(\pi^{\pr})\to\cat{Spec}(\pi)$
of pre-Grothendieck objects in $\kos{K}$
is a morphism 
$\cat{Spec}(\pi^{\pr})\to\cat{Spec}(\pi)$
of group objects in $\kos{A\!f\!f}(\kos{K})$.

A \emph{representation} of a pre-Grothendieck object $\cat{Spec}(\pi)$ in $\kos{K}$
is a representation of $\cat{Spec}(\pi)$ as a group object in $\kos{A\!f\!f}(\kos{K})$:
it is a pair $X=(x,\rho_x)$
of an object $x$ in $\CK$
and a right $\pi$-coaction morphism $x\xrightarrow{\rho_x} x\otimes \pi$
in $\CK$.
Recall the strong $\kos{K}$-tensor category
$(\kos{R\!e\!p}(\kos{K}),\kos{t}_{\pi}^*)$
of representations of $\cat{Spec}(\pi)$ in $\kos{K}$
defined in Definition~\ref{def RepGrpAff(K) Rep(pi)Ktensorcat}.

\begin{lemma} \label{lem preGroobj Rep(pi)}
  Let $\cat{Spec}(\pi)$ be a pre-Grothendieck object in $\kos{K}$
  and recall Lemma~\ref{lem RepGrpAff(K) Rep(pi)}.
  \begin{enumerate}
    \item 
    We have a connected Grothendieck $\kos{K}$-prekosmos
    $\mathfrak{Rep}(\pi)=(\kos{R\!e\!p}(\pi),\kos{t}_{\pi})$
    of representations of $\cat{Spec}(\pi)$.

    \item 
    We have a surjective pre-fiber functor
    $\varpi_{\pi}=(\omega_{\pi},\what{\omega}_{\pi}):\mathfrak{K}\to \mathfrak{Rep}(\pi)$
    whose left adjoint
    $\omega_{\pi}^*$ is the functor of forgetting the actions of $\cat{Spec}(\pi)$.
  \end{enumerate}
\end{lemma}
\begin{proof}
  Recall the symmetric monoidal category
  $\kos{R\!e\!p}(\pi)=(\cat{Rep}(\pi),\tensor\!_{\pi},\unit\!_{\pi})$
  of representations of $\cat{Spec}(\pi)$.
  The category
  $\cat{Rep}(\pi)$ has coreflexive equalizers
  since the functor $\slot\otimes\pi:\CK\to \CK$
  preserves coreflexive equalizers.
  Therefore $\kos{R\!e\!p}(\pi)$ is a Grothendieck prekosmos.
  The fully faithful strong symmetric monoidal functor
  $\kos{t}_{\pi}^*:\kos{K}\to\kos{R\!e\!p}(\pi)$
  of constructing trivial representations
  admits a right adjoint
  $\kos{t}_{\pi*}$
  which sends each representation $X=(x,\rho_x)$ of $\cat{Spec}(\pi)$
  to the following coreflexive equalizer in $\CK$.
  \begin{equation*}
    \vcenter{\hbox{
      \xymatrix@C=30pt{
        \kos{t}_{\pi*}(X)
        \ar@{^{(}->}[r]
        &x
        \ar@<0.5ex>[r]^-{\rho_x}
        \ar@<-0.5ex>[r]_-{u_{\pi}\otimes I_x}
        &\pi\otimes x
        \ar@/_2pc/@<-1ex>[l]|-{e_{\pi}\otimes I_x}
      }
    }}
  \end{equation*}
  Thus we obtain a connected Grothendieck morphism
  \begin{equation*}
    \vcenter{\hbox{
      \xymatrix{
        \kos{K}
        \ar@/_1pc/[d]_{\kos{t}^*_{\pi}}
        \\
        \kos{R\!e\!p}(\pi)
        \ar@/_1pc/[u]_{\kos{t}_{\pi*}}
      }
    }}
    \qquad\quad
    \kos{t}_{\pi}:\kos{R\!e\!p}(\pi)\to \kos{K}
  \end{equation*}
  and the pair
  $\mathfrak{Rep}(\pi):=(\kos{R\!e\!p}(\pi),\kos{t}_{\pi})$
  is a connected Grothendieck $\kos{K}$-prekosmos.
  Recall that in Lemma~\ref{lem RepGrpAff(K) Rep(pi)},
  we introduced the coreflective left-strong $\kos{K}$-tensor adjunction
  \begin{equation*}
    \vcenter{\hbox{
      \xymatrix{
        (\kos{R\!e\!p}(\pi),\kos{t}_{\pi}^*)
        \ar@/_1pc/[d]_-{(\omega_{\pi}^*,\what{\omega}_{\pi}^*\text{ }\!)}
        \\
        (\kos{K},\id_{\kos{K}})
        \ar@/_1pc/[u]_-{(\omega_{\pi*},\what{\omega}_{\pi*})}
      }
    }}
    \qquad\quad
    (\omega_{\pi},\what{\omega}_{\pi}):
    (\kos{K},\id_{\kos{K}})\to (\kos{R\!e\!p}(\pi),\kos{t}^*_{\pi})
  \end{equation*}
  satisfying the projection formula,
  such that the right adjoint
  $(\omega_{\pi*},\what{\omega}_{\pi*})$
  is a coreflective lax $\kos{K}$-tensor functor
  and the underlying functor $\omega_{\pi*}$ is conservative.
  Thus 
  \begin{equation*}
    \varpi_{\pi}=(\omega_{\pi},\what{\omega}_{\pi})
    :\mathfrak{K}\to\mathfrak{Rep}(\pi)
  \end{equation*}
  is a pre-fiber functor for $\mathfrak{Rep}(\pi)$.
  The forgetful functor
  $\omega^*$ is conservative,
  and it preserves coreflexive equalizers
  since the functor 
  $\slot\otimes\pi:\CK\to \CK$
  preserves coreflexive equalizers.
  We conclude that
  $\varpi_{\pi}:\mathfrak{K}\to\mathfrak{Rep}(\pi)$
  is a surjective pre-fiber functor.
  This completes the proof of Lemma~\ref{lem preGroobj Rep(pi)}.
\qed\end{proof}

Let $\cat{Spec}(\pi)$ be a pre-Grothendieck object in $\kos{K}$.
In (\ref{eq RepGrpAff(K) Hom(-,Spec(pi))})
we explained the presheaf of groups
$\Hom_{\cat{Aff}(\kos{K})}\bigl(\slot,\cat{Spec}(\pi)\bigr)$
on $\cat{Aff}(\kos{K})$,
where the group operation $\star$
is the convolution product.
The group
$\Hom_{\cat{Aff}(\kos{K})}\bigl(\cat{Spec}(\kappa),\cat{Spec}(\pi)\bigr)$
acts on the presheaf of groups
$\Hom_{\cat{Aff}(\kos{K})}\bigl(\slot,\cat{Spec}(\pi)\bigr)$
by conjugation.
The conjugation action of each element
$\cat{Spec}(\kappa)\xrightarrow{\theta}\cat{Spec}(\pi)$
corresponds to the automorphism 
$\sigma_{\theta}:\cat{Spec}(\pi)\xrightarrow{\cong}\cat{Spec}(\pi)$
where
\begin{equation*}
  \sigma_{\theta}:
  \xymatrix@C=50pt{
    \pi
    \ar[r]^-{\cp_{\pi}\circ (\cp_{\pi}\otimes \varsigma_{\pi})}
    &\pi\otimes \pi\otimes \pi
    \ar[r]^-{\theta\otimes I_{\pi}\otimes \theta}
    &\kappa\otimes \pi\otimes \kappa
    \ar[r]^-{(\imath^{-1}_p\otimes I_{\kappa})\circ \jmath^{-1}_p}_-{\cong}
    &\pi
    .
  }
\end{equation*}
We call $\sigma_{\theta}$ as the \emph{inner automorphism} of $\cat{Spec}(\pi)$ associated to $\theta$.
We have the relation
$\sigma_{\theta\star\tilde{\theta}}:
\cat{Spec}(\pi)\xrightarrow[\cong]{\sigma_{\tilde{\theta}}}\cat{Spec}(\pi)\xrightarrow[\cong]{\sigma_{\theta}}\cat{Spec}(\pi)$
for every pair
$\cat{Spec}(\kappa)\xrightarrow[\cong]{\theta}\cat{Spec}(\pi)$,
$\cat{Spec}(\kappa)\xrightarrow[\cong]{\tilde{\theta}}\cat{Spec}(\pi)$
and we have 
$\sigma_{e_{\pi}}=I_{\cat{Spec}(\pi)}$
where $\cat{Spec}(\kappa)\xrightarrow{e_{\pi}} \cat{Spec}(\pi)$.

Let $\cat{Spec}(\pi^{\pr})\xrightarrow{f} \cat{Spec}(\pi)$
be a morphism of pre-Grothendieck objects in $\kos{K}$.
For each $\cat{Spec}(\kappa)\xrightarrow{\theta^{\pr}}\cat{Spec}(\pi^{\pr})$,
let us denote $f\theta^{\pr}:\cat{Spec}(\kappa)\xrightarrow{\theta^{\pr}} \cat{Spec}(\pi^{\pr})\xrightarrow{f} \cat{Spec}(\pi)$.
Then the associated inner automorphisms
$\sigma^{\pr}_{\theta^{\pr}}:
\cat{Spec}(\pi^{\pr})\xrightarrow{\cong}\cat{Spec}(\pi^{\pr})$,
$\sigma_{f\theta^{\pr}}:
\cat{Spec}(\pi)\xrightarrow{\cong}\cat{Spec}(\pi)$
satisfy the following relation.
\begin{equation*}
  \vcenter{\hbox{
    \xymatrix@C=30pt{
      \cat{Spec}(\pi^{\pr})
      \ar[d]_-{\sigma^{\pr}_{\theta^{\pr}}}^-{\cong}
      \ar[r]^-{f}
      &\cat{Spec}(\pi)
      \ar[d]^-{\sigma_{f\theta^{\pr}}}_-{\cong}
      \\
      \cat{Spec}(\pi^{\pr})
      \ar[r]^-{f}
      &\cat{Spec}(\pi)
    }
  }}
\end{equation*}

\begin{definition}\label{def preGroobj GROOBJpre2cat}
  We define the $(2,1)$-category
  \begin{equation*}
    \mathbb{GROOBJ}^{\cat{pre}}(\kos{K})
  \end{equation*}
  of pre-Grothendieck objects in $\kos{K}$ as follows.
  \begin{itemize}
    \item 
    A $0$-cell is a pre-Grothendieck object $\cat{Spec}(\pi)$ in $\kos{K}$.

    \item 
    A $1$-cell $\cat{Spec}(\pi^{\pr})\xrightarrow{f} \cat{Spec}(\pi)$
    is a morphism of pre-Grothendieck objects in $\kos{K}$,
    i.e., a morphism of group objects in $\kos{A\!f\!f}(\kos{K})$.

    \item 
    A $2$-cell $\theta:f_1\Rightarrow f_2$
    between $1$-cells $f_1$, $f_2:\cat{Spec}(\pi^{\pr})\to \cat{Spec}(\pi)$
    is a morphism $\cat{Spec}(\kappa)\xrightarrow{\theta} \cat{Spec}(\pi)$
    in $\cat{Aff}(\kos{K})$
    whose associated inner automorphism
    $\sigma_{\theta}:\cat{Spec}(\pi)\xrightarrow{\cong}\cat{Spec}(\pi)$
    satisfies the relation $\sigma_{\theta}\circ f_1=f_2$.
    \begin{equation*}
      \vcenter{\hbox{
        \xymatrix@C=30pt{
          \cat{Spec}(\pi^{\pr})
          \ar@/^1pc/[r]^-{f_1}
          \ar@/_1pc/[r]_-{f_2}
          \xtwocell[r]{}<>{<0>{\text{ }\theta}}
          &\cat{Spec}(\pi)
        }
      }}
      \qquad\quad
      \vcenter{\hbox{
        \xymatrix@R=10pt@C=10pt{
          \text{ }
          &\cat{Spec}(\pi)
          \ar[dd]^-{\sigma_{\theta}}_-{\cong}
          \\
          \cat{Spec}(\pi^{\pr})
          \ar[ur]^-{f_1}
          \ar[dr]_-{f_2}
          &\text{ }
          \\
          \text{ }
          &\cat{Spec}(\pi)
        }
      }}
    \end{equation*}
    This is equivalent to saying that
    the relation $\theta\star f_1=f_2\star \theta$ holds.
    \begin{equation*}
      \vcenter{\hbox{
        \xymatrix@C=30pt{
          \kappa\otimes \pi^{\pr}
          \ar[r]^-{\imath^{-1}_{\pi^{\pr}}}_-{\cong}
          &\pi^{\pr}
          &\pi^{\pr}\otimes \kappa
          \ar[l]_-{\jmath^{-1}_{\pi^{\pr}}}^-{\cong}
          \\
          \pi\otimes \pi
          \ar[u]^-{\theta\otimes f_1}
          &\pi
          \ar[l]_-{\cp_{\pi}}
          \ar[r]^-{\cp_{\pi}}
          &\pi\otimes \pi
          \ar[u]_-{f_2\otimes \theta}
        }
      }}
    \end{equation*}

    \item 
    Identity $1$-cell of the $0$-cell $\cat{Spec}(\pi)$
    is the identity morphism $I_{\cat{Spec}(\pi)}$.
    Composition of $1$-cells is the usual composition of morphisms of pre-Grothendieck objects in $\kos{K}$.
    Identity $2$-cell of the $1$-cell
    $\cat{Spec}(\pi^{\pr})\xrightarrow{f} \cat{Spec}(\pi)$
    is $\cat{Spec}(\kappa)\xrightarrow{e_{\pi}} \cat{Spec}(\pi)$.

    \item 
    Vertical composition of $2$-cells
    $\xymatrix@C=18pt{f_1\ar@2{->}[r]^-{\theta_1}&f_2\ar@2{->}[r]^-{\theta_2}&f_3}$
    between $1$-cells $f_1$, $f_2$, $f_3:\cat{Spec}(\pi^{\pr})\to \cat{Spec}(\pi)$
    is the convolution product $\theta_2\star\theta_1:f_1\Rightarrow f_3$.
    \begin{equation*}
      \vcenter
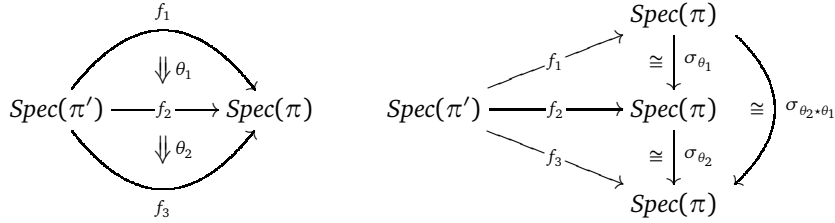
{\hbox{
        \xymatrix@C=40pt{
          \cat{Spec}(\pi^{\pr})
          \ar@/^2.5pc/[r]^-{f_1}
          \ar[r]|-{f_2}
          \ar@/_2.5pc/[r]_-{f_3}
          \xtwocell[r]{}<>{<-3>{\text{ }\text{ }\theta_1}}
          \xtwocell[r]{}<>{<3>{\text{ }\text{ }\theta_2}}
          &\cat{Spec}(\pi)
        }
      }}
      \qquad
      \vcenter{\hbox{
        \xymatrix@R=20pt@C=50pt{
          \text{ }
          &\cat{Spec}(\pi)
          \ar[d]^-{\sigma_{\theta_1}}_-{\cong}
          \ar@/^2pc/@<3ex>[dd]^-{\sigma_{\theta_2\star \theta_1}}_-{\cong}
          \\
          \cat{Spec}(\pi^{\pr})
          \ar[ur]|-{f_1}
          \ar[r]|-{f_2}
          \ar[dr]|-{f_3}
          &\cat{Spec}(\pi)
          \ar[d]^-{\sigma_{\theta_2}}_-{\cong}
          \\
          \text{ }
          &\cat{Spec}(\pi)
        }
      }}
    \end{equation*}

    \item 
    Every $2$-cell $\theta:f_1\Rightarrow f_2:\cat{Spec}(\pi^{\pr})\to \cat{Spec}(\pi)$
    is invertible, whose inverse is
    $\cat{Spec}(\kappa)\xrightarrow{\theta}\cat{Spec}(\pi)\xrightarrow[\cong]{\varsigma_{\pi}}\cat{Spec}(\pi)$.

    \item 
    Horizontal composition of $2$-cells
    $\theta:f_1\Rightarrow f_2:\cat{Spec}(\pi^{\pr})\to \cat{Spec}(\pi)$
    and
    $\theta^{\pr}:f^{\pr}_1\Rightarrow f^{\pr}_2:\cat{Spec}(\pi^{\ppr})\to \cat{Spec}(\pi^{\pr})$
    is defined as
    $\theta\centerdot\theta^{\pr}
    :=\theta\star f_1\theta^{\pr}
    =f_2\theta^{\pr}\star \theta$.
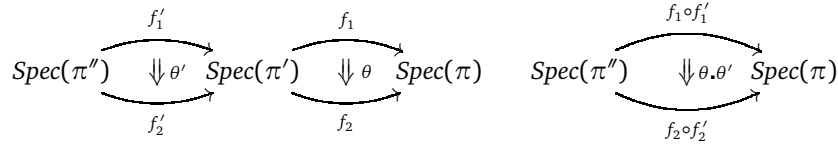
    \begin{equation*}
      \vcenter{\hbox{
        \xymatrix@C=30pt{
          \cat{Spec}(\pi^{\ppr})
          \ar@/^1pc/[r]^-{f_1^{\pr}}
          \ar@/_1pc/[r]_-{f_2^{\pr}}
          \xtwocell[r]{}<>{<0>{\text{ }\text{ }\theta^{\pr}}}
          &\cat{Spec}(\pi^{\pr})
          \ar@/^1pc/[r]^-{f_1}
          \ar@/_1pc/[r]_-{f_2}
          \xtwocell[r]{}<>{<0>{\text{ }\text{ }\theta}}
          &\cat{Spec}(\pi)
        }
      }}
      \quad
      \vcenter{\hbox{
        \xymatrix@C=40pt{
          \cat{Spec}(\pi^{\ppr})
          \ar@/^1.2pc/[r]^-{f_1\circ f_1^{\pr}}
          \ar@/_1.2pc/[r]_-{f_2\circ f_2^{\pr}}
          \xtwocell[r]{}<>{<0>{\quad\theta\centerdot \theta^{\pr}}}
          &\cat{Spec}(\pi)
        }
      }}
    \end{equation*}

    \item 
    Associators and left, right unitors are identities.
  \end{itemize}
\end{definition}

\subsection{Pre-Grothendieck $\kos{K}$-category}
\label{subsec preGroKCat}

Let $\kos{K}=(\CK,\otimes,\kappa)$ be a Grothendieck prekosmos.

\begin{definition}
  A \emph{pre-Grothendieck $\kos{K}$-category}
  is a Grothendieck $\kos{K}$-prekosmos $\mathfrak{T}=(\kos{T},\kos{t})$
  which admits a surjective pre-fiber functor
  $\varpi=(\omega,\what{\omega}):\mathfrak{K}\to \mathfrak{T}$.
\end{definition}

Let $\cat{Spec}(\pi)$ be a pre-Grothendieck object in $\kos{K}$
and recall Lemma~\ref{lem preGroobj Rep(pi)}.
Then the connected Grothendieck $\kos{K}$-prekosmos
$\mathfrak{Rep}(\pi)=(\kos{R\!e\!p}(\pi),\kos{t}_{\pi})$
of representations of $\cat{Spec}(\pi)$
is a pre-Grothendieck $\kos{K}$-category,
as it admits a surjective pre-fiber functor
$\varpi_{\pi}:\mathfrak{K}\to\mathfrak{Rep}(\pi)$.

\begin{lemma} \label{lem preGroKCat prefib is surjective}
  Let $\mathfrak{T}=(\kos{T},\kos{t})$ be a pre-Grothendieck $\kos{K}$-category.
  Then every pre-fiber functor $\mathfrak{K}\to\mathfrak{T}$ is surjective.
\end{lemma}
\begin{proof}
  From the definition of a pre-Grothendieck $\kos{K}$-category $\mathfrak{T}$,
  there exists some
  surjective pre-fiber functor
  $\varpi=(\omega,\what{\omega}):\mathfrak{K}\to \mathfrak{T}$.
  Let $\varpi^{\pr}=(\omega^{\pr},\what{\omega}^{\pr}):\mathfrak{K}\to\mathfrak{T}$
  be an arbitrary pre-fiber functor.
  We claim that $\varpi$ and $\varpi^{\pr}$ are locally isomorphic
  and conclude that $\varpi^{\pr}$ is also surjective.
  Let us explain in detail.
  By Corollary~\ref{cor Groprekosmoi summary}
  the object
  $\cat{Spec}(p)$,
  $p:=\omega^{\pr*}\omega_*(\kappa)$ in $\cat{Aff}(\kos{K})$
  represents the presheaf of invertible Grothendieck $\kos{K}$-transformations
  from $\varpi^{\pr}$ to $\varpi$,
  and the universal element is a monoidal $\kos{K}$-tensor natural isomorphism
  \begin{equation} \label{eq preGroKCat prefib is surjective}
    \xi:
    \xymatrix{
      (\kos{p}^*\omega^{\pr*},\what{\kos{p}^*\omega^{\pr*}})
      \ar@2{->}[r]^-{\cong}
      &(\kos{p}^*\omega^*,\what{\kos{p}^*\omega^*})
      :(\kos{T},\kos{t}^*)\to (\kos{K}_p,\kos{p}^*).
    }
  \end{equation}
  Note that the functor
  $\otimes p:\CK\to \CK$
  is conservative and preserves coreflexive equalizers.
  This is because we have natural isomorphisms
  \begin{equation*}
    \xymatrix@C=40pt{
      \omega^*\omega^{\pr}_*
      \ar@2{->}[r]^-{\varsigma^{\omega^{\pr},\omega}}_-{\cong}
      &\omega^{\pr*}\omega_*
      \ar@2{->}[r]^-{\bigl(\tahar{\omega^{\pr*}\omega_*}\bigr)^{-1}}_-{\cong}
      &\otimes p
      :\CK\to\CK
    }
  \end{equation*}
  where $\varsigma^{\omega^{\pr},\omega}$
  is mentioned in Corollary~\ref{cor Groprekosmoi summary},
  and the functor
  $\omega^*\omega^{\pr}_*$
  is conservative and preserves coreflexive equalizers.
  Thus the functor
  $\kos{p}^*:\kos{K}\to\kos{K}_p$
  is conservative and preserves coreflexive equalizers.
  As we have the natural isomorphism
  (\ref{eq preGroKCat prefib is surjective})
  the composition
  $\kos{p}^*\omega^{\pr*}$
  is conservative and preserves coreflexive equalizers.
  We conclude that
  $\omega^{\pr*}$
  is conservative and preserves coreflexive equalizers.
  This shows that the pre-fiber functor
  $\varpi^{\pr}:\mathfrak{K}\to\mathfrak{T}$
  is surjective.
  This completes the proof of Lemma~\ref{lem preGroKCat prefib is surjective}.
\qed\end{proof}

\begin{theorem} \label{thm preGroKCat mainThm1}
  Let $\mathfrak{T}=(\kos{T},\kos{t})$ be a pre-Grothendieck $\kos{K}$-category
  and let $\varpi:\mathfrak{K}\to\mathfrak{T}$ be a pre-fiber functor for $\mathfrak{T}$.
  \begin{enumerate}
    \item 
    We have a pre-Grothendieck object
    $\cat{Spec}(\omega^*\omega_*(\kappa))$ in $\kos{K}$
    which represents the presheaf of groups of
    invertible Grothendieck $\kos{K}$-transformations
    from $\varpi$ to $\varpi$.
    \begin{equation*}
      \vcenter{\hbox{
        \xymatrix{
          \Hom_{\cat{Aff}(\kos{K})}\bigl(\slot,\cat{Spec}(\omega^*\omega_*(\kappa))\bigr)
          \ar@2{->}[r]^-{\cong}
          &\underline{\Aut}_{\mathbb{GRO}^{\cat{pre}}_{\kos{K}}}(\varpi)
          :\cat{Aff}(\kos{K})^{\op}\to\cat{Grp}
        }
      }}
    \end{equation*}

    \item
    The given pre-fiber functor
    $\varpi:\mathfrak{K}\to \mathfrak{T}$
    is surjective,
    and factors through as an equivalence of Grothendieck $\kos{K}$-prekosmoi
    $\widebreve{\varpi}:\mathfrak{Rep}(\omega^*\omega_*(\kappa))\xrightarrow{\simeq} \mathfrak{T}$.
    \begin{equation*}
      \vcenter{\hbox{
        \xymatrix@C=40pt{
          \mathfrak{T}
          &\mathfrak{Rep}(\pi)
          \ar[l]_-{\widebreve{\varpi}}^-{\simeq}
          \\
          \text{ }
          &\mathfrak{K}
          \ar@/^1pc/[ul]^-{\varpi}
          \ar[u]_-{\varpi_{\pi}}
        }
      }}
      \qquad\quad
      \pi=\omega^*\omega_*(\kappa)
    \end{equation*}
  \end{enumerate}
\end{theorem}
\begin{proof}
  In Proposition~\ref{prop RepGrpAff(K) comparison functor}
  we showed that we have a group object
  $\cat{Spec}(\omega^*\omega_*(\kappa))$ in $\kos{A\!f\!f}(\kos{K})$
  which represents the presheaf of groups 
  $\underline{\Aut}_{\mathbb{GRO}^{\cat{pre}}_{\kos{K}}}(\varpi)$.
  The product, unit and antipode morphisms of $\cat{Spec}(\omega^*\omega_*(\kappa))$
  are explicitly given as follows.
  \begin{equation*}
    \begin{aligned}
      \cp_{\omega^*\omega_*(\kappa)}
      &:
      \xymatrix@C=40pt{
        \omega^*\omega_*(\kappa)
        \ar[r]^-{\omega^*(\eta_{\omega_*(\kappa)})}
        &\omega^*\omega_*\omega^*\omega_*(\kappa)
        \ar[r]^-{(\tahar{\omega^*\omega_*}_{\omega^*\omega_*(\kappa)})^{-1}}_-{\cong}
        &\omega^*\omega_*(\kappa)
        \otimes
        \omega^*\omega_*(\kappa)
      }
      \\
      e_{\omega^*\omega_*(\kappa)}
      &:
      \xymatrix{
        \omega^*\omega_*(\kappa)
        \ar[r]^-{\epsilon_{\kappa}}
        &\kappa
      }
      \\
      \varsigma_{\omega^*\omega_*(\kappa)}
      &:
      \xymatrix{
        \omega^*\omega_*(\kappa)
        \ar[r]^-{\varsigma^{\omega,\omega}_{\kappa}}_-{\cong}
        &\omega^*\omega_*(\kappa)
      }
    \end{aligned}
  \end{equation*}
  By Lemma~\ref{lem preGroKCat prefib is surjective},
  the given pre-fiber functor $\varpi:\mathfrak{K}\to \mathfrak{T}$
  is surjective.
  Therefore
  $\omega^*\omega_*:\kos{K}\to \kos{K}$
  is conservative and preserves coreflexive equalizers.
  As we have the natural isomorphism
  $\xymatrix@C=15pt{
    \tahar{\omega^*\omega_*}:
    \slot\otimes\omega^*\omega_*(\kappa)
    \ar@2{->}[r]^-{\cong}
    &\omega^*\omega_*
  }$
  the functor
  $\slot\otimes \omega^*\omega_*(\kappa):\CK\to \CK$
  also preserves coreflexive equalizers.
  Thus $\cat{Spec}(\omega^*\omega_*(\kappa))$
  is a pre-Grothendieck object in $\kos{K}$.
  This proves statement 1.

  Next we prove statement 2.
  Let us denote $\pi:=\omega^*\omega_*(\kappa)$.
  By Proposition~\ref{prop RepGrpAff(K) comparison functor},
  the strong $\kos{K}$-tensor functor
  $(\omega^*,\what{\omega}^*)
  :(\kos{T},\kos{t}^*)\to (\kos{K},\id_{\kos{K}})$
  factors through as a strong $\kos{K}$-tensor functor
  $(\widebreve{\omega}\!^*,\what{\widebreve{\omega}}^*)
  :(\kos{T},\kos{t}^*)\to (\kos{R\!e\!p}(\pi),\kos{t}_{\pi}^*)$.
  \begin{equation*}
    \vcenter{\hbox{
      \xymatrix@C=40pt{
        (\kos{T},\kos{t}^*)
        \ar[r]^-{(\widebreve{\omega}\!^*,\what{\widebreve{\omega}}^*)}
        \ar@/_1pc/[dr]_-{(\omega^*,\what{\omega}^*)}
        &(\kos{R\!e\!p}(\pi),\kos{t}_{\pi}^*)
        \ar[d]^-{(\omega_{\pi}^*,\what{\omega}_{\pi}^*)}
        \\
        \text{ }
        &(\kos{K},\id_{\kos{K}})
      }
    }}
  \end{equation*}  
  To conclude statement 2,
  it suffices to show that
  $\widebreve{\omega}\!^*$
  is an equivalence of categories.
  This follows from the crude comonadicity theorem.
  Since the category $\CT$ has coreflexive equalizers,
  the functor
  $\widebreve{\omega}\!^*$
  admits a right adjoint
  $\widebreve{\omega}\!_*$
  which sends each object $(x,\rho_x)$ in $\cat{Rep}(\pi)$
  to the following coreflexive equalizer in $\CT$.
  \begin{equation*}
    \vcenter{\hbox{
      \xymatrix@C=35pt{
        \widebreve{\omega}\!_*(x,\rho_x)
        \ar@{^{(}->}[r]^-{\widebreve{\textit{eq}}_{(x,\rho_x)}}
        &\omega_*(x)
        \ar@<0.5ex>[rr]^-{\eta_{\omega_*(x)}}
        \ar@<-0.5ex>[rr]_-{\omega_*\bigl(\rho_x\text{ }\!\circ\text{ }\! \tahar{\omega^*\omega_*}_x\bigr)}
        &\text{ }
        &\omega_*\omega^*\omega_*(x)
        \ar@<-1ex>@/_2pc/[ll]_-{\omega_*(\epsilon_x)}
      }
    }}
  \end{equation*}
  The component 
  $\widebreve{\epsilon}\!_{(x,\rho_x)}:
  \widebreve{\omega}\!^*\!\widebreve{\omega}\!_*(x,\rho_x)
  \to(x,\rho_x)$
  of the adjunction counit $\widebreve{\epsilon}$
  at each object $(x,\rho_x)$ in $\cat{Rep}(\pi)$
  is given by
  \begin{equation*}
    \xymatrix@C=50pt{
      \widebreve{\epsilon}\!_{(x,\rho_x)}:
      \omega^*\!\widebreve{\omega}\!_!(x,\rho_x)
      \ar@{^{(}->}[r]^-{\omega^*(\widebreve{\textit{eq}}_{(x,\rho_x)})}
      &\omega^*\omega_*(x)
      \ar[r]^-{\epsilon_x}
      &x
      .
    }
  \end{equation*}
  The component 
  $\widebreve{\eta}\!_Y:
  Y\to \widebreve{\omega}\!_*\!\widebreve{\omega}\!^*(Y)$
  of the adjunction unit $\widebreve{\eta}$
  at each object $Y$ in $\CT$ is the
  unique morphism in $\CT$ which satisfies the relation below.
  Note that we have
  $\omega^*(\eta_Y)=
  \hatar{\omega^*\omega_*}_{\omega^*(Y)}\circ \rho_{\omega^*(Y)}:\omega^*(Y)\to\omega^*\omega_*\omega^*(Y)$.
  \begin{equation*}
    \vcenter{\hbox{
      \xymatrix@C=50pt{
        Y
        \ar@{.>}[d]_-{\widebreve{\eta}\!_Y}^-{\exists!}
        \ar@/^1pc/[dr]^-{\eta_Y}
        &\text{ }
        &\text{ }
        \\
        \widebreve{\omega}\!_*\!\widebreve{\omega}\!^*(Y)
        \ar@{^{(}->}[r]^-{\widebreve{\textit{eq}}_{\widebreve{\omega}\!^*(Y)}}
        &\omega_*\omega^*(Y)
        \ar@<0.5ex>[r]^-{\eta_{\omega_*\omega^*(Y)}}
        \ar@<-0.5ex>[r]_-{\omega_*\omega^*(\eta_Y)}
        &\omega_!\omega^*\omega_!\omega^*(Y)
      }
    }}
  \end{equation*}
  One can check that $\widebreve{\eta}$, $\widebreve{\epsilon}$
  satisfy the triangle identities.
  Since $\omega^*$ preserves coreflexive equalizers,
  one can show that
  $\widebreve{\epsilon}$ is a natural isomorphism.
  As the functor $\omega^*$ is also conservative,
  we further obtain that
  $\widebreve{\eta}$
  is a natural isomorphism.
  This shows that $\widebreve{\omega}\!^*\dashv \widebreve{\omega}\!_*$
  is an adjoint equivalence of categories.
  This completes the proof of Theorem~\ref{thm preGroKCat mainThm1}.
\qed\end{proof}

\begin{corollary}
  Every pre-Grothendieck $\kos{K}$-category
  $\mathfrak{T}=(\kos{T},\kos{t})$
  is connected as a Grothendieck $\kos{K}$-prekosmos.
\end{corollary}
\begin{proof}
  This is an immediate consequence of 
  Lemma~\ref{lem preGroobj Rep(pi)}
  and 
  Theorem~\ref{thm preGroKCat mainThm1}.
\qed\end{proof}

\subsection{Pre-fundamental functors}
\label{subsec Groprefundamentalfunct}

Let $\kos{K}=(\CK,\otimes,\kappa)$ be a Grothendieck prekosmos.

In this subsection, we introduce morphisms between pre-Grothendieck $\kos{K}$-categories
and define 
the $(2,1)$-category $\mathbb{GROCAT}^{\cat{pre}}_*(\kos{K})$
of pointed pre-Grothendieck $\kos{K}$-categories.
The main goal is to establish an equivalence
of $2$-categories
between 
$\mathbb{GROOBJ}^{\cat{pre}}(\kos{K})$,
$\mathbb{GROCAT}^{\cat{pre}}_*(\kos{K})$.

\begin{definition}
  Let $\mathfrak{T}=(\kos{T},\kos{t})$, $\mathfrak{T}^{\pr}=(\kos{T}^{\pr},\kos{t}^{\pr})$
  be pre-Grothendieck $\kos{K}$-categories.
  A \emph{pre-fundamental functor} from $\mathfrak{T}^{\pr}$ to $\mathfrak{T}$
  is a Grothendieck $\kos{K}$-morphism
  $\mathfrak{f}:\mathfrak{T}^{\pr}\to\mathfrak{T}$
  which has the following property:
  there exists some pre-fiber functor
  $\varpi^{\pr}:\mathfrak{K}\to\mathfrak{T}^{\pr}$
  such that the composition 
  $\mathfrak{f}\varpi^{\pr}:\mathfrak{K}\xrightarrow{\varpi^{\pr}}\mathfrak{T}^{\pr}\xrightarrow{\mathfrak{f}}\mathfrak{T}$
  is also a pre-fiber functor.
\end{definition}

\begin{lemma} \label{lem Groprefundamentalfunct Hom=Isom}
  Let $\mathfrak{T}=(\kos{T},\kos{t})$, $\mathfrak{T}^{\pr}=(\kos{T}^{\pr},\kos{t}^{\pr})$
  be pre-Grothendieck $\kos{K}$-categories
  and let
  $\mathfrak{f}_1$, $\mathfrak{f}_2:\mathfrak{T}^{\pr}\to\mathfrak{T}$
  be pre-fundamental functors.
  Every Grothendieck $\kos{K}$-transformation
  $\vartheta:\mathfrak{f}_1\Rightarrow\mathfrak{f}_2$
  is invertible.
\end{lemma}
\begin{proof}
  Choose a pre-fiber functor $\varpi^{\pr}:\mathfrak{K}\to\mathfrak{T}^{\pr}$
  such that the composition
  $\mathfrak{f}_2\varpi^{\pr}:\mathfrak{K}\to\mathfrak{T}$
  is also a pre-fiber functor.
  Consider an arbitrary Grothendieck $\kos{K}$-transformation
  $\vartheta:\mathfrak{f}_1\Rightarrow\mathfrak{f}_2:\mathfrak{T}^{\pr}\to\mathfrak{T}$
  which is a monoidal $\kos{K}$-tensor natural transformation
  \begin{equation*}
    \vartheta:(\kos{f}_1^*,\what{\kos{f}_1^*})\Rightarrow (\kos{f}_2^*,\what{\kos{f}_2^*})
    :(\kos{T},\kos{t}^*)\to (\kos{T}^{\pr},\kos{t}^{\pr*})
    .
  \end{equation*}
  By Lemma~\ref{lem preGroKCat prefib is surjective},
  the pre-fiber functor $\varpi^{\pr}$ is surjective.
  In particular, the left adjoint $\omega^{\pr*}$ is conservative.
  Thus to conclude that $\vartheta$ is invertible,
  it suffices to show that
  \begin{equation*}
    \omega^{\pr*}\vartheta:
    (\omega^{\pr*}\kos{f}_1^*,\what{\omega^{\pr*}\kos{f}_1^*})\Rightarrow (\omega^{\pr*}\kos{f}_2^*,\what{\omega^{\pr*}\kos{f}_2^*})
    :(\kos{T},\kos{t}^*)\to (\kos{K},\id_{\kos{K}})
  \end{equation*}
  is invertible.
  This follows from Proposition~\ref{prop LSKtensoradj HomRep,Hom=Isom},
  since 
  $(\kos{f}_2\omega^{\pr},\what{\kos{f}_2\omega^{\pr}}):(\kos{K},\id_{\kos{K}})\to (\kos{T},\kos{t}^*)$
  is a coreflective left-strong $\kos{K}$-tensor adjunction
  satisfying the projection formula,
  and
  $(\omega^{\pr*}\kos{f}_1^*,\what{\omega^{\pr*}\kos{f}_1^*}):(\kos{T},\kos{t}^*)\to (\kos{K},\id_{\kos{K}})$
  is a strong $\kos{K}$-tensor functor.
  \begin{equation*}
    \vcenter{\hbox{
      \xymatrix{
        (\kos{T},\kos{t}^*)
        \ar[d]_-{(\omega^{\pr*}\kos{f}_1^*\text{ }\!,\text{ }\!\what{\omega^{\pr*}\kos{f}_1^*})}
        \\
        (\kos{K},\id_{\kos{K}})
      }
    }}
    \qquad\qquad\quad
    \vcenter{\hbox{
      \xymatrix{
        (\kos{T},\kos{t}^*)
        \ar@/_1pc/[d]_-{(\omega^{\pr*}\kos{f}_2^*\text{ }\!,\text{ }\!\what{\omega^{\pr*}\kos{f}_2^*})}
        \\
        (\kos{K},\id_{\kos{K}})
        \ar@/_1pc/[u]_-{(\kos{f}_{2*}\omega_*^{\pr}\text{ }\!,\text{ }\!\what{\kos{f}_{2*}\omega_*^{\pr}})}
      }
    }}
  \end{equation*}
  This completes the proof of Lemma~\ref{lem Groprefundamentalfunct Hom=Isom}.
\qed\end{proof}

\begin{definition}
  We define the $(2,1)$-category 
  \begin{equation*}
    \mathbb{GROCAT}^{\cat{pre}}_*(\kos{K})
  \end{equation*}
  of pre-Grothendieck $\kos{K}$-categories pointed with pre-fiber functors
  as follows.
  \begin{itemize}
    \item 
    A $0$-cell is a pair $(\mathfrak{T},\varpi)$
    of a pre-Grothendieck $\kos{K}$-category $\mathfrak{T}$
    and a pre-fiber functor $\varpi:\mathfrak{K}\to\mathfrak{T}$.

    \item 
    A $1$-cell $(\mathfrak{f},\alpha):(\mathfrak{T}^{\pr},\varpi^{\pr})\to (\mathfrak{T},\varpi)$
    is a pair of a pre-fundamental functor $\mathfrak{f}:\mathfrak{T}^{\pr}\to\mathfrak{T}$
    and an invertible Grothendieck $\kos{K}$-transformation
    $\xymatrix@C=15pt{\alpha:\mathfrak{f}\varpi^{\pr}\ar@2{->}[r]^-{\cong}&\varpi.}$
    \begin{equation*}
      \vcenter{\hbox{
        \xymatrix@C=15pt{
          \mathfrak{T}
          &\text{ }
          \xtwocell[d]{}<>{<0>{\alpha}}
          &\mathfrak{T}^{\pr}
          \ar[ll]_-{\mathfrak{f}}
          \\
          \text{ }
          &\mathfrak{K}
          \ar[ur]_-{\varpi^{\pr}}
          \ar[ul]^-{\varpi}
        }
      }}
      \qquad\quad
      \alpha:
      \xymatrix@C=18pt{
        \omega^{\pr*}\kos{f}^*
        \ar@2{->}[r]^-{\cong}
        &\omega^*
        :\kos{T}\to\kos{K}
      }
      \qquad
      \vcenter{\hbox{
        \xymatrix{
          \kos{T}
          \ar[d]_-{\kos{f}^*}
          \ar[r]^-{\omega^*}
          &\kos{K}
          \ar@{=}[d]
          \\
          \kos{T}^{\pr}
          \ar[r]_-{\omega^{\pr*}}
          &\kos{K}
          \xtwocell[ul]{}<>{<0>{\alpha\text{ }\text{ }}}
        }
      }}
    \end{equation*}
    We denote the mate associated to $\alpha$ as
    \begin{equation*}
      \vcenter{\hbox{
        \xymatrix{
          \kos{T}
          \ar[d]_-{\kos{f}^*}
          &\kos{K}
          \ar[l]_-{\omega_*}
          \ar@{=}[d]
          \\
          \kos{T}^{\pr}
          \xtwocell[ur]{}<>{<0>{\text{ }\text{ }\alpha_*}}
          &\kos{K}
          \ar[l]^-{\omega^{\pr}_*}
        }
      }}
      \qquad
      \alpha_*:
      \xymatrix{
        \kos{f}^*\omega_*
        \ar@2{->}[r]^-{\eta^{\pr}\kos{f}^*\omega_*}
        &\omega^{\pr}_*\omega^{\pr*}\kos{f}^*\omega_*
        \ar@2{->}[r]^-{\omega^{\pr}_*\alpha\omega_*}_-{\cong}
        &\omega^{\pr}_*\omega^*\omega_*
        \ar@2{->}[r]^-{\omega^{\pr}_*\epsilon}
        &\omega^{\pr}_*
        .
      }
    \end{equation*}

    \item
    A $2$-cell $\vartheta:(\mathfrak{f}_1,\alpha_1)\Rightarrow (\mathfrak{f}_2,\alpha_2)$
    between $1$-cells
    $(\mathfrak{f}_1,\alpha_1)$, $(\mathfrak{f}_2,\alpha_2):(\mathfrak{T}^{\pr},\varpi^{\pr})\to (\mathfrak{T},\varpi)$
    is a Grothendieck $\kos{K}$-transformation
    $\vartheta:\mathfrak{f}_1\Rightarrow\mathfrak{f}_2$
    between pre-fundamental functors.
    It is automatically an invertible Grothendieck $\kos{K}$-transformation
    by Lemma~\ref{lem Groprefundamentalfunct Hom=Isom}.
    \begin{equation*}
      \vcenter{\hbox{
        \xymatrix@C=35pt{
          (\mathfrak{T},\varpi)
          \xtwocell[r]{}<>{<0>{\text{ }\text{ }\vartheta}}
          &(\mathfrak{T}^{\pr},\varpi^{\pr})
          \ar@/_1pc/@<-0.5ex>[l]_-{(\mathfrak{f}_1,\alpha_1)}
          \ar@/^1pc/@<0.5ex>[l]^-{(\mathfrak{f}_2,\alpha_2)}
        }
      }}
      \qquad\quad
      \vartheta:
      \xymatrix@C=15pt{
        \mathfrak{f}_1
        \ar@2{->}[r]^-{\cong}
        &\mathfrak{f}_2
      }
    \end{equation*}

    \item
    Identity $1$-cell of a $0$-cell $(\mathfrak{T},\varpi)$
    is a pair of an identity fundamental functor $I_{\mathfrak{T}}:\mathfrak{T}\to\mathfrak{T}$
    and the identity Grothendieck $\kos{K}$-transformation
    $\varpi=I_{\mathfrak{T}}\varpi$.
    Identity $2$-cell of a $1$-cell
    $(\mathfrak{f},\alpha):(\mathfrak{T}^{\pr},\varpi^{\pr})\to (\mathfrak{T},\varpi)$
    is the identity Grothendieck $\kos{K}$-transformation
    $I_{\mathfrak{f}}:\mathfrak{f}=\mathfrak{f}$.

    \item 
    Composition of $1$-cells is given as follows.
    \begin{equation*}
      \vcenter{\hbox{
        \xymatrix@C=20pt{
          (\mathfrak{T},\varpi)
          &(\mathfrak{T}^{\pr},\varpi^{\pr})
          \ar[l]_-{(\mathfrak{f},\alpha)}
          &(\mathfrak{T}^{\ppr},\varpi^{\ppr})
          \ar[l]_-{(\mathfrak{f}^{\pr},\alpha^{\pr})}
          \ar@/^1pc/@<0.5ex>[ll]^-{(\mathfrak{f},\alpha)\circ(\mathfrak{f}^{\pr},\alpha^{\pr})=\bigl(\mathfrak{f}\mathfrak{f}^{\pr},\alpha\circ (\mathfrak{f}\alpha^{\pr})\bigr)}
        }
      }}
      \quad
      \vcenter{\hbox{
        \xymatrix@C=15pt{
          \alpha\circ (\mathfrak{f}\alpha^{\pr}):
          \mathfrak{f}\mathfrak{f}^{\pr}\varpi^{\ppr}
          \ar@2{->}[r]^-{\mathfrak{f}\alpha^{\pr}}_-{\cong}
          &\mathfrak{f}\varpi^{\pr}
          \ar@2{->}[r]^-{\alpha}_-{\cong}
          &\varpi
        }
      }}
    \end{equation*}

    \item
    Vertical composition of $2$-cells
    is given by vertical composition as Grothendieck $\kos{K}$-transformations.
    \begin{equation*}
      \vcenter{\hbox{
        \xymatrix@C=50pt{
          (\mathfrak{T},\varpi)
          \xtwocell[r]{}<>{<-3>{\text{ }\text{ }\vartheta_1}}
          \xtwocell[r]{}<>{<3>{\text{ }\text{ }\vartheta_2}}
          &(\mathfrak{T}^{\pr},\varpi^{\pr})
          \ar@/_2pc/@<-0.5ex>[l]_-{(\mathfrak{f}_1,\alpha_1)}
          \ar[l]|-{(\mathfrak{f}_2,\alpha_2)}
          \ar@/^2pc/@<0.5ex>[l]^-{(\mathfrak{f}_3,\alpha_3)}
        }
      }}
      \qquad\quad
      \begin{aligned}
        \vartheta_2\circ\vartheta_1
        &:
        \xymatrix@C=18pt{
          \mathfrak{f}_1
          \ar@2{->}[r]^-{\vartheta_1}_-{\cong}
          &\mathfrak{f}_2
          \ar@2{->}[r]^-{\vartheta_2}_-{\cong}
          &\mathfrak{f}_3
        }
      \end{aligned}
    \end{equation*}

    \item 
    Horizontal composition of $2$-cells 
    is given by horizontal composition
    as Grothendieck $\kos{K}$-transformations.
    \begin{equation*}
      \vcenter{\hbox{
        \xymatrix@C=35pt{
          (\mathfrak{T},\varpi)
          \xtwocell[r]{}<>{<0>{\text{ }\text{ }\vartheta}}
          &(\mathfrak{T}^{\pr},\varpi^{\pr})
          \ar@/_1pc/@<-0.5ex>[l]_-{(\mathfrak{f}_1,\alpha_1)}
          \ar@/^1pc/@<0.5ex>[l]^-{(\mathfrak{f}_2,\alpha_2)}
          \xtwocell[r]{}<>{<0>{\text{ }\text{ }\text{ }\vartheta^{\pr}}}
          &(\mathfrak{T}^{\ppr},\varpi^{\ppr})
          \ar@/_1pc/@<-0.5ex>[l]_-{(\mathfrak{f}^{\pr}_1,\alpha^{\pr}_1)}
          \ar@/^1pc/@<0.5ex>[l]^-{(\mathfrak{f}^{\pr}_2,\alpha^{\pr}_2)}
        }
      }}
    \end{equation*}
    \begin{equation*}
      \vcenter{\hbox{
        \xymatrix@C=50pt{
          (\mathfrak{T},\varpi)
          \xtwocell[r]{}<>{<0>{\quad\vartheta\vartheta^{\pr}}}
          &(\mathfrak{T}^{\ppr},\varpi^{\ppr})
          \ar@/_1pc/@<-0.5ex>[l]_-{(\mathfrak{f}_1,\alpha_1)\circ(\mathfrak{f}^{\pr}_1,\alpha^{\pr}_1)}
          \ar@/^1pc/@<0.5ex>[l]^-{(\mathfrak{f}_2,\alpha_2)\circ(\mathfrak{f}^{\pr}_2,\alpha^{\pr}_2)}
        }
      }}
      \qquad\quad
      \begin{aligned}
        \vartheta\vartheta^{\pr}
        &:
        \mathfrak{f}_1\mathfrak{f}_1^{\pr}
        \Rightarrow
        \mathfrak{f}_2\mathfrak{f}_2^{\pr}
        :\mathfrak{T}^{\ppr}\to\mathfrak{T}
        \\
        \vartheta^{\pr}\vartheta
        &:
        \kos{f}_1^{\pr*}\kos{f}_1^*
        \Rightarrow
        \kos{f}_2^{\pr*}\kos{f}_2^*
        :\kos{T}\to\kos{T}^{\ppr}
      \end{aligned}
    \end{equation*}
  \end{itemize}
\end{definition}

\begin{proposition}\label{prop GroObjtoGroCat}
  Let $\kos{K}$ be a Grothendieck prekosmos $\kos{K}$.
  We have a $2$-functor
  \begin{equation*}
    \vcenter{\hbox{
      \xymatrix@C=30pt{
        \mathbb{GROOBJ}^{\cat{pre}}(\kos{K})
        \ar[r]
        &\mathbb{GROCAT}^{\cat{pre}}_*(\kos{K}).
      }
    }}
  \end{equation*}
\end{proposition}
\begin{proof}
  First, we describe how the $2$-functor
  $\mathbb{GROOBJ}^{\cat{pre}}(\kos{K})
  \to \mathbb{GROCAT}^{\cat{pre}}_*(\kos{K})$
  sends $0$-cells to $0$-cells.
  For each pre-Grothendieck object $\cat{Spec}(\pi)$ in $\kos{K}$,
  we have the pointed pre-Grothendieck $\kos{K}$-category
  $(\mathfrak{Rep}(\pi),\varpi_{\pi})$
  which we described in Lemma~\ref{lem preGroobj Rep(pi)}.
  
  Next, we describe how the $2$-functor
  $\mathbb{GROOBJ}^{\cat{pre}}(\kos{K})
  \to \mathbb{GROCAT}^{\cat{pre}}_*(\kos{K})$
  sends
  $1$-cells to $1$-cells.
  Let $\cat{Spec}(\pi^{\pr})\xrightarrow{f}\cat{Spec}(\pi)$
  be a morphism of pre-Grothendieck objects in $\kos{K}$.
  Then we have a Grothendieck $\kos{K}$-morphism
  $\mathfrak{f}=(\kos{f},\what{\kos{f}}):
  \mathfrak{Rep}(\pi^{\pr})\to \mathfrak{Rep}(\pi)$
  and an invertible Grothendieck $\kos{K}$-transformation
  $I:\mathfrak{f}\varpi_{\pi^{\pr}}\cong \varpi_{\pi}$.
  \begin{equation*}
    \vcenter{\hbox{
      \xymatrix@C=15pt{
        \mathfrak{Rep}(\pi)
        &\text{ }
        \xtwocell[d]{}<>{<0>{I}}
        &\mathfrak{Rep}(\pi^{\pr})
        \ar[ll]_-{\mathfrak{f}}
        \\
        \text{ }
        &\mathfrak{K}
        \ar[ul]^-{\varpi_{\pi}}
        \ar[ur]_-{\varpi_{\pi^{\pr}}}
      }
    }}
  \end{equation*}
  Let us explain the underlying Grothendieck morphism
  $\kos{f}:\kos{R\!e\!p}(\pi^{\pr})\to \kos{R\!e\!p}(\pi)$.
  The left adjoint $\kos{f}^*:\kos{R\!e\!p}(\pi)\to \kos{R\!e\!p}(\pi^{\pr})$
  sends each object $X=(x,\rho_x)$ in $\cat{Rep}(\pi)$ to
  $\kos{f}^*(X)=
  (x,\xymatrix@C=15pt{
    x
    \ar[r]^-{\rho_x}
    &x\otimes \pi
    \ar[r]^-{I_x\otimes f}
    &x\otimes \pi^{\pr}
  }\!)$.
  The monoidal coherence isomorphisms of $\kos{f}^*$ are given by identity morphisms.
  The right adjoint $\kos{f}_*$ sends
  each object
  $Z=(z,\rho^{\pr}_z)$ in $\cat{Rep}(\pi^{\pr})$ to
  $\kos{f}_*(Z)=(z\otimes_{\pi^{\pr}}\pi,\rho_{z\otimes_{\pi^{\pr}}\pi})$
  whose underlying object in $\CK$ is the coreflexive equalizer in $\CK$
  \begin{equation*}
    \vcenter{\hbox{
      \xymatrix@C=30pt{
        z\otimes_{\pi^{\pr}}\pi
        \ar@{^{(}->}[r]^-{\textit{eq}^f_Z}
        &z\otimes \pi
        \ar@<0.5ex>[rr]^-{\rho^{\pr}_z\otimes I_{\pi}}
        \ar@<-0.5ex>[rr]_-{(I_z\otimes f\otimes I_{\pi})\circ (I_z\otimes \cp_{\pi})}
        &\text{ }
        &z\otimes \pi^{\pr}\otimes \pi
        \ar@/_2pc/@<-1ex>[ll]|-{I_{\pi}\otimes e_{\pi^{\pr}}\otimes I_z}
      }
    }}
  \end{equation*}
  and the right $\pi$-coaction
  $\rho_{z\otimes_{\pi^{\pr}}\pi}$
  is the unique morphism in $\CK$ satisfying the relation below.
  \begin{equation*}
    \vcenter{\hbox{
      \xymatrix@C=50pt{
        z\otimes_{\pi^{\pr}}\pi
        \ar@{^{(}->}[r]^-{\textit{eq}^f_Z}
        \ar@{.>}[d]^-{\rho_{z\otimes_{\pi^{\pr}}\pi}}_-{\exists!}
        &z\otimes \pi
        \ar[d]^-{I_z\otimes \cp_{\pi}}
        \\
        (z\otimes_{\pi^{\pr}}\pi)\otimes \pi
        \ar@{^{(}->}[r]^-{\textit{eq}^f_Z\otimes I_{\pi}}
        &z\otimes \pi\otimes \pi
      }
    }}
  \end{equation*}
  The component of the adjunction counit $\epsilon^f$ at
  each object $Z=(z,\rho^{\pr}_z)$ in $\cat{Rep}(\pi^{\pr})$ is
  the morphism $\epsilon^f_Z:\kos{f}^*\kos{f}_*(Z)\to Z$ in $\cat{Rep}(\pi^{\pr})$
  whose underlying morphism in $\CK$ is
  \begin{equation*}
    \xymatrix{
      \epsilon^f_Z:
      z\otimes_{\pi^{\pr}}\pi
      \ar@{^{(}->}[r]^-{\textit{eq}^f_Z}
      &z\otimes \pi
      \ar[r]^-{I_z\otimes e_{\pi}}
      &z\otimes \kappa
      \ar[r]^-{\jmath_z^{-1}}_-{\cong}
      &z
      .
    }
  \end{equation*}
  The component of the adjunction unit $\eta^f$
  at each object $X=(x,\gamma_x)$ in $\cat{Rep}(\pi)$
  is the unique morphism
  $\eta^f_X:X\to \kos{f}_*\kos{f}^*(X)$
  in $\cat{Rep}(\pi)$
  whose underlying morphism in $\CK$ satisfies the following relation.
  \begin{equation*}
    \vcenter{\hbox{
      \xymatrix@C=40pt{
        x
        \ar@{.>}[d]^-{\eta^f_X}_-{\exists!}
        \ar@/^1pc/[dr]^-{\rho_x}
        &\text{ }
        \\
        x\otimes \pi
        \ar@{^{(}->}[r]^-{\textit{eq}^f_{\kos{f}^*(X)}}
        &x\otimes_{\pi^{\pr}}\pi
      }
    }}
  \end{equation*}
  This is the description of
  the Grothendieck morphism
  $\kos{f}:\kos{R\!e\!p}(\pi^{\pr})\to \kos{R\!e\!p}(\pi)$.
  The invertible Grothendieck morphism
  $\what{\kos{f}}^*:
  \kos{t}_{\pi^{\pr}}\cong \kos{t}_{\pi}\kos{f}$
  is given by identity natural transformation
  $\what{\kos{f}}^*:\kos{t}^*_{\pi^{\pr}}=\kos{f}^*\kos{t}^*_{\pi}$.
  This is the description of the Grothendieck $\kos{K}$-morphism
  $\mathfrak{f}=(\kos{f},\what{\kos{f}}):\mathfrak{Rep}(\pi^{\pr})\to \mathfrak{Rep}(\pi)$.
  Finally, the invertible Grothendieck $\kos{K}$-transformation
  $I:\mathfrak{f}\varpi_{\pi^{\pr}}\cong \varpi_{\pi}$
  is also given by identity natural transformation
  $I:
  (\omega_{\pi^{\pr}}^*\kos{f}^*,\what{\omega_{\pi^{\pr}}^*\kos{f}^*})
  =(\omega_{\pi}^*,\what{\omega}_{\pi}^*)$.

  Now we describe how the $2$-functor
  $\mathbb{GROOBJ}^{\cat{pre}}(\kos{K})
  \to \mathbb{GROCAT}^{\cat{pre}}_*(\kos{K})$
  sends $2$-cells to $2$-cells.
  Let $\theta:f_1\Rightarrow f_2$ be a $2$-cell between
  morphisms $f_1$, $f_2:\cat{Spec}(\pi^{\pr})\to\cat{Spec}(\pi)$
  of pre-Grothendieck objects $\cat{Spec}(\pi)$, $\cat{Spec}(\pi^{\pr})$ in $\kos{K}$.
  We have an invertible Grothendieck $\kos{K}$-transformation
  \begin{equation*}
    \vcenter{\hbox{
      \xymatrix@C=40pt{
        (\mathfrak{Rep}(\pi),\varpi_{\pi})
        \xtwocell[r]{}<>{<0>{\text{ }\text{ }\vartheta}}
        &(\mathfrak{Rep}(\pi^{\pr}),\varpi_{\pi^{\pr}})
        \ar@/_1pc/[l]_-{(\mathfrak{f}_1,I)}
        \ar@/^1pc/[l]^-{(\mathfrak{f}_2,I)}
      }
    }}
    \qquad\quad
    \vartheta:
    \xymatrix@C=18pt{
      \mathfrak{f}_1
      \ar@2{->}[r]^-{\cong}
      &\mathfrak{f}_2
    }
  \end{equation*}
  whose component $\vartheta_X:\kos{f}_1^*(X)\xrightarrow[]{\cong} \kos{f}_2^*(X)$
  at each object $X=(x,\rho_x)$ in $\cat{Rep}(\pi)$ is
  \begin{equation*}
    \vartheta_X:
    \xymatrix{
      x
      \ar[r]^-{\rho_x}
      &x\otimes \pi
      \ar[r]^-{I_x\otimes \theta}
      &x\otimes \kappa
      \ar[r]^-{\jmath^{-1}_x}_-{\cong}
      &x
      .
    }
  \end{equation*}
  We leave for the readers to check that
  the correspondence
  $\mathbb{GROOBJ}^{\cat{pre}}(\kos{K})
  \to \mathbb{GROCAT}^{\cat{pre}}_*(\kos{K})$
  preserves identity $1$-cells, identity $2$-cells, horizontal composition of $1$-cells,
  vertical composition of $2$-cells 
  and 
  horizontal composition of $2$-cells.
  This completes the proof of Proposition~\ref{prop GroObjtoGroCat}.
\qed\end{proof}

\begin{proposition}
  \label{prop GroCattoGroObj}
  Let $\kos{K}$ be a Grothendieck prekosmos $\kos{K}$.
  We have a $2$-functor
  \begin{equation*}
    \vcenter{\hbox{
      \xymatrix@C=30pt{
        \mathbb{GROCAT}^{\cat{pre}}_*(\kos{K})
        \ar[r]
        &\mathbb{GROOBJ}^{\cat{pre}}(\kos{K}).
      }
    }}
  \end{equation*}
\end{proposition}
\begin{proof}
  First, we describe how the $2$-functor
  $\mathbb{GROCAT}^{\cat{pre}}_*(\kos{K})
  \to \mathbb{GROOBJ}^{\cat{pre}}(\kos{K})$
  sends $0$-cells to $0$-cells.
  For each pointed pre-Grothendieck $\kos{K}$-category
  $(\mathfrak{T},\varpi)$
  we have a pre-Grothendieck object
  $\cat{Spec}(\omega^*\omega_*(\kappa))$ in $\kos{K}$
  as explained in Theorem~\ref{thm preGroKCat mainThm1}.
  
  Next, we describe how the $2$-functor
  $\mathbb{GROCAT}^{\cat{pre}}_*(\kos{K})
  \to \mathbb{GROOBJ}^{\cat{pre}}(\kos{K})$
  sends $1$-cells to $1$-cells.
  Let 
  $(\mathfrak{f},\alpha):(\mathfrak{T}^{\pr},\varpi^{\pr})\to (\mathfrak{T},\varpi)$
  be a $1$-cell between pointed pre-Grothendieck $\kos{K}$-categories.
  \begin{equation*}
    \vcenter{\hbox{
      \xymatrix@C=15pt{
        \mathfrak{T}
        &\text{ }
        \xtwocell[d]{}<>{<0>{\alpha}}
        &\mathfrak{T}^{\pr}
        \ar[ll]_-{\mathfrak{f}}
        \\
        \text{ }
        &\mathfrak{K}
        \ar[ur]_-{\varpi^{\pr}}
        \ar[ul]^-{\varpi}
      }
    }}
    \qquad\quad
    \alpha:
    \xymatrix@C=18pt{
      \omega^{\pr*}\kos{f}^*
      \ar@2{->}[r]^-{\cong}
      &\omega^*
      :\kos{T}\to\kos{K}
    }
    \qquad
    \vcenter{\hbox{
      \xymatrix{
        \kos{T}
        \ar[d]_-{\kos{f}^*}
        \ar[r]^-{\omega^*}
        &\kos{K}
        \ar@{=}[d]
        \\
        \kos{T}^{\pr}
        \ar[r]_-{\omega^{\pr*}}
        &\kos{K}
        \xtwocell[ul]{}<>{<0>{\alpha\text{ }\text{ }}}
      }
    }}
  \end{equation*}
  We denote the mate associated to $\alpha$ as
  \begin{equation*}
    \vcenter{\hbox{
      \xymatrix{
        \kos{T}
        \ar[d]_-{\kos{f}^*}
        &\kos{K}
        \ar[l]_-{\omega_*}
        \ar@{=}[d]
        \\
        \kos{T}^{\pr}
        \xtwocell[ur]{}<>{<0>{\text{ }\text{ }\alpha_*}}
        &\kos{K}
        \ar[l]^-{\omega^{\pr}_*}
      }
    }}
    \qquad
    \alpha_*:
    \xymatrix{
      \kos{f}^*\omega_*
      \ar@2{->}[r]^-{\eta^{\pr}\kos{f}^*\omega_*}
      &\omega^{\pr}_*\omega^{\pr*}\kos{f}^*\omega_*
      \ar@2{->}[r]^-{\omega^{\pr}_*\alpha\omega_*}_-{\cong}
      &\omega^{\pr}_*\omega^*\omega_*
      \ar@2{->}[r]^-{\omega^{\pr}_*\epsilon}
      &\omega^{\pr}_*
      .
    }
  \end{equation*}
  If we denote
  \begin{equation*}
    \check{\alpha}:
    \xymatrix@C=30pt{
      \omega^*\omega_*
      \ar@2{->}[r]^-{\alpha^{-1}\omega_*}_-{\cong}
      &\omega^{\pr*}\kos{f}^*\omega_*
      \ar@2{->}[r]^-{\omega^{\pr*}\alpha_*}
      &\omega^{\pr*}\omega^{\pr}_*
    }
  \end{equation*}
  then
  $\check{\alpha}_{\kappa}:
  \cat{Spec}(\omega^{\pr*}\omega^{\pr}_!(\kappa))
  \to \cat{Spec}(\omega^*\omega_!(\kappa))$
  is a morphism of pre-Grothendieck objects in $\kos{K}$.

  We also describe how the $2$-functor
  $\mathbb{GROCAT}^{\cat{pre}}_*(\kos{K})
  \to \mathbb{GROOBJ}^{\cat{pre}}(\kos{K})$
  sends $2$-cells to $2$-cells.
  Let 
  $\vartheta:(\mathfrak{f}_1,\alpha_1)\Rightarrow (\mathfrak{f}_2,\alpha_2)$
  be a $2$-cell between $1$-cells
  $(\mathfrak{f}_1,\alpha_1)$, $(\mathfrak{f}_2,\alpha_2):(\mathfrak{T}^{\pr},\varpi^{\pr})\to (\mathfrak{T},\varpi)$.
  If we denote
  \begin{equation*}
    \check{\vartheta}:
    \xymatrix@C=30pt{
      \omega^*\omega_*
      \ar@2{->}[r]^-{\alpha_1^{-1}\omega_*}_-{\cong}
      &\omega^{\pr*}\kos{f}_1^*\omega_*
      \ar@2{->}[r]^-{\omega^{\pr*}\vartheta\omega_*}_-{\cong}
      &\omega^{\pr*}\kos{f}_2^*\omega_*
      \ar@2{->}[r]^-{\alpha_2\omega_*}_-{\cong}
      &\omega^*\omega_*
      \ar@2{->}[r]^-{\epsilon}
      &I_{\kos{K}}
    }
  \end{equation*}
  then 
  $\check{\vartheta}_{\kappa}:\omega^*\omega_*(\kappa)\to\kappa$
  is a $2$-cell
  $\check{\vartheta}_{\kappa}:
  \check{\alpha}_{1\kappa}\Rightarrow \check{\alpha}_{2\kappa}
  :\cat{Spec}(\omega^{\pr*}\omega^{\pr}_*(\kappa))
  \to \cat{Spec}(\omega^*\omega_*(\kappa))$
  between morphisms of pre-Grothendieck objects in $\kos{K}$.

  We leave for the readers to check that
  the correspondence
  $\mathbb{GROCAT}^{\cat{pre}}_*(\kos{K})
  \to \mathbb{GROOBJ}^{\cat{pre}}(\kos{K})$
  preserves identity $1$-cells,
  identity $2$-cells,
  vertical composition of $2$-cells,
  composition of $1$-cells 
  and horizontal composition of $2$-cells.
  This completes the proof of Proposition~\ref{prop GroCattoGroObj}.
\qed\end{proof}

\begin{theorem}
  \label{thm GROOBJbiequivGROCAT}
  Let $\kos{K}$ be a Grothendieck prekosmos $\kos{K}$.
  We have a biequivalence of $(2,1)$-categories
  \begin{equation*}
    \vcenter{\hbox{
      \xymatrix@C=30pt{
        \mathbb{GROOBJ}^{\cat{pre}}(\kos{K})
        \ar@<0.5ex>[r]^-{\simeq}
        &\mathbb{GROCAT}^{\cat{pre}}_*(\kos{K}).
        \ar@<0.5ex>[l]^-{\simeq}
      }
    }}
  \end{equation*}
\end{theorem}
\begin{proof}
  We explain why the $2$-functors described in 
  Proposition~\ref{prop inducedfundamentalfunctor}
  and
  Proposition~\ref{prop CATtoOBJ}
  are quasi-inverse to each other.
  It is clear that the composition $2$-functor
  \begin{equation*}
    \mathbb{GROOBJ}^{\cat{pre}}(\kos{K})
    \to 
    \mathbb{GROCAT}^{\cat{pre}}_*(\kos{K})
    \to
    \mathbb{GROOBJ}^{\cat{pre}}(\kos{K})
  \end{equation*}
  is $2$-naturally isomorphic to the identity $2$-functor
  of $\mathbb{GROOBJ}^{\cat{pre}}(\kos{K})$.
  Consider the composition $2$-functor
  \begin{equation}\label{eq GroCat composition 2-functor}
    \mathbb{GROCAT}^{\cat{pre}}_*(\kos{K})\to
    \mathbb{GROOBJ}^{\cat{pre}}(\kos{K})\to 
    \mathbb{GROCAT}^{\cat{pre}}_*(\kos{K})
    .
  \end{equation}
  We have the following
  weak $2$-natural transformation
  from the composition $2$-functor
  (\ref{eq GroCat composition 2-functor})
  to the identity $2$-functor of
  $\mathbb{GROCAT}^{\cat{pre}}_*(\kos{K})$.
  By the term `weak $2$-natural transformation'
  we mean that we associate
  to each $1$-cell
  $(\mathfrak{f},\alpha)$
  a constraint $2$-cell $\widebreve{\alpha}$
  which we explain below.
  \begin{itemize}
    \item 
    The component of the weak $2$-natural transformation
    at each $0$-cell $(\mathfrak{T},\varpi)$
    is the $1$-cell
    \begin{equation}\label{eq GRO strong transformation component}
      \vcenter{\hbox{
        \xymatrix@C=30pt{
          \mathfrak{T}
          &\mathfrak{Rep}(\omega^*\omega_*(\kappa))
          \ar[l]_-{\widebreve{\varpi}}^-{\simeq}
          \xtwocell[d]{}<>{<5>{I}} 
          \\
          \text{ }
          &\mathfrak{K}
          \ar@/^1pc/[ul]^-{\varpi}
          \ar[u]_-{\varpi_{\omega^*\omega_*(\kappa)}}
        }
      }}
      \quad
      \xymatrix@C=15pt{
        (\widebreve{\varpi},I):
        \bigl(\mathfrak{Rep}(\omega^*\omega_*(\kappa)),\varpi_{\omega^*\omega_*(\kappa)}\bigr)
        \ar[r]^-{\simeq}
        &(\mathfrak{T},\varpi)
      }
    \end{equation}
    where 
    $\widebreve{\varpi}$
    is the equivalence of
    Grothendieck $\kos{K}$-prekosmoi
    described in Theorem~\ref{thm preGroKCat mainThm1}.
    The invertible Grothendieck $\kos{K}$-transformation
    $I:\widebreve{\varpi}\varpi_{\omega^*\omega_*(\kappa)}\Rightarrow\varpi$
    is the identity monoidal $\kos{K}$-tensor natural isomorphism.
    
    \item
    Let 
    $(\mathfrak{f},\alpha):(\mathfrak{T}^{\pr},\varpi^{\pr})\to (\mathfrak{T},\varpi)$
    be a $1$-cell.
    Then
    the invertible Grothendieck $\kos{K}$-transformation
    $\alpha:\mathfrak{f}^*\widebreve{\varpi}\!^{\pr}\cong\widebreve{\varpi}$
    becomes an invertible $2$-cell
    $\widebreve{\alpha}:
    \mathfrak{f}\widebreve{\varpi}\!^{\pr}
    \cong
    \widebreve{\varpi}\mathfrak{f}(\check{\alpha}_{\kappa})$.
  \end{itemize}
  To each $1$-cell $(\mathfrak{f},\alpha)$,
  we associate the above constraint $2$-cell $\widebreve{\alpha}$.
  This defines a weak $2$-natural transformation
  from the composition $2$-functor
  (\ref{eq GroCat composition 2-functor})
  to the identity $2$-functor on
  $\mathbb{GROCAT}^{\cat{pre}}_*(\kos{K})$.
  It is also invertible in weak sense,
  as every component $1$-cell
  (\ref{eq GRO strong transformation component}) has a quasi-inverse $1$-cell.
  We conclude that the $2$-functors defined in 
  Proposition~\ref{prop GroObjtoGroCat}
  and
  Proposition~\ref{prop GroCattoGroObj}
  are quasi-inverse to each other.
  Thus we established the biequivalence between
  $\mathbb{GROOBJ}^{\cat{pre}}(\kos{K})$
  and
  $\mathbb{GROCAT}^{\cat{pre}}_*(\kos{K})$
  as we claimed.
  This completes the proof of Theorem~\ref{thm GROOBJbiequivGROCAT}.
\qed\end{proof}

\subsection{Torsors of a pre-Grothendieck object}
The content of this subsection is parallel to 
that of \textsection~\ref{subsec Torsors GAL}.
Thus we will omit the details when presenting the results.

Let $\kos{K}=(\CK,\otimes,\kappa)$ be a Grothendieck prekosmos.

\begin{definition}
  Let $\cat{Spec}(\pi)$ be a pre-Grothendieck object in $\kos{K}$.
  A \emph{left $\cat{Spec}(\pi)$-torsor over $\cat{Spec}(\kappa)$}
  is a pair $(\cat{Spec}(p),\lambda_p)$
  of an object $\cat{Spec}(p)$ in $\cat{Aff}(\kos{K})$
  and morphism $\cat{Spec}(\pi)\times \cat{Spec}(p)\xrightarrow{\lambda_p} \cat{Spec}(p)$
  in $\cat{Aff}(\kos{K})$
  such that
  \begin{itemize}
    \item 
    $\slot\otimes p:\CK\to \CK$ is conservative and preserves coreflexive equalizers;
  
    \item
    $\cat{Spec}(\pi)\times \cat{Spec}(p)\xrightarrow{\lambda_p} \cat{Spec}(p)$
    satisfies the left $\pi$-action relations;
  
    \item
    the morphism
    $\cat{Spec}(\pi)\times \cat{Spec}(p)
    \xrightarrow{\tau_p}
    \cat{Spec}(p)\times \cat{Spec}(p)$
    in $\cat{Aff}(\kos{K})$,
    which corresponds to the morphism
    $\tau_p:p\otimes p\xrightarrow{\lambda_p\otimes I_p}\pi\otimes p\otimes p\xrightarrow{I_{\pi}\otimes \pc_p}\pi\otimes p$
    in $\cat{Comm}(\kos{K})$,
    is an isomorphism.
  \end{itemize}
  We define the \emph{division morphism}
  $\cat{Spec}(p)\times \cat{Spec}(p)\xrightarrow{\mon{d}_p} \cat{Spec}(\pi)$
  of $(\cat{Spec}(p),\lambda_p)$ as follows.
  \begin{equation*}
    \mon{d}_p:
    \xymatrix{
      \pi
      \ar[r]^-{\jmath_{\pi}}_-{\cong}
      &\pi\otimes \kappa
      \ar[r]^-{I_{\pi}\otimes u_p}
      &\pi\otimes p
      \ar[r]^-{\tau_p^{-1}}_-{\cong}
      &p\otimes p
    }
  \end{equation*}
\end{definition}

We often omit the left $\cat{Spec}(\pi)$-action
$\lambda_p$
and simply denote a left $\cat{Spec}(\pi)$-torsor
$(\cat{Spec}(p),\lambda_p)$ over $\cat{Spec}(\kappa)$
as $\cat{Spec}(p)$.
A \emph{morphism}
$\cat{Spec}(p)\to \cat{Spec}(p^{\pr})$ of
left $\cat{Spec}(\pi)$-torsors over $\cat{Spec}(\kappa)$
is a morphism in $\cat{Aff}(\kos{K})$
which is compatible with left $\cat{Spec}(\pi)$-actions
$\lambda_p$, $\lambda_{p^{\pr}}$.
We denote
$(\cat{Spec}(\pi)\cat{-Tors})_{\kappa}$
as the category of left $\cat{Spec}(\pi)$-torsors over $\cat{Spec}(\kappa)$.

\begin{remark}
  Let $\cat{Spec}(\pi)$ be a pre-Grothendieck object in $\kos{K}$.
  The pair $(\cat{Spec}(\pi),\cp_{\pi})$ is a left $\cat{Spec}(\pi)$-torsor over $\cat{Spec}(\kappa)$.
  The composition 
  \begin{equation*}
    \varphi_{\pi}:=\tau_{\pi}:
    \pi\otimes \pi
    \xrightarrow{\cp_{\pi}\otimes I_{\pi}}
    \pi\otimes \pi\otimes \pi
    \xrightarrow{I_{\pi}\otimes \pc_{\pi}}
    \pi\otimes \pi
  \end{equation*}
  is an isomorphism
  whose inverse is
  $\varphi_{\pi}^{-1}:
  \pi\otimes \pi
  \xrightarrow{\cp_{\pi}\otimes I_{\pi}}
  \pi\otimes \pi\otimes \pi
  \xrightarrow[\cong]{I_{\pi}\otimes\varsigma_{\pi}\otimes I_{\pi}}
  \pi\otimes \pi\otimes \pi
  \xrightarrow{I_{\pi}\otimes \pc_{\pi}}
  \pi\otimes \pi$.
  The division morphism is given by
  $\mon{d}_{\pi}:
  \pi
  \xrightarrow{\cp_{\pi}}
  \pi\otimes \pi
  \xrightarrow[\cong]{I_{\pi}\otimes \varsigma_{\pi}}
  \pi\otimes \pi$.
\end{remark}

We introduce basic properties of left
$\cat{Spec}(\pi)$-torsors over $\cat{Spec}(\kappa)$.

\begin{lemma}
  Let $\cat{Spec}(\pi)$ be a pre-Grothendieck object in $\kos{K}$
  and let $\cat{Spec}(p)$ be a left $\cat{Spec}(\pi)$-torsor over $\cat{Spec}(\kappa)$.
  The isomorphism
  $\tau_p$
  satisfies the following relations.
  \begin{equation*}
    \vcenter{\hbox{
      \xymatrix@C=30pt{
        p\otimes p\otimes p
        \ar[d]_-{\tau_p\otimes I_p}^-{\cong}
        \ar[r]^-{I_p\otimes \pc_p}
        &p\otimes p
        \ar[d]^-{\tau_p}_-{\cong}
        &p
        \ar[r]^-{\jmath_p}_-{\cong}
        \ar@/_1pc/[drr]_-{\lambda_p}
        &p\otimes \kappa
        \ar[r]^-{I_p\otimes u_p}
        &p\otimes p
        \ar[d]^-{\tau_p}_-{\cong}
        \\
        \pi\otimes p\otimes p
        \ar[r]^-{I_{\pi}\otimes \pc_p}
        &\pi\otimes p
        &\text{ }
        &\text{ }
        &\pi\otimes p
        \\
        p\otimes p
        \ar[d]_-{\tau_p}^-{\cong}
        \ar[r]^-{\lambda_p\otimes I_p}
        &\pi\otimes p\otimes p
        \ar[d]^-{I_{\pi}\otimes \tau_p}_-{\cong}
        &p\otimes p
        \ar[d]_-{\tau_p}^-{\cong}
        \ar@/^1pc/[drr]^-{\cp_{\pi}}
        &\text{ }
        &\text{ }
        \\
        \pi\otimes p
        \ar[r]^-{\cp_{\pi}\otimes I_p}
        &\pi\otimes \pi\otimes p
        &\pi\otimes p
        \ar[r]^-{e_{\pi}\otimes I_p}
        &\kappa\otimes p
        \ar[r]^-{\imath_p^{-1}}_-{\cong}
        &p
      }
    }}
  \end{equation*}
\end{lemma}

\begin{lemma}
  Let $\cat{Spec}(\pi)$ be a pre-Grothendieck object in $\kos{K}$
  and let $\cat{Spec}(p)$ be a left $\cat{Spec}(\pi)$-torsor over $\cat{Spec}(\kappa)$.
  The division morphism
  $\mon{d}_p$
  satisfies the following relations.
  \begin{equation*}
    \vcenter{\hbox{
      \xymatrix@C=30pt{
        \pi
        \ar[d]_-{\cp_{\pi}}
        \ar[r]^-{\mon{d}_p}
        &p\otimes p
        \ar[d]^-{\lambda_p\otimes I_p}
        &\pi
        \ar[r]^-{\mon{d}_p}
        \ar[d]_-{e_{\pi}}
        &p\otimes p
        \ar[d]^-{\pc_{\pi}}
        &\pi\otimes p
        \ar[r]^-{\mon{d}_p\otimes I_p}
        \ar@/_1pc/[dr]_-{\tau_p^{-1}}^-{\cong}
        &p\otimes p\otimes p
        \ar[d]^-{I_p\otimes \pc_p}
        \\
        \pi\otimes \pi
        \ar[r]^-{I_{\pi}\otimes \mon{d}_p}
        &\pi\otimes p\otimes p
        &\kappa
        \ar[r]^-{u_p}
        &p
        &\text{ }
        &p\otimes p
      }
    }}
  \end{equation*}
\end{lemma}

\begin{lemma}
  Let $\cat{Spec}(\pi)$ be a pre-Grothendieck object in $\kos{K}$.
  Every morphism of 
  left $\cat{Spec}(\pi)$-torsors over $\cat{Spec}(\kappa)$
  is an isomorphism, i.e.,
  the category
  $(\cat{Spec}(\pi)\cat{-Tors})_{\kappa}$
  of left $\cat{Spec}(\pi)$-torsors over $\cat{Spec}(\kappa)$
  is a groupoid.
\end{lemma}

Let $\cat{Spec}(\pi)$ be a pre-Grothendieck object in $\kos{K}$.
A \emph{right $\cat{Spec}(\pi)$-torsor over $\cat{Spec}(\kappa)$}
is a pair $(\cat{Spec}(p),\rho_p)$
of an object $\cat{Spec}(p)$ in $\cat{Aff}(\kos{K})$
and a morphism
$\cat{Spec}(p)\times \cat{Spec}(\pi)\xrightarrow{\rho_p} \cat{Spec}(p)$ in $\cat{Aff}(\kos{K})$
such that
\begin{itemize}
  \item 
  $\slot\otimes p:\CK\to \CK$ is conservative and preserves coreflexive equalizers;

  \item 
  $\cat{Spec}(p)\times \cat{Spec}(\pi)\xrightarrow{\rho_p} \cat{Spec}(p)$
  satisfies the right $\cat{Spec}(\pi)$-action relations;

  \item 
  the induced morphism
  $\cat{Spec}(p)\times \cat{Spec}(\pi)
  \to \cat{Spec}(p)\times \cat{Spec}(p)$
  is an isomorphism.
\end{itemize}
We often omit $\rho_p$
and denote a right $\cat{Spec}(\pi)$-torsor
$(\cat{Spec}(p),\rho_p)$ over $\cat{Spec}(\kappa)$
as $\cat{Spec}(p)$.

\begin{remark}
  Let $\cat{Spec}(\pi)$ be a pre-Grothendieck object in $\kos{K}$.
  Suppose $\cat{Spec}(p)$ is a left $\cat{Spec}(\pi)$-torsor
  over $\cat{Spec}(\kappa)$ with left $\cat{Spec}(\pi)$-action
  $\lambda_p$.
  Then $\cat{Spec}(p)$ equipped with
  a right $\cat{Spec}(\pi)$-action
  \begin{equation*}
    \cat{Spec}(p)\times \cat{Spec}(\pi)\xrightarrow{\rho_p}
    \cat{Spec}(p),
    \quad
    \xymatrix{
      \rho_p:
      p
      \ar[r]^-{\lambda_p}
      &\pi\otimes p
      \ar[r]^-{s_{\pi,p}}_-{\cong}
      &p\otimes \pi
      \ar[r]^-{I_p\otimes \varsigma_{\pi}}_-{\cong}
      &p\otimes \pi
    }
  \end{equation*}
  becomes a right $\cat{Spec}(\pi)$-torsor over $\cat{Spec}(\kappa)$.
  Conversely, if $\cat{Spec}(p)$ is a right $\cat{Spec}(\pi)$-torsor over $\cat{Spec}(\kappa)$
  with right $\cat{Spec}(\pi)$-action
  $\rho_p$,
  then $\cat{Spec}(p)$ equipped with
  \begin{equation*}
    \cat{Spec}(\pi)\times \cat{Spec}(p)
    \xrightarrow{\lambda_p}
    \cat{Spec}(p),
    \quad
    \lambda_p:
    \xymatrix{
      p
      \ar[r]^-{\rho_p}
      &p\otimes \pi
      \ar[r]^-{s_{p,\pi}}_-{\cong}
      &\pi\otimes p
      \ar[r]^-{\varsigma_{\pi}\otimes I_p}_-{\cong}
      &\pi\otimes p
    }
  \end{equation*}
  becomes a left $\cat{Spec}(\pi)$-torsor over $\cat{Spec}(\kappa)$.
\end{remark}

Let $\cat{Spec}(\pi^{\pr})$ be another pre-Grothendieck object in $\kos{K}$.
By a \emph{$(\cat{Spec}(\pi),\cat{Spec}(\pi^{\pr}))$-bitorsor over $\cat{Spec}(\kappa)$}
we mean a triple
$(\cat{Spec}(p),\lambda_p,\rho^{\pr}_p)$
where
$(\cat{Spec}(p),\lambda_p)$
is a left $\cat{Spec}(\pi)$-torsor over $\cat{Spec}(\kappa)$,
$(\cat{Spec}(p),\rho^{\pr}_p)$
is a right $\cat{Spec}(\pi^{\pr})$-torsor over $\cat{Spec}(\kappa)$
and 
the left $\cat{Spec}(\pi)$-action
$\lambda_p$
is compatible with
the right $\cat{Spec}(\pi^{\pr})$-action
$\rho^{\pr}_p$.
We often omit
$\lambda_p$
as well as
$\rho^{\pr}_p$
and simply denote a
$(\cat{Spec}(\pi),\cat{Spec}(\pi^{\pr}))$-bitorsor
$(\cat{Spec}(p),\lambda_p,\rho^{\pr}_p)$
over $\cat{Spec}(\kappa)$
as $\cat{Spec}(p)$.
The
\emph{opposite
$(\cat{Spec}(\pi^{\pr}),\cat{Spec}(\pi))$-bitorsor
of $\cat{Spec}(p)$ over $\cat{Spec}(\kappa)$}
is a triple 
$\cat{Spec}(\mathring{p})
=(\cat{Spec}(p),\lambda^{\pr}_p,\rho_p)$
where
\begin{equation*}
  \begin{aligned}
    \lambda^{\pr}_p
    &:
    \xymatrix{
      p
      \ar[r]^-{\rho^{\pr}_p}
      &p\otimes \pi^{\pr}
      \ar[r]^-{s_{p,\pi^{\pr}}}_-{\cong}
      &\pi^{\pr}\otimes p
      \ar[r]^-{\varsigma_{\pi^{\pr}\otimes }I_p}_-{\cong}
      &\pi^{\pr}\otimes p
      ,
    }
    \\
    \rho_p
    &:
    \xymatrix{
      p
      \ar[r]^-{\lambda_p}
      &\pi\otimes p
      \ar[r]^-{s_{\pi,p}}_-{\cong}
      &p\otimes \pi
      \ar[r]^-{I_p\otimes \varsigma_{\pi}}_-{\cong}
      &p\otimes \pi
      .
    }
  \end{aligned}
\end{equation*}

\begin{lemma}
  \label{lem TorsorsGro pip definition}
  Let $\cat{Spec}(\pi)$ be a pre-Grothendieck object in $\kos{K}$
  and let $\cat{Spec}(p)$ be a left $\cat{Spec}(\pi)$-torsor over $\cat{Spec}(\kappa)$.
  We define an object $\pi^p$ in $\cat{Comm}(\kos{K})$
  as the following coreflexive equalizer.
  \begin{equation*}
    \vcenter{\hbox{
      \xymatrix@C=50pt{
        \pi^p
        \ar@{^{(}->}[r]^-{\textit{eq}_p}
        &p\otimes p
        \ar@<0.5ex>[rr]^-{u_p\otimes I_{p\otimes p}}
        \ar@<-0.5ex>[rr]_-{(\pc_p\otimes I_{p\otimes p})\circ (I_p\otimes \mon{d}_p\otimes I_p)\circ (I_p\otimes \lambda_p)}
        &\text{ }
        &p\otimes p\otimes p
        \ar@<-1ex>@/_2pc/[ll]|-{\pc_p\otimes I_p}
      }
    }}
  \end{equation*}
  \begin{enumerate}
    \item 
    We have a unique morphism
    $\rho^p_p:p\to p\otimes \pi^p$
    in $\cat{Comm}(\kos{K})$
    which satisfies the following relation.
    \begin{equation*}
      \vcenter{\hbox{
        \xymatrix@C=40pt{
          p
          \ar@{.>}[d]^-{\rho^p_p}_-{\exists!}
          \ar[r]^-{\lambda_p}
          &\pi\otimes p
          \ar[d]^-{\mon{d}_p\otimes I_p}
          \\
          p\otimes \pi^p
          \ar@{^{(}->}[r]^-{I_p\otimes \textit{eq}_p}
          &p\otimes p\otimes p
        }
      }}
    \end{equation*}
    Moreover,
    $\rho^p_p$
    is compatible with
    $\lambda_p$.

    \item 
    We have an isomorphism
    $\Phi:\pi^p\otimes p \xrightarrow{\cong}p\otimes p$
    in $\cat{Comm}(\kos{K})$
    defined by
    \begin{equation*}
      \Phi:
      \xymatrix@C=30pt{
        \pi^p\otimes p
        \ar@{^{(}->}[r]^-{\textit{eq}_p\otimes I_p}
        &p\otimes p\otimes p
        \ar[r]^-{s_{p,p}\otimes I_p}_-{\cong}
        &p\otimes p\otimes p
        \ar[r]^-{I_p\otimes \pc_p}
        &p\otimes p
      }
    \end{equation*}
    whose inverse is 
    \begin{equation*}
      \Phi^{-1}:
      \xymatrix@C=30pt{
        p\otimes p
        \ar[r]^-{\rho_p^p\otimes I_p}
        &p\otimes \pi^p\otimes p
        \ar[r]^-{s_{p,\pi^p}\otimes I_p}_-{\cong}
        &\pi^p\otimes p\otimes p
        \ar[r]^-{I_{\pi^p}\otimes \pc_p}
        &\pi^p\otimes p.
      }
    \end{equation*}
    Moreover, the isomorphism
    $\Phi$
    satisfies the following relation.
    \begin{equation*}
      \vcenter{\hbox{
        \xymatrix@C=30pt{
          \pi^p\otimes \pi^p\otimes p
          \ar[d]_-{I_{\pi^p}\otimes \Phi}^-{\cong}
          \ar[rrr]^-{\pc_{\pi^p}\otimes I_p}
          &\text{ }
          &\text{ }
          &\pi^p\otimes p
          \ar[d]^-{\Phi}_-{\cong}
          \\
          \pi^p\otimes p\otimes p
          \ar@{^{(}->}[r]^-{\textit{eq}_p\otimes I_{p\otimes p}}
          &p\otimes p\otimes p\otimes p
          \ar[r]^-{s_{p,p}\otimes I_{p\otimes p}}_-{\cong}
          &p\otimes p\otimes p\otimes p
          \ar[r]^-{\pc_{p\otimes p}}
          &p\otimes p
        }
      }}
    \end{equation*}
  \end{enumerate}
\end{lemma}

\begin{proposition}
  \label{prop TorsorsGro lefttors become bitors}
  Let $\cat{Spec}(\pi)$ be a pre-Grothendieck object in $\kos{K}$
  and let $\cat{Spec}(p)$ be a left $\cat{Spec}(\pi)$-torsor over $\cat{Spec}(\kappa)$.
  Recall the object $\pi^p$ in $\cat{Comm}(\kos{K})$
  and the morphism $\rho^p_p:p\to p\otimes \pi^p$
  in $\cat{Comm}(\kos{K})$
  introduced in Lemma~\ref{lem TorsorsGro pip definition}.
  \begin{enumerate}
    \item 
    We have a pre-Grothendieck object $\cat{Spec}(\pi^p)$ in $\kos{K}$,
    whose product, unit morphisms
    $\cp_{\pi^p}$, $e_{\pi^p}$
    are the unique morphisms
    in $\cat{Comm}(\kos{K})$
    satisfying the following relations.
    \begin{equation*}
      \vcenter{\hbox{
        \xymatrix@C=40pt{
          \pi^p
          \ar@{.>}[d]^-{\cp_{\pi^p}}_-{\exists!}
          \ar@{^{(}->}[r]^-{\textit{eq}_p}
          &p\otimes p
          \ar[d]^-{I_p\otimes \rho_p^p}
          &\pi^p
          \ar@{^{(}->}[r]^-{\textit{eq}_p}
          \ar@{.>}[d]^-{e_{\pi^p}}_-{\exists!}
          &p\otimes p
          \ar[d]^-{\pc_p}
          \\
          \pi^p\otimes \pi^p
          \ar@{^{(}->}[r]^-{\textit{eq}_p\otimes I_{\pi^p}}
          &p\otimes p\otimes \pi^p
          &\kappa
          \ar[r]^-{u_p}
          &p
        }
      }}
    \end{equation*}

    \item
    The triple
    $(\cat{Spec}(p),\lambda_p,\rho^p_p)$
    is a $(\cat{Spec}(\pi),\cat{Spec}(\pi^p))$-bitorsor
    over $\cat{Spec}(\kappa)$.
  \end{enumerate}
\end{proposition}

\begin{lemma}
  \label{lem TorsorsGro antipode formula}
  Let $\cat{Spec}(\pi)$
  be a pre-Grothendieck object in $\kos{K}$
  and let $\cat{Spec}(p)$
  be a left $\cat{Spec}(\pi)$-torsor over $\cat{Spec}(\kappa)$.
  \begin{enumerate}
    \item 
    The isomorphism
    $\tau_p:p\otimes p\xrightarrow{\cong}\pi\otimes p$
    satisfies the following relation.
    \begin{equation*}
      \vcenter{\hbox{
        \xymatrix@C=35pt{
          p\otimes p\otimes p
          \ar[rrr]^-{\tau_p\otimes I_p}_-{\cong}
          \ar[d]_-{\lambda_p\otimes I_{p\otimes p}}
          &\text{ }
          &\text{ }
          &\pi\otimes p\otimes p
          \ar[dd]^-{I_{\pi}\otimes \tau_p}_-{\cong}
          \\
          \pi\otimes p\otimes p\otimes p
          \ar[d]_-{I_{\pi}\otimes s_{p,p}\otimes I_p}^-{\cong}
          \ar@/^1pc/[dr]|-{(\otimes p)_{\pi,p}}
          &\text{ }
          &\text{ }
          &\text{ }
          \\
          \pi\otimes p\otimes p\otimes p
          \ar[r]^-{I_{\pi\otimes p}\otimes \pc_p}
          &\pi\otimes p\otimes p
          \ar[r]^-{I_{\pi}\otimes \tau_p}_-{\cong}
          &\pi\otimes \pi\otimes p
          \ar[r]^-{\varphi_{\pi}\otimes I_p}_-{\cong}
          &\pi\otimes \pi\otimes p
        }
      }}
    \end{equation*}
  
    \item 
    The division morphism $\mon{d}_p:\pi\to p\otimes p$
    satisfies the following relation.
    \begin{equation*}
      \vcenter{\hbox{
        \xymatrix@C=40pt{
          \pi
          \ar[d]_-{\varsigma_{\pi}}^-{\cong}
          \ar[r]^-{\mon{d}_p}
          &p\otimes p
          \ar[d]^-{s_{p,p}}_-{\cong}
          \\
          \pi
          \ar[r]^-{\mon{d}_p}
          &p\otimes p
        }
      }}
    \end{equation*}

    \item 
    The induced right $\cat{Spec}(\pi)$-action
    $\xymatrix{
      \rho_p:
      p
      \ar[r]^-{\lambda_p}
      &\pi\otimes p
      \ar[r]^-{s_{\pi,p}}_-{\cong}
      &p\otimes \pi
      \ar[r]^-{I_p\otimes \varsigma_{\pi}}_-{\cong}
      &p\otimes \pi
    }$
    satisfies the following relation.
    \begin{equation*}
      \vcenter{\hbox{
        \xymatrix@C=30pt{
          p\otimes p
          \ar[d]_-{\pc_p}
          \ar[r]^-{\rho_p\otimes I_p}
          &p\otimes \pi\otimes p
          \ar[r]^-{I_p\otimes \tau_p^{-1}}_-{\cong}
          &p\otimes p\otimes p
          \ar[d]^-{\pc_p\otimes I_p}
          \\
          p
          \ar[r]^-{\imath_p}_-{\cong}
          &\kappa\otimes p
          \ar[r]^-{u_p\otimes I_p}
          &p\otimes p
        }
      }}
    \end{equation*}
  \end{enumerate}
\end{lemma}

\subsubsection{Pre-fiber functor twisted by a torsor}

\begin{lemma}\label{lem Gro twistbyatorsor}
  Let $\cat{Spec}(\pi)$ be a pre-Grothendieck object in $\kos{K}$
  and let $\cat{Spec}(p)$ be a left 
  $\cat{Spec}(\pi)$-torsor over $\cat{Spec}(\kappa)$.
  We can twist the forgetful
  strong $\kos{K}$-tensor functor
  $(\omega_{\pi}^*,\what{\omega}_{\pi}^*):
  (\kos{R\!e\!p}(\pi),\kos{t}^*_{\pi})
  \to
  (\kos{K},\id_{\kos{K}})$
  by $\cat{Spec}(p)$ and obtain the strong $\kos{K}$-tensor functor
  \begin{equation*}
    \vcenter{\hbox{
      \xymatrix@R=30pt@C=40pt{
        \text{ }
        &\kos{R\!e\!p}(\pi)
        \ar[d]^-{\omega^{p*}}
        \\
        \kos{K}
        \ar@/^1pc/[ur]^-{\kos{t}_{\pi}^*}
        \ar[r]_-{\id_{\kos{K}}}
        &\kos{K}
        \xtwocell[l]{}<>{<3>{\what{\omega}^{p*}\quad}}
      }
    }}
    \qquad\quad
    (\omega^{p*},\what{\omega}^{p*}):
    (\kos{R\!e\!p}(\pi),\kos{t}^*_{\pi})
    \to
    (\kos{K},\id_{\kos{K}})
    .
  \end{equation*}
  \begin{itemize}
    \item 
    The underlying functor $\omega^{p*}$ sends 
    each object $X=(x,\rho_x)$ in $\cat{Rep}(\pi)$
    to the following coreflexive equalizer in $\CK$.
    \begin{equation}\label{eq Gro twistbyatorsor}
      \vcenter{\hbox{
        \xymatrix@C=40pt{
          \omega^{p*}(X)
          \ar@{^{(}->}[r]^-{\xi^p_X}
          &x\otimes p
          \ar@<0.5ex>[r]^-{\rho_x\otimes I_p}
          \ar@<-0.5ex>[r]_-{I_x\otimes \lambda_p}
          &x\otimes \pi\otimes p
          \ar@<-1ex>@/_2pc/[l]|-{I_x\otimes e_{\pi}\otimes I_p}
        }
      }}
    \end{equation}
    For every object $z$ in $\CK$,
    the functor $z\otimes \slot:\CK\to \CK$
    preserves the coreflexive equalizer (\ref{eq Gro twistbyatorsor}).
    
    \item 
    The strong symmetric monoidal coherence isomorphisms of
    $\omega^{p*}$ are unique morphisms satisfying the relations below.
    Let $Y=(y,\rho_y)$
    be another object in $\cat{Rep}(\pi)$.
    \begin{equation*}
      \vcenter{\hbox{
        \xymatrix@C=40pt{
          \omega^{p*}(X)\otimes \omega^{p*}(Y)
          \ar@{^{(}->}[d]_-{\xi^p_X\otimes \xi^p_Y}
          \ar@{.>}[r]^-{\text{ }\omega^{p*}_{X,Y}}_-{\cong}
          &\omega^{p*}(X\tensor\!_{\pi}Y)
          \ar@{^{(}->}[d]^-{\xi^p_{X\tensor\!_{\pi}Y}}
          \\
          (x\otimes p)\otimes (y\otimes p)
          \ar[r]^-{(\otimes p)_{x,y}}
          &(x\otimes y)\otimes p
        }
      }}
      \qquad
      \vcenter{\hbox{
        \xymatrix@C=40pt{
          \kappa
          \ar[d]_-{u_p}
          \ar@{.>}[r]^-{\exists!\text{ }\omega^{p*}_{\unit\!_{\pi}}}_-{\cong}
          &\omega^{p*}(\unit\!_{\pi})
          \ar@{^{(}->}[d]^-{\xi^p_{\unit\!_{\pi}}}
          \\
          p
          \ar[r]^-{\imath_p}_-{\cong}
          &\kappa \otimes p
        }
      }}
    \end{equation*}
  
    \item
    The monoidal natural isomorphism
    $\what{\omega}^{p*}:\id_{\kos{K}}\cong \omega^{p*}\kos{t}_{\pi}^*$
    is described below.
    Let $z$ be an object in $\CK$.
    \begin{equation*}
      \vcenter{\hbox{
        \xymatrix@C=50pt{
          z
          \ar@{.>}[d]_-{\exists!\text{ }\what{\omega}^{p*}_z}^-{\cong}
          \ar@/^1pc/[dr]^-{I_z\otimes u_p}
          &\text{ }
          \\
          \omega^{p*}\kos{t}_{\pi}^*(z)
          \ar@{^{(}->}[r]^-{\xi^p_{\kos{t}_{\pi}^*(z)}}
          &z\otimes p
        }
      }}
    \end{equation*}
    If we denote
    $z\otimes X
    =z\acts X
    =\kos{t}_{\pi}^*(z)\tensor\!_{\pi}X$
    as the action of $z$ on $X$,
    then the associated $\kos{K}$-equivariance 
    $\cevar{\omega}^{p*}$ of
    $(\omega^{p*},\what{\omega}^{p*})$
    is described below.
    \begin{equation*}
      \vcenter{\hbox{
        \xymatrix@C=40pt{
          z\otimes \omega^{p*}(X)
          \ar@{^{(}->}[d]_-{I_z\otimes \xi^p_X}
          \ar@{.>}[r]^-{\exists!\text{ }\cevar{\omega}^{p*}_{z,X}}_-{\cong}
          &\omega^{p*}(z\otimes X)
          \ar@{^{(}->}[d]^-{\xi^p_{z\otimes X}}
          \\
          z\otimes (x\otimes p)
          \ar[r]^-{a_{z,x,p}}_-{\cong}
          &(z\otimes x)\otimes p
        }
      }}
    \end{equation*}
  \end{itemize}
\end{lemma}

\begin{proposition} \label{prop Gro twistbyatorsor}
  Let $\cat{Spec}(\pi)$, $\cat{Spec}(\pi^{\pr})$
  be pre-Grothendieck objects in $\kos{K}$
  and $\cat{Spec}(p)$
  be a $(\cat{Spec}(\pi),\cat{Spec}(\pi^{\pr}))$-bitorsor over $\cat{Spec}(\kappa)$.
  \begin{enumerate}
    \item 
    We have a pre-fiber functor
    \begin{equation*}
      \varpi^p=(\omega^p,\what{\omega}^p)
      :\mathfrak{K}\to \mathfrak{Rep}(\pi)      
    \end{equation*}
    which is obtained by twisting the pre-fiber functor
    $\varpi_{\pi}=(\omega_{\pi},\what{\omega}_{\pi})
    :\mathfrak{K}\to\mathfrak{Rep}(\pi)$.
    The left adjoint strong $\kos{K}$-tensor functor
    $(\omega^{p*},\what{\omega}^{p*})$
    of $\varpi^p$
    is described in Lemma~\ref{lem Gro twistbyatorsor}.

    \item
    The twisted fiber functor
    $\varpi^p:\mathfrak{K}\to\mathfrak{Rep}(\pi)$
    factors through as an equivalence of Grothendieck $\kos{K}$-prekosmoi
    $\widebreve{\varpi}\!^p:
    \mathfrak{Rep}(\pi^{\pr})
    \xrightarrow{\simeq}
    \mathfrak{Rep}(\pi)$.
    \begin{equation*}
      \vcenter{\hbox{
        \xymatrix@C=30pt{
          \mathfrak{Rep}(\pi)
          &\mathfrak{Rep}(\pi^{\pr})
          \ar[l]_-{\widebreve{\varpi}\!^p}^-{\simeq}
          \\
          \text{ }
          &\mathfrak{K}
          \ar[u]_-{\varpi_{\pi^{\pr}}}
          \ar@/^1pc/[ul]^-{\varpi^p}
        }
      }}
    \end{equation*}
  \end{enumerate}
\end{proposition}

\begin{corollary} \label{cor Gro twistbyatorsor}
  Let $\cat{Spec}(\pi)$ be a pre-Grothendieck object in $\kos{K}$
  and let $\cat{Spec}(p)$ be a left $\cat{Spec}(\pi)$-torsor over $\cat{Spec}(\kappa)$.
  Recall Proposition~\ref{prop TorsorsGro lefttors become bitors}.
  We denote $\cat{Spec}(\pi^p)$ as the pre-Grothendieck object in $\kos{K}$
  obtained by twisting $\cat{Spec}(\pi)$ by $\cat{Spec}(p)$,
  and $\cat{Spec}(p)$ becomes a $(\cat{Spec}(\pi),\cat{Spec}(\pi^p))$-bitorsor over $\cat{Spec}(\kappa)$.
  \begin{enumerate}
    \item 
    We have a pre-fiber functor
    $\varpi^p:\mathfrak{K}\to \mathfrak{Rep}(\pi)$
    which is obtained by twisting the pre-fiber functor
    $\varpi_{\pi}:\mathfrak{K}\to\mathfrak{Rep}(\pi)$.
    See Lemma~\ref{lem Gro twistbyatorsor}
    and Proposition~\ref{prop Gro twistbyatorsor}.

    \item
    The twisted pre-fiber functor 
    $\varpi^p:\mathfrak{K}\to \mathfrak{Rep}(\pi)$
    factors through as an equivalence of Grothendieck $\kos{K}$-prekosmoi
    $\widebreve{\varpi}\!^p:\mathfrak{Rep}(\pi^p)\xrightarrow[]{\simeq}\mathfrak{Rep}(\pi)$.
    \begin{equation*}
      \vcenter{\hbox{
        \xymatrix@C=40pt{
          \mathfrak{Rep}(\pi)
          &\mathfrak{Rep}(\pi^p)
          \ar[l]_-{\widebreve{\varpi}\!^p}^-{\simeq}
          \\
          \text{ }
          &\mathfrak{K}
          \ar@/^0.5pc/[ul]^-{\varpi^p}
          \ar[u]_-{\varpi_{\pi^{p}}}
        }
      }}
    \end{equation*}

    \item 
    The presheaf of groups of
    invertible Grothendieck $\kos{K}$-transformations
    from $\varpi^p$ to $\varpi^p$
    is represented by the twisted pre-Grothendieck object $\cat{Spec}(\pi^p)$ in $\kos{K}$.
    \begin{equation*}
      \xymatrix{
        \Hom_{\cat{Aff}(\kos{K})}(\slot,\cat{Spec}(\pi^p))
        \ar@2{->}[r]^-{\cong}
        &\underline{\Aut}_{\mathbb{GRO}^{\cat{pre}}_{\kos{K}}}(\varpi^p)
        :\cat{Aff}(\kos{K})^{\op}\to\cat{Grp}
      }
    \end{equation*}

    \item 
    The presheaf of
    invertible Grothendieck $\kos{K}$-transformations
    from $\varpi^p$ to $\varpi_{\pi}$
    is represented by $\cat{Spec}(p)$.
    \begin{equation*}
      \xymatrix{
        \Hom_{\cat{Aff}(\kos{K})}(\slot,\cat{Spec}(p))
        \ar@2{->}[r]^-{\cong}
        &\underline{\Isom}_{\mathbb{GRO}^{\cat{pre}}_{\kos{K}}}(\varpi^p,\varpi_{\pi})
        :\cat{Aff}(\kos{K})^{\op}\to\cat{Set}
      }
    \end{equation*}
    The component of the universal element
    $\xi^p$
    of the representation 
    at each object $X=(x,\rho_x)$ in $\cat{Rep}(\pi)$
    is the coreflexive equalizer
    $\xymatrix@C=15pt{
      \xi^p_X:\omega^{p*}(X)\ar@{^{(}->}[r]&x\otimes p
    }$
    described in (\ref{eq Gro twistbyatorsor}).
  \end{enumerate}
\end{corollary}

\begin{lemma} \label{lem Gro twistbyatorsor functorial}
  Let $\cat{Spec}(\pi)$ be a pre-Grothendieck object in $\kos{K}$
  and let $f:\cat{Spec}(p)\xrightarrow{\cong}\cat{Spec}(q)$
  be an isomorphism of left $\cat{Spec}(\pi)$-torsors over $\cat{Spec}(\kappa)$.
  We have an invertible Grothendieck $\kos{K}$-transformation
  \begin{equation*}
    \vcenter{\hbox{
      \xymatrix{
        \mathfrak{Rep}(\pi)
        \\
        \mathfrak{K}
        \ar@/^1.2pc/[u]^-{\varpi^q}
        \ar@/_1.2pc/[u]_-{\varpi^p}
        \xtwocell[u]{}<>{<0>{\vartheta^f}}
      }
    }}
    \qquad\quad
    \vartheta^f:
    \xymatrix@C=15pt{
      \varpi^q
      \ar@2{->}[r]^-{\cong}
      &\varpi^p
      :\mathfrak{K}\to \mathfrak{Rep}(\pi)
    }
  \end{equation*}
  between twisted pre-fiber functors
  in the opposite direction.
  The component of
  $\vartheta^f:\omega^{q*}\cong \omega^{p*}$
  at each object $X=(x,\rho_x)$ in $\cat{Rep}(\pi)$ is
  the unique isomorphism satisfying the relation below.
  \begin{equation*}
    \vcenter{\hbox{
      \xymatrix@C=40pt{
        \omega^{q*}(X)
        \ar@{.>}[d]_-{\exists!\text{ }\vartheta^f_X}^-{\cong}
        \ar@{^{(}->}[r]^-{\xi^q_X}
        &x\otimes q
        \ar[d]^-{I_x\otimes f}_-{\cong}
        \\
        \omega^{p*}(X)
        \ar@{^{(}->}[r]^-{\xi^p_X}
        &x\otimes p
      }
    }}
  \end{equation*}
\end{lemma}

\subsubsection{Torsors and pre-fiber functors}

\begin{theorem}
  \label{thm Gro TorsFib adjequiv}
  Let $\mathfrak{T}=(\kos{T},\kos{t})$ be a pre-Grothendieck $\kos{K}$-category
  and let $\varpi=(\omega,\what{\omega}):\mathfrak{K}\to\mathfrak{T}$
  be a pre-fiber functor for $\mathfrak{T}$.
  Recall Theorem~\ref{thm preGroKCat mainThm1}.
  We denote $\cat{Spec}(\pi)$,
  $\pi:=\omega^*\omega_!(\kappa)$
  as the pre-Grothendieck object in $\kos{K}$
  which represents 
  $\underline{\Aut}_{\mathbb{GRO}^{\cat{pre}}_{\kos{K}}}(\varpi)$,
  and the given pre-fiber functor
  $\varpi$ factors through as equivalence of Grothendieck $\kos{K}$-prekosmoi
  $\widebreve{\varpi}:\mathfrak{Rep}(\pi)\xrightarrow{\simeq}\mathfrak{T}$.
  \begin{equation*}
    \vcenter{\hbox{
      \xymatrix@C=40pt{
        \mathfrak{T}
        &\mathfrak{Rep}(\pi)
        \ar[l]_-{\widebreve{\varpi}}^-{\simeq}
        \\
        \text{ }
        &\mathfrak{K}
        \ar[u]_-{\varpi_{\pi}}
        \ar@/^0.5pc/[ul]^-{\varpi}
      }
    }}
  \end{equation*}
  Then we have an adjoint equivalence of groupoids
  \begin{equation*}
    \xymatrix{
      \cat{Fib}(\mathfrak{T})^{\cat{pre}}
      \ar@<0.5ex>[r]^-{\simeq}
      &(\cat{Spec}(\pi)\cat{-Tors})_{\kappa}^{\op}
      \ar@<0.5ex>[l]^-{\simeq}
    }
  \end{equation*}
  between the groupoid of pre-fiber functors for $\mathfrak{T}$
  and the opposite groupoid of left
  $\cat{Spec}(\pi)$-torsors over $\cat{Spec}(\kappa)$.
  \begin{itemize}
    \item 
    Each pre-fiber functor $\varpi^{\pr}$
    for $\mathfrak{T}$
    is sent to the left $\cat{Spec}(\pi)$-torsor
    $\cat{Spec}(\omega^{\pr*}\omega_*(\kappa))$ over $\cat{Spec}(\kappa)$,
    which represents the presheaf
    $\underline{\Isom}_{\mathbb{GRO}^{\cat{pre}}_{\kos{K}}}(\varpi^{\pr},\varpi)$
    of invertible Grothendieck $\kos{K}$-transformations
    from $\varpi^{\pr}$ to $\varpi$.
    \begin{equation*}
      \vcenter{\hbox{
        \xymatrix{
          \Hom_{\cat{Aff}(\kos{K})}\bigl(\slot,\cat{Spec}(\omega^{\pr*}\omega_*(\kappa))\bigr)
          \ar@2{->}[r]^-{\cong}
          &\underline{\Isom}_{\mathbb{GRO}^{\cat{pre}}_{\kos{K}}}(\varpi^{\pr},\varpi)
          :\cat{Aff}(\kos{K})^{\op}\to\cat{Set}
        }
      }}
    \end{equation*}
    The left $\cat{Spec}(\pi)$-action of $\cat{Spec}(\omega^{\pr*}\omega_*(\kappa))$ is
    \begin{equation*}
      \xymatrix@C=35pt{
        \lambda_{\omega^{\pr*}\omega_*(\kappa)}\!:\!
        \omega^{\pr*}\omega_*(\kappa)
        \ar[r]^-{\omega^{\pr*}(\eta_{\omega_*(\kappa)})}
        &\omega^{\pr*}\omega_*\omega^*\omega_*(\kappa)
        \ar[r]^-{(\tahar{\omega^{\pr*}\omega_*}_{\omega^*\omega_*(\kappa)})^{-1}}_-{\cong}
        &\omega^*\omega_*(\kappa)\!\otimes\! \omega^{\pr*}\omega_*(\kappa)
        .
      }
    \end{equation*}
  
    \item 
    Each left $\cat{Spec}(\pi)$-torsor $\cat{Spec}(p)$ over $\cat{Spec}(\kappa)$
    is sent to the pre-fiber functor
    \begin{equation*}
      \widebreve{\varpi}\!\varpi^p:
      \xymatrix{
        \mathfrak{K}
        \ar[r]^-{\varpi^p}
        &\mathfrak{Rep}(\pi)
        \ar[r]^-{\widebreve{\varpi}}_-{\simeq}
        &\mathfrak{T}
      }
    \end{equation*}
    where $\varpi^p:\mathfrak{K}\to \mathfrak{Rep}(\pi)$
    is the pre-fiber functor
    twisted by the given torsor $\cat{Spec}(p)$
    which we introduced in Corollary~\ref{cor Gro twistbyatorsor}.
  \end{itemize}
\end{theorem}

\end{document}